\documentclass[10pt]{article}
\usepackage{amsmath}
\usepackage{amssymb}
\usepackage{amscd}
\usepackage{amsthm}
\usepackage{enumitem}
\usepackage{stmaryrd} 
\usepackage{hyperref}
\def\ver{arxiv-jan15.tex}
\hypersetup{
    colorlinks=true,
    linkcolor=blue,
    filecolor=magenta,      
    urlcolor=blue,
    pdftitle={Overleaf Example},
    pdfpagemode=FullScreen,
    }
\theoremstyle{plain}
\newtheorem{theorem}{Theorem}[section]
\newtheorem*{theoremone-ref}{Theorem 1.1 {\rm (Reformulation)}}
\newtheorem*{theoremtwo-ref}{Theorem 1.2 {\rm (Reformulation)}}
\newtheorem*{theoremone-dual}{Theorem 1.1$'$}
\newtheorem*{theoremone-rot}{Theorem 1.1$''$}
\newtheorem*{theoremtwo-dual}{Theorem 1.2$'$}
\newtheorem*{theoremtwo-rot}{Theorem 1.2$''$}
\newtheorem*{theoremthree-ref}{Theorem 1.3 {\rm (Reformulation)}}
\newtheorem*{theoremthree-prime}{Theorem 1.3$'$}
\newtheorem*{theoremone-again}{Theorem 1.1 {\rm (Again reformulated)}}
\newtheorem*{theoremfour-recall}{Theorem 1.4}
\newtheorem{proposition}[theorem]{Proposition}
\newtheorem{lemma}[theorem]{Lemma}
\newtheorem{corollary}[theorem]{Corollary}
\newtheorem*{corollarythree-ref}{Corollary}
\newtheorem*{corollaryfive}{Corollary}
\newtheorem*{corollaryseven}{Corollary}
\theoremstyle{definition}
\newtheorem*{definition}{Definition}
\newtheorem*{definitions}{Definitions}
\newtheorem{remark}[theorem]{Remark}
\newtheorem{remarks}[theorem]{Remarks}
\newtheorem*{remarkfour}{Remark}
\newtheorem{warning}[theorem]{Warning}
\newtheorem*{warning-four-three-cont}{Warning 4.3 {\rm (continued)}}
\newtheorem{example}[theorem]{Example}
\renewcommand{\thetheorem}{\arabic{section}.\arabic{theorem}}
\renewcommand{\thesection}{\arabic{section}}
\renewcommand{\thesubsection}{(\arabic{section}.\arabic{subsection})}

\input prepictex   \input pictex    \input postpictex
\font\kl=cmr5  
   \font\rmk=cmr8      \font\ttk=cmtt8

\font\Gross=cmbx10 scaled\magstep3

\font\Grossit=cmmib10 scaled \magstep3

\def\bs{\backslash}

\def\ind{\operatorname{ind}}

\def\mod{\operatorname{mod}}
\def\Mod{\operatorname{Mod}}
\def\MOD{\operatorname{MOD}}

\def\Hom{\operatorname{Hom}}
\def\End{\operatorname{End}}
\def\Ext{\operatorname{Ext}}

\def\rad{\operatorname{rad}}

\def\Ker{\operatorname{Ker}}
\def\Eig{\operatorname{Eig}}
\def\Cok{\operatorname{Cok}}
\def\soc{\operatorname{soc}}

\def\Im{\operatorname{Im}}

\def\bdim{\operatorname{\mathbf{dim}}}
\def\pr{\operatorname{\mathbf{pr}}}
\def\bpar{\operatorname{\mathbf{par}}}

\def\under{\operatorname{\backslash}}
\def\D{\operatorname{D}}
\def\:{\,{\rm :}\;}
\def\ss{\ssize}
\def\sss{\sssize}
\def\arr#1#2{\arrow <1.5mm> [0.25,0.75] from #1 to #2}

\def\s{\hfill \square} 
\def\GL{\operatorname{GL}}

\def\Mimo{\operatorname{Mimo}}
\def\uwbb{\operatorname{\mathbf{uwb}}}

\def\stackrel#1#2{\mathrel{\mathop{#2}\limits^{#1}}}

\def\Cal{\mathcal}
\def\ssize{\scriptstyle}
\def\sssize{\scriptscriptstyle}

\begin{document}
{\tt \ver}
\sloppypar

\medskip
\centerline{\bf Invariant Subspaces of Nilpotent Operators.}

\medskip
\centerline{\bf Level, Mean, and Colevel:
        The Triangle $\mathbb T(n)$.}

\bigskip\smallskip
\centerline{Claus Michael Ringel and Markus Schmidmeier}

\bigskip\bigskip
{\narrower 
Abstract. We consider the category $\Cal S(n)$ of pairs $X = (U,V)$, where $V$ is
a finite-dimensional vector space with a nilpotent 
operator $T$ with $T^n = 0$, and $U$ is a subspace
of $V$ such that $T(U) \subseteq U.$
For any vector space $V$, let $|V|$ denote its dimension
(or length). 
Note that $\Cal S(n)$ is just the category
of Gorenstein-projective $T_2(\Lambda)$-modules, where $\Lambda = k[T]/\langle T^n\rangle$
and $T_2(\Lambda)$ is the ring of upper triangular $(2\times 2)$-matrices with
coefficients in $\Lambda$. 
We consider three related invariants for the objects $X$ 
in $\Cal S(n)$, the {\it mean} $qX$, the {\it level} $pX$
and the {\it colevel} $rX$.  
By definition,  $qX = |V|/bV,\,$ 
$pX = |U|/bV$, and
$rX = |V/U|/bV$. Here, $bV$ denotes the dimension of the kernel of the operator $T$,
thus the number of its Jordan blocks; we call $bV$ the {\it width} of $V$.
The objects $X$ with $bX = 1$ are called
{\it pickets.}

For any $X$ in $\Cal S(n)$, both numbers $pX,\ rX$ are non-negative and 
$pX + rX = qX \le n$. It is the pr-triangle $\mathbb T(n)$ of
vectors $(p,r)$ with $p\ge 0,\ r\ge 0,\ p+r \le n,$ 
which we want to study in order to overview the category $\Cal S(n)$. If $X$ is
an indecomposable object in $\Cal S(n)$, we call $(pX,rX) \in \mathbb T(n)$ 
its {\it support.}

We use $\mathbb T(n)$ to visualize part 
  of the categorical structure of $\Cal S(n)$:
The action of the duality $\D$  and of the square $\tau_n^2$ of the Auslander-Reiten translation 
are represented on  $\mathbb T(n)$  by a reflection and by a rotation by  $120^\circ$,
respectively. Moreover for $n\geq 6$, each component of the Auslander-Reiten quiver of  $\Cal S(n)$
has support either contained in the center of $\mathbb T(n)$  or with the center as its only accumulation point.

We show that the only indecomposable objects $X$ in $\Cal S(n)$ with support having
boundary distance smaller than 1 are the pickets which lie on the boundary, whereas 
any rational vector in $\mathbb T(n)$ with boundary distance at least 2 supports 
infinitely many indecomposable objects. At present, it is not clear at all what happens
for vectors with boundary distance between 1 and 2; several partial results are included 
in the paper. 

The use of $\mathbb T(n)$ provides even in the (quite well-understood) case $n = 6$
some surprises: We will show that any indecomposable object in $\Cal S(6)$ 
lies on one of 12 central lines in $\mathbb T(6)$ and that the center of $\mathbb T(6)$ is
the only
vector which supports infinitely many indecomposables of $\Cal S(6)$. 

A further target of our investigations is to 
single out settings which are purely combinatorial: This concerns not only 
the behaviour near the boundary of the triangle $\mathbb T(n)$, 
but also sets of indecomposable objects: For example, the pickets, the bipickets,
as well as the objects $X = (U,V)$ with $U$ being cyclic.

The paper is essentially self-contained, all prerequisites which are needed 
are outlined in detail. 
\par}

\makeatletter
\renewcommand\l@subsection{\@dottedtocline{2}{1.5em}{3.5em}}
\makeatother
  
\tableofcontents

\renewcommand\listfigurename{Some relevant pictures and tables.}
\listoffigures

\vfill\eject
\section{Introduction.}
\label{sec-one}

{\leftskip=2truecm
\parindent=0pt
    {\bf Question of some undergraduate students:} {\it Ist in der Mathematik
      eigentlich nicht schon alles erforscht?} \newline
{\bf Answer:}  {\it In der Linearen Algebra von endlichdimensionalen Vektor\-r\"aumen
ist in der Tat alles erforscht.}

		\smallskip
        \hfill From an interview with G\"unther M. Ziegler (2007), see \cite{Z}.
\par}
	\bigskip\bigskip 
\subsection{Basic concepts and the triangle \protect $\mathbb T(n)$.}
\label{sec-one-one}

Let $k$ be an arbitrary field. All vector spaces which we consider will
be finite-dimensional 
$k$-spaces, the dimension (or length) of a vector space $V$ will be denoted by $|V|$.
We denote by $\Cal N(n)$ the category of nilpotent operators with nilpotence index at most $n$
(its objects are pairs $V = (V,T)$, where $V$ is a vector space and $T\:V \to V$ a linear
transformation with $T^n = 0$), or, what is the same, the category of $\Lambda$-modules,
where $\Lambda = k[T]/\langle T^n\rangle.$ Let $\Cal N = \bigcup_n \Cal N(n)$. 
For $m\ge 1,$
we denote by $[m]$ the $m$-dimensional vector space with a nilpotent 
operator with kernel of length 1; sometimes, it will be convenient to
write $[0]$ for the zero object in $\Cal N$. 
If $V$ is an object of $\Cal N(n)$, we denote by $[V]$ its isomorphism class. Note that 
the isomorphism classes of objects in $\Cal N(n)$ correspond bijectively to the
partitions bounded by $n$: If $\lambda = 
(\lambda_1,\dots,\lambda_b)$ is a partition, the corresponding 
module is $[\lambda] = 
[\lambda_1,\dots,\lambda_b] = \bigoplus_{i=1}^b [\lambda_i].$ 
The {\it height} of $V$ (or $\lambda$)
is defined to be $\lambda_1$ in case $V \neq 0$, and $0$ 
otherwise. We denote by $bV$ the length of the kernel 
of the operator $T$ and call it the {\it width} of $V$
(thus, $bV$ is the number of Jordan blocks of the operator $T$, or, what
is the same, the cardinality of a minimal generating set of $V$;
note that $bV$ is the Krull-Remak-Schmidt multiplicity of $V$ in
$\Cal N$; and also the number of parts of the partition corresponding
to $[V]$).
	
If $V$ belongs to $\Cal N(n)$ with operator $T$,
a subspace $U \subseteq V$ is said to be an {\it invariant} subspace
provided $T(U) \subseteq U.$
The objects of the category $\Cal S(n)$ are the pairs $X = (U,V),$
where $V$ is in $\Cal N(n)$, and $U\subseteq V$ is an invariant subspace of $V$;
if $X = (U,V)$, we write $VX = V$ and call it the {\it global space} of $X$ and
$UX = U,$ and call it the {\it subspace} in $X$, and we put $vX = \dim VX,$
$uX = \dim UX,$ and $wX=\dim V/U$.
If $X = (U,V)$ is an object of $\Cal S(n)$, we consider (as in \cite{RS1}) the
pair $\bdim X = (uX,vX);$ it is called the {\it dimension vector} of $X$.
But in addition, we also will look at $bX = bV.$
The maps $(U,V) \to (U',V')$ in $\Cal S(n)$ are the $\Lambda$-homomorphisms $f\:V \to V'$
such that $f(U) \subseteq U'.$
Note that $\Cal S(n)$ is an additive, but not an abelian, category.
The {\it direct sum} of objects $(U,V)$ and $(U',V')$ is defined to be
$(U,V)\oplus (U',V') = (U\oplus U', V\oplus V').$ An object $(U,V)$ is {\it indecomposable}
provided it is non-zero and not isomorphic to the direct sum of two non-zero objects.
For $n\leq m$, the category $\Cal S(n)$ is contained in $\Cal S(m)$
hence we may form the union $\Cal S = \bigcup_n \Cal S(n).$
If $X = (U,V)$ is an object in $\Cal S$, the
{\it height} of $X$ is by definition the height of $V$, thus
the smallest number $n$ such that $X$ belongs to $\Cal S(n)$.

Some special objects have to be mentioned already here:
A {\it picket} is an object $X$ in $\Cal S$ with $bX = 1$, thus it is of the form 
$([t],[m])$ with $m\ge 1$ and $0 \le t\le m$.
An indecomposable object $X$
in $\Cal S(n)$ with $bX = 2$ will be called a {\it bipicket.}

\medskip
Our aim is to provide information about the indecomposable objects in $\Cal S(n).$ 
In Sections \ref{sec-one-two} to \ref{sec-one-eight},
we are going to formulate and explain the main results of 
the paper, namely Theorems \ref{theoremone} to \ref{theoremeight}.
In order to overview the indecomposables in $\Cal S(n)$,
we introduce the reference space $\Bbb T(n)$; it is
a subset of the real plane $\Bbb R^2$, namely a triangle and its
interior. We will study it in detail in Section~\ref{sec-one-four}.

By definition, $\Bbb T(n)$ is
the set of pairs $(p,r)$ of real numbers with $p,\ r,$ as well as $p+r$
belonging to the interval $[0,n].$ Given a non-zero object $X$ in $\Cal S(n)$,
we attach to it its {\it pr-vector} $\pr X = (pX,rX)$, where $pX$ is defined as
$pX = uX/bX$ and called the {\it level} of $X$ and where $rX = wX/bX$ is
called the {\it colevel} of $X$ (finally, $qX = pX + rX$ is called the
{\it mean} of $X$; this explains the title of the paper).
As we will see easily, for any object $X$ in $\Cal S(n)$, the value
$\pr X$ belongs to $\Bbb T(n)$, see Section~\ref{sec-one-four}.

The first two results, Theorems~\ref{theoremone} and \ref{theoremtwo} and
their accompanying theorems in Sections~\ref{sec-one-two}
and \ref{sec-one-three}, describe the indecomposable objects
in $\mathcal S(n)$ with pr-vector
near  the boundary of the triangle $\mathbb T(n)$.
What we will achieve can be
\phantomsection{visualized},
\addcontentsline{lof}{subsection}{The boundary of the pr-triangle $\mathbb T(n)$.}%
say in the case $n=8$, as follows: 
$$  
{\beginpicture
   \setcoordinatesystem units <1cm,1cm>
   \setcoordinatesystem units <.2887cm,.5cm>
\multiput{} at -8 1  8 4 /
\setdots <1mm>
\plot -8 0  8 0  0 8  -8 0 /
\multiput{$\ss \bullet$} at -7 1  -6 2  -5 3  -4 4  -3 5  -2 6  -1 7  0 8
    7 1  6 2  5 3  4 4  3 5  2 6  1 7  
    -6 0  -4 0  -2 0  0 0 
    6 0  4 0  2 0  8 0
     -5 1  5 1  0 6 /
\multiput{$\ss \bigcirc$} at  -5 1  5 1  0 6 /
\setdashes <1mm>
\plot -5 1  5 1  0 6  -5 1 /
\setshadegrid span <.5mm>
\vshade -4.8 1.2 1.2 <z,z,,> 
        -4.6 1 1.4 <z,z,,>
        -0.2 1 5.8  <z,z,,> 
        0.2 1 5.8  <z,z,,> 
        4.6 1 1.4 <z,z,,>
        4.8 1.2 1.2 /
\put{$\Bbb T(8)$} at -7 7 

\setsolid 
\arr{5 -.7}{7 -.7}
\put{$r$} at 7.7 -.7
\arr{-2.8 6.5}{-1.8 7.5}
\put{$p$} at -1.3 8

\endpicture}
$$
The interior (shaded) triangle is bounded
by the lines $p = 1,$ $r=1$ and $p+r = n-1.$
We describe 
all the indecomposable objects in $\Cal S(n)$ with pr-vector outside or at a corner 
of the shaded triangle.
On the boundary of $\Bbb T(n)$ are the pr-vectors of the pickets $X = (U,V)$ with
$U = 0,$ or $W = 0,$ or with $V$ a free $\Lambda$-module
(Theorems~\ref{theoremone}, 1.1$'$, 1.1$''$, see Section~\ref{theoremone-dual}),
and there are no other indecomposable objects in
$\Cal S(n)$ outside of the shaded triangle!
The corners of the shaded
triangle have coordinates $(1,1),\ (1,n-2)$ and $(n-2,1)$ and each of them is the 
pr-vector of precisely two indecomposable objects (Theorems~\ref{theoremtwo}, 1.2$'$, 1.2$''$).

This provides essential information about
the indecomposable objects of $\Cal S(n)$ with
boundary distance at most 1. But note: Whereas for all $n$ the number
of indecomposable objects with boundary distance smaller than 1
is finite, namely equal to $3n-1$, we will see in Section~\ref{sec-seven-four}
that for $n\ge 9$, there are infinitely many isomorphism classes
of indecomposable objects $X$ with $pX = 1$, thus all with
boundary distance equal to $1$.

On the other hand,
Theorem~\ref{theoremfour} formulated in Section \ref{sec-one-five} concerns the indecomposable
objects of
$\Cal S(n)$ with boundary distance at least 2. It asserts that
any rational pr-vector with boundary distance at least 2 is
the pr-vector of infinitely many indecomposable objects.

Our investigation of the symmetries of the category $\Cal S(n)$
reveal that $\mathbb T(n)$ should be considered as an equilateral triangle. 
The symmetries of $\Cal S(n)$ are the theme of Section \ref{sec-one-three}:
The decisive Theorem \ref{theoremthree} draws the
attention to the Auslander-Reiten translation $\tau = \tau_n$ in the category $\Cal S(n)$,
as considered in our previous paper \cite{RS2}. As we will see, the square $\tau^2$ of $\tau$ 
categorifies the rotation of $\mathbb T(n)$ by $120^\circ$. This surprising fact is a main
ingredient for nearly all the results of the present paper. 

Sections \ref{sec-one-six} to \ref{sec-one-eight}
are devoted to the half-line support as well as the triangle support
of suitable classes of objects. In this way, we obtain some insight into the structure of
the categories $\Cal S(n)$, in particular about 
the relevance of the central lines in $\mathbb T(n)$
(the lines which pass through the center $z(n) = (n/3,n/3)$ of $\mathbb T(n)$).
The results mentioned until now concern all the categories $\Cal S(n)$, or at least
those of infinite type; these are the categories $\Cal S(n)$ with $n\ge 6.$ 
The final theorem which will be presented in the introduction,
namely Theorem \ref{theoremeight}, draws the attention to
the case $n = 6$. This case seems to be of special interest and 
was considered already in \cite{RS1} and \cite{S2} quite in detail. Whereas the previous
investigations put the emphasis on its large wealth of indecomposables, the present
study establishes (on the contrary) severe finiteness conditions which are valid for $\Cal S(6)$. 
In particular, we will see
in Theorem \ref{theoremeight} that all indecomposable objects of $\Cal S(6)$ lie on just 12 central lines
in $\mathbb T(6)$ --- a very remarkable feature! 

\subsection{First results.}
\label{sec-one-two}

\begin{theorem}
  \label{theoremone}
  Let $X$ be an indecomposable object of $\Cal S$. Then either $uX = 0$
  (and then $bX = 1$, thus $X$ is a picket), or else $uX \ge bX.$
\end{theorem}

\medskip
Theorem~\ref{theoremone}
describes the indecomposable objects with pr-vector on or near to the lower
boundary of $\Bbb T(n)$, as follows; here the  dashed line is the line $p = 1$.
$$  
{\beginpicture
   \setcoordinatesystem units <1cm,1cm>
   \setcoordinatesystem units <.2887cm,.5cm>
\multiput{} at -8 1  8 4 /
\setdots <1mm>
\plot -5 3  -8 0  8 0  5 3  /
\multiput{$\ss \bullet$} at 
    -6 0  -4 0  -2 0  0 0 
    8 0  6 0  4 0  2 0   /
\setdashes <1mm>
\plot -7 1  7 1  /
\setshadegrid span <.5mm>
\vshade -7 1 1 <z,z,,> 
        -5 1 3 <z,z,,>
        5 1 3 <z,z,,>
        7 1 1  /
\put{$\Bbb T(8)$} at -10 3 

\setsolid 
\arr{5 -.7}{7 -.7}
\put{$r$} at 7.7 -.7
\arr{-7.8 1.5}{-6.8 2.5}
\put{$p$} at -6.3 3

\endpicture}
$$
The only indecomposable objects $X = (U,V)$ with $uX < bX$ are pickets, namely
the $n$ pickets of the form $(0,[t])$ with $1 \le t \le n$; and the shaded part 
contains all the pr-vectors of objects $X$ with $pX \ge 1$, thus with $uX \ge bX.$

\medskip
\begin{theorem}
  \label{theoremtwo}
  Let $X$ be an indecomposable object of $\Cal S$. If $vX \le 2bX$,
  then $X$ belongs to $\Cal S(3)$ (and is a picket or a bipicket).
\end{theorem}

	\medskip 
We picture the indecomposable objects in
$\Cal S$ with $vX \le 2\cdot bX$:
There are precisely 6 isomorphism classes, namely five pickets
and the object $E^2_2$.

$$
{\beginpicture
    \setcoordinatesystem units <.3cm,.3cm>
\put{\beginpicture
\multiput{} at  0 0  1 3 /
\plot 0 0  1 0  1 1  0 1  0 0  /
\endpicture} at 0 0
\put{\beginpicture
\multiput{} at  0 0  1 3 /
\plot 0 0  1 0  1 1  0 1  0 0  /
\put{$\bullet$} at 0.5  0.5 
\endpicture} at 5 0
\put{\beginpicture
\multiput{} at  0 0  1 3 /
\plot 0 0  1 0  1 2  0 2  0 0  /
\plot 0 1  1 1 /
\endpicture} at 10 0
\put{\beginpicture
\multiput{} at  0 0  1 3 /
\plot 0 0  1 0  1 2  0 2  0 0  /
\plot 0 1  1 1 /
\put{$\bullet$} at 0.5  0.5 
\endpicture} at 15 0
\put{\beginpicture
\multiput{} at  0 0  1 3 /
\plot 0 0  1 0  1 2  0 2  0 0  /
\plot 0 1  1 1 /
\put{$\bullet$} at 0.5  1.5 
\endpicture} at 20 0
\put{\beginpicture
\multiput{} at 0 0  1 3 /
\plot 0 0  1 0   1 3  0 3  0 0  /
\plot 0 1  2 1  2 2  0 2  /
\multiput{$\bullet$} at 0.5 1.5  1.5 1.5 /
\plot 1.5 1.5  .5 1.5 /
\plot 1.5 1.55 .5 1.55 /
\plot 1.5 1.45 .5 1.45 /
\endpicture} at 25 0
\put{$([0],[1])$} at 0 -4
\put{$([1],[1])$} at 5 -4
\put{$([0],[2])$} at 10 -4
\put{$([1],[2])$} at 15 -4
\put{$([2],[2])$} at 20 -4
\put{$E_2^2$} at 25 -4

\endpicture}
$$

\medskip
The bipicket $E_2^2$ has global space $[3,1]$, subspace $[2]$ and
factor space $[2]$.
As in the previous paper \cite{RS1}, we try to visualize objects in $\Cal S$ by 
using boxes and bullets connected by lines; for an outline,
see Section~\ref{sec-two-two} below. 

\medskip
The essential assertion of Theorem~\ref{theoremtwo} can be reformulated as follows:
{\it There are only $2$ isomorphism classes of
indecomposable objects $X$ in $\Cal S$ 
with  $uX = bX = wX$, namely the picket $([1],[2])$ and 
the bipicket $E_2^2$.}

\bigskip
Also here, we provide a picture of $\Bbb T(n)$, in case $n \ge 3$. 
Theorem~\ref{theoremtwo}
concerns the lower left corner of $\Bbb T(n)$, since it describes all the
indecomposable objects $X$ in $\Cal S(n)$
with $vX \le 2bX$ thus with $qX \le 2$ (the dashed line is the line $q = 2$).
The bullets refer to pickets, the
circle to a bipicket: 
$$  
{\beginpicture
   \setcoordinatesystem units <1cm,1cm>
   \setcoordinatesystem units <.2887cm,.5cm>
\multiput{} at -8 1  8 4 /
\setdots <1mm>
\plot -5 3  -8 0  8 0  5 3  /
\plot  -5 1  -7 1  -6 0  -5 1 /
\setdashes <1mm>
\plot -6 2  -4 0  /
\setshadegrid span <.5mm>
\vshade -6 2 2 <z,z,,> 
        -5 1 3 <z,z,,>
        -4 0 3  <z,z,,>
         5 0 3  <z,z,,>
         8 0 0  /

\multiput{$\ss \bullet$} at -7 1  -6 2  
    -6 0  -4 0  
     -5 1   /
\multiput{$\ss \bigcirc$} at  -5 1   /

\setsolid 
\arr{5 -.7}{7 -.7}
\put{$r$} at 7.7 -.7
\arr{-7.8 1.5}{-6.8 2.5}
\put{$p$} at -6.3 3

\put{$\Bbb T(n)$} at -10 3 

\endpicture}
$$

\medskip
Theorems 1.1 and 1.2
may be considered as basic statements in linear algebra:

\medskip
\begin{theoremone-ref}
  \label{theoremone-ref}
  Let $V=(V,T)$ be a linear operator and $U$ a $T$-invariant subspace
  of $V$ such that the pair $(U,V)$ is indecomposable.
  Then either $ U=0 $ (and $V$ has only one Jordan block)
  or the dimension of $ U $ is at least the number of Jordan blocks of $ V $.
\end{theoremone-ref}

\smallskip
Theorem~\ref{theoremone} can also be phrased as follows: {\it Let $X = (U,V)$ be an object in
$\Cal S$. If $uX < bX,$ 
then $V$ has a direct decomposition $V = V_1\oplus V_2$ with $TV_i \subseteq V_i$
for $i = 1,2,$ such that $U \subseteq V_1$ and $V_2\neq 0.$ }
	
\medskip
\begin{theoremtwo-ref}
  \label{theoremtwo-ref}
  Let  $(U,V)$  be a non-zero pair in $\Cal S$ such that
  the average size of the Jordan blocks of  $V$  is at most two.
  Then  $(U,V)$  has a direct summand which
  is either a picket $([t],[m])$ with $m\leq 2$ or the bipicket $E_2^2$.
\end{theoremtwo-ref}

\subsection{Symmetries: duality and rotation.}
\label{sec-one-three}

We next draw the attention to two important internal symmetries of the category $\Cal S(n)$.
First of all, there is the duality functor $\D$ 
on $\Cal S(n)$, defined by $\D(U,V) = ((V/U)^*,V^*),$
where $V^* = \Hom(V,k)$. Obviously, this contravariant functor has the property that
$$
  u(\D X) = wX,\ v(\D X) = vX,\  w(\D X) = uX,\
\text{and}\ \ b(\D X) = bX,
$$ 
for all objects $X\in \Cal S(n)$. Using $\D$, we obtain:

\medskip
\begin{theoremone-dual}
  \label{theoremone-dual}
  Let $X$ be an indecomposable object of $\Cal S$. Then either $wX = 0$
  (and $X$ is a picket), or else $wX \ge bX.$
\end{theoremone-dual}

\medskip
Of course, it is not surprising that Theorem \ref{theoremone}
has this consequence, since dealing with finite-dimensional vector spaces, there 
always are such duality features. 
	
There is a second symmetry, and this symmetry is 
really unexpected and exciting. It is based on the (relative) 
Auslander-Reiten functor 
of $\Cal S(n)$, as considered in \cite{RS2}.
Let us recall the setting.
The category $\Cal S(n)$ may be identified with the category of the torsionless 
$T_2(\Lambda)$-modules, where $\Lambda = k[T]/\langle T^n\rangle$
and $T_2(\Lambda)$ is the ring of upper triangular $(2\times 2)$-matrices with
coefficients in $\Lambda.$ Since $\Cal S(n)$ is just the category of torsionless 
$T_2(\Lambda)$-modules, $\Cal S(n)$ has (relative) Auslander-Reiten sequences 
(by Auslander-Smal\o{}), and we denote by $\tau_n$ (or just $\tau$) the relative
Auslander-Reiten translation of $\Cal S(n)$. 

The pickets $(0,[n])$ and $([n],[n])$ will be said
to be the {\it projective} pickets, since they are the indecomposable projective 
$T_2(\Lambda)$-modules.
We say that $X$ in $\Cal S(n)$ is {\it reduced}  
provided $X$ has no direct summand which is a projective picket.
If  $X$ is a projective picket,
then $\tau_n X = 0$. 
{\it If $X$ in $\Cal S(n)$ is reduced, then $\tau_n X$ is also reduced 
and $\tau_n^6 X = X$,} see \cite[Corollary 6.5]{RS2}.

\medskip
\begin{theorem}
  \label{theoremthree}
  If $X$ is a reduced object of $\Cal S(n)$, 
  then we have (for $\tau = \tau_n$):
$$
u(\tau^2X) = wX,\ v(\tau^2X) =  n\cdot bX-uX ,\ w(\tau^2 X) =  n\cdot bX-vX,
\ \text{\it and}\ \ 
b(\tau^2X) = bX.
$$
\end{theorem}

\medskip	
In particular, we see that the dimension vector (as well as the width) of 
$\tau^2 X$ only depends on the dimension vector and the width of $X$ (this leads us in
Section~\ref{sec-one-four} to focus the attention to these invariants). Let us stress that 
the corresponding assertions for $\tau$ itself do not hold:
Starting with indecomposable objects $X, X'$ with equal dimension vector and equal width, 
the dimension vectors of $\tau X$ and $\tau X'$ (also the width) usually are different. 
Already in $\Cal S(5)$, we may look at the two indecomposable objects $X = (U,V),\ X'= (U',V')$, both
with dimension vector $(2,6)$ and width 2, namely with $V = [4,2]$ and $V' = [5,1]$;
the total space of $\tau_5X$ is $[5,3,1]$, its subspace $[4,2]$, thus $\bdim \tau_5X = (6,9)$,
whereas $\tau_5 X' = ([1],[4])$ is the picket with dimension vector $\bdim \tau_5X' = (1,4)$.

\medskip
Using $\tau^2$, we can reformulate Theorem \ref{theoremone} as follows:

\medskip
\begin{theoremone-rot}
  \label{theoremone-rot}
  Let $X$ be an indecomposable object of $\Cal S(n)$. Then either $vX = n\cdot bX$
  (and $X$ is a picket), or else $vX \le (n-1)bX.$
\end{theoremone-rot}

\medskip
In order to obtain Theorem 1.1$''$ from
Theorem 1.1$'$,
we start with $X$ indecomposable. Either $X$ is a projective picket, then $vX = n,$ and
$bX = 1$, thus $vX = n\cdot bX$. 
Thus, we can assume that $X$ is not a projective picket, thus $\tau^2 X$ is indecomposable 
and Theorem 1.1$'$ asserts that  $w(\tau^2 X) = 0$
and $\tau^2 X$ is a picket, or else that $w(\tau^2 X) \ge b(\tau^2 X) = bX.$
Now $w(\tau^2 X) =  n\cdot bX-vX$. If $w(\tau^2 X) = 0$, 
then $n\cdot bX = vX$ (and $bX = b(\tau^2 X)
= 1$ shows that $X$ is a picket). If $w(\tau^2 X) \ge bX,$ then 
$n\cdot bX-vX = w(\tau^2 X) \ge bX,$ and therefore $vX \le (n-1)bX.$ $\s$

\medskip
Looking back at Theorem \ref{theoremtwo}, the indecomposable objects with $uX < bX$
are known by Theorem \ref{theoremone}, those with $wX < bX$ are known by
Theorem~1.1$'$;
therefore, as we have mentioned already in Section~\ref{sec-one-two},
the essential assertion of Theorem \ref{theoremtwo} is:
{\it There are only $2$ isomorphism classes of
indecomposable objects $X$ in $\Cal S$ 
with  $uX = bX = wX$, namely the picket $([1],[2])$ and 
the bipicket $E_2^2$.} Using Theorem \ref{theoremthree}, we see:

\medskip
\begin{theoremtwo-dual}
  \label{theoremtwo-dual}
  There are only $2$ isomorphism classes of
  indecomposable objects $X$ in $\Cal S(n)$ 
  with  $uX = bX$ and $vX = (n-1)bX$, namely 
  the picket $([1],[n-1])$ and the bipicket $\tau^2E_2^2$.
\end{theoremtwo-dual}

\medskip
\begin{theoremtwo-rot}
  \label{theoremtwo-rot}
  There are only $2$ isomorphism classes of
  indecomposable objects $X$ in $\Cal S(n)$ 
  with  $wX = bX$ and $vX = (n-1)bX$, 
  the picket $([n-2],[n-1])$ and the bipicket $\tau^4E_2^2$.
\end{theoremtwo-rot}

\medskip
The four objects mentioned in Theorem~1.2$'$ and 
Theorem~1.2$''$ can easily be visualized;
for $n = 6,$ they look as follows:
$$
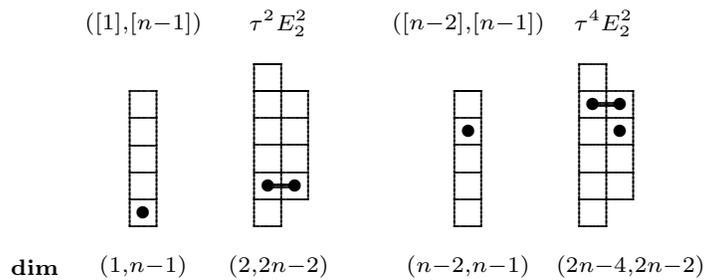
{\beginpicture
    \setcoordinatesystem units <.3cm,.3cm>
\put{\beginpicture
\multiput{} at 0 0  1 6 /
\plot 0 0  1 0  1 5  0 5  0 0  /
\plot 0 1  1 1 /
\plot 0 4  1 4 /
\plot 0 3  1 3 /
\plot 0 2  1 2 /
\put{$\bullet$} at 0.5  .5 
\endpicture} at 0 0
\put{\beginpicture
\multiput{} at 0 0  1 6 /
\plot 0 0  1 0   1 6  0 6  0 0  /
\plot 0 1  2 1  2 5  0 5  /
\plot 0 5  2 5 /
\plot 0 6  1 6 /
\plot 0 4  2 4 /
\plot 0 3  2 3 /
\plot 0 2  2 2 /
\multiput{$\bullet$} at 0.5 1.5  1.5  1.5 /
\plot 1.5 1.5  .5 1.5 /
\plot 1.5 1.55 .5 1.55 /
\plot 1.5 1.45 .5 1.45 /
\endpicture} at 5 0

\put{\beginpicture
\multiput{} at 0 0  1 6 /
\plot 0 0  1 0  1 5  0 5  0 0  /
\plot 0 1  1 1 /
\plot 0 5  1 5 /
\plot 0 4  1 4 /
\plot 0 3  1 3 /
\plot 0 2  1 2 /
\put{$\bullet$} at 0.5  3.5 
\endpicture} at 12 0
\put{\beginpicture
\multiput{} at 0 0  1 6 /
\plot 0 0  1 0   1 6  0 6  0 0  /
\plot 0 1  2 1  2 5  0 5  /
\plot 0 4  2 4 /
\plot 0 3  2 3 /
\plot 0 2  2 2 /
\multiput{$\bullet$} at 0.5  4.5  1.5 4.5  1.5 3.5 /
\plot 1.5 4.5  .5 4.5 /
\plot 1.5 4.55 .5 4.55 /
\plot 1.5 4.45 .5 4.45 /
\endpicture} at 17 0

\put{$\ssize \bdim$} at -4 -4.5
\put{$\ssize (1,n-1)$} at 0 -4.5
\put{$\ssize (2,2n-2)$} at 5 -4.5
\put{$\ssize (n-2,n-1)$} at 12 -4.5
\put{$\ssize (2n-4,2n-2)$} at 18 -4.5

\put{$\ss([1],[n-1])$} at 0 4.5
\put{$\ss\tau^2E_2^2$} at 5 4.5
\put{$\ss([n-2],[n-1])$} at 12 4.5
\put{$\ss\tau^4E_2^2$} at 17 4.5
\endpicture}
$$

\subsection{uwb-vectors and symmetries in the pr-triangle $\mathbb T(n)$.}
\label{sec-one-four}

Looking at indecomposable objects $X$ in $\Cal S(n)$, the paper \cite{RS1} has drawn the attention
to the dimension vector $\bdim X = (uX,vX)$, thus to the invariants $u$ and $v$, or
equivalently, to the invariants $u$ and $w$ (since $v = u+w)$. 
According to the results mentioned already, one
also should take into account the invariant $b$
(after all, the objects in $\Cal S$ are vector spaces with an operator 
and a subspace: The invariants $u$ and $w$ refer to the space and its subspace, and it is 
the invariant $b$ which points to the operator). To repeat: It seems important
to focus the attention not only to the invariants $u$ and $w$, but also to $b$. It is 
{\it the relationship between 
the invariants $u,$ $w$ and $b$,} which
has to be studied. Starting with these three invariants $u,\ w,\ b$,
one may build further invariants: linear combinations, quotients, and so on. 
Our main interest will lie in the quotients $u/b$ and $w/b$.
(Later, the treatment of central lines will rely on the ratio $(3u-nb)/(3w-nb)$.)

	\bigskip
Given an object $X$ in $\Cal S$, we 
call the triple $\uwbb X = (uX,wX,bX)$ the {\it uwb-vector} of $X.$ Actually, 
instead of looking at the reference space $\mathbb R^3$ with the triples $(uX,wX,bX)$, 
we will focus the attention to the
corresponding projective space which contains for a non-zero object $X$ 
the pair $\pr X = (uX/bX,wX/bX)$; we call $\pr X$ the {\it support} of $X$.

Thus, here are the main definitions for the present paper:
Let $X$ be a non-zero object in $\Cal S$. 
The {\it level} $pX$, the {\it mean} $qX$, and the {\it colevel} $rX$ of $X$ are defined as
follows:
$$
   pX = uX/bX, \quad qX = vX/bX, \quad rX = wX/bX
$$
(Since for $X = 0$, the numbers $uX$, $vX$, $wX$, and $bX$ are all
zero, the quotients used in the definition of $pX$, $qX$ and $rX$ are not defined,  thus
mean, level, and colevel can only be considered for non-zero objects.) 
For any non-zero object $X$ in $\Cal S$, we have $qX = pX+rX$,
thus any two of the invariants $pX,\ qX,\ rX$ determine the third one. 	
[A short hint to explain the terminology. 
Recall that the isomorphism class of a module in $\Cal N$ may be considered 
as a partition. The ``mean'' of an object $X = (U,V)$ in $\Cal S$ is just
the mean (or average) of the partition $[V]$, thus 
the mean (or average) of the sizes of the Jordan blocks of the given operator. 
The intuition behind the chosen word ``level'' stems from the vision of considering $U$ 
as a kind of filling of $V$, thus we measure how much of the global space is filled by the
subspace (by dividing the length $|U|$ of the subspace 
through $bV$, the length of the socle of $V$).]
	
The aim of the present paper is to consider the support $\pr X = (pX,rX)$ of the indecomposable
objects $X$ in $\Cal S(n)$. 
Of course, for any non-zero $X$, both numbers 
$pX,\ rX$ are non-negative. It is obvious that for $X$ in $\Cal S(n)$, we have
$vX \le n\cdot bX$, thus $qX \le n$.  
It is the {\bf pr-triangle} $\mathbb T(n)$ of
vectors $(p,r)$ with $p\ge 0,\ r\ge 0,\ p+r \le n,$ 
which we want to study. If $\mathbb X$ is a subset of $\mathbb T(n)$,
and $X$ is an indecomposable object in $\Cal S(n)$, we will say that $X$ {\it lives} on
$\mathbb X$ (or also that $\mathbb X$ {\it supports} $X$) 
provided $\pr X$ belongs to $\mathbb X$.
The supports of the pickets in $\Cal S(n)$ together with 
the zero vector $(0,0)$ provide the grid of the vectors in $\mathbb T(n)$ 
with integral coefficients. 

We denote by $z(n) = (\frac n3,\frac n3)$ the center of $\mathbb T(n)$ (this may be thought of
as the center of gravity of the triangle). 
An indecomposable object $X$ with $\pr X = z(n)$ is said to be
{\it central}. 

Let us present a picture of the triangle $\mathbb T(n)$. We show $\mathbb T(8)$
with its triangular grid given by the lines with $p,\ q,$ or $r$ being an integer.
The center $z(8)$ is marked by a black square $\ss \blacksquare$.
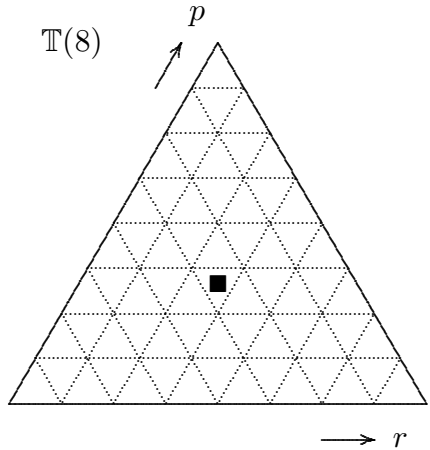
$$  
{\beginpicture
   \setcoordinatesystem units <1cm,1cm>
\put{\beginpicture
   \setcoordinatesystem units <.57735cm,1cm>
\multiput{} at -6 -.5  6 4.4 /
\plot -4 0  4 0  0 4  -4 0 /
\setdots <.5mm>
\multiput{$\ssize \blacksquare$} at 0 1.333  /
\setdots <.5mm> 

\plot -2.5 1.5  -2 1 .5 3.5  -.5 3.5  2 1  2.5 1.5  -2.5 1.5 /
\plot -2 2      -1 1  1 3  -1 3  1 1  2 2  -2 2 /
\plot -1.5 2.5  0 1  1.5 2.5  -1.5 2.5 /
\plot -3 1  3 1  /
\plot -3.5 0.5  3.5 0.5 /
\plot -3 1  -2 0  -1 1  0 0  1 1  2 0  3 1 /
\plot -3.5 .5  -3 0 -2 1  -1 0  0 1  1 0  2 1  3 0  3.5  0.5  /
\put{$\mathbb T(8)$} at -2.8 4 
\setsolid 
\arr{2 -.4}{3 -.4}
\put{$r$} at 3.5 -.4
\arr{-1.2 3.5}{-.7 4}
\put{$p$} at -.4 4.3
\endpicture} at 0 0 
\endpicture}
$$

We stress that the triangle $\mathbb T(n)$ should be considered as an equilateral one, since
the dihedral group $\Sigma_3$ of order 6, as the symmetry group of $\mathbb T(n)$,
plays a decisive role in our investigation. 
Let $\D$ be the reflection of $\mathbb T(n)$ defined by
$\D(p,r) = (r,p),$  and 
$\rho$ the rotation of $\mathbb T(n)$ by $120^\circ$ (with center $z(n)$),
so that $\rho(p,r) = (r,n-p-r).$ Thus,
the symmetry group $\Sigma_3$ of the triangle $\mathbb T(n)$ is generated by
$\rho$ and $\D$. According to Section \ref{sec-one-three}, we have:
{\it If $X$ is a non-zero object in $\Cal S(n)$, then $b(\D X) = bX$ and 
$\pr(\D X) = \D\pr X.$} And second, Theorem \ref{theoremthree} can be reformulated as follows:

\medskip
\begin{theoremthree-ref}
  \label{theoremthree-ref}
  If $X$ is a reduced non-zero object in $\Cal S(n)$, then
    $$
    \pr \tau^2 X = \rho\pr X,  \quad\text{\it and} \quad b(\tau^2X) = bX.
    $$
\end{theoremthree-ref}

\medskip 

We have shown in \cite{RS2} that $\tau^6 X = X$ for any reduced object
in $\Cal S(n)$ and that, as a consequence, all Auslander-Reiten 
components of $\Cal S(n)$ are tubes of rank 1, 2, 3 or 6.

\medskip
\begin{corollarythree-ref}
  \label{corollarythree-ref}
  Let $n\geq 6$.  If $X$ is an indecomposable object in $\Cal S(n)$ which
  occurs in a tube of rank $1$ or $2,$ then $X$ is central.
\end{corollarythree-ref}

\smallskip
Namely, if $X=\tau^2X$ then $\pr X=\rho\pr X$ so $X$ is central. $\s$

\medskip{\bf Boundary distance.} 
The boundary of $\mathbb T(n)$
consists of the vectors $(p,r)$ with $p = 0,$ or $r = 0$, or $p+r = n.$
Given a vector $(p,r)$ in $\mathbb T(n)$, its {\it boundary distance} is 
$d(p,r) = \min\{p,r,n-p-r\}.$
Thus, $(p,r)$ belongs to the boundary if and only if its boundary distance is equal to $0$, whereas
$z(n) = (n/3,n/3)$ is the only vector with boundary distance $n/3$ (and there are no vectors with
boundary distance greater than $n/3$). 
For $0\le d < n/3$, we write $\Delta_d$ for the 
set of pr-vectors with boundary distance $d,$ and call $\Delta_d$ a {\it standard
triangle.} The standard triangles and the one-element set $\{z(n)\}$
provide a partition of $\mathbb T(n)$. 

If $X$ is indecomposable in $\Cal S(n)$, the number $dX = d(\pr X)$ 
is called the {\it boundary distance} of $X$. Of course, $X$ is central if and only if 
$dX = d(z(n))=n/3$; also, $X$ lies on the standard triangle $\Delta_d$ 
if and only if $dX = d.$ The invariant $d$ (and the corresponding standard
triangles) play an important role in our investigation.
One of the reasons is of course the equality $d(\tau^2 X) = dX$ for $X$
reduced. 

A picket whose pr-vector lies on the boundary will be called a {\it boundary picket.}
There are $3n-1$ boundary pickets in $\Cal S(n)$, namely $n$ pickets of the form
$([0],[m])$ with $1\le m\le n$ (the {\it void} pickets), 
also $n$ pickets of the form $([m],[m])$ with $1\le m \le n$ (the {\it full} pickets),
and finally $n-1$ pickets of the form $([t],[n])$ with $1\le t < n$ (the non-projective
pickets of height $n$).

Section \ref{sec-one-one} provides an illustration of Theorems~\ref{theoremone}, \ref{theoremtwo}
and their accompanying theorems in Sections~\ref{sec-one-two} and \ref{sec-one-three},
in the case $n = 8.$
In the picture, the
boundary triangle $\Delta_0$ is dotted, the triangle $\Delta_1$ is dashed. 
The small bullets $\scriptstyle\bullet$ mark
the position of the boundary pickets, as well as the corners of $\Delta_1$. These corners are,
in addition, encircled, in order to indicate that they support also one of the bipickets
$E_2^2$, $\tau^2 E_2^2$ and $\tau^4E_2^2$. The support of the remaining indecomposable
objects has boundary distance at least 1 and is not a corner of $\Delta_1$, thus it lies in the
shaded region.

\medskip
Note that for $n = 3$, Theorems~\ref{theoremone}, \ref{theoremtwo}
and their accompanying theorems in Sections~\ref{sec-one-two} and \ref{sec-one-three}
recover the full classification of the
indecomposable objects in $\Cal S(n)$: There are 9 pickets and the bipicket $E_2^2.$
Below, we provide two versions of the triangle $\mathbb T(3).$ 
In the left version, the pickets are marked by bullets, whereas the circle 
at the center $z(3)$ indicates 
that $z(3)$ is also the pr-vector of the only additional indecomposable object, 
the bipicket $E^2_2$. On the right, 
all the individual objects are visualized.
$$
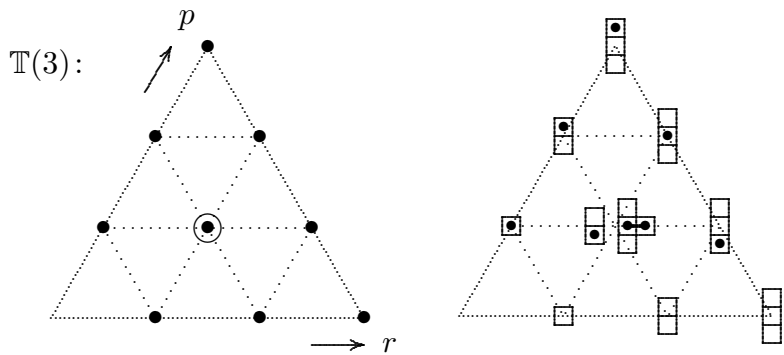
{\beginpicture
   \setcoordinatesystem units <.57735cm,1cm>
\put{\beginpicture
\put{$\mathbb T(3)\:$} at -3 2.8
\multiput{} at -3 0  3 0  0 3 /
\setdots <.5mm>
\plot -3 0  3 0  0 3  -3 0 /
\setdots <1mm>
\plot -2 1  -1 0  0 1  1 0  2 1  -2 1 /
\plot -1 2  0 1  1 2  -1 2 /
\multiput{$\bullet$} at -1 0  1 0  3 0  -2 1  0 1  2 1  -1 2  1 2  0 3 /
\put{$\bigcirc$} at 0 1 
\setsolid 
\arr{2 -.3}{3 -.3}
\put{$r$} at 3.5 -.3
\arr{-1.2 2.5}{-.7 3}
\put{$p$} at -.4 3.3
\endpicture} at 0 0 
\put{\beginpicture
\multiput{} at -3 0  3 0  0 3 /
\setdots <.5mm>
\plot -3 0  3 0  0 3  -3 0 /
\setdots <1mm>
\plot -2 1  -1 0  0 1  1 0  2 1  -2 1 /
\plot -1 2  0 1  1 2  -1 2 /

\put{\beginpicture
   \setcoordinatesystem units <.2cm,.2cm>
  \multiput{} at 0 0  1 1 /
  \setsolid
  \plot 0 0  1 0  1 1  0 1  0 0 /
  \endpicture} at -1 0
\put{\beginpicture
   \setcoordinatesystem units <.2cm,.2cm>
  \multiput{} at 0 0  1 2 /
  \setsolid
  \plot 0 0  1 0  1 2  0 2  0 0 /
  \plot 0 1  1 1 /
  \endpicture} at 1 0
\put{\beginpicture
   \setcoordinatesystem units <.2cm,.2cm>
  \multiput{} at 0 0  1 3 /
  \setsolid
  \plot 0 0  1 0  1 3  0 3  0 0 /
  \plot 0 1  1 1 /
  \plot 0 2  1 2 /
  \endpicture} at 3 0

\put{\beginpicture
   \setcoordinatesystem units <.2cm,.2cm>
  \multiput{} at 0 0  1 1 /
  \setsolid
  \plot 0 0  1 0  1 1  0 1  0 0 /
  \put{$\ssize\bullet$} at 0.5 0.5 
  \endpicture} at -2 1
\put{\beginpicture
   \setcoordinatesystem units <.2cm,.2cm>
  \multiput{} at 0 0  1 2 /
  \setsolid
  \plot 0 0  1 0  1 2  0 2  0 0 /
  \plot 0 1  1 1 /
  \put{$\ssize\bullet$} at 0.5 1.5 
  \endpicture} at -1 2
\put{\beginpicture
   \setcoordinatesystem units <.2cm,.2cm>
  \multiput{} at 0 0  1 3 /
  \setsolid
  \plot 0 0  1 0  1 3  0 3  0 0 /
  \plot 0 1  1 1 /
  \plot 0 2  1 2 /
  \put{$\ssize\bullet$} at 0.5 2.5 
  \endpicture} at 0 3
\put{\beginpicture
   \setcoordinatesystem units <.2cm,.2cm>
  \multiput{} at 0 0  1 3 /
  \setsolid
  \plot 0 0  1 0  1 3  0 3  0 0 /
  \plot 0 1  1 1 /
  \plot 0 2  1 2 /
  \put{$\ssize\bullet$} at 0.5 1.5 
  \endpicture} at 1 2
\put{\beginpicture
   \setcoordinatesystem units <.2cm,.2cm>
  \multiput{} at 0 0  1 3 /
  \setsolid
  \plot 0 0  1 0  1 3  0 3  0 0 /
  \plot 0 1  1 1 /
  \plot 0 2  1 2 /
  \put{$\ssize\bullet$} at 0.5 0.5 
  \endpicture} at 2 1

\put{\beginpicture
   \setcoordinatesystem units <.2cm,.2cm>
  \multiput{} at 0 0  1 2 /
  \setsolid
  \plot 0 0  1 0  1 2  0 2  0 0 /
  \plot 0 1  1 1 /
  \put{$\ssize\bullet$} at 0.5 .5 
  \endpicture} at -.4 1

\put{\beginpicture
   \setcoordinatesystem units <.2cm,.2cm>
  \multiput{} at 0 0  2 3 /
  \setsolid
  \plot 0 0  1 0  1 3  0 3  0 0 /
  \plot 0 1  2 1  2 2  0 2 /
  \multiput{$\ssize\bullet$} at 0.5 1.5  1.5 1.5 /
  \plot 0.5 1.5  1.5 1.5 /
  \plot 0.5 1.45  1.5 1.45 /
  \plot 0.5 1.55  1.5 1.55 /
  \endpicture} at .4 1 
\endpicture} at 8 0

\endpicture}
$$

\medskip
As we have seen in Section~\ref{sec-one-three}, the vectors in $\mathbb T(n)$ which support indecomposable
objects and which have boundary distance smaller than one
have integral coordinates. 
The coordinates of the pr-vectors $(p,r)$ of a bipicket have
denominator at most $2$, thus $(p,r)$ belongs to a grid line. But already for $n = 5$, there
are indecomposable objects which are neither pickets, nor bipickets
(still all vectors in $\mathbb T(5)$ which support indecomposable
objects lie on grid lines, see Sections~\ref{sec-fifteen-one} and \ref{sec-ten-nine}). Indecomposable
objects $X$ whose pr-vectors do not lie on a grid line occur for $n = 6$; here are two
examples $X$ and $Y$ with $\pr X = (5/3,5/3)$  and $\pr Y = (7/3,7/3).$

$$
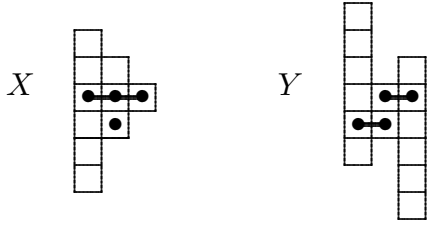
{\beginpicture
    \setcoordinatesystem units <.3cm,.3cm>
\put{\beginpicture
\multiput{} at 0 0  3 6 /
\plot 0 0  1 0   1 6  0 6  0 0  /
\plot 0 1  1 1 /
\plot 0 2  2 2  2 5  0 5 /
\plot 0 3  3 3  3 4  0 4 /
\plot 0 2  2 2 /
\multiput{$\bullet$} at 0.5 3.5  1.5  3.5  2.5 3.5  1.5 2.5 /
\plot 2.5 3.5  .5 3.5 /
\plot 2.5 3.55 .5 3.55 /
\plot 2.5 3.45 .5 3.45 /
\put{$X$} at -2 4
\endpicture} at 0 0
\put{\beginpicture
\multiput{} at 0 0  3 8 /
\plot 0 2  1 2  1 8  0 8  0 2 /
\plot 2 0  3 0  3 6  2 6  2 0 /
\plot 2 1  3 1 /
\plot 2 2  3 2 /
\plot 0 3  3 3 /
\plot 0 4  3 4 /
\plot 0 5  3 5 /
\plot 0 6  1 6 /
\plot 0 7  1 7 /
\multiput{$\bullet$} at 0.5 3.5  1.5  3.5  1.5 4.5  2.5 4.5 /
\plot 1.5 3.5  .5 3.5 /
\plot 1.5 3.55 .5 3.55 /
\plot 1.5 3.45 .5 3.45 /

\plot 2.5 4.5  1.5 4.5 /
\plot 2.5 4.55 1.5 4.55 /
\plot 2.5 4.45 1.5 4.45 /
\put{$Y$} at -2 5
\endpicture} at 10 0
\endpicture}
$$
There are many more, as we will see.

\subsection{pr-vectors with boundary distance at least 2.}
\label{sec-one-five}

As we have seen in Section~\ref{sec-one-three}, there are only finitely many indecomposable objects
which have boundary distance smaller than 1,  namely $3n-1$ pickets.
And there are precisely two indecomposable objects which lie 
on a given corner of $\Delta_1$: a picket and a bipicket. 

A {\it $\mathbb P^1$-family} $M = \{M_c\mid c\in \mathbb P^1(k)\}$  of objects is
a set of pairwise non-isomorphic 
indecomposable objects with fixed uwb-vector, indexed by the elements of the
projective line $\mathbb P^1 = \mathbb P^1(k)$
over $k$ (the elements $c = (c_0:c_1)$ of $\mathbb P^1$ are the
one-dimensional subspaces of $k^2$, we write $c = (c_0:c_1)$ for the subspace generated
by the non-zero element $(c_0,c_1) \in k^2$). 

We say that the vector $(p,r)$ in $\mathbb T(n)$ 
is a {\it BTh-vector} provided there is $a\in\mathbb N_1$ (the set
of positive integers) 
such that for any $t\in \mathbb N_1$, there is a $\mathbb P^1$-family of indecomposable objects in
$\Cal S(n)$ with
uwb-vector $(atp,atr,at)$ (the corresponding {\it BTh-family}).
(Here, the letters BTh refer to Brauer-Thrall, since the property which defines
BTh-vectors reminds of the second Brauer-Thrall conjecture in the form which seems now 
to be the accepted one: To have $\mathbb P^1$-families of indecomposable objects for all
dimensions which are a multiple of a fixed one.)
Of course: {\it If $(p,r)$ is a BTh-vector, then there are infinitely many 
isomorphism classes of indecomposable objects in $\Cal S(n)$ with pr-vector $(p,r)$,}
even if $k$ is finite. 

It is well-known that 
$\Cal S(5)$ has only finitely many (namely 50) isomorphism classes
of indecomposable objects, see for example \cite{RS1}, thus 
there are no BTh-vectors in $\mathbb T(5)$. 
On the other hand,  $(2,2) \in \mathbb T(6)$ is a BTh-vector, 
Namely, there is 
a one-parameter family of indecomposable objects $X = (U,V)$, where
$V = [6,4,2]$ and  the subspaces $U$ are suitable 6-dimensional subspaces of $V$
(see \cite{RS1}, or already \cite{Bh};
details concerning this family will be recalled in the example in Section~\ref{sec-two-seven}).
There are indecomposable objects $X[t]$ in $\Cal S(6)$ which have a filtration with $t$ factors,
all being isomorphic to $X$, where $t \ge 1.$ 
Since $\uwbb X = (6,6,3)$, we have $\uwbb X[t] = (6t,6t,3t)$
(and therefore $\pr X[t] = (2,2)$) for all $t\ge 1.$
	
Actually, $(2,2)$ is the  only BTh-vector in $\mathbb T(6)$, 
see Theorem \ref{theoremeight} below. 
In general, for $n\ge 6$, let us look at the set of pr-vectors with boundary distance at least 2
(in the case $n = 6,$ this set consists just of the single vector $(2,2)$). 

\medskip
\begin{theorem}
  \label{theoremfour}
  Any rational pr-vector with boundary distance at least $2$ is a BTh-vector.
\end{theorem}

\medskip
Let us present again the case $\mathbb T(8)$. The pr-vectors
with boundary distance at least 2 lie in the dark region bounded by the triangle $\Delta_2.$ 

$$  
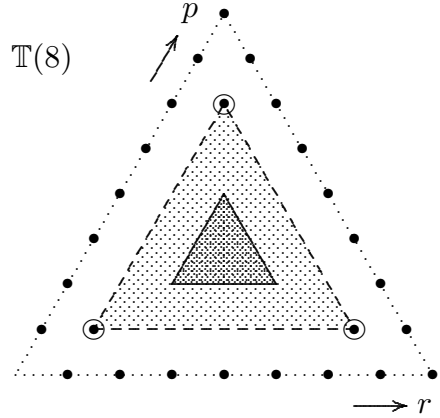
{\beginpicture
   \setcoordinatesystem units <1cm,1cm>
   \setcoordinatesystem units <.2887cm,.5cm>
\multiput{} at -8 1  8 4 /
\setdots <1mm>
\plot -8 0  8 0  0 8  -8 0 /
\multiput{$\ss \bullet$} at -7 1  -6 2  -5 3  -4 4  -3 5  -2 6  -1 7  0 8
    7 1  6 2  5 3  4 4  3 5  2 6  1 7  
    -6 0  -4 0  -2 0  0 0 
    6 0  4 0  2 0  8 0
     -5 1  5 1  0 6 /
\multiput{$\ss \bigcirc$} at  -5 1  5 1  0 6 /
\setdashes <1mm>
\plot -5 1  5 1  0 6  -5 1 /
\setshadegrid span <.5mm>
\vshade -4.8 1.2 1.2 <z,z,,> 
        -4.6 1 1.4 <z,z,,>
        -0.2 1 5.8  <z,z,,> 
        0.2 1 5.8  <z,z,,> 
        4.6 1 1.4 <z,z,,>
        4.8 1.2 1.2 /
\put{$\mathbb T(8)$} at -7 7 

\setsolid
\plot -2 2  2 2  0 4  -2 2 /
\setshadegrid span <.3mm>
\vshade -2 2 2 <z,z,,> 
         0 2 4 <z,z,,>
         2 2 2 /
\setsolid 
\arr{5 -.7}{7 -.7}
\put{$r$} at 7.7 -.7
\arr{-2.8 6.5}{-1.8 7.5}
\put{$p$} at -1.3 8

\endpicture}
$$

\medskip
The proof of Theorem \ref{theoremfour} will provide further information on the objects which we
are able to construct, namely, we will show that required indecomposable objects 
can be constructed by starting with objects which are rather similar to those in 
the one-parameter family of indecomposable objects $X = (U,V)$ in $\Cal S(6)$ with 
$V = [6,4,2]$ mentioned above.
	
In addition to the BTh-vectors provided by Theorem \ref{theoremfour}, 
already for $n = 7,$ there are further BTh-vectors see Section \ref{sec-eight};
they have boundary distance smaller than 2. 
For $n \ge 9,$ there are even BTh-vectors with boundary distance equal to 1,
see Section~\ref{sec-seven}.
For any $n$, we present a possible guess how the region of BTh-vectors may look like, see
Section~\ref{sec-sixteen-five}.

\medskip
\begin{remarkfour}
  \label{remarkfour}
  In the present paper only vectors in $\mathbb T(n)$ with rational coordinates
  play a role, since we restrict the attention to finite-dimensional vector spaces. 
  If one considers (as one should!) also infinite-dimensional vector spaces, then it will be
  useful to remember that any vector in $\mathbb T(n)$ with boundary distance at least 2 
  (and with arbitrary 
  real, not necessarily rational coordinates) is an accumulation point of BTh-vectors.
\end{remarkfour}

\subsection{Central half-lines.}
\label{sec-one-six}

Let us assume now that $n\ge 6,$ and let us draw the attention
to the central half-lines in $\mathbb T(n)$, these are the half-lines starting in $z(n)$.
Note that a central line (i.e.\ a line which passes through $z(n)$) 
is the union of two central half-lines; they are said to be 
{\it complementary.} 
The relevance of central half-lines when studying $\Cal S(n)$ is based on the following
Theorem. 

\medskip
\begin{theorem}
  \label{theoremfive}
  Let $n\ge 6.$ Let $X$ be indecomposable.
  Then there is a (uniquely determined)
  indecomposable object $X^+$ and a positive integer $\beta X$ (also uniquely determined)
  such that 
  $$
  \frac{u|w}{b}X^+ = \frac{u|w}{b}X + \beta X\cdot \frac{n|n}{3},
  $$
  and such that there is a sectional path from $X$ to $X^+$ of length $6.$
\end{theorem}

\medskip
In particular, $X$ is central if and only if $X^+$ is central. If $X$ is not central, then $X^+$ 
belongs to the central half-line $H$ through $X$. 
Starting with $X$, Theorem \ref{theoremfive}
provides countable many indecomposable objects which lie on $H$.

\medskip
If $\Cal X$ is a class of indecomposable objects in $\Cal S(n)$,
the {\it half-line support} of $\Cal X$ is defined to be the set of central 
half-lines $H$ in $\mathbb T(n)$ which have the property that there is at least
one non-central object in $\Cal X$ which lies on $H.$

\medskip 
\begin{corollaryfive}
  \label{corollaryfive}
  Let $n\ge 6.$ Assume that $H$ belongs to the half-line support of 
  an Auslander-Reiten component $\Cal C$.
  Then there are infinitely many vertices in $H$ which
  support indecomposable objects in $\Cal C$ and which converge to $z(n)$.
\end{corollaryfive}

\medskip
Namely, assume that the central half-line $H$ supports the non-central object $X$ in $\Cal C$. 
Starting with $X$, 
Theorem \ref{theoremfive} provides inductively an infinite sequence $X^{+t}$ 
indexed by $t\in \mathbb N_0$ and defined by $X^{+0} = X$ and
$X^{+(t+1)} = (X^{+t})^+.$ The objects $X^{+t}$ have pairwise different pr-vectors, all 
lie on the central half-line $H$ and the corresponding pr-vectors converge to $z(n)$.
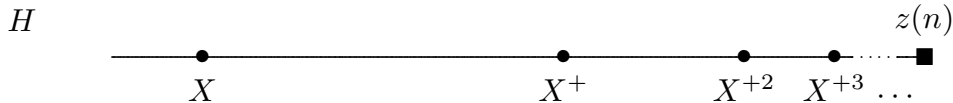
$$
{\beginpicture
    \setcoordinatesystem units <1cm,.4cm>
\multiput{} at -1 0  8 0  /
\plot -1 0  7.2 0 /
\plot 7.7 0  8 0 /
\setdots <1mm>
\plot 7.2 0  7.8 0 /
\multiput{$\bullet$} at 0 0  4 0  6 0  7 0  /
\put{$\ssize\blacksquare$} at 8 0 
\put{$X$\strut} at 0 -1
\put{$X^+$\strut} at 4 -1
\put{$X^{+2}$\strut} at 6 -1
\put{$X^{+3}$\strut} at 7 -1
\put{$\cdots$\strut} at 7.7 -1.1
\put{$z(n)$} at 8 1
\put{$H$} at -2 1 
\endpicture}
$$

\medskip 
{\bf Addition to Theorem \ref{theoremfive}.} As we have mentioned, 
the object $X^+$ as well as the number $\beta X$ provided for $X$ in Theorem \ref{theoremfive} 
are uniquely determined by $X$. In addition, as we will see,
$\beta X$ only depends on the Auslander-Reiten component $\Cal C$
to which $X$ belongs and will be called the {\it weight} of
$\Cal C$. Since with $X$ also $X^+$ belongs to the component $\Cal C$, we have 
$\beta X^+ = \beta X.$
As a consequence, the infinite sequence $X^{+t}$ given by the Corollary is very well behaved
and will be called an {\it arithmetical sequence.}

\subsection{The half-line support and the triangle support of a component.}
\label{sec-one-seven}

For any class $\Cal X$ of indecomposable objects in $\Cal S(n)$,
we have defined in Section~\ref{sec-one-six} its half-line support. Let us also define its triangle support.
Both together seem to shed a lot of light on such a class $\Cal X$.
The {\it triangle support of $\Cal X$}
is the union of the standard triangles $\Delta_d$ in $\mathbb T(n)$ which support at least one 
non-central object in $\Cal X$.  We write
$$\Psi(\Cal X)=\{0\leq d< n/3 \mid
   \Delta_d \;\text{is contained in the triangle support of}\;\Cal X\}$$
for the index set and use the abbreviation $\Psi=\Psi(\Cal S(6))$.
(By definition, both the half-line support as well as the triangle support
of a set of central indecomposables is empty.) 

We will consider now the half-line support and the triangle support for 
an Auslander-Reiten component in $\Cal S(n)$, with $n\ge 6.$ In Section~\ref{sec-one-eight}, we will consider 
the half-line support and the triangle support for the whole category $\Cal S(6).$ 

\bigskip
Let us recall the following facts, see \cite{RS1,RS2}.
Since the category $\Cal S(n)$ has
an Auslander-Reiten translation, we may look at the Auslander-Reiten quiver of $\Cal S(n)$.
The Auslander-Reiten quiver is connected in case $\Cal S(n)$ has finite type,
thus if $n\le 5.$ If $n\ge 6$, the Auslander-Reiten components are tubes of rank 1, 2, 3 or 6.
For any $n\ge 1$, 
there is a unique component $\Cal P(n)$ which is not stable, it is the component which contains the
simple object $S = (0,[1])$; we call $\Cal P(n)$ the {\it principal component}.
The isomorphism classes of the projective pickets $(0,[n])$ and $([n],[n])$
are the projective as well as the injective vertices of $\Cal P(n)$.
If we delete the projective vertices from $\Cal P(n)$, we obtain 
a stable component;
for $n\ge 6$, this stable component is a tube of rank $6$. 
We assume in Section~\ref{sec-one-seven} from now on that $n\ge 6.$ 

\medskip
\begin{theorem}
  \label{theoremsix}
  Let $\Cal C$ be an Auslander-Reiten component in $\Cal S(n)$ where $n\ge 6.$

  \begin{itemize}[leftmargin=3em]
    \item[{\rm (a)}] 
In case $\Cal C$ is stable, the half-line support of $\Cal C$ is
the union of at most $15$ pairs of complementary central half-lines.
The half-line support of $\Cal P(n)$ is the union of $18$ or $24$
half-lines, for $n = 6$, or $n \ge 7,$ respectively, and always precisely $12$
of these half-lines form complementary pairs.

\item[{\rm (b)}]  The triangle support of $\Cal C$ is the union of triangles
$\Delta_d$
with $d\in \Psi(\Cal C).$
If $\Cal C$ has at least one non-central object (so that $\Psi (\Cal C)$ 
is not empty), then
$n/3$ is an accumulation point of $\Psi(\Cal C)$, and the only one.

\item[{\rm (c)}]  Any non-central vector in $\mathbb T(n)$ is the support of at most $11$
objects in $\Cal C$. 
The center $z(n)$ is the support of countably many
objects in $\Cal C.$
  \end{itemize}
\end{theorem}

\bigskip
Some comments. 

\medskip
{\bf The half-line support of a component.}
The bound $15$ in the first assertion is 
optimal for $n = 7$, see Section~\ref{sec-ten-two};
but, as we will see in Theorem~\ref{theoremeight}, it is not optimal 
for $n =6$: According to Theorem~\ref{theoremeight}, the category $\Cal S(6)$ is supported by just 12
central lines, thus any component of $\Cal S(6)$ is supported by at most 12 lines. 

The precise description of the half-line support of $\Cal P(n)$ is as follows:
First of all, there are 
the central half-lines $\mathbb P(n)$ which are
parallel to the boundary lines. Then 
there are the central half-lines  $\mathbb D(n)$
contained in the diagonal lines.
These
are the 6 pairs of complementary half-lines in the half-line support of $\Cal P(n)$. 
Then, for all $n\ge 6,$ there is the central half-line which supports $S$;
using reflections and rotations this yields 6 half-lines $\mathbb H_\ell(n)$.
If $n \ge 7,$ the additional 6 half-lines $\mathbb K_s(n)$ are
obtained from the central half-line which supports $E_2^{n-2}$ using reflections and rotations
(the global space of $E_2^{n-2}$ is $[n-1,1],$ the subspace is $[2]$, the factor space $[n-2]$);
for $n = 6,$ the object $E_2^{n-2}$ is contained in $\mathbb P(6)$.
The half-line support of $\Cal P(6)$ is 
$\mathbb P(6)\cup\mathbb D(6)\cup\mathbb H_\ell(6)$, and the half-line support of $\Cal P(n)$ 
for $n\ge 7$ is 
$\mathbb P(n)\cup\mathbb D(n)\cup\mathbb H_\ell(n)\cup\mathbb K_s(n)$.

Let us stress the following difference between $\Cal P(n)$ and the
stable components: The half-line support of a stable component is the union of lines,
whereas there are half-lines in the half-line support of $\Cal P(n)$
such that the complementary ones are not contained in the half-line support.

Note that the non-central pickets in $\Cal S(6)$ are
supported by $\mathbb P(6)\cup\mathbb D(6)\cup \mathbb H_\ell(6)$ 
and this is the half-line support of $\Cal P(6).$ 
The half-line support of the complete category $\Cal S(6)$ will be discussed
in Section~\ref{sec-eleven}.

\bigskip
{\bf The triangle support of a component $\Cal C$.} 
For a stable component $\Cal C$, the set $\Psi(\Cal C)$ is given by the numbers of the form
$d = (u+6\beta t)/(b+3\beta t)$, where $\beta$ is the weight of $\Cal C$, with $(u,b)$
belonging to a set of at most 10 pairs. 

\medskip 
Now we consider $\Cal C = \Cal P(n).$
The increasing sequence $\Psi (\Cal P(n))$ is a subset of the interval $[0,n/3[$.  
The elements of $\Psi (\Cal P(n))$ can be distributed into 5 subsets of the form $(u+6t)/(b+3t)$ with $t\in \mathbb N_0$ and $(u,b)$ being one of the following 5 pairs, see Section~\ref{sec-ten-six}:
$$
 (0,1),\ (1,1),\ (2,2),\  (n\!-\!2,3),\  (n\!-\!1,3).
$$

\medskip
    {\bf Number of indecomposables with same pr-vector.}
    We show in Proposition~\ref{prop-nine-fourteen}
    that any
    standard triangle is the support of at most 33 objects in $\Cal C$ (and the number
33 is optimal). This
implies that any non-central vertex in $\mathbb T(n)$ is the support of at most 11
objects in $\Cal C$ (but the number 11 does not seem to be optimal).

\subsection{The case $n = 6.$}
\label{sec-one-eight}

As we have mentioned, the categories $\Cal S(n)$ with $n\le 5$ have finite type, whereas
those with $n\ge 6$ have infinite type.
The border case $n = 6$ is clearly of special interest. We have devoted
a lot of energy in \cite{RS1, S2} in order to describe the category $\Cal S(6)$ explicitely. 
It is worthwhile to have a further look at $\Cal S(6)$, now taking into account the 
invariant $b$. First, we will deal with bounds which concern $b$ (similar to the known one
for $u$ and $w$), 
then we will use 
the bounds for $u,w,b$ in order to establish several finiteness
results which are similar to the finiteness results shown in Theorem \ref{theoremsix} for single
components. 

\medskip
{\bf Bounds.}
A basic result of our previous paper \cite{RS1} asserts that for $X$ indecomposable in $\Cal S(6)$,
one has $|vX-2uX| \le 6,$ or, equivalently, that $|wX-uX| \le 6.$ We have stressed several times
that the aim of the present paper is to invoke besides the invariants $uX$ and $wX$ also
the width $bX$. There are the following three inequalities which relate the width $b$ to
$u,\ v$ and $w.$

\medskip
\begin{theorem}
  \label{theoremseven}
  If $X$ is indecomposable in $\Cal S(6)$, then there are the following inequalities
$$  
   |uX -2\cdot bX| \le 4,\quad |vX -4\cdot bX| \le 4,\quad |wX -2\cdot bX| \le 4.
$$ 
   for all indecomposable objects $X$ in $\Cal S(6)$.
\end{theorem}

\medskip
These inequalities are optimal: For example, let us look at the first inequality:
The objects $X = P[6t]$ with $t\ge 0$ have $uX = 6(t+1)$
and $bX = 3t+1,$ thus $uX-2\cdot bX = (6t+6)- 2(3t+1) = 4.$
The object $X = ([3,1],[6,4,3,1])$ has $uX = 4 = bX,$ thus $uX-2\cdot bX = -4.$
(If we consider only the objects in the principal component $\Cal P(6)$, we get
better inequalities. Let us add that we have 
$-2 \le uX - 2\cdot bX \le 4$ 
and $-3 \le vX - 4\cdot bX \le 2$ 
for $X$ in $\Cal P(6),$ see Section~\ref{sec-eleven-four}.)
	
There are indecomposable objects with $u = 2\cdot bX - 4,$ and 
also indecomposable objects with $w = 2\cdot bX - 4,$ but the middle inequality in Theorem \ref{theoremseven}
asserts that we cannot have both conditions at the same time!
	
The known symmetries assert
that the three inequalities are equivalent (thus only one of them has to be shown, 
this will be done in Section~\ref{sec-eleven-four}). Namely, the 
duality $\D$ shows that the first and the last inequality are obviously equivalent.
In order to show that the first inequality implies the second, let us assume that 
$|uY -2\cdot bY| \le 4$ for all indecomposable objects $Y$ in $\Cal S(6).$ 
Let $X$ be indecomposable in $\Cal S(6)$. If $X$ is projective, then
$vX = 6$ and $bX = 1,$ thus $vX -4\cdot bX = 2.$ 
Thus, we may assume that $X$ is not projective. Then $X = \tau^2Y$ for some indecomposable
object $Y$ in $\Cal S(6)$ and by assumption, we have $|uY -2\cdot bY| \le 4$. According to
Theorem \ref{theoremthree},  $vX = v(\tau^2Y) = 6\cdot bY - uY $ and $bX = b(\tau^2Y) = bY$.
Therefore $|vX -4\cdot bX| =  |6\cdot bY - uY-4\cdot bY| = |uY-2\cdot bY| \le 4.$
A similar argument shows that the second inequality implies the first one.
$\s$

\bigskip
We may divide the inequalities by $bX$ in order to get 
the following estimates for $pX-2,$ for $qX-4,$ and for $rX-2$. 

\medskip
\begin{corollaryseven}
  \label{corollaryseven}
  If $X$ is indecomposable in $\Cal S(6)$, then there are the following inequalities
  $$  
  |pX -2| \le \frac4{bX},\quad |qX -4| \le \frac4{bX},\quad |rX -2| \le \frac4{bX}
  $$ 
  for all indecomposable objects $X$ in $\Cal S(6)$. $\s$
\end{corollaryseven}

\medskip
{\bf Finiteness results.} 
Our previous paper \cite{RS1} has provided a complete classification of the indecomposable
objects in $\Cal S(6)$, using three kinds of invariants: rational numbers 
(called the index, or, better, the rationality index), 
the elements of the projective line $\mathbb P^1$ (as parameter inside a tubular family), 
as well as combinatorial
ones (in order to describe the position in a tube). Our present approach to look at 
the distribution of the pr-vectors is less ambitious, but it shows its power even 
in the case $n = 6$, see the following Theorem \ref{theoremeight}.

Since we consider here just the possible pr-vectors,
we cannot expect any new information concerning central objects.
Whereas the classification given in \cite{RS1} may give the impression that $\Cal S(6)$ has
an abundance of indecomposable objects, Theorem \ref{theoremeight} focuses the attention to several 
strong finiteness properties of $\Cal S(6)$. We denote by $\mathbb L(6)$ the
union of $\mathbb P(6),\ \mathbb D(6),\ \mathbb H(6)$, where $\mathbb H(6)$ is obtained from
$\mathbb H_\ell(6)$ by adding the complementary half-lines. Thus 
$\mathbb L(6)$ is the union of 12 pairs of complementary half-lines. 

\medskip 
\begin{theorem}
  \label{theoremeight}
  \phantom{m.}
  \begin{itemize}[leftmargin=3em]
    \item[{\rm (a)}] The half-line support of $\Cal S(6)$ is $\mathbb L(6)$, thus the union of
      $12$ pairs of complementary central half-lines.

    \item[{\rm (b)}]
      Write the triangle support of $\Cal S(6)$ as the disjoint union of the 
      triangles $\Delta_d$ with $d \in \Psi$.
      Then this index set $\Psi$ is an increasing sequence
      of rational numbers $0 \le d < 2$ which converge to $2.$

      \item[{\rm (c)}] Any non-central vector in $\mathbb T(6)$ is the
support of a finite number of indecomposable objects in $\Cal S(6)$;
these numbers are {\bf not} bounded. 
The center $z(6)$ is the support of infinitely many indecomposable objects of $\Cal S(6)$,
the number is $\max(\aleph_0,|k|).$
  \end{itemize}
\end{theorem}

\medskip
The last assertion in (c) follows directly from Theorem \ref{theoremfive}. The remaining assertions are
established in Section~\ref{sec-eleven}.
Here are some comments on Theorem \ref{theoremeight}.

\bigskip
{\bf The half-line support of $\Cal S(6)$.} According to Theorems \ref{theoremsix} and \ref{theoremeight}, 
we obtain the half-line support of $\Cal S(6)$ 
from  the half-line support $\mathbb P(6)\cup \mathbb D(6)\cup \mathbb H_\ell(6)$ of $\Cal P(n)$ 
by adding the half-lines complementary to the half-lines in $\mathbb H_\ell(6).$ 
As we have mentioned, all pickets of $\Cal S(6)$ live on
$\mathbb P(6)\cup \mathbb D(6)\cup \mathbb H_\ell(6)$, this holds also true for the bipickets.
The smallest non-central objects in $\Cal S(6)$
which live on $\mathbb H(6)\setminus \mathbb H_\ell(6)$ 
have width 3, they are exhibited in Section~\ref{sec-eleven-seven}, see also Section~\ref{sec-fifteen-two}--(c). 

A central line $L$ in $\mathbb T(6)$ can be described by specifying its slope 
$\phi\in \mathbb R\cup \{\infty\},$ namely $L = L_\phi,$
where $L_\phi = \{(p,r)\in \mathbb T(6)\mid (u-2b)/(w-2b) =
\phi\}$. The set 
$\mathbb L(6)$ is the set of the lines $L_\phi$ with $\phi\in \Phi$:
$$
 \Phi = \{0,\ \tfrac12,\ 1,\ 2,\ \infty,\ -3,\ -2,\ -\tfrac32,\ -1,\ 
   -\tfrac23,\ -\tfrac12,\ -\tfrac13 \},
$$    
(The 12 lines in $\mathbb L(6)$ can be written in a unified way if we refer 
also to the syzygy module $\Omega V$ of $V$ in $\Cal N(n).$
Namely, let $\omega = |\Omega V|$. Using the functions $u, w, \omega$, all
the 12 lines are given by equations of the form
$u = 2b$ (one of the functions is constant, namely equal to $2b$),
$u = w$ (two of the functions coincide), 
and finally $u = 2(w-b),$ see Section~\ref{sec-eleven-seven}.)

\medskip
{\bf The triangle support.}
The increasing sequence $\Psi$ is a subset of the interval $[0,2[$.  
Here are the first numbers in $\Psi$:
$$
 0,\ 1,\ \frac54,\ \frac43,\ \frac 75,\ \frac {10}7,\,
\frac 32,\ \frac{11}7,\ \frac85,\ \frac{13}8,\
 \frac 53,\ \cdots.
$$
All elements of $\Psi$ are of the form $2-c/m$ with $c = 3$ or $c= 4$ and $m\in \mathbb N_1$.
They can be distributed into 10 subsets of the form $(u+6t)/(b+3t)$ with $t\in \mathbb N_0$
and $(u,b)$ being one of the following pairs
$$
 (0,1),\ (1,1),\ (2,2),\ (3,2),\ (3,3),\ (4,3),\ (4,4),\ (5,3),\ (5,4),\ (7,5).
$$

\medskip
Theorem \ref{theoremeight} provides for $n = 6$ a generalization of an essential
part of Theorem \ref{theoremone}. 
Namely, Theorem~\ref{theoremone}, Theorem~1.1$'$ (in Section~\ref{theoremone-dual}),
and Theorem~1.1$''$ (in Section~\ref{theoremone-rot}) assert that there is
no indecomposable object with support in $\Delta_d$, where $d$
belongs to the open interval $]0,1[$. Using Theorem \ref{theoremeight}, we obtain 
countably many open intervals
with this property. We do not know whether a corresponding assertion may exist
for $7\le n \le 9.$
For $n\geq 10$, Corollary in Section~\ref{sec-eight-five} shows that 
for any rational number 
$d$ with $1\le d \le \frac n3,$ there is an indecomposable object $X$ in
$\Cal S(n)$ with
boundary distance $d$.

\bigskip 
{\bf The intersections of central half-lines and standard triangles.}
According to Theorem \ref{theoremeight}, any non-central indecomposable object of $\Cal S(6)$
lives on the intersection of a half-line in $\mathbb L(6)$ and a triangle
$\Delta_d$ with $d\in \Psi$.  However, not all intersection vertices arise in this way:
For all values $d\in \Psi$, the set 
$\mathbb L(6)\cap \Delta_d$ has cardinality 24, however, according to Theorem \ref{theoremone},
there are only 17 indecomposable objects with pr-vector in $\Delta_0.$ 
For $d =1,$ all vertices of $\mathbb L(6)\cap \Delta_1$ 
support indecomposable objects, see Section~\ref{sec-fifteen-two}.
On the other hand, for $d = \frac54,$ the set $\Delta_d$
is the support of only few indecomposables, see Section~\ref{sec-fifteen-two}.

\bigskip
{\bf The essence of Theorem \ref{theoremeight}.} 
Altogether, Theorem \ref{theoremeight}
combines {\bf three} different finiteness assertions for the
set of non-central indecomposable objects in $\Cal S(6),$ similar to Theorem \ref{theoremsix}. 
First of all, there are only $12$ central lines which support non-central
indecomposable objects. Second, given $a < 2,$ there are 
only finitely many numbers $d \le  a$ such that $\Delta_d$ supports
an indecomposable object. And third, for any vector $(p,r)\neq (2,2)$ in $\mathbb T(6)$,
there are only finitely many indecomposable objects with pr-vector $(p,r).$ But note that
all these finiteness conditions concern the non-central indecomposable 
objects of $\Cal S(6)$:  The central ones are collected in the black square, and
this black square may be considered as a sort of black hole of the category.

\bigskip
Let us present the pr-triangle $\mathbb T(6)$ in two ways: once together 
with the 12 lines $L_\phi$ (left), once with 
some of the triangles $\Delta_d$ (right). 
As before, 
the center $z(6) = (2,2)$ of $\mathbb T(6)$ will be marked by 
a black square $\ssize\blacksquare$, 
it is the support of the unique central picket $([2],[4])$; the support of the 
remaining pickets will be shown as small bullets $\bullet$.

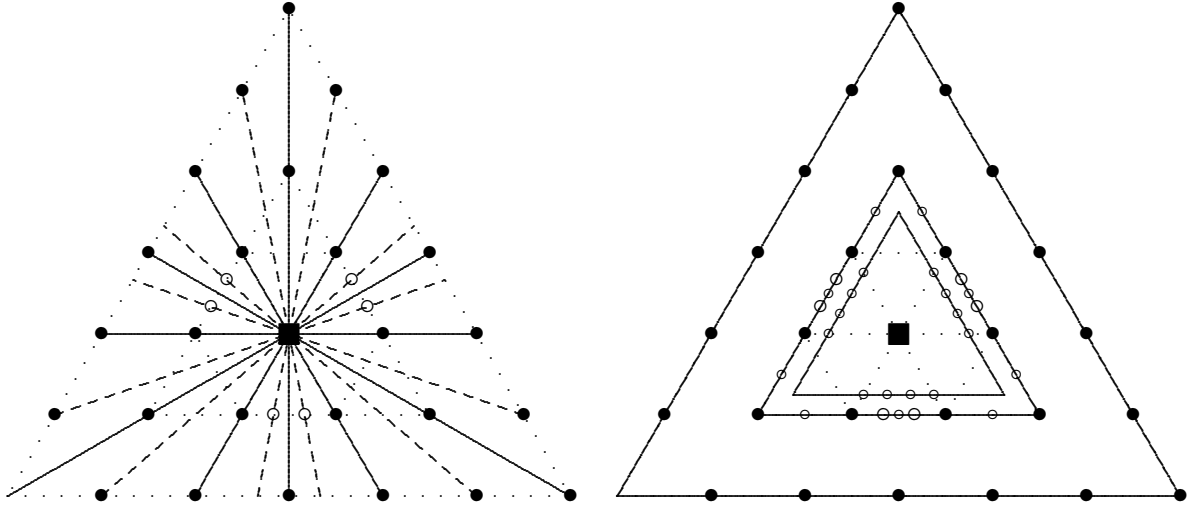
$$
{\beginpicture
   \setcoordinatesystem units <.462cm,.8cm>
\put{\beginpicture
\multiput{} at -6 -.5  6 6 /
\setdots <2mm> 
\plot -6 0  6 0  0 6  -6 0 /
\plot -3 1  3 1  0 4  -3 1 /
\plot -1 1  1 3  -1 3  1 1  2 2  -2 2  -1 1 /

\multiput{$\bullet$} at 
     -4 0  -2 0  0 0  2 0  4 0  
     -3 1  -1 1  1 1  3 1 
     -2 2  0 2  2 2
     -1 3  1 3
      0 4
       -5 1  -4 2  -3 3  -2 4  -1 5  1 5  2 4  3 3  4 2  5 1  /
\multiput{$\circ$} at -.333 1  .333 1  1.333 2.667  1.667 2.333  -1.333 2.667  -1.667 2.333 /
\put{$\blacksquare$} at 0 2
\setsolid
\plot 0 0  0 6 /
\plot -6 0  3 3 /
\plot  6 0  -3 3 /

\plot -2 0  2 4 /
\plot 2 0  -2 4 /
\plot -4 2  4 2 /

\setdashes <1mm>
\plot  -4 0  0 2  /
\plot   4 0  0 2 /
\plot -5 1  0 2 /
\plot  5 1  0 2 /
\plot 0 2  1 5 /
\plot  0 2  -1 5 /

\setdashes <1mm>
\plot  0 2  2.667  3.33 /
\plot   0 2  -2.667  3.33 /
\plot 0 2  3.333 2.667 /
\plot  0 2  -3.333 2.667 /
\plot -0.667 0  0 2 /
\plot  0.667 0  0 2 /

\setsolid
\put{$\bullet$} at  6 0  
\put{$\bullet$} at 0 6

\endpicture} at 0 0
\put{\beginpicture
\multiput{} at -6 -.5  6 6 /
\setdots <2mm> 
\plot -6 0  6 0  0 6  -6 0 /
\plot -3 1  3 1  0 4  -3 1 /
\plot -1 1  1 3  -1 3  1 1  2 2  -2 2  -1 1 /

\multiput{$\bullet$} at 
     -4 0  -2 0  0 0  2 0  4 0  
     -3 1  -1 1  1 1  3 1 
     -2 2  0 2  2 2
     -1 3  1 3
      0 4
       -5 1  -4 2  -3 3  -2 4  -1 5  1 5  2 4  3 3  4 2  5 1  /
\multiput{$\circ$} at -.333 1  .333 1  1.333 2.667  1.667 2.333  -1.333 2.667  -1.667 2.333 /
\put{$\blacksquare$} at 0 2

\setsolid
\plot -6 0  6 0  0 6  -6 0 /
\plot -3 1  3 1  0 4  -3 1 /
\plot -2.25 1.25  2.25 1.25  0 3.5  -2.25 1.25 /

\put{$\bullet$} at  6 0  
\put{$\bullet$} at 0 6
\multiput{$\ssize \circ$} at  -2 1  2 1 /
\multiput{$\ssize \circ$} at  
    -2.5 1.5  2.5 1.5  -.5 3.5  .5 3.5 /
\multiput{$\ssize \circ$} at -.25 1.25  .25 1.25   1.25 2.25  1 2.5  -1.25 2.25  -1 2.5 /
\multiput{$\ssize \circ$} at  0 1  -1.5 2.5  1.5 2.5  /

\multiput{$\ssize \circ$} at  -.75 1.25 .75 1.25  1.5 2  -1.5 2   0.75 2.75   -.75 2.75 /
\endpicture} at 13 0
\endpicture}
$$

To repeat: In the picture on the left, we have shown the 12 central lines $L_\phi$ in $\mathbb L(6)$ 
which contain a boundary picket (for the values $\phi$, see Section~\ref{sec-eleven-seven}). 
These 12 central lines are drawn in two different ways: The diagonal lines $\mathbb D(6)$
and the lines $\mathbb P(6)$ 
parallel to the boundary are drawn as solid lines, the 
remaining half lines in $\mathbb H(6)$
are dashed (for the half-lines in 
$\mathbb H(6)\setminus\mathbb H_\ell(6)$, we have indicated by a small circle $\circ$ the position of the smallest non-central object which lives on this half-line; for the shape of these objects,
see Section~\ref{sec-eleven-seven}).

In the picture on the right, we have drawn the sets $\Delta_d$ for the first
three values $d = 0,\ 1,\ \frac 54$ and the position of the indecomposable
objects living on these triangles $\Delta_d$ (see Section~\ref{sec-fifteen-two} and Appendix~\ref{app-B}).
Of course, if we insert
further sets $\Delta_d$, we obtain a nested sequence of triangles with size 
converging to zero. 

\bigskip
{\bf The relevant
\phantomsection{hexagon}
\addcontentsline{lof}{subsection}{The relevant hexagon in $\mathbb T(6)$.}%
in $\mathbb T(6)$.} Let us draw the attention 
to the convex hull of the $\Sigma_3$-orbit of the pr-vector
$(1,2)$; this is the hexagon and its interior, and is shaded in the picture below.  
Outside of the hexagon live only pickets and bipickets (the 17 boundary pickets, as well
as 3 pickets and 9 bipickets with boundary distance 1, their position is marked 
by two kinds of circles). 

The hexagon is the support of 18 indecomposable
objects with width 3, as well as of six indecomposable objects with width 4. The six objects 
with width 4 are supported by the midpoints of the sides of the hexagon, thus also by
$\mathbb D(6)$ (they will be discussed for example in Section~\ref{sec-eleven-eleven}).
Three of them belong to the principal component
$\Cal P(6)$, thus they have rationality index $0$; the remaining three 
have rationality index 1. The position of these 24 indecomposable objects with width 3 and 4
is marked by small bullets 
(actually, some of these positions support in addition also pickets and bipickets).

Looking at the hexagon, we see nicely the direction of the 12 lines in 
$\mathbb L(6) = \mathbb P(6)\cup \mathbb D(6)\cup \mathbb H(6).$ 
Namely, as we have mentioned already, the midpoints of the
sides of the hexagon support objects in $\mathbb D(6)$. The six
corners of the hexagon support objects in $\mathbb P(6)$,  The remaining 12 marks
on the hexagon support objects in $\mathbb H(6).$ 

$$
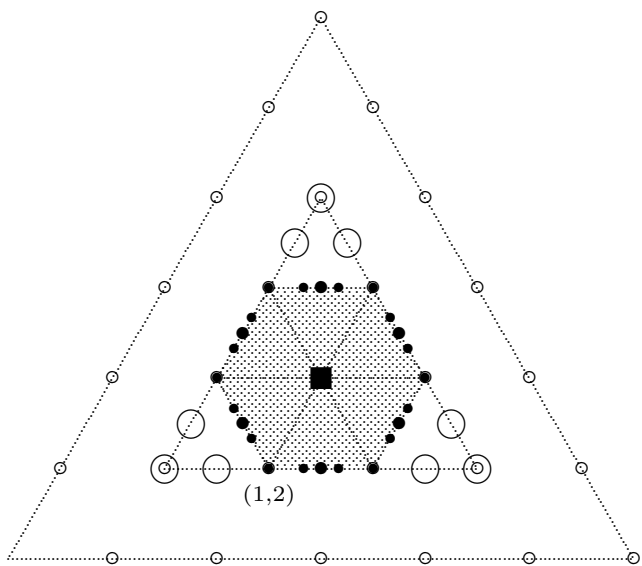
{\beginpicture
   \setcoordinatesystem units <.57735cm,1cm>
\multiput{} at -6 -.5  6 6 /
\setdots <.5mm> 
\plot -6 0  6 0  0 6  -6 0 /
\plot -3 1  3 1  0 4  -3 1 /
\plot -1 1  1 3  -1 3  1 1  2 2  -2 2  -1 1 /

\multiput{$\circ$} at 
     -4 0  -2 0  0 0  2 0  4 0  
     -3 1  -1 1  1 1  3 1 
     -2 2  0 2  2 2
     -1 3  1 3
      0 4   6 0   0 6
       -5 1  -4 2  -3 3  -2 4  -1 5  1 5  2 4  3 3  4 2  5 1  /
\multiput{$\bigcirc$} at -.5 3.5  .5 3.5  -2 1  2 1  -2.5 1.5  2.5 1.5 /
\multiput{$\bigcirc$} at -3 1  3 1  0 4  /
  
\multiput{$\ssize\bullet$} at -.333 1  .333 1  1.333 2.667  1.667 2.333  -1.333 2.667  -1.667 2.333 
     -2 2  -1 1  -1 3  1 1  1 3  2 2 -.333 3  .333 3 
      1.333 1.333   1.667 1.667  -1.333 1.333  -1.667 1.667 /
\multiput{$\bullet$} at 0 1  0 3  -1.5 1.5  1.5 1.5  -1.5 2.5  1.5 2.5  /
\put{$\blacksquare$} at 0 2
   \setshadegrid span <.4mm>
   \vshade -2 2 2  <,z,,>  -1 1 3   <z,z,,>  1 1 3   <z,z,,>  2 2 2  /
\put{$\ssize (1,2)$} at -1 0.7
\endpicture}
$$
As we have seen, there are just 83 indecomposables
which live on or outside of the hexagon. Thus, all
the interesting features of $\Cal S(6)$ happen inside
the hexagon.

\subsection{Optimal results? In no way!}
\label{sec-one-nine}

Our paper is rather long, together with \cite{RS1} and \cite{RS2}
there are well over 200 pages. 
But the results presented by us are far away from being optimal. For example, when we deal
with BTh-vectors, we focus the attention to a typical tame behaviour, in spite of being
nearly always in a wild realm. 
	
The reader will notice the great number of open questions which we have selected
in Section~\ref{sec-sixteen}.
And many of these questions are really what one may call ``Ziegler questions'': 
questions, which are completely elementary to formulate, which are easy to grasp even 
by a first year undergraduate student, but which seem to be quite difficult to solve. 
In our opinion, it is important to mention such problems to students in order to make them aware
of the complexity of mathematics, but also to advertise possible research topics. 

\subsection{The general context.}
\label{sec-one-ten}

Let us add some remarks regarding the context of the present paper.

\medskip
{\bf The Birkhoff problem.}
The Birkhoff problem concerns subgroups of finite abelian groups: To describe
the possible embeddings, up to automorphisms of the global group.
It is easy to see that one may restrict to deal with
subgroups of finite abelian $p$-groups, thus with submodules of $\mathbb Z/p^n$-modules
for some $n$. Clearly, there is no problem to describe the isomorphism
classes of the subgroups of a given finite abelian group $G$, but it turns out to be
difficult to describe the possible embeddings, or better the equivalence classes of embeddings,
where two embeddings are considered as {\it equivalent} if they are obtained from each other by
an automorphism of $G$. 

Slightly more generally, one may pose the
same problem for $\Lambda$-modules, where $\Lambda$ is an arbitrary (commutative local)
uniserial ring. It is an open problem whether the Birkhoff problem for uniserial rings
in general has different answers for different uniserial rings of the same length
$n \ge 4$ (the cases $n\le 3$ have been settled in the new preprint \cite{GKKP};
for a general survey, see \cite{S3}).
For an old discussion of the relationship
between $\mathbb Z/p^n$ and $k[T]/\langle T^n\rangle$ we may refer
to Kaplansky \cite{Kap} who was eager to stress the similarity between
abelian group theory and the theory of linear transformations.

The case of $\Lambda = k[T]/\langle T^n\rangle$ is easier to attack, since one may use
covering theory. Thus our previous papers, as well as the present one, are devoted to
submodules of $\Lambda$-modules, where $\Lambda = k[T]/\langle T^n\rangle$, and we hope
that any result concerning the Birkhoff problem for $\Lambda = k[T]/\langle T^n\rangle$
may help to understand the Birkhoff problem for a general uniserial ring $\Lambda$,
and thus for $\Lambda = \mathbb Z/p^n$.

Given a ring $\Lambda$ and any $\Lambda$-module $V$,
one of the first question may concern the
non-zero elements $x\in V$:
how are they embedded into $V$. But this
means that one studies the pair $(U,V)$, where $U$ is the cyclic
$\Lambda$-module generated by $x$.
Corresponding results of Pr\"ufer \cite{P} and Kaplansky \cite{Kap} will be
recovered for $\Lambda = k[T]/\langle n\rangle$ in our Sections~\ref{sec-thirteen-two}
through \ref{sec-thirteen-four}.

We also should mention two recent preprints which are devoted to the Birkhoff problem 
and its generalizations. First of all, there is the report \cite{Kva}
by Kvamme which presents a comprehensive survey on general monomorphism categories.
In particular, it also outlines the relationship between submodule 
categories and $p$-valued abelian groups. 
Second, the paper \cite{KSS} by Kosakowska, Schmidmeier, and Schreiner
shows that a pair $(U,V)$, where $V$ is a direct sum of cyclic $p$-groups (not necessarily bounded) and $U$ a subgroup
of $V$ such that $p$ annihilates $V/U$, is a direct sum of pickets. 

\medskip
{\bf Gorenstein-projective modules.} If $\Lambda$
is a uniserial ring,
the ring  $T_2(\Lambda)$ is a Gorenstein ring of
dimension 1, thus a $T_2(\Lambda)$-module $M$ is Gorenstein-projective if and only if $M$ is torsionless
(i.e. isomorphic to a submodule of a projective $T_2(\Lambda)$-module). 
In particular, the category $\Cal S(n)$ is the category of
Gorenstein-projective $T_2(\Lambda)$-modules, where $\Lambda = k[T]/\langle T^n\rangle$.
Note that the category of $T_2(\Lambda)$-modules is just the category of homomorphisms
of $\Lambda$-modules (thus of homomorphisms between linear operators), and the Gorenstein-projective
$T_2(\Lambda)$-modules are just the embeddings of $\Lambda$-modules: We may say the
embeddings of linear operators. 

\medskip 
{\bf Operator theory.} 
The emphasis that $\Cal S(n)$ is a nice subcategory of $\mod T_2(\Lambda(n))$, 
namely the category of Gorenstein-projective modules, could hide the fact that all 
the results presented here concern properties of the 
category $\Cal N(n) = \mod\Lambda(n)$ itself, thus 
properties of the category of nilpotent operators with a bound on the nilpotence index.
To deal with the category $\Cal N(n)$ seems, on a first look, quite easy,
since certainly all the objects, the indecomposables, but also the decomposable ones, are
easy to describe. However, as soon as one asks for a description of the maps, one
encounters severe difficulties. 
This is the message of Birkhoff, when he posed the problem of
describing possible embeddings, and our aim is to follow his challenge. 
Since all the maps are concatenations of monomorphisms and
epimorphisms, any information about the possible embeddings helps to 
understand the category $\Cal N(n)$.

The relevance of a description of $\Cal S(n)$ for an understanding of the 
category $\Cal N(n)$
can be seen also in a different way, by stressing the relationship between $\Cal S(n)$
and the stable module category $\underline\mod\Lambda(n)$ (which should be seen as a
triangulated category).
Namely, 
any morphism $\overline g$ in $\underline\mod\Lambda(n)$
can be represented by a monomorphism in $\mod\Lambda(n)$, 
thus by an object in $\Cal S(n)$,
namely by $\Mimo g,$ see \cite{RS2}. 

The study of invariant subspaces is, of course, part of
operator theory. Here: dealing with a nilpotent operator. If one deals with a single
operator, it is usually not helpful to invoke the assumption that the operator is
indecomposable: The decomposable ones are those of interest.

\medskip
{\bf A basic problem in linear algebra.} To deal with invariant subspaces of nilpotent
operators should be seen as a basic problem in linear algebra. Any introductory course
in linear algebra introduces both subspaces and linear transformations, and to look at
a combined system (just one subspace, just one transformation $T$, even with the additional 
condition on $T$ to be nilpotent) should be seen as a very natural object to study.
Lecturers in linear algebra often pretend that elementary linear algebra is a subject
which is sealed --- but actually, it should be stressed that it is and will be
open-ended. A typical result which
usually is presented is the dimension formula for pairs of subspaces of a vector space.
However the corresponding description of a vector space with three given subspaces
(presented by Dedekind in 1900) is rarely discussed, the situation of four subspaces
(considered in 1967 by Nazarova and in 1970 by Gelfand and Ponomarev) presumably never. 

\medskip
{\bf The present paper.} 
Let us compare the present investigation with our previous paper \cite{RS1}.
In \cite{RS1} we focussed the attention to the invariants $u$ and $w$ 
(actually, to $u$ and $v = u+w$),
now, as we have mentioned, we stress the additional relevance of $b$.
The paper \cite{RS1} was based on the fact that for $X$ indecomposable in $\Cal S(6)$,
the numbers $uX$ and $wX$ are roughly the same (the difference is bounded by $6$).
The present paper started with the observation that for $\Cal S(6)$ 
the numbers $uX$ and $2\cdot bX$ are roughly the same (the difference is bounded by $4$),
thus $uX$ is roughly equal to $\frac n3\cdot bX$. 
This is a general feature, for all $n$: The numbers
$uX$ and $wX$ are roughly equal to one third of $n\cdot bX$.
Note that $n\cdot bX$ is the length of the projective cover $PV$ of the global space of $V$.
Whereas in \cite{RS1} we looked at the global space $V$ of an object $X$ as being an
extension of the {\bf two} vector spaces $UX$ and $WX$, we now fix the attention to the projective
cover $PV$ of the global space $V$ of $X$ as a $\Lambda$-module, 
and its filtration with the {\bf three} factors
$\Omega V,\ UX,\ WX$. The position of $\pr X$ in $\mathbb T(n)$ using triangular coordinates
is given by the dimension of these three vector spaces and $\tau^2$ provides a cyclic
rotation of these numbers. In this way, we see clearly 
the relevance of the  pr-triangle with center of gravity $z(n) = (\frac n3,\frac n3)$:
The rotation $\rho$ of the triangle $\mathbb T(n)$ is realized (for all
reduced objects) by the functor $\tau^2.$ 

\subsection{Outline of the paper.} 
\label{sec-one-eleven}

Some further definitions and basic concepts will be mentioned in Section~\ref{sec-two}.
Our investigation starts in Section~\ref{sec-three}, where we look at the Auslander-Reiten
translation $\tau$. In particular, we show Theorem~\ref{theoremthree} which
asserts that $\tau^2$ induces the rotation of $\mathbb T(n)$ by $120^\circ$. 
This result will be used frequently. 

Section~\ref{sec-four} is devoted to the proof of Theorem \ref{theoremone},
Section~\ref{sec-five} to the proof of Theorem \ref{theoremtwo}.
Both Sections~\ref{sec-four} and \ref{sec-five} deal with suitable filtrations of objects in $\Cal S$.
Such filtrations may be of further interest. Let us draw the attention
to the additional {\bf Theorem~\ref{theoremnine}} at the end of Section~\ref{sec-five}. 

Sections~\ref{sec-six}, \ref{sec-seven} and \ref{sec-eight}
are devoted to provide BTh-vectors. Section~\ref{sec-six} deals with 
pr-vectors with boundary distance at least 2 and presents a proof of Theorem \ref{theoremfour}.
Note that in Section~\ref{sec-sixteen-five}, we formulate a guess about the shape of the region of all
BTh-vectors. It is a guess, not a conjecture: After all, whereas 
we know some methods for showing that a given pr-vector is a BTh-vector, at present no
effective way seems to be available to decide that a given pr-vector is not a BTh-vector. 
Sections~\ref{sec-seven} and \ref{sec-eight}
deal with pr-vectors with boundary distance smaller than 2.
In particular, in Section~\ref{sec-seven} we draw the attention to pr-vectors with $p = 1,$
thus to the indecomposable objects $X$ in $\Cal S$ with $uX = bX$. 
We have to admit that we were surprised when we realized that already this subcategory
is wild! In view of our interest in BTh-vectors, we are going to show that 
for $n\ge 9,$ the pr-vector $(1,4)$ is a BTh-vector: There is an infinite family 
of indecomposable objects $X$ with $uX = bX = 6$ and $wX = 24$. 
In Section~\ref{sec-eight-five} one finds the additional {\bf Theorem~\ref{theoremten}} which 
provides for $n\ge 10$ many BTh-vectors with boundary distance between 1 and 2.
	
Sections~\ref{sec-nine} to \ref{sec-eleven} are devoted to the half-line and the triangle support. 
We show Theorems \ref{theoremfive}, \ref{theoremsix}, \ref{theoremseven} and
\ref{theoremeight}. Sections~\ref{sec-nine} and \ref{sec-ten}
deal with properties of $\Cal S(n)$, where $n\ge 6$ is not further specified.
Section~\ref{sec-nine} provides the proof of Theorem \ref{theoremfive}.
The discussion of 
Auslander-Reiten components will be found in Section~\ref{sec-ten}.
Special assertions which are valid only in the case $n = 6$ are discussed in Section~\ref{sec-eleven}. 
These considerations have to be seen as an expansion of those presented in \cite{RS1},
where we have discussed relations between $uX$ and  $wX$, now it is the triple
$uX,\ wX,\ bX$ which we study.
In particular, Theorem \ref{theoremseven} shows that for an indecomposable object $X$ in $\Cal S(6),$
we have $|uX- 2bX| \le 4,$  and $|vX- 4bX| \le 4.$
This is the basis for the proof of Theorem \ref{theoremeight}.
In this way, we show that for $n = 6,$ the pr-vector $(2,2)$ is the only 
BTh-vector. 
	
The remaining sections are devoted, on the one hand, to special objects, or, on the other hand,
to the weird (but actually prevalent) non-gradable objects. Some non-gradable objects are 
discussed in Section~\ref{sec-fourteen}.
	
Section~\ref{sec-twelve} deals with pickets and bipickets. 
Section~\ref{sec-thirteen} devotes a lot of
attention to the objects $X = (U,V)$ in $\Cal S$ with $bU = 1;$ here we can 
follow and extend previous investigations by Pr\"ufer and Kaplansky. 
	
Section~\ref{sec-fifteen} will be a sort of gallery: Here, 
a lot of further examples are exhibited in order
to draw the attention to special features of the categories $\Cal S(n)$. 	
In Section~\ref{sec-sixteen} we draw the attention to open questions. Indeed, dealing
with invariant subspaces, there are presently more questions than results! We add some
general comments in Section~\ref{sec-seventeen}.

Appendix~\ref{app-A} by Schmidmeier
lists the positive roots of the Dynkin diagram $\mathbf E_8$ with some additional
combinatorial data used
in the paper.  
In Appendix~\ref{app-B}, Schmidmeier presents the remaining tripickets in $\Cal S(6)$ as well as 
some tetrapickets and some pentapickets. 
In Appendix~\ref{app-C}, Schmidmeier describes briefly some applications of invariant subspaces.
Appendix~\ref{app-D} by Ringel is a 
slightly revised reprint from the Izmir lectures containing a report (with full proof)
on the covering theorem.

\medskip 
Most parts of the paper are self-contained.
We will use two standard techniques from representation theory: coverings and simplification.
Both methods are explained in Section 2. 
We will use frequently the one-parameter family of indecomposable objects 
$X = (U,V)$ in $\Cal S(6)$
where $V = [6,4,2]$. In addition, another one-parameter 
family of indecomposable objects $X = (U,V)$
will play a role: It belongs to $\Cal S(7)$ and $U$ has height $3$; it was exhibited by
Schmidmeier in \cite{S1}, see Section~\ref{sec-eight-two}.

According to \cite{AS}, the category $\Cal S(n)$ has Auslander-Reiten sequences. The 
Auslander-Reiten translation $\tau = \tau_n$ was explicitly described in \cite{RS2},
see Section~\ref{sec-three} below. 
In particular, it has been shown in \cite{RS2} that $\tau_n$ has period 6, thus stable
(Auslander-Reiten) components are tubes of rank $1,2,3$ or $6$.

The case $n=6$ is of special interest. We will use the classification of the
indecomposable objects of $\Cal S(6)$ as obtained in \cite{RS1} with the help of 
a tubular algebra, see the beginning of Section~\ref{sec-eleven}.

\medskip 
We include various remarks throughout the paper. These remarks may be skipped at a first
reading, but we hope that they are helpful for a better understanding of the setting. 

\bigskip
\section{Preliminaries.}
\label{sec-two}

We are going to provide some further definitions which we will need. In Section~\ref{sec-two-two},
we recall from \cite{RS1} 
in which way we try to visualize at least some of the objects of $\Cal S$.
Then, we  
report about methods which have been developed quite a long time ago and which will
be used throughout the paper.

\subsection{Further definitions which are needed.}
\label{sec-two-one}

Given an object $V$ in $\Cal N(n)$, let $PV$ be a projective cover of $V$ and 
$\Omega V$ the (first) syzygy module of $V$ (the kernel of an epimorphism $PV \to V$).
For $1\le m \le n-1$,
we have $\Omega [m] = [n-m]$ (and, of course, $\Omega [n] = 0$). For any $V$ in $\Cal N(n)$, 
we have $|\Omega V| = n\cdot bV -|V|.$

\bigskip
In dealing with objects $X = (U,V)$
in $\Cal S(n)$, it seems appropriate to introduce the partition triple
$\bpar X$ as
$$
 \bpar X = ([U],[V],[V/U])
$$
(recall that for $V\in \Cal N(n)$, we denote by $[V]$ its isomorphism class, or,
equivalently, the corresponding partition). Thus, $\bpar X$ is a triple of partitions.
By abuse of notation, we may write $X = (U,V,V/U)$ in case 
$\bpar X = ([U],[V],[V/U])$. Thus, if an object $X$ of $\Cal S(n)$ is presented by
a triple, then actually only its partition vector is given (and there may be several
different isomorphism classes of objects $X$ with these data). 

\medskip
As we have mentioned, if  $X$ is an object in $\Cal S(n)$, the triple
$(uX,wX,bX)$ is called the uwb-vector of $X$. Actually, since we usually are interested
in the corresponding pr-vector $(uX/bX,wX/bX)$ (provided $X$ is non-zero), 
we sometimes will write the triple $(uX,wX,bX)$ in the form
$\frac{u|w}bX = \frac{uX|wX}{bX}$. Note that this stresses our intention to
consider $\mathbb T(n)$ 
as the projective space for the vector space $\mathbb R^3$ (with coordinates
in the form $\frac{u|w}b$).

\subsection{Visualization of objects in $\Cal S$.}
\label{sec-two-two}

Let $X = (U,V)$
in $\Cal S$. The isomorphism class  $[V]$ of $V$ can be identified with a partition 
$\lambda = (\lambda_1,\dots,\lambda_b)$  say with $b$ parts. We visualize such
a partition by using boxes. Considering
the part with index $i$, there are $\lambda_i$ boxes, and we arrange these boxes 
vertically (and not, as it is quite common, horizontally): This corresponds to the concept of a
composition series: The top factor should appear at the top, the socle factor at the bottom.
In this way, we obtain for any part a vertical strip of boxes. These strips are arranged from
left to right, but not necessarily in the decreasing sequence stressed by the partition.
Let us repeat: The vertical collection of boxes show the parts of the partition $\lambda$,
(thus they show $[V]$), but the parts are arranged in a suitable order. The
various parts are usually vertically adjusted in order to accommodate generators of a subspace,
as described below.

We may consider the boxes as providing a basis of $V$. Namely, there is a generating
set $\{v_1,\dots,v_b\}$ of $V$ with $T^{\lambda_i}v_i = 0,$ for $1\le i \le b$.
A vertical 
collection of boxes corresponds to such a generator $v_i$, 
and the boxes, going from the top down, to the elements $v_i,\ Tv_i,\ T^2v_i,\ \dots, 
T^{\lambda_i-1}v_i.$
An element $u$ of $V$ is a linear 
combination of the elements $T^jv_i$. Of particular interest will be the elements of $V$
which are sums of elements of the box basis, that is sums of some of the 
elements $T^jv_i$. We will mark the elements $T^jv_i$ by inserting a bullet into the box
and connecting the bullets by a black curve $c(u)$.
(Actually, we will mostly deal with
gradable objects where we are able to obtain horizontal lines;
here the vertical adjustments come into play.  For examples of
non-gradable objects see Section~\ref{sec-fourteen}.)

Now, we look at the subspace $U$. 
If there exists a generating set of $U$ where all the elements $u$ of this set
are sums of elements of the box basis, then we obtain  a clear visualization
using the boxes which describe $V$ and the bullet-lines which describe the
generating set of $U$. 

Here is an example $(U,V)$. We start with $[V] = [6,6,3]$; there is a
corresponding generating set $\{v_1,v_2,v_3\}$ with $T^6v_1 = T^3v_2 = T^6v_3 = 0$;
and we consider the submodule $U$ of $V$ generated by $T^4v_1+T^2v_2$ and $Tv_2+T^2v_3$.
$$
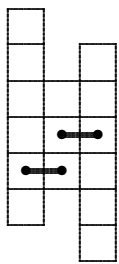
{\beginpicture
    \setcoordinatesystem units <.4cm,.4cm>
\setsolid
\multiput{} at 0 -1  3 6 /
\plot 0 0  1 0  1 6  0 6  0 0 /
\plot 2 -1  3 -1  3 5  2 5  2 -1 /
\plot 0 1  3 1  /
\plot 0 2  3 2  /
\plot 0 3  3 3 /
\plot 0 4  3 4 /
\plot 0 5  1 5 /
\plot 2 0  3 0 /
\multiput{$\bullet$} at 0.5 1.5  1.5 1.5  1.5 2.5  2.5 2.5 /
\plot  0.5 1.5   1.5 1.5 /
\plot  0.5 1.55  1.5 1.55 /
\plot  0.5 1.45  1.5 1.45 /
\plot  1.5 2.5   2.5 2.5 /
\plot  1.5 2.55  2.5 2.55 /
\plot  1.5 2.45  2.5 2.45 /
\endpicture} 
$$

The general situation (where $U$ is not generated by sums of elements of the box basis)
cannot be handled so easily. In some situations, we may be able
to use parameters which we inscribe into the boxes in order 
to provide a description of $U$, see the examples in Sections~\ref{sec-six},
\ref{sec-seven}, \ref{sec-eight}.
The reader should be aware that in general, there does not seem to exist
a satisfying strategy.

\subsection{Covering theory. {\rm \normalsize (A report)}.}
\label{sec-two-three}

  We have introduced objects in the category $\Cal S(n)$ as
  pairs $(U,V)$ where $V$ is a vector space with a nilpotent
  linear operator $T:V\to V$
  of nilpotency index at most $n$ and where $U$ is an invariant subspace,
  hence they can be characterized as representations of the following quiver $Q$
  subject to the commutativity relation $\beta\alpha'=\alpha\beta$,
  the nilpotency relations $\alpha^n=0=(\alpha')^n$ and the condition that
  $\beta$ is a monomorphism.  
$$
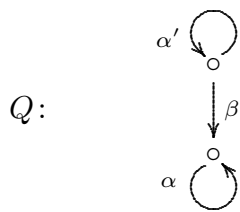
\hbox{\beginpicture
\setcoordinatesystem units <0.5cm,0.5cm>
\put{} at -4 -1
\put{} at 1 3.2
\put{$\circ$} at 0 0
\put{$\circ$} at 0 2
\circulararc 150 degrees from 0 3.2 center at 0 2.7
\circulararc -160 degrees from 0 3.2 center at 0 2.7
\circulararc 150 degrees from 0 -1.2 center at 0 -.7
\circulararc -160 degrees from 0 -1.2 center at 0 -.7
\arr{0 1.6}{0 0.4}
\put{$\ssize \alpha'$} at -1 2.6
\put{$\ssize \alpha$} at -1 -.6
\put{$\ssize \beta$} at 0.4 1
\arr{-0.3 2.31} {-0.2 2.2}  
\arr{0.3 -.31} {0.2 -.2}
\put{$Q\:$} at -4 1
\endpicture}
$$

The category $\widetilde{\Cal S}(n)$ consists of representations of the
universal covering of the quiver $Q$ which is the following locally
finite quiver $\widetilde Q$ subject to commutativity relations
$\beta_{i-1}\alpha'_i=\alpha_i\beta_i$, nilpotency relations
$\alpha_i\cdots\alpha_{i+n-1}=0=\alpha'_i\cdots\alpha'_{i+n-1}$ and the
condition that the vertical maps $\beta_i$ are monomorphisms.
$$
{\beginpicture
    \setcoordinatesystem units <1.2cm,1.2cm>
\multiput{} at -2.5 0  3 1 /
\put{$\widetilde Q\:$} at -2.5 .5
\multiput{$\circ$} at 0 0  1 0  2 0  
     0 1  1 1  2 1   /
\put{$\ssize 1$} at 0 -.25
\put{$\ssize 2$} at 1 -.25 
\put{$\ssize 3$} at 2 -.25
\put{$\ssize \alpha_2$} at .5 -.2
\put{$\ssize \alpha_3$} at 1.5 -.2
\put{$\ssize 1'$} at 0 1.27
\put{$\ssize 2'$} at 1 1.27
\put{$\ssize 3'$} at 2 1.27
\put{$\ssize \alpha'_2$} at .5 1.2
\put{$\ssize \alpha'_3$} at 1.5 1.2
\arr{-.3 0}{-.7 0}
\arr{0.7 0}{0.3 0 }
\arr{1.7 0}{1.3 0 }
\arr{2.7 0}{2.3 0 }
\arr{-.3 1}{-.7 1}
\arr{0.7 1}{0.3 1 }
\arr{1.7 1}{1.3 1 }
\arr{2.7 1}{2.3 1 }
\arr{0 0.7}{0 0.3}
\arr{1 0.7}{1 0.3}
\arr{2 0.7}{2 0.3}
\put{$\ssize \beta_1$} at .2 .5
\put{$\ssize \beta_2$} at 1.2 .5
\put{$\ssize \beta_3$} at 2.2 .5
\setdots <1mm>
\plot -1.5 0  -1 0 /
\plot -1.5 1  -1 1 /
\plot 3 0  3.5 0 /
\plot 3 1  3.5 1 /
\endpicture}
$$
The categories $\widetilde{\Cal S}(n)$ and $\Cal S(n)$ are related
via the {\it covering functor} or {\it push-down functor} $\pi_\lambda$, it maps a representation
$M = ((M_x)_{x\in\widetilde Q_0},(M_\gamma)_{\gamma\in\widetilde Q_1})$ to
the pair $\pi_\lambda M = (U,V)$ with $U=\bigoplus_{i\in\mathbb Z}M_{i'}$,
$V=\bigoplus_{i\in\mathbb Z}M_i$ where the embedding is given by
the diagonal map $\bigoplus_{i\in\mathbb Z}\beta_i$ and the action of
$T$ on $V$ is given by the maps $\alpha_i$.
The group $G=\mathbb Z$ acts on the quiver $\widetilde Q$ and on
representations of $\widetilde Q$ as index shift (so $Q=\widetilde Q/G$).
Note that the quiver $\widetilde Q$ is locally finite and the group $G = \mathbb Z$ is torsionfree and
acts freely on $\widetilde Q_0.$ 
	
The union of the categories $\widetilde{\Cal S}(n)$ with $n\in \mathbb N$ will be denoted by
$\widetilde{\Cal S}.$ If $M$ belongs to $\widetilde{\Cal S}(n)$, we call $\frac{u|w}bM =
\frac{u|w}b\pi_\lambda M$ the {\it uwb-vector,} and $\pr M = \pr \pi_\lambda M$ the
{\it pr-vector of $M$}.

\medskip
\begin{theorem}[Covering Theorem (Gabriel, Dowbor--Skowro\'nski)]
  Let $\widetilde Q$ be a locally finite quiver and $G$
  a torsionfree group of automorphisms of $\widetilde Q$ which acts freely on $\widetilde Q_0$.
  Let $Q = \widetilde Q/G$.
  The covering functor $\pi_\lambda$ provides an injective map from the set of \ $G$-orbits
  of isomorphism classes of indecomposable $k\widetilde Q$-modules to the set of
  isomorphism classes of indecomposable $kQ$-modules. $\s$
\end{theorem}

\medskip
In general, the injective map induced by $\pi_\lambda$ on the set of indecomposable modules
is not surjective.
But we do get a bijection if the category of $k\widetilde Q$-modules is locally
support finite, see \cite{DLS}). Regarding invariant subspaces for $n\leq 6$, we deal with the
subcategory $\widetilde{\Cal S}(n)$ of 
$\mod k\widetilde Q$ which is locally support finite, and it similarly follows that 
the restriction of $\pi_\lambda$ to $\widetilde{\Cal S}(n)$ induces
a bijection between the $G$-orbits of isomorphism classes of indecomposable
objects in $\widetilde{\Cal S}(n)$ and the isomorphism classes of indecomposables in
$\Cal S(n)$ \cite[(2.1) and (2.2)]{RS1}.

An object in $\Cal S(n)$ is said to be {\it gradable} provided it is in the image
of the push-down functor $\pi_\lambda$. All pickets are, of course, gradable, as are all
bipickets, see Section~\ref{sec-twelve-two}. As we have mentioned, all objects in $\Cal S(n)$ for $n\le 6$
are gradable. This is no longer true for $n = 7,$ for a typical example of a 
non-gradable object in $\Cal S(7)$ as well as further remarks on non-gradable
objects in $\Cal S(n)$, we refer to Section~\ref{sec-fourteen}. Note that
all the objects of $\Cal S(n)$ considered in the remaining parts of the paper are gradable.

\medskip
\begin{remark}
  Covering theory has been developed by Gabriel, Bongartz and Riedtmann,
  and in parallel by Gordan and Green, dealing with group-graded algebras and the
  corresponding graded modules. The formulation of the covering theorem 
  is quoted from \cite{R3}, reprinted here, slightly revised, as Appendix~\ref{app-D},
  where a short proof is given. Note that there is no assumption on the base field. 
\end{remark}

\medskip
For the purpose of exhibiting special objects, we usually restrict 
the domain of $\pi_\lambda$:
For $\ell \le m \in \mathbb Z$, we denote by $Q[\ell,m]$
the following fully commutative quiver 
with $2(m-\ell+1)$ vertices,
$$
{\beginpicture
    \setcoordinatesystem units <1cm,1cm>
\multiput{} at 0 0  5 1 /
\put{$Q[\ell,m]:$} at -2 .5
\multiput{$\circ$} at 0 0  1 0  2 0  4 0  5 0 
     0 1  1 1  2 1    5 1  4 1   /
\put{$\ssize \ell$} at 0 -.25
\put{$\ssize \ell+1$} at 1 -.25 
\put{$\ssize \ell+2$} at 2.05 -.25
\put{$\ssize m-1$} at 4 -.25
\put{$\ssize m$} at 5 -.25
\put{$\ssize \ell'$} at 0 1.27
\put{$\ssize (\ell+1)'$} at .95 1.27
\put{$\ssize (\ell+2)'$} at 2.1 1.27
\put{$\ssize (m-1)'$} at 4 1.27
\put{$\ssize m'$} at 5 1.27
\arr{0.7 0}{0.3 0 }
\arr{1.7 0}{1.3 0 }
\arr{2.6 0}{2.3 0 }
\arr{4.7 0}{4.3 0 }
\plot 3.75 0  3.5 0 /

\arr{0.7 1}{0.3 1 }
\arr{1.7 1}{1.3 1 }
\arr{2.6 1}{2.3 1 }
\arr{4.7 1}{4.3 1 }
\plot 3.75 1  3.5 1 /
\arr{0 0.7}{0 0.3}
\arr{1 0.7}{1 0.3}
\arr{2 0.7}{2 0.3}
\arr{5 0.7}{5 0.3}
\arr{4 0.7}{4 0.3}
\setdots <1mm>
\plot 3 0  3.5 0 /
\plot 3 1  3.5 1 /
\endpicture}
$$
and denote by $A=A(\ell,m)$ the factor of its path algebra 
modulo all the commutativity relations.

Similarly, we write $\Cal S[\ell,m]$ for the category of all representations of 
$\widetilde Q$ with all maps $\beta_i$ (with $\ell \le i \le m$) being monomorphisms. 

\medskip
\begin{enumerate}
\item[$\bullet$] For $\ell<m$, the algebra $A$ has global dimension 2.
\item[$\bullet$] For tuples $x=(x_i)_{i\in Q_0}$, $y=(y_i)_{i\in Q_0}$
the homological bilinear form for $A$ is given by
$$\langle x,y\rangle=\sum_{i\in Q_0} x_iy_i-\sum_{(\alpha:i\to j)\in Q_1} x_iy_j
+\sum_{i=\ell+1}^m x_{i'}y_{i-1}.$$
\item[$\bullet$] The objects in $\Cal S[\ell,m]$ have projective
dimension at most one in $\mod A$.  Namely, the simple $A$-modules $S(i)$
as well as the non-split extensions of $S(i')$ by $S(i)$ (where
$\ell\leq i\leq m$), have this property, and any object in 
$\Cal S[\ell,m]$ has a filtration with such factors. 
\item[$\bullet$] For objects $X,Y\in \Cal S[\ell,m]$,
the bilinear form satisfies \cite[2.4 Lemma]{R2}
$$
  \langle \bdim X,\bdim Y\rangle=
  \dim \Hom_A(X,Y)-\dim \Ext^1_A(X,Y).
$$
\end{enumerate}

\bigskip
Let $X$ be an indecomposable gradable object. The {\it standard grading} is the grading 
such that there do exist non-zero elements of degree 1 and such that all non-zero elements have positive degree.

\subsection{The $d$-Kronecker algebras $K(d)$. {\rm \normalsize (A report)}.}
\label{sec-two-four}

We denote by $K(d)$ the $d$-Kronecker algebra with coefficients in $k$ and $d$ arrows
(it is the path algebra of the quiver with 2 vertices, say labeled 1 and 2, and $d$
arrows $1\leftarrow 2$). The {\it Kronecker algebra} is the $2$-Kronecker algebra. 

There are two simple $K(d)$-modules, $S(1)$ and $S(2)$,
corresponding to the vertices of $K(d)$. The module 
$S(1)$ is projective, the module $S(2)$ is injective. For $d = 0$, the algebra $K(0)$ is
semisimple, the modules $S(1)$ and $S(2)$ are the only indecomposable modules. For $d = 1$,
there is a unique indecomposable module which is not simple, say $I$, with a non-split
exact sequence $0 \to S(1) \to I \to S(2) \to 0.$

For $d \ge 2,$ the algebra $K(d)$ is representation infinite, with a preprojective component
(containing $S(1)$), a preinjective component (containing $S(2)$), as well as 
with infinitely many stable
Auslander-Reiten components. The $K(d)$-modules without a non-zero direct summand which is preprojective or preinjective are said to be {\it regular.}

The case $d=2$ is of special interest. In this case, the full subcategory $\Cal R$ of regular modules
is an exact abelian subcategory. The simple objects in $\Cal R$ are said to be {\it quasi-simple}
in $\mod K(2)$. Any indecomposable object in $\Cal R$ has a unique composition series 
in $\Cal R$, 
say of length $\ell$, called its {\it quasi-length}
and its composition factors in $\Cal R$ are isomorphic, say isomorphic to $R$: In this
case we will denote the object by $R[\ell].$ In this way we obtain a bijection between the
(isomorphism classes of the) indecomposable objects in $\Cal R$ and the set of pairs $(R,\ell)$, where
$R$ is (the isomorphism class of) a simple object in $\Cal R$ and $\ell\in \mathbb N_1.$
For any simple object $R$ in $\Cal R$, the objects of the form $R[\ell]$ with $\ell\in \mathbb N_1$ form
an Auslander-Reiten component in $\mod K(2)$ and this component is a stable tube of rank 1.

\medskip
Recall that we denote by $\mathbb P^1$ the projective line over $k$; its elements are written
in the form $c = (c_0:c_1) = k(c_0,c_1)$, where $(c_0,c_1)$ is a non-zero element of $k^2$.
Any element $c$ of $\mathbb P^1$ gives rise to an indecomposable representation
$R_c = (k\leftleftarrows^{c_0}_{c_1} k)$ of length 2 in $\mod K(2)$
(and any indecomposable representation of $K(2)$ of length 2 is obtained in this way). 
The representation $R_c$ is a simple object in $\Cal R$.
The further simple objects of $\Cal R$ are related to finite field extensions of $k$;
in particular, {\it all simple objects of $\Cal R$ have length 2 in $\mod K(2)$ if and only if $k$ is
algebraically closed.} In order to avoid field extensions, we usually will restrict to look
at the representations $R_c$, instead of taking all simple objects in $\Cal R$ into account. 
In this way, we focus the attention to the $\mathbb P^1$-family $R_c$, and, more generally, 
to the corresponding $\mathbb P^1$-families $R_c[\ell],$ with fixed $\ell\in \mathbb N_1.$

\subsection{Simplification. {\rm\normalsize (A report)}.}
\label{sec-two-five}

For the construction of modules in $\widetilde{\Cal S}$, we often will use simplification.
The following result is a special case of \cite[(1.5) Lemma]{R1}.
Given a finite dimensional $k$-algebra $A$, we say that the $A$-modules $X,Y$
are an {\it orthogonal pair} provided
$\End X=\End Y=k$, $\Hom(X,Y)=\Hom(Y,X)=0$.
Let $d$ be a natural number. A {\it $d$-Kronecker pair} $X,Y$ is an (ordered) orthogonal pair
such that $\dim_k\Ext^1(X,Y) = d.$

\medskip
\begin{theorem}[Simplification Lemma]
  Let $A$ be a $k$-algebra. 
  Let $X,Y$ be a $d$-Kronecker pair in $\mod A$.
  The full subcategory 
  of all objects $Z$ with a subobject $Z'$ isomorphic to a direct sum of copies of $Y$
  such that $Z/Z'$ is isomorphic to a direct sum of copies of $X$, is equivalent to the
  category of $K(d)$-modules.
\end{theorem}

\medskip
We may reformulate the Simplification Lemma as follows:
Let $A$ be a $k$-algebra and assume that 
$X,Y$ is a $d$-Kronecker pair. Then there is a full exact embedding
$$F_{X,Y}\:\mod K(d)\to \mod A$$ such that
$F_{X,Y}(S(2))=X$ and $F_{X,Y}(S(1))=Y$.

The image of the functor $F_{X,Y}$ is the full subcategory of $\mod A$ 
consisting of all modules which have a submodule which is the direct sum of copies of $Y$
such that the corresponding factor module is the direct sum of copies of $X$.

\subsection{Kronecker subcategories. Kronecker families. BTh-vectors.}
\label{sec-two-six}

An exact subcategory $\Cal K$ of $\mod \widetilde Q$ 
will be called a {\it Kronecker subcategory} provided
$\Cal K$ is equivalent (as an exact subcategory) to $\mod K(2).$ Of course, this means that there
are given two orthogonal objects $X,Y$ in $\Cal K$ (namely the simple objects of $\Cal K$) 
with $\End(X) = 
\End(Y) = k$, and an indecomposable object $Z$ with a short exact sequence 
$0 \to Y^2\to Z \to X \to  0$ (the projective cover of $X$ inside $\Cal K$). 

If $\Cal K$ is a Kronecker subcategory of $\widetilde {\Cal S},$ the set of indecomposable
objects in $\Cal K$ which have length 2 in $\Cal K$, will be called a {\it Kronecker family} in 
$\widetilde {\Cal S}.$ Note that a Kronecker family is a $\mathbb P^1$-family of indecomposable 
and pairwise non-isomorphic objects $M_c$ in $\widetilde {\Cal S},$ with 
$c = (c_0:c_1)\in \mathbb P^1.$ 
For any $\ell\in \mathbb N_1,$ there exists the object $M_c[\ell]$ in $\Cal K$ (it has
a filtration with $\ell$ factors $M_c$), and the set of objects $M_c[\ell]$ 
with $c = (c_0:c_1)\in \mathbb P^1$ is a $\mathbb P^1$-family of indecomposable 
and pairwise non-isomorphic objects in $\widetilde {\Cal S}$. 

\medskip
\begin{warning}
  The notation may be 
  misleading. The construction of the object $M_c[\ell]$ uses not only $M_c$, but actually 
  the choice of $X$ and $Y$, thus the subcategory $\Cal K$.
  As we will see, an object of a Kronecker family in $\widetilde{\Cal S}$ may belong to several
  different Kronecker subcategories (see Remark~\ref{rem-two-seven-two} below).
  But we hope the 
  reader is not irritated in this way. 
\end{warning}

\subsection{Example: The standard family in $\Cal S(6)$.}
\label{sec-two-seven}

We present here in detail a construction for a BTh-vector
which will be used repeatedly in this manuscript. This BTh-family in 
$\Cal S(6)$ was presented in \cite[(2.3)]{RS1}, but actually has
already been observed by Birkhoff \cite{Bh}.

\medskip
Given the picket $X=([4],[6])$ and the bipicket $Y=([2],[4,2],[3,1])$,
consider the indecomposable representations for $Q[1,6]$
of dimension type
$\bdim \widetilde X=\smallmatrix 111100\\ 111111\endsmallmatrix$
and $\bdim \widetilde Y=
\smallmatrix 011000\\ 012210\endsmallmatrix$.
Clearly, $\widetilde X$ and $\widetilde Y$ is an orthogonal pair in $\mod A(1,6)$
and $X=\pi_\lambda \widetilde X$, $Y=\pi_\lambda \widetilde Y$. 

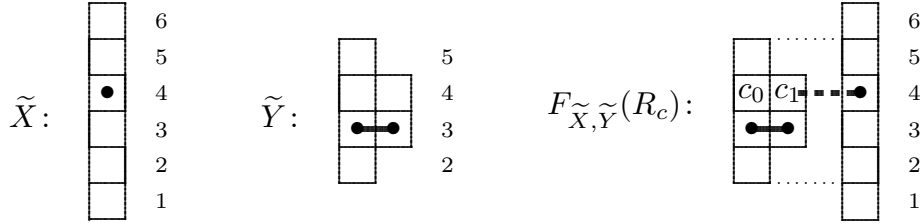
\addcontentsline{lof}{subsection}{The standard family in $\mathcal S(6)$.}

$$
{\beginpicture
    \setcoordinatesystem units <.4cm,.4cm>
\multiput{} at -1.5 0   3  6 /
\put{$\widetilde X\:$} at -1.5 3
\plot 0 0  1 0  1 6   0 6  0 0 /
\plot 0 1  1 1  1 2  0 2  0 3  1 3  1 4  0 4  0 5  1 5 /
\put{$\scriptstyle 1$} at 2 0.5
\put{$\scriptstyle 2$} at 2 1.5
\put{$\scriptstyle 3$} at 2 2.5
\put{$\scriptstyle 4$} at 2 3.5
\put{$\scriptstyle 5$} at 2 4.5
\put{$\scriptstyle 6$} at 2 5.5
\multiput{$\bullet$} at 0.5 3.5 /
\endpicture
}
\qquad
{\beginpicture
    \setcoordinatesystem units <.4cm,.4cm>
\multiput{} at -1.5 0   4  6 /
\put{$\widetilde Y\:$} at -1.5 3
\plot 0 1  1 1  1 5  0 5  0 1 /
\plot 0 2  2 2  2 4  0 4  0 3  2 3 /
\put{$\scriptstyle 2$} at 3 1.5
\put{$\scriptstyle 3$} at 3 2.5
\put{$\scriptstyle 4$} at 3 3.5
\put{$\scriptstyle 5$} at 3 4.5
\multiput{$\bullet$} at 0.5 2.5  1.5 2.5 /
\plot 0.5 2.5  1.5 2.5 /
\plot 0.5 2.45  1.5 2.45 /
\plot 0.5 2.55  1.5 2.55 /
\endpicture
}\qquad
{\beginpicture
    \setcoordinatesystem units <.4cm,.4cm>
\multiput{} at -7 0   1  6 /
\put{$F_{\widetilde X,\widetilde Y}(R_c)\:$} at -6 3
\plot 0 0  1 0  1 6   0 6  0 0 /
\plot 0 1  1 1  1 2  0 2  0 3  1 3  1 4  0 4  0 5  1 5 /
\put{$\scriptstyle 1$} at 2 0.5
\put{$\scriptstyle 2$} at 2 1.5
\put{$\scriptstyle 3$} at 2 2.5
\put{$\scriptstyle 4$} at 2 3.5
\put{$\scriptstyle 5$} at 2 4.5
\put{$\scriptstyle 6$} at 2 5.5
\multiput{$\bullet$} at 0.5 3.5 /
\plot -3 1  -2 1  -2 5  -3 5  -3 1 /
\plot -3 2  -1 2  -1 4  -3 4  -3 3  -1 3 /
\multiput{$\bullet$} at -1.5 2.5  -2.5 2.5 /
\plot -1.5 2.5  -2.5 2.5 /
\plot -1.5 2.55 -2.5 2.55 /
\plot -1.5 2.45 -2.5 2.45 /
\put{$c_0$} at -2.5 3.5 
\put{$c_1$} at -1.5 3.5 
\setdashes <.9mm>
\plot -1.2 3.5   0.5 3.5 /
\plot -1.2 3.55  0.5 3.55 /
\plot -1.2 3.45  0.5 3.45 /

\setdots <1mm>
\plot -2 5  0 5 /
\plot -2 1  0 1 /
\endpicture}
$$
Using the formulas from Section~\ref{sec-two-three}, we verify that 
$\langle\bdim \widetilde X,\bdim \widetilde Y\rangle=-2$,
hence 
$\dim\Ext^1(\widetilde X,\widetilde Y)=2$ and the embedding
$F_{\widetilde X,\widetilde Y}\:\mod K(2)\to \mod A(1,6)$
yields a BTh-family in $\widetilde{\Cal S}(6)$ with uwb-vector $\frac{6\,|\,6}3$ and, via the covering functor,
in $\Cal S(6)$.  For $c=(c_0:c_1)\in\mathbb P^1(k)$, the objects
$F_{\widetilde X,\widetilde Y}(R_c)$ are sketched above.
Note that $\uwbb X = \frac{4|2}1$ and  $\uwbb Y = \frac{2|4}2$:

\medskip
We describe the functor $F_{\widetilde X,\widetilde Y}$ in detail.
Let $J=\big(J_1\leftleftarrows^\zeta_\eta J_2\big)$ be a representation for the
Kronecker algebra $K(2)$.  The corresponding representation of the quiver $Q[1,6]$
is as follows.

$$
\beginpicture 
    \setcoordinatesystem units <2cm,1.5cm>
\multiput{} at 0 0  5 1 /
\multiput{$0$} at 4 1  5 1 /
\multiput{$J_2$} at 0 0  0 1  3 1  5 0 /
\multiput{$J_1\oplus J_2$} at 1 0  4 0  1 1  2 1 /
\multiput{$J_1^2\oplus J_2$} at 2 0  3 0 /
\arr{.65 0}{.3 0}
\arr{1.7 0}{1.35 0}
\arr{3.7 0}{3.35 0}
\arr{4.75 0}{4.35 0}
\arr{.65 1}{.3 1}
\plot 1.35 1  1.65 1 /
\plot 1.35 1.04 1.65 1.04 /

\arr{2.7 1}{2.35 1}
\arr{3.7 1}{3.3 1}
\arr{4.7 1}{4.3 1}
\arr{2 .7}{2 .3}
\arr{3 .7}{3 .3}
\arr{4 .7}{4 .3}
\arr{5 .7}{5 .3}

\plot 0 .7  0 .3 /
\plot .03 .7  .03 .3 /

\plot 1 .7  1 .3 /
\plot 1.03 .7  1.03 .3 /

\plot 2.35 0  2.65 0 /
\plot 2.35 .04 2.65 .04 /

\multiput{$\iota$} at 4.55 .15  3.55 .15  2.55 1.15 /
\multiput{$\pi$} at  .5 .15  1.5 .15  .5 1.15 /
\put{$\left[\smallmatrix \zeta \cr
                         \eta \cr
                         1\endsmallmatrix
      \right]$} at 3.25 .5
\put{$\left[\smallmatrix 
                         1&\zeta \cr
                         1&\eta \cr 
                         0&1\endsmallmatrix
      \right]
     $} at 2.3 .5

\endpicture
$$
By
$\iota$ and $\pi$ we denote the canonical inclusion into the last component(s) and the
canonical projection modulo the first component.

\medskip        
\begin{lemma}
  \label{lem-two-seven}
  \begin{itemize}[leftmargin=3em]
  \item[\rm(a)] The functor $F_{\widetilde X,\widetilde Y}$ is a full and exact embedding.
  \item[\rm(b)] For $J\in\mod K(2)$
    of dimension vector $(j_1,j_2)$, 
    the corresponding object $F_{\widetilde X,\widetilde Y}J$
    in $\widetilde{\Cal S}(6)$ has uwb-vector $\frac{2j_1+4j_2\;|\;4j_1+2j_2}{2j_1+j_2}$.
  \item[\rm(c)] The images of the modules in homogeneous tubes in $\mod K(2)$ under 
    $F_{\widetilde X,\widetilde Y}$ provide a BTh-family in $\widetilde{\Cal S}(6)$ with pr-vector $(2,2)$.
    $\s$
  \end{itemize}
\end{lemma}

	\medskip
Under $F_{\widetilde X,\widetilde Y}$, the maps in $\mod K(2)$ give rise to maps between
the corresponding objects in $\widetilde{\Cal S}(6)$. 

\medskip
\begin{remark}[A correction to \cite{RS1}]
  \label{rem-two-seven-one}
  Let $\Cal K$ be the image of the functor 
  $F_{\widetilde X,\widetilde Y}$. This is a Kronecker category of $\widetilde{\Cal S}(6)$,
  thus equivalent to $\mod K(2).$ In particular, $\Cal K$ has relative Auslander-Reiten
  sequences. Let $c = (c_0:c_1) \in \mathbb P^1(k)$. The Auslander-Reiten component 
  of $F_{\widetilde X,\widetilde Y}(R_c)$ in $\Cal K$ is a stable tube of rank 1.  

We are going to present explicitly the first three objects in this tube, with 
the (relative) irreducible maps between them (of course, these maps are not
necessarily irreducible morphisms in $\widetilde{\Cal S}(6)$).

First, let us consider the case that $c_0\neq 0$, thus, we can assume that $c_0 = 1$
(and $c_1$ is an arbitrary element of $k$).
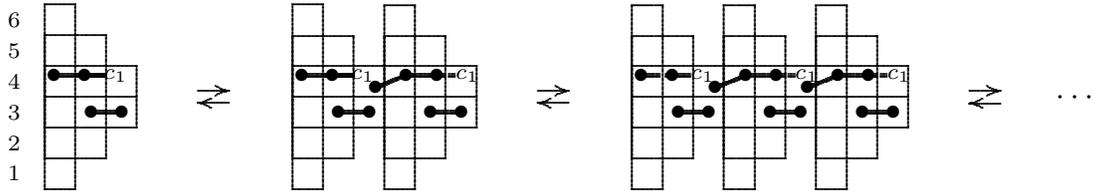
$$
\beginpicture
    \setcoordinatesystem units <.34cm,.34cm>
\multiput{} at -1 0   6 6 /
\plot 0 3  0 6  1 6  1 0  0 0  0 5  2 5  2 1  0 1  0 4  3 4  3 2  0 2  0 3  3 3 /
\put{$\scriptstyle 1$} at -1 0.5
\put{$\scriptstyle 2$} at -1 1.5
\put{$\scriptstyle 3$} at -1 2.5
\put{$\scriptstyle 4$} at -1 3.5
\put{$\scriptstyle 5$} at -1 4.5
\put{$\scriptstyle 6$} at -1 5.5
\multiput{$\bullet$} at 0.3 3.7  1.5 2.5  2.5 2.5 /
\plot 1.5 2.5  2.5 2.5 /
\plot 1.5 2.55  2.5 2.55 /
\plot 1.5 2.45  2.5 2.45 /
\put{$\bullet$} at 1.3 3.7 
\put{$\ss c_1$} at 2.3 3.7

\plot .3 3.7  2 3.7 /  \plot 1.1 3.7  2 3.7 /
\plot .3 3.75  2 3.75 /  \plot 1.1 3.75  2 3.75 /
\plot .3 3.65  2 3.65 /  \plot 1.1 3.65  2 3.65 /

\arr{5 3.2}{6 3.2}
\arr{6 2.8}{5 2.8}
\endpicture
\qquad
\beginpicture
    \setcoordinatesystem units <.34cm,.34cm>
\multiput{} at 0 0   9 6 /
\plot 0 3  0 6  1 6  1 0  0 0  0 5  2 5  2 1  0 1  0 4  6 4  6 2  0 2  0 3  6 3 /
\plot 3 3  3 6  4 6  4 0  3 0  3 5  5 5  5 1  3 1  3 3 /
\multiput{$\bullet$} at 0.3 3.7  2.7 3.3  3.7 3.7  1.5 2.5  2.5 2.5  4.5 2.5  5.5 2.5 /
\plot 1.5 2.5  2.5 2.5 /
\plot 1.5 2.55  2.5 2.55 /
\plot 1.5 2.45  2.5 2.45 /
\plot 4.5 2.5  5.5 2.5 /
\plot 4.5 2.55  5.5 2.55 /
\plot 4.5 2.45  5.5 2.45 /
\multiput{$\bullet$} at 1.3 3.7  4.7 3.7  / 
\multiput{$\ss c_1$} at 2.3 3.7  5.7 3.7  /

\plot .3 3.7  2 3.7 /  \plot 1.1 3.7  2 3.7 /
\plot .3 3.75  2 3.75 /  \plot 1.1 3.75  2 3.75 /
\plot .3 3.65  2 3.65 /  \plot 1.1 3.65  2 3.65 /

\plot 5.1 3.7    5.3 3.7 /  \plot  2.7 3.3  3.7 3.7  5  3.7 /
\plot 5.1 3.75   5.3 3.75 /  \plot 2.7 3.35  3.7 3.75  5  3.75 /
\plot 5.1 3.65   5.3 3.65 /  \plot 2.7 3.25  3.7 3.65  5  3.65 /

\arr{8 3.2}{9 3.2}
\arr{9 2.8}{8 2.8}
\endpicture
\qquad
\beginpicture
    \setcoordinatesystem units <.34cm,.34cm>
\multiput{} at 0 0   9 6 /
\plot 0 3  0 6  1 6  1 0  0 0  0 5  2 5  2 1  0 1  0 4  9 4  9 2  0 2  0 3  9 3 /
\plot 3 3  3 6  4 6  4 0  3 0  3 5  5 5  5 1  3 1  3 3 /
\plot 6 3  6 6  7 6  7 0  6 0  6 5  8 5  8 1  6 1  6 3 /
\multiput{$\bullet$} at 0.3 3.7  2.7 3.3  3.7 3.7  1.5 2.5  2.5 2.5  4.5 2.5  5.5 2.5 
                     6.7 3.7  5.7 3.3  7.5 2.5  8.5 2.5 /
\plot 1.5 2.5  2.5 2.5 /
\plot 1.5 2.55  2.5 2.55 /
\plot 1.5 2.45  2.5 2.45 /
\plot 4.5 2.5  5.5 2.5 /
\plot 4.5 2.55  5.5 2.55 /
\plot 4.5 2.45  5.5 2.45 /
\plot 7.5 2.5   8.5 2.5 /
\plot 7.5 2.55  8.5 2.55 /
\plot 7.5 2.45  8.5 2.45 /
\multiput{$\bullet$} at 1.3 3.7  4.7 3.7  7.7 3.7 / 
\multiput{$\ss c_1$} at 2.3 3.7  5.7 3.7  8.7 3.7 / 

\plot .3 3.7  .9 3.7 /  \plot 1.1 3.7  1.9 3.7 /
\plot .3 3.75  .9 3.75 /  \plot 1.1 3.75  1.9 3.75 /
\plot .3 3.65  .9 3.65 /  \plot 1.1 3.65  1.9 3.65 /

\plot 5.1 3.7    5.3 3.7 /  \plot  2.7 3.3  3.7 3.7  5  3.7 /
\plot 5.1 3.75   5.3 3.75 /  \plot 2.7 3.35  3.7 3.75  5  3.75 /
\plot 5.1 3.65   5.3 3.65 /  \plot 2.7 3.25  3.7 3.65  5  3.65 /

\plot 8.1 3.7    8.3 3.7 /  \plot  5.7 3.3  6.7 3.7  8  3.7 /
\plot 8.1 3.75   8.3 3.75 /  \plot 5.7 3.35  6.7 3.75  8  3.75 /
\plot 8.1 3.65   8.3 3.65 /  \plot 5.7 3.25  6.7 3.65  8  3.65 /

\arr{11 3.2}{12 3.2}
\arr{12 2.8}{11 2.8}
\put{$\cdots$} at 14.5 3
\endpicture
$$
If $c_0 = 0$, then the corresponding first three objects look as follows:
$$
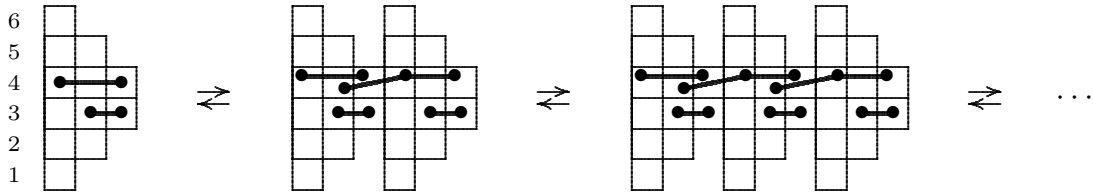
\beginpicture
    \setcoordinatesystem units <.34cm,.34cm>
\multiput{} at -1 0   6 6 /
\plot 0 3  0 6  1 6  1 0  0 0  0 5  2 5  2 1  0 1  0 4  3 4  3 2  0 2  0 3  3 3 /
\put{$\scriptstyle 1$} at -1 0.5
\put{$\scriptstyle 2$} at -1 1.5
\put{$\scriptstyle 3$} at -1 2.5
\put{$\scriptstyle 4$} at -1 3.5
\put{$\scriptstyle 5$} at -1 4.5
\put{$\scriptstyle 6$} at -1 5.5
\multiput{$\bullet$} at 0.5 3.5  1.5 2.5  2.5 2.5 /
\plot 1.5 2.5  2.5 2.5 /
\plot 1.5 2.55  2.5 2.55 /
\plot 1.5 2.45  2.5 2.45 /
\put{$\bullet$} at 2.5 3.5 
\plot .5 3.5   2.5 3.5 /
\plot .5 3.55  2.5 3.55 /
\plot .5 3.45  2.5 3.45 /
\arr{5 3.2}{6 3.2}
\arr{6 2.8}{5 2.8}
\endpicture
\qquad
\beginpicture
    \setcoordinatesystem units <.34cm,.34cm>
\multiput{} at 0 0   9 6 /
\plot 0 3  0 6  1 6  1 0  0 0  0 5  2 5  2 1  0 1  0 4  6 4  6 2  0 2  0 3  6 3 /
\plot 3 3  3 6  4 6  4 0  3 0  3 5  5 5  5 1  3 1  3 3 /
\multiput{$\bullet$} at 0.3 3.7  1.7 3.3  3.7 3.7  1.5 2.5  2.5 2.5  4.5 2.5  5.5 2.5 /
\plot 1.5 2.5  2.5 2.5 /
\plot 1.5 2.55  2.5 2.55 /
\plot 1.5 2.45  2.5 2.45 /
\plot 4.5 2.5  5.5 2.5 /
\plot 4.5 2.55  5.5 2.55 /
\plot 4.5 2.45  5.5 2.45 /
\multiput{$\bullet$} at 2.3 3.7  5.3 3.7  /

\plot .3 3.7   2.3 3.7 /
\plot .3 3.75  2.3 3.75 /
\plot .3 3.65  2.3 3.65 /
\plot  1.7 3.3   3.7  3.7  5.3 3.7 /
\plot 1.7 3.35  3.7  3.75  5.3 3.75 /
\plot 1.7 3.25  3.7  3.65  5.3 3.65 /
\arr{8 3.2}{9 3.2}
\arr{9 2.8}{8 2.8}
\endpicture
\qquad
\beginpicture
    \setcoordinatesystem units <.34cm,.34cm>
\multiput{} at 0 0   9 6 /
\plot 0 3  0 6  1 6  1 0  0 0  0 5  2 5  2 1  0 1  0 4  9 4  9 2  0 2  0 3  9 3 /
\plot 3 3  3 6  4 6  4 0  3 0  3 5  5 5  5 1  3 1  3 3 /
\plot 6 3  6 6  7 6  7 0  6 0  6 5  8 5  8 1  6 1  6 3 /
\multiput{$\bullet$} at 0.3 3.7  1.7 3.3  3.7 3.7  1.5 2.5  2.5 2.5  4.5 2.5  5.5 2.5 
                     6.7 3.7  5.3 3.7  7.5 2.5  8.5 2.5 /
\plot 1.5 2.5  2.5 2.5 /
\plot 1.5 2.55  2.5 2.55 /
\plot 1.5 2.45  2.5 2.45 /
\plot 4.5 2.5  5.5 2.5 /
\plot 4.5 2.55  5.5 2.55 /
\plot 4.5 2.45  5.5 2.45 /
\plot 7.5 2.5   8.5 2.5 /
\plot 7.5 2.55  8.5 2.55 /
\plot 7.5 2.45  8.5 2.45 /
\multiput{$\bullet$} at 2.3 3.7  4.7 3.3  8.3 3.7 / 

\plot .3 3.7   2.3 3.7 /
\plot .3 3.75  2.3 3.75 /
\plot .3 3.65  2.3 3.65 /
\plot  1.7 3.3   3.7  3.7  5.3 3.7 /
\plot 1.7 3.35  3.7  3.75  5.3 3.75 /
\plot 1.7 3.25  3.7  3.65  5.3 3.65 /
\plot  4.7 3.3   6.7  3.7  8.3 3.7 /
\plot 4.7 3.35  6.7  3.75  8.3 3.75 /
\plot 4.7 3.25  6.7  3.65  8.3 3.65 /

\arr{11 3.2}{12 3.2}
\arr{12 2.8}{11 2.8}
\put{$\cdots$} at 14.5 3
\endpicture
$$
(Note that this corrects the presentation in \cite{RS1}
given at the end of Section~2.3.)

We should mention that we deal with indecomposable objects which belong to the tubular family
of $\widetilde{\Cal S}(6)$ which contains the principal component (the family with
rationality index $0$). 
If $c_1=0$, the objects occur in the principal component;
if $c_0 = 0$, they occur in the 2-tube; for  $c_0=c_1$, they lie in the 3-tube.
\end{remark}

\medskip
\begin{remark}[A second $\mathbb P^1$-family in $\Cal S(6)$ with width $3$]
  \label{rem-two-seven-two}
  Let us now start with the Kronecker pair $\widetilde X'$, $\widetilde Y'$
with $\pi_\lambda \widetilde X' = ([2],[2],[0])$ and
$\pi_\lambda\widetilde Y' = ([3,1],[6,4],[4,2])$, and look at the corresponding Kronecker family 
$$
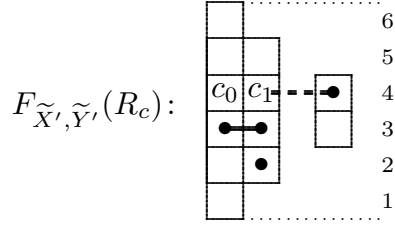
{\beginpicture
    \setcoordinatesystem units <.4cm,.4cm>
\multiput{} at -7 0   1  6 /
\put{$F_{\widetilde X',\widetilde Y'}(R_c)\:$} at -6 3
\plot 0 2  1 2  1 4   0 4  0 2 /
\plot 0 3  1 3  /
\put{$\scriptstyle 1$} at 2 0.5
\put{$\scriptstyle 2$} at 2 1.5
\put{$\scriptstyle 3$} at 2 2.5
\put{$\scriptstyle 4$} at 2 3.5
\put{$\scriptstyle 5$} at 2 4.5
\put{$\scriptstyle 6$} at 2 5.5
\multiput{$\bullet$} at 0.5 3.5 /
\plot -3 0  -2 0  -2 6  -3 6  -3 0 /
\plot -3 1  -1 1  -1 5  -3 5  -3 1  /
\plot -3 2  -1 2 /
\plot -3 3  -1 3 /
\plot -3 4  -1 4 /

\multiput{$\bullet$} at -1.5 2.5  -2.5 2.5  -1.5 1.5 /
\plot -1.5 2.5  -2.5 2.5 /
\plot -1.5 2.55 -2.5 2.55 /
\plot -1.5 2.45 -2.5 2.45 /
\put{$c_0$} at -2.5 3.5 
\put{$c_1$} at -1.5 3.5 
\setdashes <.9mm>
\plot -1.2 3.5   0.5 3.5 /
\plot -1.2 3.55  0.5 3.55 /
\plot -1.2 3.45  0.5 3.45 /

\setdots <1mm>
\plot -2 6  2 6 /
\plot -2 0  2 0 /
\endpicture}
$$
The modules which we obtain and which belong to homogeneous tubes are 
those of the standard family
$F_{\widetilde X,\widetilde Y}(R_{c})$.

However, let us look at the 3-tube $\widetilde{\Cal C}$ with rationality index 0. 
The module $M = F_{\widetilde X,\widetilde Y}(R_{(1,1)})$
does not occur in the image of the functor $F_{\widetilde X',\widetilde Y'}$,
and, conversely, the module $M' = F_{\widetilde X',\widetilde Y'}(R_{(0,1)})$
is not in the image of the functor $F_{\widetilde X,\widetilde Y}$. 
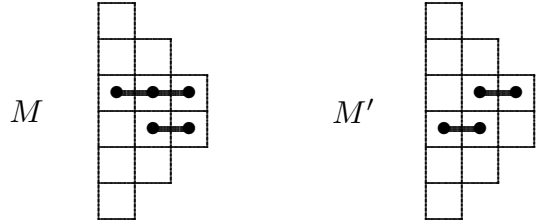
$$
{\beginpicture
    \setcoordinatesystem units <.4cm,.4cm>
\put{\beginpicture
\multiput{} at  0 0   3  6 /
\put{$M$} at -2 3
\plot 0 0  1 0  1 6  0 6  0 0 /
\plot 0 1  2 1  2 5  0 5 /
\plot 0 2  3 2  3 4  0 4 /
\plot 0 3  3 3 /

\multiput{$\bullet$} at 0.5 3.5  1.5 3.5  2.5  3.5  1.5 2.5  2.5  2.5 /
\plot 0.5 3.5  2.5 3.5 /
\plot 0.5 3.55 2.5 3.55 /
\plot 0.5 3.45 2.5 3.45 /
\plot 2.5 2.5   1.5 2.5 /
\plot 2.5 2.55  1.5 2.55 /
\plot 2.5 2.45  1.5 2.45 /
\endpicture} at 0 0 
\put{\beginpicture
\multiput{} at  0 0   3  6 /
\put{$M'$} at -2 3
\plot 0 0  1 0  1 6  0 6  0 0 /
\plot 0 1  2 1  2 5  0 5 /
\plot 0 2  3 2  3 4  0 4 /
\plot 0 3  3 3 /

\multiput{$\bullet$} at 1.5 3.5  2.5  3.5  0.5 2.5  1.5  2.5 /
\plot 1.5 3.5  2.5 3.5 /
\plot 1.5 3.55 2.5 3.55 /
\plot 1.5 3.45 2.5 3.45 /
\plot .5 2.5   1.5 2.5 /
\plot .5 2.55  1.5 2.55 /
\plot .5 2.45  1.5 2.45 /
\endpicture} at 9 0 
\endpicture}
$$
Note that both $M$ and $M'$ belong to $\widetilde{\Cal C}$.

Namely, the Auslander-Reiten component of $\widetilde{\Cal S}(6)$
which contains $M$ is a stable 3-tube $\widetilde{\Cal C}$, the 
quasi-socle of $M$ in $\widetilde{\Cal C}$ maps under the covering functor $\pi$ to $([3],[6],[3])$.
Also $M'$ belongs to this component $\widetilde{\Cal C}$, but 
its quasi-socle in $\widetilde{\Cal C}$ maps under $\pi$ to $([3],[3],[0])$, 
thus $M$ and $M'$ are not isomorphic (of course, one may also check directly that 
the modules $M$ and $M'$ are not isomorphic, but belong to the same Auslander-Reiten
component $\widetilde{\Cal C}$ of $\widetilde{\Cal S}(6)$).
Let us stress that given a Kronecker subcategory $\Cal K$ of $\widetilde {\Cal S}(6)$, the
simple regular objects of $\Cal K$ always belong to pairwise different 
Auslander-Reiten components of $\widetilde {\Cal S}(6)$. 
\end{remark}

\subsection{Some $\mathbb P^1$-families in $\Cal S(6)$ with width 6.}
\label{sec-two-eight}

Let us end this section with a slightly more complex situation: Simplification of a
subcategory $\Cal L$ of $\widetilde{\Cal S}(6)$
which contains four pairwise othogonal modules with endomorphism ring $k$. 
	\medskip

We consider the following four pairwise orthogonal modules $\widetilde A, 
\widetilde B, \widetilde B', \widetilde C$ in $\widetilde{\Cal S}(6)$
with endomorphism ring $k$:
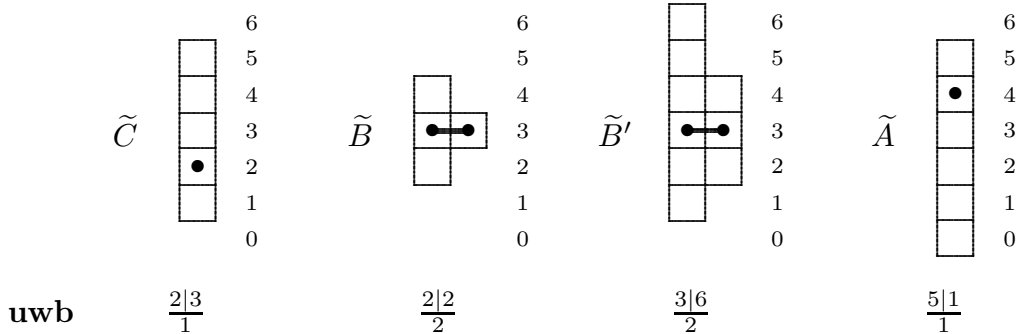
$$ {\beginpicture
    \setcoordinatesystem units <.4cm,.4cm>
\put{\beginpicture
\multiput{} at 0 0  2 7 /
\multiput{$\bullet$} at 0.5 2.5 /
\put{$\widetilde C$} at -1.5 3.5
\plot 0 1  1 1  1 6  0 6  0 1 /
\plot 0 2  1 2 /
\plot 0 3  1 3 /
\plot 0 4  1 4 /
\plot 0 5  1 5 /

\put{$\scriptstyle 0$} at 2 0.5 \put{$\scriptstyle 1$} at 2 1.5
\put{$\scriptstyle 2$} at 2 2.5 \put{$\scriptstyle 3$} at 2 3.5
\put{$\scriptstyle 4$} at 2 4.5 \put{$\scriptstyle 5$} at 2 5.5
\put{$\scriptstyle 6$} at 2 6.5
\endpicture} at 0 0
\put{\beginpicture
\multiput{} at 0 0  2 7 /
\plot 0 2  1 2  1 5  0 5  0 2 /
\plot 0 3  2 3  2 4  0 4 /
\multiput{$\bullet$} at 0.5 3.5  1.5  3.5 /
\plot 0.5 3.5  1.5  3.5 /
\plot 0.5 3.55  1.5  3.55 /
\plot 0.5 3.45  1.5  3.45 /
\put{$\widetilde B$} at -1.5 3.5

\put{$\scriptstyle 0$} at 3 0.5 \put{$\scriptstyle 1$} at 3 1.5
\put{$\scriptstyle 2$} at 3 2.5 \put{$\scriptstyle 3$} at 3 3.5
\put{$\scriptstyle 4$} at 3 4.5 \put{$\scriptstyle 5$} at 3 5.5
\put{$\scriptstyle 6$} at 3 6.5
\endpicture} at 7 0
\put{\beginpicture
\multiput{} at 0 0  2 7 /
\plot 0 1  1 1  1 7  0 7  0 1 /
\plot 0 2  2 2  2 5  0 5 /
\plot 0 3  2 3 /
\plot 0 4  2 4 /
\plot 0 6  1 6 /
\multiput{$\bullet$} at 0.5 3.5  1.5  3.5 /
\plot 0.5 3.5  1.5  3.5 /
\plot 0.5 3.55  1.5  3.55 /
\plot 0.5 3.45  1.5  3.45 /
\put{$\widetilde B'$} at -1.5 3.5

\put{$\scriptstyle 0$} at 3 0.5 \put{$\scriptstyle 1$} at 3 1.5
\put{$\scriptstyle 2$} at 3 2.5 \put{$\scriptstyle 3$} at 3 3.5
\put{$\scriptstyle 4$} at 3 4.5 \put{$\scriptstyle 5$} at 3 5.5
\put{$\scriptstyle 6$} at 3 6.5
\endpicture} at 14 0 
\put{\beginpicture
\multiput{} at 0 0  2 7 /
\multiput{$\bullet$} at 0.5 4.5 /
\put{$\widetilde A$} at -1.5 3.5
\plot 0 0  1 0  1 6  0 6  0 0 /
\plot 0 1  1 1 /
\plot 0 2  1 2 /
\plot 0 3  1 3 /
\plot 0 4  1 4 /
\plot 0 5  1 5 /

\put{$\scriptstyle 0$} at 2 0.5 \put{$\scriptstyle 1$} at 2 1.5
\put{$\scriptstyle 2$} at 2 2.5 \put{$\scriptstyle 3$} at 2 3.5
\put{$\scriptstyle 4$} at 2 4.5 \put{$\scriptstyle 5$} at 2 5.5
\put{$\scriptstyle 6$} at 2 6.5
\endpicture} at 21 0

\put{$\uwbb$} at -4 -5
\put{$\frac{2|3}1$} at 0 -5
\put{$\frac{2|2}2$} at 7 -5
\put{$\frac{3|6}2$} at 14 -5
\put{$\frac{5|1}1$} at 21 -5
\endpicture}
$$
The $\Ext$-quiver looks as follows:
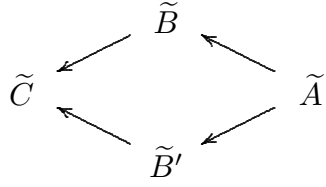
$$ {\beginpicture
    \setcoordinatesystem units <.4cm,.4cm>
\put{$\widetilde C$} at 0 0
\put{$\widetilde B$} at 4 2
\put{$\widetilde B'$} at 4 -2
\put{$\widetilde A$} at 8 0
\arr{3 1.5}{1 0.5}
\arr{7 0.5}{5 1.5}
\arr{3 -1.5}{1 -0.5}
\arr{7 -0.5}{5 -1.5}
\endpicture}
$$
and all extensions are g-split. (All these assertions are easy to verify.) Here are 
the corresponding non-trivial extensions:
$$ 
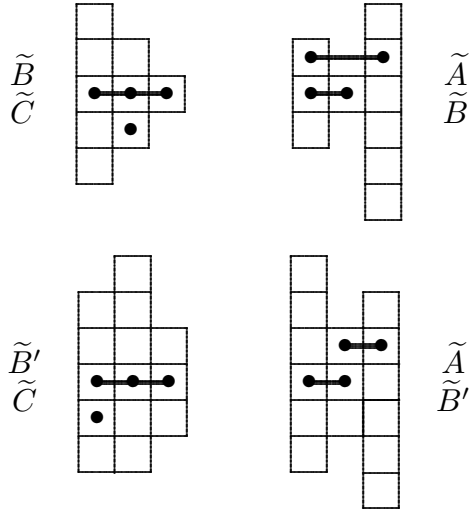
{\beginpicture
    \setcoordinatesystem units <.4cm,.4cm>
\put{\beginpicture
\multiput{} at 0 0  3 7 /
\put{$\begin{matrix} \widetilde B\cr \widetilde C\end{matrix}$} at -1.5 3.5
\plot 0 1  1 1  1 6  0 6  0 1 /
\plot 0 2  2 2  2 5  0 5 /
\plot 0 3  3 3  3 4  0 4 /
\multiput{$\bullet$} at 1.5 2.5  0.5 3.5  1.5 3.5  2.5 3.5 /
\plot 0.5 3.55  2.5 3.55 /
\plot 0.5 3.5  2.5 3.5 /
\plot 0.5 3.45  2.5 3.45 /

\endpicture} at 3 0
\put{\beginpicture
\multiput{} at 0 0  3 7 /
\plot 0 1  0 6  2 6 /
\plot 0 1  2 1  2 7  1 7  1 1 /
\plot 0 2  3 2  3 5  0 5 /
\plot 0 3  3 3 /
\plot 0 4  3 4 /
\multiput{$\bullet$} at 0.5 3.5  1.5  3.5  2.5 3.5  0.5 2.5 /
\plot 0.5 3.5  2.5  3.5 /
\plot 0.5 3.55  2.5  3.55 /
\plot 0.5 3.45  2.5  3.45 /
\put{$\begin{matrix} \widetilde B'\cr \widetilde C\end{matrix}$} at -1.5 3.5

\endpicture} at 3 -8
\put{\beginpicture
\multiput{} at 0 0  3 7 /
\plot 2 1  2 0  3 0  3 1 /
\plot 0 2  1 2  1 5  0 5  0 2 /
\plot 0 3  3 3 /
\plot 0 4  3 4 /
\plot 2 1  3 1  3 6  2 6  2 1 /
\plot 2 2  3 2 /
\plot 2 5  3 5 /
\multiput{$\bullet$} at 0.5 3.5  1.5 3.5  0.5 4.5  2.5 4.5 /
\plot 0.5 3.5  1.5  3.5 /
\plot 0.5 3.55  1.5  3.55 /
\plot 0.5 3.45  1.5  3.45 /
\plot 0.5 4.5  2.5  4.5 /
\plot 0.5 4.55  2.5  4.55 /
\plot 0.5 4.45  2.5  4.45 /
\put{$\begin{matrix} \widetilde A\cr \widetilde B\end{matrix}$} at 4.5 3.5

\endpicture} at 11 0 

\put{\beginpicture
\multiput{} at 0 0  3 7 /

\multiput{} at 0 0  2 7 /
\plot 0 1  1 1  1 7  0 7  0 1 /
\plot 0 2  2 2  2 5  0 5 /
\plot 0 3  3 3 /
\plot 0 4  3 4 /
\plot 0 6  1 6 /
\plot 3 3  2 3  2 0  3 0  3 6  2 6  2 5  3 5 /
\plot 2 1  3 1 /
\plot 2 2  3 2 /
\multiput{$\bullet$} at 0.5 3.5  1.5  3.5  1.5 4.5 2.5 4.5 /
\plot 0.5 3.5  1.5  3.5 /
\plot 0.5 3.55  1.5  3.55 /
\plot 0.5 3.45  1.5  3.45 /
\plot 1.5 4.5  2.5  4.5 /
\plot 1.5 4.55  2.5  4.55 /
\plot 1.5 4.45  2.5  4.45 /
\put{$\begin{matrix} \widetilde A\cr \widetilde B'\end{matrix}$} at 4.5 3.5
\endpicture} at 11 -8
\endpicture}
$$
(If $X, Y$ are objects and there is a unique indecomposable $M$
with a subobject $Y$ such that $M/Y = X,$ then we may write
$ M = \begin{matrix} X \cr Y \end{matrix}$.)

Let $\Cal L$ be the full subcategory of all objects which have a filtration with factors
of the form $\widetilde A, 
\widetilde B, \widetilde B',
\widetilde C$; the $\Ext$-quiver shows that $\Cal L$ is equivalent to an abelian 
full subcategory of
the category of representations of the quiver $\widetilde{\mathbf A}_{2,2}.$ 
Actually, as we will see, $\Cal L$ is hereditary, thus $\Cal L$ is equivalent to
the category of representations of the quiver $\widetilde{\mathbf A}_{2,2}.$ 

The object $\widetilde B$ belongs to the 6-tube $\widetilde{\Cal C}$ with rationality index 1, as does
the indecomposable object $\widetilde M$ with a filtration with  factors (going downwards) $\widetilde A,
\widetilde B', \widetilde C.$ It is easy to see that $\widetilde{\Cal C}$ contains infinitely many
objects which have a filtration with factors of the form $\widetilde B$ and $\widetilde M$. As a 
consequence, the category $\Cal L$ has infinitely many isomorphism classes of indecomposable objects. It follows that $\Cal L$ has to be hereditary. 

The object $\widetilde B'$ belongs to the 3-tube with rationality index 1, as does
the indecomposable object with a filtration with  factors (going downwards) $\widetilde A,
\widetilde B, \widetilde C.$
	\medskip

The $\Ext$-quiver shows that there are the following Kronecker pairs (below we note
the uwb-vectors):
$$ 
{\beginpicture
    \setcoordinatesystem units <.4cm,.5cm>
\put{$
\begin{matrix} \widetilde B\oplus \widetilde B' \cr \widetilde C\end{matrix} \ \leftleftarrows \ 
 \widetilde A$} at 0 0
\put{$
\begin{matrix} \widetilde B \cr \widetilde C\end{matrix} \ \leftleftarrows \
   \begin{matrix} \widetilde A \cr \widetilde B'\end{matrix} $} at 7 0
\put{$
\begin{matrix} \widetilde B' \cr \widetilde C\end{matrix} \ \leftleftarrows \  
   \begin{matrix} \widetilde A \cr \widetilde B\end{matrix}\ $}
 at 14 0
\put{$
\widetilde C \ \leftleftarrows \ 
\begin{matrix} \widetilde A\cr \widetilde B\oplus \widetilde B' \end{matrix} $} at 21 0

\put{$\uwbb$} at -5 -2
\put{$\frac{7|11}5$} at -1 -2
\put{$\frac{5|1}1$} at  2 -2

\put{$\frac{4|5}3$} at 6 -2
\put{$\frac{8|7}3$} at  8 -2

\put{$\frac{5|9}3$} at 13 -2
\put{$\frac{7|3}3$} at  15 -2

\put{$\frac{2|3}1$} at 19 -2
\put{$\frac{10|9}5$} at  22 -2

\endpicture}
$$
They give rise to four Kronecker subcategories
of $\Cal L.$ In this way, we have 
constructed four $\mathbb P^1$-families of indecomposable objects
with rationality index 1; the objects which we obtain and which belong to homogeneous tubes coincide.

For example, for the second pair $\widetilde X =  \begin{matrix} \widetilde A \cr \widetilde B'\end{matrix}$ and
$\widetilde Y = \begin{matrix} \widetilde B \cr \widetilde C\end{matrix}$,
we deal with the following 
$\mathbb P^1$-family:
$$
{\beginpicture
    \setcoordinatesystem units <.4cm,.4cm> \multiput{} at -2 0 7 7 /
\put{$F_{\widetilde X,\widetilde Y}(R_c)\:$\strut} at -3 3.5
\put{$\scriptstyle 0$} at 8 0.5 \put{$\scriptstyle 1$} at 8 1.5
\put{$\scriptstyle 2$} at 8 2.5 \put{$\scriptstyle 3$} at 8 3.5
\put{$\scriptstyle 4$} at 8 4.5 \put{$\scriptstyle 5$} at 8 5.5
\put{$\scriptstyle 6$} at 8 6.5

\plot 0 1 1 1 1 5 0 5 0 1 / 
\plot 0 2 2 2 2 4 0 4 0 3 2 3 / 
\plot 0 5 0 6 1 6 1 5 2 5 2 4 3 4 3 3 2 3 / 
\multiput{$\bullet$} at 0.3 3.3 1.3 3.3 2.3 3.3
  1.5 2.5 / 
\plot 0.3 3.3 2.3 3.3 / 
\plot 0.3 3.25 2.3 3.25 / 
\plot 0.3 3.35 2.3 3.35 /

\plot 4 1 5 1 5 7 4 7 4 1 / 
\plot 4 6 5 6 / 
\plot 6 0 7 0 7 6 6 6 6 0 /
\plot 6 1 7 1 / 
\plot 4 2 7 2 / 
\plot 4 3 7 3 / 
\plot 4 4 7 4 / 
\plot 4 5 7 5 / 
\multiput{$\bullet$} at 6.5 4.5 5.5 4.5 5.5 3.5 4.5 3.5 / 
\plot 6.5 4.5 5.5 4.5 / 
\plot 6.5 4.45 5.5 4.45 / 
\plot 6.5 4.55 5.5 4.55 /

\plot 5.5 3.5 4.5 3.5 / 
\plot 5.5 3.45 4.5 3.45 / 
\plot 5.5 3.55 4.5 3.55 /
\put{$\ss c_0$} at 1.5 4.5 
\put{$\ss c_1$} at 0.5 3.7 
\setdashes <.9mm>
\plot 1.8 4.5 5.5 4.5 / 
\plot 1.8 4.55 5.5 4.55 / 
\plot 1.8 4.45 5.5 4.45 /
\plot 0.8 3.7  2.8 3.7 4.5 3.5 / 
\plot 0.8 3.75 2.8 3.75 4.5 3.55 / 
\plot 0.8 3.65 2.8 3.65 4.5 3.45 /

\put{$\widetilde A$} at  6.5 -1
\put{$\widetilde B'$} at  4.9 -1
\put{$\widetilde B$} at 2 -1 
\put{$\widetilde C$} at  0.5 -1

\endpicture}
$$

\bigskip
\section{The functor $\tau^2$.}
\label{sec-three}

In this section
we work in $\Cal S(n)$ with $n$ fixed, thus we write $\tau = \tau_n.$

\medskip
The aim is to show Theorem \ref{theoremthree}, but also some consequences, and
related observations. Theorem \ref{theoremthree} will be derived
from Theorem~1.3$'$ in Section~\ref{theoremthree-prime}
which sheds more light
on the rotation $\rho$. Instead of looking at the triangle $\mathbb T(n)$ with its rotation $\rho$,
we will consider the following set $\mathbb E(n)$.

\subsection{The set $\mathbb E(n)$.}
\label{sec-three-one}

We consider triples $E = ([E_0],[E_1],[E_2])$ 
of isomorphism classes of modules $E_0, E_1,E_2$ in $\Cal N(n)$, and write
$|E| = |E_0|+|E_1|+|E_2|$. 

The set $\mathbb E(n)$ consists of the triples $E = ([E_0],[E_1],[E_2])$ 
of isomorphism classes of modules $E_0, E_1,E_2$ in $\Cal N(n)$ which
satisfy the condition that $n$ divides $|E|$. Instead of 
$([E_0],[E_1],[E_2])$, we usually will write $E = (E_0\backslash E_1\backslash E_2)$; the
number $bE = |E|/n$ will be called the {\it width} of $E$ (by assumption, 
this is a natural number).
Note that {\it for any $n$ and $b$, the set of elements $E$ in $\mathbb E(n)$ 
with fixed width $b$ is finite.}
Namely, if $E = (E_0\backslash E_1\backslash E_2)$
belongs to $\mathbb E(n)$ and has width $b$, then
$[E_i]$ is a partition of $|E_i|$, and we have $|E_i| \le bn$, for $0\le i \le 2$; 
of course, there are only finitely
many partitions of numbers bounded by $bn$.  

Let 
$$
 \rho (E_0\backslash E_1\backslash E_2) = (E_1\backslash E_2\backslash E_0),
$$ 
also, let 
$$
   \Omega (E_0\backslash E_1\backslash E_2) = 
  (\Omega E_0\backslash \Omega E_1\backslash \Omega E_2);
$$ 
If $E = (E_0\backslash E_1\backslash E_2)\in \mathbb E(n),$ then, of course, $\rho E$ 
belongs to $\mathbb E(n)$, but also $\Omega E$ belongs to $\mathbb E(n)$ (namely, for every
$V$ in $\Cal N(n)$, the number $|\Omega V|+|V| = |PV|$ is divisible by $n$, since
$PV$ is projective in $\Cal N(n)$).

Always, we have 

\begin{equation*}
  \label{equ-one}
   b(\rho E) = bE, \tag{1}
\end{equation*}
(but $b(\Omega E)$ usually is different from $bE$). 

To every object $X = (U,V) = (U,V,W)$ in $\Cal S(n)$, we attach the triple
$$
    EX = (\Omega V\backslash U\backslash W) \in \mathbb E(n);
$$
these three modules $\Omega V, U, W$ are just the factors of the filtration
$$
   0 \subseteq \Omega V \subseteq \widetilde U \subseteq PV 
$$
of $PV$, where $\widetilde U$ is defined by 
$\widetilde U/\Omega V = U \, \subseteq \, V = PV/\Omega V$
(the triple $EX$ may be called the {\it hidden filtration factors} of $X$). 
$$
\CD 
@. 0 @. 0 \\
@. @VVV @VVV @. @. \\
@. \Omega V @= \Omega V @. @. \\
@. @VVV @VVV @. @. \\
0 @>>> \widetilde U @>>> PV @>>> W @>>> 0 \\
@. @VVV @VVV @| @. \\
0 @>>> U @>>> V @>>> W @>>> 0 \\
@. @VVV @VVV @. @. \\
@. 0 @. 0 
\endCD
$$
This filtration of $PV$
shows that $|EX| = |\Omega V|+|U|+|W| = |PV|$ is divisible by $n$,
therefore $EX \in \mathbb E(n)$,   and that 
\begin{equation*}
  \label{equ-two}
 b(EX) = bV = bX. \tag{2}
\end{equation*}

\medskip
Theorem 1.3$'$ describes the hidden filtration factors of $\tau X$ and $\tau^2 X$ 
in terms of the hidden filtration factors of $X$. 

\medskip
\begin{theoremthree-prime}
  \label{theoremthree-prime}
  Let $X$ be a reduced object of $\Cal S(n)$. Then
  $$ 
  E(\tau^2 X) = \rho EX\quad \text{\it and}\quad E(\tau X) = \Omega\rho^2 EX.
  $$
\end{theoremthree-prime}

\medskip
It is the first assertion $E(\tau^2 X) = \rho EX$ which is crucial for all
our investigations and which immediately implies Theorem~\ref{theoremthree},
see Section~\ref{sec-three-two}. This is what the 
reader should keep in mind: If we apply $\tau^2$ to
any reduced object in $\Cal S(n)$, the three hidden 
filtration factors of $X$ are just rotated. A first look at an object 
$X = (U,V)$ in $\Cal S(n)$ (and its trace in $\mathbb T(n)$) concentrates the
attention to the two visible factors $U$ and $V/U$, but for a proper
understanding one should take into account also the additional factor $\Omega V$
of $PV$.  

\smallskip
\begin{proof}[Proof of Theorem~1.3$'$]
Given $V$ in $\Cal N(n)$, let $V'$ be obtained from
$V$ by deleting all direct summands of the form $[n].$ Thus
$V = V'\oplus P$, where $P = [n]^{bP}$ and $V'$ has Loewy length at most $n-1$.

\medskip 
We assume that $X$ is a reduced object in $\Cal S(n)$. Let $X = (U,V,W)$. 
Since we assume that $X$ is reduced, both $U$ and $W$
belong to $\Cal N(n-1).$ We decompose $V = V'\oplus P$ with $V'\in \Cal N(n-1)$
and $P$ a free $\Lambda$-module.
We recall from \cite[Theorem 5.1]{RS2}
the description of the Auslander-Reiten translation $\tau$ in $\Cal S(n)$:
$$
  \tau X = (V',W\oplus Q,\Omega U),
$$
where $Q$ is a suitable free $\Lambda$-module.
Since $\Omega P = 0,$
we have $\Omega V = \Omega (V'\oplus P) = \Omega(V').$
Applying $\tau$ twice, we obtain a corresponding description of $\tau^2 X.$

Altogether, we look at the following objects in $\Cal S(n)$:
\begin{align*}
   X &= (U,V,W) = (U,V'\oplus P,W), \cr
  \tau X &= (V',W\oplus Q,\Omega U), \cr
 \tau^2 X &= (W,\Omega U\oplus R, \Omega(V')) = (W,\Omega U\oplus R, \Omega V). \cr
\end{align*}
with $U,V',W\in \Cal N(n-1)$, 
and free $\Lambda$-modules $P,\ Q,\ R.$  
Since $U, V', W$ are in $\Cal N(n-1)$, we have 
$\Omega^2 U = U,\ \Omega^2 (V') = V' = \Omega^2V,\ \Omega^2 W = W$. 

\medskip
This provides the first assertion $E (\tau^2 X) = (E_1,E_2,E_0) = \rho EX.$
Namely, the subspace in $\tau^2X$ is $W$,
the corresponding factor space is $\Omega V$; the syzygy module $\Omega V(\tau^2X)$
of $V(\tau^2X)$ is $\Omega(\Omega U\oplus R) = \Omega^2U \oplus 0 \simeq U.$
	
It follows that $\rho^2 EX = 
(E_2\backslash E_0\backslash E_1) = (W\backslash \Omega V\backslash U).$
On the other hand, since $V' = \Omega^2V,$ we have  
$E(\tau X) = (\Omega W\backslash V'\backslash \Omega U) = 
(\Omega W\backslash \Omega^2V\backslash \Omega U)$.
Thus 
$E(\tau X) = \Omega\rho^2 EX.$
This is the second assertion. 
\end{proof}

\subsection{Proof of Theorem \ref{theoremthree}.}
\label{sec-three-two}

Given a non-zero element $E = (E_0\backslash E_1\backslash E_2)$
in $\mathbb E(n)$, let 
$$  
   pE = |E_1|/bE,\quad rE = |E_2|/bE,\quad \text{and also}\quad qE = pE+rE
$$
(note that $bE$ is non-zero, since we assume that $E$ is non-zero). Of course, $\pr E =
(pE,rE)$ will be called the {\it pr-vector} of $E$.
We have 
\begin{equation*}
 \pr \rho E = \rho \pr E, \tag{3}
\end{equation*}
for all $E\in \mathbb E(n).$

\smallskip
\begin{proof}[Proof of (3)]
  Let $E = (E_0\backslash E_1\backslash E_2).$ 
  Since $|E_0|+|E_1|+|E_2| = n\cdot bE,$ we have $|E_0|/bE = n-|E_1|/bE - |E_2|/bE = n -pE-rE.$
  Therefore $\pr\rho E 
  = \pr (E_1\backslash E_2\backslash E_0) = (|E_2|/bE,|E_0|/bE) = (rE,n-pE-rE) = \rho \pr E.$ 
\end{proof}

\smallskip
Of course, for any non-zero
object $X$ in $\Cal S(n)$, we have 
$pEX = pX,\  rEX = rX$ (and therefore also
$qEX = pEX + rEX = pX + rX = qX$), thus
\begin{equation*}
 \pr EX = \pr X .\tag{4}
\end{equation*}
	
	\medskip 
        Now, let $X$ be reduced in $\Cal S(n)$ and non-zero. Then Theorem~1.3$'$
        in Section~\ref{theoremthree-prime} shows that 
$$
 \pr \tau^2X \underset{(4)}= \pr E(\tau^2X) = \pr \rho EX 
  \underset{(3)}= \rho\pr EX \underset{(4)}= \rho \pr X.
$$
Similarly, we have 
$$
   b (\tau^2X) \underset{(2)}= b E(\tau^2X) = b \rho EX 
  \underset{(1)}= b EX \underset{(2)}= b X.
$$
This completes the proof of Theorem~\ref{theoremthree}. 
$\s$

\subsection{Remark.}
\label{sec-three-three}

Two special cases for the assertion $b(\tau^2 X) = bX$ for $X$ reduced 
should be kept in mind: {\it If $X$ is a
reduced picket, then $\tau^2 X$ is a reduced picket.} Namely, $X$ is a picket
if and only if $bX = 1$. Similarly, {\it if $X$ is a bipicket, then $\tau^2 X$ is a bipicket.} 
Section~\ref{sec-twelve} will provide further information on pickets and bipickets; in particular,
it describes the $\tau$-orbits of pickets.

\subsection{Average of the values of $p$, $q$ and $r$.}
\label{sec-three-four}

\begin{corollary}
  \label{cor-three-four}
  If $X$ is a reduced object of $\mathcal S(n)$, then
  $$
  pX+p(\tau_n^2X)+p(\tau_n^4X) = n, \quad\text{\it and}\quad 
  rX+r(\tau_n^2X)+r(\tau_n^4X) = n.
  $$
\end{corollary}

Note that the  Corollary asserts that the average of both of the values of $p$ and $r$
on the $\tau^2$-orbit of any reduced object of $\Cal S(n)$ is $\frac n3$.
As a consequence, the average of $q = p+r$
on the $\tau^2$-orbit of any reduced object is $\frac {2n}3$.
Thus, these numbers are also the average 
of the values of $p,\ q,\ r$ for the corresponding $\tau$-orbits.	

\smallskip
\begin{proof}
  We have $p(\tau^2X) = rX$ and
$r(\tau^2X) = n-pX-rX$. Therefore $p(\tau^4X) = r(\tau^2X) = n-pX-rX$
and $r(\tau^4X) = n-p(\tau^2X)-r(\tau^2X) = n-rX-n+pX+rX = pX.$
It follows that
\begin{align*}
 pX+p(\tau^2X)+p(\tau^4X) &= pX+rX+n-pX-rX = n, \cr
 rX+r(\tau^2X)+r(\tau^4X) &= rX+n-pX-rX+pX = n.
\end{align*}
\end{proof}

\subsection{The uwb-vectors and the triples $(|E_0|,|E_1|,|E_2|).$}
\label{sec-three-five}

Let us stress that our use of uwb-vectors as introduced in Section~\ref{sec-one-four} refers 
to the triples $E = (E_0,E_1,E_2).$ Namely, being interested in such triples,
it is natural to consider the corresponding dimension vectors $(|E_0|,|E_1|,|E_2|)$
in $\mathbb R^3$ and the uwb-vectors are just obtained by the following base change:
$$ 
  u = |E_1|,\ w = |E_2|,\ b = \tfrac1n(|E_0|+|E_1|+|E_2|).
$$ 

\medskip
Why do we prefer to work with uwb-vectors instead of the triples
$(|E_0|,|E_1|,|E_2|)$? Starting with an object $X$ in  $\Cal S(n)$, we may (and often will)
consider $X$ also as an object in $\Cal S(m)$ with $m\ge n$ and its uwb-vector
will remain the same!

\medskip
{\bf The boundary distance.}
Given $E = (E_0,E_1,E_2)$ in $\mathbb E(n)$, let 
$$
 mE = \min(|E_0|,|E_1|,|E_2|).
$$
Given $X\in \Cal S(n)$, let $mX = m(EX);$ then {\it $dX = mX/bX$ is just
the boundary distance of $\pr X$} in $\mathbb T(n)$, as defined in Section~\ref{sec-one-four}.

We may say that $E = (E_0,E_1,E_2)\in \mathbb E(n)$ is {\it central} provided 
$|E_0| = |E_1| = |E_2|$, or, equivalently, provided $mE = n/3.$
We recall that an indecomposable object $X$ in $\Cal S(n)$
is said to be central provided its pr-vector is $z(n) = (n/3,n/3),$ thus $X$
is central if and only if $EX$ is central.

\medskip 
{\bf $u$-minimal objects.}
If $E = (E_0,E_1,E_2)$ belongs to $\mathbb E(n)$, we say that $E$ is 
{\it $u$-minimal} provided $|E_1| \le |E_2|$ and $|E_1| < |E_0|.$ Thus 
$E$ is $u$-minimal if and only if the pr-vector $\pr E$ of $E$ satisfies the inequalities
$pE \le rE$ and $2pE < n-rE$ (so that $\pr E$ lies in the following  
shaded part of $\mathbb T(n)$; the black square marks the center $z(n)$):
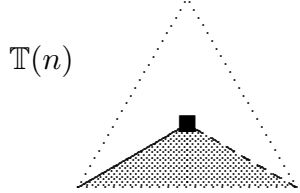
$$  
{\beginpicture
   \setcoordinatesystem units <.404cm,.7cm>
\multiput{} at -3 0  3 3 /
\setdots <1mm>
\plot -3 0  3 0  0 3  -3 0 /
\put{$\ss\blacksquare$} at 0 1
\setsolid
\plot -3 0  0 1 /
\setdashes <1mm>
\plot 3 0  0 1 /
\setshadegrid span <.3mm>
\vshade -3 0 0 <z,z,,> 0 0 1 <z,z,,> 3 0 0 /
\put{$\mathbb T(n)\:$} at -4 2
\endpicture}
$$

If $X$ in $\Cal S(n)$ is indecomposable, we say that $X$ is {\it $u$-minimal} provided
$EX$ is $u$-minimal. 
Of course, if $X$ is $u$-minimal, then $X$ is not central;
also, $X$ is not projective.

\medskip
\begin{lemma}
  \label{lem-three-one}
        If $X$ in $\mathcal S(n)$ is indecomposable, and neither central nor projective, then
          precisely one of the objects $X$, $\tau^2X$, $\tau^4X$ is $u$-minimal.
\end{lemma}

\begin{proof}
This follows immediately from Theorem~1.3$'$ in Section~\ref{theoremthree-prime}.
\end{proof}

\medskip The next result states that the indecomposable objects $X$ with $mX=0$
are exactly the boundary pickets.

\medskip
\begin{lemma}
  \label{lem-three-two}
  Let $X$ be indecomposable and $EX = (E_0\backslash E_1\backslash E_2)$.
  One of the $E_i$ is zero, if and only if $X$ is a boundary picket.
\end{lemma}

\begin{proof}
  If $X$ is a boundary picket, then either $uX = 0$, thus $E_1 = 0;$
  or $w = 0,$ thus $E_2 = 0$, or $u+w = n$, and then $E_0 = 0.$ Conversely, assume that
  $E_i = 0$ for some $i$. If $E_1 = 0,$ then $uX = 0$ and then clearly $X$ has to be 
  a picket. Dually, if $E_2 = 0,$ then $wX = 0$ and $X$ is a picket. Finally, assume
  that $E_0 = 0.$ Then $V$ is a projective $\Lambda$-module. But it is well-known that 
  the only 
  indecomposable objects $(U,V)$ with $V$ being projective are the pickets.
\end{proof}

\subsection{Barycentric coordinates.}
\label{sec-three-six}
  
This seems to be the right place to insert an important, however obvious, remark. 
The functions $p$ and $r$, the level and the colevel, are (of course) part of the barycentric
coordinate system for the triangle $\mathbb T(n)$,
the additional coordinate is the function $n-q = n-(p+r)$
(thus corresponds to the mean).
$$  
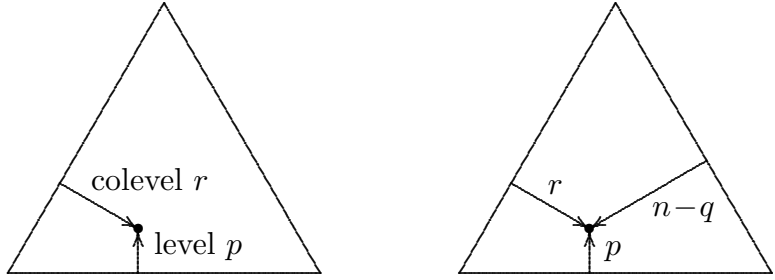
{\beginpicture
   \setcoordinatesystem units <1cm,1cm>
\put{\beginpicture
   \setcoordinatesystem units <.57735cm,1cm>
\multiput{} at -3 0  3 3 /
\plot -3 0  3 0  0 3  -3 0 /
\put{$\ss\bullet$} at -.5 .5 
\arr{-.5 0}{-.5 .45}
\arr{-2 1}{-.54 .51}
\put{colevel $r$} [l] at -1.4 1.05
\put{level $p$} [l] at -.2 0.3
\endpicture} at 0 0
\put{\beginpicture
   \setcoordinatesystem units <.57735cm,1cm>
\multiput{} at -3 0  3 3 /
\plot -3 0  3 0  0 3  -3 0 /
\put{$\ss\bullet$} at -.5 .5 
\arr{-.5 0}{-.5 .45}
\arr{-2 1}{-.54 .51}
\arr{1.75 1.25}{-.46 .51}
\put{$r$} [l] at -1.3 .95
\put{$p$} [l] at -.2 0.25
\put{$n\!-\!q$} [l] at .7 .7
\endpicture} at 5 0

\endpicture}
$$
In terms of $\mathbb E$, the three functions $p,\ r,\ n-q$ are just the functions $|E_1|/b,\
|E_2|/b$ and $|E_0|/b.$	
Given an object $X = (U,V) = (U,V,W)$ in $\Cal S(n)$, we have attached the triple
$$
    EX = (\Omega V\backslash U\backslash W) \in \mathbb E(n);
$$
Of course, the last two coordinates $p$ and $r$ are given by $p(EX) = pX,$ and $r(EX) = rX$.
The first coordinate $n-q = n-p-r$ is just $\omega/b$, where $\omega X = |\Omega V|.$
The relevant formula to have in mind is:
$$
 \omega+u+w = nb,
$$
which is valid for all objects $X$ in $\Cal S(n)$. 

\medskip
Using $\omega$ as additional coordinate, Theorem~\ref{theoremone} and its accompanying theorems
in Sections~\ref{sec-one-two} and \ref{sec-one-three} can be rewritten as follows.

\medskip
\begin{theoremone-again}
  \label{theoremone-again}
    Let $X$ be an indecomposable object in $\Cal S(n)$.
  Then either $X$ is a boundary picket, or else $u \ge b$, $w\ge b$, and $\omega \ge b$.$\s$
\end{theoremone-again}

\subsection{Central objects.}
\label{sec-three-seven}

\begin{proposition}
  \label{prop-three-one}
  Let $X = (U,V)$ be indecomposable. Then the following conditions are 
  equivalent:
  \begin{itemize}[leftmargin=3em]
  \item[\rm(i)] $X$ is central.
  \item[\rm(ii)] $uX = wX = |\Omega V|$.
  \item[\rm(iii)] $uX = u(\tau^2X),\ wX = w(\tau^2X).$
  \end{itemize}
\end{proposition}

\begin{proof}
  First, assume that $X$ is projective. Then $X$ cannot be central. Also, $\Omega V = 0,$
thus condition (ii) is not valid ($uX = wX = 0$ implies that $X = 0$); similarly,
(iii) is not valid. 

Thus, we can assume that $X$ is not projective and therefore
$\tau^2X = (W,V(\tau^2X),\Omega V).$ Therefore $u(\tau^2X) = wX$ and 
$w(\tau^2X) = |\Omega V|.$

Now $X$ is central if and only if $EX = (\Omega V, U, W)$ is central, thus if and only if condition (ii)
is satisfied. It remains to show the equivalence of (ii) and (iii).
Since $u(\tau^2X) = wX$, the equality $uX = wX$ is the same as $uX = u(\tau^2X).$
Since $w(\tau^2X) = |\Omega V|,$ the equality $wX = |\Omega V|$ is the same as 
$wX = w(\tau^2X).$ 
\end{proof}

\medskip
We add here some remarks concerning the Auslander-Reiten components of $\Cal S(n).$ 

\medskip
As a consequence of Proposition~\ref{prop-three-one}, we see that any tube $\Cal C$ 
contains infinitely many central objects.
Namely, for $\Cal C$ being stable, there is the following assertion: 

\medskip
\begin{corollary}
  \label{cor-three-one}
  Let $X$ be indecomposable. If $\tau^2 X = X,$ then $X$ is central
  (thus, if $X$ belongs to a tube of rank $1$ or $2$, then $X$ is central). 
  If $X$ belongs to a stable tube of rank $r = 3$ or $r = 6,$ and
  has quasi-length divisible by $r$, then $X$ is central.
\end{corollary}

\begin{proof}
  We use the implication (iii) $\implies$ (i). The first assertion is a direct
  consequence. For the second assertion, assume that $\Cal C$ has rank $r = 3$ or $6$.
  Let $X$ be in $\Cal C$ with quasi-length $rt$. Let $Z$ be quasi-simple in $\Cal C$.
  Then $X$ has a filtration with factors $\tau^iZ$, where $0\le i < rt$ and $Z$
  is quasi-simple. But $\tau^2X$ has a similar filtration, with factors being just
  rotated. Thus $uX = u(\tau^2X)$ and $wX = w(\tau^2X).$
\end{proof}

\medskip
The discussion of $\Cal P(n)$ in Section~\ref{sec-nine-five} shows that for $\Cal P(n)$ with $n\ge 6$,
and $t\in \mathbb N_1$, precisely one of the two $\tau^2$-orbits of objects of quasi-length
$6t$ consists of central objects, and these are the only ones in $\Cal P(n).$

\medskip
\begin{remark}
  There are additional central objects.
\end{remark}

\smallskip
The bipicket $X$
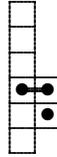
$$
{\beginpicture
  \setcoordinatesystem units <.28cm,.28cm>
\multiput{} at 0 0  2 6 /
\plot 0 0  1 0  1 6  0 6  0 0 /
\plot 0 1  2 1  2 3  0 3 /
\plot 0 4  1 4 /
\plot 0 2  2 2 /
\plot 0 5  1 5 /
\multiput{$\bullet$} at  0.5 2.5  1.5 2.5  1.5 1.5  /
\plot 0.5 2.5  1.5 2.5 /
\plot 0.5 2.45  1.5 2.45 /
\plot 0.5 2.55  1.5 2.55 /
\endpicture}
$$
has uwb-vector $\frac{4|4}2$, thus $X$ is central in $\Cal S(6)$. We have $X = \tau_6 E_2^2$
(and $E_2^2$ is not central in $\Cal S(6)$). Here, we deal with objects of quasi-length 2
in the 6-tube 
of $\Cal S(6)$ with rationality index $1$: Three of these objects are central, the
other three not. 

\bigskip
\centerline{$*\ *\ *$}

\medskip 
The main assertions of this section concern $\tau^2.$
Let us end by looking at the action of $\tau$ itself.
First, let us point out that 
there are some pecularities.

\subsection{Warnings.}
\label{sec-three-eight}

\phantom m
\vspace{-3mm}

{\bf (1)} 
{\it If $X, X'$ are indecomposable objects in $\Cal S(n)$ with $\pr X =
\pr X'$, we usually have $\pr \tau_n X \neq \pr \tau_n X'.$}

\medskip
\begin{example}[for $n= 4$]
  There are two indecomposable objects in $\Cal S(4)$
  with pr-vector $(1,1),$
  namely $X = ([1],[2])$ and $X' = E_2^2.$ We have $\tau_4 X = E_2^3$ with 
  $\pr \tau_4 X = (1,3/2),$
  whereas the global space of $\tau_4 X'$ is $[4,2]$, its subspace is $[3,1]$, its factor space 
  is $[2],$
  thus $\pr\tau_4 X' = (2,1).$
\end{example}

\medskip
{\bf (2)} {\it If $X, X'$ are indecomposable objects in $\Cal S(n)$ with $\bdim X =
\bdim X'$, and $bX = bX',$ we may have $\bdim \tau_n X \neq \bdim \tau_n X'$, $b\tau_n X \neq b\tau_n X'$
as well as $\pr \tau_n X \neq \pr \tau_n X'.$}

\medskip
\begin{example}[for $n = 6$]
  We consider three indecomposable objects with the same
  global space $V = [6,4,2],$ using subspaces $U, U', U''$ of dimension 6. 
$$
{\beginpicture
    \setcoordinatesystem units <.28cm,.28cm>
\put{\beginpicture
\multiput{} at 0 0  3 8 /
\plot 0 0  1 0  1 6  0 6  0 0 /
\plot 0 1  2 1  2 5  0 5 /
\plot 0 2  3 2  3 4  0 4 /
\plot 0 3  3 3 /
\multiput{$\bullet$} at  0.5 3.5  2.5 3.5  1.5 2.5  2.5 2.5  /
\plot .5 3.5  2.5 3.5 /
\plot .5 3.45  2.5 3.45 /
\plot .5 3.55  2.5 3.55 /
\plot 1.5 2.5  2.5 2.5 /
\plot 1.5 2.45  2.5 2.45 /
\plot 1.5 2.55  2.5 2.55 /

\put{$X = (U,V,W)$} at 1 -1.5
\endpicture} at 0 0
\put{\beginpicture
\multiput{} at 0 0  3 8 /
\plot 0 0  1 0  1 6  0 6  0 0 /
\plot 0 1  2 1  2 5  0 5 /
\plot 0 2  3 2  3 4  0 4 /
\plot 0 3  3 3 /
\multiput{$\bullet$} at  0.5 3.5  1.5 3.5  1.5 2.5  2.5 2.5  /
\plot .5 3.5  1.5 3.5 /
\plot .5 3.45  1.5 3.45 /
\plot .5 3.55  1.5 3.55 /
\plot 1.5 2.5  2.5 2.5 /
\plot 1.5 2.45  2.5 2.45 /
\plot 1.5 2.55  2.5 2.55 /

\put{$X' = (U',V,W')$} at 1 -1.5
\endpicture} at 11 0
\put{\beginpicture
\multiput{} at 0 0  3 8 /
\plot 0 0  1 0  1 6  0 6  0 0 /
\plot 0 1  2 1  2 5  0 5 /
\plot 0 2  3 2  3 4  0 4 /
\plot 0 3  3 3 /
\multiput{$\bullet$} at  0.5 3.5  1.5 3.5  2.5 3.5  1.5 1.5  2.5 2.5  /
\plot .5 3.5  2.5 3.5 /
\plot .5 3.45  2.5 3.45 /
\plot .5 3.55  2.5 3.55 /

\put{$X'' = (U'',V,W'')$} at 1 -1.5
\endpicture} at 22 0
\endpicture}
$$
Here, $U \simeq U' \simeq [4,2] \simeq W \simeq W'',$ and $U'' \simeq W' \simeq [4,1,1].$
Obviously, we have $\bdim X = \bdim X' = \bdim X'' = (6,12),$ and $bX = bX' = bX'' = 3.$

We get 
$ \bpar\tau_6 X = ([4,2],[6,4,2],[4,2]),$
  $\bpar\tau_6 X' = ([4,2],[6,4,1,1],[4,2]),$ and finally
  $\bpar\tau_6 X'' = ([4,2],[6,6,4,2],[5,5,2]).$
It follows that $\bdim \tau_6 X =  \bdim\tau_6 X_n' = (6,12),$ whereas  $\bdim \tau X_6'' = (6,18).$
Also, $bX = 3,$ but 
$bX' = bX'' = 4.$ Thus $\pr \tau_6 X = (2,2),$ $\pr \tau_6 X' = (3/2,3/2)$
and $\pr \tau_6 X'' = (3/2,3).$ 
\end{example}

\medskip
\begin{remark}
  \label{rem-three-one}
  The partition vector $\bpar(\tau_n X)$ is determined by $\bpar X$.
\end{remark}

\begin{proof}
  Let $X$ be an object in $\Cal S(n)$ with $X = (U,V'\oplus P,W)$,
  where $V'$ has height at most $n-1$ and $P$ is a free $\Lambda$-module. 
  Thus $\bpar X = ([U],[V'\oplus P],[W])$, with $[V']$ and $[P]$ both being
  determined by $[V]$ (and $n$).
  As we know, $\bpar(\tau_nX) = ([V'],[W]\oplus [Q],[\Omega U])$, thus it remains
  to see that $[Q]$ is determined by $\bpar X.$ But we have
  $|Q|=|V'|+|\Omega U| - |W|=|V'|+|\Omega U|-|V|+|U|=n\cdot bU-|P|.$
\end{proof}

\medskip
\begin{remark}
  \label{rem-three-two}
  If the indecomposable object $X$ occurs in a stable tube and has even quasi-length,
  then $b\tau_nX=bX$.
\end{remark}

\begin{proof}
  Write $X=Z[\ell]$ where $Z$ is on the mouth of the tube and
  the quasi-length $\ell$ of $X$ is even.
  Using the additivity of $b$ on Auslander-Reiten sequences in stable tubes (see 9.4),
  we have $bX=\sum_{i=0}^{\ell-1}b\tau_n^{-i}Z=\sum_{i=0}^{\ell-1} b\tau_n^{1-i}Z=b\tau_n X$
  since $b\tau_n^{-\ell+1}Z=b\tau_nZ$ by Theorem~\ref{theoremthree}, using that $\ell$ is even.
\end{proof}

\bigskip
There is the following general observation concerning $\tau_n X$. 

\subsection{The multiplicity of $\Lambda$ in $V$.}
\label{sec-three-nine}

\begin{proposition}
  \label{prop-three-nine}
  Let $X = (U,V,W)$ be any object of $\Cal S(n)$. Then
$$
  b(\tau_n X) \le bU + bW,
$$
  with equality if and only if $X$ has height at most $n-1.$
\end{proposition}

\medskip
The difference between $bU+bW$ and $b(\tau_nX)$ is of interest (it is the Krull-Remak-Schmidt
multiplicity of $\Lambda$ in the $\Lambda$-module $V$); this will be further discussed
in Section~\ref{sec-three-twelve}.

\smallskip
\begin{proof}
  We can assume that $X$ is indecomposable. If $X$ is projective, then $\tau_nX = 0,$
  thus nothing has to be shown.  

  Therefore, we can assume that $X$ is reduced. As in the proof of Theorem~1.3$'$
  in Section~\ref{theoremthree-prime}, we have 
$V = V'\oplus P$, and 
$\tau_nX = (V',W\oplus Q,\Omega U)$, 
where $V'$ has height at most $n\!-\!1$ and 
where $P$ and $Q$ are free $\Lambda$-modules. 

The exact sequence $0 \to U \to V'\oplus P \to W \to 0$ yields
$|U| = |V'| + |P| - |W|$. The exact sequence $0 \to V' \to W\oplus  Q \to 
\Omega U \to 0$ yields
$|\Omega U| = |W| + |Q| - |V'|$. Thus we obtain 
$$
 n\cdot bU = |U|+|\Omega U| = |P|+|Q| = n\cdot bP + n\cdot bQ.
$$
Therefore 
\begin{equation*}
  bU = bP+bQ. \tag{$*$}
\end{equation*}
Since $bP \ge 0$, we have $bQ \le bU.$ Since the global space of $\tau_n X$
is $W\oplus Q,$ we see that $b(\tau_nX) = bW+bQ \le bW+bU.$ 

We have
\begin{equation*}
 bU+bW-b(\tau_nX) = bU+bW - bW-bQ = bU - bQ = bP. \tag{$**$}
\end{equation*}
Of course, $bP = 0$ if and only if $X$ belongs to $\Cal S(n-1).$
\end{proof}

\subsection{The width of $\tau_nX$.}
\label{sec-three-ten}

\begin{proposition}
  If $X$ in $\Cal S(n)$ is indecomposable and not projective, then
$$
 \tfrac12 bX \le b(\tau_n X) \le 2 bX.
 $$
 \end{proposition}

\begin{proof}
  Let $X = (U,V,W).$ According to Section~\ref{sec-three-nine}, we have $b(\tau_n X) \le bU+bW$. Since both $U$ and $W$
are subquotients of $V$, we have $bU \le bV$ and $bW \le bV$. Therefore
$$
 b(\tau_n X) \le  bU+bW \le  2bX.
$$

If we consider $\tau_n X$ instead of $X$, we see that $b(\tau_n^2 X) \le 2b(\tau_n X).$ But 
since $X$ is indecomposable and not projective, we have $bX = b(\tau_n^2X)$, thus
$bX = b(\tau_n^2X ) \le 2b(\tau_n X)$, and therefore $\frac12 bX \le b(\tau_n X).$
\end{proof}

\medskip
\begin{example}
{\it The inequality is optimal.} For $n\ge 3$ 
it is easy to construct indecomposable objects $X = (U,V,W)$
in $\Cal S(n-1)$ with $bU = bW = bX$ and then $b(\tau_{n}X) = bU+bW = 2bX.$

For example, there is the following object $X$ in $\Cal S(8)$ with $bX = 3$
$$
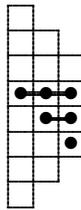
{\beginpicture
    \setcoordinatesystem units <.28cm,.28cm>
\multiput{} at 0 0  3 8 /
\plot 0 0  1 0  1 8  0 8  0 0 /
\plot 0 1  2 1  2 7  0 7 /
\plot 0 2  3 2  3 6  0 6 /
\plot 0 3  3 3 /
\plot 0 4  3 4 /
\plot 0 5  3 5 /
\multiput{$\bullet$} at  0.5 4.5 
                  2.5 4.5  1.5 3.5  2.5  3.5  2.5 2.5 /
\plot 0.5 4.5  2.5 4.5 /
\plot 0.5 4.45  2.5 4.45 /
\plot 0.5 4.55  2.5 4.55 /

\plot 1.5 3.5  2.5 3.5 /
\plot 1.5 3.45  2.5 3.45 /
\plot 1.5 3.55  2.5 3.55 /
\endpicture}
$$
with the partition vector
$$
 [531]\quad [864]\quad [531],
$$
and $\tau_9 X$ has the partition vector
$$
 [864]\quad [999531]\quad[864].
$$
\end{example}

\begin{remark}
If $X$ is a bipicket in $\Cal S(n)$ with $n\ge 6,$ then all
cases $b(\tau_n X) = 1,2,3,4$ are possible.
\end{remark}

\subsection{A formula for $bX+b(\tau_nX)$.}
\label{sec-three-eleven}

The following formula for $bX+b(\tau_n X)$ is valid for all objects $X$ in
$\Cal S(n)$.

\medskip
\begin{proposition}
  \label{prop-three-eleven}
  Let $X = (U,V,W)$ be an object in $\Cal S(n)$. Then
  $$
  bX+b(\tau_n X) = b\Omega V+bU+bW.
  $$
\end{proposition}

\begin{proof}
  Here we use for $X$ in $\Cal S(n)$ the definition of $EX$
  and the formula $E(\tau_n X) = \Omega\rho^2EX$ given in Theorem~1.3$'$
  in Section~\ref{theoremthree-prime}.
  \begin{align*}
    n(b\Omega V+bU+bW) &= |P\Omega V|+|PU|+|PW| \cr
    &= \left(
    |\Omega V| + |U| + |W|\right) + 
    \left(|\Omega^2V| + |\Omega U| + |\Omega W|\right)\cr
    & = |EX| + |E(\tau_n X)| \cr
    &=
    n\cdot bX + n\cdot b(\tau_n X).
  \end{align*}
\end{proof}

\bigskip
\centerline{$*\ *\ *$}

\medskip 
It seems to be a good place to introduce here a further invariant for the objects in $\Cal S(n)$.
We denote it by $c_n$. In contrast to the invariants $u,v,w,b$, this new
one depends not only on the object $X \in \Cal S$, but also on the height $n$ of $X$. 
As $u,v,w,b$, the values of $c_n$ are non-negative integers. Whereas $u,v,w,b$ vanish only on
the zero object, we have $c_nX = 0$ if and only if $X$ has height at most $n-1$. 

\subsection{The invariant $c_n$.}
\label{sec-three-twelve}

The new invariant $c_n$ will be defined first for the objects
$V$ in $\Cal N(n)$ (as we did when we defined the invariant $b$). 
Let $V$ be a module in $\Cal N(n)$. Then $c_nV$ is defined to be 
the dimension of $T^{n-1}V$, thus $c_nV$ 
is the number of Jordan blocks of size $n$ of the operator $T$, or, equivalently, the
Krull-Remak-Schmidt multiplicity of $\Lambda$ in $V$. Thus, if we
write $V = V'\oplus P$, where $V'$ has height at most $n-1$ and $P$
is a free $\Lambda$-module, then $c_n$ is the rank of $P$.
If $X = (U,V)$ is an object in $\Cal S(n)$, we set $c_nX = c_nV$.

\medskip
The following formulae rewrite some previous observation.

\medskip 
\begin{proposition}
  \label{prop-three-twelve-one}
  Let $X$ be any object in $\Cal S(n)$. Then
$$
 b(\tau_nX)= bU+bW- c_nX. 
 $$
\end{proposition}

\begin{proof}
  The formula holds for projective objects and has been shown for
  reduced objects in the proof of Proposition~\ref{prop-three-nine}, see  $(**)$.
\end{proof}

\medskip
\begin{corollary}
  Let $X = (U,V,W)$ be a reduced object in $\Cal S(n)$ of height $n$. 
  Then for all $t\ge 1$, we have
$$
  b(\tau_{n}X) <  bU+bW = b(\tau_{n+t}X).
  $$
\end{corollary}

\begin{proof}
  Since $X$ has height $n$, we have $c_nX > 0$, thus
  Proposition~\ref{prop-three-twelve-one} asserts that $b(\tau_{n}X) < bU+bW.$
  For every $t\ge 1$, the object $X$ is reduced in $\Cal S(n+t)$. Since its height is
  smaller than $n+t$, Proposition~\ref{prop-three-twelve-one}
  asserts that $b(\tau_{n+t}X) = bU+bW.$
\end{proof}

\medskip
There is also the following assertion.

\medskip
\begin{proposition}
  \label{prop-three-twelve-two}
  If $X= (U,V)$ is a reduced object in $\Cal S(n)$, then
$$
 bU = c_nX + c_n(\tau_nX).
 $$
\end{proposition}

\begin{proof}
  See the proof of Proposition~\ref{prop-three-nine}, $(*)$.
\end{proof}

\medskip
\begin{corollary}
  If $X= (U,V)$ is a reduced object in $\Cal S(n)$, then
  $$
  c_nX \le bU,\quad\text{\it and}\quad c_nX \le bW.
  $$
\end{corollary}

\begin{proof}
  The first assertion is a direct consequence of Proposition~\ref{prop-three-twelve-two},
  since $c_n(\tau_nX) \ge 0.$
  The second assertion follows by duality.
\end{proof}

\medskip
\begin{proposition}
  \label{prop-three-twelve-three}
  If $X= (U,V)$ is a reduced object in $\Cal S(n)$, then
  $$
  c_n(\tau_n^2X) = bX - bU.
  $$
\end{proposition}

\begin{proof}
  We have mentioned in the proof of Theorem~1.3$'$ in Section~\ref{theoremthree-prime}
  that for $X = (U,V,W)$
  reduced, we have $\tau_n^2 X = (W,\Omega U\oplus R,\Omega V)$ with $R$ a free $\Lambda$-module.
  Since $\Omega U$ has height at most $n-1$, we see that $c_n(\tau_n^2 X)$ is the rank of $R$,
  therefore $b X = b(\tau_n^2 X) = b\Omega U + bR = bU + c_n(\tau_n^2 X)$ (where we also use that
  $b U = b\Omega U$). 
\end{proof}

\medskip
\begin{corollary} If  $X$ is any object in $\Cal S(n)$, then
$$
   bX = c_nX + c_n(\tau_nX) + c_n(\tau_n^2 X).
$$
\end{corollary}

\begin{proof}
  We can assume that $X$ is indecomposable. If $X$ is projective, $bX = 1 = c_nX$
  yield the assertion.
  
  Thus, we can assume that $X$ is reduced. Proposition~\ref{prop-three-twelve-two}
  asserts that $bU = c_nX + c_n(\tau_nX).$
  Proposition~\ref{prop-three-twelve-three} asserts that $bX - bU = c_n(\tau_n^2X).$
  The sum of the two equalities is the required one.
\end{proof}

\vfill\eject
\centerline{\Gross First part: Sparsity.}

\addcontentsline{toc}{part}{First part: Sparsity}

\section{Telescope filtrations.}
\label{sec-four}

In this section we show Theorem~\ref{theoremone}.  For this, we introduce
telescope filtrations and show that for every object in $\mathcal S(n)$
a telescope filtration exists.  The key Proposition~\ref{prop-four-one} shows that
in each telescope filtration, factors which are $0$-pickets split off as direct
summands.  As a consequence, we obtain Theorem~\ref{theoremone}.

\subsection{Admissible submodules and telescope filtrations.}
\label{sec-four-one}

We always will interpret $\Cal S(n)$ as a full subcategory of $\mod T_2(\Lambda)$. In particular,
when we need homological invariants such as $\Ext$-groups, the reference
category is $\mod T_2(\Lambda)$. Thus, if $X,Z$ are objects in $\Cal S(n)$,
the group $\Ext^1(Z,X)$ may be considered
as the set of equivalence classes of short exact sequences
$0\to X \to Y \to Z \to 0$ in $\mod T_2(\Lambda)$. 

\medskip
The subcategory $\Cal S(n)$ of $\mod T_2(\Lambda)$ is closed under submodules.
If $(U,V)$ is an object of $\Cal S(n)$ and $(U',V')$ is a $T_2(\Lambda)$-submodule of $(U,V)$,
the factor module $(U,V)/(U',V')$ usually will not belong to $\Cal S(n)$ (since
the induced map $U/U' \to V/V'$ does not have to be a monomorphism). We say
that $(U',V')$ is an {\it admissible} submodule of $(U,V)$, provided
$U' = U\cap V'$, thus provided 
$(U,V)/(U',V')$ belongs to $\Cal S(n).$ 

If $(U,V)$ is an object of $\Cal S(n)$ and $(U',V')$ is an admissible 
submodule of $(U,V)$,
we will say that
$(U',V')$ is a {\it g-split} subobject of $(U,V)$, provided the embedding $V' \to V$
is a split monomorphism (g-split means: ``globally split''). 
The equivalence classes of g-split exact sequences form a
subgroup of $\Ext^1(Z,X)$ which we denote by $\Ext^1_{\text g}(Z,X).$

\medskip
Let $0\to X_1 \to X \to X_2 \to 0$ be an exact sequence in $\Cal S(n)$, with
$X_i = (U_i,V_i)$ for $i=1,2$ and $X = (U,V).$ Then, obviously,
$$
 |U| = |U_1|+|U_2|,  \quad\text{and} \quad  
  |V/U| = |V_1/U_1|+|V_2/U_2|, 
$$
If the sequence is, in addition, g-split, then also
$$
 bV = bV_1+bV_2,
$$
so that 
$$
 pX = \frac{|U_1|+|U_2|}{bV_1+bV_2}, \quad\text{and}\quad 
 rX = \frac{|V_1/U_1|+|V_2/U_2|}{bV_1+bV_2}.
$$

\medskip
If $X$ has a g-split filtration with factors $F_i$, $1\le i \le t$, then
\begin{align*}
 &\min \{pF_i\mid 1\le i \le t\} \le pX \le
 \max \{pF_i\mid 1\le i \le t\}, \cr
 &\min \{rF_i\mid 1\le i \le t\} \le rX \le
 \max \{rF_i\mid 1\le i \le t\},
\end{align*}
with proper inequalities in case $\min < \max.$ 

\medskip 
\begin{definition}
  A {\it telescope} filtration of $X \in \Cal S(n)$ is a filtration 
  $$
  0 = X_0 \subset X_1 \subset \cdots \subset X_{b-1} \subset X_b = X
  $$
  consisting of g-split admissible 
  subobjects, such that each factor $X_i/X_{i-1}$ is a picket,
  say equal to $([t_i],[m_i])$, 
  and such that, in addition, we have $m_1 \ge m_2 \ge \cdots \ge m_b.$
  (Since we deal with g-split subobjects, we see that the global space of $X$ is
  equal to the direct sum $\bigoplus_{i=1}^b [m_i],$ as a consequence, $bX = b$,
  and the subspace has a filtration with factors $[t_i],$ thus its length is equal to
  $\sum_{i=1}^b t_i$.)  
\end{definition}

\medskip
\begin{lemma}
  \label{lem-four-one}
  Any object $X$ in $\Cal S(n)$ has a telescope filtration.
\end{lemma}

\begin{proof}
  Let $X = (U,V)$. Let $b = bV$. Let $V = \bigoplus_{i=1}^b [\lambda_i],$ where
  $(\lambda_1,\dots,\lambda_b)$ is a partition. For $0\le t\le b$, let $V_t 
  = \bigoplus_{i=1}^t [\lambda_i],$ and $X_t = (V_t\cap U,V_t)$. Then $(X_t)_t$ is
  a telescope filtration of $X$.
\end{proof}

\subsection{Splitting off $0$-pickets.}
\label{sec-four-two}

A picket of the form $(0,[m])$ will be called a {\it $0$-picket.} 

\medskip
\begin{lemma}
  \label{lem-four-two}
  Let $0 \to X \to Y \to Z \to 0$ be an exact sequence in $\Cal S$.
  Let $F$ be a $0$-picket.
  \begin{itemize}[leftmargin=3em]
  \item[\rm(a)] If $X = F$ and the height of $X$ is greater or equal to the height of $Y$, 
    then the
    sequence splits (in particular, $Y$ has a direct summand isomorphic to $F$).
  \item[\rm(b)] If $Z = F$ and the sequence is g-split, then the sequence splits
    (in particular, $Y$ has a direct summand isomorphic to $F$).
  \end{itemize}
\end{lemma}

\begin{proof} (a) Let $h$ denote the height of $X$.
By assumption, the exact sequence lies in $\Cal S(h)$, and $\Cal S(h)$ is 
the category of Gorenstein-projective
$T_2(\Lambda)$-modules, where $\Lambda = k[T]/T^h.$ Since $(0,[h])$ is a projective 
$T_2(\Lambda)$-module, we have $\Ext^1(Z,(0,[h]))=0$ for all objects in $\Cal S(h).$

	\smallskip
(b) Really, nothing has to be shown: 
Let $0 \to X \to Y \to Z \to 0$ be a g-split exact sequence with
$X = (U,V)$ and $Z = (0,[h]).$
Then the global space of $Y$ is $V\oplus [h]$, its subspace is 
$U\oplus 0 = U$ with its given embedding into $V$, thus $Y = X\oplus(0,[h])$.
\end{proof}

\medskip
\begin{warning}
  In order to understand the setting of Lemma~\ref{lem-four-two}
  let us exhibit several examples.
	\smallskip
        
In (a), we need the height condition. In order to see this, let $X = (0,[2])$ 
be the $0$-picket of height
$h = 2$. We provide two non-split exact
sequences $0 \to X \to Y \to Z\to 0\ $.
First of all, let $Y = (0,[3])$ and $Z = Y/X = (0,[1]).$
Here, $Z$ has height at most $h$ (but note that $Y$ is again a $0$-picket). 
Second, let $Y = ([1],[2])\oplus (0,[3])$
with global space generated by $x,y$ such that $T^2x = T^3y = 0$ and with subspace
generated by $Tx$. Let $X$ be the admissible subobject generated by $x+Ty$. Then 
$X$ is isomorphic to $([0],[2])$ and $Z = Y/X$ is isomorphic to $([1],[3])$. Note that 
$Z$ has height 3 and that the sequence is g-split. 

	\smallskip
In (b), we need the g-splitting, as the case $X = ([1],[1]),\
Y = ([1],[2]),\ Z = (0,[1])$ shows. In this example, 
$bY < bX+bZ.$ An example with $bY = bX+bZ$
is given by $Y = E_2^2$, where the global space of $E_2^2$ is
generated by $x,y$ with $T^3x = 0 = Ty$ and the subspace in $E_2^2$ is
generated by $Tx+y.$  Let $X$ be the admissible subobject of $E_2^2$
generated by $Tx+y.$ 
Then $X$ is isomorphic to $([2],[2])$ and $Z = Y/X$ is isomorphic to $(0,[2]).$
	\smallskip
        
Conversely, there are exact sequences $0 \to X \to Y \to Z \to 0$ such that 
$Y$ has a $0$-picket $F$ as a direct summand, whereas neither $X$, nor $Z$,
have $F$ as a direct summand. For example, 
let $Y = ([0,[2])\oplus ([2],[2])$, say with 
global space generated by $x$ and $y$, and subspace generated by $y$. Let $X$ be the admissible
subobject of $Y$ generated by $Tx+y$. Then both $X$ and $Y/X$ are isomorphic to $([1],[2])$.
\end{warning}

\medskip
\begin{proposition}
  \label{prop-four-one}
  Let $F$ be a $0$-picket. 
  If $X$ has a telescope filtration where at least one of the factors is isomorphic to $F$, 
  then $X$ has a direct summand isomorphic to $F$.
\end{proposition}

\begin{proof}
  Let $(X_i)_i$ be a telescope filtration of $X$ and assume that $F = X_t/X_{t-1}$ is 
  a $0$-picket of height $h$.
  Then $X_{t-1} \subset X_t$ is a g-split embedding, and $X_t/X_{t-1}$,
  and hence $X/X_{t-1}$, have height at most $h$. 

  We apply part (a) of Lemma~\ref{lem-four-two} to the embedding of
  $X_t/X_{t-1}$ into $X/X_{t-1}$, so there exists a submodule $Z$ of $X$
  containing $X_{t-1}$ such 
  that $X/X_{t-1} = X_t/X_{t-1}\oplus Z/X_{t-1}$
  (thus $X_t + Z = X,$ and $X_t\cap Z=  X_{t-1}$).
  We apply part (b) of Lemma~\ref{lem-four-two} to the embedding
  of $X_{t-1}$ into $X_t$ and see that 
  $X_t = C\oplus X_{t-1}$ (thus $C+X_{t-1} = X_t$ and $C\cap X_{t-1} = 0$).  
  $$
  
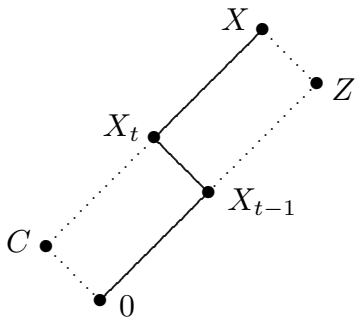
{\beginpicture
    \setcoordinatesystem units <.6cm,.6cm>
    \multiput{} at 0 0  5 4 /
    \multiput{$\bullet$} at 0 1  1 0  2 3  3 2  4 5  5 4 /
    \plot 1 0  3 2  2 3  4 5 /
    \setdots <1mm>
    \plot 1 0  0 1  2 3 /
    \plot 3 2  5 4  4 5 /
    \put{$0$} at 1.5 -0.1
    \put{$C$} at -.5 1.1
    \put{$X_{t-1}$} at 4 1.8
    \put{$X_{t}$} at 1.4 3.2
    \put{$X$} at 3.5 5.2
    \put{$Z$} at 5.5 3.9 
    \endpicture}
  $$
  We have $C+Z = C+X_{t-1}+Z = X_t+Z = X$ and $C\cap Z = C\cap X_t\cap Z = C\cap X_{t-1} = 0.$
  Therefore $X = C\oplus Z.$ Of course, $C$ is isomorphic to $X_t/X_{t-1} = F.$
\end{proof}

\subsection{Proof of Theorem~\ref{theoremone}.}
\label{sec-four-three}

\begin{lemma}
  \label{lem-four-three}
  Let $X$ be a non-zero object in $\Cal S(n)$.
  If $X$ has a telescope filtration, such that no factor is a $0$-picket,
  then $pX \ge 1$.
\end{lemma}

\begin{proof}
  Let $X$ be a non-zero object in $\Cal S(n)$ and let
  $$
  0 = X_0 \subset X_1 \subset \cdots \subset X_{b-1} \subset X_b
  $$
be a telescope filtration of $X$ with factors $F_i = X_i/X_{i-1}.$
Assume 
that none of the factors $F_i$, $1\le i\le b$ is a $0$-picket.
Let $F_i = [h_i,m_i]$, where $1 \le h_i \le m_i.$
The subspace $U$ in $X$ has length 
$|U| = \sum_{i=1}^b h_i \ge b$,
and we have $bX = b.$ Thus $|U| \ge bX,$ therefore 
$pX = |U|/bX \ge 1.$ 
\end{proof}

\medskip
\begin{warning-four-three-cont}
  Note that an object $X$ with a telescope filtration such 
  that no factor is a $0$-picket, may have a direct summand which is a $0$-picket.

\smallskip 
{\bf Example.}
  The object $X = ([2],[2]) \oplus (0,[1])$ has an admissible  g-split subobject 
of the form $X'$ isomorphic to $([1],[2])$ and the corresponding factor object is isomorphic to 
$([1],[1]).$ Namely, let $x,y$ be generators of the global space of $X$, say
with $T^2x = 0 = Ty$ and such that the subspace in $X$ is generated by $x$.
Let $X'$ be the admissible subobject
of $X$ whose total space is generated by $x+y.$ Then $X'$ is isomorphic to $([1],[2])$
and $X/X'$ is isomorphic to $([1],[1]).$ $\s$
\end{warning-four-three-cont}

	\smallskip
        
In particular, we see: {\it The class of objects without direct summands which are $0$-pickets
is not closed under extensions.}

\smallskip
\begin{proof}[Proof of Theorem~\ref{theoremone}]
  We show in addition the statements Theorem~1.1$'$ and Theorem~1.1$''$
  in Section~\ref{theoremone-dual} and Theorem~1.1 (reformulated) in Section~\ref{theoremone-again}.
  Assume that $X$ is a non-zero object in $\Cal S(n)$.

\medskip
(a) {\it If $X$ has no direct summand which is a $0$-picket, then $pX \ge 1.$}
	\smallskip
        
Proof. Let $X$ be a non-zero object and assume that $X$ has no direct summand which is a
$0$-picket. Take any telescope filtration of $X$, say with factors $F_i$. According to
Proposition~\ref{prop-four-one}, no $F_i$ is a $0$-picket. Lemma~\ref{lem-four-three} yields
that $pX \ge 1.$ $\s$

\medskip
(b) {\it Let $X$ be a non-zero object in $\Cal S(n).$
If $X$ has no direct summand which is a picket of height $n$, 
then $pX+rX \le n-1.$}
	\smallskip
        
Proof. We assume that $X$ is non-zero and has no direct summand which is a picket
of height $n$. In particular, $X$ is reduced. Let $Y = \tau^4X.$ Then $Y$ is reduced and
$\tau^2Y = \tau^6 X = X.$ If $Y$ would have a direct summand which is a $0$-picket, then
$X = \tau^2Y$ has a direct summand which is a picket of height $n$
a contradiction (here we use that $\tau^2(0,[m]) = ([m],[n])$).
Thus, according to (a), we have $pY \ge 1.$ According to Theorem~\ref{theoremthree},
we have 
$pX+rX = p(\tau^2Y) + r(\tau^2Y) = rY + n-pY-rY = n-pY$. Since $pY \ge 1,$ it follows
that $pX+rX = n-pY \le n-1.$ $\s$

\medskip
(c) {\it 
If $X$ has no direct summand which is a picket of the form $([m],[m])$, then $rX \ge 1.$}
	\smallskip
        
Proof. We could use $\tau^2$ and (b), or $\tau^4$ and (a). But actually, it
is easier to use duality. We assume that $X$ has no direct summand which is a picket
of the form $([m],[m])$. Then $\D X$ has no direct summand which is a picket of
the form $(0,[m])$, thus a $0$-picket. According to (a) we have $p(\D X) \ge 1.$
Since $rX = p(\D X),$ we see that $rX \ge 1.$ 
\end{proof}

\bigskip

\section{Special filtrations and nice filtrations.}
\label{sec-five}

We want to prove Theorem~\ref{theoremtwo}; it follows from Corollary~\ref{cor-five-one}
for the proof of which we introduce special filtrations and nice filtrations.
We conclude this section with the slightly stronger statement Theorem~\ref{theoremnine}.

\subsection{Extended pickets and the proof of Theorem~\ref{theoremtwo}.}
\label{sec-five-one}

Theorem~\ref{theoremtwo} concerns the corners of the triangle $\Delta_1$.
Let us first look at the corner $(1,1).$ We claim:

\medskip
\begin{theorem}
  \label{thm-five-one}
  Let $X$ be an indecomposable object in $\Cal S(n)$ with
  pr-vector $(1,1)$. Then $X$ is the picket $([1],[2])$ or the unique bipicket in $\Cal S(3)$.
\end{theorem}

\medskip 
Note that there is just one bipicket in $\Cal S(3)$, namely $E_2^2 = (U,V)$ with
$V = [3,1]$ and $U$ and $V/U$ both isomorphic to $[2].$
The essential information of Theorem~\ref{thm-five-one} is the following: {\it 
Let $X$ be indecomposable in $\Cal S(n)$. If $qX \le 2$, then $X$ belongs to $\Cal S(3)$.}
Namely, the indecomposable objects $X$ of height at most 3 are well-known and $qX$ 
is easy to calculate for these objects:  There are just six indecomposable 
objects in $\Cal S(3)$ with $qX \le 2,$ namely the five pickets $(0,[1]),\ ([1],[1]),\
(0,[2]),\ ([1],[2]),\ ([2],[2]),$ and in addition $E_2^2$.

\bigskip
A bipicket $X = (U,V)$ is said to be an {\it extended picket} provided $V$ has a part equal to $[1].$
The extended pickets are the objects of the form $E_u^w = ([u],[u+w-1,1],[w])$ with $u,w\ge 2$.
Note that if  $\bpar X = ([u],[u+w-1,1],[w]),$ then $X$ has to be
indecomposable and we have for $n = u+w-1$ that 
$\tau_nX = (\soc [w],[w])$, thus $X$ is uniquely determined by $\bpar X.$
Here are the extended pickets
$E_u^w$ with $u+w \le 6$:
$$
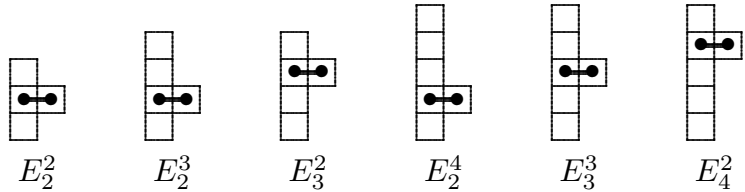
{\beginpicture
    \setcoordinatesystem units <.3cm,.3cm>
\put{\beginpicture
\multiput{} at 0 0  2 5 /
\plot 0 0  1 0  1 3  0 3  0 0 /
\plot 0 1  2 1  2 2  0 2 /
\multiput{$\bullet$} at 0.5 1.5  1.5 1.5 /
\plot .5 1.5   1.5 1.5 /
\plot .5 1.55  1.5 1.55 /
\plot .5 1.45  1.5 1.45 /

\endpicture} at  0 0 
\put{\beginpicture
\multiput{} at 0 0  2 5 /
\plot 0 0  1 0  1 4  0 4  0 0 /
\plot 0 1  2 1  2 2  0 2 /
\plot 0 3  1 3 /
\multiput{$\bullet$} at 0.5 1.5  1.5 1.5 /
\plot .5 1.5   1.5 1.5 /
\plot .5 1.55  1.5 1.55 /
\plot .5 1.45  1.5 1.45 /

\endpicture} at  5 0 
\put{\beginpicture
\multiput{} at 0 0  2 5 /
\plot 0 0  1 0  1 4  0 4  0 0 /
\plot 0 2  2 2  2 3  0 3 /
\plot 0 1  1 1 /
\multiput{$\bullet$} at 0.5 2.5  1.5 2.5 /
\plot .5 2.5   1.5 2.5 /
\plot .5 2.55  1.5 2.55 /
\plot .5 2.45  1.5 2.45 /

\endpicture} at  10 0 
\put{\beginpicture
\multiput{} at 0 0  2 5 /
\plot 0 0  1 0  1 5  0 5  0 0 /
\plot 0 1  2 1  2 2  0 2 /
\plot 0 3  1 3 /
\plot 0 4  1 4 /
\multiput{$\bullet$} at 0.5 1.5  1.5 1.5 /
\plot .5 1.5   1.5 1.5 /
\plot .5 1.55  1.5 1.55 /
\plot .5 1.45  1.5 1.45 /

\endpicture} at  15 0 
\put{\beginpicture
\multiput{} at 0 0  2 5 /
\plot 0 0  1 0  1 5  0 5  0 0 /
\plot 0 2  2 2  2 3  0 3 /
\plot 0 1  1 1 /
\plot 0 4  1 4 /
\multiput{$\bullet$} at 0.5 2.5  1.5 2.5 /
\plot .5 2.5   1.5 2.5 /
\plot .5 2.55  1.5 2.55 /
\plot .5 2.45  1.5 2.45 /

\endpicture} at  20 0 
\put{\beginpicture
\multiput{} at 0 0  2 5 /
\plot 0 0  1 0  1 5  0 5  0 0 /
\plot 0 3  2 3  2 4  0 4 /
\plot 0 1  1 1 /
\plot 0 2  1 2 /
\multiput{$\bullet$} at 0.5 3.5  1.5 3.5 /
\plot .5 3.5   1.5 3.5 /
\plot .5 3.55  1.5 3.55 /
\plot .5 3.45  1.5 3.45 /

\endpicture} at  25 0 
\put{$E_2^2$} at 0 -3.7
\put{$E_2^3$} at 5 -3.7
\put{$E_3^2$} at 10 -3.7
\put{$E_2^4$} at 15 -3.7
\put{$E_3^3$} at 20 -3.7
\put{$E_4^2$} at 25 -3.7
\endpicture}
$$
Of course, $E_2^2$ is just the unique indecomposable object of $\Cal S(3)$ which
is not a picket. 

\medskip 
\begin{proposition}
  \label{prop-five-one}
  Any object $X$ is a direct sum $X = X'\oplus X''$, where $X'$ has a g-split 
  filtration whose factors are pickets of height at least $2$ and extended pickets,
  and where $X''$ has height at most $1.$
\end{proposition}

\medskip 
In particular, if $X$ is indecomposable and not of height $1$, then $X$ 
has a g-split filtration whose factors are pickets of height at least $2$ and extended pickets.

\medskip
\begin{corollary}
  \label{cor-five-one}
  If $X$ is indecomposable with $qX \le 2$, then $X$ belongs to $\Cal S(3)$.
\end{corollary}

\begin{proof}
  For pickets, the mean is just the height.
  For the extended picket $Y = ([u],[u+w-1,1],[w])$, we have $qY = (u+w)/2 \ge 2$.
  Thus, if $Y$ is a picket or an extended picket with $qY \le 2$, then $Y$ is a picket of 
  height at most $2$ or else $Y = E_2^2$; always, $Y$ belongs to $\Cal S(3)$.

  Assume now that $X$ is indecomposable with $qX \le 2$.
  According to Proposition~\ref{prop-five-one}, 
  $X$ has either height 1 or has a g-split filtration with factors in $\Cal S(3)$.
  Therefore $X$ belongs to $\Cal S(3)$.
\end{proof}

\subsection{Special filtrations and nice filtrations.}
\label{sec-five-two}

Proposition~\ref{prop-five-one}
will follow directly  from Lemma~\ref{lem-five-one} which provides
a slightly stronger statement, namely further details 
on the direct decomposition $X = X'\oplus X''.$ 
We need some definitions.

\medskip 
A subobject $X'$ of $X$ is said to be a  
{\it special subobject} provided the following two conditions are satisfied: First, 
$X'$ is a picket of height equal to 
the height of $X$ (as a consequence, $X'$ will be a g-split subobject),
and second, the dimension vector of $X'$ is maximal (that means: Let $h$ be the height
of $X$, so that $\bdim X' = (t,h)$ for some $0 \le t \le h$; if now $(i,h)$ is the dimension vector
of any picket which is a subobject of $X$, then $i \le t$). Note that a special subobject $X'$ is
always an admissible subobject.  
	
A {\it special filtration} of $X$ is defined to be 
a sequence $0 = X_0 \subseteq X_1 \subseteq \cdots
\subseteq X_b = X$ of (admissible g-split) subobjects such that $X_i/X_{i-1}$ is a special subobject
of $X/X_{i-1}$, for $1\le i \le b.$
	
Of course, any object has a special filtration. Any special filtration is a telescope filtration,
but the converse is 
not true: Take $X = ([2],[2])\oplus(0,[1])$. Thus, the global space of $X$ is
$V = [2,1]$, say with generators $x,y$, with $T^2x = 0$ and $Ty = 0$. Let $V_1$
be the submodule of $V$ generated by $x+y$, then $U_1 = U \cap V_1$ is generated by $Tx$. Thus
we have the telescope filtration $0 = X_0 \subset X_1 \subset X_2 = X$ with $X_1 = (U_1,V_1)$
isomorphic to $([1],[2])$ and $X/X_1$ isomorphic to $([1],[1]).$ This filtration is not special,
since $X_1$ is not a special subobject of $X$. 

Finally, a {\it nice filtration} is by definition a g-split
filtration whose factors are pickets of height at least $2$ and extended pickets.

\medskip
\begin{remark}
  The set of factors of a nice filtration is not uniquely determined.
\end{remark}

\medskip
The principal component $\Cal P(6)$ contains a (unique) indecomposable object $X$ 
with $\bpar X = ([4,2],[6,4,1,1],[4,2])$ (note that there is a sectional path
from $S = (0,[1])$ to $X$ of length $5$). 

There are two different nice filtrations of $X$, given by g-split exact sequences:
$$
 0 \to E_3^4 \to X \to E_3^2 \to 0, \quad \text{and}\quad
 0 \to E_2^3 \to X \to E_4^3 \to 0. 
$$
Here are pictures:
$$
{\beginpicture
    \setcoordinatesystem units <.4cm,.4cm>
\put{\beginpicture
\multiput{} at -3 1   2 6 /
\plot 0 2  1 2  1 6   0 6  0 2 /
\plot 0 4  2 4  2 5  0 5 /
\plot 0 3  1 3 /
\multiput{$\bullet$} at 0.5 4.5  1.5 4.5   /
\plot 0.5 4.5  1.5 4.5 /
\plot 0.5 4.55  1.5 4.55 /
\plot 0.5 4.45  1.5 4.45 /

\plot -3 3  -1 3  -1 4  -2 4  -3 4 /
\plot -3 1  -2 1  -2 7  -3 7  -3 1  /
\plot -3 2  -2 2 /
\plot -3 5  -2 5 /
\plot -3 6  -2 6 /
\multiput{$\bullet$} at -1.5 3.5  -2.5 3.5 /
\plot -1.5 3.5  -2.5 3.5 /
\plot -1.5 3.55 -2.5 3.55 /
\plot -1.5 3.45 -2.5 3.45 /
\put{$\circ$} at -2.5 4.5 
\setdashes <.9mm>
\plot -2.2 4.5   0.8 4.5 /
\plot -2.2 4.55  0.8 4.55 /
\plot -2.2 4.45  0.8 4.45 /

\setdots <1mm>
\plot -2 6  0 6 /
\plot -1 4  0 4 /
\plot -1 3  0 3 /
\plot -2 2  0 2 /
\endpicture} at 0 0 
\put{\beginpicture
\multiput{} at -3 1   2 6 /
\plot 0 2  0 1  1 1  1 2 /
\plot 0 2  1 2  1 6   0 6  0 2 /
\plot 0 4  2 4  2 5  0 5 /
\plot 0 3  1 3 /
\plot 0 6  0 7  1 7  1 6 /
\multiput{$\bullet$} at 0.5 4.5  1.5 4.5   /
\plot 0.5 4.5  1.5 4.5 /
\plot 0.5 4.55  1.5 4.55 /
\plot 0.5 4.45  1.5 4.45 /

\plot -3 3  -1 3  -1 4  -2 4  -3 4 /
\plot -3 2  -2 2  -2 6  -3 6  -3 2  /
\plot -3 2  -2 2 /
\plot -3 5  -2 5 /
\plot -3 6  -2 6 /
\multiput{$\bullet$} at -1.5 3.5  -2.5 3.5 /
\plot -1.5 3.5  -2.5 3.5 /
\plot -1.5 3.55 -2.5 3.55 /
\plot -1.5 3.45 -2.5 3.45 /
\put{$\circ$} at -2.5 4.5 
\setdashes <.9mm>
\plot -2.2 4.5   0.8 4.5 /
\plot -2.2 4.55  0.8 4.55 /
\plot -2.2 4.45  0.8 4.45 /

\setdots <1mm>
\plot -2 6  0 6 /
\plot -1 4  0 4 /
\plot -1 3  0 3 /
\plot -2 2  0 2 /
\endpicture} at 13 0 
\endpicture}
$$
This shows: {\it The set of factors of a nice filtration
  is not uniquely determined.}

\subsection{Splitting off summands of height one.}
\label{sec-five-three}

We are going to show:

\medskip
\begin{lemma}
  \label{lem-five-one}
  Let $(X_i)_i$ with $0 \le i \le b$ be a special filtration of $X$ with factors $F_i$. 
  Choose $0 \le t\le b$ minimal such that all factors $F_i$ with $i > t$ have height $1$.
  Then there is a direct decomposition $X = X'\oplus X''$, with $X_t \subseteq X'$ such that
  $X'$ has a nice filtration (and $X''$ has height at most $1$).
\end{lemma}

\medskip
Note that $X''$ has height at most 1, since 
$X'' = X/X'$ is a factor object of $X/X_t$.

\smallskip 
\begin{proof}[Proof of the lemma]
  
{\bf (1)} First, we consider {\bf the case $t = 1$,} thus $X = (U,V)$ with $[V] = [h,1,\dots,1]$. 
We claim that
$X = Y\oplus Z$, where $Y$ is a picket or an extended picket and $Z$ has height at
most $1$; we call such a direct decomposition a {\it good} decomposition. 
To produce such a decomposition is slightly tedious, but not difficult, and is left to the reader;
(note that in the cases $2\le h \le 5$, this decomposition is well-known and there is no further
complication in the general case).

It remains to determine the special subobjects $X_1$ of $X$ and to show that for
any $X_1$, there is a good decomposition $X = Y'\oplus Z$ with $X_1 \subseteq Y'.$
We distinguish three cases. 

First, let $Y$ be a picket isomorphic to $([h],[h])$.
Let $y$ be a generator of the global space of $Y$. 
Then the global space of an arbitrary special subobject of $X$ is generated by $y+z$ for some
$z\in U\cap Z,$ thus $X_1 = (\Lambda(y+z),\Lambda(y+z))$
and $X = X_1\oplus Z$ is a good
decomposition.

Second, let $Y$ be a picket of type $([s],[h])$ with $s < h.$
Again, we assume that the global space of $Y$ is generated by $y$. 
Then the global space of any special subobject $X_1$ of $X$ is generated by $y+z$, now with
an arbitrary element $z\in Z.$ 
Then  $X_1 = (U\cap\Lambda(y+z),\Lambda(y+z))$ and again, $X = X_1\oplus Z$
is a good decomposition. 

Third, let $Y$ be an extended picket, say generated by $y$ and $w$, where $|\Lambda w| = 1.$ 
Then the global space of an arbitrary special subobject $X_1$ of $X$ is generated by $y+w'+z$ for some
$w'\in \Lambda w$ and $z\in Z.$ Then $X = Y'\oplus Z$, where the global space of $Y'$ is generated by 
$y+w'+z$ and $w$, thus by $y+z$ and $w$. Again, this is a good decomposition and $X_1 \subseteq Y'.$

\medskip
{\bf (2) The induction.} We start with an arbitrary object $X$ with a special filtration $(X_i)_i$ with 
factors $F_i$ such that the factors $F_1,\dots, F_t$  have height at least $2$, whereas the
remaining factors $F_i$ have height $1$. 

If $t = 0$, then $X$ has height at most $1$; thus, take $X' = 0$ and $X'' = X.$ 
The case $t = 1$ has been considered already. Thus, we assume now that $t\ge 2.$

We look at $Y = X/X_{t-1}.$ By assumption, there is given a special filtration
of $Y$, namely by the subobjects $Y_i = X_{t+i-1}/X_{t-1}$, and just one factor has height at least $2$,
namely $Y_1/Y_0 = X_t/X_{t-1}.$ According to (1), 
we obtain a direct decomposition $Y = Y'\oplus Y''$, where
$Y'$ is a picket of height at least 2 or an extended picket, and where $Y''$ has height at most $1$.
Let $Y' = W/X_{t-1}$ and $Y'' = Z/X_{t-1}$, where $X_{t-1} \subseteq W,Z\subseteq X$.
Then $X_{t-1} = W\cap Z$ and $W+Z = X.$ 

Now we consider $Z$. Since $X_{t-1} \subseteq Z \subseteq X$, we can construct a special 
filtration $(Z_i)_i$ of $Z$
starting with the subobjects $Z_i = X_i$ with $0 \le i \le t-1$ (note that for $1\le i \le t-1$,
$Z_i/Z_{i-1}$ is a special subobject of $Z/Z_{i-1}$, since $Z/Z_{i-1}$ is a subobject of $X/X_{i-1}$).
We complete the sequence $0 = Z_0 \subseteq \cdots \subseteq Z_{t-1}$ to a special filtration of $Z$
by taking suitable subobjects $Z_i$ with $i \ge t.$ Note that the factors $Z_i/Z_{i-1}$ with $i \ge t$
have height $1$, whereas $Z_{t-1}/Z_{t-2} = X_{t-1}/X_{t-2}$ has height at least $2.$

Now we apply induction: $Z$ has the special filtration $(Z_i)_i$ 
with precisely $t-1$ factors of height at least $2$,
therefore we have a direct decomposition $Z = Z'\oplus Z''$, with $Z_{t-1} (= X_{t-1}) \subseteq Z'$,
such that $Z'$ has a nice
filtration, whereas
$Z''$ has height at most $1$.
	\medskip

Let $X' = Z'+W.$ 
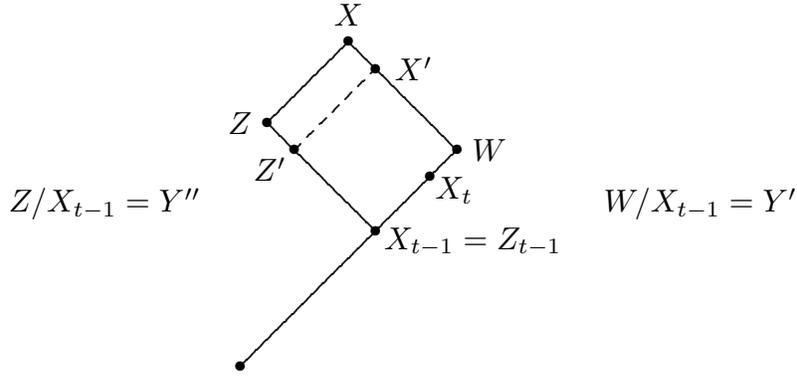
$$
{\beginpicture
    \setcoordinatesystem units <.3cm,.3cm>
\multiput{} at 0 0  8 12 / 
\plot 0 0  8 8  4 12  1 9  5 5 / 
\multiput{$\ssize\bullet$} at 0 0  5 5  8 8  4 12  1 9  7 7  2 8  5 11 /
\setdashes <1mm>
\plot 2 8  5 11 /
\put{$X$} at  4 13
\put{$X'$} at 6.4  11 
\put{$W$} at 9.2  8
\put{$X_{t}$} at 7.9 6.5 
\put{$X_{t-1} = Z_{t-1}$} at 8.6 4.6 

\put{$Z$} at 0 9
\put{$Z'$} at 1.1 7.3

\put{$Z/X_{t-1} = Y''$} at -5 6
\put{$W/X_{t-1} = Y'$} at 17 6
\endpicture}
$$

Then $X'$ is a g-split extension of $Z'$ by $X'/Z' \simeq W/X_{t-1} = Y'$.
Now $Z'$ has a nice filtration, and $Y'$ is 
a picket of height at least $2$ or an extended picket. 
Thus $X'$ has a nice filtration.

On the other hand, $Z = Z'\oplus Z''$ implies that $X = X'\oplus Z''$. 
The decomposition $X = X'\oplus Z''$ is a decomposition as required:
$X'$ has a nice filtration, and $Z''$ has height at most $1.$
\end{proof}

\subsection{A slightly stronger statement.}
\label{sec-five-four}

The proof of Lemma~\ref{lem-five-one} yields also the following assertion:

\medskip
\begin{theorem}
  \label{theoremnine}
  Let $X$ be an object without any direct summand of height $1$ and without any direct
  summand which is a $0$-picket. Then $X$ has a g-split
  filtration with factors $F_i,$ $1\le i \le t$, which are pickets or extended
  pickets of height $h_i$, but not $0$-pickets, 
  such that $h_1 \ge h_2 \ge \cdots \ge h_t \ge 2.$
\end{theorem}

\begin{proof}
  The proof of Lemma~\ref{lem-five-one} produces a nice filtration $(X'_i)_i$ such that
  $F'_i = X'_i/X'_{i-1}$ contains a special subobject of $X/X'_{i-1}$. The
  height $h_i$ of $F'_i$ is just the height of $X/X'_{i-1}$, and we have $h_i \ge h_{i+1}$
  for all $i$. 
  
  Assume that $F'_i$ is a $0$-picket. Then the height of $X/X'_{i-1}$ is
  equal to $h_i$. Since $X/X'_i$ is
  Gorenstein-projective in $\Cal S(h_i)$, we see that $F'_i$ is a direct summand of
  $X/X'_{i-1}.$ Since $X'_{i-1}$ is a g-split subobject of $X'_i$,
  the inclusion $X'_{i-1} \to X'_i$ splits (see Lemma~\ref{lem-four-two}),
  thus $F'_i$ is a direct summand
  of $X$. Since $F'_i$ is a $0$-picket, we obtain a contradiction.
\end{proof}

\vfill\eject
\centerline{\Gross Second part: Density.}
\addcontentsline{toc}{part}{Second part: Density}

\section{BTh-vectors in $\mathbb T(n)$.}
\label{sec-six}

We are going to show  Theorem \ref{theoremfour}.
Note that in this section, we always will consider graded objects, 
thus objects in $\widetilde{\Cal S}(n).$ Recall that a vector $(p,r) \in \mathbb T(n)$
is called a BTh-vector provided there is a positive natural number $a$
such that for any $t\in \mathbb N_1$, there is a $\mathbb P^1$-family of indecomposable objects with
uwb-vector $(tap,tar,ta).$  We recall

\medskip
\begin{theoremfour-recall}
  \label{theoremfour-recall}
  Any rational pr-vector $(p,r)$ in $\mathbb T(n)$ 
  with boundary distance at least $2$ is a BTh-vector.
\end{theoremfour-recall}

\medskip
Our constructions of BTh-families will rely on the standard $\mathbb P^1$-family 
$M_c$ in $\widetilde{\Cal S}(6)$ with $VM_c = [6,4,2]$ and $UM_c = [4,2].$ 
Namely, we will construct some larger objects which can functorially be reduced to members
of the standard family. In this way, we will obtain objects which are indecomposable and pairwise non-isomorphic.

\subsection{Some endofunctors of $\widetilde{\Cal S}(n)$.}
\label{sec-six-one}

First, the endofunctor $G_z$ of $\widetilde{\Cal S}(n)$ for $z\in \mathbb Z.$
For any $z\in \mathbb Z$, the endofunctor $G_z$ of $\mod \widetilde Q$ is defined
for $M = (M_x,M_{x'})_{x,x'}$ in $\mod \widetilde Q$  by deleting the
vector spaces with index $z$ and $z'$, relabeling the remaining vector spaces
as follows: For $x < z,$ we define  $(G_zM)_x = M_x$
and $(G_zM)_{x'} = M_{x'},$ whereas for $z \le x,$ we define  $(G_zM)_x = M_{x+1}$
and $(G_zM)_{x'} = M_{(x+1)'};$ finally, we use as map $(G_zM)_{z-1} \leftarrow (G_zM)_{z}$
the composition $M_{z-1} \leftarrow M_z \leftarrow M_{z+1},$ 
and similarly as map $(G_zM)_{(z-1)'} \leftarrow (G_zM)_{z'}$
the composition $M_{(z-1)'} \leftarrow M_{z'} \leftarrow M_{(z+1)'}.$ 

We stress the following obvious fact: {\it If $G_z(M) = 0,$ then $M_x = 0$ for all $x\neq z.$}
In particular, if $G_z(M)$ is indecomposable, then $M = M'\oplus M''$,
where $M'$ is indecomposable and $(M'')_x = 0$ for all $x\neq z.$

\medskip 
Second, the endofunctor $H_z$ of $\widetilde{\Cal S}(n)$ for $z\in \mathbb Z.$
If $M$ is a representation of $\widetilde Q,$ and $z\in \mathbb Z$, let $H_zM = M/M'$, where
$M' = M_z\cap \soc M.$ If all the maps $M_{x'} \to M_x$ are inclusion maps (so that 
$M$ is Gorenstein-projective), then $H_zM$ is Gorenstein-projective if and only if $M_{z'}\cap \soc M = 0.$

\subsection{Solid objects and strongly solid objects.}
\label{sec-six-two}

We say that an object $X = (U,V)$ in $\widetilde{\Cal S}$ 
is {\it $st$-solid} where $s\le t$ are integers, provided first, 
$V$ is a direct sum of indecomposables of the form $[i,j]$ with $i\le s \le t \le j$, 
second $U_i = V_i$ for $i< s$, and third, $U_j = 0$ for $t<j.$
Typical examples of $34$-solid objects are the objects in the standard $\mathbb P^1$-family of
$\widetilde{\Cal S}(6).$ 

An $st$-solid object $X = (U,V)$ 
is said to be {\it strongly $st$-solid,} provided $X$ is $st$-solid
and the indecomposable direct summands of $V$ are even of the form
$[i,j]$ with $i < s \le t < j$. 

If $X = (U,V)$ is strongly $st$-solid, then $\soc V \subseteq U \subseteq TV$, and we may form
$X^-_- = (U/\soc V,TV/\soc V)$; the object $X^-_-$ 
is $st$-solid, and we have $bX^-_- = bX.$ 

\medskip
\begin{lemma}
  \label{lem-six-two}
  The construction 
  $X \mapsto X^-_-$ provides an equivalence between the full subcategory of strongly
  $st$-solid objects and the full subcategory of $st$-solid objects.
\end{lemma}

\medskip 
The inverse functor will be denoted by $Z \mapsto Z^+_+.$

\smallskip
\begin{proof}
  If $X = (U,V)$ and $X' = (U',V')$ are $st$-solid, and $f\:X \to X'$ is a
  homomorphism, then by definition $f$ can be considered as a homomorphism $V \to V'$
  which maps $U$ into $U'$. Since the indecomposable direct summands of $V$ and $V'$
  are of the form $[i,j]$ with $i\le s \le j$, we have $f = 0$ if and only if $f_s = 0.$
  
  Now, let $X = (U,V)$ and $X' = (U',V')$ 
  be strongly $st$-solid, and $f\:X \to X'$ a homomorphism. Then $f$ maps $\soc V$ into $\soc V'$
  and $TV$ into $TV'$, thus it yields a map $f^-_-\:X^-_- \to (X')^-_-.$ Assume that $f^-_- = 0.$ 
  Then $f_s = (f^-_-)_s = 0,$ thus $f = 0.$ 

  For any $st$-solid object $X = (U,V)$ we have to construct an $st$-solid object $X^+_+$ with
  $(X^+_+)^-_- = X.$ We write the global space $V$ of $X$ as the direct sum of indecomposable
  objects: They are of the form $[i,j]$ with $i\le s\le t \le j$. Note that 
  $[i,j]$ is in a unique way
  a subquotient of $[i-1,j+1]$. The direct sum $V^+_+$ 
  of the corresponding objects $[i-1,j+1]$ will be
  the global space of $X^+_+$. In this way, we have $(V^+_+)_a = V_a$ for $s\le a \le t,$ 
  thus we may consider for $s\le a \le t$ the subspace $U_a$ as a subspace of $(V^+_+)_a$.
  We take as $U(X^+_+)$ the space
  $\bigoplus_{s\le a \le t} U_a\oplus \bigoplus_{a'<s}(V^+_+)_{a'}.$ Clearly, if
  $X = (U,V)$ and $X' = (U',V')$ are $st$-solid, any homomorphism $f\:X \to X'$ can be
  extended in a unique way to a homomorphism $f^+_+\:X^+_+ \to (X')^+_+.$
\end{proof}

\subsection{Expansions and coexpansions of solid Kronecker families.}
\label{sec-six-three}

\begin{proposition}
  \label{prop-six-three}
  Let $M = M_c$ be a solid Kronecker family in $\widetilde{\Cal S}(n)$,
  and let $0 \le e,f \le \ell\cdot bM$ be integers. 
  Then there is a $\mathbb P^1$-family of indecomposable and pairwise non-isomorphic 
  solid objects $M_c[\ell;e,f]$ in
  $\widetilde{\Cal S}(n+2)$ with uwb-vector $(\ell\cdot uM + e, \ell\cdot wM + f, \ell\cdot bM).$
  The objects $M_c[\ell;0,f]$  and $M_c[\ell;e,0]$ 
  belong to $\widetilde{\Cal S}(n+1)$.
\end{proposition}

\medskip 
Note that $\pr M_c[\ell;e,f] = (pM+\frac e{\ell\cdot bM},rM+\frac f{\ell\cdot bM}).$

\medskip
\begin{corollary}
  \label{cor-six-three}
  Let $M = M_c$ be a solid Kronecker family in $\widetilde{\Cal S}(n)$.
  Let $0 \le e,f \le \ell\cdot bM.$ Then 
  $(pM+\frac e{\ell\cdot bM},rM+\frac f{\ell\cdot bM})$ is a BTh-vector in $\mathbb T(n+2)$
  and $(pM+\frac e{\ell\cdot bM},rM)$ 
  and $(pM,rM+\frac f{\ell\cdot bM})$ are BTh-vectors in $\mathbb T(n+1).$
\end{corollary}

\medskip
The pr-vectors which we obtain in this way are the vectors in the convex hull of 
$(p,r) = \pr M,$  $(p,r+1),$ $(p+1,r+1)$, and $(p+1,r).$
$$  
{\beginpicture
   \setcoordinatesystem units <.404cm,.7cm>
   \setcoordinatesystem units <.808cm,1.4cm>
\multiput{$\bullet$} at 0 0  1 1  2 0  3 1 /
\setsolid
\plot 0 0  1 1  2 0  0 0 /
\plot 1 1  3 1  2 0 /
\plot 2 .01  -0.01 .01  0.99 1.01 /
\plot 2 -.01  .01 -.01  1.01 .99 /
\put{$M$} at -.4 -.1 
\put{$\ssize (p,r+1)$} at 2.7 -.1
\put{$\ssize (p+1,r+1)$} at 3.9 1.1 
\put{$\ssize (p+1,r)$} at 0.3 1.1 
\setshadegrid span <.3mm>
\vshade 0 0 0 <z,z,,> 1 0 1 <z,z,,> 2 0 1 <z,z,,> 3 1 1 /

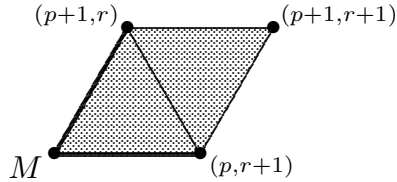
\endpicture}
$$

\medskip
We will show in Section~\ref{sec-six-six}
that the corollary implies Theorem~\ref{theoremfour} for $n \ge 8$, using the
symmetries of $\mathbb T(n)$.
As we will see in Section~\ref{sec-eight-five}, the corollary also yields
Theorem~\ref{theoremten} which 
concerns pr-vectors with boundary distance between 1 and 2.

\medskip
\begin{remark}
  If $(U,V)$ is an object and $(U \subseteq V' \subseteq V)$ with $bV' = bV$, 
  then we say that $(U,V)$ is an {\it expansion of} $(U,V')$.  
  Note that we have
  $p(U,V') = p(U,V)$ and $r(U,V') \le r(U,V)$ (looking at an expansion, we keep the
  subspace and its level, but increase the factor space $W = V/U$ and the colevel).

  The dual concept is that of a coexpansion: Here we increase the subspace and its level, but keep 
  the factor space $W$ and its colevel: Here, let $(U,V)$ be an object and $U' \subseteq U$ a
  subspace. If $b(V/U') = bV,$ then $(U,V)$ will be called a {\it coexpansion of} $(U/U',V/U').$
\end{remark}

\medskip
The objects constructed in the proof of the proposition will be expansions and coexpansions
of the given Kronecker family. 

\bigskip\bigskip 
For the proof of Proposition~\ref{prop-six-three},
we start with an $st$-solid Kronecker family $M = M_c$, say given by a Kronecker pair $X,Y$.
Then $N = M^+_+$ is a strongly solid Kronecker family, given by the Kronecker pair $X^+_+, Y^+_+.$
The construction of $N$ shows that $M$ can be recovered from $N$ 
in the following way: Let $U = U(N), V = V(N).$ 
We define $\soc N$
to be the object 
$(\soc V,\soc V).$ Then there is a canonical embedding $M \to N/\soc N$ with image
$(U,TV)/\soc N$ (thus, $M = (U/\soc V,TV/\soc V)$).

Now we are going to consider subquotients of the form $N'/N''$, where 
$N'' \subseteq \soc N$ and where $N'$ is a subobject of $N$ which contains $(U,TV).$

For $i>t$ we call a subobject $N' = (U,V')$ of $N = (U,V)$ {\it $i$-admissible} provided 
$$
 TV \subseteq V' \subseteq V,\quad V_j \subseteq V'\quad\text{for $j>i$,}\quad
\text{and} \quad V'_j\subseteq TV\quad\text{for $j<i$.}
$$
In this case we also call the subspace $V'$ of $V$ $i$-admissible.

\medskip
\begin{lemma}
  A subspace $V'$ of $V$ is admissible for some $i>t$ provided $V'$
  is comparable with the following filtration of $V$.  $\s$
  $$
  TV \subseteq TV+V_{n+1} \subseteq TV+V_{n+1}+V_{n} \subseteq \cdots 
  \subseteq  TV+V_{n+1}+\cdots+V_{t} = V.
  $$
\end{lemma}

Since the index of $TV$ in $V$ is equal to $bV$, we see: {\it For any natural number $f$
with $0 \le f \le bV$, there is an admissible subspace $V'$ of $V$ with dimension 
$\dim V - bV + f.$}

\medskip 
\begin{lemma}
  If $N'$ is $i$-admissible, then 
  the inclusion $M\to N/\soc N$ induces an isomorphism $M\to G_i(N')/\soc N$. $\s$
\end{lemma}

\bigskip
Dually, for $i<s$ we define a subobject $N''=(V'',V'')$ of $N$ to be $i$-{\it coadmissible} if
$$V''\subset\soc V,\qquad (\soc V)_j\subset V'' \quad\text{for}\; j>i,
\quad \text{and}\quad V_j''=0 \quad\text{for}\; j<i.$$
In this case, the subspace $V''$ of $V$ is also called $i$-coadmissible.

\medskip
\begin{lemma}
  A subspace $V''$ of $V$ is $i$-coadmissible for some $i<s$ provided $V''$ is comparable
  with the following filtration of $\soc V$. $\s$
$$0\subset \soc V\cap V_{s-1}\subset \soc V\cap(V_{s-1}+V_{s-2})\subset\cdots
\subset \soc V\cap(V_{s-1}+V_{s-2}+\cdots+V_0)=\soc V$$
\end{lemma}

           \medskip
           Since $\dim \soc V=bV$, for any natural number $e$ with $0\leq e\leq bV$
           there is a coadmissible subspace $V''$
of $V$ of dimension $bV-e$.

\medskip
\begin{lemma}
  If $N''$ is $i$-coadmissible, then the canonical map $TN\to M$ induces an isomorphism
  $G_i(TN/N'')[1]\to M$. $\s$
\end{lemma}

\smallskip
\begin{proof}[Proof of Proposition~\ref{prop-six-three}]
  We start with $M = M_c$, we form $M^+_+.$ 
  In this way, we obtain a Kronecker family and we consider the corresponding extensions
  $N = M^+_+[\ell]$ with $\ell$ factors isomorphic to $M^+_+$. We have $bN = \ell\cdot bM.$
  Recall that we denoted the subobject $(\soc V,\soc V)$ by $\soc N$. Correspondingly, let
  $TN = (U,TV).$
  
  In addition, there are given integers $0 \le e,f\le \ell\cdot bM.$ 
  The previous considerations show that there are submodules 
  $N',N''$ of $N$ with
  $$
  0 \subseteq N'' \subseteq \soc N \subset TN \subseteq N' \subseteq N
  $$
  such that $N'$ is $i$-admissible for some $i>t$, and $N''$ is $j$-coadmissible for some $j< s$,
  and such that $\soc N/N''$ has dimension $e$, whereas $N'/TN$ has dimension $f$. 
  
  Using the functors first $G_i$, then $G_j$, we see that the objects 
  $N'/N''$ are indecomposable and (in reference to the parameter $c$) pairwise non-isomorphic.
  This completes the proof.
\end{proof}

\subsection{A triangle, for $n = 7.$}
\label{sec-six-four}

Here, $n = 7.$ We want to construct a BTh-family for any rational vertex in 
the triangle with corners $(7/3,8/3),\ (2,2),\ (8/3,7/3)$ (see the illustration at the end
of Section~\ref{sec-six-five}). 

\medskip
We start with the $Q[1,7]$-modules $D = D_c$ with dimension vector
$\smallmatrix 1 & 2 & 2 & 2 & 1 \cr 1 & 2 & 3 & 4 & 3 & 2 & 1\endsmallmatrix$; say with
global space generated by $x_1,x_2,x_3,x_4$ in degrees $7,6,5,4,$ and with $\Lambda x_i$ of 
length $7,5,3,1$, respectively; with subspace generated by $T^2x_1+c_0Tx_2 + c_1x_3$ and
$T^2x_2+Tx_3+x_4$; here $c = (c_0:c_1)\in \mathbb P^1(k)$. 
It is the Kronecker family for the Kronecker subcategory of $\widetilde{\Cal S}$ 
given by the following Kronecker pair $X,Y$:

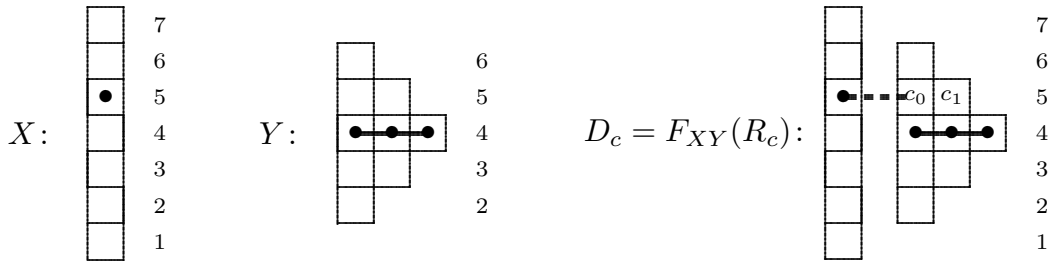
$$
{\beginpicture
    \setcoordinatesystem units <.4cm,.4cm>
\multiput{} at -1.5 0   3  7 /
\put{$X\:$} at -1.5 3.5
\plot 0 0  1 0  1 6   0 6  0 0 /
\plot 0 1  1 1  1 2  0 2  0 3  1 3  1 4  0 4  0 5  1 5 /
\plot 0 6  0 7  1 7  1 6 /
\put{$\scriptstyle 1$} at 2 0.5
\put{$\scriptstyle 2$} at 2 1.5
\put{$\scriptstyle 3$} at 2 2.5
\put{$\scriptstyle 4$} at 2 3.5
\put{$\scriptstyle 5$} at 2 4.5
\put{$\scriptstyle 6$} at 2 5.5
\put{$\scriptstyle 7$} at 2 6.5
\multiput{$\bullet$} at 0.5 4.5 /
\endpicture
}
\qquad
{\beginpicture
    \setcoordinatesystem units <.4cm,.4cm>
\multiput{} at -1.5 0   4  7 /
\put{$Y\:$} at -1.5 3.5
\plot 0 1  1 1  1 5  0 5  0 1 /
\plot 0 2  2 2  2 4  0 4  0 3  2 3 /
\plot 0 5  0 6  1 6  1 5  2 5  2 4  3 4  3 3  2 3 /
\put{$\scriptstyle 2$} at 4 1.5
\put{$\scriptstyle 3$} at 4 2.5
\put{$\scriptstyle 4$} at 4 3.5
\put{$\scriptstyle 5$} at 4 4.5
\put{$\scriptstyle 6$} at 4 5.5
\multiput{$\bullet$} at 0.5 3.5  1.5 3.5  2.5  3.5 /
\plot 0.5 3.5  2.5 3.5 /
\plot 0.5 3.45  2.5 3.45 /
\plot 0.5 3.55  2.5 3.55 /
\endpicture
}
\qquad\quad
{\beginpicture
    \setcoordinatesystem units <.4cm,.4cm>
\multiput{} at -2 0   6  7 /
\put{$D_c = F_{XY}(R_c)\:$} at -3.5 3.5
\plot 0 0  1 0  1 7  0 7  0 0 /
\plot  0 1  1 1  1 2  0 2  0 3  1 3  1 4  0 4  0 5  1 5  1 6  0 6 /
\put{$\scriptstyle 1$} at 6 0.5
\put{$\scriptstyle 2$} at 6 1.5
\put{$\scriptstyle 3$} at 6 2.5
\put{$\scriptstyle 4$} at 6 3.5
\put{$\scriptstyle 5$} at 6 4.5
\put{$\scriptstyle 6$} at 6 5.5
\put{$\scriptstyle 7$} at 6 6.5
\plot 2 1  3 1  3 6  2 6  2 1 /
\plot 2 2  4 2  4 5  2 5 /
\plot 2 3  5 3  5 4  2 4 /
\multiput{$\bullet$} at 0.5 4.5  2.5 3.5  3.5  3.5  4.5  3.5 /
\plot 2.5 3.5  4.5 3.5 /
\plot 2.5 3.55 4.5 3.55 /
\plot 2.5 3.45 4.5 3.45 /
\put{$\ss c_0$} at 2.5 4.5 
\put{$\ss c_1$} at 3.5 4.5 
\setdashes <.9mm>
\plot 2.2 4.5   0.5 4.5 /
\plot 2.2 4.55  0.5 4.55 /
\plot 2.2 4.45  0.5 4.45 /

\endpicture}
$$

Thus we deal with the following module $D$ shown in the middle:
$$
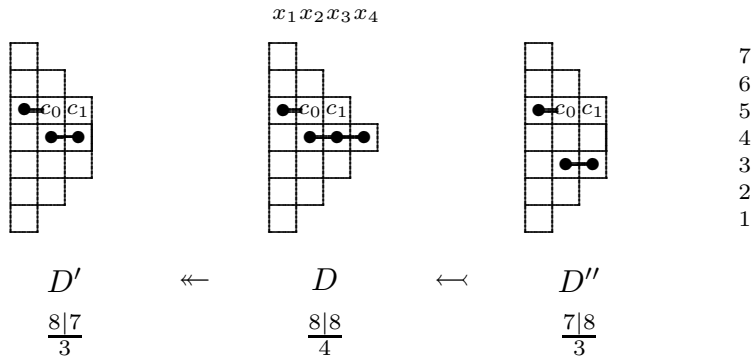
{\beginpicture
   \setcoordinatesystem units <.3cm,.3cm>


\put{\beginpicture
\multiput{} at 0 0  4 8.5 /
\plot 0 0  1 0  1 7  0 7  0 0 /
\plot 0 1  2 1  2 6  0 6 /
\plot 0 2  3 2  3 5  0 5 /
\plot 0 3  3 3  3 4  0 4 /
\multiput{$\bullet$} at 0.5 4.5  1.5 3.5  2.5 3.5 /
\put{$\ssize c_0$} at 1.5 4.5 
\put{$\ssize c_1$} at 2.5 4.5 
\plot 1.5 3.5   2.5 3.5 /
\plot 1.5 3.45   2.5 3.55 /
\plot 1.5 3.55   2.5 3.55 /

\plot 0.5 4.5   1.2 4.5 /
\plot 0.5 4.45   1.2 4.45 /
\plot 0.5 4.55   1.2 4.55 /
\endpicture} at 0 0

\put{\beginpicture
\multiput{} at 0 0  4 8.5 /
\plot 0 0  1 0  1 7  0 7  0 0 /
\plot 0 1  2 1  2 6  0 6 /
\plot 0 2  3 2  3 5  0 5 /
\plot 0 3  4 3  4 4  0 4 /
\multiput{$\bullet$} at 0.5 4.5  1.5 3.5  2.5 3.5  3.5 3.5 /
\put{$\ssize c_0$} at 1.5 4.5 
\put{$\ssize c_1$} at 2.5 4.5 
\plot 1.5 3.5   3.5 3.5 /
\plot 1.5 3.45   3.5 3.55 /
\plot 1.5 3.55   3.5 3.55 /

\plot 0.5 4.5   1.2 4.5 /
\plot 0.5 4.45   1.2 4.45 /
\plot 0.5 4.55   1.2 4.55 /

\put{$\ssize x_1$} at 0.6 8
\put{$\ssize x_2$} at 1.6 8
\put{$\ssize x_3$} at 2.6 8
\put{$\ssize x_4$} at 3.6 8
\endpicture} at 9,6 0 
 
\put{\beginpicture
\multiput{} at 0 0  4 8.5 /
\plot 0 0  1 0  1 7  0 7  0 0 /
\plot 0 1  2 1  2 6  0 6 /
\plot 0 2  3 2  3 5  0 5 /
\plot 0 3  3 3  3 4  0 4 /
\multiput{$\bullet$} at 0.5 4.5  1.5 2.5  2.5 2.5 /
\put{$\ssize c_0$} at 1.5 4.5 
\put{$\ssize c_1$} at 2.5 4.5 
\plot 1.5 2.5   2.5 2.5 /
\plot 1.5 2.45  2.5 2.55 /
\plot 1.5 2.55  2.5 2.55 /

\plot 0.5 4.5   1.2 4.5 /
\plot 0.5 4.45   1.2 4.45 /
\plot 0.5 4.55   1.2 4.55 /
\endpicture} at 19 0 
\put{\beginpicture
\multiput{} at 0 0  4 8.5 /
\put{$\ssize 1$} at 0 0.5
\put{$\ssize 2$} at 0 1.5
\put{$\ssize 3$} at 0 2.5
\put{$\ssize 4$} at 0 3.5
\put{$\ssize 5$} at 0 4.5
\put{$\ssize 6$} at 0 5.5
\put{$\ssize 7$} at 0 6.5
\endpicture} at 27 0 
\put{$D'$} at 0 -6
\put{$\twoheadleftarrow$} at 4.8 -6
\put{$D$} at 9.6 -6
\put{$\leftarrowtail$} at 14.3 -6
\put{$D''$} at 19 -6

\put{$\frac{8|7}3$} at 0 -8
\put{$\frac{8|8}4$} at 9.6 -8
\put{$\frac{7|8}3$} at 19 -8

\endpicture}
$$
There is the factor module $D'$ of $D$, it is 
shown on the left (the kernel of $D \to D'$ is generated by $x_4\in VD$ and 
is of the form $S_4$).
And there is the submodule $D''$ of $D$, with $V(D'')$
generated by the elements $x_1,x_2,x_3$ of $VD$; it is  
shown on the right (note that $D/D''$ is of the form $\overline S_4$, with global space
generated by the residue class of $x_4$). We have added the uwb-vectors of $D',D,D''.$

\medskip
Since $D_c$ is a Kronecker family, there are given the $\mathbb P^1$-families $D[\ell]$ for
all $\ell\in \mathbb N_1.$ Let $d',d,d''$ be natural numbers (not all zero). 
We start with the self-extension $B = D[d'+d+d'']$ using $d'+d+d''$
copies of $D$. Let $B'$ be the submodule of $B$ generated by the first $d'$ copies of $x_4$.
Let $B''$ be the submodule of $B$ generated by all copies of $x_1,x_2,x_3$ and the first $d'+d$
copies of $x_4$, thus $B' \subseteq B''\subseteq B$. The module we are interested in is $B''/B'.$
By construction, it has a filtration, going upwards, with first 
$d'$ copies of $D'$, then $d$ copies of $D$, and finally, $d''$ copies of
$D'',$ thus its uwb-vector is
$$
  \uwbb B''/B' = d'(8,7,3)+d(8,8,4)+d''(7,8,3).
$$

\medskip
We use the functor $G_4$ introduced in 6.1. 
The functor $G_4$ sends all three objects $D',$ $D$ and $D''$
to $M_c;$ it sends $D[d'+d+d'']$ and also $B''/B'$ to $M_c[d'+d+d''].$ 
The functor $G_4$ shows that all the modules $B''/B'$ are indecomposable,
and also, for pairwise different elements $c\in \mathbb P^1(k),$ 
pairwise non-isomorphic.  This shows that we 
obtain a BTh-family $B''/B'.$

\subsection{The trapezoid $T$ in $\mathbb T(7).$}
\label{sec-six-five}

We start with the $Q[1,7]$-modules $E = E_c$ with dimension vector
$\smallmatrix 1 & 2 & 3 & 2 & 1 \cr 1 & 2 & 3 & 4 & 3 & 2 & 1\endsmallmatrix$; say with
global space generated by $x_1,x_2,x_3,x_4$ in degrees $7,6,5,4,$ and with $\Lambda x_i$ of 
length $7,5,3,1$, respectively; with subspace generated by $T^2x_1+c_0Tx_2 + c_1x_3,$ 
$T^2x_2+Tx_3+x_4$; and $T^2x_3$, here $c = (c_0:c_1)\in \mathbb P^1(k)$. 
It is the Kronecker family for the Kronecker subcategory of $\widetilde{\Cal S}$ 
given by the following Kronecker pair $X,Y'$:
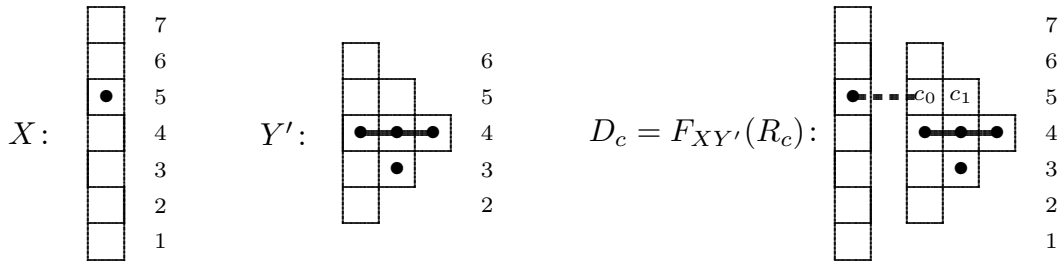
$$
{\beginpicture
    \setcoordinatesystem units <.4cm,.4cm>
\multiput{} at -1.5 0   3  7 /
\put{$X\:$} at -1.5 3.5
\plot 0 0  1 0  1 6   0 6  0 0 /
\plot 0 1  1 1  1 2  0 2  0 3  1 3  1 4  0 4  0 5  1 5 /
\plot 0 6  0 7  1 7  1 6 /
\put{$\scriptstyle 1$} at 2 0.5
\put{$\scriptstyle 2$} at 2 1.5
\put{$\scriptstyle 3$} at 2 2.5
\put{$\scriptstyle 4$} at 2 3.5
\put{$\scriptstyle 5$} at 2 4.5
\put{$\scriptstyle 6$} at 2 5.5
\put{$\scriptstyle 7$} at 2 6.5
\multiput{$\bullet$} at 0.5 4.5 /
\endpicture
}
\qquad
{\beginpicture
    \setcoordinatesystem units <.4cm,.4cm>
\multiput{} at -1.5 0   4  7 /
\put{$Y'\:$} at -1.5 3.5
\plot 0 1  1 1  1 5  0 5  0 1 /
\plot 0 2  2 2  2 4  0 4  0 3  2 3 /
\plot 0 5  0 6  1 6  1 5  2 5  2 4  3 4  3 3  2 3 /
\put{$\scriptstyle 2$} at 4 1.5
\put{$\scriptstyle 3$} at 4 2.5
\put{$\scriptstyle 4$} at 4 3.5
\put{$\scriptstyle 5$} at 4 4.5
\put{$\scriptstyle 6$} at 4 5.5
\multiput{$\bullet$} at 0.5 3.5  1.5 3.5  2.5  3.5 /
\multiput{$\bullet$} at 1.5 2.5   /
\plot 0.5 3.5  2.5 3.5 /
\plot 0.5 3.45  2.5 3.45 /
\plot 0.5 3.55  2.5 3.55 /
\endpicture
}
\qquad\quad
{\beginpicture
    \setcoordinatesystem units <.4cm,.4cm>
\multiput{} at -2 0   6  7 /
\put{$E_c = F_{XY'}(R_c)\:$} at -3.5 3.5
\plot 0 0  1 0  1 7  0 7  0 0 /
\plot  0 1  1 1  1 2  0 2  0 3  1 3  1 4  0 4  0 5  1 5  1 6  0 6 /
\put{$\scriptstyle 1$} at 6 0.5
\put{$\scriptstyle 2$} at 6 1.5
\put{$\scriptstyle 3$} at 6 2.5
\put{$\scriptstyle 4$} at 6 3.5
\put{$\scriptstyle 5$} at 6 4.5
\put{$\scriptstyle 6$} at 6 5.5
\put{$\scriptstyle 7$} at 6 6.5
\plot 2 1  3 1  3 6  2 6  2 1 /
\plot 2 2  4 2  4 5  2 5 /
\plot 2 3  5 3  5 4  2 4 /
\multiput{$\bullet$} at 0.5 4.5  2.5 3.5  3.5  3.5  4.5  3.5 /
\multiput{$\bullet$} at 3.5 2.5   /
\plot 2.5 3.5  4.5 3.5 /
\plot 2.5 3.55 4.5 3.55 /
\plot 2.5 3.45 4.5 3.45 /
\put{$\ss c_0$} at 2.5 4.5 
\put{$\ss c_1$} at 3.5 4.5 
\setdashes <.9mm>
\plot 2.2 4.5   0.5 4.5 /
\plot 2.2 4.55  0.5 4.55 /
\plot 2.2 4.45  0.5 4.45 /

\endpicture}
$$

We consider also the submodule $D$ of $E$:
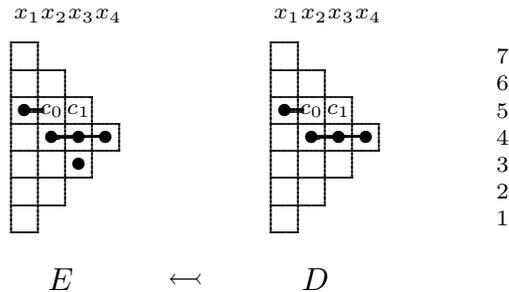
$$
{\beginpicture
   \setcoordinatesystem units <.3cm,.3cm>
\put{\beginpicture
\multiput{} at 0 0  4 8.5 /
\plot 0 0  1 0  1 7  0 7  0 0 /
\plot 0 1  2 1  2 6  0 6 /
\plot 0 2  3 2  3 5  0 5 /
\plot 0 3  4 3  4 4  0 4 /
\multiput{$\bullet$} at 0.5 4.5  1.5 3.5  2.5 3.5  3.5 3.5  2.5 2.5 /
\put{$\ssize c_0$} at 1.5 4.5 
\put{$\ssize c_1$} at 2.5 4.5 
\plot 1.5 3.5   3.5 3.5 /
\plot 1.5 3.45   3.5 3.55 /
\plot 1.5 3.55   3.5 3.55 /

\plot 0.5 4.5   1.2 4.5 /
\plot 0.5 4.45   1.2 4.45 /
\plot 0.5 4.55   1.2 4.55 /

\put{$\ssize x_1$} at 0.6 8
\put{$\ssize x_2$} at 1.6 8
\put{$\ssize x_3$} at 2.6 8
\put{$\ssize x_4$} at 3.6 8
\endpicture} at 0 0

\put{\beginpicture
\multiput{} at 0 0  4 8.5 /
\plot 0 0  1 0  1 7  0 7  0 0 /
\plot 0 1  2 1  2 6  0 6 /
\plot 0 2  3 2  3 5  0 5 /
\plot 0 3  4 3  4 4  0 4 /
\multiput{$\bullet$} at 0.5 4.5  1.5 3.5  2.5 3.5  3.5 3.5 /
\put{$\ssize c_0$} at 1.5 4.5 
\put{$\ssize c_1$} at 2.5 4.5 
\plot 1.5 3.5   3.5 3.5 /
\plot 1.5 3.45   3.5 3.55 /
\plot 1.5 3.55   3.5 3.55 /

\plot 0.5 4.5   1.2 4.5 /
\plot 0.5 4.45   1.2 4.45 /
\plot 0.5 4.55   1.2 4.55 /

\put{$\ssize x_1$} at 0.6 8
\put{$\ssize x_2$} at 1.6 8
\put{$\ssize x_3$} at 2.6 8
\put{$\ssize x_4$} at 3.6 8
\endpicture} at 9.6 0 
 
\put{\beginpicture
\multiput{} at 0 0  4 8.5 /
\put{$\ssize 1$} at 0 0.5
\put{$\ssize 2$} at 0 1.5
\put{$\ssize 3$} at 0 2.5
\put{$\ssize 4$} at 0 3.5
\put{$\ssize 5$} at 0 4.5
\put{$\ssize 6$} at 0 5.5
\put{$\ssize 7$} at 0 6.5
\endpicture} at 18 0 
\put{$E$} at -.2 -6
\put{$\leftarrowtail$} at 4.4 -6
\put{$D$} at 9.2 -6
\endpicture}
$$

Note that we have $D_{4'}\cap \soc D = 0$ and $E_{4'}\cap \soc E = 0,$ 
thus the modules $D' = H_4D$ and $E' = H_4E$ are again Gorenstein-projective. 

\bigskip
We look at the trapezoid $T$ with corners $D,\ D',\  E',\ E.$ 

\medskip 
Using the canonical embedding $\mu\:D \to E$, 
we can interpolate between $D$ and $E$ by looking at extensions
induced by the extension $E[d+e],$ where $d,e$ are non-negative integers, 
namely
$$
\CD
 0 @>>> E[e] @>>> E[e]\under D[d] @>>> D[d] @>>> 0, \cr
 @.     @|            @VVV             @VV\mu[d] V   \cr
 0 @>>> E[e] @>>> E[d+e]          @>>> E[d] @>>> 0.
\endCD
$$

Then we can factor out from $E[e]\under D[d]$
copies of $x_4$, say $f$ copies, where $0 \le f \le d+e.$ 
In this way, we obtain
an object $B$ with pr-vector inside the trapezoid $T$ (it lies on the line which
connects $(0,1)$ with the pr-vector of $E[e]\under D[d]$; but for $f> 0,$ 
the pr-vector obtained will lie outside of this interval).

Note that {\it any pr-vector inside the trapezoid $T$ is obtained in this way!}
Let us outline how we find numbers $d,\ e,\ f,$ when we 
start with a pr-vector $(p,r)$ in the trapezoid $T$.
Recall that our aim is to find modules $M$ with pr-vector $(p,r)$ 
thus the uwb-vector of $M$ has to be of the form
$$
 d\frac{u|w}b D' + e\frac{u|w}b E' + f'\frac{u|w}b S_4
$$
for some non-negative integers $d,e,f',$ therefore equal to
\begin{align*}
 d\frac{8|7}3  + e\frac{9|6}3 + f'\frac{0|1}1 &= \frac {8d+9e\, |\, 7d+6e+f'}{3d+3e+f'} \cr
    &=  d\frac{8|8}4  + e\frac{9|7}4 - (d+e-f')\frac{0|1}1 \cr
    &=   d\frac{u|w}b D + e\frac{u|w}b E - (d+e-f')\frac{u|w}b S_4.
\end{align*}
If $(p,r)$ belongs to the triangle with corners $S_4, D, E$, then it lies on
the segment $D, E.$ As a consequence, we see 
that $f = d+e-f'$ cannot be negative. 
Also, $d+e-f = f' \ge 0$ shows that  $f\le d+e,$ thus, altogether we have $0 \le f \le d+e.$
This provides the required 
non-negative numbers $d,e,f$ (and not all can be zero). 

\medskip
It remains to be shown that all the modules $B = B_c$ obtained 
from $E[e]\under D[d] = E_c[e]\under D_c[d]$ by factoring out
copies of $x_4$ are indecomposable and for different $c$ non-isomorphic. 
We use the functor $G_z$ and $H_z$ defined in 6.1.

The functor $G_3H_4$ sends both objects $D,$ $E$ to $M_c;$ 
it sends $D[d+e],$ $E[e]\under D[d],$ $E[d+e],$ but also any $B$ (obtained from 
$E[e]\under D[d]$ by factoring out copies of $x_4$) to $M_c[d+e].$
In this way, the functor $G_3H_4$ shows that all the modules $B$ are indecomposable,
and also, for pairwise different elements $c\in \mathbb P^1(k),$ 
pairwise non-isomorphic.  This shows that we 
obtain a BTh-family. 

Altogether, we obtain BTh-families inside the trapezoid $T$ with corners
$D, D' = H_4D, E' = H_4E, E.$ 

\bigskip
Let us present in $\mathbb T(7)$ the triangle with corners $D,\ D' = H_4D,\ D''$
 as well as the trapezoid $T$ with corners
$D,\ D' = H_4D,\ E' = H_4E,\ E.$

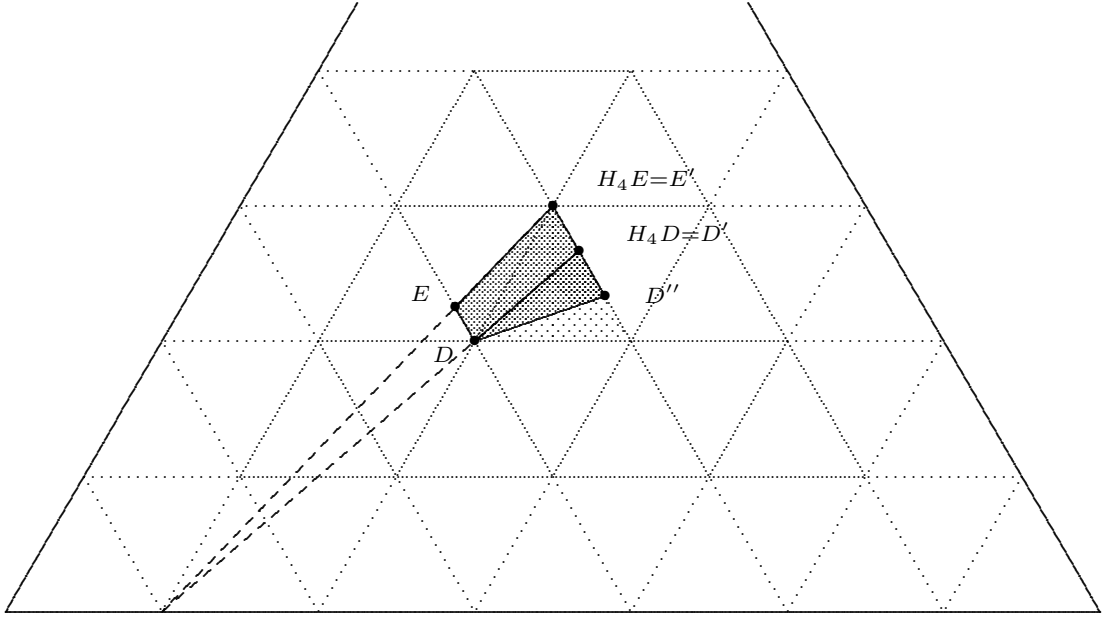
\medskip
$$  
{\beginpicture
  \setcoordinatesystem units <.8371cm,1.45cm>
\multiput{} at -6 -3  6 -7 /
\plot  -2.5 -2.5  -7 -7   /
\plot  2.5 -2.5  7 -7 /
\plot -7 -7  7 -7 /
\setdots <1mm>
\plot -6 -6  -4 -6  -5 -5  -3 -5  -4 -4  -2 -4  -3 -3  -1 -3 /
\plot  6 -6   4 -6   5 -5   3 -5   4 -4   2 -4   3 -3   1 -3 /
\plot -6 -6  -5 -7  -4 -6  -3 -7  -2 -6  -1 -7  0 -6  1 -7  2 -6  3 -7  4 -6  5 -7  6 -6 /
\setdots <.5mm>
\plot -4 -6  4 -6  1 -3  -2 -6  -3 -5  3 -5  2 -6  -1 -3  -4 -6  /
\plot -2 -4  2 -4  0 -6  -2 -4 /
\plot -1 -3  1 -3 /

\setsolid
\plot -1 -5  -1.25 -4.75  0 -4  0.333  -4.333 /

\setshadegrid span <.6mm>
\vshade -1 -5 -5 <z,z,,> 0 -5 -4 <z,z,,> 1 -5 -5 /

\setshadegrid span <.3mm>
\vshade -1 -5 -5 <z,z,,> 0.333 -4.75  -4.333 <z,z,,> 0.667 -4.667 -4.667  /
\vshade -1.25 -4.75 -4.75 <z,z,,> -1 -5 -4.6 <z,z,,> 0 -4.6 -4 <z,z,,>
    0.333 -4.75  -4.333 <z,z,,> 0.667 -4.667 -4.667  /

\put{$\ssize E$} at -1.7 -4.65
\put{$\ssize D$} at  -1.4 -5.1

\put{$\ssize H_4E = E'$} at 1.2  -3.8
\put{$\ssize H_4D = D'$} at  1.6 -4.2
\put{$\ssize D''$} at 1.4 -4.65

\multiput{$\ssize \bullet$} at 0.333 -4.333  0.667 -4.667  0 -4  -1 -5  -1.25 -4.75 /

\setsolid
\plot 0.333 -4.333  0.667 -4.667  /
\plot 0.333 -4.333 -1 -5 /
\plot -1 -5  0.667 -4.667  /

\setdashes <1mm>
\plot 0.333 -4.333  -5 -7 /
\plot 0 -4  -5 -7 /
\endpicture}
$$

\subsection{Proof of Theorem~\ref{theoremfour}.}
\label{sec-six-six}

We are going to use the previous results in order to prove Theorem~\ref{theoremfour}.
As we will see, we only have to invoke the symmetries for $\mathbb T(m)$,
with $m \le n.$

First, let $n = 7.$ In Section~\ref{sec-six-five}, we have found a trapezoid of BTh-vectors.
Using duality, we see that also the reflected trapezoid consists of BTh-vectors.
These two trapezoids, together with the triangle exhibited in Section~\ref{sec-six-four}, provide  
BTh-vectors with boundary distance at least 2.

Next, $n = 8.$ According to Corollary~\ref{cor-six-three},
we have BTh-vectors inside the parallelogram
given by $(2,2)$, $(3,2)$, $(3,3)$, $(2,3)$.
The rotation $\rho$ shows that all vectors with boundary distance at least 2 are BTh-vectors.

Finally, for $n\ge 9$, we use induction on $n$. As we have shown already, if $n = 8$,
any pr-vector in $\mathbb T(n)$ 
with boundary distance at least 2, is a BTh-vector. 
Thus, let $n \ge 9.$ By induction, 
any pr-vector in $\mathbb T(n-1)$ 
with boundary distance at least 2, is a BTh-vector. 
If we deal with a vector $\mathbf x$ 
in $\mathbb T(n)$, then its image under $\rho$ or $\rho^2$
will belong to $\mathbb T(n-1)$. Since $\rho(\mathbf x)$ or $\rho^2(\mathbf x)$ 
is a BTh-vector, also $\mathbf x$ itself is a BTh-vector.
$\s$

\subsection{BTh vectors in $\mathbb T(7)$.}
\label{sec-six-seven}

\begin{remark}
  Let us look again at $\mathbb T(7)$. Starting with Sections~\ref{sec-six-two} and \ref{sec-six-three}
  and using duality and rotation, we cover the following shaded region which includes properly
  the triangle of all pr-vectors of boundary distance at least 2.
\end{remark}

$$  
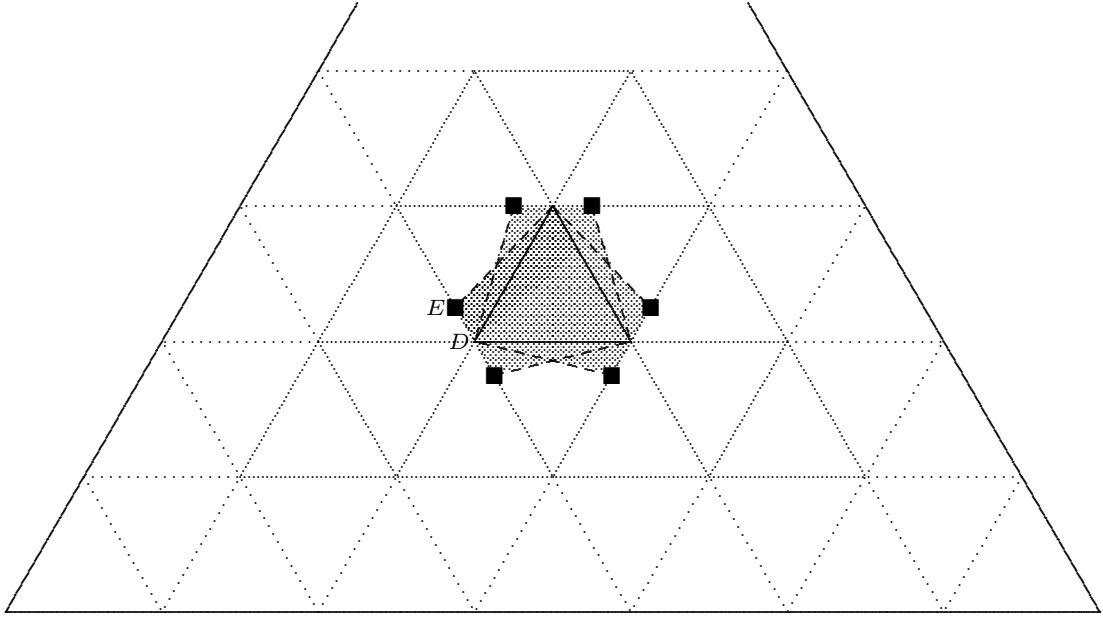
{\beginpicture
    \setcoordinatesystem units <.57735cm,1cm>
  \setcoordinatesystem units <.83713cm,1.45cm>
\multiput{} at -6 -3  6 -7 /
\plot  -2.5 -2.5  -7 -7   /
\plot  2.5 -2.5  7 -7 /
\plot -7 -7  7 -7 /
\setdots <1mm>
\plot -6 -6  -4 -6  -5 -5  -3 -5  -4 -4  -2 -4  -3 -3  -1 -3 /
\plot  6 -6   4 -6   5 -5   3 -5   4 -4   2 -4   3 -3   1 -3 /
\plot -6 -6  -5 -7  -4 -6  -3 -7  -2 -6  -1 -7  0 -6  1 -7  2 -6  3 -7  4 -6  5 -7  6 -6 /
\setdots <.5mm>
\plot -4 -6  4 -6  1 -3  -2 -6  -3 -5  3 -5  2 -6  -1 -3  -4 -6  /
\plot -2 -4  2 -4  0 -6  -2 -4 /
\plot -1 -3  1 -3 /

\setshadegrid span <.6mm>
\vshade -1 -5 -5 <z,z,,> 0 -5 -4 <z,z,,> 1 -5 -5 /

\setshadegrid span <.3mm>
\vshade -1.25 -4.75 -4.75 <z,z,,> -1 -5 -4.6 <z,z,,> 0 -4.5 -4 <z,z,,> 0.333 -4.333 -4.333 /

\put{$\ssize D$} at -1.2 -5
\put{$\ssize E$} at -1.5 -4.75

\multiput{$\ssize \blacksquare$} at  -1.25 -4.75  1.25 -4.75  .5 -4  -.5 -4  .75 -5.25  -.75 -5.25 /

\setdashes <1mm>
\plot 0 -4  -1.25 -4.75 /
\plot -.5 -4  -1 -5 /
\plot 0 -4  1.25 -4.75 /
\plot .5 -4  1 -5 /
\plot -.75 -5.25  1 -5 /
\plot .75 -5.25  -1 -5 /

\setshadegrid span <.3mm>
\vshade -1.25 -4.75 -4.75 <z,z,,> 
    -.75 -5.25 -4.45 <z,z,,>
    -.5 -5.2  -4 <z,z,,> 
    0 -5.15 -4 <z,z,,> 
    .5 -5.2  -4 <z,z,,>  
    .75  -5.25 -4.45 <z,z,,>   
    1.25 -4.75 -4.75 /

\setsolid
\plot 0 -4  1 -5  -1 -5  0 -4 /

\endpicture}
$$

\bigskip

\section{Objects with level $1$.}
\label{sec-seven}

Here we want to consider some examples of objects in $\Cal S(n)$ with level 1.

\subsection{Objects without $0$-pickets as direct summands.}
\label{sec-seven-one}

We return to the assertion of Proposition~\ref{prop-four-one}:

\medskip
\begin{proposition}
  Let $X$ be a non-zero object in $\Cal S$ which has no direct summands
  which are a $0$-pickets. Then the following assertions hold:
  \begin{itemize}[leftmargin=3em]
  \item[\rm(a)] $pX \ge 1.$
  \item[\rm(a$'$)] $uX \ge bX.$
  \item[\rm(b)] All factors $F_i$ of any telescope filtration satisfy $pF_i \ge 1.$
  \item[\rm(b$'$)] All factors $F_i$ of any telescope filtration satisfy $uF_i \ge 1.$
  \item[\rm(b$''$)] No factor of a telescope filtration of $X$ is a $0$-picket.
  \end{itemize}
\end{proposition}

\begin{proof}
  If $X$ is non-zero, the assertions (a) and (a$'$) are obviously  
  equivalent, as are the assertions (b), (b$'$) and (b$''$). 
  Thus, assume now that $X$ is non-zero and has no direct summand which
  is a $0$-picket. According to Proposition~\ref{prop-four-one}, the assertion (b$'$) holds true. 
  Finally, we show that (b$')$ implies (a$'$): The telescope filtration $(F_i)_i$ has
  $bX$ factors, thus, if (b$'$) is satisfied, then $uX = \sum_{i=1}^{bX} uF_i \ge bX.$
\end{proof}

\medskip 
(Actually, our interest in the level 
of the indecomposable objects in $\Cal S$ started with this observation!)

\bigskip
There is the following consequence:

\medskip
\begin{corollary}
  Let $X$ be a non-zero object in $\Cal S$ which has no direct summand
  which is a $0$-picket. Then the following assertions are equivalent:
  \begin{itemize}[leftmargin=3em]
  \item[\rm(i)] $pX = 1.$
  \item[\rm(ii)] All factors $F_i$ of any telescope filtration satisfy $pF_i = 1.$
  \item[\rm(iii)] There exists a telescope filtration whose factors $F_i$ satisfy $pF_i = 1.$
  \end{itemize}
\end{corollary}

\medskip
\begin{remark}
  {\bf Two examples.} Let us look at objects $X$ in $\Cal S$ which have 
  a filtration whose factors are level 1 pickets. Here are two examples which one should have in
  mind. First, there is the exact sequence 
  $$
  0\to ([1],[1])\to ([2],[3])\to ([1],[2])\to 0
  $$
  which shows that we may have $pX > 1.$ Of course, the sequence is not g-split. 
  
  Second there is an exact sequence
  $$
  0\to ([1],[2])\to ([2],[2])\oplus (0,[1]) \to ([1],[1])\to 0
  $$
  (mentioned already in Section~\ref{sec-four}). This sequence is g-split. The middle term $X$ has a
  direct summand which is a $0$-picket. 
\end{remark}

\subsection{The indecomposable objects $X = (U,V)$ with $bU = 1 = pX$.}
\label{sec-seven-two}

We refer to Section~\ref{sec-thirteen} for a discussion of indecomposable objects
with $bU=1$.
We say that a partition $(\lambda_1,\dots,\lambda_b)$
is {\it strongly decreasing}, provided $\lambda_i-\lambda_{i+1} \ge 2$ for $1\le i < b.$

\medskip
\begin{proposition}
  There is a bijection between, on the one hand,
  the isomorphism classes of the
  indecomposable objects
  $X = (U,V)$ in $\Cal S(n)$ with $pX = 1$ and $bU = 1,$ and,  on the other hand, the
  strongly decreasing partitions bounded by $n$, 
  by sending  $X = (U,V)$ to $[V].$
\end{proposition}

\medskip
If $X$ corresponds under this bijection to the partition $\lambda$, 
we write $X = C_\lambda.$

\smallskip
\begin{proof}
  Let $X = (U,V)$ be indecomposable in $\Cal S(n)$
  with $bU = 1.$ We decompose $V = \bigoplus_{i=1}^b
  V_i$ with all $V_i$ indecomposable in $\Cal N(n)$. Let $|V_i| = \lambda_i$, and we can
  assume that $\lambda_1 \ge \lambda_2 \ge \cdots \ge \lambda_b.$ 
  Using the projection maps $V \to V_i$,
  the embedding $u\:U \to V$ yields (non-zero) maps $u_i:U \to V_i.$ Since $bU = 1$, the 
  $\Lambda$-module $U$ is generated by an element $x$. 
  For any $i$, we can choose a generator $y_i\in V_i$
  such that $u_i(x) = T^{m_i}y_i$, with $0 \le m_i < \lambda_i.$
  
  Looking at the map $U \to V_i\oplus V_{i+1}$, we see that we must have 
  $m_i > m_{i+1}$ and $\lambda_i-m_i > \lambda_{i+1}-m_{i+1},$
  otherwise $X$ will be decomposable. In particular, it follows that the sequence
  $(\lambda_1,\lambda_2,\dots,\lambda_b)$ has to be strongly decreasing. 
  
  Since $u$ is injective, $u_1$ has to be injective, thus $U$ has length $\lambda_1-m_1$
  and there is the sequence
  \begin{equation*}
    |U| = \lambda_1-m_1 > \lambda_2-m_2 > \cdots > \lambda_b-m_b \ge 1. \tag{$*$}
  \end{equation*}
  
  Now, let us assume, in addition, that $pX = 1,$ thus  $b = bV = |U|.$
  The sequence $(*)$ is strictly decreasing, starts with $b$, and has $b$ terms,   
  thus it has to be the sequence $b, b-1, b-2,\dots, 2,1$. 
  Thus $\lambda_i-m_i = b-i+1.$ But this means that $m_i = \lambda_i-b+i-1,$
  thus $m_i$ is determined by the partition $\lambda$. 
  \smallskip
  
  Conversely, starting with a strongly decreasing partition $\lambda$, 
  let $m_i = \lambda_i-b+i-1.$ Let $V_i$ be the indecomposable $\Lambda$-module
  of length $\lambda_i$ generated by an element $x_i$.
  Let $V$ be the direct sum of the $\Lambda$-modules $V_i$.
  Let $U$ be the submodule of $V$ generated by $(T^{m_1}x_1,\dots,T^{m_b}x_b).$
  Then, clearly, $X = (U,V)$ is indecomposable.
  By construction, $bU = 1$ and $|U| = b = bV,$ thus also $pX = 1.$
\end{proof}

\medskip
\begin{example} For the partition $(7,4,2)$, the object $C_{(7,4,2)}$
may be visualized as follows:
$$
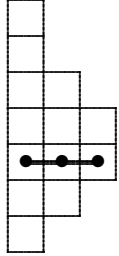
{\beginpicture
    \setcoordinatesystem units <.4cm,.4cm>
\multiput{} at 0 0   3  7 /
\plot 0 0  1 0  1 7   0 7  0 0 /
\plot 0 1  2 1  2 5  0 5 /
\plot 0 2  3 2  3 4  0 4 /
\plot 0 3  3 3 /
\plot 0 6  1 6 /
\multiput{$\bullet$} at 0.5 2.5  1.5 2.5  2.5 2.5  /
\plot 0.5 2.5  2.5 2.5 /
\plot 0.5 2.55  2.5 2.55 /
\plot 0.5 2.45  2.5 2.45 /
\endpicture}
$$
(Of course, the object $E_2^2$ mentioned already in the introduction, the only
indecomposable object of $\Cal S(3)$ which is not a picket, has to be mentioned here:
We have $E_2^2 = C_{(3,1)}$.)
\end{example}

\bigskip 
In general, if $X$ is indecomposable in $\Cal S$ and $U(X)$ belongs to
$\Cal N(2),$ then $X$ is a picket or a bipicket, see \cite[Proposition 3.4]{S1}.
This shows that {\it an
indecomposable object $X = (U,V)$ with $bX \le 3$ and $pX = 1$ satisfies $bU = 1.$}

\bigskip
For $n\le 5$, all indecomposable objects $X$ in $\Cal S(n)$ with
$pX = 1$ are of the form $C_\lambda.$ 
In $\Cal S(6)$, there is an indecomposable object $Y$ with $pY = 1$ and $bY = 4,$
which we will present now.

\subsection{An example from $\mathcal S_3(6)$.}
\label{sec-seven-three}

Here is an example of an indecomposable object $Y = (U,V)$ with 
$pY = 1$ and $bU > 1.$ It belongs to $\Cal S(6)$ and has $pY = 1$ and $bY = 4.$ 
Note that an indecomposable object $Y\in \Cal S$ with $pY = 1 = bU$ and $bY = 4$ has height at least 7.

\medskip
\begin{example}
  The indecomposable object $Y\in \Cal S(6)$ with $pY = 1$,
  and $bY = 4.$ 
$$
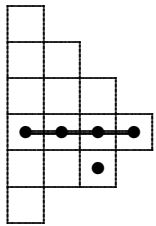
{\beginpicture
    \setcoordinatesystem units <.4cm,.4cm>
\multiput{} at 0 0   4  6 /
\plot 0 0  1 0  1 6   0 6  0 0 /
\plot 0 1  3 1  3 4  0 4 /
\plot 0 2  4 2  4 3  0 3 /
\plot 2 1  2 5  0 5 /
\multiput{$\bullet$} at 0.5 2.5  1.5 2.5  2.5 2.5  3.5 2.5  2.5  1.5 /
\plot 0.5 2.5  3.5 2.5 /
\plot 0.5 2.55  3.5 2.55 /
\plot 0.5 2.45  3.5 2.45 /
\endpicture}
$$
We may consider $Y$ as a g-split extension of $C_{(3)}$ by $C_{(6,4,1)}$.
In general, for any $c\in k$, there is the following 
g-split extension of $C_{(3)}$ by $C_{(6,4,1)}$:
$$
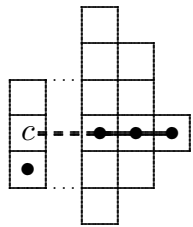
{\beginpicture
    \setcoordinatesystem units <.4cm,.4cm>
\multiput{} at -2 0   3 6 /
\plot 0 0  1 0  1 6   0 6  0 0 /
\plot 0 1  2 1  2 5  0 5 /
\plot 0 2  3 2  3 3  0 3 /
\plot 0 4  2 4 /
\multiput{$\bullet$} at 0.5 2.5  1.5 2.5  2.5 2.5 /
\plot 0.5 2.5  2.5 2.5 /
\plot 0.5 2.55  2.5 2.55 /
\plot 0.5 2.45  2.5 2.45 /

\plot -2 1  -1 1  -1 4  -2 4  -2 1 /
\plot -2 2  -1 2 /
\plot -2 3  -1 3 /
\multiput{$\bullet$} at -1.5 1.5  /
\put{$c$} at -1.5 2.5 
\setdashes <.9mm>
\plot -1.2 2.5  0.8 2.5 /
\plot -1.2 2.55  0.8 2.55 /
\plot -1.2 2.45  0.8 2.45 /

\setdots <1mm>
\plot -1 4  0 4 /
\plot -2 1  0 1 /

\endpicture}
$$
To be precise, we consider $X_c = (U,V)$, with $V$ being the direct sum of $\Lambda$-modules
$\Lambda x_i$, where $0 \le i \le 3$; with $|\Lambda x_i|$ equal to 
$3,6,4,1$ for $i = 0,\dots,3,$ and with $U$ being generated by the elements
$$
  u_1 = T^2x_0,\quad u_2 = cTx_0+T^3x_1+T^2x_2+x_3,
$$
with $c \in k.$ Note that $X_1 = Y.$ 
We have $X_0 = C_{(3)}\oplus C_{(6,4,1)}$, whereas, all objects $X_c$ with $c\neq 0$ 
are isomorphic (thus isomorphic to $Y$).
\end{example}

\subsection{A $\mathbb P^1$-family of indecomposable objects $X$ in
$\Cal S(9)$ with $pX = 1$.}
\label{sec-seven-four}

Similar to the previous object, let us now look at the following 
g-split extensions of $C_{(6,3)}$ by $C_{(9,7,4,1)}$
$$
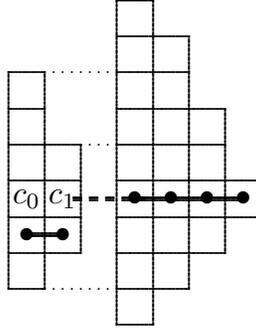
{\beginpicture
    \setcoordinatesystem units <.4cm,.4cm>
\multiput{} at -3 0   4  9 /
\plot 0 0  1 0  1 9   0 9  0 0 /
\plot 0 1  2 1  2 8  0 8 /
\plot 0 2  3 2  3 6  0 6 /
\plot 0 3  4 3  4 4  0 4 /
\plot 0 5  3 5 /
\plot 0 7  2 7 /
\multiput{$\bullet$} at 0.5 3.5  1.5 3.5  2.5 3.5  3.5 3.5 /
\plot 0.5 3.5  3.5 3.5 /
\plot 0.5 3.55  3.5 3.55 /
\plot 0.5 3.45  3.5 3.45 /

\plot -3 2  -1 2  -1 5  -2 5  -3 5 /
\plot -3 1  -2 1  -2 7  -3 7  -3 1  /
\plot -3 3  -1 3 /
\plot -3 4  -1 4 /
\plot -3 6  -2 6 /
\multiput{$\bullet$} at -1.5 2.5  -2.5 2.5 /
\plot -1.5 2.5  -2.5 2.5 /
\plot -1.5 2.55 -2.5 2.55 /
\plot -1.5 2.45 -2.5 2.45 /
\put{$c_0$} at -2.5 3.5 
\put{$c_1$} at -1.5 3.5 
\setdashes <.9mm>
\plot -1.2 3.5   0.8 3.5 /
\plot -1.2 3.55  0.8 3.55 /
\plot -1.2 3.45  0.8 3.45 /

\setdots <1mm>
\plot -2 7  0 7 /
\plot -1 5  0 5 /
\plot -2 1  0 1 /
\endpicture}
$$
To be precise, we consider $X = (U,V)$, with $V$ being the direct sum of $\Lambda$-modules
$\Lambda x_i$, where $0 \le i \le 5$; with $|\Lambda x_i|$ equal to 
$6,3,9,7,4,1$ for $i = 0,\dots,5,$ and with $U$ being generated by the elements
$$
  u_1 = T^4x_0+T^2x_1,\quad u_2 = c_0T^3x_0+c_1Tx_1+T^5x_2+T^4x_3+T^2x_4+x_5,
$$
with $c = (c_0,c_1)\in k^2.$ 
The object $X = X_c$ turns out to be indecomposable if and only if $c \neq 0.$
And indecomposable objects $X_c,X_{c'}$ are isomorphic if and only if $c$ and $c'$ are
multiples of each other (thus we obtain a $\mathbb P^1$-family of indecomposable objects indexed
by the projective line $\mathbb P^1(k)$). 

For any $c\in k^2,$ we have $|U| = 6$ and $bV = 6,$ thus $pX = 1.$

\medskip
For dealing with this example, we use the quiver $Q[1,9]$ as introduced in Section 2.3.
There are indecomposable representations $\widetilde Y,\ \widetilde Z$ of $Q[1,9]$
with dimension vector
$$
 \bdim \widetilde Y = \smallmatrix 0& 1& 1& 0& 0 &0 &0 &0 &0 \cr
                                     0& 1& 2& 2& 2 &1 &1 &0 &0 \endsmallmatrix, \qquad
 \bdim \widetilde Z = \smallmatrix 1& 1& 1& 1& 0 &0 &0 &0 &0 \cr
                                     1& 2& 3& 4& 3 &3 &2 &2 &1 \endsmallmatrix,
$$
and $\widetilde Y$ is sent under $\pi$ to $C_{(6,3)}$, whereas 
$\widetilde Z$ is sent to $C_{(9,7,4,1)}$. 
One has to verify that  
$\widetilde Y$ and $\widetilde Z$ is an orthogonal pair  with
$\dim_k \Ext^1(\widetilde Z,\widetilde Y) = 2.$ 

Let $\Cal C$ be 
the full subcategory of all representations of $Q[1,9]$ which have a submodule which is
isomorphic to a direct sum of copies of $\widetilde Y$, such that the 
factor module is isomorphic to a direct sum of copies of $\widetilde Z$. 
The category $\Cal C$ 
is a length category which is equivalent to the category of Kronecker modules.
Thus, in $\Cal C$ we find a $\mathbb P^1$-family of indecomposable objects $\widetilde X_c$, 
indexed by $c\in \mathbb P^1(k)$, such that $\widetilde X_c$ has 
a submodule of the form $\widetilde Y$, with factor module of the form $\widetilde Z$.

Since $\widetilde Y$ is sent under $\pi$ to $C_{(6,3)}$ and 
$\widetilde Z$ to $C_{(9,7,4,1)}$, the modules $\widetilde X_c$ are sent under $\pi$
to objects in $\Cal S(9)$ which are g-split extensions of $C_{(6,3)}$ by
$C_{(9,7,4,1)}:$ these are the required objects $X_c$. 
$\s$

\subsection{Wildness. An example.}
\label{sec-seven-five}

\begin{example}
  There is a 
  $\mathbb P^2$-family of indecomposable objects $X = (U,V)$ in $\Cal S(12)$
  with $|U| = [5,3],$ $|V| = [12,10,9,7,6,4,3,1]$, therefore $|U| = 8,$
  $bV = 8$ and $pX = 1$.
\end{example}

\medskip
Here is the picture: We deal with g-split extensions of $C_{(9,6,3)}$
by $C_{(12,10,7,4,1)}$:
$$
{\beginpicture
    \setcoordinatesystem units <.4cm,.4cm>
\multiput{} at -4 0   5  12 /
\plot 0 0  1 0  1 12   0 12  0 0 /
\plot 0 1  2 1  2 11  0 11 /
\plot 0 2  3 2  3 9  0 9 /
\plot 0 3  4 3  4 7  0 7 /
\plot 0 4  5 4  5 5  0 5 /
\plot 0 6  4 6 /
\plot 0 8  3 8 /
\plot 0 10 2 10 /
\multiput{$\bullet$} at 0.5 4.5  1.5 4.5  2.5 4.5  3.5 4.5  4.5 4.5 /
\plot 0.5 4.5  4.5 4.5 /
\plot 0.5 4.55  4.5 4.55 /
\plot 0.5 4.45  4.5 4.45 /

\plot -4 3  -1 3  -1 6  -4 6  /
\plot -4 2  -2 2  -2 8  -4 8 /
\plot -4 1  -3 1  -3 10 -4 10  -4 1 /
\plot -4 4  -1 4 /
\plot -4 5  -1 5 /
\plot -4 7  -2 7 /
\plot -4 9  -3 9 /
\multiput{$\bullet$} at -1.5 3.5  -2.5 3.5 -3.5 3.5 /
\plot -1.5 3.5  -3.5 3.5 /
\plot -1.5 3.55 -3.5 3.55 /
\plot -1.5 3.45 -3.5 3.45 /
\put{$c_0$} at -3.5 4.5 
\put{$c_1$} at -2.5 4.5 
\put{$c_2$} at -1.5 4.5 
\setdashes <.9mm>
\plot -1.2 4.5   0.8 4.5 /
\plot -1.2 4.55  0.8 4.55 /
\plot -1.2 4.45  0.8 4.45 /

\setdots <1mm>
\plot -3 10  0 10 /
\plot -2 8  0 8 /
\plot -1 6  0 6 /
\plot -2 3  0 3 /
\plot -2 2  0 2 /
\plot -3 1  0 1 /
\endpicture}
$$

The proof is similar to the proof for the examples in Section~\ref{sec-seven-four}.
This time, we start with the fully commutative quiver $Q[1,12]$,
and with representations $\widetilde Y,\ 
\widetilde Z$ of $Q[1,12]$ such that the pushdown functor sends $\widetilde Y$ to
$C_{(9,6,3)}$ and $\widetilde Z$ to $C_{(12,10,7,4,1)}.$
Again, we have to show that $\widetilde Y$ and $\widetilde Z$ 
is an orthogonal pair, but now we need that 
$\dim_k \Ext^1(\widetilde Z,\widetilde Y) = 3.$
$\s$

\bigskip
\section{Examples of BTh-vectors with boundary distance $1 < d < 2$.}
\label{sec-eight}

We have seen in Section~\ref{sec-six}
that all vectors in $\mathbb T(n)$ with boundary distance at least 2
are BTh-vectors. In Section \ref{sec-seven},
we have found a BTh-vector in $\mathbb T(9)$ with boundary
distance 1. Here we look at the cases $n = 7$ and $n = 8$ and exhibit some BTh-vectors
with boundary distance greater than 1 and smaller than 2.

\subsection{A Kronecker pair with uwb-vector $\frac{6|10}4$.}
\label{sec-eight-one}

For $n = 7,$ the uwb-vector $\frac{6|10}4$ is given by the Kronecker pair
shown on the left. 

$$
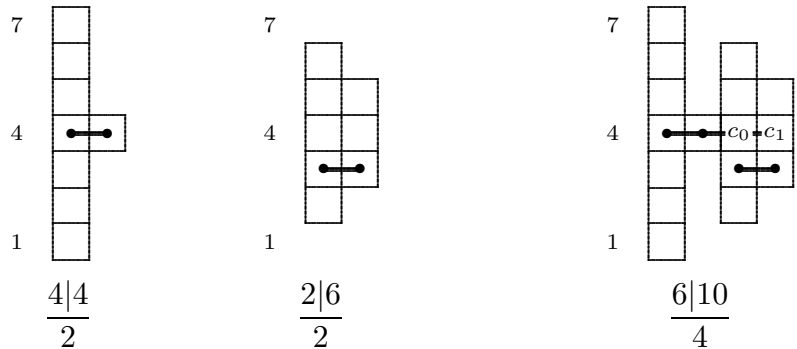
{\beginpicture
   \setcoordinatesystem units <.7cm,1cm>
\put{\beginpicture
   \setcoordinatesystem units <.4cm,0.4cm>
\setsolid
\multiput{} at 0 0  2 7 /
\plot 0 0  1 0  1 7  0 7  0 0 /
\plot 0 3  2 3  2 4  0 4 /
\plot 0 1  1 1 /
\plot 0 2  1 2 /
\plot 0 5  1 5 /
\plot 0 6  1 6 /
\multiput{$\ssize \bullet$} at 0.5 3.5  1.5 3.5 /
\plot 0.5 3.5  1.5 3.5 /
\plot 0.5 3.45  1.5 3.45 /
\plot 0.5 3.55  1.5 3.55 /
\put{$\ssize 1$} at -1 0.5
\put{$\ssize 4$} at -1 3.5
\put{$\ssize 7$} at -1 6.5
\endpicture} at 0 0
\put{$\dfrac{4|4}2$} at 0 -2

\put{\beginpicture
   \setcoordinatesystem units <.4cm,0.4cm>
\setsolid
\multiput{} at 0 0  2 5 /
\plot 0 0  1 0  1 5  0 5  0 0 /
\plot 0 1  2 1  2 4  0 4 /
\plot 0 2  2 2 /
\plot 0 3  2 3 /
\multiput{$\ssize \bullet$} at 0.5 1.5  1.5 1.5  /
\plot 0.5 1.5  1.5 1.5 /
\plot 0.5 1.45  1.5 1.45 /
\plot 0.5 1.55  1.5 1.55 /
\put{$\ssize 1$} at -1 -.5
\put{$\ssize 4$} at -1 2.5
\put{$\ssize 7$} at -1 5.5
\endpicture} at 4 0
\put{$\dfrac{2|6}2$} at 4 -2

\put{\beginpicture
   \setcoordinatesystem units <.4cm,0.4cm>
\setsolid
\multiput{} at 0 0  4 7 /
\plot 0 0  1 0  1 7  0 7  0 0 /
\plot 0 3  2 3  2 4  0 4 /
\plot 0 1  1 1 /
\plot 0 2  1 2 /
\plot 0 5  1 5 /
\plot 0 6  1 6 /
\multiput{$\ssize \bullet$} at 0.5 3.5  1.5 3.5 /
\plot 0.5 3.5  2.1 3.5 /
\plot 0.5 3.45  2.1 3.45 /
\plot 0.5 3.55  2.1 3.55 /

\plot 2.9 3.5  3.1 3.5 /
\plot 2.9 3.45  3.1 3.45 /
\plot 2.9 3.55  3.1 3.55 /

\plot 2 1  3 1  3 6  2 6  2 1 /
\plot 2 2  4 2  4 5  2 5 /
\plot 2 3  4 3 /
\plot 2 4  4 4 /
\multiput{$\ssize \bullet$} at 2.5 2.5  3.5 2.5 /
\plot 2.5 2.5  3.5 2.5 /
\plot 2.5 2.45  3.5 2.45 /
\plot 2.5 2.55  3.5 2.55 /
\put{$\ssize c_0$} at 2.5 3.5  
\put{$\ssize c_1$} at 3.5 3.5 
\put{$\ssize 1$} at -1 .5
\put{$\ssize 4$} at -1 3.5
\put{$\ssize 7$} at -1 6.5
\endpicture} at 10 0 
\put{$\dfrac{6|10}4$} at 10 -2

\endpicture}
$$

\subsection{A BTh-family in $\widetilde{\Cal S}_3(7)$ with uwb-vector $\frac{7|14}5.$}
\label{sec-eight-two}

Let $\widetilde{\Cal S}_3(7)$ be the category of pairs $(\widetilde U,\widetilde V)$ in
$\widetilde{\Cal S}(7)$ such that the height of $\widetilde U$ is bounded by 3. The category
$\widetilde{\Cal S}_3(7)$ was described completely by one of the authors in \cite{S1},
it is tame
with a unique one-parameter family of components (up to the shift). Here is a BTh-family in 
$\widetilde{\Cal S}_3(7)$.

\medskip
$$
{\beginpicture
    \setcoordinatesystem units <.4cm,.4cm>
\multiput{} at -3 0   3  7 /
\plot 0 0  1 0  1 7   0 7  0 0 /
\plot 0 1  2 1  2 5  0 5 /
\plot 0 2  3 2  3 3  0 3 /
\plot 0 4  2 4 /
\plot 0 5  1 5 /
\plot 0 6  1 6 /
\multiput{$\bullet$} at -1.3 2.3  0.7 2.3  1.7 2.3  2.7 2.3 /
\plot 0.7 2.3  2.7 2.3 /

\plot -3 0  -2 0  -2 6  -3 6  -3 0 /
\plot -3 1  -1 1  -1 4  -3 4 /
\plot -3 2  -1 2 /
\plot -3 3  -1 3 /
\plot -3 5  -2 5 /
\multiput{$\bullet$} at -1.7 2.7  -2.7 2.7 /
\plot -1.7 2.7  -2.7 2.7 /
\put{$c_0$} at -1.5 1.5 
\put{$c_1$} at 1.5 1.5
\setdashes <.9mm>
\plot -1.2 1.5  1.2 1.5 /
\plot -1.3 2.3  0.7 2.3 /

\setdots <1mm>
\plot -2 6  0 6 /
\plot -2 5  0 5 /
\plot -1 4  0 4 /
\plot -1 3  0 3 /
\plot -1 2  0 2 /
\plot -1 1  0 1 /
\plot -2 0  0 0 /
\endpicture}
$$
The uwb-vector is $\frac{7|14}5.$ 

\medskip
\begin{proposition}
  The $\widetilde{\Cal S}_3(7)$-family is 
  given by a Kronecker pair in $\mod T_2(\widetilde\Lambda).$
\end{proposition}

\begin{proof}
  Here is the Kronecker pair:
$$
X = \smallmatrix 0&1&0\cr
             0&0&0&0&0&0&0 \endsmallmatrix,
\quad
Y = \smallmatrix 2&2&2\cr
             2&4&5&4&3&2&1 \endsmallmatrix
$$
(note that $X = S_{2'}$, and, of course, $X$ does not belong to $\widetilde{\Cal S}$).
The module $Y$ is given as follows: We start with the $\mathbf E_8$-module (with arrows going down and
going left) with dimension vector
$\smallmatrix &&2\cr
             2&4&5&4&3&2&1 \endsmallmatrix$;
since the map $2 \to 5 \to 4 \to 2$ is generic, we can assume that it is the identity
map and we factor it as a sequence of identity maps 
$2 \to 2 \to 2 \to 2$ in order to obtain $Y$.
Obviously, these modules $X, Y$ are orthogonal modules with endomorphism ring $k$.

There is an exact sequence
$$
 0 \to Y \to Z \to X^2 \to 0,
$$
with $Z$ indecomposable, and we have
$$
 Z = \smallmatrix 2&4&2\cr
             2&4&5&4&3&2&1 \endsmallmatrix,
$$
where as map $4 \to 4$ we may take the identity map.
The sequence shows that  $Z$ belongs to the subcategory $\Cal C$ of all
$T_2(\widetilde\Lambda)$-modules which have a submodule which is the direct sum
of copies of $Y$ such that the quotient module is a direct sum of copies of $X$.

We may consider $Z$ and $X$ as the indecomposable injective objects
in an exact abelian subcategory $\Cal C'$ inside $\Cal C$ (then $X,Y$ are still 
the simple objects in 
$\Cal C'$, with $Y$ projective and $X$ injective).
(Actually, the
subcategory $\Cal C'$ is all of $\Cal C$; equivalently: We have $\dim \Ext^1(X,Y) = 2;$
but this does not matter at present.) 
The category $\Cal C'$ looks as follows:
$$
{\beginpicture
   \setcoordinatesystem units <.7cm,.7cm>
\multiput{} at 0 0  9 1 /
\plot 3.8 1  1 1  0 0  3.8 0 /
\plot 4 1.5  4 -.5  5 -.5  5 1.5 /
\plot 5.2 1  9 1  8 0  5.2 0 /
\put{$X$} at 9.4 1
\put{$Y$} at -.4 -.2
\put{$Z$} at 8.4  -0.2 
\put{$M_c$} at 4.3 -0.9
\multiput{$\bullet$} at 0 0  1 1  9 1  8 0  4.3 -.5 /
\setdashes <1mm> 
\plot 7 1  8 0 /
\setshadegrid span <.7mm>
\vshade 0 0 0 <z,z,,> 1 0 1 <z,z,,> 3.8  0 1 /
\vshade 4 -.5 1.5  <z,z,,> 5 -.5 1.5  /
\vshade 5.2 0 1 <z,z,,> 7 0 1  <z,z,,> 8 0 0 /
\endpicture}
$$
and 
the $\widetilde{\Cal S}_3(7)$-family $M_c$ is given by the modules $M_c$ with proper
inclusions
$$
 Y \subset M_c \subset Z.
 $$
\end{proof}

\medskip
\begin{remark}
  As we have mentioned, $X$ does not lie in 
  $\widetilde{\Cal S}(7)$, thus $\Cal C'$ does not lie inside $\widetilde{\Cal S}(7)$;
  but
  all objects in $\Cal C'$ without a direct summand of the form $X$ lie inside $\widetilde{\Cal S}(7)$;
  this is the shaded part in the picture.
\end{remark}

\subsection{A BTh-family in $\widetilde{\Cal S}(8)$ with uwb-vector $\frac{6|17}5$.}
\label{sec-eight-three}

Again, the family is given by a Kronecker pair.

$$
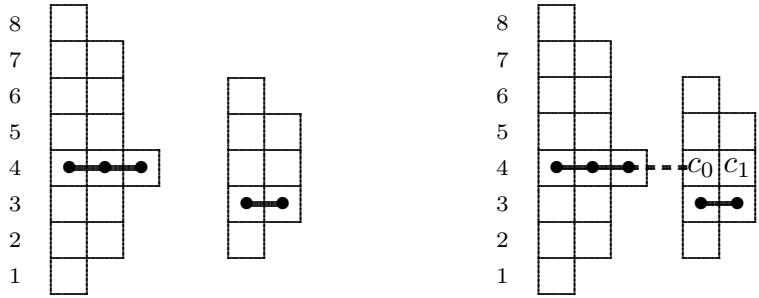
{\beginpicture
    \setcoordinatesystem units <.4cm,.4cm>
\put{\beginpicture
\multiput{} at -1 0   3  8 /
\plot 0 0  1 0  1 8  0 8  0 0 /
\plot 0 1  2 1  2 7  0 7 /
\plot 0 2  2 2 /
\plot 0 3  3 3  3 4  0 4 /
\plot 0 5  2 5 /
\plot 0 6  2 6 /
\multiput{$\bullet$} at 0.5 3.5  1.5 3.5  2.5 3.5 /
\plot 0.5 3.5  2.5 3.5 /
\plot 0.5 3.55  2.5 3.55 /
\plot 0.5 3.45  2.5 3.45 /
\put{$\ssize 1$} at -1 0.5
\put{$\ssize 2$} at -1 1.5
\put{$\ssize 3$} at -1 2.5
\put{$\ssize 4$} at -1 3.5
\put{$\ssize 5$} at -1 4.5
\put{$\ssize 6$} at -1 5.5
\put{$\ssize 7$} at -1 6.5
\put{$\ssize 8$} at -1 7.5
\endpicture} at 0 0
\put{\beginpicture
\multiput{} at 0 0   2  8 /
\plot 0 1  1 1  1 6  0 6  0 1 /
\plot 0 2  2 2  2 5  0 5 /
\plot 0 3  2 3 /
\plot 0 4  2 4 /
\multiput{$\bullet$} at 0.5 2.5  1.5 2.5  /
\plot 0.5 2.5  1.5 2.5 /
\plot 0.5 2.55  1.5 2.55 /
\plot 0.5 2.45  1.5 2.45 /
\endpicture} at 5 0

\put{\beginpicture
\multiput{} at -1 0   6  8 /
\plot 0 0  1 0  1 8   0 8  0 0 /
\plot 0 1  2 1  2 7  0 7 /
\plot 0 3  3 3  3 4  0 4 /
\plot 0 5  2 5 /
\plot 0 6  2 6 /
\plot 0 2  2 2 /
\multiput{$\bullet$} at 0.5 3.5  1.5 3.5  2.5 3.5 /
\plot 0.5 3.5  2.5 3.5 /
\plot 0.5 3.55  2.5 3.55 /
\plot 0.5 3.45  2.5 3.45 /

\plot 4 1  5 1  5 6  4 6  4 1 /
\plot 4 2  6 2  6 5  4 5 /
\plot 4 3  6 3 /
\plot 4 4  6 4 /
\multiput{$\bullet$} at 4.5 2.5  5.5 2.5 /
\plot 4.5 2.5  5.5 2.5 /
\plot 4.5 2.55  5.5 2.55 /
\plot 4.5 2.45  5.5 2.45 /

\setdashes <.9mm>
\plot 4.1 3.5   2.5 3.5 /
\plot 4.1 3.55  2.5 3.55 /
\plot 4.1 3.45  2.5 3.45 /

\put{$c_0$} at 4.5 3.5
\put{$c_1$} at 5.5 3.5
\put{$\ssize 1$} at -1 0.5
\put{$\ssize 2$} at -1 1.5
\put{$\ssize 3$} at -1 2.5
\put{$\ssize 4$} at -1 3.5
\put{$\ssize 5$} at -1 4.5
\put{$\ssize 6$} at -1 5.5
\put{$\ssize 7$} at -1 6.5
\put{$\ssize 8$} at -1 7.5
\endpicture} at 15 0
\endpicture}
$$
The uwb-vector is $\frac{6|17}5.$

\subsection{Interpolation between $\frac{6|17}5$ and $\frac{6|16}4$.}
\label{sec-eight-four}

We start with the family $D$ given in \ref{sec-eight-three} (with uwb-vector $\frac{6|17}5$) and
form $D' = H_4 D$ (with uwb-vector $\frac{6|16}4$).
$$
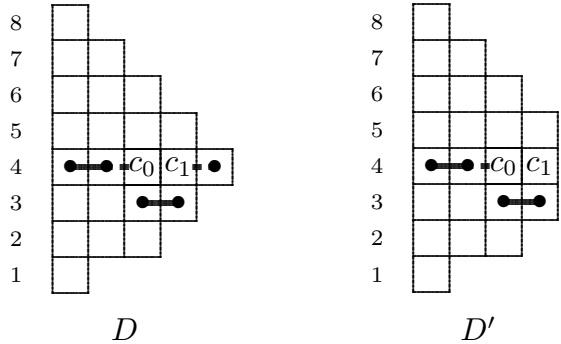
{\beginpicture
    \setcoordinatesystem units <.4cm,.4cm>
\put{\beginpicture
\multiput{} at -1 0   6  8 /
\plot 0 0  1 0  1 8   0 8  0 0 /
\plot 0 1  2 1  2 7  0 7 /
\plot 0 3  3 3  3 4  0 4 /
\plot 0 5  2 5 /
\plot 0 6  2 6 /
\plot 0 2  2 2 /
\multiput{$\bullet$} at 0.5 3.5  1.5 3.5  4.5 3.5 /
\plot 0.5 3.5  1.5 3.5 /
\plot 0.5 3.55  1.5 3.55 /
\plot 0.5 3.45  1.5 3.45 /

\plot 2 1  3 1  3 6  2 6 /
\plot 2 2  4 2  4 5  2 5 /
\plot 2 3  4 3 /
\plot 2 4  4 4 /

\plot 4 3  5 3  5 4  4 4 /
\multiput{$\bullet$} at 2.5 2.5  3.5 2.5 /
\plot 2.5 2.5  3.5 2.5 /
\plot 2.5 2.55  3.5 2.55 /
\plot 2.5 2.45  3.5 2.45 /

\setdashes <.9mm>
\plot 2.1 3.5   1.5 3.5 /
\plot 2.1 3.55  1.5 3.55 /
\plot 2.1 3.45  1.5 3.45 /

\plot 3.9 3.5   4.5 3.5 /
\plot 3.9 3.55  4.5 3.55 /
\plot 3.9 3.45  4.5 3.45 /

\put{$c_0$} at 2.5 3.5
\put{$c_1$} at 3.5 3.5
\put{$\ssize 1$} at -1 0.5
\put{$\ssize 2$} at -1 1.5
\put{$\ssize 3$} at -1 2.5
\put{$\ssize 4$} at -1 3.5
\put{$\ssize 5$} at -1 4.5
\put{$\ssize 6$} at -1 5.5
\put{$\ssize 7$} at -1 6.5
\put{$\ssize 8$} at -1 7.5
\put{$D$} at 2 -1
\endpicture} at 0 0
\put{\beginpicture
\multiput{} at -1 0   6  8 /
\plot 0 0  1 0  1 8   0 8  0 0 /
\plot 0 1  2 1  2 7  0 7 /
\plot 0 3  3 3  3 4  0 4 /
\plot 0 5  2 5 /
\plot 0 6  2 6 /
\plot 0 2  2 2 /
\multiput{$\bullet$} at 0.5 3.5  1.5 3.5  /
\plot 0.5 3.5  1.5 3.5 /
\plot 0.5 3.55  1.5 3.55 /
\plot 0.5 3.45  1.5 3.45 /

\plot 2 1  3 1  3 6  2 6 /
\plot 2 2  4 2  4 5  2 5 /
\plot 2 3  4 3 /
\plot 2 4  4 4 /

\multiput{$\bullet$} at 2.5 2.5  3.5 2.5 /
\plot 2.5 2.5  3.5 2.5 /
\plot 2.5 2.55  3.5 2.55 /
\plot 2.5 2.45  3.5 2.45 /

\setdashes <.9mm>
\plot 2.1 3.5   1.5 3.5 /
\plot 2.1 3.55  1.5 3.55 /
\plot 2.1 3.45  1.5 3.45 /

\put{$c_0$} at 2.5 3.5
\put{$c_1$} at 3.5 3.5
\put{$\ssize 1$} at -1 0.5
\put{$\ssize 2$} at -1 1.5
\put{$\ssize 3$} at -1 2.5
\put{$\ssize 4$} at -1 3.5
\put{$\ssize 5$} at -1 4.5
\put{$\ssize 6$} at -1 5.5
\put{$\ssize 7$} at -1 6.5
\put{$\ssize 8$} at -1 7.5
\put{$D'$} at 1.8 -1
\endpicture} at 10 0
\endpicture}
$$
We should mention that we may obtain $D'$ also from the BTh-family constructed in
Section~\ref{sec-eight-one} 
by applying duality and the rotation $\rho = \tau_8^2$ in $\mathbb T(8).$ 

\medskip
In order to obtain a BTh-family with uwb-vector $d \frac{6|17}5 + e\frac{6|16}4,$
we form the Jordan extension $D[d+e]$ and factor out $e$ copies of $S_4$ (similar to our
procedure in Sections \ref{sec-six-three} and \ref{sec-six-four}). 

$$
{\beginpicture
   \setcoordinatesystem units <0.404cm,0.7cm>
\multiput{} at -8 0  2 4 /
\setdots <1mm>
\plot -4 4  -8 0  2 0 /
\plot -7 1  -6 0  -5 1  -4 0  -3 1  -2 0  -1 1  0 0  1 1  2 0 /
\plot -6 2  -5 1  -4 2  -3 1  -2 2  -1 1  0 2  1 1  2 2 /
\plot -5 3  -4 2  -3 3  -2 2  -1 3  0 2  1 3  2 2 /
\plot -7 1  2 1 /
\plot -6 2  2 2 /
\setshadegrid span <.5mm>
\vshade -2 2 2 <z,z,,>  0 2 4 <z,z,,> 2 2 2   /
\multiput{$\bullet$} at -6 0  /
\multiput{$\ssize\blacksquare$} at  1.5 1.5  0 1.2  2.5 2.5 /
\setdashes <1mm>
\plot 1.5 1.5  2.5 2.5 /
\setsolid
\plot 0 1.2   1.5 1.5 /
\setdots <0.5mm>
\plot -6 0  0 1.2 /
\put{$\frac{0|1}1$} at -6 -.5
\put{$\frac{6|17}5$} at 0 .6
\put{$\frac{6|16}4$} at 1.8 1
\put{$\frac{10|16}4$} at 3.2 2.8
\endpicture}
$$

\subsection{BTh-vectors with arbitrary boundary distance between 1 and 2.}
\label{sec-eight-five}

The following result provides BTh-vectors with arbitrary boundary distance between 1 and 2.

\medskip
\begin{theorem}
  \label{theoremten}
  Let $n = 10.$ 
  Any rational pr-vector $(p,r)$ in $\mathbb T(n)$ which belongs to one of the intervals
  between $(1,4)$ and $(1,5),$ or between
  $(1,4)$ and $(2,4),$ or between $(1,5)$ and $(2,4)$ is a BTh-vector.

  Let $n\ge 11$. 
  Any rational pr-vector $(p,r)$ in $\mathbb T(n)$ which belongs to the
  convex hull of $(1,4),\ (1,n-5),\ (2,n-6),\ (2,4),$ is a BTh-vector.
\end{theorem}

\medskip
Here is the case $n = 10.$
$$  
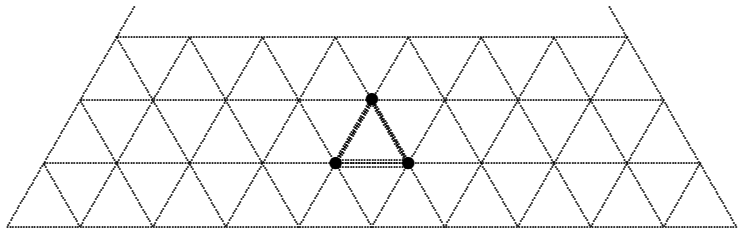
{\beginpicture
   \setcoordinatesystem units <.404cm,.7cm>
\multiput{} at -10 0  10 3 /
\multiput{$\bullet$} at -1 1  1 1  0 2  /
\setdots <.3mm>
\plot -7 3  -10 0  10 0  7 3 /
\plot 9 1  -9 1  -8 0  -5 3  -2 0  1 3  4 0  7 3 /
\plot 8 2  -8 2  -6 0  -3 3  0 0  3 3  6 0  8 2 /
\plot 7 3  -7 3  -4 0  -1 3  2 0  5 3  8 0  9 1 /
\plot 7 3  6.5 3.5 /
\plot -7 3  -6.5 3.5 /
\setshadegrid span <.3mm>
\plot -1 1   0 2   1 1   -1 1 /
\plot -1 .95   0 1.95  1 0.95  -1 0.95 /
\plot -1 1.05   0 2.05  1 1.05  -1 1.05 /
\endpicture}
$$

And the case $n = 12.$
$$  
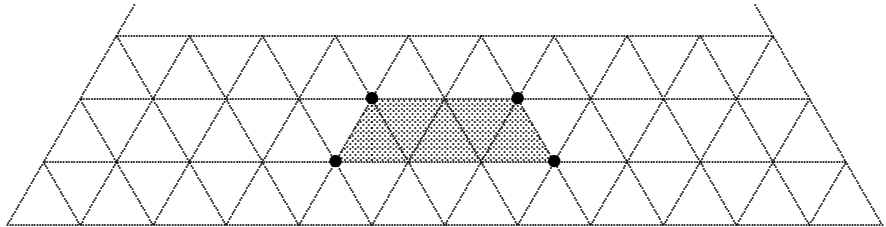
{\beginpicture
   \setcoordinatesystem units <.404cm,.7cm>
\multiput{} at -10 0  10 3 /
\multiput{$\bullet$} at -1 1  5 1  0 2   4 2 /
\setdots <.3mm>
\plot -7 3  -10 0  14 0  11 3 /
\plot 9 1  -9 1  -8 0  -5 3  -2 0  1 3  4 0  7 3 /
\plot 8 2  -8 2  -6 0  -3 3  0 0  3 3  6 0  8 2 /
\plot 7 3  -7 3  -4 0  -1 3  2 0  5 3  8 0  9 1 /
\plot 11 3  10.5 3.5 /
\plot -7 3  -6.5 3.5 /
\plot 10 0  7 3  11 3 /
\plot 8 2  9 3  12 0  13 1  9 1  11 3 /
\plot 8 2  12 2  10 0  /  
\setshadegrid span <.3mm>
\vshade -1 1 1 <z,z,,> 0 1 2 <z,z,,> 4 1 2 <z,z,,> 5 1 1 /
\endpicture}
$$

\smallskip
\begin{proof}[Proof of Theorem~\ref{theoremten}]
  We use expansion and coexpansion as in Section~\ref{sec-six-three}
  but now based on the Kronecker families
  in Sections~\ref{sec-seven-four} and \ref{sec-eight-one}.
  Our goal is to construct the convex regions of
  BTh-vectors in $\mathbb T(n)$ for $n\geq 10$ claimed to exist in Theorem~\ref{theoremten}
  However, the families mentioned are not st-solid, so the regions obtained
  will be somewhat smaller, compared to Corollary~\ref{cor-six-three}.

  More precisely, both families $M=M_c$ with $M=(U,V)\in\widetilde{\Cal S}$
  are {\it almost} 44-solid in the sense that $V$ is a direct sum of indecomposables of the form
  $[i,j]$ with $i\leq4\leq j$ and $U_i=0$ holds for $i>4$. But the condition $U_i=V_i$ is only
  satisfied for $i\leq 1$ and $i\leq 2$ for the families in Sections~\ref{sec-seven-four}
  and \ref{sec-eight-one}, respectively.

  As a consequence, expansions as in Section~\ref{sec-six-three} are possible,
  but for coexpansions the range of the
  parameter $e$ is limited.

\medskip  
\begin{lemma}
  \label{lem-eight-five}
  \begin{itemize}[leftmargin=3em]
    \item[{\rm (a)}]
    Let $M = M_c$ be the Kronecker family in $\widetilde{\Cal S}(7)$
    presented in Section~\ref{sec-eight-one}, it satisfies $pM=\frac32$, $rM=\frac52$ and $bM=4$.
    Let $0 \le e\le \frac12\ell\cdot bM$ and $0\le f \le \ell\cdot bM.$
    Then 
    $(pM+\frac e{\ell\cdot bM},rM+\frac f{\ell\cdot bM})$ is a BTh-vector in $\mathbb T(9)$
    and $(pM+\frac e{\ell\cdot bM},rM)$ 
    and $(pM,rM+\frac f{\ell\cdot bM})$ are BTh-vectors in $\mathbb T(8).$
    \item[{\rm (b)}]
    Let $M = M_c$ be the Kronecker family in $\widetilde{\Cal S}(9)$
    presented in Section~\ref{sec-seven-four} with $pM=1$, $rM=4$ and $bM=6$. 
    Let $0 \le e\le \frac12\ell\cdot bM$ and $0\le f \le \ell\cdot bM.$
    Then 
    $(pM+\frac e{\ell\cdot bM},rM+\frac f{\ell\cdot bM})$ is a BTh-vector in $\mathbb T(11)$
    and $(pM+\frac e{\ell\cdot bM},rM)$ 
    and $(pM,rM+\frac f{\ell\cdot bM})$ are BTh-vectors in $\mathbb T(10).$
  \end{itemize}
\end{lemma}

\medskip
The pr-vectors which we obtain in this way are the vectors in the following region
which we picture in the case $n = 12.$  It is the union of the convex region bounded by
$(\frac32,\frac52)$, $(2,\frac52)$, $(2,\frac72)$, $(\frac32,\frac72)$, given by Part (a),
and the convex region bounded by $(1,4)$, $(\frac32,4)$, $(\frac32,5)$, $(1,5)$,
given by Part (b).

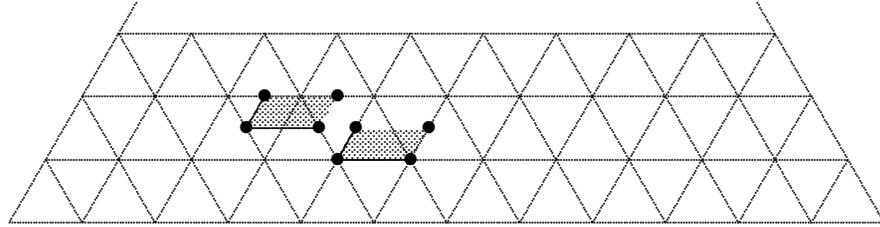
$$  
{\beginpicture
   \setcoordinatesystem units <.404cm,.7cm>
\multiput{} at -10 0  10 3 /
\multiput{$\bullet$} at -3.5 1.5  -1.5 1.5  -1 2  -3 2 /
\multiput{$\bullet$} at -1 1  -.5 1.5  1.5 1.5  1 1 /
\setdots <.3mm>
\plot -7 3  -10 0  14 0  11 3 /
\plot 9 1  -9 1  -8 0  -5 3  -2 0  1 3  4 0  7 3 /
\plot 8 2  -8 2  -6 0  -3 3  0 0  3 3  6 0  8 2 /
\plot 7 3  -7 3  -4 0  -1 3  2 0  5 3  8 0  9 1 /
\plot 11 3  10.5 3.5 /
\plot -7 3  -6.5 3.5 /
\plot 10 0  7 3  11 3 /
\plot 8 2  9 3  12 0  13 1  9 1  11 3 /
\plot 8 2  12 2  10 0  /  
\setshadegrid span <.3mm>
\vshade -3.5 1.5 1.5  <z,z,,> -3 1.5 2  <z,z,,> -1.5 1.5 2 <z,z,,> -1 2 2 /
\vshade -1 1 1 <z,z,,> -.5 1 1.5   <z,z,,> 1 1 1.5 <z,z,,> 1.5 1.5 1.5  /
\setsolid
\plot -3 2  -3.5 1.5  -1.5 1.5 /
\plot -.5 1.5  -1 1  1 1 /
\endpicture}
$$

We sketch here the objects for the interpolation, each position is one of the marked vertices,
the order is from left to right.  First, the family based on Section~\ref{sec-eight-one}.

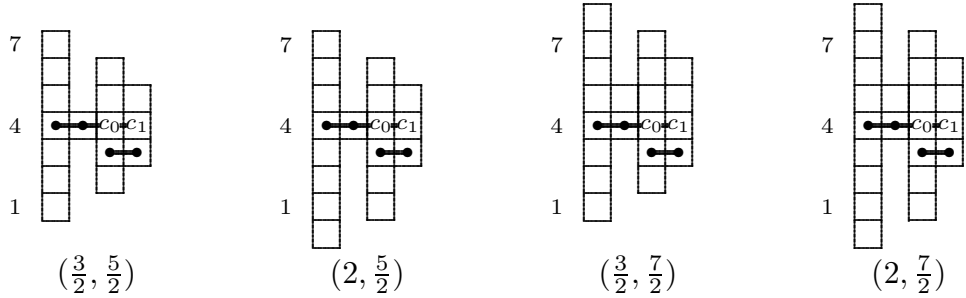
$$\beginpicture   \setcoordinatesystem units <.3cm,0.3cm>
\put{\beginpicture
\setsolid
\multiput{} at 0 0  4 7 /
\plot 0 0  1 0  1 7  0 7  0 0 /
\plot 0 3  2 3  2 4  0 4 /
\plot 0 1  1 1 /
\plot 0 2  1 2 /
\plot 0 5  1 5 /
\plot 0 6  1 6 /
\multiput{$\ssize \bullet$} at 0.5 3.5  1.5 3.5 /
\plot 0.5 3.5  2.1 3.5 /
\plot 0.5 3.45  2.1 3.45 /
\plot 0.5 3.55  2.1 3.55 /

\plot 2.9 3.5  3.1 3.5 /
\plot 2.9 3.45  3.1 3.45 /
\plot 2.9 3.55  3.1 3.55 /

\plot 2 1  3 1  3 6  2 6  2 1 /
\plot 2 2  4 2  4 5  2 5 /
\plot 2 3  4 3 /
\plot 2 4  4 4 /
\multiput{$\ssize \bullet$} at 2.5 2.5  3.5 2.5 /
\plot 2.5 2.5  3.5 2.5 /
\plot 2.5 2.45  3.5 2.45 /
\plot 2.5 2.55  3.5 2.55 /
\put{$\ssize c_0$} at 2.5 3.5  
\put{$\ssize c_1$} at 3.5 3.5 
\put{$\ssize 1$} at -1 .5
\put{$\ssize 4$} at -1 3.5
\put{$\ssize 7$} at -1 6.5
\put{$(\frac32,\frac52)$} at 2 -2
\endpicture} at 0 0 

\put{\beginpicture
\setsolid
\multiput{} at 0 0  4 7 /
\plot 0 -1  1 -1  1 7  0 7  0 -1 /
\plot 0 3  2 3  2 4  0 4 /
\plot 0 0  1 0 /
\plot 0 1  1 1 /
\plot 0 2  1 2 /
\plot 0 5  1 5 /
\plot 0 6  1 6 /
\multiput{$\ssize \bullet$} at 0.5 3.5  1.5 3.5 /
\plot 0.5 3.5  2.1 3.5 /
\plot 0.5 3.45  2.1 3.45 /
\plot 0.5 3.55  2.1 3.55 /

\plot 2.9 3.5  3.1 3.5 /
\plot 2.9 3.45  3.1 3.45 /
\plot 2.9 3.55  3.1 3.55 /

\plot 2 0  3 0  3 6  2 6  2 0 /
\plot 2 2  4 2  4 5  2 5 /
\plot 2 3  4 3 /
\plot 2 4  4 4 /
\plot 2 1  3 1 /
\multiput{$\ssize \bullet$} at 2.5 2.5  3.5 2.5 /
\plot 2.5 2.5  3.5 2.5 /
\plot 2.5 2.45  3.5 2.45 /
\plot 2.5 2.55  3.5 2.55 /
\put{$\ssize c_0$} at 2.5 3.5  
\put{$\ssize c_1$} at 3.5 3.5 
\put{$\ssize 1$} at -1 .5
\put{$\ssize 4$} at -1 3.5
\put{$\ssize 7$} at -1 6.5
\put{$(2,\frac52)$} at 2 -2
\endpicture} at 10 0 

\put{\beginpicture
\setsolid
\multiput{} at 0 0  4 7 /
\plot 0 0  1 0  1 8  0 8  0 0 /
\plot 0 3  2 3  2 5  0 5 /
\plot 0 1  1 1 /
\plot 0 2  1 2 /
\plot 0 4  2 4 /
\plot 0 6  1 6 /
\plot 0 7  1 7 /
\multiput{$\ssize \bullet$} at 0.5 3.5  1.5 3.5 /
\plot 0.5 3.5  2.1 3.5 /
\plot 0.5 3.45  2.1 3.45 /
\plot 0.5 3.55  2.1 3.55 /

\plot 2.9 3.5  3.1 3.5 /
\plot 2.9 3.45  3.1 3.45 /
\plot 2.9 3.55  3.1 3.55 /

\plot 2 1  3 1  3 7  2 7  2 1 /
\plot 2 2  4 2  4 6  2 6 /
\plot 2 3  4 3 /
\plot 2 4  4 4 /
\plot 2 5  4 5 /
\multiput{$\ssize \bullet$} at 2.5 2.5  3.5 2.5 /
\plot 2.5 2.5  3.5 2.5 /
\plot 2.5 2.45  3.5 2.45 /
\plot 2.5 2.55  3.5 2.55 /
\put{$\ssize c_0$} at 2.5 3.5  
\put{$\ssize c_1$} at 3.5 3.5 
\put{$\ssize 1$} at -1 .5
\put{$\ssize 4$} at -1 3.5
\put{$\ssize 7$} at -1 6.5
\put{$(\frac32,\frac72)$} at 2 -2 
\endpicture} at 20 0 

\put{\beginpicture
\setsolid
\multiput{} at 0 0  4 7 /
\plot 0 -1  1 -1  1 8  0 8  0 -1 /
\plot 0 3  2 3  2 5  0 5 /
\plot 0 0  1 0 /
\plot 0 1  1 1 /
\plot 0 2  1 2 /
\plot 0 4  2 4 /
\plot 0 7  1 7 /
\plot 0 6  1 6 /
\multiput{$\ssize \bullet$} at 0.5 3.5  1.5 3.5 /
\plot 0.5 3.5  2.1 3.5 /
\plot 0.5 3.45  2.1 3.45 /
\plot 0.5 3.55  2.1 3.55 /

\plot 2.9 3.5  3.1 3.5 /
\plot 2.9 3.45  3.1 3.45 /
\plot 2.9 3.55  3.1 3.55 /

\plot 2 0  3 0  3 7  2 7  2 0 /
\plot 2 2  4 2  4 6  2 6 /
\plot 2 5  4 5 /
\plot 2 3  4 3 /
\plot 2 4  4 4 /
\plot 2 1  3 1 /
\multiput{$\ssize \bullet$} at 2.5 2.5  3.5 2.5 /
\plot 2.5 2.5  3.5 2.5 /
\plot 2.5 2.45  3.5 2.45 /
\plot 2.5 2.55  3.5 2.55 /
\put{$\ssize c_0$} at 2.5 3.5  
\put{$\ssize c_1$} at 3.5 3.5 
\put{$\ssize 1$} at -1 .5
\put{$\ssize 4$} at -1 3.5
\put{$\ssize 7$} at -1 6.5
\put{$(2,\frac72)$} at 2 -2
\endpicture} at 30 0 
\endpicture
$$

Next, the objects obtained from the family in Section~\ref{sec-seven-five}.
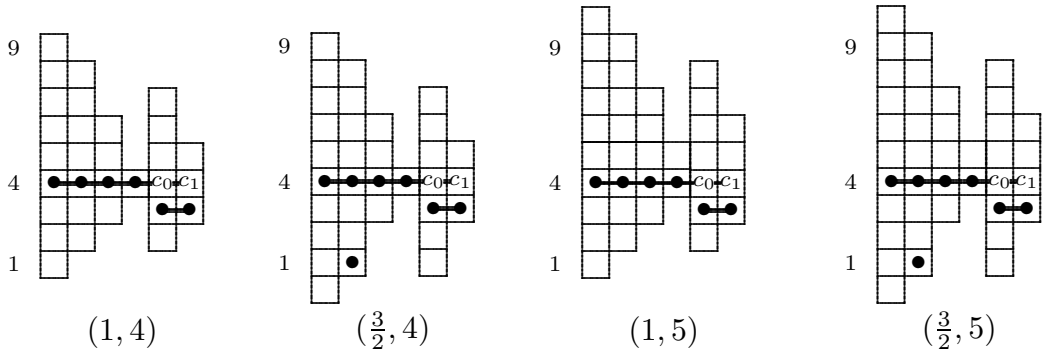
$$
\beginpicture
    \setcoordinatesystem units <.3cm,.3cm>
\put{\beginpicture
\multiput{} at -1 -2   6 10 /
\plot 0 0  1 0  1 9   0 9  0 0 /
\plot 0 1  2 1  2 8  0 8 /
\plot 0 2  3 2  3 6  0 6 /
\plot 0 3  4 3  4 4  0 4 /
\plot 0 5  3 5 /
\plot 0 7  2 7 /
\multiput{$\bullet$} at 0.5 3.5  1.5 3.5  2.5 3.5  3.5 3.5 /
\plot 0.5 3.5   4.1 3.5 /
\plot 0.5 3.55  4.1 3.55 /
\plot 0.5 3.45  4.1 3.45 /
\plot 4 2  6 2  6 5  5 5  4 5 /
\plot 4 1  5 1  5 7  4 7  4 1  /
\plot 4 3  6 3 /
\plot 4 4  6 4 /
\plot 4 6  5 6 /
\multiput{$\bullet$} at 5.5 2.5  4.5 2.5 /
\plot 4.5 2.5  5.5 2.5 /
\plot 4.5 2.55 5.5 2.55 /
\plot 4.5 2.45 5.5 2.45 /
\put{$\scriptstyle c_0$} at 4.5 3.5 
\put{$\scriptstyle c_1$} at 5.5 3.5 
\setsolid
\plot 4.9 3.5  5.1 3.5 /
\plot 4.9 3.55 5.1 3.55 /
\plot 4.9 3.45 5.1 3.45 /
\put{$\scriptstyle 1$} at -1 .5
\put{$\scriptstyle 4$} at -1 3.5
\put{$\scriptstyle 9$} at -1 8.5
\put{$(1,4)$} at 3 -2
\endpicture} at 0 0

\put{\beginpicture
\multiput{} at -1 -2  6 10 /
\plot 0 -1  1 -1  1 9   0 9  0 -1 /
\plot 0 0  2 0  2 8  0 8 /
\plot 0 1  2 1 /
\plot 0 2  3 2  3 6  0 6 /
\plot 0 3  4 3  4 4  0 4 /
\plot 0 5  3 5 /
\plot 0 7  2 7 /
\multiput{$\bullet$} at 0.5 3.5  1.5 3.5  2.5 3.5  3.5 3.5  1.5 .5 /
\plot 0.5 3.5   4.1 3.5 /
\plot 0.5 3.55  4.1 3.55 /
\plot 0.5 3.45  4.1 3.45 /
\plot 4 2  6 2  6 5  5 5  4 5 /
\plot 4 0  5 0  5 7  4 7  4 0  /
\plot 4 1  5 1 /
\plot 4 3  6 3 /
\plot 4 4  6 4 /
\plot 4 6  5 6 /
\multiput{$\bullet$} at 5.5 2.5  4.5 2.5 /
\plot 4.5 2.5  5.5 2.5 /
\plot 4.5 2.55 5.5 2.55 /
\plot 4.5 2.45 5.5 2.45 /
\put{$\scriptstyle c_0$} at 4.5 3.5 
\put{$\scriptstyle c_1$} at 5.5 3.5 
\setsolid
\plot 4.9 3.5  5.1 3.5 /
\plot 4.9 3.55 5.1 3.55 /
\plot 4.9 3.45 5.1 3.45 /
\put{$\scriptstyle 1$} at -1 .5
\put{$\scriptstyle 4$} at -1 3.5
\put{$\scriptstyle 9$} at -1 8.5
\put{$(\frac32,4)$} at 3 -2
\endpicture} at 10 0

\put{\beginpicture
\multiput{} at -1 -2  6 10 /
\plot 0 0  1 0  1 10   0 10  0 0 /
\plot 0 1  2 1  2 9  0 9 /
\plot 0 2  3 2  3 7  0 7 /
\plot 0 3  4 3  4 5  0 5 /
\plot 0 5  3 5 /
\plot 0 6  3 6 /
\plot 0 8  2 8 /
\plot 0 4  4 4 /
\multiput{$\bullet$} at 0.5 3.5  1.5 3.5  2.5 3.5  3.5 3.5 /
\plot 0.5 3.5   4.1 3.5 /
\plot 0.5 3.55  4.1 3.55 /
\plot 0.5 3.45  4.1 3.45 /
\plot 4 2  6 2  6 6  4 6 /
\plot 4 1  5 1  5 8  4 8  4 1  /
\plot 4 3  6 3 /
\plot 4 4  6 4 /
\plot 4 5  6 5 /
\plot 4 7  5 7 /
\multiput{$\bullet$} at 5.5 2.5  4.5 2.5 /
\plot 4.5 2.5  5.5 2.5 /
\plot 4.5 2.55 5.5 2.55 /
\plot 4.5 2.45 5.5 2.45 /
\put{$\scriptstyle c_0$} at 4.5 3.5 
\put{$\scriptstyle c_1$} at 5.5 3.5 
\setsolid
\plot 4.9 3.5  5.1 3.5 /
\plot 4.9 3.55 5.1 3.55 /
\plot 4.9 3.45 5.1 3.45 /
\put{$\scriptstyle 1$} at -1 .5
\put{$\scriptstyle 4$} at -1 3.5
\put{$\scriptstyle 9$} at -1 8.5
\put{$(1,5)$} at 3 -2
\endpicture} at 20 0

\put{\beginpicture
\multiput{} at -3 -1   4  10 /
\plot 0 -1  1 -1  1 10   0 10  0 -1 /
\plot 0 0  2 0  2 9  0 9 /
\plot 0 1  2 1 /
\plot 0 2  3 2  3 7  0 7 /
\plot 0 3  4 3  4 5  0 5 /
\plot 0 6  3 6 /
\plot 0 8  2 8 /
\plot 0 4  4 4 /
\multiput{$\bullet$} at 0.5 3.5  1.5 3.5  2.5 3.5  3.5 3.5  1.5 .5 /
\plot 0.5 3.5   4.1 3.5 /
\plot 0.5 3.55  4.1 3.55 /
\plot 0.5 3.45  4.1 3.45 /
\plot 4 2  6 2  6 6  4 6 /
\plot 4 0  5 0  5 8  4 8  4 0  /
\plot 4 1  5 1 /
\plot 4 3  6 3 /
\plot 4 4  6 4 /
\plot 4 5  6 5 /
\plot 4 7  5 7 /
\multiput{$\bullet$} at 5.5 2.5  4.5 2.5 /
\plot 4.5 2.5  5.5 2.5 /
\plot 4.5 2.55 5.5 2.55 /
\plot 4.5 2.45 5.5 2.45 /
\put{$\scriptstyle c_0$} at 4.5 3.5 
\put{$\scriptstyle c_1$} at 5.5 3.5 
\setsolid
\plot 4.9 3.5  5.1 3.5 /
\plot 4.9 3.55 5.1 3.55 /
\plot 4.9 3.45 5.1 3.45 /
\put{$\scriptstyle 1$} at -1 .5
\put{$\scriptstyle 4$} at -1 3.5
\put{$\scriptstyle 9$} at -1 8.5
\put{$(\frac32,5)$} at 3 -2
\endpicture} at 30 0
\endpicture
$$

In order to complete the proof of Theorem~\ref{theoremten},
it suffices to observe that the region
of BTh-vectors claimed to exist can be covered by the regions
given by Lemma~\ref{lem-eight-five},
and their images under reflections on vertical lines through vertices
$(0,\frac\ell2)\in\mathbb T(n)$
for suitable values $\ell \le n$.  Each such reflection can be realized by
applying the functor $\tau_\ell^2\D$ to objects in $\Cal S(\ell)$. 
\end{proof}

\medskip
\begin{corollary}
  Let $n\ge 10$ and let $d$ be a rational number with
  $1\le d\le \frac n3.$ There exists an
  indecomposable object in $\Cal S(n)$ with boundary distance $d.$ $\s$
\end{corollary}

\vfill\eject
\centerline{\Gross Third part: Half-line and triangle support.}
\addcontentsline{toc}{part}{Third part: Half-line support.  Triangle support.}

\section{Rays and central half-lines: the plus-construction.}
\label{sec-nine}

We consider now Auslander-Reiten components of $\Cal S(n)$. 
The main result of this Section is Theorem \ref{theoremfive}.

\subsection{Rays, corays, and the quasi-length of indecomposable objects.}
\label{sec-nine-one}

We assume that $n\ge 6.$ 
We say that a vertex $z$ of a tube is {\it ray-simple} provided only one arrow ends in
$z$. Dually, a vertex $z$ of a tube is {\it coray-simple} provided only one arrow starts in
$z$. 
Since $n \ge 6,$ there is a (unique) infinite sectional path which starts in a 
given ray-simple vertex, we  will call it a
{\it (categorical) ray.} Dually, 
there is a (unique) infinite sectional path which ends in a 
given coray-simple vertex, we  will call it a {\it (categorical) coray.} 

In the stable components, the ray-simple vertices are just the vertices at the mouth.
In the principal component $\Cal P(n)$, the projective vertices are ray-simple, and 
there are four ray-simple vertices which are not projective: 
the immediate predecessors
of the projective vertices, as well as the isomorphism classes of $S = (0,[1])$ and of 
$([1],[1]).$
Dually, in $\Cal P(n)$, the projective vertices are coray-simple, and 
there are four coray-simple vertices which are not projective: 
the immediate successors
of the projective vertices, as well as again $S = (0,[1])$ and 
$([1],[1]).$
Here are two sketches of the mouth of the principle component $\Cal P(n)$
of $\Cal S(n)$. On the left, the ray-simple vertices are marked by bullets, on the right,
the coray-simple objects are marked by bullets:
$$
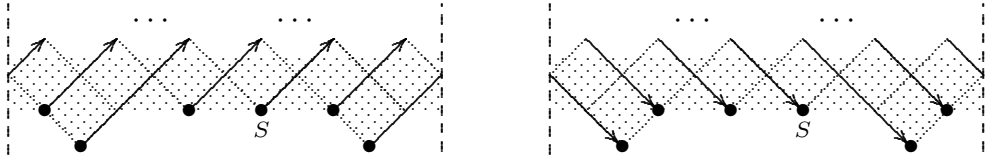
{\beginpicture  
    \setcoordinatesystem units <.4cm,.4cm>
\put{\beginpicture
\multiput{} at 0 0  12 4  /
\setdashes <1mm>
\plot 0 -.2  0 4 /
\plot 12 -.1  12 4 /
\multiput{$\bullet$} at  1 1  2 0  5 1  7 1  9 1  10 0  /
\setdots <.5mm>
\plot 0 2  2 0  5 3  7 1  9 3  11 1 /
\plot 1 3  3 1  /
\plot 3 3  5 1 /
\plot 7 3  10 0 /
\plot 9 3  11 1 /
\plot 11 3  12 2 /
\setsolid
\arr{0 2}{1 3}
\arr{1 1}{3 3}
\arr{2 0}{5 3}
\arr{5 1}{7 3}
\arr{7 1}{9 3}
\arr{9 1}{11 3}
\plot 10 0  12 2 /
\put{$\ss S$} at 7 .5
\setshadegrid span <.7mm>
\vshade 0 1 2.5  <,z,,> 1 1 2.5  <z,z,,> 2 0  2.5  <z,z,,> 3 1 2.5  <z,z,,> 8 1 2.5 
    <z,z,,> 9 1 2.5   <z,z,,> 10 0 2.5 <z,z,,> 11 1 2.5   <z,z,,> 12 1 2.5  /
\multiput{$\cdots$} at 4 3.5  8 3.5 /
\endpicture} at 0 0
\put{\beginpicture
\multiput{} at 0 0  12 4  /
\setdashes <1mm>
\plot 0 -.2  0 4 /
\plot 12 -.1  12 4 /
\multiput{$\bullet$} at  3 1  2 0  5 1  7 1  11 1  10 0  /
\setdots <.5mm>
\plot 0 2  2 0  5 3  7 1  9 3  11 1 /
\plot 0 2  1 3 /
\plot 1 1  3 3 /
\plot 5 1  7 3 /
\plot 9 1  11 3 /
\plot 10 0  12 2 /
\setsolid
\arr{0 2}{1.9 0.1}
\arr{1 3}{2.9 1.1}
\arr{3 3}{4.9 1.1}
\arr{5 3}{6.9 1.1}
\arr{7 3}{9.9 0.1}
\arr{9 3}{10.9 1.1}
\plot 11 3  12 2 /
\put{$\ss S$} at 7 .5
\setshadegrid span <.7mm>
\vshade 0 1 2.5  <,z,,> 1 1 2.5  <z,z,,> 2 0  2.5  <z,z,,> 3 1 2.5  <z,z,,> 8 1 2.5 
    <z,z,,> 9 1 2.5   <z,z,,> 10 0 2.5 <z,z,,> 11 1 2.5   <z,z,,> 12 1 2.5  /
\multiput{$\cdots$} at 4 3.5  8 3.5 /
\endpicture} at 15 0
\endpicture}
$$
(on the left, for every ray-simple object $z$, we have drawn a solid arrow which indicates
the beginning of the ray starting at $z$; similarly, on the right we indicate the corays 
by arrows). 

As we have mentioned, 
the two projective vertices $z$ in $\Cal P(n)$ are ray-simple, and the ray starting in the
projective vertex $z$ may be written in the form
$$
 z = z[0] \to z[1] \to \cdots,
$$
with vertices $z[\ell]$, where $\ell\in \mathbb N_0.$ 
If $z$ is a ray-simple and non-projective vertex in the Auslander-Reiten quiver of $\Cal S(n)$,
the ray starting with $z$ may be written in the form 
$$
 z = z[1] \to z[2] \to \cdots,
$$
with vertices $z[\ell]$, where $\ell\in \mathbb N_1.$ 
Always, $\ell$ will be called the {\it quasi-length} of $z[\ell]$.
(Note that this corresponds to the established use of the word 
dealing with stable tubes, or more generally, with stable components of
tree type $\mathbf A_\infty$: If the indecomposable object $X$ is not projective,
then its quasi-length is just the quasi-length of the object $[X]$ in the
stable Auslander-Reiten quiver; if $X$ is projective, then $[X]$ vanishes
in the stable Auslander-Reiten quiver and its quasi-length is set to be zero.)

\medskip
{\bf The plus-construction.} 
For any vertex $x = z[\ell]$, let $x^+ = z[\ell+6]$. If $X$ is an indecomposable object,
let $X^+$ be a representative in the isomorphism class $[X]^+.$ 
	\smallskip
        
If $\Cal C$ is a stable component, then any plus-orbit of
non-central objects in $\Cal C$ has a representative which has quasi-length at most 5,
thus the number of plus-orbits is at most $30.$
Any plus-orbit of non-central objects in $\Cal P(n)$ has a representative 
which has quasi-length $0\le \ell\le 5$ or is equal to $S[6],$
thus the number of plus-orbits is at most $33$, see the diagram in Section~\ref{sec-nine-five}.

\subsection{The weight of a component.}
\label{sec-nine-two}

Let $n\ge 6$. For any Auslander-Reiten component $\Cal C$, 
we are going to define the {\it weight} $\beta \Cal C$ of $\Cal C$ as follows.
If $X$ belongs to a stable tube $\Cal C$, let $\beta \Cal C$ be the minimum of the values
$bX+b(\tau X)$ where $X$ is an object in $\Cal C$; obviously,  
$\beta \Cal C = bZ+b(\tau Z)$, where $Z$ is a ray-simple object in $\Cal C$
(according to Theorem \ref{theoremthree},
all the numbers $bZ+b(\tau Z)$ with $Z$ ray-simple in $\Cal C$ 
are equal). For the principal component $\Cal P(n)$, we set $\beta \Cal P(n) = 1$.
A characterization of $\beta \Cal C$ using central objects of quasi-length 6 will be given 
in Section~\ref{sec-nine-six}.

\subsection{The main result.}
\label{sec-nine-three}

\begin{theorem}
  Let $X$ be indecomposable and $\Cal C$ the component it belongs to.
  Then
$$
  \dfrac{u|w}b X^+ =  \dfrac{u|w}b X + \beta \Cal C\cdot \dfrac{n|n}3.
  $$
\end{theorem}

\medskip
The proof will be given in Sections~\ref{sec-nine-eight} and \ref{sec-nine-ten},
first for $\Cal C$ being stable, then for $\Cal C = \Cal P(n).$ 

\medskip 
We will use additive functions on translation quivers. The functions $u$ and $w$
are additive on the Auslander-Reiten quiver. Also, $b$ is nearly always additive on Auslander-Reiten
sequences, the only exception is the Auslander-Reiten sequence ending in $S$, thus an
Auslander-Reiten sequence in the principal component $\Cal P(n)$, see the following lemma. 

\subsection{The additivity of the function $b$.}
\label{sec-nine-four}

\begin{lemma}
  The function
  $b$ is additive on all Auslander-Reiten sequences but one, the
  only exception is the Auslander-Reiten sequence ending in $(0,[1])$.
\end{lemma}

\begin{proof}
  The function $b$ is nothing else
  but $\dim_k\Hom((0,[1]),-)$. For any $A$-module $M$, where $A$ is an artin algebra, 
  the functor $\Hom(M,-)$ is exact on an Auslander-Reiten sequence $0 \to X \to Y \to Z \to 0$
  provided $Z$ is not a direct summand of $M$. Thus, $b = \dim_k\Hom((0,[1]),-)$ is exact on
  all Auslander-Reiten sequences which do not end in $(0,[1]).$ 
\end{proof}

\medskip
Thus, $b$ (as well as $u,v$) are additive functions on the translation quiver which is obtained
from the Auslander-Reiten quiver of $\Cal S(n)$ by deleting the isomorphism class $[S] =
[(0,[1])]$ from the domain of the translation $\tau$. The translation quiver obtained
in this way will be said to be the {\it adjusted} translation quiver. 

\medskip
Thus, {\it the uwb-vectors are additive on the adjusted Auslander-Reiten quiver of $\Cal S(n)$.} 
In particular, the uwb-vectors for the objects in a tube are known as soon as those
on the boundary of the adjusted component are known.

\subsection{The principal component.}
\label{sec-nine-five}

In our paper, we will need a lot of information about
the principal component $\Cal P(n)$. It will be necessary to know
the uwb-triples $\dfrac{uX|wX}{bX}$ for all $X$ in $\Cal P(n)$ with quasi-length 
at most 7. The following picture shows these uwb-vectors, and even those
of the objects with quasi-length 8.

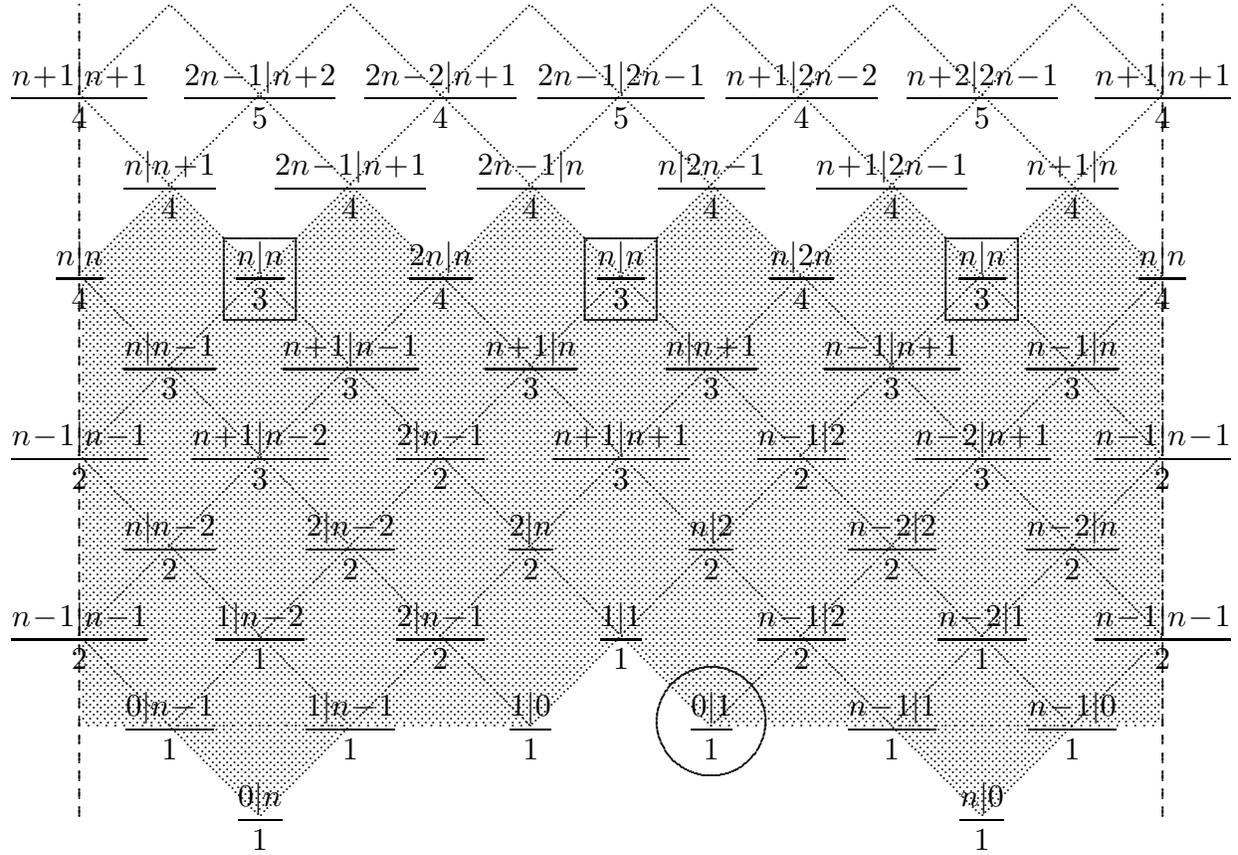
\addcontentsline{lof}{subsection}{The uwb-vectors of the objects in $\mathcal P(n)$.}

$$
{\beginpicture
    \setcoordinatesystem units <.85cm,1cm>
\multiput{} at 0 -.5  12 9.2  /
\setdashes <1mm>
\plot 0 0  0 9 /
\plot 12 0  12 9 /
\setdots <1mm>
\plot 0 1  5 1 /
\plot 7 1  12  1 /
\setdots <.5mm>
\plot 0 4  3 1  /
\plot 0 2  6 8  12 2 /
\plot 0 2  2 0  10 8  12 6 /
\plot 0 4  4 8  11 1 /
\plot 0 6  2 8  10 0  12 2 /
\plot 0 8  7 1  12 6 /
\plot 1 1  8 8  12 4  9 1 /
\plot 0 6  5 1  12 8 /
\plot 0 8  1 9  2 8  3 9  4 8  5 9  6 8  7 9  8 8  9 9  10 8  11 9  12 8 /  
  \setshadegrid span <.4mm>
\vshade 0 1 6 <,z,,> 1 1 7 <z,z,,> 2 0 6 <z,z,,> 3 1 7  <z,z,,> 4 1 6 
    <z,z,,> 5 1 7  <z,z,,> 6 2 6  <z,z,,> 7 1 7   <z,z,,>  8 1 6  <z,z,,>  9 1 7
    <z,z,,>  10 0 6   <z,z,,>  11 1 7   <z,z,,> 12 1 6 /

\put{$\dfrac{0|n}{1}$} at 2 0
\put{$\dfrac{n|0}{1}$} at 10 0

\put{$\dfrac{0|n\!-\!1}{1}$} at 1 1
\put{$\dfrac{1|n\!-\!1}{1}$} at 3 1
\put{$\dfrac{1|0}{1}$} at 5 1
\put{$\dfrac{0|1}{1}$} at 7 1
\put{$\dfrac{n\!-\!1|1}{1}$} at 9 1
\put{$\dfrac{n\!-\!1|0}{1}$} at 11 1

\put{$\dfrac{n\!-\!1|n\!-\!1}{2}$} at 0 2
\put{$\dfrac{1|n\!-\!2}{1}$} at 2 2
\put{$\dfrac{2|n\!-\!1}{2}$} at 4 2
\put{$\dfrac{1|1}{1}$} at 6 2
\put{$\dfrac{n\!-\!1|2}{2}$} at 8 2
\put{$\dfrac{n\!-\!2|1}{1}$} at 10 2
\put{$\dfrac{n\!-\!1|n\!-\!1}{2}$} at 12 2

\put{$\dfrac{n|n\!-\!2}{2}$} at 1 3
\put{$\dfrac{2|n\!-\!2}{2}$} at 3 3
\put{$\dfrac{2|n}{2}$} at 5 3
\put{$\dfrac{n|2}{2}$} at 7 3
\put{$\dfrac{n\!-\!2|2}{2}$} at 9 3
\put{$\dfrac{n\!-\!2|n}{2}$} at 11 3

\put{$\dfrac{n\!-\!1|n\!-\!1}{2}$} at 0 4
\put{$\dfrac{n\!+\!1|n\!-\!2}{3}$} at 2 4
\put{$\dfrac{2|n\!-\!1}{2}$} at 4 4
\put{$\dfrac{n\!+\!1|n\!+\!1}{3}$} at 6 4
\put{$\dfrac{n\!-\!1|2}{2}$} at 8 4
\put{$\dfrac{n\!-\!2|n\!+\!1}{3}$} at 10 4
\put{$\dfrac{n\!-\!1|n\!-\!1}{2}$} at 12 4

\put{$\dfrac{n|n\!-\!1}{3}$} at 1 5
\put{$\dfrac{n\!+\!1|n\!-\!1}{3}$} at 3 5
\put{$\dfrac{n\!+\!1|n}{3}$} at 5 5
\put{$\dfrac{n|n\!+\!1}{3}$} at 7 5
\put{$\dfrac{n\!-\!1|n\!+\!1}{3}$} at 9 5
\put{$\dfrac{n\!-\!1|n}{3}$} at 11 5

\put{$\dfrac{n|n}{4}$} at 0 6
\put{$\dfrac{n|n}{3}$} at 2 6
\put{$\dfrac{2n|n}{4}$} at 4 6
\put{$\dfrac{n|n}{3}$} at 6 6
\put{$\dfrac{n|2n}{4}$} at 8 6
\put{$\dfrac{n|n}{3}$} at 10 6
\put{$\dfrac{n|n}{4}$} at 12 6

\put{$\dfrac{n|n\!+\!1}{4}$} at 1 7
\put{$\dfrac{2n\!-\!1|n\!+\!1}{4}$} at 3 7
\put{$\dfrac{2n\!-\!1|n}{4}$} at 5 7
\put{$\dfrac{n|2n\!-\!1}{4}$} at 7 7
\put{$\dfrac{n\!+\!1|2n\!-\!1}{4}$} at 9 7
\put{$\dfrac{n\!+\!1|n}{4}$} at 11 7

\put{$\frac{n\!+\!1|n\!+\!1}{4}$} at 0 8
\put{$\frac{2n\!-\!1|n\!+\!2}{5}$} at 2 8
\put{$\frac{2n\!-\!2|n\!+\!1}{4}$} at 4 8
\put{$\frac{2n\!-\!1|2n\!-\!1}{5}$} at 6 8
\put{$\frac{n\!+\!1|2n\!-\!2}{4}$} at 8 8
\put{$\frac{n\!+\!2|2n\!-\!1}{5}$} at 10 8
\put{$\frac{n\!+\!1|n\!+\!1}{4}$} at 12 8

\setsolid
\circulararc 360 degrees from 7.6 1  center at 7 1.05 

\plot  1.6 5.5  2.4 5.5  2.4 6.4  1.6 6.4  1.6 5.5  / 
\plot  5.6 5.5  6.4 5.5  6.4 6.4  5.6 6.4  5.6 5.5  / 
\plot  9.6 5.5  10.4 5.5  10.4 6.4  9.6 6.4  9.6 5.5  / 
\endpicture}
$$
We have mentioned that the functions $u,w,b$ are additive on the adjusted
Auslander-Reiten quiver. In our picture, the shaded part consists of meshes 
of the adjusted Auslander-Reiten quiver; these are the meshes which are needed in order
to show the formula in Theorem~\ref{theoremfive} for the ray-simple and the coray-simple
objects of $\Cal P(n)$, see Section~\ref{sec-nine-ten}.
The objects enclosed in square boxes are the central ones. The encircled object is $S$.

In order to verify all the exhibited values, we only have to look
at the boundary of the adjusted principal component and to check the additivity, going upwards.
(Of course, all the objects in the principal component are known, see \cite{RS1}, and one
easily may calculate their uwb-vectors; 
but we want to stress that for the present discussion, only the knowledge of the uwb-vectors
at the boundary is relevant, and these objects all are pickets.)

\subsection{Central objects of quasi-length 6.}
\label{sec-nine-six}

\begin{proposition}
  \label{prop-nine-six}
  Let $Y$ be a central indecomposable object with quasi-length $6$ in the
  Auslander-Reiten component $\Cal C$. Then
  $$
  uY = wY = \beta \Cal C\cdot n,\quad bY = \beta\Cal C\cdot 3,
  $$
  thus
  $ \dfrac{u|w}b Y = \beta\Cal C \cdot\dfrac{n|n}3.$
\end{proposition}

\begin{proof}
  It is sufficient to show that $bY = 3\cdot\beta Y$. Namely, we assume that $Y$ is central,
  thus the pr-vector of $Y$ is $(n/3,n/3)$, and consequently $uY = wY = n/3\cdot bY =
  n\cdot \beta \Cal C.$ 
  
  First, let us assume that $Y$ belongs to a stable component $\Cal C$. Let $Z$ be quasi-simple
  in $\Cal C$ so that $Y = Z[6]$. Then $X$ has a filtration with factors 
  $Z,\ \tau^-Z,\  \dots,\, 
  \tau^{-5}Z.$ As we know, $b$ is additive on $\Cal C$ (see Section~\ref{sec-nine-four}), thus 
  $$
  bY = \sum_{i=0}^5 b(\tau^{-i}Z) = 3(bZ + b(\tau^-Z)) = 3\cdot\beta \Cal C,
  $$
  where we use that $b(\tau^2 X) = bX$ for any reduced object $X$. 
  
  Second, let $\Cal C = \Cal P(n)$. According to Section~\ref{sec-nine-five},
  there are three possibilities for $Y$
  (they are enclosed in a square box),
  and the uwb-vector of $Y$ is always $\dfrac{n|n}3$, in particular $bY = 3.$ 
  On the other hand, $\beta \Cal P(n) = 1.$
  Thus $bY = 3 = 3\cdot \beta \Cal C.$
\end{proof}

\subsection{Proof of Theorem~\ref{theoremfive} for stable tubes.}
\label{sec-nine-seven}

There is an exact sequence
$$
 0 \to X \to X^+ \to X^+/X \to 0.
$$
 We apply Proposition~\ref{prop-nine-six}
 to $Y = X^+/X$ which has quasi-length 6 (and is central as
 all objects in stable tubes of quasi-length 6) and use the additivity of the functions
 $u,w,b$. $\s$
 
\medskip
It remains to look at the principal component. We start with the following lemma
which is valid for all tubes (we need it for the principal component $\Cal P(n)$, see
Section~\ref{sec-nine-nine} but also in Section~\ref{sec-nine-eleven}). 

\subsection{The induction step.}
\label{sec-nine-eight}

\begin{lemma}
  Let $\Cal C$ be an Auslander-Reiten component of 
  $\Cal S(n)$ with $n\ge 6.$ If 
  $$
  \dfrac{u|w}b X^+ =  \dfrac{u|w}b X + \beta \Cal C\cdot \dfrac{n|n}3.
  $$
  holds for $X$ coray-simple in $\Cal C$, then it holds for all
  $X$ in $\Cal C$.
\end{lemma}

\begin{proof}
  If we have an additive function $f$ on a translation quiver, then for any translation
  subquiver of the form
  $$
  
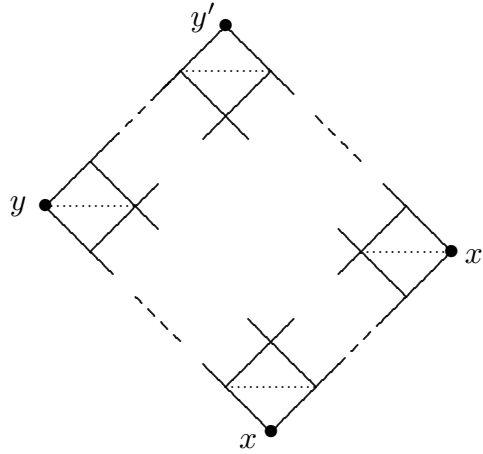
{\beginpicture
    \setcoordinatesystem units <.5cm,.5cm>
\multiput{} at 0 0  9 9 /
\plot 1.5 6.5  0 5  1.5 3.5 /
\plot 3.5 1.5  5 0  6.5 1.5 /
\plot 7.5 2.5  9 4  7.5 5.5 /
\plot 2.5 7.5 4 9  5.5 7.5 /
\plot 1 4  2.5 5.5  /
\plot 1 6  2.5 4.5 /
\plot 4 1  5.5 2.5 /
\plot 4.5 2.5  6 1 /
\plot 3 8  4.5 6.5 /
\plot 3.5 6.5  5 8 /
\plot 6.5 3.5  8 5 /
\plot 6.5 4.5  8 3 /
\multiput{$\bullet$} at 0 5  5 0  9 4  4 9 /
\setdashes <1mm>
\plot 2 3  3 2  /
\plot 1.5 6.5  3 8  /
\plot 6.5 1.5  7.5  2.5 /
\plot 6 7  7 6 /
\setdots <.7mm>
\plot 0 5  2 5  /
\plot 4 1  6 1  /
\plot 7 4  9 4 /
\plot 3 8  5 8   /
\put{$x$} at 4.5 -.2
\put{$x'$} at 9.6 4

\put{$y$} at -.6 5
\put{$y'$} at 3.5 9.2 
\endpicture}
  $$
  (with arrows going from left to right)
  the additivity of $f$ implies that 
  \begin{equation*}
    f(y')-f(y) = f(x')-f(x)
    \tag{$*$}
  \end{equation*}
  
  Assume that the assertion of Theorem~\ref{theoremfive} holds for 
  all coray-simple objects $X$ in $\Cal C$.
  Let $Y$ be in $\Cal C$ and $Y' = Y^+$. Then we find a coray-simple object $X$ and 
  a rectangle as exhibited above with 
  $y = [Y],\ y' = [Y^+]$ and $x = [X],\ x' = [X^+]$. 
  According to $(*)$, we see that the assertion of Theorem~\ref{theoremfive}  holds for $Y$.
\end{proof}

\subsection{Proof of Theorem~\ref{theoremfive} for the principal component $\Cal P(n)$.}
\label{sec-nine-nine}

In Section~\ref{sec-nine-five},
we have shaded the meshes which are needed in order to calculate the ubw-vectors of
the objects $X^+$ with $X$ being coray-simple. For $X$ coray-simple, we have to compare its
uwb-vector with the uwb-vector of $X^+$. We 
see
$$
 \dfrac{u|w}{b} X^+ - \dfrac{u|w}{b} X = \dfrac{n|n}{3}.
$$
for all coray-simple objects $X$ in $\Cal P(n)$. 
Since $\beta X = 1$, the formula in Theorem~\ref{theoremfive}  is satisfied for the 
coray-simple objects in $\Cal P(n)$, thus, according to Section~\ref{sec-nine-five},  the formula 
is valid in general for the principal component. 
This concludes the proof of Theorem~\ref{theoremfive} for 
the principal component, thus for all components. 
$\s$

\subsection{The global space of the objects in $\Cal P(n)$.}
\label{sec-nine-ten}

We have shown in Section~\ref{sec-nine-nine}: If $X$ belongs to $\Cal P(n)$, then
$u(X^+) = uX + n,\ w(X^+) = wX + n$ (therefore
$v(X^+) = vX + 2n$), and $b(X^+) = bX +3.$ These assertions can be refined as
follows:

\medskip
\begin{proposition}
  \label{prop-nine-ten}
  If $X = (U,V,W)$ belongs to $\Cal P(n)$, then:
  \begin{itemize}[leftmargin=3em]
    \item[{\rm (a)}] The global space of $X^+$ is 
      $$
      V(X^+) = VX\oplus[n,n-2,2]
      $$ 
      with subspace $U(X^+) = UX\oplus [n-2,2]$ and factor space 
      $W(X^+) = WX\oplus [n-2,2]$.
    \item[{\rm (b)}]
      The parts occurring in $U, V, W$ are $[1], [2], [n-2], [n-1], [n].$
  \end{itemize}
\end{proposition}

\medskip
Note that (a) implies immediately 
that $u(X^+) = uX+n$, $v(X^+) = vX+2n$, $w(X^+) = wX+n$, and
$b(X^+) = bX+3.$

\medskip
Outline of a proof of Proposition~\ref{prop-nine-ten}:
We denote by $[\Cal N(n)]$ the set of isomorphism classes $[N]$
of the objects of $\Cal N(n)$, or, what is the same, the set of partitions of height 
at most $n$.
For $1\le t \le n,$ let $[\Cal N(n)]_t$ be the subset of all $[N]$ in $[\Cal N(n)]$
such that $[t]$ is a direct summand of $N$, thus the set of partitions of height 
at most $n$ with at least one part equal to $t$. We define functions $\kappa_i$ with
$i\in \mathbb N_1$ and with $\kappa_i = \kappa_j$ for $i\equiv j \mod\ 6$, all with values in
$[\Cal N(n)]$, as follows: The functions $\kappa_1, \kappa_4, \kappa_5$ are defined on all
of $[\Cal N(n)]$, the functions $\kappa_2, \kappa_3,\kappa_6$ are defined on 
$[\Cal N(n)]_n, [\Cal N(n)]_{n-1},$ and $[\Cal N(n)]_1$, respectively, namely as follows:
\begin{align*}
 \kappa_1[N] &= [N\oplus [n]],\cr
 \kappa_2[N\oplus[n]] &= [N\oplus[n\!-\!1]],\cr
 \kappa_3[N\oplus[n\!-\!1]] & = [N\oplus[n,n\!-\!2]],\cr
 \kappa_4[N] &= [N],\cr
 \kappa_5[N] &= [N\oplus [1]],\cr
 \kappa_6[N\oplus[1]] &= [N\oplus[2]].
\end{align*}
It is easy to check that for all $i\in \mathbb N_1$, we have
\begin{equation*}
 \kappa_{i+5}\kappa_{i+4}\cdots \kappa_{i+1}\kappa_{i}[N] = [N\oplus[n,n-2,2]], \tag{1}
\end{equation*}
provided $\kappa_{i}[N]$ is defined. 

\medskip\medskip
Let $Z_1 = S,$ $Z_2=([n],[n])$,
$Z_3 = ([n-1],[n-1])$, $Z_4 = (0,[n]),$ $Z_5 = ([1],[n]),$ $Z_6 = ([1],[1]),$
and $Z_j = Z_i$ provided $i\equiv j\mod \ 6.$
Thus, $Z_1,\dots,Z_6$ are the coray-simple objects.
$$
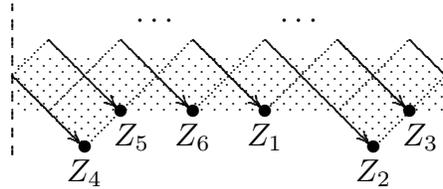
{\beginpicture  
    \setcoordinatesystem units <.4cm,.4cm>
\multiput{} at 0 0  12 4  /
\setdashes <1mm>
\plot 0 -.2  0 4 /
\plot 12 -.1  12 4 /
\put{$Z_1$} at 7 .3
\put{$Z_2$} at 10 -.7
\put{$Z_3$} at 11.3 .3
\put{$Z_4$} at 2 -.7 
\put{$Z_5$} at 3.3 .3
\put{$Z_6$} at 5 .3
\multiput{$\bullet$} at  3 1  2 0  5 1  7 1  11 1  10 0  /
\setdots <.5mm>
\plot 0 2  2 0  5 3  7 1  9 3  11 1 /
\plot 0 2  1 3 /
\plot 1 1  3 3 /
\plot 5 1  7 3 /
\plot 9 1  11 3 /
\plot 10 0  12 2 /
\setsolid
\arr{0 2}{1.9 0.1}
\arr{1 3}{2.9 1.1}
\arr{3 3}{4.9 1.1}
\arr{5 3}{6.9 1.1}
\arr{7 3}{9.9 0.1}
\arr{9 3}{10.9 1.1}
\plot 11 3  12 2 /
\setshadegrid span <.7mm>
\vshade 0 1 2.5  <,z,,> 1 1 2.5  <z,z,,> 2 0  2.5  <z,z,,> 3 1 2.5  <z,z,,> 8 1 2.5 
    <z,z,,> 9 1 2.5   <z,z,,> 10 0 2.5 <z,z,,> 11 1 2.5   <z,z,,> 12 1 2.5  /
\multiput{$\cdots$} at 4 3.5  8 3.5 /
\endpicture}
$$
We denote by $[t]Z$ the object of quasi-length $t$ in the coray ending in $Z$, thus
any object in $\Cal P(n)$ is of the form $[t]Z$, and
a sectional path $X \to X^+$ is of the form
$$
 X = [t]Z_i \to [t\!+\!1]Z_{i+1}  \to \cdots \to [t\!+\!5]Z_{i+5}  \to [t\!+\!6]Z_{i+6} = X^+.
$$
There is the following formula 
$$
  [V([t+1]Z_{i+1})] = \kappa_i(V([t]Z_i))
$$
which can be shown by induction on $t$.

\medskip
It follows that 
\begin{equation*}
 [V(X^+)] = \kappa_{i+5}\cdots \kappa_{i+1}\kappa_{i}[VX]. \tag{2}
\end{equation*}

\medskip
The equalities (2) and (1) show:
$$
 [V(X^+)] = \kappa_{i+5}\cdots \kappa_{i+1}\kappa_{i}[VX] = [VX\oplus [n,n-2,2]].
$$

\medskip
Proof of (b). It is sufficient to consider $V$. Let $\mathbf N$ be the subset of $|\Cal N(n)|$
given by the objects $N$ whose indecomposable direct summands are of the form $[1],[2],[n-2],
[n-1],[n].$ Then we have $\kappa_i(\mathbf N) \subseteq \mathbf N,$ for $1\le i \le 6.$ Also,
if $X$ is ray-simple in $\Cal P(n)$, then $VX$ belongs to $\mathbf N.$ It follows by induction on
the quasi-length that for all objects $X$ in $\Cal P(n)$, we have $[VX] \in \mathbf N.$
$\s$

\medskip\medskip
Here are the
\phantomsection{partitions}
\addcontentsline{lof}{subsection}{The partitions $V X$ of the objects $X$ in $\mathcal P(n)$.}
$[VX]$ for the objects $X$ in $\Cal P(n)$ 
with quasi-length at most 7. 
For the objects of the form $X^+$ (they have quasi-length 6 or 7), we have encircled
$[VX]$ (it occurs as a part of the
partition $[V(X^+)]$).

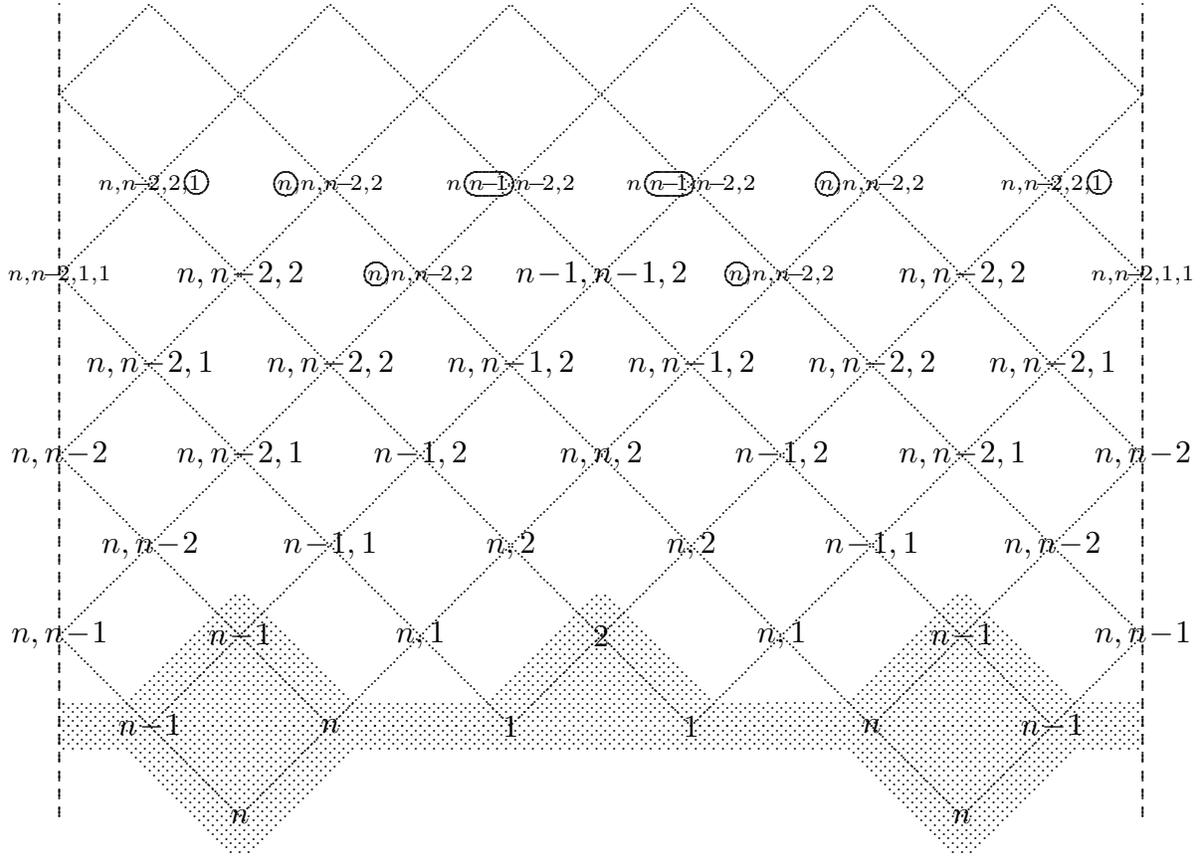
$$
{\beginpicture
    \setcoordinatesystem units <.85cm,1cm>
\multiput{} at 0 -.5  12 8.5  /
\setdashes <1mm>
\plot 0 0  0 9 /
\plot 12 0  12 9 /
\setdots <1mm>
\plot 0 1  5 1 /
\plot 7 1  12  1 /
\setdots <.5mm>
\plot 0 4  3 1 /
\plot 0 2  6 8  12 2 /
\plot 0 2  2 0  10 8  12 6 /
\plot 0 4  4 8  11 1 /
\plot 0 6  2 8  10 0  12 2 /
\plot 0 8  7 1  12 6 /
\plot 1 1  8 8  12 4  9 1 /
\plot 0 6  5 1  12 8 /
\plot 0 8  1 9  2 8  3 9  4 8  5 9  6 8  7 9  8 8  9 9  10 8  11 9  12 8 /  

\multiput{$n$} at 2 0  3 1  9 1  10 0 /
\multiput{$n\!-\!1$} at 1 1  2 2  10 2  11 1  /
\multiput{$1$} at 5 1  7 1 /
\multiput{$2$} at 6 2  /
\multiput{$n,n\!-\!1$} at 0 2  12 2  /
\multiput{$n,1$} at 4 2  8 2  /
\multiput{$n,n\!-\!2$} at 1 3  11 3  0 4  12 4 /
\multiput{$n\!-\!1,1$} at 3 3  9 3   /
\multiput{$n,2$} at 5 3  7 3   /
\multiput{$n\!-\!1,2$} at 4 4  8 4  /
\multiput{$n,n\!-\!2,1$} at 2 4  10 4  1 5  11 5 /
\multiput{$n,n,2$} at 6 4 /
\multiput{$n,n\!-\!1,2$} at 5 5  7 5  /
\multiput{$n,n\!-\!2,2$} at 3 5  9 5  2 6  10 6 /
\multiput{$n\!-\!1,n\!-\!1,2$} at 6 6  /
\multiput{$\ss n,n\!-\!2,1,1$} at 0 6  12 6  /
\multiput{$\ss n,n,n\!-\!2,2$} at 4 6  8 6   3 7  9 7 /
\multiput{$\ss n\;n\!-\!1\;n\!-\!2,2$} at 5 7  7 7 /
\multiput{$\ss n,n\!-\!2,2,1$} at 1 7  11 7 /
\setsolid

\circulararc 360 degrees from 1.65 7.03  center at 1.52 7.03
\circulararc 360 degrees from 11.65 7.03  center at 11.52 7.03
\circulararc 360 degrees from 2.64 7.01  center at 2.51 7.01
\circulararc 360 degrees from 8.64 7.01  center at 8.51 7.01
\circulararc 360 degrees from 3.64 6.01  center at 3.51 6.01
\circulararc 360 degrees from 7.64 6.01  center at 7.51 6.01

\circulararc 180 degrees from 6.62 7.14  center at 6.62 7.01
\circulararc -180 degrees from 6.9 7.14  center at 6.9 7.01
\plot 6.9 7.14  6.62 7.14 /
\plot 6.9 6.88  6.62 6.88 /

\circulararc 180 degrees from 4.62 7.14  center at 4.62 7.01
\circulararc -180 degrees from 4.9 7.14  center at 4.9 7.01
\plot 4.9 7.14  4.62 7.14 /
\plot 4.9 6.88  4.62 6.88 /

\setshadegrid span <.5mm>
\vshade 0 0.7 1.3 <z,z,,> 0.7  0.7 1.3   <z,z,,> 2  -.5 2.5 <z,z,,>  3.3  .7 1.3 
                   <z,z,,>  4.7  0.7 1.3   <z,z,,> 6  .7 2.5 <z,z,,>  7.3  .7 1.3 
                   <z,z,,>  8.7  0.7 1.3   <z,z,,> 10  -.5 2.5 <z,z,,>  11.3  .7 1.3 
                   <z,z,,> 12 0.7 1.3 /
\endpicture}
$$

\medskip
\begin{remark}
  Theorem~\ref{theoremfive} provides a numerical comparison between $X$ and $X^+$,
  for any indecomposable object $X$ in $\Cal S(n)$. Also, we know that there is a 
  sectional path from $X$ to $X^+$ and one may analyse this path.
  In case $X$ belongs to a stable tube, this path yields a monomorphism (whose cokernel
  belongs to the same component as $X$); in particular $X$ is a subobject of $X^+$.
  In case $X$ belongs to the principal component $\Cal P(n)$, we may get a different
  behaviour.
\end{remark}

\medskip
\begin{proposition}
  Let $X$ be indecomposable in $\Cal S(n)$. The following conditions are
  equivalent:
  \begin{itemize}[leftmargin=3em]
  \item[\rm(i)] $X$ is a subobject of $X^+$.
  \item[\rm(ii)] A sectional path $X\to \cdots \to X^+$ yields a monomorphism $X \to X^+$.
  \item[\rm(iii)] $X$ does not belong to the coray ending in $([n],[n]).$
  \item[\rm(iv)] There is no sectional path from $X$ to  $([n],[n]).$
  \end{itemize}
\end{proposition}

\begin{proof}[Proof (outline)]
  We may assume that $X$ belongs to $\Cal P(n)$. 
  If $X$ does not belong to the coray ending in $([n],[n])$, then it is easy to see
  that again any sectional path $X\to \cdots \to X^+$ yields a monomorphism $X \to X^+$,
  thus $X$ is again a subobject of $X^+.$
  
  On the other hand, we have to consider the case that
  $X$ belongs to the coray ending in $([n],[n])$.
  Then the Auslander-Reiten quiver exhibits a short exact sequence
  $$0\to X\to X[1]\oplus([n],[n])\to ([n-1],[n-1])\to 0$$ hence the map $X\to X[1]$
  cannot be a monomorphism.
  But this implies that $\mu$ has a non-trivial kernel. A contradiction. 
\end{proof}

\subsection{Half-lines.}
\label{sec-nine-eleven}

As the title of this section indicates, we also will deal with half-lines: 
These are half-lines in the pr-plane $\mathbb T(n)$, 
in the usual sense: If we start with a line in $\mathbb T(n)$, any vertex $x$ of the line
decomposes the line into two half-lines with initial vertex $x$ (it is common to consider
$x$ as an element of both half-lines, and the two half-lines which are obtained 
will be said to be {\it complementary}). If $y\neq x$ are different vertices in $\mathbb T(n)$,
we write $[x,y\rangle$ for the half-line with initial vertex $x$ passing through $y$.

We recall that the center of $\mathbb T(n)$ is $z(n) = (\frac n3,\frac n3)$ and that 
a line in $\mathbb T(n)$ is said to be central, provided $z(n)$ belongs to it; a half-line is
said to be central provided its initial vertex is $z(n)$.
Also, an indecomposable object $X$ with $\pr X = z(n) = (\frac n3,\frac n3)$ is said to
be central.
Note that {\it an indecomposable object $X$ with $\tau^2X = X$ is central.}

\subsection{Arithmetical sequences.}
\label{sec-nine-twelve}

We are going to prove the Corollary of Theorem~\ref{theoremfive} in Section~\ref{sec-one-six}.
Starting with a non-central 
object $X$, the central half-line which passes through $\pr X$ contains infinitely many vertices of the
form $\pr X_i$, with $X_0 = X$ such that all the objects $X_i$ belong to the 
same (categorical) ray. 

\medskip
There is the following obvious generalization of Theorem~\ref{theoremfive}
(see also Section~\ref{sec-nine-eight}). We define
inductively $X^{+t}$ for $t\in \mathbb N_0$ 
as follows: $X^{+0} = X$, and if $X^{+t}$ is already defined, then
$X^{+(t+1)} = (X^{+t})^+.$

\medskip
\begin{theorem}[Theorem~\ref{theoremfive} (general form)]
  Let $X$ be indecomposable in the Auslander-Reiten component
  $\Cal C$.
  Then for all $t\in \mathbb N_0$
$$
\dfrac{u|w}b X^{+t} =  \dfrac{u|w}b X + \beta \Cal C\cdot t\cdot \dfrac{n|n}3.
$$
\end{theorem}

\medskip
\begin{corollary}
  \label{cor-nine-twelve}
  Let $X$ be indecomposable.
  If $X$ is central, all objects $X^{+t}$ are central. 
  Otherwise, the sequence $X = X^{+0}, X^{+1}, X^{+2},\ldots$ 
  lies on the half-line $[z(n),\pr X\rangle$
    in $\mathbb T(n)$.
    The distance of $\pr X^{+t}$ to the center $z(n)$ is strictly
    decreasing and converges to zero. 
\end{corollary}

\medskip
(In particular, the vertices $\pr X_t$ are pairwise
  different.)

\medskip
If $X$ is not central, we may visualize the setting as follows: On
the half-line $[z(n),\pr X\rangle$, we mark the (pr-vectors of the)
objects $X_t$ by bullets.
$$
{\beginpicture
    \setcoordinatesystem units <1cm,.5cm>
\multiput{} at -2 0  /
\plot -2 0  7.2 0 /
\plot 7.7 0  8 0 /
\setdots <1mm>
\plot 7.2 0  7.8 0 /
\multiput{$\bullet$} at 0 0  4 0  6 0  7 0 /
\put{$\blacksquare$} at 8 0 
\put{$X = X^{+0}$} at 0 -1
\put{$X^{+1}$} at 4 -1
\put{$X^{+2}$} at 6 -1
\put{$X^{+3}$} at 7 -1
\put{$\cdots$} at 7.7 -1.2
\put{$z(n)$} at 8 1

\endpicture}
$$

\medskip        
\begin{remark}
  We may add to Theorem~\ref{theoremfive} and Corollary~\ref{cor-nine-twelve}
  the following observation:
  If $X$ is indecomposable, then 
  $$
  \dfrac{u|w}b X^{+t} =  
  \dfrac{u|w}b (X \oplus Y^t),
  $$
  where $Y = X^+/X$ in case $X$ belongs to a stable tube,
  and where $Y$ is one of the indecomposable
  objects in $\Cal P(n)$ with $\dfrac{u|w}b Y = \dfrac{n|n}3$, in case $X$ belongs to $\Cal P(n)$.
  Thus, always $Y$ is indecomposable and belongs to the same Auslander-Reiten component as $X$. 
  $\s$
\end{remark}

\subsection{The sequence $d(X^{+t})$ converges.}
\label{sec-nine-thirteen}

The boundary distance $dX$ of an indecomposable object $X\in\Cal S(n)$ has been introduced in
Section~\ref{sec-one-four}.

\medskip
\begin{proposition}
  \label{prop-nine-thirteen}
  Let $X$ be indecomposable in $\Cal S(n)$. Then $X$ is central if and only if
  $X^+$ is central.
  Also, $X$ is $u$-minimal if and only if $X^+$ is $u$-minimal. If $X$ is not central,
  then $dX < d(X^+)$ and the sequence $d(X^{+t})$ converges to $n/3$.
\end{proposition}

\begin{proof}
  The first two assertions follow directly from Theorem~\ref{theoremfive}.
  It remains to show the last
  assertions. Thus, assume that $X$ is not central. Assume that $X$ belongs to the component 
  $\Cal C$ and let $\beta = \beta\Cal C$. We may assume that $X$ is $u$-minimal,
  thus $uX/bX < n/3,$ therefore $3uX < nbX.$ It follows that
  $$
  dX = uX/bX < (uX+n\beta)/(bX+3\beta) = d(X^+).
  $$
  It is easy to verify that 
  the sequence $d(X^{+t}) = (uX+n\beta\cdot t)/(bX+3\beta\cdot t)$ converges to
  $n/3.$
\end{proof}

\subsection{Objects with fixed boundary distance.}
\label{sec-nine-fourteen}

\begin{proposition}
  \label{prop-nine-fourteen}
  If $\Cal C$ is a stable
  component, the number of objects in $\Cal C$ with fixed boundary distance $d < n/3$ 
  is at most $30$.
  The number of objects in $\Cal P(n)$ with fixed boundary distance $d < n/3$ is at most $33$.
\end{proposition}

\begin{proof}
  This follows from the inequality $dX < d(X^+)$ given in Section~\ref{sec-nine-twelve}
  and the following obvious assertions (see the end of Section \ref{sec-nine-one}).
  If
  $\Cal C$ is a stable component, then any plus-orbit of
  non-central objects in $\Cal C$ has a representative which has quasi-length at most 5,
  thus the number of plus-orbits is at most $30.$
  Any plus-orbit of non-central objects in $\Cal P(n)$ has a representative 
  which has quasi-length $0\le \ell\le 5$ or is equal to $S[6],$
  thus the number of plus-orbits is at most $33.$
\end{proof}

\medskip
\begin{remark}
  We will see at the beginning of Section~\ref{sec-ten-seven}
  that the number $33$ for $\Cal P(n)$ is optimal. 
\end{remark}

\bigskip
\section{Half-line support and triangle support of a component.}
\label{sec-ten}

In this section, we usually suppose that $n\ge 6,$ but in Sections \ref{sec-ten-nine}
and \ref{sec-ten-ten}, we
deal with all values $n\ge 1.$

\medskip
We have attached to every indecomposable non-central object $X$ the central half-line line
$H(X)$ which starts at the center $z(n)$ and contains the pr-vector of $X$. 
Since we have $H(X^+) = H(X)$, we see that for any Auslander-Reiten component $\Cal C$, we 
may obtain at most 36 different half-lines. Actually, since for $\Cal C$ being stable, 
the objects of quasi-length 6 are central, we see that the half-line support of
a stable component consists of at most 30 half-lines.
Similarly, the three central objects in $\Cal P(n)$ with quasi-length 6 show that the
half-line support of $\Cal P(n)$ can contain at most 33 central half-lines. 
As we will see below, the number 30 for stable components is optimal.
However, the half-line support of $\Cal P(n)$ consists of only 18 half-lines in case $n = 6,$
and of  24  half-lines in case $n\ge 7.$ 

\bigskip
First, we look at the half-line support of a stable tube.

\subsection{Complementary half-lines.}
\label{sec-ten-one}

We look at a stable tube $\Cal C$ and
want to show: If a central half-line supports $\Cal C$, also the complementary 
half-line supports $\Cal C$. There is the following Proposition.

\medskip
\begin{proposition}
  \label{prop-ten-one}
  Let $X$ be indecomposable object in some  
  stable component $\Cal C$. Then there is $X'$ in $\Cal C$ such that
  $$
  \dfrac{u|w}b X  +  \dfrac{u|w}b X' 
  $$
  is an integral multiple of $\beta \Cal C\cdot \dfrac{n|n}3$.
\end{proposition}

\begin{proof}
  We assume that $X = Z[6t-i]$ for some ray-simple object $Z$, with $t\ge 1$ and
  $0 \le i \le 5.$ If $X$ is not central, the quasi-length $6t-i$ of $X$ is not divisible by 6,
  thus $i > 0.$ Let $X' = Z[6t]/Z[6t-i]$. Since $i \ge 1$, $X'$ is indecomposable and belongs to $\Cal C.$ Also, we have
  $$
  \dfrac{u|w}b X  +  \dfrac{u|w}b X'  = \dfrac {u|w}b Z[6t]
  $$ 
  and $\dfrac {u|w}b Z[6t]$ is a multiple of $\beta \Cal C\cdot \dfrac{n|n}3$, see
  Section~\ref{sec-nine-seven}.
\end{proof}

\medskip
There is the following consequence:

\medskip
\begin{corollary}
  If $X$ is indecomposable, non-central and belongs to a stable tube
  $\Cal C$, then there is $X'$ in $\Cal C$ such that $X$ and $X'$ live on 
  complementary half-lines.
  $$
  
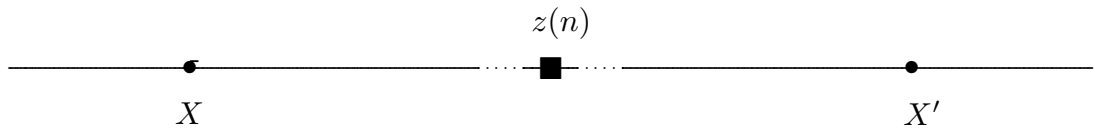
{\beginpicture
    \setcoordinatesystem units <1cm,.5cm>
\multiput{} at -2 0  10 0 /
\plot -2 0  3.2 0 /
\plot 10 0  4.8 0 /
\plot 3.7 0  4.3 0 /
\setdots <1mm>
\plot 3.2 0  3.8 0 /
\plot 4.2 0  4.8 0 /
\multiput{$\bullet$} at 0 0  8 0 /
\put{$\blacksquare$} at 4 0 
\put{$X$} at 0 -1
-\put{$X'$} at 8 -1
\put{$z(n)$} at 4 1

\endpicture}
$$
\end{corollary}

\begin{proof}
  Let $X = Z[\ell]$ with $Z$ ray-simple and $\ell \ge 1.$ According to
  Proposition~\ref{prop-ten-one},
  $X$ lives on the same central half-ray as some $Z[t]$ with $1\le t \le 5.$
  Let $X' = \tau^{-t}Z[6-t]$. The formula in Proposition~\ref{prop-nine-six}
  asserts that $X'$ lives on the
  complementary central half-line.
\end{proof}

\medskip
\begin{remark}
  Proposition~\ref{prop-ten-one} provides a recipe how to relate complementary half-lines
  in the support of a stable tube. We may refine these considerations by taking into
  account the precise $\tau$-period of the objects in question.
  
  Namely, let $X$ be a non-central indecomposable in a stable tube. 
  Then the $\tau$-period of $X$ is $3$ or $6$. 
\end{remark}

\medskip
\begin{corollary}
  \label{cor-ten-one-two}
  The half-line support of a stable tube $\Cal C$ consists of $g$  pairs
  of complementary central half-lines, where $0 \le g \le 15$ is divisible by $3$.
  
  If $\Cal C$ has rank $3$, the half-line support of $\Cal C$
  consists of $g$ pairs of complementary half-lines, where $g = 0$ or $g = 3.$
  $\s$
\end{corollary}

\subsection{Example: $15$ complementary pairs of half-lines.}
\label{sec-ten-two}

Here we deal with the
Auslander-Reiten component in $\Cal S(7)$ which contains $E_2^3.$

\smallskip
We start with the bipicket $X = ([4],[7,2],[4,1])$, as shown on the right:
$$
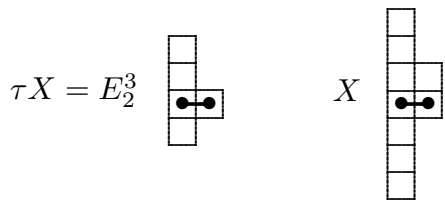
\beginpicture
   \setcoordinatesystem units <0.3cm,0.3cm>
\put{\beginpicture
\multiput{} at 0 0  2 4 /
\plot 0 0  1 0  1 4  0 4  0 0 /
\plot 0 1  2 1  2 2  0 2 /
\plot 0 3  1 3 /

\multiput{$\bullet$} at 0.5 1.5  1.5 1.5 /
\plot 0.5 1.5  1.5 1.5 / 
\plot 0.5 1.45  1.5 1.45 / 
\plot 0.5 1.55  1.5 1.55 / 
\put{$\tau X = E_2^3$} at -3.5 2
\endpicture} at 0 .5
\put{\beginpicture
\multiput{} at 0 0  2 7 /
\plot 0 0  1 0  1 7  0 7  0 0 /
\plot 0 1  1 1 /
\plot 0 2  1 2 /
\plot 0 3  2 3  2 5  0 5 /
\plot 0 4  2 4 /
\plot 0 6  1 6 /
\multiput{$\bullet$} at 0.5 3.5  1.5 3.5 /
\plot 0.5 3.5  1.5 3.5 / 
\plot 0.5 3.45  1.5 3.45 / 
\plot 0.5 3.55  1.5 3.55 / 
\put{$X$} at -1.5 4
\endpicture} at 10 0 
\endpicture
$$

According to Section~\ref{sec-three}, we have $\tau X = ([2],[4,1],[3]) = E_2^3$ and $\tau^2 X = 
([4,1],[7,3],[5])$. Put $Y=\tau X[2]$ and $Z=\tau^2X[3]$.  Then:
$$
\beginpicture
   \setcoordinatesystem units <2.5cm,0.6cm>
\multiput{} at 0 0  2 5 /
\put{$X$}        at 0 5
\put{$\tau X$}   at 0 4
\put{$\tau^2 X$} at 0 3
\put{$Y$}        at 0 2
\put{$\tau Y$}   at 0 1
\put{$Z$}        at 0 0

\put{uwb-vector} at 1 6
\put{$\frac{4|5}2$}   at 1 5
\put{$\frac{2|3}2$}   at 1 4
\put{$\frac{5|5}2$}   at 1 3
\put{$\frac{6|8}4$}   at 1 2
\put{$\frac{7|8}4$}   at 1 1
\put{$\frac{11|13}6$}   at 1 0

\put{$\phi = \frac{3u-7b}{3w-7b}$} at 2 6
\put{$-2$}   at 2 5
\put{$\frac85$}   at 2 4
\put{$1$}   at 2 3
\put{$\frac52$}   at 2 2
\put{$\frac74$}   at 2 1
\put{$3$}   at 2 0

\endpicture
$$
The number $\phi$ in the last column
is the ratio for the line through $z(7)$ and the pr-vector $(\frac ub,\frac wb)$, see 
Section~\ref{sec-ten-three}.

\medskip
Here are the pr-vectors of $X,\ \tau X,\ \tau^2X,\ Y,\ \tau Y,\ Z$ in $\mathbb T(7)$:
$$
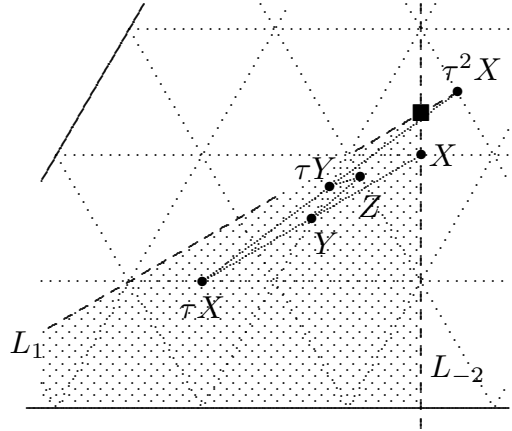
\beginpicture
   \setcoordinatesystem units <.808cm,1.4cm>
\multiput{} at 0 0  6.2 3.2 /
\plot -.4 0  6.2 0 /
\plot -.2 1.8  1.2 3.2 /
\setdots <1mm>
\plot -.2 .2  0 0  3.2 3.2 /
\plot 0 2  2 0  5.2 3.2 /
\plot 1 3  4 0  6.2 2.2 /
\plot 2.8 3.2  6 0  6.2 0.2 /
\plot 4.8 3.2  6.2 1.8 /
\plot -.2 1  6.2 1 /
\plot 0 2  6.2 2 /
\plot 1 3  6.2 3 /
\put{$\ssize \blacksquare$} at 5 2.333 
\multiput{$\ssize \bullet$} at 5 2  2 1  5.5 2.5  3.5 1.5  3.75 1.75  4.1667  1.833 /
\put{$X$}        at 5.3 2 
\put{$\tau X$}   at 2 .8
\put{$\tau^2 X$} at 5.7 2.7 
\put{$Y$}        at 3.7 1.3
\put{$\tau Y$}   at 3.55 1.9 
\put{$Z$}        at 4.3 1.6 

\setdots <.4mm>
\plot 5 2  2 1  /
\plot  2 1  5.5 2.5  /
\plot   3.75 1.75  4.1667  1.833  3.5 1.5  /
\setdashes <1mm>
\plot 5 3.2  5 -.2 /
\plot 5.5 2.5  4 2 /
\plot 3 1.667 -.2 0.6 /
\put{$L_1$} at -.4 0.5 
\put{$L_{-2}$} at 5.5 0.3 
\setshadegrid span <.7mm>
\vshade -.2 0 0.6  <z,z,,> 5 0 2.333  /

\endpicture
$$

We see that the uwb-vectors of $X,\ \tau X,\ Y,\ \tau Y,\ Z$ lie in the 
(shaded) fundamental
region bounded by the lines $L_1,\ L_{-2}$ and $p = 0$ (and not on $L_1$).  If we apply 
the rotations $\rho$ and $\rho^2$ (thus 
$\tau^2$ and $\tau^4$), the 5 central half-lines which pass through
$X,\ \tau X,\ Y,\ \tau Y,\ Z$ yield 
a set of 15 half-lines which belong to the half-line support of $\Cal C$, and which contains
no pair of complementary half-lines. Thus the half-line support
of $\Cal C$ consists of 15 pairs of complementary half-lines.

\subsection{Some central lines and half-lines.}
\label{sec-ten-three}

A central line $L$ which is not parallel to a boundary line, is the union of two
central half-lines which are denoted by $L_s$ (the short one) and $L_{\ell}$ (the long
one); here ''short'' and ''long'' refer to the intersection with $\mathbb T(n)$; for
example, 
in case $L$ contains a vertex of the form $(0,r)$ with $0\le r\le n,$ then $(0,r)$
belongs to 
$L_s$ if and only if $\frac n3< r < \frac {2n}3.$).

We denote by $\mathbb P(n)$ the union of the
\phantomsection{central}
\addcontentsline{lof}{subsection}{The lines $\Bbb P(n), \Bbb D(n), \Bbb H_\ell(n), \Bbb K_s(n)$.}%
lines 
which are parallel to the three coordinate axes, to the boundary lines:
These are the lines $p = \frac n3,\ r = \frac n3,\ n-q = \frac n3.$ 
We denote by $\mathbb D(n)$ the union of the diagonal lines (the reflection lines):
These are the lines $p = r,\ r = n-q,\ p = n-q.$ 

We denote by $\mathbb H(n)$ the union of the central lines which pass throught a vertex in
the $\Sigma_3$-orbit of $\pr S = (0,1)$, and by $\mathbb K(n)$, 
for $n\ge 5$, the central lines which pass throught a vertex in the $\Sigma_3$-orbit of $\pr E_2^{n-2} = (1,(n-2)/2)$. 
Note that for $n = 6$, $\mathbb K(n) = \mathbb P(n)$, since in this case $\pr E_2^{n-2}= (1,2)$. 
	

Here are the lines $\mathbb P(n), \mathbb D(n)$ and the half-lines $\mathbb H_\ell(n)$ 
(for $n\ge 4$), as well as the half-lines $\mathbb K_s(n)$ (for $n\ge 7$):
$$  
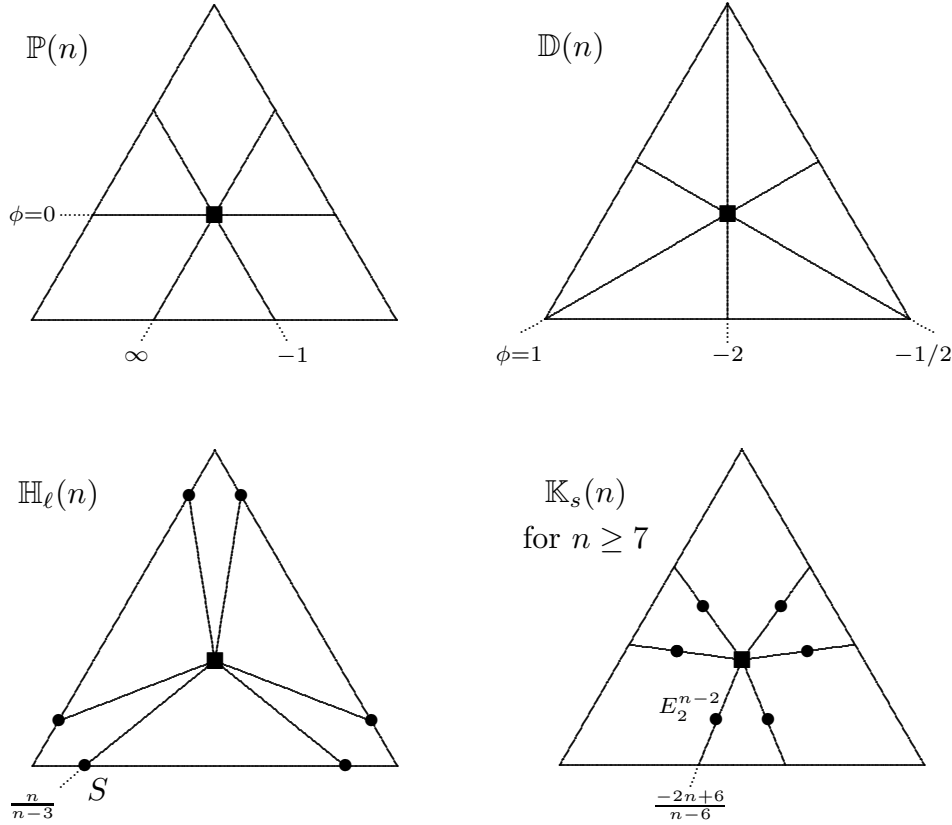
{\beginpicture
 \setcoordinatesystem units <.28868cm,.5cm>
\put{\beginpicture
\multiput{} at -7 0  7 0  0 7  /
\plot  -7 0  7 0  0 7  -7 0 /
\plot -2.333 0  2.333 4.667  /
\plot  2.333 0  -2.333 4.667  /
\plot  -4.667 2.333  4.667 2.333 /

\setdots <.5mm>
\plot -2.333 0  -2.85 -.5 /
\plot  2.333 0   2.85 -.5 /
\plot  -4.667 2.333  -6 2.333 /

\put{$\ss \blacksquare$} at 0 2.333 
\put{$\mathbb P(n)$} at -6 6 
\put{$\ss\infty$} at -3 -.8
\put{$\ss-1$} at 3 -.8
\put{$\ss \phi = 0$} at -7 2.333
\endpicture} at 0 0
\put{\beginpicture
\multiput{} at -7 0  7 0  0 7  /
\plot  -7 0  7 0  0 7  -7 0 /
\plot 0 2.333  -7 0 /
\plot 0 2.333  7 0 /
\plot 0 2.333  0 7 /
\plot 0 2.333  3.5 3.5 /
\plot 0 2.333  -3.5 3.5 /
\plot 0 2.333  0 0 /

\setdots <.5mm>
\plot 0 0  0 -.5 /
\plot 7 0  8 -.33 /
\plot -7 0  -8 -.33 /

\put{$\ss \blacksquare$} at 0 2.333 
\put{$\mathbb D(n)$} at -6 6 
\put{$\ss\phi = 1$} at -8 -.8
\put{$\ss-2$} at 0 -.8
\put{$\ss-1/2$} at 7.5 -.8
\endpicture} at 20 0

\put{\beginpicture
\multiput{} at -7 0  7 0  0 7  0 -.6 /
\plot  -7 0  7 0  0 7  -7 0 /
\multiput{$\bullet$} at -5 0  5 0  -6 1  6 1  -1 6  1 6 /
\plot 0 2.333  -5 0 /
\plot 0 2.333   5 0 /
\plot 0 2.333  -6 1 /
\plot 0 2.333   6 1 /
\plot 0 2.333  -1 6 /
\plot 0 2.333   1 6  /
\setdots <.5mm>
\plot  -5 0  -6 -.5 /
\put{$\ss \frac n{n-3}$\strut} at -7 -.8 
\put{$\ss \blacksquare$} at 0 2.333 
\put{$S$} at -4.5 -.5  
\put{$\mathbb H_\ell(n)$} at -6 6 
\endpicture} at 0 -10
\put{\beginpicture
\multiput{} at -7 0  7 0  0 7  0 -.6 /
\plot  -7 0  7 0  0 7  -7 0 /
\multiput{$\bullet$} at -1 1  1 1  -2.5 2.5  2.5 2.5  -2.5 2.5  1.5 3.5  -1.5 3.5 /
\plot 0 2.333   -1.667 0  /
\plot 0 2.333   1.667 0  /
\plot 0 2.333   -4.3 2.667 /
\plot 0 2.333    4.3 2.667  /
\plot 0 2.333  -2.6 4.4 /
\plot 0 2.333   2.6 4.4 /
\put{$\ss \blacksquare$} at 0 2.333 
\setdots <.5mm>
\plot  -1.667 0  -2 -.5 /
\put{$\ss \frac{-2n+6}{n-6}$\strut} at -2 -.85
\put{$\ss E_2^{n-2}$} at -2 1.3  
\put{$\mathbb K_s(n)$} at -6 6 
\put{for $n\ge 7$} at -6 5 
\endpicture} at 20 -10
\endpicture}
$$

\medskip
For $\phi\in \mathbb R\cup\{\infty\}$, let $L_\phi$ be the set of elements of $\mathbb T(n)$
with $p-\frac n3 = \phi\cdot (r-\frac n3).$ This is a line which contains 
$z(n) = (\frac n3,\frac n3)$, thus a central line. We call $\phi$ the {\it slope} of 
$L_\phi$. (It is the slope of the line through $z(n)$ and $(p,r)$
in the cartesian coordinate system with axes $x=r$, $y=p$.) We stress that
{\it the non-central pr-vector $(p,r)$ belongs to $L_\phi$ with $\phi = \frac{3p-n}{3r-n}.$}

The central lines $p=\frac n3,\ r=\frac n3,\ n-q=\frac n3$ in $\mathbb P(n)$ 
have slope $0,\ \infty,\ -1,$ respectively. 
The diagonal lines $p=r,\ r = n-q,\ p = n-q$ in 
$\mathbb D(n)$ have slope $1,\ -2,\ -1/2,$ respectively.

In the lower row, we mention the slope $\frac n{n-3}$
of the central line through $S$ and the slope $\frac{-2n+6}{n-6}$ of the
central line through $E_2^{n-2}.$

\medskip
\begin{lemma}
  \label{lem-ten-three}
  The union $\mathbb P(6)\cup\mathbb D(6)
  \cup \mathbb H_{\ell}(6)$ consists of $18$ central half-lines.
  For $n\ge 7,$ the union
  $\mathbb P(n)\cup\mathbb D(n)\cup \mathbb H_\ell(n) \cup \mathbb K_s(n)$ 
  consists of $24$ central half-lines.
  
  For any $n\ge 6,$ precisely $12$ of these central 
  half-lines form complementary pairs, namely the
  half-lines in $\mathbb P(n)\cup\mathbb D(n)$.
\end{lemma}

\medskip
The proof is easy. In order to show that for $n \ge 7$, no half-line in $\mathbb H_\ell(n)$ is
complementary to a half-line in $\mathbb K_s(n)$, it is sufficient to show that the
central half-line which contains $\D\tau^2 S$ is not complementary to the central
half-line which contains $E_2^{n-2}$. Now $\pr \D\tau^2 S = (n-1,1)$ and $(n-1,1)$ does
not lie on the line $L_\phi$ with $\phi = \frac {-2n+6}{n-6}.$
$\s$

\subsection{The half-line support of $\Cal P(n)$ for $n\ge 6$.}
\label{sec-ten-four}

When dealing with the principal component $\Cal P(n)$, where $n\ge 6,$ we may
use the following reference system: Besides the projective objects $([0],[n])$ and
$([n],[n])$, there are the objects of the form $\tau^iS[\ell]$ 
with $\ell \in \mathbb N_1$ and $0 \le i \le 5$ (see the picture at the beginning of
Section~\ref{sec-ten-six}. 

Note that there is a single object of quasi-length 6 which is neither of the form $X^+$
with $X\in \Cal P(n)$, nor central, namely $S[6]$. The uwb-vector of $S[6]$ is $\frac{6|6}4$
(the central objects of quasi-length 6 have uwb-vector $\frac{6|6}3$).

\medskip
First, let us determine the half-line support of $\Cal P(n)$. 

\medskip
\begin{proposition}
  \label{prop-ten-four}
  The half-line support of $\Cal P(6)$ is 
  $\mathbb P(6)\cup\mathbb D(6)\cup\mathbb H_\ell(6)$. For $n\ge 7,$
  the half-line support of $\Cal P(n)$ 
  is $\mathbb P(n)\cup\mathbb D(n)\cup\mathbb H_\ell(n)\cup\mathbb K_s(n)$.
\end{proposition}

\begin{proof}
  {\it The objects with even quasi-length are supported by $\mathbb D(n)$.}
  In Section~\ref{sec-ten-seven} we will
  show a picture of $\Cal P(n)$; there, the objects on the dashed vertical lines 
  are invariant under $\D$, and all objects of even quasi-length belong to their 
  $\tau^2$-orbits. 
  The two objects of quasi-length 2 which are invariant under $\D$ are $\tau S[2]$
  with pr-vector $(1,1)$ and $\tau^4 S[2]$ with pr-vector $((n-1)/2,(n-1)/2)$,
  thus they belong to complementary half-lines. 

  Now we look at the objects of quasi-length 1,\;3,\;5, and, in addition, at $S[6]$. 
  {\it The objects with quasi-length $1$ have half-line support $\mathbb H_\ell(n)$,}
  since $S$ is one of them. 
  {\it The objects with quasi-length $3$ have half-line support $\mathbb K_s(n)$, if $n\ge 7$,
    and half-line support $\mathbb P(6)$ for $n = 6.$} This follows from the fact that 
  $E_2^{n-2} = \tau^3 S[3]$. 
  {\it The objects with quasi-length $5$ have half-line support $\mathbb P(n)$,} since 
  $\tau^2 S[5], \tau^5 S[5]$ have quasi-length 5 and belong to complementary half-lines in
  $\mathbb P(n)$. Finally, $S[6]$ belongs to $\mathbb D(n).$
  
  Altogether, we see: The half-line support of the class of objects in $\Cal P(n)$ 
  which have quasi-length at most 5
  is $\mathbb P(n)\cup \mathbb D(n)\cup \mathbb H_\ell(n)\cup \mathbb K_s(n)$
  in case $n \ge 7$, and
  is $\mathbb P(n)\cup \mathbb D(n)\cup \mathbb H_\ell(n)$ in case $n = 6$.
  Also, $S[6]$ belongs to a half-line in $\mathbb D(n)$.
  
  It remains to use Theorem~\ref{theoremfive}.
\end{proof}

\subsection{The triangle support of a component and
  proof of Theorem~\ref{theoremsix}.}
\label{sec-ten-five}

For the proof of Theorem~\ref{theoremsix} it remains to deal with
part (b).

\medskip
\begin{proposition}
  \label{prop-ten-five}
  Let $\Cal C$ be an Auslander-Reiten component of $\Cal S(n)$,
  with $n\ge 6.$ 
  The triangle support of $\Cal C$ is the union of the triangles $\Delta_d$ with
  $d\in \Psi(\Cal C),$ where $\Psi(\Cal C)$ is a set of rational numbers $0 \le d < n/3$.
  Moreover, either $\Psi(\Cal C)$ is empty, or else $n/3$ is the only accumulation point of 
  $\Psi(\Cal C)$. Any triangle $\Delta_d$ is the support of only finitely many 
  indecomposable objects of $\Cal C.$
\end{proposition}

\begin{proof}
  Let $\Cal C$ be a component, and let $\beta = \beta\Cal C$. 
  We have seen in Section~\ref{sec-nine}: If $X$ is indecomposable, and $Y = X^{+t}$, then 
$X$ is $u$-minimal if and only if $Y$ is $u$-minimal; also, 
$uY = uX+n\beta t$, $bY = bX+3\beta t$.

\medskip
We call a pair $(u,b)$ of natural numbers 
a {\it primitive pair for} $\Cal C$ provided
there is $X$ in $\Cal C$ which is not central, $u$-minimal, and has quasi-length at most 5,
such that $u = uX$ and $b = bX.$ 

\medskip
(1) {\it For any component $\Cal C$, there are at most $10$ primitive pairs.} 

\medskip
Proof. First of all, the projective indecomposable objects are not $u$-minimal, thus
the $u$-minimal objects in $\Cal C$ of quasi-length at most 5
correspond bijectively to the non-central $\tau^2$-orbits of
objects of quasi-length between 1 and 5
(namely, any $\tau^2$-orbit of non-central objects contains 
precisely one $u$-minimal object). The assertion follows from the fact
that there are at most 10 $\tau^2$-orbits of 
objects of quasi-length between 1 and 5. $\s$

\medskip
If $(u,b)$ is a primitive pair for $\Cal C$, and $\beta = \beta\Cal C$, let
$$
 d_{ub}^{\beta}(t) = (u+n\beta\cdot t)/(b+3\beta\cdot t).
$$

\medskip
Let $\Cal C$ be a component. Let $\Psi(\Cal C)$ be the
set of numbers $dX$ with $X\in \Cal C$ non-central. 

\bigskip
(2) {\it Let $\Cal C$ be a component. The set $\Psi(\Cal C)$ is the set of all 
values $d_{ub}^{\beta\Cal C}(t),$ where $(u,b)$ is a primitive pair for $\Cal C$
and $t\in \mathbb N_0.$}

\medskip
Proof. First, let $(u,b)$ be a primitive pair for $\Cal C$, 
thus there is $X$ in $\Cal C$
with $uX = u$ and $bX = b.$ Let $t\in \mathbb N_0.$ 
According to Theorem~\ref{theoremfive}, we have $d(X^{+t}) = d_{ub}^{\beta\Cal C}(t).$
This shows that the values $d_{ub}^{\beta\Cal C}(t)$ belong to $\Psi(\Cal C).$

Conversely, let  $X\in \Cal C$ be non-central. We have to show that $dX = d_{ub}^{\beta\Cal C}(t)$
for some primitive pair $(u,b)$ and some $t\in \mathbb N_0.$

First, assume that $X = P$ is projective or that $X = S[6].$ 
Then $\Cal C = \Cal P(n)$. Clearly, $dP = 0$ and $d(S[6]) = n/4$ 
(since $u(S[6]) = n,$ and $b(S[6]) = 4$). 
On the other hand, $S$ is $u$-minimal in $\Cal C$ with $uS = 0$ and
$b = 1$, thus $(0,1)$ is a primitive pair for $\Cal C$ and 
$d_{01}^{\beta\Cal C}(0) = 0,$ whereas
$d_{01}^{\beta\Cal C}(1) = n/4.$ This shows that in these three cases, $dX$ has the
required form. 

Next, assume that $X$ has quasi-length between 1 and 5. Then the $\tau^2$-orbit of $X$
contains some object $X'$ which is $u$-minimal and $dX = dX'.$ Let $u = u(X')$, and
$b = b(X').$ Then $d(X') = d_{ub}^{\beta\Cal C}(0)$ has the required form. 

Finally, we can assume that the quasi-length of $X$ is at least 6 and that $X$ 
is different from $S[6].$ Then $X = Y^{+t}$, where $Y$ in $\Cal C$ has quasi-length
at most 5 or is equal to $S[6].$ As we have seen already, $dY = 
d_{ub}^{\beta\Cal C}(0)$ for some primitive pair $(u,b)$, and Theorem~\ref{theoremfive} now yields that
$dX = d(Y^{+t}) = d_{ub}^{\beta\Cal C}(t)$.
$\s$

\medskip
(3) {\it The triangle support of $\Cal C$ consists of the triangles $\Delta_d$ 
with $d \in \Psi(\Cal C)$, where $\Psi(\Cal C)$ is the 
set of rational numbers $d_{ub}^\beta(t) = 
(u+n\beta t)/(b+3\beta t)$, with $(u,b)$ is 
a primitive pair for $\Cal C$, and $t\in \mathbb N_0.$}

\medskip
Proof. This is an immediate consequence of (2). $\s$

\medskip
(4) {\it The sequence  $d_{ub}^\beta(t)$ is strictly increasing and converges to $n/3.$}

\medskip
Proof. This can be checked easily.
$\s$

\medskip
(5) {\it If the set $\Psi(\Cal C)$ is non-empty, then $n/3$ is its only accumulation
point.}

\medskip
Proof. According to (2), $\Psi(\Cal C)$ is the union of the 
sequences $d_{ub}^\beta(t)$, and according to (4), 
these sequences are strictly increasing and converge to
$n/3,$ The set of sequences is indexed by the set of primitive pairs. According to
(1), there are only finitely many primitive pairs. 
Since $\Psi(\Cal C)$ is the union of
finitely many sequences which converge to $n/3$, the assertion follows. 
$\s$.

\medskip
(6) {\it Any standard 
triangle is the support of only finitely many indecomposables in $\Cal C$.}

\medskip
Proof. See Section~\ref{sec-nine-thirteen}. $\s$

\medskip
This completes the proof of Proposition~\ref{prop-ten-five}.
\end{proof}

\begin{proof}[Proof of Theorem~\ref{theoremsix}]
  We have just shown part (b) of Theorem~\ref{theoremsix}.
  Part (a) follows from Corollary~\ref{cor-ten-one-two}
  and Proposition~\ref{prop-ten-four}.  Part (c)
  is a consequence of Proposition~\ref{prop-nine-fourteen}.
\end{proof}

\subsection{The primitive pairs for $\Cal P(n)$.}
\label{sec-ten-six}

\begin{proposition}
  \label{prop-ten-six}
  The primitive pairs for $\Cal P(n)$ are
$
(0,1),\; (1,1),\; (2,2),\; (n-2,3),\; (n-1,3).
$
\end{proposition}

\begin{proof}
  We determine the $u$-minimal objects of quasi-length at most 5. 
  For quasi-length 0, there is no $u$-minimal object.
  For quasi-length 1, there are two $u$-minimal object, namely $S$ and $\tau^3S$,
  both yield the primitive pair $(0,1).$ 
  For quasi-length 2, there are the $u$-minimal objects $\tau^iS[2]$ with $i = 1,2;$
  the corresponding primitive pairs are  $(1,1),\ (2,2),$ respectively.
  There are two $u$-minimal objects of quasi-length 3, namely $\tau^2S[3]$ and $\tau^3S[3]$, both
  with  primitive pair  $(2,2).$
  There are two $u$-minimal objects of quasi-length 4, namely $S[4]$ with primitive pair $(n-2,3)$
  and $\tau^3S[4]$ with primitive pair $(2,2).$ Finally, $S[5]$ and $\tau S[5]$ are the
  $u$-minimal objects of quasi-length 5, both yield the primitive pair $(n-1,3)$.
  Here are the $u$-minimal objects in $\Cal P(n)$ which have quasi-length at most 5:
  
  \medskip

$$
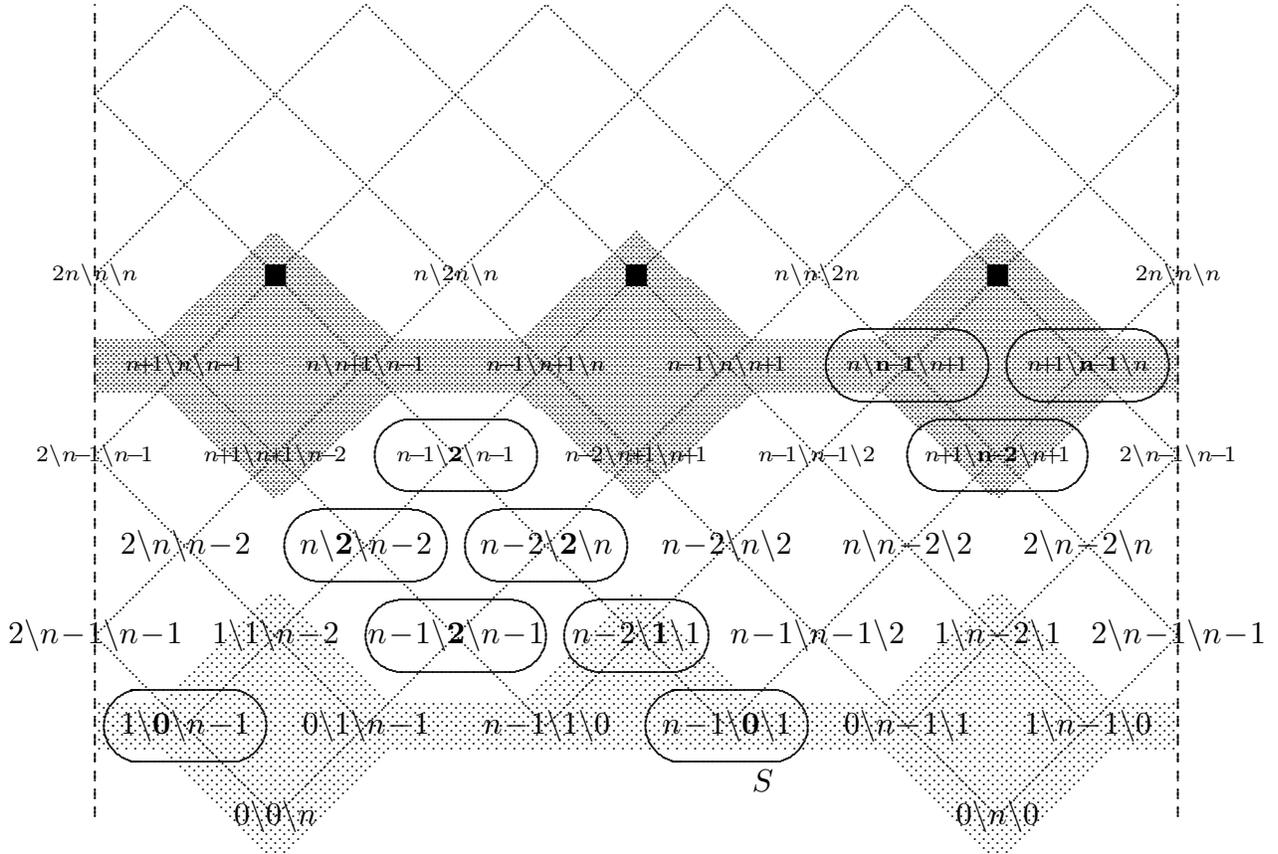
{\beginpicture
    \setcoordinatesystem units <.85cm,1cm>
\multiput{} at 0 -.5  12 8.5  /
\setdashes <1mm>
\plot 0 0  0 9 /
\plot 12 0  12 9 /
\setdots <1mm>
\plot 0 1  5 1 /
\plot 7 1  12  1 /
\setdots <.5mm>
\plot 0 4  3 1 /
\plot 0 2  6 8  12 2 /
\plot 0 2  2 0  10 8  12 6 /
\plot 0 4  4 8  11 1 /
\plot 0 6  2 8  10 0  12 2 /
\plot 0 8  7 1  12 6 /
\plot 1 1  8 8  12 4  9 1 /
\plot 0 6  5 1  12 8 /
\plot 0 8  1 9  2 8  3 9  4 8  5 9  6 8  7 9  8 8  9 9  10 8  11 9  12 8 /  

\put{$\ss 0\bs 0\bs n$} at 2 0
\put{$\ss 0\bs n\bs 0$} at 10 0

\put{$\ss 1\bs \mathbf 0\bs n\!-\!1$} at 1 1
\put{$\ss 0\bs 1\bs n\!-\!1$} at 3 1
\put{$\ss n\!-\!1\bs 1\bs 0$} at 5 1
\put{$\ss n\!-\!1\bs \mathbf 0\bs 1$} at 7 1
\put{$\ss 0\bs n\!-\!1\bs 1$} at 9 1
\put{$\ss 1\bs n\!-\!1\bs 0$} at 11 1

\put{$\ss 2\bs n\!-\!1\bs n\!-\!1$} at 0 2
\put{$\ss 1\bs 1\bs n\!-\!2$} at 2 2
\put{$\ss n\!-\!1\bs \mathbf 2\bs n\!-\!1$} at 4 2
\put{$\ss n\!-\!2\bs \mathbf 1\bs 1$} at 6 2
\put{$\ss n\!-\!1\bs n\!-\!1\bs 2$} at 8 2
\put{$\ss 1\bs n\!-\!2\bs 1$} at 10 2
\put{$\ss 2\bs n\!-\!1\bs n\!-\!1$} at 12 2

\put{$\ss 2\bs n\bs n\!-\!2$} at 1 3
\put{$\ss n\bs \mathbf 2\bs n\!-\!2$} at 3 3
\put{$\ss n\!-\!2\bs \mathbf 2\bs n$} at 5 3
\put{$\ss n\!-\!2\bs n\bs 2$} at 7 3
\put{$\ss n\bs n\!-\!2\bs 2$} at 9 3
\put{$\ss 2\bs n\!-\!2\bs n$} at 11 3

\put{$\ss 2\bs n\!-\!1\bs n\!-\!1$} at 0 4
\put{$\ss n\!+\!1\bs n\!+\!1\bs n\!-\!2$} at 2 4
\put{$\ss n\!-\!1\bs \mathbf 2\bs n\!-\!1$} at 4 4
\put{$\ss n\!-\!2\bs n\!+\!1\bs n\!+\!1$} at 6 4
\put{$\ss n\!-\!1\bs n\!-\!1\bs 2$} at 8 4
\put{$\ss n\!+\!1\bs \mathbf{n\!-\!2}\bs n\!+\!1$} at 10 4
\put{$\ss 2\bs n\!-\!1\bs n\!-\!1$} at 12 4

\put{$\ss n\!+\!1\bs n\bs n\!-\!1$} at 1 5
\put{$\ss n\bs n\!+\!1\bs n\!-\!1$} at 3 5
\put{$\ss n\!-\!1\bs n\!+\!1\bs n$} at 5 5
\put{$\ss n\!-\!1\bs n\bs n\!+\!1$} at 7 5
\put{$\ss n\bs \mathbf{n\!-\!1}\bs n\!+\!1$} at 9 5
\put{$\ss n\!+\!1\bs \mathbf{n\!-\!1}\bs n$} at 11 5

\put{$\ss 2n\bs n\bs n$} at 0 6
\put{$\ss n\bs 2n\bs n$} at 4 6
\put{$\ss n\bs n\bs 2n$} at 8 6
\put{$\ss 2n\bs n\bs n$} at 12 6

\setsolid
\circulararc 180 degrees from 7.5 0.6  center at 7.5 1
\circulararc -180 degrees from 6.5 0.6  center at 6.5 1
\plot 6.4 1.4  7.5 1.4 /
\plot 6.4 0.6  7.5 0.6 /

\put{$S$} at 7.4 0.4 

\circulararc 180 degrees from 1.5 0.6  center at 1.5 1
\circulararc -180 degrees from 0.5 0.6  center at 0.5 1
\plot .4 1.4  1.5 1.4 /
\plot .4 0.6  1.5 0.6 /

\circulararc 180 degrees from 4.6 1.6  center at 4.6 2
\circulararc -180 degrees from 3.4 1.6  center at 3.4 2
\plot 3.4 1.6  4.6 1.6 /
\plot 3.4 2.4  4.6 2.4 /

\circulararc 180 degrees from 6.4 1.6  center at 6.4 2
\circulararc -180 degrees from 5.6 1.6  center at 5.6 2
\plot 5.6 1.6  6.4 1.6 /
\plot 5.6 2.4  6.4 2.4 /

\circulararc 180 degrees from 3.5 2.6  center at 3.5 3
\circulararc -180 degrees from 2.5 2.6  center at 2.5 3
\plot 2.5 2.6  3.5 2.6 /
\plot 2.5 3.4  3.5 3.4 /

\circulararc 180 degrees from 5.5 2.6  center at 5.5 3
\circulararc -180 degrees from 4.5 2.6  center at 4.5 3
\plot 4.5 2.6  5.5 2.6 /
\plot 4.5 3.4  5.5 3.4 /

\circulararc 180 degrees from 4.5 3.6  center at 4.5 4
\circulararc -180 degrees from 3.5 3.6  center at 3.5 4
\plot 3.5 3.6  4.5 3.6 /
\plot 3.5 4.4  4.5 4.4 /

\circulararc 180 degrees from 10.6 3.6  center at 10.6 4
\circulararc -180 degrees from 9.4 3.6  center at 9.4 4
\plot 9.4 3.6  10.6 3.6 /
\plot 9.4 4.4  10.6 4.4 /

\circulararc 180 degrees from 9.5 4.6  center at 9.5 5
\circulararc -180 degrees from 8.5 4.6  center at 8.5 5
\plot 8.5 4.6  9.5 4.6 /
\plot 8.5 5.4  9.5 5.4 /

\circulararc 180 degrees from 11.5 4.6  center at 11.5 5
\circulararc -180 degrees from 10.5 4.6  center at 10.5 5
\plot 10.5 4.6  11.5 4.6 /
\plot 10.5 5.4  11.5 5.4 /

\setshadegrid span <.5mm>
\vshade 0 0.7 1.3 <z,z,,> 0.7  0.7 1.3   <z,z,,> 2  -.5 2.5 <z,z,,>  3.3  .7 1.3 
                   <z,z,,>  4.7  0.7 1.3   <z,z,,> 6  .7 2.5 <z,z,,>  7.3  .7 1.3 
                   <z,z,,>  8.7  0.7 1.3   <z,z,,> 10  -.5 2.5 <z,z,,>  11.3  .7 1.3 
                   <z,z,,> 12 0.7 1.3 /

\setshadegrid span <.3mm>
\vshade 0 4.7 5.3 <z,z,,> 0.7  4.7 5.3   <z,z,,> 2  3.5 6.5 <z,z,,>  3.3  4.7 5.3 
                   <z,z,,>  4.7  4.7 5.3   <z,z,,> 6  3.5 6.5 <z,z,,>  7.3  4.7 5.3 
                   <z,z,,>  8.7  4.7 5.3   <z,z,,> 10  3.5 6.5 <z,z,,>  11.3  4.7 5.3 
                   <z,z,,> 12 4.7 5.3 /

\multiput{$\blacksquare$} at 2 6  6 6  10 6  /
\endpicture}
$$
\end{proof}

\subsection{The $m$-partition and the width partition of $\Cal P(n)$, for $n \ge 6.$}
\label{sec-ten-seven}

We denote by $\Cal P_m$ the class of objects $X$ in $\Cal P(n)$ with $mX = m.$ 
(Recall from Section~\ref{sec-three-six}
that $mX=\min(|\Omega V|,|U|,|V/U|)$ for $X=(U,V)\in\Cal S(n)$.)
We are going to study the $m$-partition in detail. Before we do this, let us 
insert the following remark which concerns Proposition~\ref{prop-nine-fourteen}.

\medskip
\begin{remark}
  The bound $33$ in Proposition~\ref{prop-nine-fourteen}
  is optimal: There are $33$ 
  objects in $\Cal P(7)$ with boundary distance $d = 2$.
\end{remark}

\medskip
Namely consider the objects in
$\Cal P_6,\ \Cal P_8,\ \Cal P_{12},\ \Cal P_{14}$ and $\Cal P_{16}$.
For $n = 7,$ the objects in $\Cal P_m$ with $m=6,\ 8,\ 12,\  14,\ 16$
have width $3,\ 4,\ 6,\ 7,\ 8$ respectively, thus $d = m/b = 2.$ And we have
$|\Cal P_m| = 6,\ 3,\ 3,\ 9,\ 12,$ respectively. 

\medskip
\begin{lemma}
  \label{lem-ten-seven}
  Let $n\ge 6.$ 
  Let $X$ be an object in $\Cal P_m$ with 
  $mX = u+nt$, where $0\le u \le n-1,$ then $u = 0,\ 1,\ 2,\ n-2,\ n-1$ and the width of
  $X$ is $b+3t$ with $b = 1,\ 1,\ 2,\ 3,\ 3$, respectively. Thus, 
  the $m$-partition refines the width-partition of the class of non-central objects in
  $\Cal P(n).$
\end{lemma}

\begin{proof}
  First, assume that $X$ has quasi-length at most $5$ and is $u$-minimal.
  Then $(uX,bX)$ is just one of the 5 primitive pairs $(0,1),\ (1,1),\ (2,2),\ (n-2,3),\
  (n-1,3)$ and the width of $X$ is as asserted. If $X$ is not necessarily $u$-minimal, 
  we look at its $\tau^2$-orbit: Both $mX$ and the width of $X$ are constant
  on the $\tau^2$-orbit.
  
  We also have to look at $X = S[6]$. We have $m(S[6]) = n = 0+nt$ with $t = 1$, and
  $b(S[6]) = 4 = 1+3t,$ as asserted. 
  
  Now any plus-orbit of non-central objects in $\Cal P(n)$ 
  contains either an object of quasi-length at most 5 or else $S[6]$.
  Theorem~\ref{theoremfive} asserts that $m(X^+) = mX+n,$ and $b(X^+) = bX+3.$
\end{proof}

\medskip
Let us denote by $\Cal P\langle b \rangle$ the class of objects in $\Cal P(n)$
with width $b$. Then Lemma~\ref{lem-ten-seven} asserts that for $t\ge 0$, we have
\begin{align*}
  \Cal P\langle 3t+1\rangle &= \Cal P_{nt}\cup \Cal P_{nt+1} \cr
  \Cal P\langle 3t+2\rangle &= \Cal P_{nt+2} \cr
  \Cal P\langle 3t+3\rangle &= \Cal P_{nt+n-2}\cup \Cal P_{nt+n-1}\cup \Cal Z(6t+6) \cr
\end{align*}
where $\Cal Z(\ell)$ is the class of central objects in $\Cal P(n)$ with quasi-length $\ell$.

The following picture shows parts $\Cal P_m$
of the $m$-partition (and the cardinalities $|\Cal P_m|$), as well as (using shading) parts of
the width partition.
$$
{\beginpicture
    \setcoordinatesystem units <.6cm,.7cm>
\multiput{} at -6 0  6 14.5 /

\put{width} at -8 14.5
\put{$b = 1$} at  -8 1
\put{$b = 2$} at  -8 3
\put{$b = 3$} at  -8 5
\put{$b = 4$} at  -8 7
\put{$b = 5$} at  -8 9
\put{$b = 6$} at  -8 11
\put{$b = 7$} at  -8 13

\setdashes <2mm>
\plot -6 0  -6 14.5 /
\plot  6 0   6 14.5 /
\plot  0 0   0 14.5 /

\setdots <.5mm>
\plot -6 2  -4 0  6 10  4 12 /
\plot -6 4  -3 1   /
\plot -6 6  -1 1  6 8  2 12 /
\plot -6 8  1 1   6 6   0 12 /
\plot -6 10  4 0  6 2 /
\plot -6 12  5 1 /
\plot 3 1  6 4  -2 12 /
\plot -4 12  6 2 /
\plot -5 1  6 12 /
\plot -6 2  4 12 /
\plot -6 4  2 12 /
\plot -6 6  0 12 /
\plot -6 8  -2 12 /
\plot -6 10  -4 12 /
\plot -6 12  -4.5 13.5 /
\plot -5.5 13.5  -4 12  -2.5 13.5 /
\plot -3.5 13.5  -2 12  -.5 13.5 /
\plot -6 12  -4.5 13.5 /
\plot  5.5 13.5   4 12  2.5 13.5 /
\plot  3.5 13.5   2 12  .5 13.5 /
\plot  1.5 13.5   0 12  -1.5 13.5 /
\plot  6 12  4.5  13.5 /
\setshadegrid span <.5mm>
\vshade -6 0.7 1.3 <z,z,,> -5.3  0.7 1.3   <z,z,,> -4  -.5 2.5 <z,z,,>  -2.7  .7 1.3 
                   <z,z,,>  -1.3  0.7 1.3   <z,z,,> 0  .7 2.5 <z,z,,>  1.3  .7 1.3 
                   <z,z,,>  2.7  0.7 1.3   <z,z,,> 4  -.5 2.5 <z,z,,>  5.3  .7 1.3 
                   <z,z,,> 6 0.7 1.3 /
\vshade -6 5.5 8.5 <z,z,,> -4.7  6.7 7.3  <z,z,,> -3.3  6.7 7.3  <z,z,,>  -2 5.5 8.5 
                   <z,z,,> -0.7  6.7 7.3  <z,z,,> 0.7  6.7 7.3  <z,z,,>  2 5.5 8.5 
                    <z,z,,> 3.3  6.7 7.3  <z,z,,> 4.7  6.7 7.3  <z,z,,>  6 5.5 8.5 /
\vshade -6 12.7 13.3 <z,z,,> -5.3  12.7 13.3   <z,z,,> -4  11.5 14.5 <z,z,,>  -2.7  12.7 13.3 
                   <z,z,,>  -1.3  12.7 13.3   <z,z,,> 0  11.5 14.5 <z,z,,>  1.3  12.7 13.3 
                   <z,z,,>  2.7  12.7 13.3   <z,z,,> 4  11.5 14.5 <z,z,,>  5.3  12.7 13.3 
                   <z,z,,> 6 12.7 13.3 /

\vshade -9 .2 1.8 <,,,> -7 .2 1.8 /
\vshade -9 6.2 7.8 <,,,> -7 6.2 7.8 /
\vshade -9 12.2 13.8 <,,,> -7 12.2 13.8 /

\setshadegrid span <.3mm>
\vshade -6 4.7 5.3 <z,z,,> -5.3  4.7 5.3   <z,z,,> -4  3.5 6.5 <z,z,,>  -2.7  4.7 5.3 
                   <z,z,,>  -1.3  4.7 5.3   <z,z,,> 0  3.5 6.5 <z,z,,>  1.3  4.7 5.3 
                   <z,z,,>  2.7  4.7 5.3   <z,z,,> 4  3.5 6.5 <z,z,,>  5.3  4.7 5.3 
                   <z,z,,> 6 4.7 5.3 /
\vshade -6 9.5 12.5 <z,z,,> -4.7  10.7 11.3  <z,z,,> -3.3  10.7 11.3  <z,z,,>  -2 9.5 12.5
                   <z,z,,> -0.7  10.7 11.3  <z,z,,> 0.7  10.7 11.3  <z,z,,>  2 9.5 12.5 
                    <z,z,,> 3.3  10.7 11.3  <z,z,,> 4.7  10.7 11.3  <z,z,,>  6 9.5 12.5 /

\vshade -9 4.2 5.8 <,,,> -7 4.2 5.8 /
\vshade -9 10.2 11.8 <,,,> -7 10.2 11.8 /

\multiput{$\blacksquare$} at -4 6  0 6  4 6  -6 12  -2 12  2 12  6 12 /
\setdashes <1mm>
\setdashes <.5mm>
\plot -6 1  -5 1 /
\plot -3 1  -1 1 /
\plot 5 1  6 1 /
\plot 3 1  1 1 /

\put{$S$} at 1 1

\put{$\Cal P_m$} [l] at 8 14.5
\put{$\Cal P_{2n}$} [l] at 8 12.75
\put{$\Cal P_{2n-1}$} [l] at 8 11
\put{$\Cal P_{2n-2}$} [l] at 8 10
\put{$\Cal P_{n+2}$} [l] at 8 9
\put{$\Cal P_{n+1}$} [l] at 8 8
\put{$\Cal P_n$} [l] at 8 6.75
\put{$\Cal P_{n-1}$} [l] at 8 5
\put{$\Cal P_{n-2}$} [l] at 8 4
\put{$\Cal P_2$} [l] at 8 3
\put{$\Cal P_1$} [l] at 8 2
\put{$\Cal P_0$} [l] at 8 0.75

\put{$|\Cal P_m|$} at 10 14.5

\put{$9$}  at 10 12.75
\put{$6$} at 10 11
\put{$3$}  at 10 10
\put{$12$} at 10 9
\put{$3$}  at 10 8
\put{$9$}  at 10 6.75
\put{$6$}  at 10 5
\put{$3$}  at 10 4
\put{$12$} at 10 3
\put{$3$}  at 10 2
\put{$8$}  at 10 0.75

\setsolid
\plot -6 1.5  10.5 1.5 /
\plot -6 4.5  10.5 4.5 /
\plot -6 5.5  10.5 5.5 /
\plot -6 7.5  10.5 7.5 /
\plot -6 10.5  10.5 10.5 /
\plot -6 11.5  10.5 11.5 /
\plot -6 13.5  10.5 13.5 /

\circulararc 360 degrees from -3.5 2 center at -4 2
\circulararc 360 degrees from   .5 2 center at  0 2
\circulararc 360 degrees from  3.5 2 center at  4 2

\circulararc 360 degrees from -3.5 4 center at -4 4
\circulararc 360 degrees from   .5 4 center at  0 4
\circulararc 360 degrees from  3.5 4 center at  4 4

\circulararc 180 degrees from -6 7.5 center at -6 8
\circulararc 360 degrees from  -2.5 8 center at  -2 8
\circulararc 360 degrees from  2.5 8 center at  2 8
\circulararc -180 degrees from  6 7.5 center at  6 8

\circulararc 180 degrees from -6 9.5 center at -6 10
\circulararc 360 degrees from  -2.5 10 center at  -2 10
\circulararc 360 degrees from  2.5 10 center at  2 10
\circulararc -180 degrees from  6 9.5 center at  6 10

\circulararc 360 degrees from -3.5 14 center at -4 14
\circulararc 360 degrees from   .5 14 center at  0 14
\circulararc 360 degrees from  3.5 14 center at  4 14

\endpicture}
$$
The left scale describes the width of the various shaded areas, the right
scale names the parts $\Cal P_m$
(the objects in $\Cal P_m$ with $m = 1,\ n-2,\ n+1,\ 2n-2,\ 2n+1$
are encircled); the black squares $\blacksquare$
are the central objects, they do not belong to
any class $\Cal P_m.$ 

\bigskip
Let us return to the triangle support of $\Cal P(n)$. As we know, any class $\Cal P_m$
lives on a single standard triangle, thus the triangle support of the classes
$\Cal P\langle b\rangle$ with $b \equiv 1,2,3\ \mod n$ consists of 2,\ 1, 2, 
standard triangles, respectively:

Here are the corresponding pictures for all $b\ge 1$ 
(for $b = 1$, one has to delete the small circle at the left
lower corner):
$$
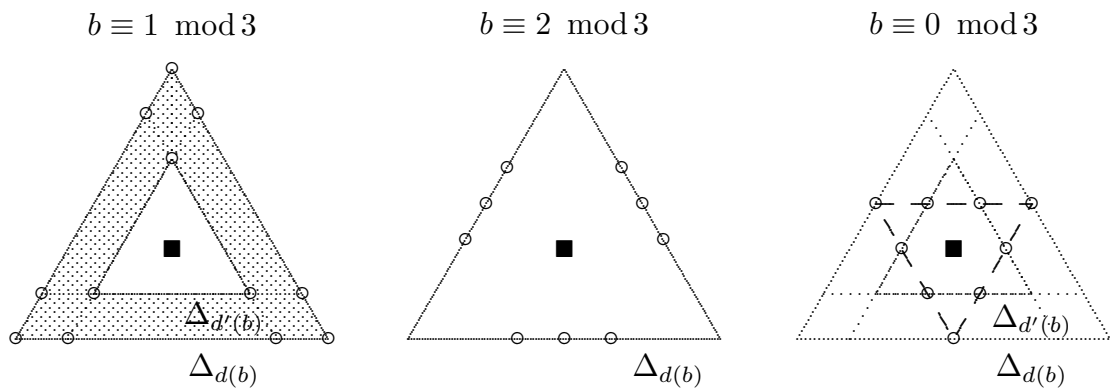
{\beginpicture
   \setcoordinatesystem units <.288675cm,.5cm>

\put{\beginpicture
\multiput{} at 0 6  6 0  0 5 /
\put{$b \equiv 1 \ \mod 3$} at 0 7
\setdots <.3mm>
\plot -6 0  0 6  6 0  -6 0 /
\plot -3 1  3 1  0 4  -3 1 /
\setdots <1mm>
\plot -4 0  1 5 /
\plot 4 0  -1 5 /
\plot -5 1  5 1 /
\put{$\circ$} at -6 0 
\multiput{$\circ$} at -5 1  -3 1  -4 0 
         -1 5  0 6  0 4  1 5 
         3 1  5 1  4 0  6 0 /
\put{$\ss \blacksquare$} at 0 2
\setdots <.5mm>
\put{$\Delta_{d'(b)}$} at 2 0.45 
\put{$\Delta_{d(b)}$} at 2 -.7  
\setshadegrid span <.5mm>
\vshade -6 0 0 <,z,,> -3 0 3 <z,z,,> -3.001 1 3 <z,z,,> 0 4 6  <z,z,,> 2.999 1 3 <z,z,,> 
         3 0 3 <z,z,,>  6 0 0 /
\vshade -3 0 1 <z,z,,> 3 0 1 / 
\endpicture} at 0 0 
\put{\beginpicture
\multiput{} at 0 6  6 0  0 5 /
\put{$b \equiv 2 \ \mod 3$} at 0 7

\setdots <.3mm>
\plot -6 0  0 6  6 0  -6 0 /

\setdots <1mm>
\multiput{$\bullet$} at /
\put{$\ss \blacksquare$} at 0 2
\put{$\Delta_{d(b)}$} at 4 -.7  

\multiput{$\circ$} at -3.8 2.2  3.8 2.2 -3 3  -2.2 3.8  3 3   2.2 3.8 
     -1.8 0  0 0  1.8 0  /

\setdashes <2mm>
\endpicture} at 16 0 

\put{\beginpicture
\multiput{} at -6 6  6 0  0 5 /
\put{$b \equiv 0 \ \mod 3$} at 0 7
\setdots <.5mm>
\plot -6 0  0 6  6 0  -6 0 /
\plot -3 1  3 1  0 4  -3 1 /
\setdots <1mm>
\plot -4 0  1 5 /
\plot 4 0  -1 5 /
\plot -5 1  5 1 /
\put{$\ss \blacksquare$} at 0 2
\put{$\Delta_{d'(b)}$} at 3 0.4 
\put{$\Delta_{d(b)}$} at 3 -.7  

\multiput{$\circ$} at -3 3  3 3  0 0  -1 3  1 3  -2 2  2 2  -1 1  1 1   /
\setshadegrid span <.5mm>

\setdashes <2mm>
\plot -3 3  3 3  0 0  -3 3 /
 
\endpicture} at 30 0 

\endpicture}
$$
Here we write $d(3t+1) = \frac{nt}{3t+1},\ d'(3t+1) = \frac{nt+1}{3t+1},\
d(3t+2) = \frac{nt+2}{3t+2},\ d(3t+3) = \frac{n(t+1)-2}{3t+3}\;\text{and}\;\
d'(3t+3) = \frac{n(t+1)-1}{3t+3}.$

\bigskip
{\bf Nablas.} The last picture suggests to look not only at standard triangles, but also at
costandard ones (``nablas''): The {\it costandard triangle} $\nabla_d$ is
defined for $n/3 < d \le n$ as follows: 
one of its sides is given by the line $p = d$, the remaining two
are obtained using the rotations $\rho$ and $\rho^2$. 
One checks easily that for $b = 3t+3$, the class $\Cal P\langle b\rangle$ lies on the
costandard triangle $\nabla_{d''(b)}$, with $d'' = n/3+1/b.$

\bigskip
{\bf Circles.} The objects $X$ in $\Cal P(n)$ with fixed odd quasi-length have equal width: If 
the quasi-length of $X$ is $2b-1,$ then $bX = b.$ It follows that the $\tau$-orbit of $X$ 
lies on a circle (the dashed circles in the following picture).
The remaining objects are all supported by $\mathbb D(n)$.
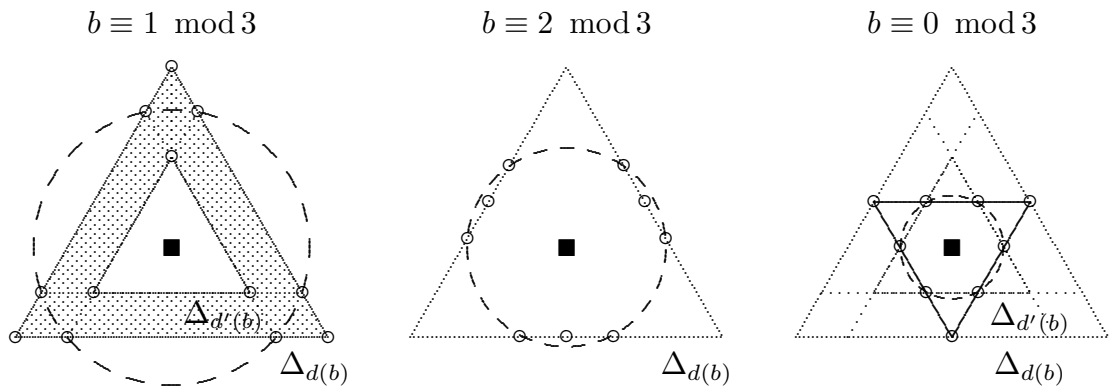
$$
{\beginpicture
   \setcoordinatesystem units <.288675cm,.5cm>

\put{\beginpicture
\multiput{} at 0 6  6 0  0 5 /
\put{$b \equiv 1 \ \mod 3$} at 0 7
\setdots <.3mm>
\plot -6 0  0 6  6 0  -6 0 /
\plot -3 1  3 1  0 4  -3 1 /
\setdots <1mm>
\plot -4 0  1 5 /
\plot 4 0  -1 5 /
\plot -5 1  5 1 /
\put{$\circ$} at -6 0 
\multiput{$\circ$} at -5 1  -3 1  -4 0 
         -1 5  0 6  0 4  1 5 
         3 1  5 1  4 0  6 0 /
\put{$\ss \blacksquare$} at 0 2
\setdots <.5mm>
\put{$\Delta_{d'(b)}$} at 2 0.45 
\put{$\Delta_{d(b)}$} at 5.5 -.7  
\setshadegrid span <.5mm>
\vshade -6 0 0 <,z,,> -3 0 3 <z,z,,> -3.001 1 3 <z,z,,> 0 4 6  <z,z,,> 2.999 1 3 <z,z,,> 
         3 0 3 <z,z,,>  6 0 0 /
\vshade -3 0 1 <z,z,,> 3 0 1 / 
\setdashes <2mm>
\circulararc 360 degrees from 5 1 center at 0 2
\endpicture} at 0 0 
\put{\beginpicture
\multiput{} at 0 6  6 0  0 5 /
\put{$b \equiv 2 \ \mod 3$} at 0 7

\setdots <.5mm>
\plot -6 0  0 6  6 0  -6 0 /

\setdots <1mm>
\multiput{$\bullet$} at /
\put{$\ss \blacksquare$} at 0 2
\put{$\Delta_{d(b)}$} at 5 -.7  

\multiput{$\circ$} at -3.8 2.2  3.8 2.2 -3 3  -2.2 3.8  3 3   2.2 3.8 
     -1.8 0  0 0  1.8 0  /

\setdashes <1.5mm>
\circulararc 360 degrees from  2.2 3.8  center at 0 2
\endpicture} at 16 0 

\put{\beginpicture
\multiput{} at 0 6  6 0  0 5 /
\put{$b \equiv 0 \ \mod 3$} at 0 7
\setdots <.5mm>
\plot -6 0  0 6  6 0  -6 0 /
\plot -3 1  3 1  0 4  -3 1 /
\setdots <1mm>
\plot -4 0  1 5 /
\plot 4 0  -1 5 /
\plot -5 1  5 1 /
\put{$\ss \blacksquare$} at 0 2
\put{$\Delta_{d'(b)}$} at 3 0.4 
\put{$\Delta_{d(b)}$} at 3 -.7  

\multiput{$\circ$} at -3 3  3 3  0 0  -1 3  1 3  -2 2  2 2  -1 1  1 1   /
\setshadegrid span <.5mm>

\setsolid
\plot -3 3  3 3  0 0  -3 3 /

\setdashes <1mm>
\circulararc 360 degrees from 0 3.14 center at 0 2
 
\endpicture} at 31 0 

\endpicture}
$$

However, one should be aware that these circles just 
indicate that the {\bf set} of the corresponding $\tau$-orbit is contained in the circle,
whereas the circle usually does {\bf not} describe the action of $\tau$ on the orbit! 
The graph which describes the action of $\tau$ is topologically also a circle, 
but its embedding into the pr-triangle $\mathbb T(n)$ may be different, namely it looks as follows:
$$
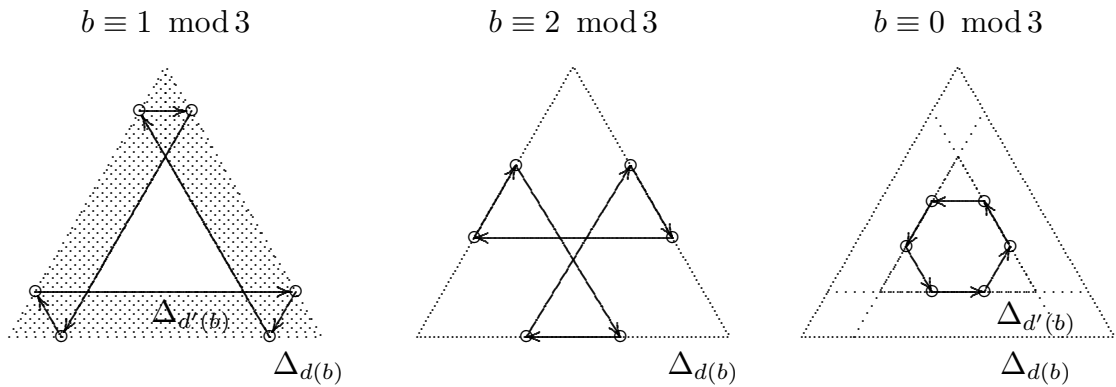
{\beginpicture
   \setcoordinatesystem units <.288675cm,.5cm>

\put{\beginpicture
\multiput{} at 0 6  6 0  0 5 /
\put{$b \equiv 1 \ \mod 3$} at 0 7
\setdots <1mm>
\plot -6 0  0 6  6 0  -6 0 /
\plot -3 1  3 1  0 4  -3 1 /
\plot -4 0  1 5 /
\plot 4 0  -1 5 /
\plot -5 1  5 1 /
\multiput{$\circ$} at -5 1   
         -1 5     1 5 
         5 1   -4 0  4 0   /
\setdots <.5mm>
\put{$\Delta_{d'(b)}$} at 1 0.45 
\put{$\Delta_{d(b)}$} at 5.5 -.7  
\setshadegrid span <.5mm>
\vshade -6 0 0 <,z,,> -3 0 3 <z,z,,> -3.001 1 3 <z,z,,> 0 4 6  <z,z,,> 2.999 1 3 <z,z,,> 
         3 0 3 <z,z,,>  6 0 0 /
\vshade -3 0 1 <z,z,,> 3 0 1 / 

\setsolid
\arr{-4 0}{-4.9 .9}
\arr{-5 1}{4.8 1}
\arr{5 1}{4.1 0.1}
\arr{4 0}{-.9 4.9}
\arr{-1 5}{.8 5}
\arr{1 5}{-3.9 0.1}

\endpicture} at 0 0 
\put{\beginpicture
\multiput{} at 0 6  6 0  0 5 /
\put{$b \equiv 2 \ \mod 3$} at 0 7

\setdots <.5mm>
\plot -6 0  0 6  6 0  -6 0 /

\setdots <1mm>
\multiput{$\bullet$} at /
\put{$\Delta_{d(b)}$} at 5 -.7  

\multiput{$\circ$} at -3.8 2.2  3.8 2.2  -2.2 3.8     2.2 3.8 
     -1.8 0   1.8 0  /

\setsolid
\arr{1.8 0}{-1.8 0}
\arr{-1.8 0}{2.2 3.8}
\arr{2.2 3.8}{3.8 2.2}
\arr{3.8 2.2}{-3.8 2.2}
\arr{-3.8 2.2}{-2.2 3.8}
\arr{-2.2 3.8}{1.8 0}
\endpicture} at 16 0 

\put{\beginpicture
\multiput{} at 0 6  6 0  0 5 /
\put{$b \equiv 0 \ \mod 3$} at 0 7
\setdots <.5mm>
\plot -6 0  0 6  6 0  -6 0 /
\plot -3 1  3 1  0 4  -3 1 /
\setdots <1mm>
\plot -4 0  1 5 /
\plot 4 0  -1 5 /
\plot -5 1  5 1 /
\put{$\Delta_{d'(b)}$} at 3 0.4 
\put{$\Delta_{d(b)}$} at 3 -.7  

\multiput{$\circ$} at   -1 3  1 3  -2 2  2 2  -1 1  1 1   /
\setshadegrid span <.5mm>

\setsolid
\arr{-1 1}{1 1}
\arr{1 1}{2 2}
\arr{2 2}{1 3}
\arr{1 3}{-1 3}
\arr{-1 3}{-2 2}
\arr{-2 2}{-1 1}

\endpicture} at 31 0 

\endpicture}
$$

\medskip
{\bf Summary.} 
The following picture presents the
\phantomsection{pr-vectors}
\addcontentsline{lof}{subsection}{The small pr-vectors of $\Cal P(6)$.}%
of all the objects in $\Cal P\langle b\rangle$ with $b \le 4.$ 
$$
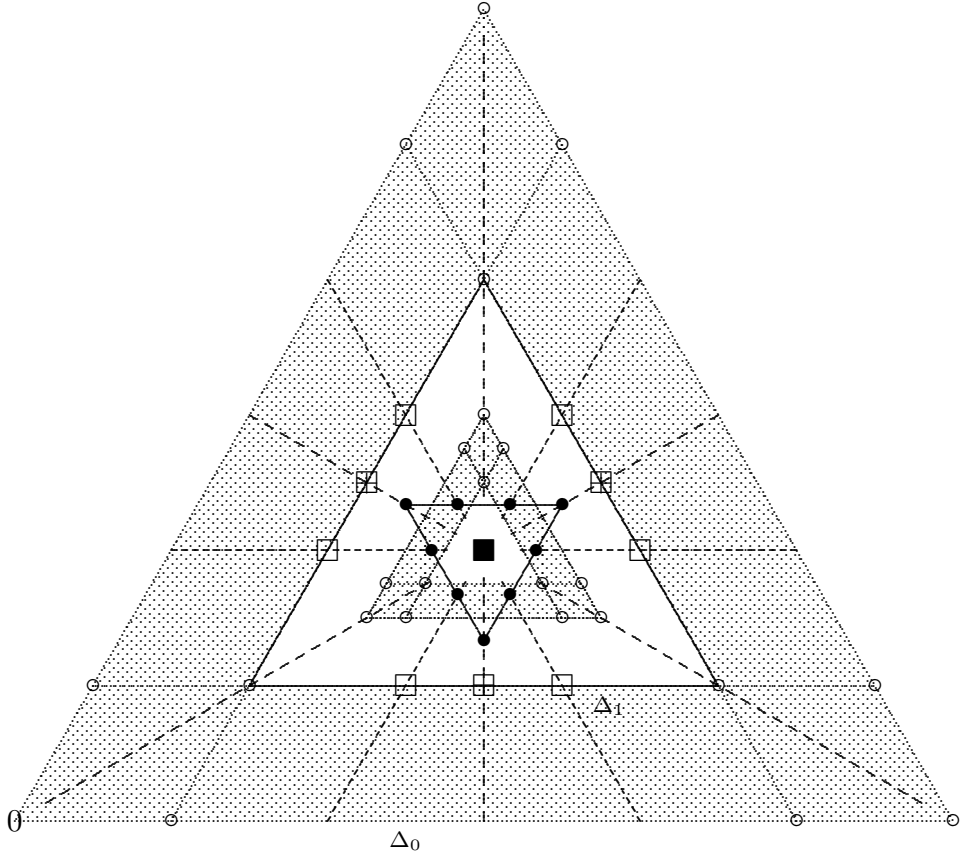
{\beginpicture
   \setcoordinatesystem units <.866025cm,1.5cm>
\multiput{} at -6 -.4  6 6 /

\setdots <1mm>

\multiput{$\bullet$} at -0.333 1.667  0.333 1.667 
                        -0.667 2      .667  2 
                        -0.333 2.333  0.333 2.333  /   
\multiput{$\square$} at -1 1  1 1
                        -2 2  2 2
                        -1 3  1 3 / 
\multiput{$\circ$} at -4 0  4 0
                        -5 1  5 1
                        -1 5  1 5 
                        0 6  6 0 / 
\multiput{$\circ$} at  -1 1.5  1 1.5 
                       -1.25 1.75   1.25 1.75    
                          -.25 2.75  .25 2.75                    
                         / 

\multiput{$\square$} at 0 1  -1.5 2.5  1.5 2.5  /
\multiput{$\circ$} at    0 4  -3 1  3 1 /

\multiput{$+$} at 0 1  -1.5 2.5  1.5 2.5 /
\multiput{$\bullet$} at 0 1.333   -1 2.333  1 2.333  /

\multiput{$\circ$} at 0 3  -1.5 1,5  1.5 1.5 /  
\multiput{$\circ$} at 0 2.5  -0.75 1.75  0.75 1.75  /  

\setdots <.5mm> 
\plot -6 0  6 0  0 6  -6 0 /    
\plot -4 0  1 5 /
\plot  4 0  -1 5 /
\plot -5 1  5 1 /     

\plot 0 3  -1.5 1,5  1.5 1.5  0 3 /   
\plot -1.25 1.75  1.25 1.75 /
\plot  -.25 2.75 1 1.5 /
\plot   -1 1.5   .25 2.75 /
 
\setsolid

\plot -3 1  3 1  0 4  -3 1 /    
\plot  0 1.333   -1 2.333  1 2.333   0 1.333 /   
                                       
\put{$\blacksquare$} at 0 2

\setdots <.8mm>  

\setdashes <1mm> 
\plot 0 0    0 1.8   /
\plot 0 5.8  0 2.4 /
\plot -5.6 0.1333  -.6 1.8 /
\plot  5.6 0.1333   .6 1.8 /
\plot    3 3  0.35 2.12 /
\plot   -3 3  -.35 2.12 /

\setdashes <.6mm>  

\plot  -2 0  -.2 1.8 /
\plot   2 0   .2 1.8 /
\plot   2 4  0.2 2.2 /
\plot  -2 4  -.2 2.2 /
\plot -4 2  -.4 2 /
\plot  4 2   .4 2 /

\setdashes <1.5mm>  


\setsolid
\put{$0$} at -6 0  

\put{$\ss \Delta_0$} at -1 -.15
\put{$\ss \Delta_1$} at 1.6 .85

\setshadegrid span <.5mm>
\vshade -6 0 0 <,z,,> -3 0 3 <z,z,,> -3.001 1 3 <z,z,,> 0 4 6  <z,z,,> 2.999 1 3 <z,z,,> 
         3 0 3 <z,z,,>  6 0 0 /
\vshade -3 0 1 <z,z,,> 3 0 1 / 
 
\vshade  -1.5 1.5 1.5 <,z,,>  -0.75 1.5 2.25 
                      <z,z,,> -0.7701 1.75 2.25 <z,z,,> 0 2.5  3 <z,z,,>
                      0.7699 1.75 2.25 <z,z,,> .75 1.5 2.25  <z,z,,> 1.5 1.5 1.5  /
\vshade  -0.75 1.5 1.75  <z,z,,>  0.75 1.5 1.75  / 

\setdots <.4mm> 
\plot 0 3   1.5 1.5  -1.5 1.5  0 3 /   
\plot  -.25 2.75  1 1.5 /
\plot   -1 1.5   .25 2.75 /

\endpicture}
$$
There is 
$\Cal P\langle1\rangle$, it is part of the outer dotted region with boundary
$\Delta_{d(1)}\cup\Delta_{d'(1)} = \Delta_0\cup \Delta_1;$ the vertices are marked by small circles. ---
Then there is $\Cal P \langle2\rangle,$ it is part of the solid triangle 
$\Delta_{d(2)} = \Delta_1$
(note that $\Delta_1$ plays a role both for $\Cal P\langle 1\rangle$ and $\Cal P\langle 2\rangle$).
The vertices for $\Cal P\langle 2\rangle$ are drawn as square boxes; the boxes on the
diagonal lines are endowed with a plus sign in order to indicate that they are the pr-vectors of two isomorphism classes.
---
Then there is 
$\Cal P\langle3\rangle$, lying on the (solid) costandard triangle $\nabla_{d''(3)}$,
marked by small bullets. 
---
Finally, $\Cal P\langle4\rangle$ is part of the inner dotted region with boundary
$\Delta_{d(4)}\cup\Delta_{d'(4)};$ the vertices are marked again by small circles. 
---
Next, one should insert 
$\Cal P\langle5\rangle$ as part of the triangle $\Delta_{d(5)}$, then 
$\Cal P\langle6\rangle$ as part of $\nabla_{d''(6)}$, then  
$\Cal P\langle7\rangle$, and so on.

We also have indicated the central lines in the half-line support of $\Cal P(n)$, namely the 
lines in $\mathbb P(n)\cup\mathbb D(n)$ (but not the remaining half-lines). 
In the diagram, note that for $n>6$ the six objects of quasi-length $3$ in $\Cal P\langle2\rangle$ lie on
$\mathbb K_s(n)$, not on $\mathbb P(n)$, see Sections \ref{sec-ten-three}, \ref{sec-ten-four}. 

\subsection{The central objects in $\Cal P(n)$.}
\label{sec-ten-eight}

Let us finish our study of the principal components $\Cal P(n)$ with $n\ge 6$ 
by looking at the
central objects in $\Cal P(n)$. 
The central objects have quasi-length divisible by $6$, thus, let us
add  a description of the objects in 
$\Cal P(n)$ of quasi-length $6t$ with 
$t\in \mathbb N_0.$ For $t = 0,$ we deal with the two
indecomposable projective objects; they belong to $\Cal P\langle 1\rangle.$ 
For any $t\ge 1$, there are 3 central objects, they
have uwb-vector $(nt,nt;3t)$;
the remaining 3 objects belong to $\Cal P\langle 3t+1 \rangle$; they
have uwb-vectors
$(nt,nt,3t+1),\ (nt,nt+n,3t+1)$ and $(nt+n,nt,3t+1)$; they are 
supported by $\mathbb D(n)$, and 
their pr-vectors are the corners of the triangle $\Delta_{nt/(3t+1)}.$ 

\subsection{The half-line support of $\Cal P(n)$, for $1 \le n \le 5.$}
\label{sec-ten-nine}

Up to now, we have considered the cases $n\ge 6.$ Of course, one may also ask
for the half-line support of $\Cal P(n)$ for $n\le 5$ (considering, as we did up to now,
{\bf central} half-lines). Here is the statement for $n = 5.$

\medskip
\begin{proposition}
  \label{prop-ten-nine}
  The number of central half-lines which contain the pr-vector of a
  non-central object in $\Cal P(5)$ is $30$. First of all, there are the $12$
  half-lines in $\mathbb P(5)\cup\mathbb D(5)$. The additional
  $18$ half-lines pass through the $\tau$-orbits of $S$,
  $E_2^3$  and $(0,[2]).$
  The half-lines which pass through the objects in the $\tau$-orbit
  of $E_2^3$ are complementary to the half-lines which pass through the objects 
  in the $\tau$-orbit of $(0,[2])$.
\end{proposition}

\begin{proof}
  Parts of $\Cal P(5)$ can be constructed 
  in the same way as $\Cal P(n)$ is constructed,
  for $n\ge 6$. The indecomposable projective objects are said to have quasi-length 0.
  Then we consider the $\tau$-orbit
of $S = (0,[1])$ and we call these the objects of quasi-length 1. 
The non-projective neighbors of objects of quasi-length 1 are said to have
quasi-length 2. We
proceed in this way in order to construct the objects of 
quasi-length 2 to 5. When we try to 
construct objects of quasi-length
 6, we obtain no longer indecomposable objects: What we get are
two $\tau$-orbits, namely the $\tau$-orbit of $([1],[3])$, 
and the $\tau$-orbit of $E_2^2;$ both $\tau$-orbits live on the diagonal lines.
Finally, there is a further $\tau$-orbit, namely the $\tau$-orbit of $(0,[2]).$ 

As we see: Most of the indecomposables live on the diagonal lines and on the 
central lines parallel to the boundary, namely (as in the cases $n\ge 7$) the 
objects of quasi-length $0,\ 2,\ 4,$ and $5$, but in addition also those in 
two of the three new orbits. 
The only exceptions are  the  $\tau$-orbit of $S$, 
the objects of quasi-length 3 (this is the $\tau$-orbit of $E_2^3$), 
as well as 
one of the new orbits, namely the orbit of $(0,[2])$. 
In order to see that we obtain in this way 18 additional 
half-lines, we look at the triangles $\mathbb F_1$ and $\mathbb F_2$ of vectors $(p,r)$ with 
$p \le r \le \frac13n$ and $\frac13n\le r \le \frac12(n-p),$  respectively.
The object $S$ lies in the interior of $\mathbb F_1$ (since $0 < 1 < \frac 53$, we have
$p < r < \frac13n$); the object $E_2^3$ also lies in the interior of $\mathbb F_1$ 
(since $1 <\frac 32 < \frac 53$, we again have $p < r < \frac13n$). 
Note that $S$ and $E_2^3$  live on different central half-lines (since there is no $a$ such that
$(1-a)(\frac 53,\frac 53) + a(0,1) = (1,\frac 32)$). 
The object $(0,[2])$ lives in the interior of $\mathbb F_2$ 
(since $\frac 53 < 2 < \frac 52$ we
have $\frac 13 n < r < \frac12(n-p)$).
Altogether
we obtain three different central half-lines inside $\mathbb F_1\cup \mathbb F_2.$

It remains to observe that $\tau^5 (0,[2]) = ([3],[5])$ and that the central half-lines
passing through $E_2^3$ and $([3],[5])$ are complementary.
\end{proof}

\medskip
{\bf Summary.} 
Let $g$ be the number of central lines with
\phantomsection{non-central}
\addcontentsline{lof}{subsection}{The number of central half lines and lines for $\Cal P(n)$.}%
objects of $\Cal P(n)$ living
on both of its central half-lines (for $n\neq 5,$ these lines 
turn out to be diagonal lines or central lines
parallel to the boundary, thus lines in $\mathbb P(n)\cup \mathbb D(n)$). Let 
$h$ be the number of half-lines in the half-line support of $\Cal P(n)$ such that the
complementary half-line does not belong to the half-line support. 
Thus $g+h$ is the number of central {\bf lines} which pass through 
non-central objects in $\Cal P(n)$, whereas
$2g+h$ is the number of central {\bf half-lines} which pass through 
non-central objects in $\Cal P(n)$.
$$
{\beginpicture
   \setcoordinatesystem units <1cm,.5cm>
\put{$n$} at -.7 0
\put{$1$} at 1 0
\put{$2$} at 2 0
\put{$3$} at 3 0
\put{$4$} at 4 0
\put{$5$} at 5 0
\put{$6$} at 6 0
\put{$7$} at 7 0
\put{$8$} at 8 0
\put{$9$} at 9 0
\put{$\cdots$} at 10 0

\put{$g$} at -.7 -1
\put{$0$} at 1 -1
\put{$2$} at 2 -1
\put{$3$} at 3 -1
\put{$3$} at 4 -1
\put{$12$} at 5 -1
\put{$6$} at 6 -1
\put{$6$} at 7 -1
\put{$6$} at 8 -1
\put{$6$} at 9 -1
\put{$\cdots$} at 10 -1

\put{$h$} at -.7 -2
\put{$2$} at 1 -2
\put{$1$} at 2 -2
\put{$2$} at 3 -2
\put{$6$} at 4 -2
\put{$6$} at 5 -2
\put{$6$} at 6 -2
\put{$12$} at 7 -2
\put{$12$} at 8 -2
\put{$12$} at 9 -2
\put{$\cdots$} at 10 -2

\put{$g+h$} at -.7 -3
\put{$2$} at 1 -3
\put{$3$} at 2 -3
\put{$5$} at 3 -3
\put{$9$} at 4 -3
\put{$18$} at 5 -3
\put{$12$} at 6 -3
\put{$18$} at 7 -3
\put{$18$} at 8 -3
\put{$18$} at 9 -3
\put{$\cdots$} at 10 -3

\put{$2g+h$} at -.7 -4
\put{$2$} at 1 -4
\put{$5$} at 2 -4
\put{$8$} at 3 -4
\put{$12$} at 4 -4
\put{$30$} at 5 -4
\put{$18$} at 6 -4
\put{$24$} at 7 -4
\put{$24$} at 8 -4
\put{$24$} at 9 -4
\put{$\cdots$} at 10 -4

\plot 0.2 0.5  0.2 -4.5 /
\plot -1.4 -.5  10.5 -.5 /
\setdots <.5mm>
\plot 5.5 0.5  5.5 -4.5 /
\plot 6.5 0.5  6.5 -4.5 /
\setshadegrid span <.6mm>
\vshade 5.5 -4.5 0.5  <z,z,,> 6.5 -4.5 0.5 /

\endpicture}
$$
We have separated the representation types: finite, tame, wild.
For $n\le 5$,
we deal with finite type, the dotted case $n = 6$ is tame, the cases $n\ge 7$
are wild. Note that for the wild cases, the numbers do not change. It comes as a 
surprise that all the numbers $g,\ g+h$ and $2g+h$ in the tame case are much smaller 
than the corresponding ones in the finite case $n = 5.$ 

\medskip
These numbers concern {\bf central} half-lines and lines. If we allow 
also non-central lines, we get different numbers, at least in the finite type cases.
There is the following interesting observation: {\it If $X$ is an indecomposable object
in $\Cal S(5)$, then at least one of the numbers $pX,\ qX,\ rX$ is an integer.}
This may be reformulated as follows: If $X$ is an indecomposable object in $\Cal S(n)$, then
$\pr X$ lies on one of the 15 lines $p = 0,\ 1,\ 2,\ 3,\ 4,$ or
$q = 1,\ 2,\ 3,\ 4,\ 5,$ or $r = 0,\ 1,\ 2,\ 3,\ 4,$ But actually, one checks easily
that $X$ lies on one of the {\bf seven} lines $p = 1,\ 2,$ or $q = 3,\ 4,\ 5$, or $r = 1,\ 2.$

\subsection{The triangle support of $\Cal P(n)$, for $1\le n \le 5.$}
\label{sec-ten-ten}

We also look at the triangle support of $\Cal S(n)$, where $1\le n \le 5.$ 
We denote by $\Psi(\Cal S(n))$ the set of numbers $d$ of the form $d = u/b$, where 
$u = uX,$ and $b = bX$ for some indecomposable object $X$ in $\Cal S(n)$ which is
$u$-minimal and not central. 

It is easy to see that
$$
 \Psi(\Cal S(1)) =  \Psi(\Cal S(2)) =  \Psi(\Cal S(3)) = \{0\},
$$
and that 
$$
 \Psi(\Cal S(4)) =  \{0,1\}.
$$
For $n = 5,$ we get (see for example Section~\ref{sec-fifteen-one}) that
$$
 \Psi(\Cal S(5)) =  \{0,1,4/3,3/2\}.
$$

Actually, looking at the case $n = 5$, we see that it is very convenient to work 
with costandard triangles, since it turns out that the indecomposable objects in $\Cal S(5)$
which live on $\Delta_{4/3}$ and $\Delta_{3/2}$ all have support in $\nabla_2.$

$$
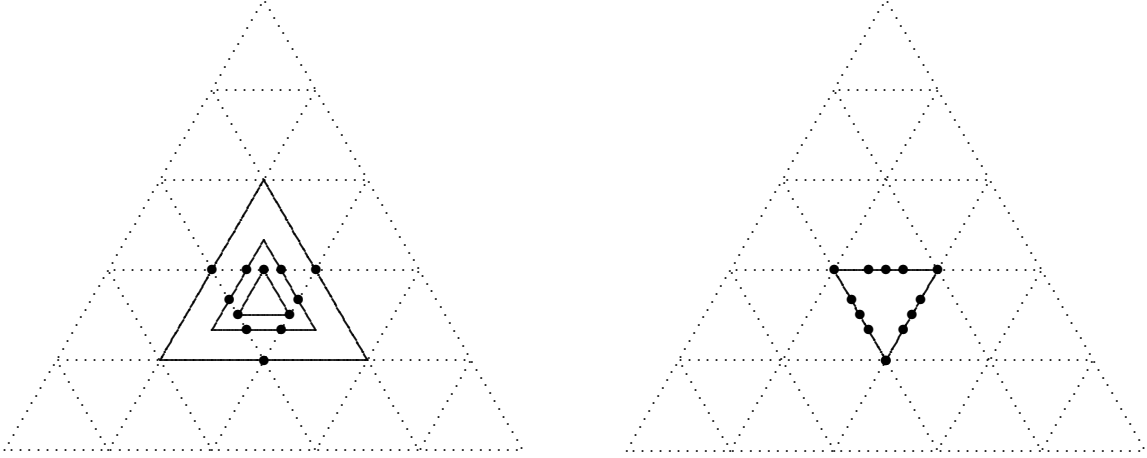
{\beginpicture
   \setcoordinatesystem units <.5169cm,.9cm>
\put{\beginpicture
\multiput{} at -5 0  5 5 /
\setdots <1mm>
\plot -5 0  0 5  5 0  -5 0 /
\plot -3 0  1 4  -1 4  3 0  4 1  -4 1  -3 0 /
\plot -1 0  2 3  -2 3  1 0  3 2  -3 2  -1 0 /
\multiput{$\ss\bullet$} at 0 1  -1 2  1 2 
    -.333 1.333  .333 1.333  -.6667 1.6667  .6667 1.6667 
   -0.333 2  .333 2   /
\multiput{$\ss\bullet$} at .5 1.5  -.5 1.5  0 2 / 
\setsolid
\plot -.5 1.5  .5 1.5  0 2  -.5 1.5 /
\plot -1 1.333   1 1.333  0 2.3333  -1 1.333 /
\plot -2 1  2 1  0 3  -2 1 /

\endpicture} at 0 0 

\put{\beginpicture
\multiput{} at -5 0  5 5 /
\setdots <1mm>
\plot -5 0  0 5  5 0  -5 0 /
\plot -3 0  1 4  -1 4  3 0  4 1  -4 1  -3 0 /
\plot -1 0  2 3  -2 3  1 0  3 2  -3 2  -1 0 /
\multiput{$\ss\bullet$} at 0 1  -1 2  1 2 /
\multiput{$\ss\bullet$} at 0 1  -.333 1.333  .333 1.333  -.6667 1.6667  .6667 1.6667 
   -0.333 2  .333 2   /
\multiput{$\ss\bullet$} at .5 1.5  -.5 1.5  0 2 / 
\setsolid
\plot -1 2  1 2  0 1  -1 2 /
\endpicture} at 12 0 

\endpicture}
$$
We see: 
A single nabla takes care of a lot of indecomposables which live on 
different triangles. In particular, we can describe $\Cal S(5)$ by using
(central half-lines as well as) just
two standard triangles and one costandard triangle.

\medskip
\centerline{$*\ *\ *$}

\medskip
The principal component $\Cal P(n)$ is closed under the duality $\D$. Components which
are closed under duality have some special properties which we are going to mention
at the end of this Section.

\subsection{Components closed under duality.}
\label{sec-ten-eleven}

As we will see in Section~\ref{sec-ten-twelve},
the $\tau$-orbits of components closed under duality have a quite
restricted behaviour. In this section, we describe basic properties of components
which are closed under duality.

\medskip
\begin{lemma}
  \label{lem-ten-eleven}
  \begin{itemize}
    \item[\bf (a)] Let $\Cal C$ be an Auslander-Reiten component of $\Cal S(n)$, where 
      $n \ge 2.$ The following conditions are equivalent:
      \begin{itemize}[leftmargin=3em]
      \item[\rm(i)] $\Cal C$ is closed under duality.
      \item[\rm(ii)] All $\tau$-orbits in $\Cal C$ are closed under duality.
      \item[\rm(iii)] There is a $\tau$-orbit in $\Cal C$ which is closed under duality.
      \item[\rm(iv)] There is an object $X$ in $\Cal C$ such that $\D X$ belongs to $\Cal C$.
      \item[\rm(v)] There exists an object $X$ in $\Cal C$ which is self-dual.
      \end{itemize}
    \item[\bf (b)] If $n\ge 6$ there is the following additional equivalent condition.
      \begin{itemize}[leftmargin=3em]
      \item[\rm(vi)] There are infinitely many objects $X$ in $\Cal C$ which are self-dual.
      \end{itemize}
  \end{itemize}
\end{lemma}

\begin{proof}
  The duality $\D$ is a contravariant equivalence, thus it sends an irreducible map
  $X \to Y$ to an irreducible map $\D Y \to \D X$. It follows that it sends Auslander-Reiten
  components to Auslander-Reiten components, thus (iv) implies (i). For $n\ge 6$, the duality 
  preserves the quasi-length of an
  indecomposable object $X$. This shows that (i) implies (ii) in case $n \ge 6.$ 
  For the cases $n \le 5$, the implication (i) implies (ii) can be checked easily.
  
  Of course, (ii) implies (iii), and (iii) implies (iv). Also, (v) implies (iv). 
  
  Finally, let us show that (ii) implies (v). We assume (ii). 
  
  First, let us show that there is at least one object in $\Cal C$ which is self-dual.
  The component $\Cal C = \Cal P(n)$ for $n \ge 2$ satisfies (v): Namely, there is 
  $X = ([1],[2])$ in $\Cal P(n)$.  
  Thus, we can
  assume that $n \ge 6$ and that 
  $\Cal C$ is a stable tube. Let $Z$ belong to the boundary of $\Cal C$, 
  thus, according to (ii), $\D Z = \tau^iZ$ for some $i \ge 0$. Let $X = (\tau^i Z)[i+1].$ 
  We claim that $\D X$ is isomorphic to $X$. Namely, the chain of irreducible epimorphism
  $\tau^i(Z[i+1]) \to \cdots \to \tau(Z[2]) \to Z$ is sent under $\D$ 
  to the chain of irreducible monomorphism 
  $\D\tau^i(Z[i+1]) \leftarrow \cdots \leftarrow \D\tau(Z[2]) \leftarrow \D Z$. Since 
  $\D Z = \tau^i Z$, we see that $\D\tau^i(Z[i+1]) \simeq (\tau^i Z)[i+1],$ thus $X$
  is self-dual. 
  
  For part (b), it remains to observe: If $X$ is self-dual, then also $X^+$ is self-dual. 
\end{proof}

\medskip
\begin{remark}
  If $X$ belongs to a component which is closed under duality, then (ii)
  asserts that $\D X = \tau^iX$ for some $0 \le i \le 5.$ Note that all cases do occur:
  The $\tau$-orbit of $S$ yields (for $n\ge 3$) objects $X$ with $i = 1,3,5$
  (and the minimal $\tau$-period of $S$ is $6$).
  The $\tau$-orbit of $([1],[2])$ yields (for $n\ge 4$) objects $X$ with 
  $i = 0,2,4$ (and the minimal $\tau$-period of $([1],[2])$ is $6$).
  Of course, the case $i = 0$ has already been addressed in (v). 
\end{remark}

\subsection{The $\tau$-orbits closed under duality.}
\label{sec-ten-twelve}

\begin{proposition}
  Let $\Cal O$ be a $\tau$-orbit in $\Cal S(n)$ which is closed
  under duality and which contains an object which does not lie on a diagonal line. Then 
  $\Cal O$ has cardinality $6$, lies on a central circle, and 
  no object in $\Cal O$ belongs to a diagonal line. There are the following two
  possibilies:
  \begin{itemize}[leftmargin=3em]
  \item[\rm (a)] The half-line support of $\Cal O$ is $\mathbb P(n)$.
  \item[\rm (b)] The half-line support of $\Cal O$ consists of six central half-lines 
    outside of $\mathbb P(n)$ and does not contain half-lines which are complementary.
  \end{itemize}
\end{proposition}

\begin{proof}
  We assume that $X \in \Cal O$ does not lie on a diagonal line. In particular,
  $X$ is not central, thus $X, \tau^2X,$ and $\tau^4X$ have pairwise different pr-vectors.
  Also, by assumption, $\D$ maps $\Cal O$ into itself. Since $X$ does not lie on a diagonal
  line, the pr-vectors of the objects in $\Cal O$ are six points on a central circle --- this
  set is the orbit of $\pr X$ under the symmetry group $\Sigma_3$ of $\mathbb T(n)$.
  
  If $\Cal O$ lies inside $\mathbb P(n)$, all objects in $\Cal O$
  lie on the central lines parallel to the boundary, and 
  the half-line support of $\Cal O$ are the six central half-lines parallel to the boundary.
  This is case (a).
  
  If $\Cal O$ does not lie inside $\mathbb P(n)$, then we consider as in Section~\ref{sec-ten-nine}
  the triangles $\mathbb F_1$ and $\mathbb F_2$
  of all vectors $(p,r)$ with $p\leq r\leq \frac13n$ and
  $\frac13n\leq r\leq\frac12(n-p)$, respectively. 
  Since one point in $\Cal O$ is contained in the union of the interiors
  $\mathbb F_1^\circ\cup \mathbb F_2^\circ$, the remaining ones are obtained by applying $\Sigma_3$. 
  Note that for each $i$, the union 
  $\cup_{\sigma\in\Sigma_3}\; \sigma\cdot \mathbb F_i^\circ$ does not contain points on
  complementary half-lines. This yields case (b).
\end{proof}

\medskip
\begin{remark}
  A typical example of case (a) are the objects of
  quasi-length 3 in $\Cal P(6)$.
  A typical example of case (b) are the objects of quasi-length 1 in $\Cal P(6)$.
\end{remark}

\medskip
\begin{remark}
  If the $\tau$-orbit $\Cal O$ is a closed under duality, and contains
  an object which lies on a diagonal line, there are many possibilities:
  \begin{itemize}[leftmargin=3em]
  \item[$\bullet$] $|\Cal O| \le 2,$ and then all objects in $\Cal O$ are
    central; see for $n = 6$ the objects in any stable tube of rank 1 or 2. Thus, the
    half-line support of $\Cal O$ is empty.
  \item[$\bullet$] $|\Cal O| = 3,$ and all objects in $\Cal O$ are
    central: see for $n = 6$ the $\tau$-orbit of $([4,2],[6,3,3]).$ Again, the
    half-line support of $\Cal O$ is empty.  The objects in this and in the following two
    examples occur in the tube of rank 3 with rationality index 0.
  \item[$\bullet$] $|\Cal O| = 3,$ no object in $\Cal O$ is
    central, and the half-line support of $\Cal O$ consists of the three diagonal half-lines
    which contain the corners of $\mathbb T(6)$:
    see for $n = 6$ the $\tau$-orbit of $([3],[4,2])$.
  \item[$\bullet$] $|\Cal O| = 3,$ no object in $\Cal O$ is
    central, and the half-line support of $\Cal O$ consists of the three diagonal half-lines
    which do not contain the corners of $\mathbb T(6)$:
    see for $n = 6$ the $\tau$-orbit of $([3],[6])$.
  \item[$\bullet$] $|\Cal O| = 6,$ and all objects in $\Cal O$ are
    central: see for $n = 6$ the objects of quasi-length 6 in any stable tube. 
    Again, the
    half-line support of $\Cal O$ is empty.
  \item[$\bullet$] $|\Cal O| = 6,$ and  precisely 3 objects in $\Cal O$ are
    central: for $n = 6$, see  the objects of quasi-length 6 in $\Cal P(6)$,
    or see the objects of quasi-length 2 or 4 in the tube of rank 6 with 
    rationality index 1
    (for quasi-length 2, the half-lines contain
    the corners of $\mathbb T(6)$; for quasi-length 4, we obtain the complementary
    half-lines).
  \item[$\bullet$] $|\Cal O| = 6,$ no object in $\Cal O$ is central, and
    the half-line support of $\Cal O$ are three central half-lines: see the objects of
    quasi-length 4 in $\Cal P(6)$ (here, the half-lines do not contain
    the corners of $\mathbb T(6)$; is there a similar example where the half-lines 
    contain the corners? There is no such example for $n=6$; after all $\Cal S(6)$
    has only two components of rank 6 which are self-dual.)
  \item[$\bullet$] $|\Cal O| = 6,$ no object in $\Cal O$ is central, and
    the half-line support of $\Cal O$ consists of six half-lines: see the objects of
    quasi-length 2 in $\Cal P(6).$ 
  \end{itemize}
\end{remark}

\medskip
\begin{remark}
  If $X$ belongs to a component which is {\bf not} closed under duality, 
  $X$ may lie on a diagonal line, whereas $\tau X$ does not lie on a diagonal line. See
  the example in Section~\ref{sec-ten-two}.
\end{remark}

\medskip
\begin{remark}
  For $n = 6,$ there are just $6$ components which are closed
  under duality and which are not homogeneous tubes.
\end{remark}

\medskip
Namely the non-homogeneous components
with rationality index $0$ and $1$.  
(In Section~\ref{sec-twelve} we will consider the case $n = 6$ in detail, based on
our previous investigations in \cite{RS1}. In particular, we will recall the rationality
index of any component in $\Cal S(6)$, it is a non-negative rational number.)

It follows from Lemma~\ref{lem-ten-eleven} that 
there are infinitely many indecomposable objects $X$ in $\Cal S(6)$ which are
self-dual. If $n < n'$, any self-dual indecomposable object in $\Cal S(n)$ yields a 
component in $\Cal S(n')$ which is closed under duality,  thus, using again Lemma 10.11,
we obtain many new indecomposable objects which are self-dual, 
and, in this way, components which are closed under duality.
It follows: {\it For $n\ge 7,$ there are infinitely many components which are closed under
  duality and which are not homogeneous.} 

\vfill\eject

\centerline{\Gross Fourth part: The case {\Grossit n}\, =\, 6.}
\addcontentsline{toc}{part}{Fourth part: The case $n=6$.}

\section{The case $n = 6.$}
\label{sec-eleven}

Here we consider the case $n = 6.$ We are going to determine 
the half-line support as well
as the triangle support of the whole category $\Cal S(6)$. The main results to be established
are Theorems \ref{theoremseven} and \ref{theoremeight}. For the proofs, we need 
details about the category $\Cal S(6)$ which have been obtained in our previous
paper \cite{RS1}, a report will be given in Section~\ref{sec-eleven-one}.
In particular, the root system for the 
covering $\widetilde{\Cal S}(6)$ will play an essential role. 
Some of the considerations in this section will be parallel
to arguments presented in \cite{RS1}
in order to compare the invariants $u$ and $v$ (or $u$ and $w$). As we mentioned already, in
the new paper 
we focus the attention to $b$ as an additional relevant invariant, and the formulae in 
Theorem~\ref{theoremseven} provide a corresponding comparison of $b$ with $u,v,w$.

\subsection{The structure of the category $\Cal S(6).$ {\rm \normalsize (A report)}.}
\label{sec-eleven-one}

In order to deal with $\Cal S(6)$, we follow the previous paper \cite{RS1}. In particular,
we use again covering theory as explained in Section~\ref{sec-two-three}.
In this way, the classification
problem for $\Cal S(6)$ was reduced in \cite{RS1} to deal with tubular algebras.

Since the principal component $\Cal P(6)$ is known, 
we may concentrate on the objects of $\Cal S(6)$ which lie in stable tubes.
In order to deal with the objects which lie in stable tubes, we use as in \cite{RS1} 
the following algebra $\Theta$, with vertices $0,1,\dots,6$ and $2',3',4'$,
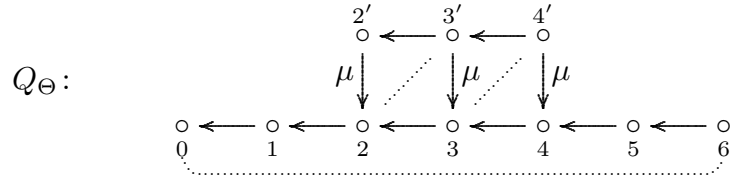
$$
\hbox{\beginpicture
\setcoordinatesystem units <0.5cm,0.5cm>
\put{} at -3 0
\put{} at 5 2
\put{$Q_\Theta\:$} at -3 1
\put{$\circ$} at 0 0
\put{$\circ$} at 2 0
\put{$\circ$} at 4 0
\put{$\circ$} at 6 0
\put{$\circ$} at 8 0
\put{$\circ$} at 10 0
\put{$\circ$} at 12 0
\put{$\circ$} at 4 2
\put{$\circ$} at 6 2
\put{$\circ$} at 8 2
\arr{1.6 0}{0.4 0}
\arr{3.6 0}{2.4 0}
\arr{5.6 0}{4.4 0}
\arr{7.6 0}{6.4 0}
\arr{9.6 0}{8.4 0}
\arr{11.6 0}{10.4 0}

\arr{5.6 2}{4.4 2}
\arr{7.6 2}{6.4 2}

\arr{4 1.6}{4 0.4}
\arr{6 1.6}{6 0.4}
\arr{8 1.6}{8 0.4}


\put{$\ssize 2'$} at  4 2.5
\put{$\ssize 3'$} at  6 2.5
\put{$\ssize 4'$} at  8 2.5

\put{$\ssize 0$} at 0 -.5
\put{$\ssize 1$} at 2 -.5
\put{$\ssize 2$} at 4 -.5
\put{$\ssize 3$} at 6 -.5
\put{$\ssize 4$} at 8 -.5
\put{$\ssize 5$} at 10 -.5
\put{$\ssize 6$} at 12 -.5

\setdots<2pt>
\plot 5.5 1.5  4.5 0.5 /
\plot 7.5 1.5  6.5 0.5 /
\plot 0 -0.7  0.1 -0.9  0.2 -1  0.3 -1.05  11.7 -1.05  11.8 -1  11.9 -0.9  12 -0.7 /

\multiput{$\mu$} at 3.6 1.1  6.4 1.1  8.4 1.1 /
\endpicture}
$$
with two commutativity relations and one zero relation. Given a representation $M$
of $\Theta$, let $VM$ be the restriction of $M$ to the full subquiver with vertices
$0,1,\dots,6$ (this is a $\Theta$-submodule of $M$), and let $UM$ be the submodule of $VM$
generated by $M_0, M_1$ and the images of the three maps $\mu$. 

Note that $VM$ is a representation of the quiver 
$\mathbb Z$ (with vertex set $\mathbb Z$ and with arrows $(i\!-\!1) \leftarrow i$ for all $i\in \mathbb Z$).
Starting with a finite-dimensional 
representation $V$ of the quiver $\mathbb Z$, let $\pi V$ be the $k[T]$-module with
underlying vector space $\bigoplus_{i\in \mathbb Z} V_i$, and with $T$ operating on this
vector space via the given maps $V_{i-1} \leftarrow V_i.$ 
Starting with a representation $M$ of $\Theta$, the zero relation of $\Theta$ shows that
the $k[T]$-module $\pi(VM)$ belongs to $\Cal N(6).$ Due to the commutativity relations, 
$\pi(UM)$ is an invariant subspace
of $\pi(VM)$. Altogether, we obtain a functor $\pi\:\mod \Theta \to \Cal S(6)$,
namely $\pi M = (\pi(UM),\pi(VM))$.

\medskip
The algebra $\Theta$ is a tubular algebra, see \cite{RS1}.
The shape of the category of $\Theta$-modules is as follows:
$$
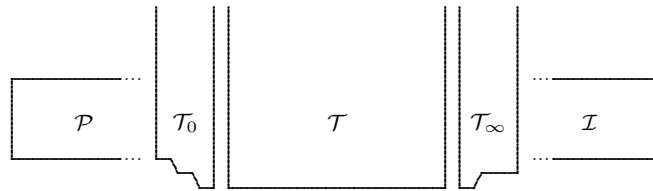
\hbox{\beginpicture
\setcoordinatesystem units <0.8cm,0.8cm>
\put{} at 0 0
\put{} at 9 2.5 
\plot 1.5 0.4  0 0.4  0 1.5  1.5 1.5 /
\setdots<2pt>
\plot 1.5 0.4  1.8 0.4 /
\plot 1.5 1.5  1.8 1.5 /
\setsolid
\plot 2 2.5  2 0.4  2.2 0.4  2.3 .2  2.5 .2  2.6 0  2.8 0  2.8 2.5  /
\plot 3 2.5  3 0  6 0  6 2.5 /
\plot 6.2 2.5  6.2 0  6.4 0  6.5  0.2  7 .2  7 2.5 /
\plot 7.5 0.4  9 0.4  9 1.5  7.5 1.5 /
\setdots<2pt>
\plot 7.5 0.4  7.2 0.4 /
\plot 7.5 1.5  7.2 1.5 /

\setsolid

\put{$\ssize \Cal P$} at 1 .9
\put{$\ssize \Cal T_0$} at 2.4 .9   
\put{$\ssize \Cal T$} at 4.5 .9
\put{$\ssize \Cal T_\infty$} at 6.6 .9
\put{$\ssize \Cal I$} at 8 .9
\endpicture}
$$
There are a preprojective component $\Cal P$ of type $\widetilde {\mathbf E}_7$,
a preinjective component $\Cal I$  of type $\widetilde {\mathbf E}_8$,
as well as tubular families $\Cal T_\gamma$
with $\gamma $ a non-negative rational number or the symbol $\infty$.
In \cite{RS1}, the number $\gamma$ was called the index of the modules in $\Cal T_\gamma$;
in order to be more specific, we now prefer to call it the {\it rationality index}.
Let $\Cal T = \bigcup_{0 < \gamma <\infty} \Cal T_\gamma$. Note that 
the indecomposable projective $\Theta$-modules corresponding to the vertices different from
$4'$ and $6$ belong to $\Cal P$, the remaining two belong to a tube of rank 6 in $\Cal T_0.$
The indecomposable injective $\Theta$-modules corresponding to the vertices different from
$0$ belong to $\Cal I$, the remaining one belongs to a tube of rank 6 in $\Cal T_\infty.$

We denote by $'\Cal D$ the union of $\Cal T$ and the stable tubes in $\Cal T_0$.
Note that $'\Cal D$ can be considered as part of the fundamental region 
of the universal covering  $\widetilde{\Cal S}(6)$ of $\Cal S(6)$
under the shift operation (the fundamental region was labelled $\Cal D$ in \cite{RS1}),
and $\pi$ can be considered as  part of the pushdown functor from $\widetilde{\Cal S}(6)$
to $\Cal S(6)$. 
It has been shown in \cite{RS1} that {\it $\pi$ provides a bijection between $'\Cal D$ and
the class of those objects in $\Cal S(6)$ which lie in stable tubes.} 
For $X$ in a stable component, we define $\widetilde X$ as the $\Theta$-module $\widetilde X = M$
in $'\Cal D$ with $\pi M = X$ (according to \cite{RS1}, $M$ exists and is uniquely determined).
In this way, we obtain a bijection $\pi$ between the class of $\Theta$-modules in $'\Cal D$ and the
objects in $\Cal S(n)$ which belong to stable tubes. 

Using this bijection, we attach to any indecomposable object in $\Cal S(6)$
its {\it rationality index}; it is a non-negative rational number.
If $X$ in $\Cal S(6)$ belongs to a stable component, the
{\it rationality index} of $X$ is just the rationality index of $\widetilde X$.
If $X$ belongs to $\Cal P(6)$, then, by definition, its rationality index is $0$
(of course, if $X$ belongs to $\pi (\Cal T_0)$, its rationality index is also $0$,  whereas
if $X$ is in $\pi \Cal T$, then its rationality index is positive). 

\medskip
We also need the algebra 
$\Xi$, obtained from $\Theta$ by deleting the vertices $0$ and $4'$.
The algebra $\Xi$ is a tilted algebra of type $\mathbf E_8.$ 
If $M$ is an indecomposable $\Theta$-module, then $\bdim M$
is a radical vector or a 
positive root of $\chi_\Theta$, and we can write
$$
  \bdim M =  \mathbf r + a_0\mathbf h_0 + a_\infty\mathbf h_\infty,
$$
where $\mathbf r = \mathbf r(M)$ 
is either zero or a (positive or negative) root of $\chi_\Xi$, where 
$$
 \mathbf h_0 = \smallmatrix  &  & 2 & 1 & 0 \cr
             1 & 2 & 3 & 3 & 2 & 1& 0 \endsmallmatrix, \quad 
 \mathbf h_\infty = \smallmatrix  &  & 2 & 2 & 1  \cr
             0 & 1 & 2 & 3 & 3 & 2 & 1 \endsmallmatrix,
$$
and where $a_0 = \dim M_0$ and $a_\infty = \dim M_{4'}$. In case 
$\mathbf r = 0,$
{\it the fraction $a_\infty/a_0$ is just the rationality index of $M$.}
(For further information about the root system of $\Xi$, we refer to Appendix~\ref{app-A}.)
If $X$ belongs to a stable tube of $\Cal S(6)$, let $\mathbf r(X) = \mathbf r(\widetilde X).$
In this way, we attach to any indecomposable object of $\Cal S(6)$ which belongs to
a stable tube an element of $\mathbb R^8$ which is a root of $\chi_\Xi$ or the zero vector.

\medskip
Starting with an indecomposable object $X$ in $\Cal S(6)$, it often is of interest to
know its rationality index. The following Lemma shows that for some objects $X$ in $\Cal S(6)$,
the width of $X$ provides at least some partial information:
It shows that for all objects $X$ in stable tubes
with $\mathbf r(X) = 0,$ there are coprime positive integers $i,j$ with
$bX = 3(i+j)$, such that the rationality index of $X$ is $i/j.$ 

\medskip
\begin{lemma}
  Let $M$ be an indecomposable $\Theta$-module which belongs to 
  $\Cal T_\gamma$. where $\gamma = i/j$ with coprime positive integers $i,j$. Assume that
  the endomorphism ring of $M$ is $k$ and that $M$ has non-trivial self-extensions. 
  Then
  $$ 
  bM = 3(i+j).
  $$
\end{lemma}

\begin{proof}
  By assumption, $M$ belongs to a stable tube, say of rank $r$.
  Since $\End(M) = k,$
  the quasi-length $\ell$ of $M$ is at most $r$, since $\Ext^1(M,M) \neq 0,$ we have $\ell \ge r.$
  Thus, $\ell = r$. Using again that $\End(M) = k,$ we see that $\bdim M = j \mathbf h_0+i\mathbf h_\infty.$
  Since $b\mathbf h_0 = b\mathbf h_\infty = 3,$ we have $b M = 3(i+j).$
\end{proof}
  
\medskip
Let us also mention that the rationality index of any 
indecomposable non-projective object $X\in\Cal S(6)$ determines 
the indices of the objects in its $\Sigma_3$-orbit:
The objects $X$, $\tau^2X$ and $\tau^4X$  have
the same rationality index, since they belong to the same tube as $X$.
If $X$ has rationality index $\gamma =0$ or $1$, the object $DX$ belongs to the same
tube as $X$, thus $DX$ has the same rationality index as $X$. 
Finally, if $X$ has rationality index $\gamma \neq0, 1$, the object
$DX$ has rationality index $\gamma^{-1}$. (It is enough to verify this for the objects $X$
of quasi-length $6$. In this case, we have $\mathbf r(X) = 0,$ thus 
$\bdim \widetilde X = a_0\mathbf h_0+a_\infty\mathbf h_\infty$, and the rationality index
is $a_\infty/a_0.$ For $Y = DX$, we have 
$\bdim \widetilde Y = a_\infty\mathbf h_0+a_0\mathbf h_\infty$, thus the rationality index
of $Y$ is $a_0/a_\infty.$)

\subsection{The width.}
\label{sec-eleven-two}

\begin{lemma}
  \label{lem-eleven-two}
  If $M$ is an indecomposable $\Theta$-module which belongs to a stable tube, then
  $$
  b(\pi M) = \dim M_3.
  $$
\end{lemma}

\begin{proof}
  Let $\Cal R$ be the class of indecomposable $\Theta$-modules which belong
  to stable tubes. We denote by $\overline
  {S(2)}$ the 2-dimensional indecomposable $\Theta$-module with 
  $\overline{S(2)}_2 = k = \overline{S(2)}_{2'}$. The essential observation for the proof of
  Lemma is the fact that $\overline{S(2)}$ belongs to the non-stable tube
  in $\Cal T_0$ and that $S(4)$ belongs to the non-stable tube in $\Cal T_\infty$.

  \medskip
  Given any algebra $A$, and indecomposable $A$-modules $X,Y$, we say that $X$ is a {\it 
    predecessor} of $Y$ provided there is a sequence $X = X_0, X_1,\dots, X_t = Y$
  of $A$-modules such that $\Hom_A(X_{i-1},X_i) \neq 0.$ 
  
  Note that $\Ext^1(Y,X) \neq 0$
  implies that $X$ is a predecessor of $Y$. 
  In the case of $A = \Theta$, we see in this way 
  that $S(0)$ is a predecessor of $S(1),$ that $S(1)$ is a
  predecessor of $\overline {S(2)}$, that $S(4)$ is a predecessor of $S(5)$, and that$
  S(5)$ is a predecessor of $S(6).$
  
  On the other hand, since $\Theta$ is a tubular algebra, we know that no module in $\Cal R$
  is a predecessor of $\overline{S(2)}$ (since $\overline{S(2)}$ belongs to the non-stable tube
  in $\Cal T_0$),
  and $S(4)$ is not a predecessor of any module in $\Cal R$ (since $S(4)$ belongs to
  the non-stable tube in $\Cal T_\infty$). 
  
  \medskip 
  Let $M$ be a $\Theta$-module in $\Cal R$.
  We claim that the maps $M_{i-1} \leftarrow M_i$ are injective, for $4\le i \le 6$
  (thus, equivalently, that $S(i)$ with $i=4,5,6$ is not a submodule of $VM$) and
  that the maps $M_{i-1} \leftarrow M_i$ are surjective, for $1\le i \le 3$
  (thus, equivalently, that $S(i)$ with $i=0,1,2$ is not a factor module of $VM$).
  
  Namely, if $S(i)$ is a submodule of $VM$, then $0 \neq \Hom(S(i),VM) \subseteq \Hom(S(i),M),$
  therefore $S(i)$ is a predecessor of $M.$ As we have mentioned, this is not the case
  for $M\in \Cal R$ and $i= 4,5,6.$ 
  Similarly, if $S(0)$ or $S(1)$ is a factor module of $VM$, then $S(0)$ or $S(1)$ is
  a factor module of $M$ itself, thus $M$ is a predecessor of $S(0)$ or $S(1),$ respectively.
  Again, for $M\in \Cal R$, this is not the case. The situation is slightly more complicated
  for $i = 2$: If $S(2)$ is a factor module of $VM$, then we may have 
  $\Hom(M,{S(2)}) = 0,$ but at least we have $\Hom(M,\overline{S(2)}) \neq 0,$
  thus $M$ is a predecessor of $\overline{S(2)}.$ Again, for $M\in \Cal R$ this is not the
  case.
  
  This shows that for $M\in \Cal R,$ the maps $M_{i-1} \leftarrow M_i$ are injective, 
  for $4\le i \le 6,$ and the maps $M_{i-1} \leftarrow M_i$ are surjective,
  for $1 \le i \le 3.$ But this implies that $b\pi(VM) = \dim M_3$ (and, of course,
  $b\pi(VM) = b(\pi M)$). 
\end{proof}

\medskip
\begin{remark}
  The formula given by the lemma is very important for dealing with
  those indecomposable objects in $\Cal S(6)$ which belong to stable tubes. 
  The assumption that we deal with stable tubes is essential.
  We stress
  that there is no corresponding formula for the remaining indecomposable
  objects of $\Cal S(6),$ those belonging to the principal component
  $\Cal P(6)$. 

  There are two facts which one has to be aware of. First of
  all, there are many objects in $\Cal P(6)$ which are not of the form 
  $\pi M$ where $M$ is a $\Theta$-module. This concerns not only the two projective
  objects but also all the objects of
  $\Cal P(6)$ which occur in the intersection of the ray starting in the projective object
  $(0,[6])$ with the coray ending in the projective object $([6],[6])$. 
  Second, even if $X$ in $\Cal P(6)$ is of the form $X = \pi M$ for some
  $\Theta$-module $M$ in $\Cal T_0$, we may have $bX \neq \dim M_3$. A typical example
  is the object $X = S[6]$ with $bX = 4$, but $\dim M_i \le 3$ for all vertices $i$. 

  Of course, this really does not matter, since 
  all the objects of $\Cal P(6)$, in particular also their width, are known,
  see Sections~\ref{sec-nine-five}, \ref{sec-nine-ten}, as well as \ref{sec-ten-six},
  \ref{sec-ten-seven}.
\end{remark}

\subsection{Central objects.}
\label{sec-eleven-three}

\begin{lemma}
  If $M$ is an indecomposable $\Theta$-module which belongs to 
  a stable tube, 
  and $\bdim M$ is a radical vector, then $\pi M$ is a central object.
\end{lemma}

\begin{proof}
  We have $\bdim M = a_0\mathbf h_0+ a_\infty\mathbf h_\infty$, therefore
  $v(\pi M)=12 a_0+12 a_\infty$.
  According to Lemma~\ref{lem-eleven-two},
  we have $b(\pi M)=\dim M_3=3 a_0+3 a_\infty$.
  Thus $q(\pi M)=v(\pi M)/b(\pi M)=4$.
\end{proof}

\medskip
\begin{remark}
  The assumption in the lemma that $M$ belongs to a stable tube is essential.
  Namely, the object $S[6]$ in $\Cal P(n)$ is of the form $\pi M$ with $\bdim M = \mathbf h_0$,
  thus $\bdim M$ is a radical vector, whereas $\pr S[6] = (6/4,6/4)$ is not central. 
\end{remark}

\subsection{Proof of Theorem~\ref{theoremseven}.}
\label{sec-eleven-four}

\begin{lemma}
  \label{lem-eleven-four-one}
  Any indecomposable $\Theta$-module $M$ satisfies the inequality
  $$
  |vM - 4 \dim M_3| \le 4.
  $$
\end{lemma}

\begin{proof}
  We have $\bdim M =  \mathbf r + a_0\mathbf h_0 + a_\infty\mathbf h_\infty,$
  where $\mathbf r$ is zero or a root of $\chi_\Xi$. The function 
  $v(M) - 4\dim M_3 $ is additive on the Grothendieck group and vanishes on $\mathbf h_0$ as
  well as on $\mathbf h_\infty$, thus we only have to look at its values on the roots
  of $\chi_\Xi$. As in the proof of (2.4.1) in \cite{RS1}, we use the fact that the 
  roots of $\chi_\Xi$ are explicitly known (they are related to the
  roots of $\mathbf E_8$, since $\Xi$
  is a tilted algebra of type $\mathbf E_8$, see Appendix~\ref{app-A}). 
  It is quite easy to verify the desired inequality for the roots $\mathbf r$.
\end{proof}

\medskip
\begin{lemma}
  \label{lem-eleven-four-two}
  For all indecomposable $X$ in stable tubes of 
  $\Cal S(6)$, we have $|vX - 4\cdot bX| \le 4$.
\end{lemma}

\begin{proof}
  Let $X = \pi M$, where $M$ is an indecomposable $\Theta$-module which belongs to
  a stable tube. According to Section~\ref{sec-eleven-two},
  we have $bX = b(\pi M) = \dim M_3.$ Of course, we 
  have $vX = vM.$ Thus Lemma~\ref{lem-eleven-four-one}
  shows that $|4bX-vX| = |4\dim M_3-vM| \le 4.$
\end{proof}

\medskip
\begin{proof}[Proof of Theorem~\ref{theoremseven}]
  Let us write $\eta X = vX-4\cdot bX.$ 
  According to Lemma~\ref{lem-eleven-four-two}, we have $|\eta| \le 4$ in case
  $X$ belongs to a stable tube. Thus, it remains to show the same in case $X$ belongs
  to the principal tube. 
  
  The values of $b,u,w$, thus also of $v = u+w$, on the 
  principal tube are exhibited in Section~\ref{sec-nine-five}.
  As a consequence, we have the following 
  values of $\eta$
  on the objects $Z[\ell]$ with $Z$ ray-simple and $0\le \ell \le 8$ in the principal component:

$$
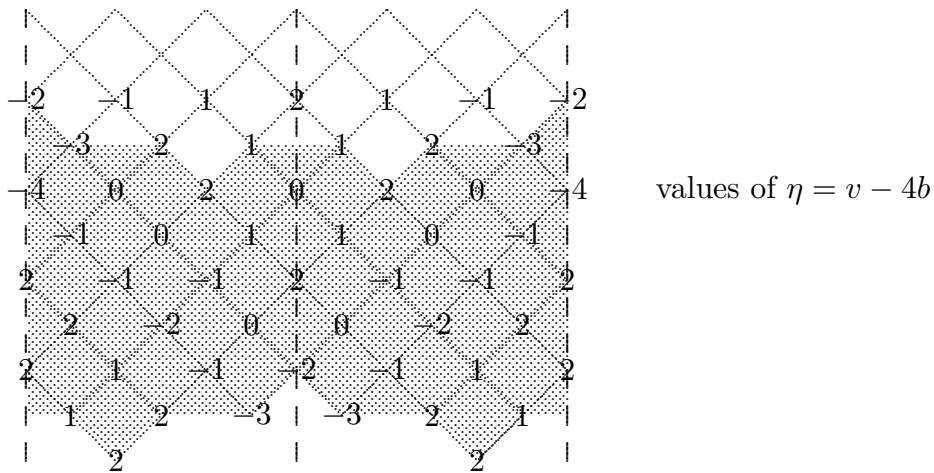
{\beginpicture
    \setcoordinatesystem units <.5cm,.5cm>
\multiput{} at -6 0  6 10 /
\setdashes <2mm>
\plot -6 0  -6 10 /
\plot  6 0   6 10 /
\plot  0 0   0 10 /

\setdots <.5mm>
\plot -6 2  -4 0  6 10  /
\plot -6 4  -3 1   /
\plot -6 6  -1 1  6 8  4 10 /
\plot -6 8  1 1   6 6   2 10 /
\plot -6 2  2 10  /
\plot -6 4  0 10  6 4  3 1 /
\plot -6 6  -2 10  6 2  4 0 /
\plot -6 8  -4 10  5 1 /
\plot -6 10  4 0  5 1 /
\plot -5 1  4 10 /
\setshadegrid span <.4mm>
\vshade -6 1 8 <,z,,> -5 1 7 <z,z,,> -4 0 7 <z,z,,> -3 1 7  <z,z,,> -2 1 6 
    <z,z,,> -1 1 7  <z,z,,> 0 2 7   <z,z,,> 1 1 7  <z,z,,>  2 1 6  
    <,z,,> 3 1 7 <z,z,,> 4 0 7 <z,z,,> 5 1 7  <z,z,,> 6 1 8 /
\multiput{$2$} at -6 2  -6  4  -5 3  -3 1  -3 7  -2 6  0 4 
                   6 2  6  4  5 3  3 1  3 7  2 6  -4 0  4 0  0 8 /
\multiput{$1$} at -5 1  -4 2 -1 7  -1 5  -2 8  5 1  4 2  1 7  1 5  2 8 /
\multiput{$0$} at -4 6  -3 5 -1 3  0 6  4 6  3 5  1 3 /
\multiput{$-1$} at  -5 5  -4 4  -2 2  -2 4  5 5  4 4  2 2  2 4 
                     -4 8  4 8 /
\multiput{$-2$} at    -3 3     3 3  0 2  -6 8  6 8 /
\multiput{$-3$} at  -5 7  -1 1  5 7 1 1 /
\multiput{$-4$} at  -6 6  6 6 /
\put{values of $\eta = v-4b$} at 11 6 
\endpicture}
$$
(Actually, we do not have to refer to Section~\ref{sec-nine-five}.
Since $\eta$ is additive on the adjusted tube, we only have to determine
the values of $\eta$ at the boundary of the adjusted tube
(the corresponding objects all are pickets), and use the
additivity of $\eta$.) 

As we see, we have $\eta Z^+ = \eta Z$ for all objects $Z$ on the boundary of the
adjusted translation quiver. 
Again, using the additivity of $\eta$, it follows that $\eta X^+ = \eta X$ for all $X$
in the principal component. We conclude that 
$-4 \le \eta X \le 2$ for all $X$ in the principal component; a fortiori, $|\eta X| \le 4$.
\end{proof}

\medskip
\begin{remark}
  As we have seen, we have the inequalites
  $-4 \le v-4b \le 2$ for all $X \in \Cal P(6)$. 
  Let us add that we also have the inequalities $-2 \le u-2b \le 4$.
  Here are the values of $u-2b$:
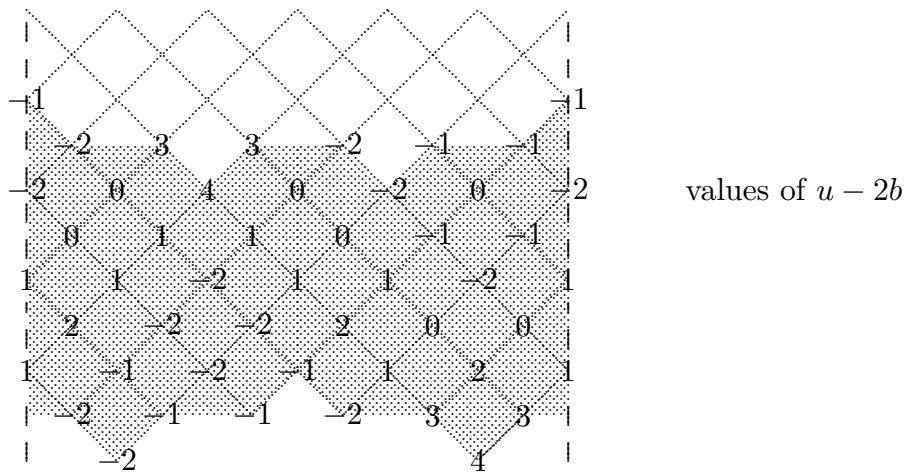
$$
{\beginpicture
    \setcoordinatesystem units <.5cm,.5cm>
\multiput{} at -6 0  14 10 /
\setdashes <2mm>
\plot -6 0  -6 1.5 /
\plot -6 2.5  -6 3.5 /
\plot -6 4.5 -6 10 /
\plot  6 0     6 1.5 /
\plot  6 2.5   6 3.5 /
\plot  6 4.5   6 10 /

\setdots <.5mm>
\plot -6 2  -4 0  6 10  /
\plot -6 4  -3 1   /
\plot -6 6  -1 1  6 8  4 10 /
\plot -6 8  1 1   6 6   2 10 /
\plot -6 2  2 10  /
\plot -6 4  0 10  6 4  3 1 /
\plot -6 6  -2 10  6 2  4 0 /
\plot -6 8  -4 10  5 1 /
\plot -6 10  4 0  5 1 /
\plot -5 1  4 10 /
\setshadegrid span <.4mm>
\vshade -6 1 8 <,z,,> -5 1 7 <z,z,,> -4 0 7 <z,z,,> -3 1 7  <z,z,,> -2 1 6
    <z,z,,> -1 1 7  <z,z,,> 0 2 7   <z,z,,> 1 1 7  <z,z,,>  2 1 6  
    <,z,,> 3 1 7 <z,z,,> 4 0 7 <z,z,,> 5 1 7  <z,z,,> 6 1 8 /
\multiput{$1$} at -6 2  -6 4  -4 4  -3 5  -1 5  0 4  2 2  2 4  6 2  6 4
                /
\multiput{$0$} at  -5 5  1 5  3 3  5 3  -4 6  0 6  4 6 /
\multiput{$-1$} at  -4 2  -3 1  0 2  3 5  5 5  -1 1  3 7  5 7  -6 8  6 8 /
\multiput{$-2$} at -5 1  -4 0  -3 3  -2 2  -2 4  -1 3  1 1  2 6  4 4  6 6  -6 6  1 7  -5 7 /
\multiput{$2$} at  4 2  -5 3  1 3 /
\multiput{$3$} at  3 1  5 1  -3 7 -1 7  /
\multiput{$4$} at 4 0  -2 6  /
\put{values of $u-2b$} at 11 6 
\endpicture}
$$
By duality, we similarly have $-2 \le w-2b \le 4$ on $\Cal P(6)$. 
\end{remark}

\subsection{A finiteness result.}
\label{sec-eleven-five}

\begin{proposition}
  \label{prop-eleven-five} 
  For any $a < 2,$ there are only finitely many 
  indecomposable objects $X\in \Cal S(6)$ 
  with $dX \le a.$
\end{proposition}

\medskip
For the proof, we need some preparation. Let 
$\eta X = vX-4\cdot bX.$ The inequality $|\eta Y| \le 4$ for $Y$ indecomposable in $\Cal S(6)$ 
(established in Section~\ref{sec-eleven-four}) has the following consequence:

\medskip
\begin{lemma}
  \label{lem-eleven-five-one}
  Let $Y$ be an indecomposable object in $\Cal S(6)$. Then
  $$
  \left |qY-4 \right | \le \frac{24}{vY}.
  $$
\end{lemma}

\begin{proof}
  We have
  $$
  |qY-4| =  \Bigl|\frac{\eta Y}{bY}\Bigr| \le
  \frac4{bY} \le \frac{24}{vY},
  $$
  where the first inequality sign is implied by $|\eta Y| \le 4$.
  The second inequality sign follows from 
  the inequality $vY \le 6\!\cdot\!bY,$ which is valid for all $Y$ in $\Cal S(6).$
\end{proof}

\medskip
\begin{lemma}
  \label{lem-eleven-five-two}
  If $|qY-4| \ge \epsilon > 0$, then
  $vY \le \frac{24}\epsilon.$
\end{lemma}

\begin{proof}
  Assume, for the contrary, that $|qY-4| \ge \epsilon$ and 
  $vY > \frac{24}\epsilon.$ Then Lemma~\ref{lem-eleven-five-one} yields 
  $|qY-4| \le \frac{24}{vY} < \epsilon,$ a contradiction.
\end{proof}

\medskip
\begin{lemma}
  \label{lem-eleven-five-three}
  For any number $c < 4$, there are only finitely many indecomposable objects $X$
  in $\Cal S(6)$ with $q(X) \le c$.
\end{lemma}

\begin{proof}
  Let $0 \le c < 4.$
  There are 
  only finitely many objects $X$ in the non-stable tube with $vX\le 24/(4-c)$;
  objects with $vX>24/(4-c)$ satisfy $|qX-4|\le 24/vX<4-c$, by
  Lemma~\ref{lem-eleven-five-one}, hence $q(X)> c.$ 

  Thus, let us look at the objects $X$ in stable tubes. According to Section~\ref{sec-eleven-one},
  we can
  choose $M \in {'\Cal D}$ with $\pi M  = X$. Since the universal covering is controlled by the
  quadratic form, $\bdim M$ is either a root or a radical vector. 
  According to Section~\ref{sec-eleven-three}, we know that
  $\bdim M$ cannot be a radical vector, thus $\bdim M$ is a root.
  There are only finitely many roots with $v(X) \le 24/(4-c)$; as above,
  the remaining objects in stable tubes satisfy $q(X)> c$.
\end{proof}

\medskip
\begin{proof}[Proof of Proposition~\ref{prop-eleven-five}]
  Using rotation, we obtain from Lemma~\ref{lem-eleven-five-three}:
  {\it For any number $c < 4,$
    there are only finitely many indecomposable objects $X$ in $\Cal S(6)$ with $pX \ge 6-c.$}
  The union of the sets of pr-vectors with $qX \le c$ and with $pX \ge 6-c$ contains all
  the pr-vectors with $r\le 2c-6.$
  
  Now assume that there is given $0 \le a < 2.$ Let $c = a/2+3.$ Then $c < 4.$ As we have seen,
  there are only finitely many indecomposables with pr-vector $(p,r)$ such that $r \le 2c-6 = a.$
  Using rotation, there are also only finitely many indecomposables with pr-vector $(p,r)$ such that
  $p \le a$ or $6-q\le a$.
  In particular, there are only finitely many indecomposable objects $X$ with $dX \le a.$
\end{proof}

\subsection{The half-line support of $\Cal S(6)$.
  {\rm\normalsize (Proof of Theorem~\ref{theoremeight}(a))}.}
\label{sec-eleven-six}

We have to show: 
{\it Let $X$ be indecomposable in $\Cal S(6)$.
Then $X$ lies on some line $L_\phi$, with $\phi\in \Phi.$}
We may assume
that $X$ belongs to a stable tube, thus there is $M$ which belongs to a stable tube of
$\mod \Theta$ with $X = \pi M$. According to Lemma~\ref{lem-eleven-two}, we have $bM = \dim M_3.$
We also can assume that $X$ is not central (since central objects lie on all
lines $L_\phi$, for example on $L_1$, since $pX = rX$).

For any element 
$$
  \mathbf z = \left(\smallmatrix & & z_{2'} & z_{3'} & z_{4'} \cr
                       z_0 & z_1 & z_2 & z_3 & z_4 & z_5 & z_6 \endsmallmatrix\right)
$$
in $K_0(\mod \Theta)$, let 
$$
   u\mathbf z = z_0+z_1+z_{2'}+z_{3'}+z_{4'}, \quad 
   v\mathbf z = \sum\nolimits_{i=0}^6 z_i, \quad w\mathbf z = v\mathbf z - u\mathbf z,\quad
   \dim_3 \mathbf z = z_3.
$$
We say that $\mathbf z$ is {\it central} provided $u\mathbf z = 2\dim_3 z = w\mathbf z$, and 
we define,
for $\mathbf z$ not central, 
$$
  \phi'\mathbf z = \frac{u\mathbf z - 2\dim_3\mathbf z}{w\mathbf z - 2\dim_3\mathbf z},
$$
this is an element of $\mathbb Q\cup\{\infty\}.$

\medskip
Any element $\mathbf z\in K_0(\mod \Theta)$ can be written as
$$
 \mathbf z = \mathbf r(\mathbf z) + z_0\mathbf h_0 + z_{4'}\mathbf h_\infty,
$$ 
thus $\mathbf r(\mathbf z)$ belongs to $K_0(\mod \Xi)$. If $M$ is an indecomposable
$\Theta$-module, then $\mathbf r(\bdim M)$ is a 
(not necessarily positive) root of the quadratic form $\chi_\Xi$. 

Since $u(\mathbf h_0) = 2\dim_3 \mathbf h_0 = w(\mathbf h_0)$, 
we see that 
\begin{align*}
   (u-2\dim_3)(\mathbf h_0) &= (w-2\dim_3)\mathbf h_0) = 0,\quad\text{and} \cr
   (u-2\dim_3)(\mathbf h_\infty) &= (w-2\dim_3)\mathbf h_\infty) = 0,
\end{align*}
thus $\mathbf z$ is central if and only if $\mathbf r(\mathbf z)$ is, and if $\mathbf z$ is not central,
then 
$$
 \phi'(\mathbf z) = \phi'(\mathbf r(\mathbf z)).
$$
 A case-by-case calculation (presented in column (5) in Appendix~\ref{app-A})
 shows that for all non-central
 roots $\mathbf r$ of $\chi_\Xi$, we have $\phi'(\mathbf r) \in \Phi.$
Actually, we only have to look at the roots of $\chi_\Xi,$ which correspond to positive roots
of $\mathbf E_8$, 
since for $\mathbf z$ not central, we have $\phi'(\mathbf z) = \phi'(-\mathbf z)$.

\medskip 
Let us summarize the considerations. Recall that we want to show that any
non-central $X$ which belongs to a stable tube 
in $\Cal S(6)$ satisfies $\phi(X) \in \Phi$ (so that $X$
lies on one of the 12 lines $L_\phi$ with $\phi\in \Phi$).
Now $X = \pi M$ for some non-central $\Theta$-module
which belongs to a stable tube in $\mod\Theta$ and we have 
$$
  \phi X = \phi M = \phi'(\bdim M) = \phi'(\mathbf r),
$$
where $\mathbf r = \mathbf r(\bdim M)$ is a non-central root of $\chi_\Xi$. Finally, we have 
$\phi'(\mathbf r) \in \Phi.$
$\s$

\subsection{More about $\mathbb L(6)$.}
\label{sec-eleven-seven}

Theorem~\ref{theoremsix} (shown in Section~\ref{sec-ten})
asserts that the principal component $\Cal P(6)$
lives on the union $\mathbb P(6)\cup\mathbb D(6)\cup\mathbb H_\ell(6)$. Recall that we denote 
by $\mathbb H_s(6)$ the 
union of the half-lines complementary to those in $\mathbb H_\ell(6)$.
There is no non-central indecomposable object $Y$ in $\Cal P(n)$ supported by $\mathbb H_s(6)$.
But there are such  objects in $\Cal S(6).$ For any
half line $H$ in $\mathbb H_s(6)$, let us exhibit a non-central indecomposable object $Y$ 
in $\Cal S(6)$ with $bY = 3$ which is supported by $H$.
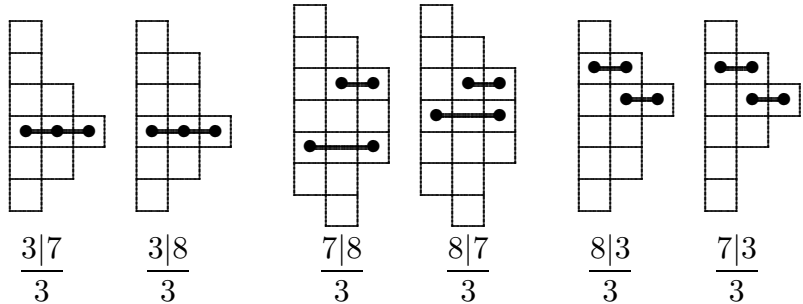
$$
{\beginpicture
    \setcoordinatesystem units <.35cm,.35cm>
\put{\beginpicture
\multiput{} at 0 0   3 6 /
\plot 0 0  1 0  1 6  0 6  0 0 /
\plot 0 1  2 1  2 4  0 4 /
\plot 0 2  3 2  3 3  0 3 /
\plot 0 5  1 5 /
\multiput{$\bullet$} at 0.5 2.5  1.5 2.5  2.5 2.5  /
\plot 0.5 2.5  2.5 2.5 /
\plot 0.5 2.55  2.5 2.55 /
\plot 0.5 2.45  2.5 2.45 /

\endpicture} at 0 0
\put{$\dfrac{3|7}3$} at -.5 -4.8
\put{\beginpicture
\multiput{} at 0 0   3 6 /
\plot 0 0  1 0  1 6  0 6  0 0 /
\plot 0 1  2 1  2 5  0 5 /
\plot 0 2  3 2  3 3  0 3 /
\plot 0 4  2 4 /
\multiput{$\bullet$} at 0.5 2.5  1.5 2.5  2.5 2.5  /
\plot 0.5 2.5  2.5 2.5 /
\plot 0.5 2.55  2.5 2.55 /
\plot 0.5 2.45  2.5 2.45 /

\endpicture} at 4 0
\put{$\dfrac{3|8}3$} at 3.5 -4.8

\put{\beginpicture
\multiput{} at 0 -1   3 6 /
\plot 0 0  1 0  1 6  0 6  0 0 /
\plot 1 0  1 -1  2 -1  2 5  0 5 /
\plot 1 0  2 0 /
\plot 0 1  3 1  3 4  0 4 /
\plot 0 2  3 2 /
\plot 0 3  3 3 /
\multiput{$\bullet$} at 0.5 1.5  2.5 1.5  1.5 3.5  2.5  3.5 /
\plot 0.5 1.5  2.5 1.5 /
\plot 0.5 1.55  2.5 1.55 /
\plot 0.5 1.45  2.5 1.45 /
\plot 1.5 3.5  2.5 3.5 /
\plot 1.5 3.55  2.5 3.55 /
\plot 1.5 3.45  2.5 3.45 /
\endpicture} at 9 0
\put{$\dfrac{7|8}3$} at 9 -4.8
\put{\beginpicture
\multiput{} at 0 -1   3 6 /
\plot 0 0  1 0  1 6  0 6  0 0 /
\plot 1 0  1 -1  2 -1  2 5  0 5 /
\plot 1 0  2 0 /
\plot 0 1  3 1  3 4  0 4 /
\plot 0 2  3 2 /
\plot 0 3  3 3 /
\multiput{$\bullet$} at 0.5 2.5  2.5 2.5  1.5 3.5  2.5  3.5 /
\plot 0.5 2.5  2.5 2.5 /
\plot 0.5 2.55  2.5 2.55 /
\plot 0.5 2.45  2.5 2.45 /
\plot 1.5 3.5  2.5 3.5 /
\plot 1.5 3.55  2.5 3.55 /
\plot 1.5 3.45  2.5 3.45 /
\endpicture} at 13 0
\put{$\dfrac{8|7}3$} at 13 -4.8

\put{\beginpicture
\multiput{} at 0 0   3 6 /
\plot 0 0  1 0  1 6  0 6  0 0 /
\plot 0 1  2 1  2 5  0 5 /
\plot 0 3  3 3  3 4  0 4 /
\plot 0 2  2 2 /
\multiput{$\bullet$} at 0.5 4.5  1.5 4.5  1.5 3.5  2.5  3.5 /
\plot 0.5 4.5  1.5 4.5 /
\plot 0.5 4.55  1.5 4.55 /
\plot 0.5 4.45  1.5 4.45 /
\plot 1.5 3.5  2.5 3.5 /
\plot 1.5 3.55  2.5 3.55 /
\plot 1.5 3.45  2.5 3.45 /
\endpicture} at 18 0
\put{$\dfrac{8|3}3$} at 17.5 -4.8

\put{$\dfrac{7|3}3$} at 21.5 -4.8
\put{\beginpicture
\multiput{} at 0 0   3 6 /
\plot 0 0  1 0  1 6  0 6  0 0 /
\plot 0 2  2 2  2 5  0 5 /
\plot 0 3  3 3  3 4  0 4 /
\plot 0 1  1 1 /
\multiput{$\bullet$} at 0.5 4.5  1.5 4.5  1.5 3.5  2.5  3.5 /
\plot 0.5 4.5  1.5 4.5 /
\plot 0.5 4.55  1.5 4.55 /
\plot 0.5 4.45  1.5 4.45 /
\plot 1.5 3.5  2.5 3.5 /
\plot 1.5 3.55  2.5 3.55 /
\plot 1.5 3.45  2.5 3.45 /
\endpicture} at 22 0

\endpicture}
$$

\bigskip
The following picture shows all the
\phantomsection{lines}
\addcontentsline{lof}{subsection}{The 12 lines $\mathbb L(6)$ and their slopes.}
$L = L_\phi$ in $\mathbb L(6) = \mathbb P(6)\cup\mathbb D(6)\cup
\mathbb H(6)$, 
with the
corresponding label $\phi$.
The lines in
$\mathbb P(6)$ and $\mathbb D(6)$ are solid, the remaining ones in $\mathbb H(6)$ are dashed. 
For every central half-line, the non-central objects with smallest possible width are indicated:
pickets by bullets $\bullet$, non-pickets by small circles $\circ$ (these are the six
objects of width 3
mentioned above).

$$
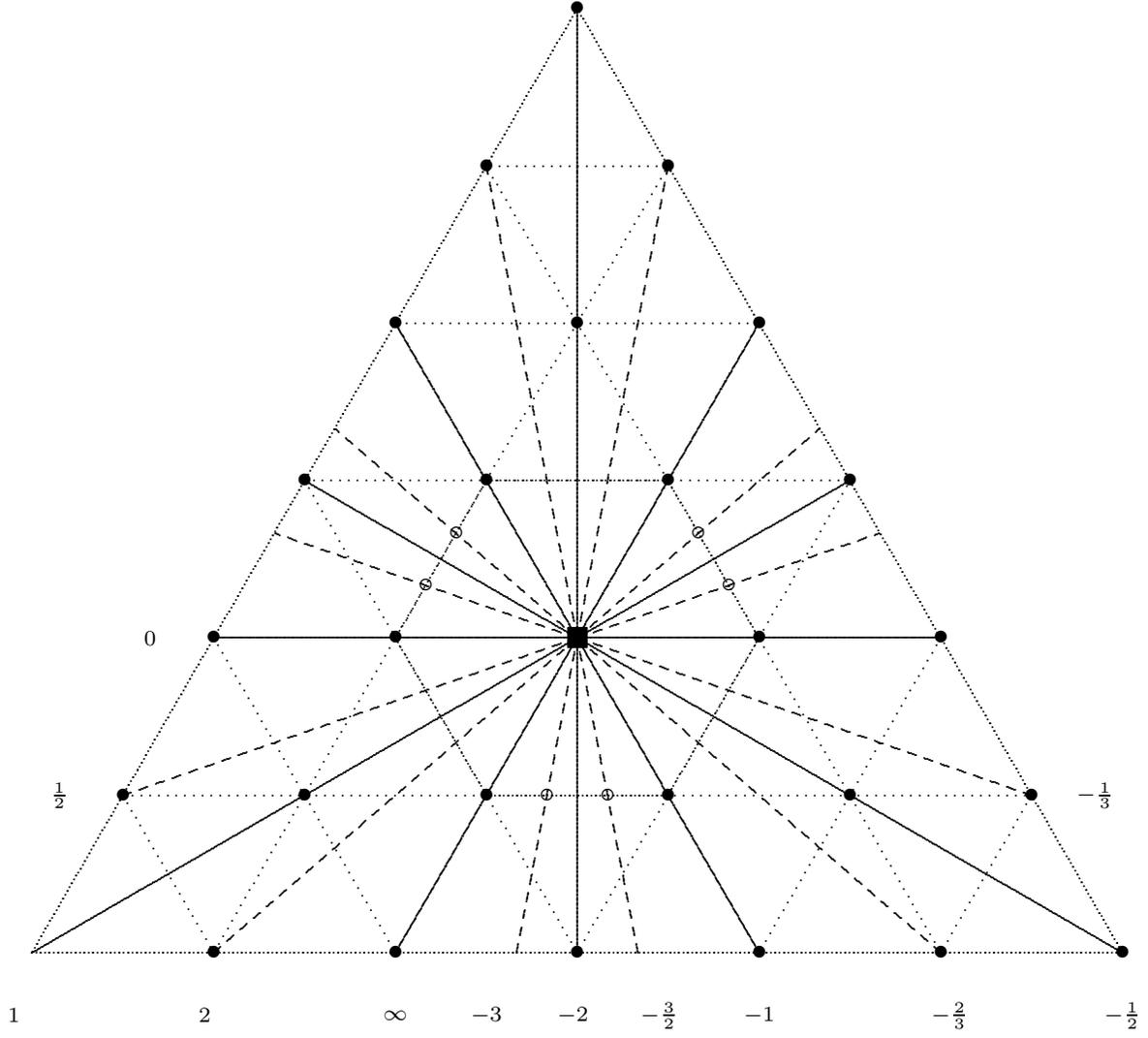
{\beginpicture
\setcoordinatesystem units <.92355cm, 1.6cm>
\multiput{} at -6 -.5  6 6 /
\setdots <1mm> 
\plot -6 0  6 0  0 6  -6 0 /
\plot -5 1  -4 0  1 5  -1 5  4 0  5 1  -5 1 /
\plot -4 2  -2 0  2 4  -2 4  2 0  4 2  -4 2 /
\plot -3 3  0 0  3 3  -3 3 /

\setdots <.5mm> 
\plot -6 0  6 0  0 6  -6 0 /
\plot -1 1  1 1  2 2  1 3  -1 3  -2 2  -1 1 /

\multiput{$\bullet$} at 
     -4 0  -2 0  0 0  2 0  4 0  
     -3 1  -1 1  1 1  3 1 
     -2 2  0 2  2 2
     -1 3  1 3
      0 4
       -5 1  -4 2  -3 3  -2 4  -1 5  1 5  2 4  3 3  4 2  5 1  /
\multiput{$\circ$} at -.333 1  .333 1  1.333 2.667  1.667 2.333  -1.333 2.667  -1.667 2.333 /
\put{$\blacksquare$} at 0 2
\setsolid
\plot 0 0  0 6 /
\plot -6 0  3 3 /
\plot  6 0  -3 3 /

\plot -2 0  2 4 /
\plot 2 0  -2 4 /
\plot -4 2  4 2 /

\setdashes <1mm>
\plot  -4 0  0 2  /
\plot   4 0  0 2 /
\plot -5 1  0 2 /
\plot  5 1  0 2 /
\plot 0 2  1 5 /
\plot  0 2  -1 5 /

\setdashes <1mm>
\plot  0 2  2.667  3.33 /
\plot   0 2  -2.667  3.33 /
\plot 0 2  3.333 2.667 /
\plot  0 2  -3.333 2.667 /
\plot -0.667 0  0 2 /
\plot  0.667 0  0 2 /

\setsolid
\put{$\bullet$} at  6 0  
\put{$\bullet$} at 0 6
 
\put{$\ss-\frac13$} at 5.7 1
\put{$\ss0$} at -4.7 2
\put{$\ss\frac12$} at -5.7 1

\put{$\ss1$} at -6.2 -.4
\put{$\ss2$} at -4.1 -.4
\put{$\ss\infty$} at -2 -.4
\put{$\ss-3$} at -1 -.4
\put{$\ss-2$} at -0.05 -.4
\put{$\ss-\frac32$} at .9 -.4
\put{$\ss-1$} at 2 -.4
\put{$\ss-\frac23$} at 4.1 -.4
\put{$\ss-\frac12$} at 6 -.4

\endpicture}
$$

Actually, the 12 lines in $\mathbb L(6)$ 
can be described in several ways. The following table mentions three slightly 
different possibilities.
The first column lists the slope $\phi$, where $L_\phi = \frac{u-2b}{w-2b}$, 
the second column rewrites this description by multiplying with the denominators.
Thus, the first two columns describe the lines by refering to $u,w,b$.  
The third column uses in addition $\omega = |\Omega V|$.
We obtain the third column from the second one by adding or subtracting the equality
$\omega+u+w = 6b$ mentioned in Section~\ref{sec-three-six}. 

$$
{\beginpicture
   \setcoordinatesystem units <3cm,.4cm>
\put{$\ \phi$} at 0.2 1.2
\put{using $u,w,b$} at 1.15 1.2
\put{using also $\omega$} at 2.08 1.2
\put{type} at 2.8 1.2

\put{$\ \;0/1$} at 0.2 0
\put{$\ \;1/2$} at 0.2 -1
\put{$\ \;1/1$} at 0.2 -2
\put{$\ \;2/1$} at 0.2 -3
\put{$\ \;1/0$} at 0.2 -4
\put{$-3/1$} at 0.2 -5
\put{$-2/1$} at 0.2 -6
\put{$-3/2$} at 0.2 -7
\put{$-1/1$} at 0.2 -8
\put{$-2/3$} at 0.2 -9
\put{$-1/2$} at 0.2 -10
\put{$-1/3$} at 0.2 -11

\put{$u$\strut} [r] at 0.89 0
\put{$$\strut} [r] at 1 0
\put{$$\strut} [r] at 1.15 0 
\put{$=$\strut} [r] at 1.27 0
\put{$2b$\strut} [r] at 1.5 0 

\put{$2u$\strut} [r] at 0.89 -1
\put{$-$\strut} [r] at 1 -1
\put{$w$\strut} [r] at 1.15 -1 
\put{$=$\strut} [r] at 1.27 -1
\put{$2b$\strut} [r] at 1.5 -1 

\put{$u$\strut} [r] at 0.89 -2
\put{$-$\strut} [r] at 1 -2
\put{$w$\strut} [r] at 1.15 -2
\put{$=$\strut} [r] at 1.27 -2
\put{$0$\strut} [r] at 1.5 -2 

\put{$u$\strut} [r] at 0.89 -3
\put{$-$\strut} [r] at 1 -3
\put{$2w$\strut} [r] at 1.15 -3
\put{$=$\strut} [r] at 1.27 -3
\put{$-2b$\strut} [r] at 1.5 -3 

\put{$$\strut} [r] at 0.89 -4
\put{$$\strut} [r] at 1 -4
\put{$w$\strut} [r] at 1.15 -4
\put{$=$\strut} [r] at 1.27 -4
\put{$2b$\strut} [r] at 1.5 -4

\put{$u$\strut} [r] at 0.89 -5
\put{$+$\strut} [r] at 1 -5
\put{$3w$\strut} [r] at 1.15 -5
\put{$=$\strut} [r] at 1.27 -5
\put{$8b$\strut} [r] at 1.5 -5

\put{$u$\strut} [r] at 0.89 -6
\put{$+$\strut} [r] at 1 -6
\put{$2w$\strut} [r] at 1.15 -6
\put{$=$\strut} [r] at 1.27 -6
\put{$6b$\strut} [r] at 1.5 -6

\put{$2u$\strut} [r] at 0.89 -7
\put{$+$\strut} [r] at 1 -7
\put{$3w$\strut} [r] at 1.15 -7
\put{$=$\strut} [r] at 1.27 -7
\put{$10b$\strut} [r] at 1.5 -7

\put{$u$\strut} [r] at 0.89 -8
\put{$+$\strut} [r] at 1 -8
\put{$w$\strut} [r] at 1.15 -8
\put{$=$\strut} [r] at 1.27 -8
\put{$4b$\strut} [r] at 1.5 -8

\put{$3u$\strut} [r] at 0.89 -9
\put{$+$\strut} [r] at 1 -9
\put{$2w$\strut} [r] at 1.15 -9
\put{$=$\strut} [r] at 1.27 -9
\put{$10b$\strut} [r] at 1.5 -9

\put{$2u$\strut} [r] at 0.89 -10
\put{$+$\strut} [r] at 1 -10
\put{$w$\strut} [r] at 1.15 -10
\put{$=$\strut} [r] at 1.27 -10
\put{$6b$\strut} [r] at 1.5 -10

\put{$3u$\strut} [r] at 0.89 -11
\put{$+$\strut} [r] at 1 -11
\put{$w$\strut} [r] at 1.15 -11
\put{$=$\strut} [r] at 1.27 -11
\put{$8b$\strut} [r] at 1.5 -11

\put{$u =  2b$} [l] at 1.75 0
\put{$w = 2(u-b)$} [l] at 1.75 -1
\put{$w = u$} [l] at  1.75 -2
\put{$u = 2(w-b)$} [l] at  1.75 -3
\put{$w = 2b$} [l] at  1.75 -4
\put{$\omega = 2(w-b)$} [l] at  1.75 -5
\put{$\omega = w $} [l] at  1.75 -6
\put{$w = 2(\omega-b)$} [l] at  1.75 -7
\put{$\omega = 2b$} [l] at  1.75 -8
\put{$u = 2(\omega-b) $} [l] at  1.75 -9
\put{$u = \omega$} [l] at  1.75 -10
\put{$\omega = 2(u-b)$} [l] at  1.75 -11

\put{$\mathbb P$} at 2.8 0
\put{$\mathbb H$} at 2.8 -1
\put{$\mathbb D$} at 2.8 -2
\put{$\mathbb H$} at 2.8 -3
\put{$\mathbb P$} at 2.8 -4
\put{$\mathbb H$} at 2.8 -5
\put{$\mathbb D$} at 2.8 -6
\put{$\mathbb H$} at 2.8 -7
\put{$\mathbb P$} at 2.8 -8
\put{$\mathbb H$} at 2.8 -9
\put{$\mathbb D$} at 2.8 -10
\put{$\mathbb H$} at 2.8 -11

\endpicture}
$$

\medskip
Note that the strange numbers which are needed to describe the lines in $\mathbb L(6)$
when we use the first or the second column, are replaced in the third column by a concise
denomination. 
It is remarkable that the third column provides {\bf just three} kinds of equations
which define the lines in $\mathbb L(6)$,
one for each of the types $\mathbb P, \mathbb D, \mathbb H$: 
For $\mathbb P$, one of the coordinates has to have
the constant value $2b$.
For $\mathbb D$, the equation asserts that 
two of the three coordinates
$u,w,\omega$ coincide. Finally for $\mathbb H$, there are the six equations of the form $w = 2(u-b)$
(clearly, these latter equations deserve some further interest).  
Recall that Theorem~\ref{theoremone}
asserts that if $X$ is an indecomposable object in $\Cal S(n)$, which is
not a boundary picket, then $u \ge b$, $w\ge b$, and $\omega \ge b.$ Thus, it seems
reasonable to work (as we do here) with the differences $u-b$, $w-b$, $\omega-b$. Looking 
at the lines of type $\mathbb H$, the equations of the form $w = 2(u-b)$ show in which way
one of the coordinate value, here $w$, is determined by another coordinate value, here $u$,
namely $w$ is the double of $u-b.$ 

Let us stress the following consequence of the table:

\medskip
\begin{proposition}
  \label{prop-eleven-seven}
  Let $X$ be indecomposable in $\Cal S(6)$. If $u < w < 2b$, or if $2b < w < u,$
  then $u = 2(w-b),$ thus $w = \frac12(u+2b)$ is the average of $u$ and $2b$.
\end{proposition}

\medskip
This concerns the pr-vectors in the interior of the following two shaded triangles:

$$
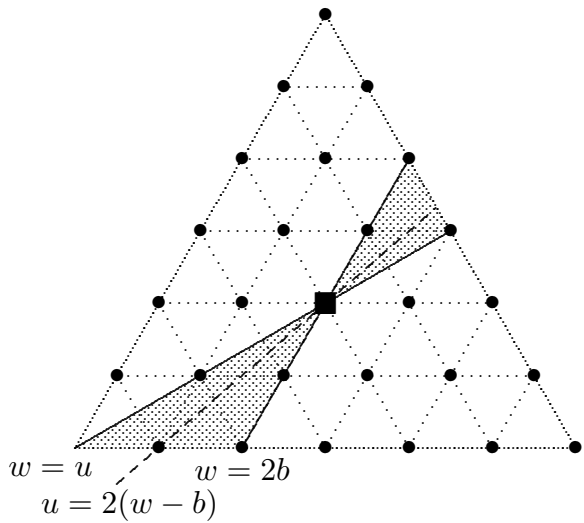
{\beginpicture
   \setcoordinatesystem units <.57735cm,1cm>
   \setcoordinatesystem units <.46188cm,.8cm>
\multiput{} at -6 -.5  6 6 /
\setdots <1mm> 
\plot -6 0  6 0  0 6  -6 0 /
\plot -5 1  -4 0  1 5  -1 5  4 0  5 1  -5 1 /
\plot -4 2  -2 0  2 4  -2 4  2 0  4 2  -4 2 /
\plot -3 3  0 0  3 3  -3 3 /

\setdots <.5mm> 
\plot -6 0  6 0  0 6  -6 0 /

\multiput{$\bullet$} at 
     -4 0  -2 0  0 0  2 0  4 0  
     -3 1  -1 1  1 1  3 1 
     -2 2  0 2  2 2
     -1 3  1 3
      0 4
       -5 1  -4 2  -3 3  -2 4  -1 5  1 5  2 4  3 3  4 2  5 1  /
\put{$\blacksquare$} at 0 2
\setsolid
\plot -6 0  3 3 /
\plot -2 0  2 4 /

\setdashes <1mm>
\plot  0 2  2.667  3.33 /

\plot -5 -.5  0 2 /
\setshadegrid span <.4mm>
\vshade -6 0 0  <,z,,>  -2 0 1.333   <z,z,,>  0 2 2  <z,z,,> 2 2.667 4 <z,z,,> 3 3 3 /
\setsolid
\put{$\bullet$} at  6 0  
\put{$\bullet$} at 0 6
 
\put{$w=u$} at -6.6 -.3
\put{$u = 2(w-b)$} at -4.7 -0.8
\put{$w = 2b$} at -2 -.3

\endpicture}
$$

\medskip
There are five similar assertions, obtained via the 
$\Sigma_3$-symmetries (or, of course, using again the table).

\subsection{The triangle support of $\Cal S(6)$.
  {\rm\normalsize (Proof of Theorem~\ref{theoremeight}(b))}.}
\label{sec-eleven-eight}

We want to show: {\it The set $\Psi=\Psi(\Cal S(6))$ is a monoton increasing sequence in $[0,2[$
which  converges to $2$.}

\medskip 
We are going to define the primitive pairs $(u,b)$ for $\Cal S(6)$.

\medskip
The algebra $\Xi$ is a tilted algebra of type $\mathbf E_8.$ 
If $M$ is an indecomposable $\Theta$-module, then $\bdim M$
is a radical vector or a 
positive root of $\chi_\Theta$, and we can write
$$
  \bdim M =  \mathbf r(M) + a_0\mathbf h_0 + a_\infty\mathbf h_\infty,
$$
where $\mathbf r(M)$ is either zero or a (not necessarily positive) root of $\chi_\Xi$, and where 
$$
 \mathbf h_0 = \smallmatrix  &  & 2 & 1 & 0 \cr
             1 & 2 & 3 & 3 & 2 & 1& 0 \endsmallmatrix, \quad 
 \mathbf h_\infty = \smallmatrix  &  & 2 & 2 & 1  \cr
             0 & 1 & 2 & 3 & 3 & 2 & 1 \endsmallmatrix
$$
(here, $a_0 = \dim M_0$ and $a_\infty = \dim M_{4'}$).
Recall that $\pi$ provides a bijection between $'\Cal D$ and the class $\Cal X$ of
indecomposable objects which belong to stable tubes. 

\medskip
For any root $\mathbf r$ of $\chi_\Xi,$ let $'\Cal D(\mathbf r)$ be the class of
indecomposable modules $M$ in $'\Cal D$ such that $\mathbf r(M) = \mathbf r.$ 
Let $M(\mathbf r)$ be an indecomposable module in $'\Cal D(\mathbf r)$ of minimal length.
Let 
$$
    u'(\mathbf r) = u(\pi M(\mathbf r)), \quad b'(\mathbf r) = b(\pi M(\mathbf r)),
   \quad a(\mathbf r) = a_0+a_\infty,
$$ 
clearly, $u'(\mathbf r),b'(\mathbf r),a(\mathbf r)$ only depend on $\mathbf r$, and 
for any object $Y$ in $'\Cal D(\mathbf r)$, there is $t\in \mathbb N_0$ such that
$$
  u(\pi Y) = u'(\mathbf r)+6t, \quad\text{and}\quad b(\pi Y) = b'(\mathbf r)+3t;
$$ 
namely, if $\bdim Y = \mathbf r + a'_0\mathbf h_0 + a'_\infty\mathbf h_\infty,$
then  $t = a'_0+a'_\infty-a(\mathbf r).$ 

\bigskip
If $M(\mathbf r)$ is $u$-minimal, let
$$
 d_{\mathbf r}(t) = (u'(\mathbf r)+6t)/(b'(\mathbf r)+3t);
$$
this is a function $\mathbb N_0 \to \mathbb R.$ 

\medskip
\begin{lemma}
  For any indecomposable object $X$ in $\Cal S(n)$, there is a root $\mathbf r$
  such that $M(\mathbf r)$ is $u$-minimal and some $t\in \mathbb N_0$
  such that $dX = d_{\mathbf r}(t).$
\end{lemma}

\begin{proof}
  We know this already for $X$ belonging to a stable component. 
  We claim that the primitive pairs $(u,b)$ for $\Cal P(6)$ are primitive pairs for $\Cal S(6).$
  There are the following five primitive pairs for $\Cal P(6)$,
  see Proposition~\ref{prop-ten-six}.
  Below each
  pair $(u,b)$, we show an object $M(\mathbf r)$ which is $u$-minimal and such that
  $u = u'(\mathbf r)$ and $b = b'(\mathbf r).$ 
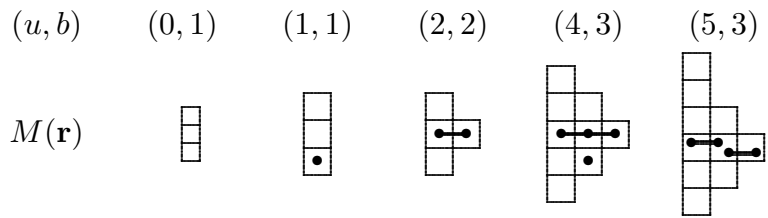
  $$
  \hbox{\beginpicture
    \setcoordinatesystem units <1.5cm,1.2cm>
    \put{$(u,b)$} at 0 0 
    \put{$(0,1)$} at 1 0 
    \put{$(1,1)$} at 2 0 
    \put{$(2,2)$} at 3 0 
    \put{$(4,3)$} at 4 0 
    \put{$(5,3)$} at 5 0 
    \put{$M(\mathbf r)$} at 0 -1  
    \put{\beginpicture
      \setcoordinatesystem units <.3cm,0.3cm>
      \multiput{} at 0 0  0 3 /
      \plot 0 0  1 0  1 3  0 3  0 0 /
      \plot 0 1  1 1 /
      \plot 0 2  1 2 /
      \endpicture} at 1 -1
    \put{\beginpicture
      \setcoordinatesystem units <.3cm,0.3cm>
      \multiput{} at 0 0  1 3 /
      \plot 0 0  1 0  1 3  0 3  0 0 /
      \plot 0 1  1 1 /
      \plot 0 2  1 2 /
      \put{$\ss\bullet$} at 0.5 0.5
      \endpicture} at 2 -1
    \put{\beginpicture
      \setcoordinatesystem units <.3cm,0.3cm>
      \multiput{} at 0 0  2 3 /
      \plot 0 0  1 0  1 3  0 3  0 0 /
      \plot 0 1  2 1  2 2   0 2  /
      \multiput{$\ss\bullet$} at 0.5 1.5  1.5 1.5 /
      \plot   0.5 1.5  1.5 1.5 /
      \plot   0.5 1.45  1.5 1.45 /
      \plot   0.5 1.55  1.5 1.55 /
      \endpicture} at 3 -1
    \put{\beginpicture
      \setcoordinatesystem units <.3cm,0.3cm>
      \multiput{} at 0 0  3 5 /
      \plot 0 0  1 0  1 5  0 5  0 0 /
      \plot 0 1  2 1  2 4  0 4 /
      \plot 0 2  3 2  3 3  0 3  /
      \multiput{$\ss\bullet$} at 0.5 2.5  1.5 2.5  2.5 2.5  1.5 1.5 /
      \plot   0.5 2.5  2.5 2.5 /
      \plot   0.5 2.45  2.5 2.45 /
      \plot   0.5 2.55  2.5 2.55 /
      \endpicture} at 4 -1
    \put{\beginpicture
      \setcoordinatesystem units <.3cm,0.3cm>
      \multiput{} at 0 0  3 6 /
      \plot 0 0  1 0  1 5  0 5  0 0 /
      \plot 0 1  2 1  2 4  0 4 /
      \plot 0 2  3 2  3 3  0 3  /
      \plot 0 5  0 6  1 6  1 5 /
      \multiput{$\ss\bullet$} at 0.3 2.7  1.3 2.7  2.7 2.3  1.7 2.3 /
      \plot   0.3 2.7  1.3 2.7  /
      \plot   0.3 2.75  1.3 2.75  /
      \plot   0.3 2.65  1.3 2.65  /
      \plot    2.7 2.3  1.7 2.3  /
      \plot    2.7 2.35  1.7 2.35  /
      \plot    2.7 2.25  1.7 2.25  /
      \endpicture} at 5 -1
    \endpicture}
  $$
  All the objects $M(\mathbf r)$ belong to stable components, since 
  the global space of an object in $\Cal P(6)$ has no part of the form $[3]$.
\end{proof}

\medskip
\begin{corollary}
  The index set $\Psi$ is contained in the set of values of the
  functions $d_{\mathbf r}(t)$, where $\mathbf r$ is a root and $M(\mathbf r)$ is $u$-minimal.
  
  For any root $\mathbf r$ with $M(\mathbf r)$ being minimal, 
  the set of values $d_{\mathbf r}(t)$ with $t\in \mathbb N_0$
  is a strictly increasing sequence which converges to $2$. $\s$
\end{corollary}

\medskip 
Theorem~\ref{theoremeight}(b)
follows from the Corollary, since the number of roots $\mathbf r$ is finite. 

\subsection{More on the set $\Psi$.}
\label{sec-eleven-nine}

\begin{proposition}
  If $d\in\Psi$, then $d = 2-c/b$, where
  $c$ is one of the numbers $1,2,3,4$ and $b\in \mathbb N_1$.
\end{proposition}

\begin{proof}
  Let $d = u/b$, where $u = uX,\ b = bX$ for some $u$-minimal object 
  $X$. Let $c = 2b-u.$ 
  The $u$-minimality implies that $d = u/b < 2,$ therefore $u < 2b,$
  thus $c > 0.$ 
  
  Now, $d = u/b = (2b-c)/b = 2-c/b.$ 
  According to Theorem~\ref{theoremseven}, we have $-4 \le u-2b \le 4,$ thus $-4 \le c \le 4.$
  Altogether we see that $0 < c \le 4.$ Since $c$ is an integer, 
  $c\in\{1,2,3,4\}.$
\end{proof}

\medskip
Here is the list of the numbers  $d = 2-c/b = (2b-c)/b$ with $1\le d\le 5/3.$ It is
sufficient to look at the cases $c = 3$ and $c = 4.$
	\smallskip
        
With $c = 3$: \quad $3/3,\ 5/4,\ 7/5,\ 9/6,\ 11/7,\ 13/8,\ 15/9.$
	\smallskip
        
With  $c = 4$: \quad $4/4,\ 6/5,\ 8/6,\ 10/7,\ 12/8,\ 14/9,\ 16/10,\ 18/11,\ 20/12.$
	\smallskip
        
Not all these numbers belong to $\Psi$. For example, there is no
indecomposable object $X\in\Cal S(6)$ on $\Delta_d$ with $d = 6/5,$ since $X$ would belong 
to a stable tube and we would have $\bdim \widetilde X = \mathbf r + \mathbf h$ for a root
$\mathbf r$ and a radical vector $\mathbf h$ with $u\mathbf h = 6,\ b\mathbf h = 3,$ thus
$u\mathbf r = 0$ and $b\mathbf r = 2$. But such a root $\mathbf r$ does not exist. 

On the other hand, we note that $d = 10/7$ does belong to $\Psi$, 
since $dX^+ = 10/7$, where $X = ([3,1],[6,4,3,1])$ (we have $uX = 4,\ bX = 4,$ thus
$uX^+ = 10,\ bX^+ = 7$).

\subsection{Proof of Theorem~\ref{theoremeight}(c): The number is unbounded.}
\label{sec-eleven-ten}

For $a_0,\ a_\infty \in \mathbb N_0,$ there is an indecomposable $\Xi$-module $M(a_0,a_\infty)$ with
$$
  \bdim M(a_0,a_\infty) =  \mathbf r + a_0\mathbf h_0 + a_\infty\mathbf h_\infty,
$$
where 
$$
 \mathbf r =  \smallmatrix  &  & 0 & 0 & 0 \cr
             0 & 0 & 0 & 1 & 0 & 0& 0 \endsmallmatrix,
$$
and we let $X(a_0,a_\infty) = \pi M(a_0,a_\infty)$.
For example, the module $M(0,0) = S_3$ is the simple module corresponding to the vertex $3$,
thus $X(0,0) = S$ and $X(1,0) = S[7].$ 

Let $t = a_0+a_\infty$. The uwb-vector of $X(a_0,a_\infty)$ is $\dfrac{6t|6t+1}{3t+1}.$ 
We obtain in this way $t+1$ indecomposable objects with the same pr-vector 
$(6t/(3t+1),(6t+1)/(3t+1)).$

\subsection{Small objects in $\Cal S(6)$.}
\label{sec-eleven-eleven}

We are going to show that Theorem~\ref{theoremseven}
implies strong restrictions on the existence of 
indecomposable objects. In particular, it provides an interesting characterization of the
object $S[6]$ in $\Cal S(6)$. We recall that $S[6]$ is the unique indecomposable object 
with a sectional path from $S = (0,[1])$
to $S[6]$ of length 5. Of course, $\Cal S[6]$ belongs to the principal component $\Cal P(6)$; its
partition vector is $([4,2],[6,4,1,1],[4,2])$ and $qS[6]=3$.

\medskip
\begin{proposition}
  Let $X$ be indecomposable in $\Cal S(6).$ 
  \smallskip
        
  If $qX < 16/5$, then: Either $X$ has height at most $5$ or else $X = S[6]$.
  \smallskip
        
  If $qX < 10/3$, then $bX \le 5$ and $vX \le 16.$
\end{proposition}

\medskip
All the indecomposable objects in $\Cal S(6)$ with $bX \le 5$ and $qX < 10/3$
will be exhibited in Section~\ref{sec-fifteen-two} and Appendix~\ref{app-B}.
In this way, we see that the bounds are optimal. 

\smallskip
\begin{proof}
  According to Theorem~\ref{theoremseven}, we have $|vX-4\cdot bX|\le 4,$ thus
$|qX - 4| \le 4/bX$ and thus $4/bX \ge 4 -qX $.

First, let $qX < 10/3,$ thus $4-qX > 4- 10/3 = 2/3$. We get 
$4/bX \ge 4 -qX > 2/3,$ so that $bX < 6$. Since $bX$ is
an integer, we have $bX \le 5.$ 
We get $vX = qX\cdot bX < 10/3\cdot 5 = 50/3,$ thus $vX \le 16.$
This finishes the proof of the second assertion. 

\medskip 
In order to show the first 
assertion, we assume now that $qX < 16/5,$ thus $4-qX > 4- 16/5 = 4/5$. We get 
$4/bX \ge 4 -qX > 4/5,$ so that $bX < 5$, and therefore $bX \le 4.$ 

Let us show that for $bX \le 3$, the object $X$ has height at most 5.
If $bX = 1,\ 2,$ or $3$, we have $vX = qX\cdot bX < 16/5\cdot bX,$ thus $vX \le 3,\ 6,$ or $9$,
respectively. If $bX = 1$ or $2$, then clearly $X$ has height at most 5. 
Thus let us consider the case $bX = 3$, so that $vX \le 9.$ Let us assume that
$X$ has height $6$, thus $VX = [6,2,1]$ or $[6,1,1].$ But this means that $U(\tau X)$
is equal to $[2,1]$ or to $[1,1],$ thus $U(\tau X)$ 
has height at most $2$. 
But for any indecomposable object in $\Cal S$ whose subspace has height at most 2,
the subspace has to be indecomposable. 

Thus, we can assume that $bX = 4.$ 
We get $v = qX\cdot bX < 16/5\cdot 4 = 65/5,$ thus $vX \le 12.$
We assume now again that $X$ has height $6$, thus either $VX = [6,4,1,1],$ or else
$[6,x,y,z]$ with $3\ge x \ge y \ge z \ge 1$.
We are going to show that only the case $[6,4,1,1]$ is possible. 
Thus, assume that $VX = [6,x,y,z]$ with $3\ge x \ge y \ge z \ge 1$ and let
$Y = \tau X$. According to Section~\ref{sec-three}, we have $UY = [x,y,z],$ thus $Y$ belongs to
$\Cal S_3(6)$. Again we use the argument that $UY$ must have height 3, thus 
$3 = x > y \ge z \ge 1.$
According to Schmidmeier \cite{S1}, the category $\Cal S_3(6)$ is known 
(it is representation-finite with
84 indecomposable objects, see Proposition 12 in \cite{S1}) and 
there are only two indecomposable objects $Y$ in $\Cal S_3(6)$
with subspace $UY = [x,y,z],$ where $3 = x > y \ge z \ge 1;$ 
in both cases $UY = [3,2,1]$. In order to
get a contradiction, we consider now $\tau Y$ (where again $\tau = \tau_6$).
It is easy to see that $b(\tau Y) = 5$. 
(For example, we can use Proposition~\ref{prop-three-twelve-one}:
Since $UY=[3,2,1]$, we have $bUY=3$.
From the objects in \cite{S1} we read off that in each case $bWY=3$ and $c_6Y=1$,
so the proposition yields $b(\tau Y)=bUY+bWY-c_6Y=5$.)
But this contradicts the fact that $b(\tau Y) = b(\tau^2X) = bX = 4$,
according to Theorem~\ref{theoremthree}.
As a consequence, we conclude that $V = [6,4,1,1].$

There are several ways to show that the only indecomposable object $X$ in $\Cal S$ with
$VX = [6,4,1,1]$ is $X = S[6].$ Let us outline two approaches. Using the root system of $\Theta$,
one easily sees that there is no indecomposable object in a stable tube with global space
$[4,2,1,1]$. Thus, it remains to look at the principal component $\Cal P(6)$. There are just
12 indecomposable objects in $\Cal P(6)$ which have width 4, and just one of them, namely
$S[6]$ has global space $[6,4,1,1].$

A second possible proof invokes Proposition~\ref{prop-five-one}: Given an indecomposable object
$X$ with $VX = [6,4,1,1]$, we obtain a g-split filtration of $X$ with two factors both
being extended pickets. This allows to identify $X$ with $S[6].$
\end{proof}

\medskip
\begin{remark}
  The Proposition and its proof provide the following characterizations of 
  $S[6]$.
  {\it Assume that $X$ is indecomposable in $\Cal S(6)$ and does not belong to $\Cal S(5)$.
    Then $X = S[6]$ if and only if $vX = 12,\ bX = 4,$ if and only if $qX = 3$, if and only if 
    $qX < 16/5,$ if and only if $VX = [6,4,1,1].$} $\s$
\end{remark}

\medskip
\begin{remark}
  There is a second interesting indecomposable object $X$ with width $4$,
  namely the object $X$ with $\bpar X = ([3,1],[6,4,3,1],[5,3,2]).$ It can be characterized as 
  follows: 
  {\it Assume that $X=(U,V)$ is indecomposable in $\Cal S(6)$ and $0\neq U$ is decomposable.
    Then $\bpar X = ([3,1],[6,4,3,1],[5,3,2])$ 
    if and only if $uX = 4,\ bX = 4,$ if and only if $pX = 1$, if and only if 
    $pX < 5/4,$ if and only if $UX = [3,1]$ and $bX=4$.} 
  For the proof, one uses Theorem~\ref{theoremeight}(b) and the fact
  that for $d\in \Psi$, we have $d \le 1$ or $d \ge 5/4,$ see Section~\ref{sec-eleven-nine}.
  See the line $p=1$ in Section \ref{sec-fifteen-two}--(g). $\s$
\end{remark}

\subsection{The function $\eta_n$ for $n\ge 7$.}
\label{sec-eleven-twelve}

We have introduced in Section~\ref{sec-eleven-four}
the function $\eta X = vX-4\cdot bX$ for $X$ on $\Cal S(6)$.
For arbitrary $n$, there is the corresponding function $\eta_n$ defined by
$\eta_n = vX-\frac{2n}3\cdot bX$ for $X$ in $\Cal S(n)$. 

\medskip
\begin{proposition}
  For $n\ge 7,$ the function 
  $\eta_n = vX-\frac{2n}3\cdot bX$ is not bounded on the indecomposable
  objects of $\Cal S(n).$
\end{proposition}

\begin{proof}
  Let $n\ge 7.$
  Any indecomposable object $X$ in $\Cal S(6)$ with $VX = [6,4,2]^t$ can be
  modified in order to produce indecomposable objects $X'$ and $X''$ in $\Cal S(n)$
  with $V(X') = [n,4,2]^t$ and $V(X'') = [n,n-2,n-4]^t.$
  
  Now $b(X') = b(X'') = 3t$ and $v(X') = (n+6)t$, whereas $v(X'') = (3n-6)t.$
  Thus
  \begin{align*}
    \eta_n(X') &= (n+6)t - \tfrac23n\cdot 3t  = (-n+6)t \ \le\;  -t, \cr
    \eta_n(X'') &=  (3n-6)t - \tfrac23n\cdot 3t= (n-6)t \ \ \ge\ t.
  \end{align*}
\end{proof}

\subsection{Central lines and Kronecker subcategories.}
\label{sec-eleven-thirteen}

As we have seen in Section \ref{sec-ten-one},
any non-central object in a stable tube $\Cal C$ of $\Cal S(n)$
gives rise to a pair of complementary central half-lines in the half-line support of
$\Cal C$. For $n = 6$, there is an additional way that complementary central half-lines
appear, namely as the half-line support of some Kronecker subcategories of 
$\widetilde{\Cal S}(6)$. 

\medskip
We say that a Kronecker-pair $X,Y$ in $\widetilde{\Cal S}(6)$ is {\it g-split}
provided $\Ext^1(X,Y) = \Ext^1_{\text g}(X,Y).$ 

\medskip
\begin{proposition}
  \label{prop-eleven-thirteen}
  Let $X,Y$ be a Kronecker pair in $\widetilde{\Cal S}(6)$ which
  is not central and g-split. Then the half-line support of the corresponding Kronecker
  subcategory $\Cal K$ is a pair of complementary central lines.
\end{proposition}

\medskip
The examples presented in Sections~\ref{sec-two-seven} and \ref{sec-two-eight}
are non-central and g-split Kronecker pairs.
For the example in Section~\ref{sec-two-seven}, with $\uwbb X = \frac{4|2}1$ and  $\uwbb Y = \frac{2|4}2$,
we obtain as half-line support the line $L_\infty$,
for the example in Section~\ref{sec-two-eight},
with $\uwbb X =\uwbb{\widetilde A\atop{\widetilde B}'} = \frac{8|7}3$
and  $\uwbb Y=\uwbb{\widetilde B\atop\widetilde C} = \frac{4|5}3$, we obtain  
the line $L_2$. 
$$  
{\beginpicture
   \setcoordinatesystem units <1cm,1cm>
\put{\beginpicture
   \setcoordinatesystem units <.57735cm,1cm>
\put{Example 2.7} at -3 4
\multiput{} at -6 .4  6 4.4 /
\setdots <.3mm>
\plot -3 1  3 1  0 4  -3 1 /
\multiput{$\ssize \blacksquare$} at 0 2  /
\setdots <1mm> 

\plot -2.5 1.5  -2 1 .5 3.5  -.5 3.5  2 1  2.5 1.5  -2.5 1.5 /
\plot -2 2      -1 1  1 3  -1 3  1 1  2 2  -2 2 /
\plot -1.5 2.5  0 1  1.5 2.5  -1.5 2.5 /
\multiput{$\bullet$} at 1 3  -.5 1.5 /
\put{$X$} at 1.5 3 
\put{$Y$} at 0 1.5 
\setsolid  
\plot -1.5 .5  1.5 3.5 /
\put{$L_\infty$} at -.8 .5 
\endpicture} at 0 0 
\put{\beginpicture
   \setcoordinatesystem units <.57735cm,1cm>
\multiput{} at -6 .4  6 4.4 /
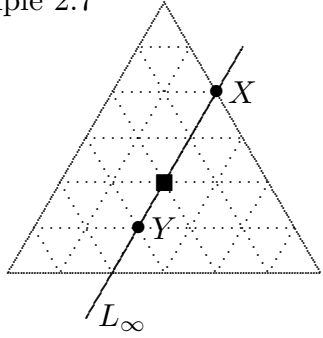
\put{Example 2.8} at -3 4
\setdots <.3mm>
\plot -3 1  3 1  0 4  -3 1 /
\multiput{$\ssize \blacksquare$} at 0 2  /
\setdots <1mm> 

\plot -2.5 1.5  -2 1 .5 3.5  -.5 3.5  2 1  2.5 1.5  -2.5 1.5 /
\plot -2 2      -1 1  1 3  -1 3  1 1  2 2  -2 2 /
\plot -1.5 2.5  0 1  1.5 2.5  -1.5 2.5 /
\multiput{$\bullet$} at 0.588 2.28  -.65 1.65 /
\put{$X$} at  1.1 2.28 
\put{$Y$} at -.2 1.5 
\setsolid  
\plot -2.5 0.75  2 3 /
\put{$L_2$} at -1.8 .5

\endpicture} at 5 0 
\endpicture}
$$

\vfill\eject

\centerline{\Gross Fifth part: Some objects.}
\addcontentsline{toc}{part}{Fifth part: Some objects.}

\section{Pickets and bipickets.}
\label{sec-twelve}

\subsection{The pickets (and their $\tau$-orbits).}
\label{sec-twelve-one}

We recall that an object $X$ in $\Cal S$ is a picket if and only if $bX = 1.$

\medskip
{\it The pickets in $\Cal S(n)$ correspond bijectively to} (and may be identified with) 
{\it the triples $(a_0,a_1,a_2)$ of
non-negative numbers $a_0,a_1,a_2$ with $a_0+a_1+a_2 = n$ and $(a_1,a_2) \neq (0,0)$;} 
the picket corresponding to the triple $(a_0,a_1,a_2)$ is $([a_1],[a_1+a_2])$ and it will
be denoted now also by $a_0\backslash a_1\backslash a_2.$ This denomination 
corresponds to the general frame established in Section~\ref{sec-three-one}; with the abbreviation
$a_0\backslash a_1\backslash a_2 = ([a_0]\backslash [a_1]\backslash [a_2]) = E(([a_1],[a_1+a_2])).$
The reduced pickets correspond to the triples $a_0\backslash a_1\backslash a_2$
different from  $0\backslash 0\backslash n$  and  $0\backslash n\backslash 0,$ 
thus to the triples $a_0,a_1,a_2$ of non-negative numbers with $a_0+a_1+a_2 = n$, and such that 
at most one of the numbers is zero. 

{\it If $a_0\backslash a_1\backslash a_2$ is a reduced picket, then}
$$
\tau_n^2 (a_0\backslash a_1\backslash a_2) = a_1\backslash a_2\backslash a_0,
$$
as asserted by Theorem~1.3$'$ in Section~\ref{theoremthree-prime}
(but also easily verified). 
(Note that Theorem~\ref{theoremthree} shows: If $X$ is a reduced picket, so that $X$ is
reduced and $bX = 1$, then also $b(\tau_n^2 X) = 1$, thus $\tau_n^2 X$ is again a picket; 
here is the precise formula.)
	
{\it If one and only one of the numbers $a_0,a_1,a_2$ is zero, then 
$\tau (a_0\backslash a_1\backslash a_2)$ is again a picket.}
Namely, $\tau_n (0\backslash a_1\backslash a_2) = ([a_2],[a_2])$; next
$\tau_n^2 (a_0\backslash 0\backslash a_2) = ([a_2],[n])$; 
and finally, $\tau (a_0\backslash a_1\backslash 0) = (0,[a_0]).$

If all numbers $a_0,a_1,a_2$ are non-zero, 
the object $\tau_n (a_0\backslash a_1\backslash a_2)$
is no longer a picket, but a bipicket $(U,V,W)$ with the following
properties: One part of the total space $V$ is $[n]$,
and both $U$ and $W$ are indecomposable $\Lambda$-modules, namely 
$\tau_n (a_0\backslash a_1\backslash a_2) = (U,V,W) = ([a_1+a_2],[n,a_2],[a_0+a_2])$,
and is constructed as follows:
There is a canonical embedding $\mu\:[a_1+a_2] \to [n]$ and a canonical projection 
$\pi\:[a_1+a_2] \to [a_2]$, in $\Cal N(n)$, 
thus there is the map $(\mu,\pi)\:[a_1+a_2] \to [n,a_2]$.
The map $(\mu,\pi)$ is a monomorphism, say with image $U$, and the cokernel is 
$W = [a_0+a_2]$.
$$
{\beginpicture
    \setcoordinatesystem units <1.5cm,1cm>
\put{$[a_1\!+\!a_2]$} at 0 0 
\put{$[n]$} at 1 1
\put{$[a_2]$} at 1 -1 
\put{$[a_0\!+\!a_2]$} at 2 0 
\arr{0.3 0.3}{0.7 0.7} 
\arr{0.3 -.3}{0.7 -.7} 
\arr{1.3 0.7}{1.7 0.3} 
\arr{1.3 -.7}{1.7 -.3} 
\put{$\mu$} at .3 .7
\put{$\pi$} at .3 -.7
\endpicture}
$$
Then $\tau_n (a_0\backslash a_1\backslash a_2)  = (U,[n,a_2],W)$.

\medskip
\begin{proposition}
  \label{prop-twelve-one}         
  The number of pickets of height $n$ is $n+1.$ The number of pickets
  of height at most $n$ is $\binom{n+2}2 -1.$
\end{proposition}

\begin{proof}
  The pickets of height $n$ are the objects $([t],[n])$ with $0 \le t \le n,$ thus the
  number of these pickets is $n+1.$ The number of pickets of height at most $n$ is therefore 
  $\sum_{i=1}^n (t+1) =
  \binom{n+2}2 -1.$
\end{proof}

\subsection{Pickets and irreducible module varieties.}
\label{sec-twelve-one-b}

It is known that pickets and their $\tau$-translates play a particular
role in geometric representation theory:
Given $n$ (nilpotency index), $v$ ($=\dim V$) and $u\leq v$ ($u=\dim U$),
let $\mathbb V_n(u,v)$ be the affine variety consisting of all triples
$M=(M_0,M_1,h_M)$ where $M_0$ is a $u\times u$-matrix, $M_1$ a $v\times v$-matrix
and $h_M$ a $u\times v$-matrix such that $M_0^n=0$, $M_1^n=0$ and
$M_0h_M=h_MM_1$.  (By convention, if $d',d''\geq 0$ and either $d'=0$ or $d''=0$
then there is a unique $d'\times d''$-matrix which behaves like zero
with respect to multiplication.)
The group $\GL_{u,v}=\GL_u(k)\times\GL_v(k)$ acts on $\mathbb V_n(u,v)$ by
$$
  (g_0,g_1)\cdot(M_0,M_1,h_M)=(g_0M_0g_0^{-1},g_1M_1g_1^{-1},g_0h_Mg_1^{-1}).
$$
It has been shown by Bobinski in \cite[Theorem 1.1]{Bob} that $\mathbb V_n(u,v)$ is an
irreducible variety for each $n$, $u$ and $v$.  More precisely,
$\mathbb V_n(u,v)$ has a unique dense $\GL_{u,v}$-orbit.
In our notation, this orbit corresponds to the following embedding
where we write $u=an+b$ and $v=cn+d$ with $a,b,c,d\in\mathbb Z$ such
that $0\leq b,d<n$ hold. The pickets
$P=([n],[n])$ and $P'=(0,[n])$ are the indecomposable projective objects.

\medskip
$P^a\oplus {P'}^{c-a} \quad \text{if} \quad 0= b=d$

\smallskip
$P^a\oplus {P'}^{c-a}\oplus (0,[d]) \quad \text{if} \quad 0=b<d$

\smallskip
$P^a\oplus {P'}^{c-a}\oplus ([b],[b]) \quad \text{if} \quad 0<b=d$

\smallskip
$P^a\oplus {P'}^{c-a-1}\oplus ([b],[n]) \quad \text{if} \quad 0=d<b$

\smallskip
$P^a\oplus {P'}^{c-a}\oplus ([b],[d]) \quad \text{if} \quad 0<b<d$

\smallskip
$P^a\oplus {P'}^{c-a-1}\oplus ([b],[n,d],[n+d-b]) \quad \text{if} \quad 0<d< b$

\medskip
Note that the bipicket $([b],[n,d],[n+d-b])$ occurs as the $\tau_n$-translate
of the picket $n\!-\!b\backslash b\!-\!d\backslash d=([b\!-\!d],[b])$.
Thus, the pickets and their $\tau_n$-translates define exactly the
indecomposable direct
summands of embeddings which give rise to the dense orbits
in representation spaces.

\subsection{Bipickets.}
\label{sec-twelve-two}

Recall that an object $X$ in $\Cal S$ is called a bipicket
if and only if $X$ is indecomposable and $bX = 2.$

\medskip
\begin{proposition}
  \label{prop-twelve-two}
  Let $X$ be a bipicket.

  \begin{itemize}[leftmargin=3em]
  \item[\rm(a)] The object $X$ is gradable. If $X = \pi (\widetilde U,\widetilde V)$,
    then $\widetilde V = \widetilde V_1\oplus \widetilde V_2$ with
    $\widetilde V_1,$ and $\widetilde V_2$ indecomposable and such that 
    $\widetilde V_2$ is a subquotient of $\rad \widetilde V_1/\soc \widetilde V_1$;  
    in particular, we have $\End \widetilde V = k\times k.$
    
  \item[\rm(b)] The bipicket $X$ is uniquely determined by $\bpar X.$
  \end{itemize}
\end{proposition}

\medskip

Let $c_1,\dots,c_5$ be natural numbers, with $c_1$, $c_3$ and $c_5$ positive.
For 
$1\le i < j \le 5,$ write $c_{ij} = \sum_{t=i}^j c_t.$ We define $B(c_1,\dots,c_5) = (U,V)$
as follows: $V = [c_{15},c_{24}]$, with generators $v_1,v_2$ (such that $T^{c_{15}}v_1 = 0
= T^{c_{24}}v_2 = 0$). Let $u_1 = T^{c_{45}}v_1+T^{c_4}v_2$ and $u_2 = T^{c_{34}}v_2$
(thus $u_2 = 0$ if and only if $c_2 = 0$), and $U$ the submodule of $V$ generated by $u_1,u_2$.
$$
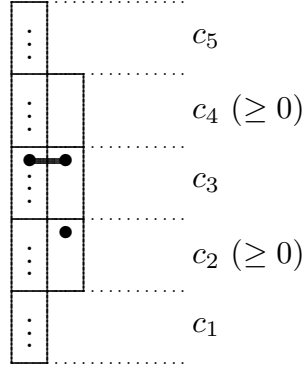
{\beginpicture
   \setcoordinatesystem units <.4cm,.4cm>
\multiput{} at 0 0  2 10 /
\plot 0 0  1 0  1 10  0 10  0 0 /
\plot 0 2  2 2  2 8  0 8 /
\plot 0 4  2 4 /
\plot 0 6  2 6 /
\multiput{$\vdots$} at 0.5 1.1  0.5 3.1  0.5 5.1  0.5  7.1  0.5 9.1 /
\multiput{$\bullet$} at 1.5 3.6  0.5 5.6  1.5  5.6 /
\plot  0.5 5.6  1.5  5.6 /
\plot  0.5 5.55  1.5  5.55 /
\plot  0.5 5.65  1.5  5.65 /
\put{$c_1$} [l] at 5 1
\put{$c_2\ (\ge 0)$} [l] at 5 3
\put{$c_3$} [l] at 5 5
\put{$c_4\ (\ge 0)$} [l] at 5 7
\put{$c_5$} [l] at 5 9
\setdots <1mm>
\plot 1 0  5 0 /
\plot 1 2  5 2 /
\plot 1 4  5 4 /
\plot 1 6  5 6 /
\plot 1 8  5 8 /
\plot 1 10  5 10 /
\endpicture}
$$
Note that $U$ is isomorphic to $[c_{13},c_2]$ and 
$W = V/U$ is isomorphic to $[c_{35},c_4]$,
so that
$$
   \bpar B(c_1,\dots,c_5) = ([c_{13},c_2],[c_{15},c_{24}],[c_{35},c_4]).
$$

\medskip
Of course, {\it all the objects $B(c_1,\dots,c_5)$ are gradable.}
The standard grading is the following:
The element $v_1$ has degree $c_{15},$ the element $v_2$ has degree $c_{14}$ (thus, $u_1$ has
degree $c_{13}$ and $u_2$ is zero or has degree $c_{12}$).

\medskip
\begin{proposition}
  \label{prop-twelve-two-two}
  The objects $B(c_1,\dots,c_5)$ are bipickets and any
  bipicket is of this form, for a unique sequence $c_1,\dots,c_5.$
\end{proposition}

\begin{proof}[Proof of Proposition~\ref{prop-twelve-two-two}]
First, let us assume that $(U,V,W)$ is a bipicket with both $U$ and $W$ being
cyclic. Let $V = V_1\oplus V_2$
with $V_1,V_2$ indecomposable. Let $\mu = \left(\smallmatrix \mu_1\cr\mu_2\endsmallmatrix
\right)\:U \to V_1\oplus V_2$ be the inclusion map. Since $U$ is not contained in
$\rad V_1\oplus \rad V_2$, at least one of the maps $\mu_1,\mu_2$ has to be surjective,
say $\mu_2.$ If $\mu_2$ would be even bijective, then $(U,V)$ would be decomposable.
Thus $c_1 = \dim \Ker \mu_2 \ge 1.$ 

Since $\mu_2$ is not injective, $\mu_1$ has to be injective. Again, $\mu_1$ cannot be 
bijective, since $(U,V)$ is indecomposable. Thus the dimension $c_5$ of the cokernel of
$\mu_1$ is at least 1. Let $v_1$ be a generator of $V_1.$ Then $\mu_1(U)$ is generated by
$u = T^{c_5}v_1$, and the element $v_2 = \mu_2(u)$ generates $V_2$. It follows that $(U,V)$ 
is isomorphic to $(\Lambda (u+v_2),\Lambda v_1\oplus \Lambda v_2) = B(c_1,0,c_3,0,c_5)$,
where $c_3$ is the dimension of $V_2$; in particular, we have  $c_3\ge 1.$

In this way, we 
see that the bipickets $(U,V,W)$ with $U$ and $W$ cyclic are the objects of the
form $B(c_1,0,c_3,0,c_5).$ 
	\smallskip
        
In general, the isomorphism classes of the indecomposable 
objects $(U,V,W)$ in $\Cal S(n)$ which satisfy
$U \subseteq \rad V$ and which are different from $(0,k)$ correspond
bijectively to the isomorphism classes of the indecomposable 
objects in $\Cal S(n-1)$,
by sending $(U,V)$ to $(U,\rad V)$ (and we have $bV = b(\rad V)$).
Under this bijection, the bipicket $B(c_1,c_2,c_3,c_4,c_5)$ with $c_4\ge 1$ corresponds 
to the bipicket $B(c_1,c_2,c_3,c_4-1,c_5).$ Note that we have 
$U \not\subseteq \rad V$ if and only if $bW < bV.$ 

Dually, the isomorphism classes of the indecomposable 
objects $(U,V,W)$ in $\Cal S(n)$ which satisfy
$\soc V \subseteq U$ and which are different from  $(k,k)$ correspond
bijectively to the isomorphism classes of the indecomposable 
objects in $\Cal S(n-1)$,
by sending $(U,V)$ to $(U/\soc V,V/\soc V)$ (and we  have $bV = b(V/\soc V)$). 
Under this bijection, the bipicket $B(c_1,c_2,c_3,c_4,c_5)$ with $c_2\ge 1$ corresponds
to the bipicket $B(c_1,c_2-1,c_3,c_4,c_5).$ Note that we have 
$\soc V \not\subseteq U$ if and only if $bU < bV.$ 
	\smallskip
        
Thus, we see that the bipickets $(U,V,W)$ are the objects of the
form $B(c_1,c_2,c_3,c_4,c_5);$ the number $c_4$ is the maximal number $t$ with $U\subseteq
\rad^tV$; the number $c_2$ is the maximal number $t$ with $\soc_t V \subseteq U.$ 
\end{proof}

\medskip
Proposition~\ref{prop-twelve-two} follows immediately from Proposition~\ref{prop-twelve-two-two}.
Note that
we easily recover the numbers $c_1,\dots,c_5$ from $\bpar B(c_1,\dots,c_5).$
$\s$

\medskip
\begin{remarks}
  We have seen in Proposition~\ref{prop-twelve-one}(b) that a bipicket $X = (U,V,W)$
is uniquely determined by $\bpar X = ([U],[V],[W])$. But {\it $X$ is not determined by two of the
three partitions $[U],[V],[W]$.} For example, the bipickets $B(1,0,2,0,2)$ and $B(2,0,1,1,1)$
have subspace $[3]$ and global space $[5,2]$, but the factor spaces are different. 
The bipickets $B(2,0,1,0,2,)$ and $B(1,0,2,0,1)$ have subspace and factor space of the form $[3]$,
but the global spaces are different. 

Also: {\it A bipicket $X = (U,V,W)$ in $\Cal S(n)$ 
is uniquely determined by $EX = ([\Omega V],[U],[W])$,}
since $[\Omega V]$ and $n$ together determine $[V]$.
\end{remarks}

\medskip
\begin{proposition}
  \label{prop-twelve-two-three}
  The number of bipickets in $\Cal S(n)$ 
  is $\binom{n+2}5$.
\end{proposition}

\begin{proof}
  The binomial coefficient $n+2\choose 5$ counts the number of five element subsets
  of $\{1,\ldots,n+2\}$, thus the number of sequences
  $1\le d_1 <  \cdots < d_5\le n+2$ of cardinality 5. 
  Given such a sequence $(d_i)_i$, let 
  $$
  c_1 = d_1,\quad c_2 = d_2-d_1-1,\quad c_3 = d_3-d_2,\quad c_4 = d_4-d_3-1, \quad c_5 = d_5-d_4.
  $$
  Then $c_i\ge 0$ for all $i$, and $c_i \ge 1$ for $i = 1,3,5.$ Also, $\sum_i c_i = (\sum_i d_i)-2 \le n.$
  We see that we can use the sequence $c_1,\dots,c_5$ in order to construct the bipicket $B(c_1,\dots,c_5)$
  and we obtain all the bipickets of height at most $n$ in this way.
\end{proof}

\medskip 
\begin{corollary}
  The number of bipickets of height $n$ is $\binom{n+1}4$.
\end{corollary}

\begin{proof} According to Proposition~\ref{prop-twelve-two-three},
  the number of bipickets in $\Cal S(n)$ 
  is $\binom{n+2}5$, thus the number of bipickets in $\Cal S(n-1)$ 
  is $\binom{n+1}5$. The assertion follows from the equality
  $$
  \tbinom{n+1}5 + \tbinom{n+1}4 = \tbinom{n+2}5.
  $$
\end{proof}

The distribution of the bipickets for $3\le n \le 7$ :
$$  
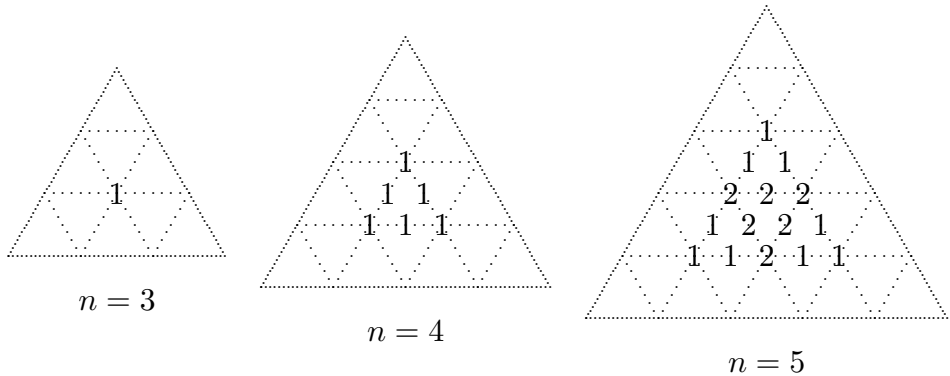
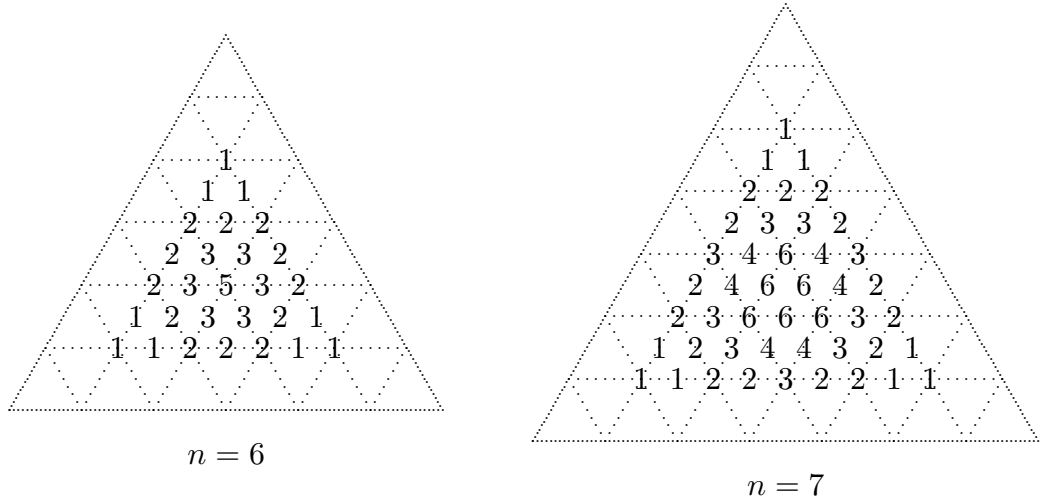
{\beginpicture
   \setcoordinatesystem units <.4cm,.69282cm>
\put{\beginpicture
\put{1} at 0 1
\multiput{} at -3 0  3 3 /
\setdots <.5mm>
\plot -3 0  3 0  0 3  -3 0 /
\setdots <1mm>
\plot -1 0  1 2  -1 2  1 0  2 1  -2 1  -1 0 /
\put{$n = 3$} at 0 -.7
\endpicture} at 0 0
\put{\beginpicture
\multiput{} at -4 0  4 4 /
\setdots <.5mm>
\plot -4 0  4 0  0 4  -4 0 /
\setdots <1mm>
\plot -2 0  1 3  -1 3  2 0  3 1  -3 1 -2 0 /
\plot 0 0  2 2  -2 2  0 0 /
\multiput{1} at -1 1  0 1  1 1  -.5 1.5  .5 1.5  0 2   /
\put{$n = 4$} at 0 -.7
\endpicture} at 8 0
\put{\beginpicture
\multiput{} at -5 0  5 5 /
\setdots <.5mm>
\plot -5 0  5 0  0 5  -5 0 /
\setdots <1mm>
\plot -3 0  1 4  -1 4  3 0  4 1  -4 1 -3 0 /
\plot -1 0  -3 2  3 2  1 0  -2 3  2 3  -1 0  /
\multiput{1} at -2 1  -1 1  1 1  2 1  -1.5 1.5  1.5 1.5  -.5 2.5  .5 2.5  0 3 
    /
\multiput{2} at 0 1  -.5 1.5  .5 1.5  -1 2  1 2  0 2 / 
\put{$n = 5$} at 0 -.7

\endpicture} at 18 0
\endpicture}
$$

$$  
{\beginpicture
   \setcoordinatesystem units <.4cm,.69282cm>
\put{\beginpicture
\multiput{} at -6 0  6 6 /
\setdots <.5mm>
\plot -6 0  6 0  0 6  -6 0 /
\setdots <1mm>
\plot -4 0  1 5  -1 5  4 0  5 1  -5 1 -4 0 /
\plot -2 0  -4 2  4 2  2 0  -2 4  2 4  -2 0  /
\plot 0 0  3 3  -3 3  0 0 /
\multiput{1} at -3 1  -2 1  2 1  3 1  -2.5 1.5  2.5 1.5  -.5 3.5  .5 3.5  0 4 
    /
\multiput{2} at -1 1  0 1  1 1  -1.5 1.5  1.5 1.5  -2 2  2 2  -1.5 2.5  1.5 2.5 
    -1 3  0 3  1 3 / 
\multiput{3} at -.5 1.5  .5 1.5  -1 2  1 2 -.5 2.5  .5 2.5 /
\multiput{5} at 0 2 / 
\put{$n = 6$} at 0 -.7
\endpicture} at 2 0
\put{\beginpicture
\multiput{} at -7 0  6 7 /
\setdots <.5mm>
\plot -7 0  7 0  0 7  -7 0 /
\setdots <1mm>
\plot -5 0  1 6  -1 6  5 0  6 1  -6 1 -5 0 /
\plot -3 0  -5 2  5 2  3 0  -2 5  2 5  -3 0  /
\plot -1 0  -4 3  4 3  1 0  -3 4  3 4  -1 0 /
\multiput{1} at -4 1  -3 1  3 1  4 1  -3.5 1.5  3.5 1.5  -.5 4.5  .5 4.5  0 5 /
\multiput{2} at -2 1  -1 1  1 1  2 1  
    -2.5 1.5  2.5 1.5  -3 2  3 2  
    -1.5 3.5  1.5 3.5 
    -1 4  0 4  1 4 
     -2.5 2.5  2.5 2.5  / 
\multiput{3} at -1.5 1.5  1.5 1.5  -2 2  2 2 -.5 3.5  .5 3.5 
     0 1  -2 3  2 3 /
\multiput{4} at -.5 1.5  .5 1.5  -1.5 2.5  1.5 2.5  -1 3  1 3 / 
\multiput{6} at -1 2  0 2  1 2  -.5 2.5  .5 2.5  0 3 / 
\put{$n = 7$} at 0 -.7
\endpicture} at 17 0
\endpicture}
$$

\subsection{Indecomposables $X$ with $bX \ge 3$.}
\label{sec-twelve-three}

{\it Let $n = 6$. If $b \equiv 0\ \mod 3,$
  then  there is a BTh-family $X$ in $\Cal S(n)$ with $bX = b.$
  If $b \not \equiv 0\ \mod 3,$ then the indecomposables $X$ with $bX = b$ are combinatorial.}

\smallskip
\begin{proof}
  For the first assertation, let $M$ be the standard family and $X = M[\ell],$ where $\ell = b/3.$
  For the second assertion, look at the root system.
\end{proof}

\medskip
{\it Let $n \ge 7.$ If $b \ge 3,$ there is a BTh-family $X$ in $\Cal S(n)$ with $bX = b.$.}

\smallskip
\begin{proof}
  If $b \equiv 0 \ \mod 3,$ we have such a family already in $\Cal S(6).$ 
  \smallskip
  
  Thus, we have to consider the cases $b \equiv 1 \ \mod 3$ and $b \equiv 2 \ \mod 3.$
  First, let $b \equiv 1 \ \mod 3.$ 
  If $b = 4,$ take $D$ with $V = [7,5,3,1]$ and $U = [4,2]$, so that $G_5H_4D = M.$ 
  For $b = 4+ 3\ell$, take the interpolation: one copy of $D$ and $\ell$ copies of $M$.
  Finally, let $b \equiv 2 \ \mod 3.$ 
  For $b = 8+ 3\ell$, take the interpolation: two copy of $D$, $\ell$ copies of $M$.
  It remains $b = 5.$ Here we take the family $\Cal S_3(7)$.
\end{proof}

\bigskip
\section{The indecomposable objects $X = (U,V)$ with $bU=1$.}
\label{sec-thirteen}

In this section, we consider invariant subspaces $(U,V)$ with $bU =1$,
thus $U$ considered as an object in $\Cal N$ is indecomposable and therefore cyclic,
but also non-zero.  (Note that the case where $U=0$ is covered by Theorem~\ref{theoremone}
as we are dealing with a boundary picket of type $(0,[\ell])$ for some $\ell\in\mathbb N$.)

The authors are grateful to William Crawley-Boevey for pointing out the
bijection between the set of isomorphism
classes of the indecomposable objects $X = (U, V )$ in
$\Cal S(n)$ with $bU = 1$ and the set of partitions with
perimeter between $1$ and $n$, see Corollary~\ref{cor-thirteen-one-two}.

\subsection{The first combinatorial descriptions.}
\label{sec-thirteen-one}

Let us attach to
any non-empty subset $E$ of $\{1,\dots,n\}$ an indecomposable module $M(E) =
(U,V)$ in $\Cal S(n)$ with $bU = 1.$ 
Also, we will attach to $E$ a partition $\lambda(E)$. In this way, we will obtain
for $n\ge 1$
bijections between first, the isomorphism classes of indecomposable objects $X = (U,V)$
in $\Cal S(n)$ with $bU = 1,$  the non-empty subsets $E$ of $\{1,2,\dots,n\}$,
and third, the partitions with perimeter at most $n$. 
\medskip

We start with a non-empty set $E =\{e_1,\dots,e_m\}$  of positive numbers, of
cardinality $m$. 
If $m$ is even, let $b = m/2,$ otherwise $b = (m+1)/2,$ thus $b = \lceil m/2 \rceil.$ 
For $1\le i \le b,$ let $d_i = e_{i+m-b}-e_b$. 
Thus, there are given the two sequences $1\le e_1 < \cdots < e_b$ and $0 \le d_1 < 
\cdots < d_b$ (and we have $d_1 = 0$ if and only if $m$ is odd).

Let $\widetilde V(i)$ be the indecomposable $\widetilde
\Lambda$-module of length $e_i+d_i,$ generated in degree $e_b+d_i$, say
with generator $x_i$, and let $\widetilde V = \bigoplus_i \widetilde V(i).$
Let $y_i = T^{d_i}x_i$, this is an element of degree $e_b$, for all 
$1\le i \le b.$
Let $\widetilde U$ be generated by the element 
$y = \sum y_i$ (which, as we know, has degree $e_b$) and $\widetilde W =
\widetilde V/\widetilde U.$ In this way, we have defined
$\widetilde M(E) = (\widetilde U,\widetilde V)$. Finally, let
$M(E) = (U,V) = \pi\widetilde M(E)$, where $\pi$ is the covering functor,
and $W = V/U.$ 

\medskip
The construction shows (for $E$ having cardinality $m$):
\begin{itemize}[leftmargin=3em]
  \item[$\bullet$]
    {\it The height of $M(E)$ is the maximum of the elements in $E$.}
  \item[$\bullet$]
    {\it The width of $M(E)$ is  $\lceil m/2 \rceil$.}
  \item[$\bullet$]
    {\it If $m$ is even, then $U$ is contained in the radical of $V$,
      therefore $b W = b V.$ }
  \item[$\bullet$]
    {\it If $m$ is odd, then $U$ is not contained in the radical of $V$,
      therefore $bW = bV -1.$}
\end{itemize}

	\bigskip
Here are two examples: We start with $E$ equal to $\{2,3,5,6,8,9\}$ or $\{2,3,5,8,9\}$.
In both cases, $b = 3$ and $\{e_1,e_2,e_3\} = \{2,3,5\};$ in particular, we have 
$e_b = 5.$ 
In the first case, $\{d_1,d_2,d_3\} = \{1,3,4\}$ (and $m$ is even), in the second case,  
$\{d_1,d_2,d_3\} = \{0,3,4\}$ (and $m$ is odd).  
The third column shows $\widetilde M(E)$.
In the last column we visualize already here
the partitions $\lambda(E)$; but $\lambda(E)$ 
will be defined only at the end of this section.

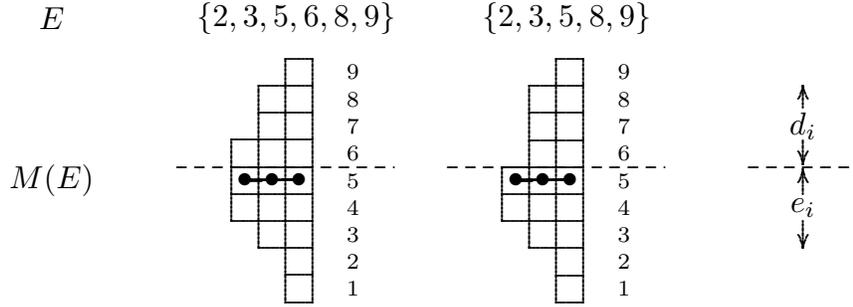
$$  
{\beginpicture
    \setcoordinatesystem units <.3cm,.3cm>
\put{\beginpicture
   \put{$E$} at 0 0 
\endpicture} at -18 7
\put{\beginpicture
   \put{$(d_1,\dots,d_b)$} at 0 0.8 
   \put{$(e_1,\dots,e_b)$} at 0 -.8 
\endpicture} at -9 7

\put{\beginpicture
   \put{$\widetilde M(E)$} at 0 0 
\endpicture} at 0 7

\put{\beginpicture
   \put{$\lambda(E)$} at 0 0 
\endpicture} at 15 7

\put{\beginpicture
   \put{$(1,3,4)$} at 0 1 
   \put{$(2,3,5)$} at 0 -1 
\endpicture} at -9 0

\put{\beginpicture
   \put{$(0,3,4)$} at 0 1 
   \put{$(2,3,5)$} at 0 -1 
\endpicture} at -9 -10
\put{\beginpicture
\put{$\{2,3,5,6,8,9\}$} at 0 0 
\endpicture} at -18 0

\put{\beginpicture
\put{$\{2,3,5,8,9\}$} at 0 0 
\endpicture} at  -18 -10

\put{\beginpicture
\multiput{} at 0 0  3 9 /
\plot 2 0  3 0  3 9  2 9  2 0 /
\plot 3 2  1 2  1 8  3 8 /
\plot 3 3  0 3  0 6  3 6 /
\plot 2 1  3 1 /
\plot 0 4  3 4 /
\plot 0 5  3 5  /
\plot 1 7  3 7 /
\multiput{$\bullet$} at .5 4.5  1.5 4.5  2.5 4.5 /
\plot 0.5 4.5  2.5 4.5 / 
\plot 0.5 4.45  2.5 4.55 / 
\plot 0.5 4.55  2.5 4.55 / 
\setdashes <1mm>
\plot -2 5  6 5 /
\put{$\ss 1$} at  4.5 0.5 
\put{$\ss 2$} at  4.5 1.5 
\put{$\ss 3$} at  4.5 2.5 
\put{$\ss 4$} at  4.5 3.5 
\put{$\ss 5$} at  4.5 4.5 
\put{$\ss 6$} at  4.5 5.5
\put{$\ss 7$} at  4.5 6.5 
\put{$\ss 8$} at  4.5 7.5 
\put{$\ss 9$} at  4.5 8.5 
\endpicture} at 0 0
\put{\beginpicture
\multiput{} at 0 0  3 9 /
\plot 2 0  3 0  3 9  2 9  2 0 /
\plot 3 2  1 2  1 8  3 8 /
\plot 3 3  0 3  0 5  1 5 /
\plot 1 6  3 6 /
\plot 2 1  3 1 /
\plot 0 4  3 4 /
\plot 0 5  3 5  /
\plot 1 7  3 7 /
\multiput{$\bullet$} at .5 4.5  1.5 4.5  2.5 4.5 /
\plot 0.5 4.5  2.5 4.5 / 
\plot 0.5 4.45  2.5 4.55 / 
\plot 0.5 4.55  2.5 4.55 / 
\setdashes <1mm>
\plot -2 5  6 5 /
\put{$\ss 1$} at  4.5 0.5 
\put{$\ss 2$} at  4.5 1.5 
\put{$\ss 3$} at  4.5 2.5 
\put{$\ss 4$} at  4.5 3.5 
\put{$\ss 5$} at  4.5 4.5 
\put{$\ss 6$} at  4.5 5.5
\put{$\ss 7$} at  4.5 6.5 
\put{$\ss 8$} at  4.5 7.5 
\put{$\ss 9$} at  4.5 8.5 
\endpicture} at 0 -11

\put{\beginpicture
\multiput{} at 0 0  3 9 /
\setdashes <1mm>
\plot -2 5  2 5 /
\setsolid
\arr{0 6}{0 5.1}
\put{$d_i$} at 0 6.5
\arr{0 7.2}{0 8}
\arr{0 4}{0 4.9}
\put{$e_i$} at 0 3.5
\arr{0 3}{0 2}
\endpicture} at 7 0
\put{\beginpicture
\multiput{} at 0 0  3 9 /
\setdashes <1mm>
\plot -2 5  2 5 /
\setsolid
\arr{0 6}{0 5.1}
\put{$d_i$} at 0 6.5
\arr{0 7.2}{0 8}
\arr{0 4}{0 4.9}
\put{$e_i$} at 0 3.5
\arr{0 3}{0 2}
\endpicture} at 7 -11

\put{\beginpicture
\multiput{} at 0 0  5 5 /
\plot 0 0  0 5  1 5  1 0  0 0  /
\plot 1 1  2 1  2 4  1 4 /
\plot 2 1  3 1  3 3  2 3 /
\plot 1 5  5 5  5 4  2 4 /
\plot 3 3  5 3  5 4 /
\plot 3 2  4 2  4 3 /
\setdots <.5mm>
\plot 1 5  1 0  0 0  0 5  5 5  5 4  0 4 /
\plot 0 3  5 3  5 4 /
\plot 0 1  3 1  3 5 /
\plot 0 2  3 2 /
\plot 2 1  2 5 /
\plot 4 3  4 5 /
\plot 3 2  4 2  4 3 /
\put{$(5,5,4,3,1)$} [l] at -.5 -1
\endpicture} at 15 0
\put{\beginpicture
\multiput{} at 0 0  5 5 /
\multiput{} at 0 0  5 5 /
\plot 0 0  0 5  1 5  1 0  0 0  /
\plot 1 1  2 1  2 4  1 4 /
\plot 2 1  3 1  3 3  2 3 /
\plot 1 5  5 5  5 4  2 4 /
\plot 3 3  5 3  5 4 /
\setdots <.5mm>
\plot 1 5  1 0  0 0  0 5  5 5  5 4  0 4 /
\plot 0 3  5 3  5 4 /
\plot 0 1  3 1  3 5 /
\plot 0 2  3 2 /
\plot 2 1  2 5 /
\plot 4 3  4 5 /
\put{$(5,5,3,3,1)$} [l] at -.5 -1
\endpicture} at 15 -11

\endpicture}
$$
The pictures in the third column (with label $\widetilde M(E)$)
show (as usual) a direct decomposition of the global space as a direct sum
of $\widetilde \Lambda$-modules $\widetilde V(i)$, with $1\le i \le b,$ as well as the generator $y$ of 
$\widetilde U.$ 
We have inserted a horizontal dashed line: For $1\le i \le b,$ the number $e_i$ is the length of 
the image of the map $\widetilde U \to \widetilde V(i)$, thus the number of boxes below the
dashed line. the number $d_i$ is the length of 
the cokernel of the map $\widetilde U \to \widetilde V(i)$, thus the number of boxes above the dashed line. 

In the presentation of the partition $\lambda(E)$ in the last column, 
one finds both the Young diagram, as well
as the sequence of numbers $(\lambda(E)_1,\dots, \lambda(E)_b).$
In the Young diagrams, we have indicated by solid lines 
sequences of boxes of length $e_*$ (vertically), and of length $d_*$ (horizontally).
Note that the indices of $e_*$ and $d_*$ occur in reverse order --- but still we 
recover the numbers $e_i$ and $d_i$, thus also $E$, from $\lambda(E).$

\medskip
\begin{proposition}
  \label{prop-thirteen-one}
  \begin{itemize}[leftmargin=3em]
  \item[\rm(a)] Any indecomposable object $X = (U,V)$ in $\Cal S(n)$ with $bU = 1$
    is gradable.
  \item[\rm(b)] The construction $E \mapsto M(E)$ provides a bijection between
    the non-empty subsets of $\{1,\dots,n\}$ and the indecomposable
    objects
    $X = (U,V)$ in $\Cal S(n)$ with $bU = 1$.
  \end{itemize}
\end{proposition}

\begin{proof}
  (1) First, let us show that the modules 
$\widetilde M(E)$ are indecomposable and pairwise non-isomorphic.
For the indecomposability, 
we show inductively that $\widetilde M(E)$ has endomorphism ring $k$: 
If $E$ has cardinality at most 2, $\widetilde M(E)$ is a
picket. If $E$ has cardinality at least 3, then 
$\widetilde M(E)$ is an extension of an orthogonal pair, where one of the modules is a  picket, the other one a module of the form 
$\widetilde M(E')$ such that the cardinality of $E'$ is smaller than the
cardinality of $E$.

In order to see that the objects $M(E) = \widetilde M(E)$ are pairwise non-isomorphic, 
we show that we can recover $E$ from $\widetilde  M(E) = (U,V)$. The global space $V$ 
is the direct sum of indecomposable $\Lambda$-modules
$V(i)$ 
of length $e_i+d_i$ for $1 \le i \le b$, and the inclusion map  $\mu:U \to V$ 
yields maps from $U$ to
$V(i)$ with image of length $e_i.$ 
	\smallskip 
(2) We show that any indecomposable object $X = (U,V)$ with $bU = 1$ is of the form 
$M(E)$ for some subset $E$ of $\{1,\dots,n\}$.
Let $\widetilde X = (\widetilde U,\widetilde V)$. We decompose 
$\widetilde V = \bigoplus \widetilde V(i)$ with $\widetilde V(i)$ indecomposable.
The inclusion map $\mu:\widetilde U \to \widetilde V$ gives maps 
$\mu_i:\widetilde U \to \widetilde V(i)$. 
Let $y$ be a generator of $\widetilde U$
and $x_i$ a generator of $\widetilde V(i)$. 
We can assume that $\mu(y) = T^{d_i}x_i$ with natural
numbers $d_i \ge 0$ such that $d_1 \le d_2 \le \cdots \le d_b.$

Let $e_i$ be the length of the image of $\mu_i.$ We claim that we have $e_i < e_{i+1}$.
Assume for the contrary that $e_t \ge e_{t+1}$ for some $1\le t < b.$
Then 
$\widetilde V = \widetilde V'\oplus \widetilde V(t+1)$, 
where $\widetilde V'$ is generated by the elements $x_i$ with $i\neq t,t+1$ and
the element $x_t+T^fx_{t+1}$, where $f = d_{t+1}-d_t,$ 
and we have $\widetilde U \subseteq \widetilde V'$. Thus $\widetilde X$ is
decomposable, a contradiction. 

Second, we claim that $d_i < d_{i+1}$ for all $1 \le i <b.$
Assume for the contrary that $d_t = d_{t+1}$ for some $1\le t < b.$
Then $\widetilde V = \widetilde V'\oplus \widetilde V(t)$, 
where $\widetilde V'$ is generated by the elements $x_i$ with $i\neq t,t+1$ and
the element $x_t+x_{t+1}$, and we have $\widetilde U \subseteq 
\widetilde V'$. Thus $\widetilde X$ is
decomposable, a contradiction. 

Let $m = 2b-1,$ if $d_1 = 0,$ and $m = 2b$ otherwise. For $1\le i \le b$, let
$e_{i+m-b} = d_i+e_b,$ and let $E = \{e_1,\dots,e_m\}$.
Then  $X$ is isomorphic to $M(E)$.
	\smallskip
(a) and (b) are direct consequences of (1) and (2).  
\end{proof}

\medskip
\begin{corollary}
  \label{cor-thirteen-one}
  Let $X = (U,V)$ be indecomposable in $\Cal S(n)$ with $bU = 1$
  Then we have $bX \le uX$ and $bX \le \frac{n+1}2.$
\end{corollary}

\begin{proof}
  According to Proposition, $X = M(E)$ for a non-empty subset 
  $E = \{e_1 < e_2 < \cdots < e_m\}$ of $\{1,\dots,n\};$ of course, we have  
  $m\le n$. The construction of $\widetilde M(E)$ shows that $b = bX$ is equal
  to $\frac m2$ or to $\frac {m+1}2$. Thus $b \le \frac{m+1}2 \le \frac{n+1}2.$
  Also, $uX$ is the length of $U$, thus equal to $e_{b}$. Since $b \le e_b$, we have 
  $b \le e_b = uX.$  (Of course, we have seen this inequality already in Theorem~\ref{theoremone}.)
\end{proof}

	\bigskip
{\bf Partitions.} Now, let us turn the attention to partitions. If $\lambda = 
(\lambda_1,\dots,\lambda_u)$ is a partition with $u\ge 1$ parts, its {\it perimeter}
is defined to be $\lambda_1+u-1$, see \cite{Str}. 

In order to visualize partitions, we have used throughout the paper diagrams with boxes,
but we have deviated from the usual arrangement of these boxes.
In this section, it will be convenient to 
work also with the usual arrangement, namely with the {\it Young diagram} of the
partition $\lambda = (\lambda_1,\dots,\lambda_u)$. It consists of the set of pairs 
$(i,j)$ of positive integers with $1 \le j \le \lambda_i$. The Young diagram of
$\lambda$ will usually by displayed by drawing $u$ rows of unit boxes, 
starting with a row with $\lambda_1$ boxes, below of it a second row consisting
of $\lambda_2$ boxes, and so on, with all rows being aligned on the left.
Here is the display of the partition $(5,4,4,1)$:
$$  
{\beginpicture
    \setcoordinatesystem units <.3cm,.3cm>
\multiput{} at 0 0  5 4 /
\plot 1 4  1 0  0 0  0 4  5 4  5 3  0 3 /
\plot 0 1  4 1  4 4 /
\plot 0 2  4 2 /
\plot 2 1  2 4 /
\plot 3 1  3 4 /
\endpicture}
$$
The boxes with labels of the form $(i,i)$ are called the {\it diagonal} boxes.
The (Durfee) {\it rank} of the partition
$\lambda$ is the number of diagonal boxes (the displayed
partition has rank 3). We denote by $\lambda'$ the partition conjugate to $\lambda$,
it is defined as follows: $\lambda'_j$ is the number of parts $\lambda_i$ 
with $\lambda_i \ge j$ (thus, for the displayed partition, the conjugate partition is
$(4,3,3,3,1)$). Finally, for every box $(i,j)$ there is defined a corresponding {\it hook}:
It consists of the boxes $(i,j')$ with $j' \ge j$ and the boxes $(i',j)$ with $i' \ge i.$
The hooks for the diagonal boxes will be called {\it diagonal hooks.}
Here are the three diagonal hooks of the partition $(5,4,4,1)$:
$$  
{\beginpicture
    \setcoordinatesystem units <.3cm,.3cm>
\multiput{} at 0 0  5 4 /
\plot 1 3  1 0  0 0  0 4  5 4  5 3  1 3 /
\plot 1 1  2 1  2 2  4 2  4 3  1 3  1 1 /
\plot 2 1  4 1  4 2 /
\setshadegrid span <.3mm>
\vshade 0 3 4 <z,z,,> 1 3 4 /
\vshade 1 2 3 <z,z,,> 2 2 3 /
\vshade 2 1 2 <z,z,,> 3 1 2 /

\endpicture}
$$
(and the diagonal boxes are shaded).

	\medskip
Starting with a non-empty subset $E$ of $\{1,2,\dots,n\}$ with $b = \lceil m/2 \rceil,$ 
we want to attach a partition $\lambda(E)$ with rank $b$. We will do this by specifying
the first numbers $\lambda_1,\dots,\lambda_b$ of the partition, as well as the 
first numbers $\lambda'_1,\dots,\lambda'_b$ of the conjugate partition:
Let $\lambda_i = d_{b-i+1}+i$, and
let $\lambda'_i = e_{b-i+1}+i-1$, for $1 \le i \le b$.
Note that our construction of $\lambda(E)$ specifies the diagonal hooks;
in particular, these hooks have length $e_b+d_b, \dots, e_1+d_1$ (in this order).  
It may be convenient to look at the examples $\lambda(\{2,3,5,6,8,9\})$ and 
$\lambda(\{2,3,5,8,9\})$ presented near the beginning of the section!

\medskip
\begin{lemma}
  \label{lem-thirteen-one-two}
  Let $n\ge 1.$ The map $\lambda$ provides a bijection between the 
  non-empty subsets of $\{1,2,\dots,n\}$ 
  and the partitions with perimeter between $1$ and $n$.
\end{lemma}

\begin{proof}\
  It is clear that for $E \subseteq \{1,2,\dots,n\},$ the
  partition $\lambda(E)$ has perimeter at most $n$, since $e_b+d_b = e_m \le n$ 
  is the largest length of a diagonal hook, thus the perimeter of $\lambda(E)$.  
  
  Conversely, given a partition $\lambda$ with perimeter between $1$ and $n$, the 
  diagonal hooks can be used to define 
  sequences $(e_1,\dots,e_b)$ and $(d_1,\dots,d_b)$ of numbers and in this way to
  define a set $E = \{e_1,\dots,e_b,d_1+e_b,\dots,d_b+e_b\}$ such that $E$ 
  has cardinality $2b-1$
  (in case $d_1 = 0$) or $2b$ (in case $d_1 > 0$), and such that $\lambda(E) = \lambda$.
\end{proof}

\medskip
\begin{corollary}
  \label{cor-thirteen-one-two}
  Let $n\ge 1$. The correspondence 
  $M(E) \leftrightarrow \lambda(E)$, where $E$ is a non-empty subset of $\{1,2,\dots,n\},$
  provides a bijection between the set of isomorphism classes of the
  indecomposable objects $X = (U,V)$ in $\Cal S(n)$ with $bU = 1$  and the set of
  partitions with perimeter between $1$ and $n.$
\end{corollary}

\medskip
Let $M(E) = (U,V)$. One has to be aware that both $V$ and $\lambda(E)$ are partitions of 
the same number $v$: This number $v$ 
is the length of $V$.
The correspondence $M(E) \leftrightarrow \lambda(E)$ is easily visualized,
when we rotate the partition $\lambda(E)$ by $45^\circ$. We have to fold (or to kink) the
parts of the
partition $V$ and obtain in this way the central hooks of $\lambda(E)$.  
As an example, let $E = \{2,3,5,6,8,9\}$ as considered at the beginning of this
section.
$$  
{\beginpicture
    \setcoordinatesystem units <.3cm,.3cm>
\put{\beginpicture
   \put{$M(E)$} at 0 0 
\endpicture} at 0 7
\put{\beginpicture
   \put{$\lambda(E)$} at 0 0 
\endpicture} at 9 7

\put{\beginpicture
\multiput{} at 0 0  3 9 /
\plot 0 0  1 0  1 9  0 9  0 0 /
\plot 0 1  1 1 /
\plot 0 2  2 2  2 8  0 8 /
\plot 0 7  2 7 /
\plot 0 3  3 3  3 6  0 6 /
\plot 0 4  3 4 /
\plot 0 5  3 5 /
\multiput{$\bullet$} at .5 4.5  1.5 4.5  2.5 4.5 /
\plot 0.5 4.5  2.5 4.5 / 
\plot 0.5 4.45  2.5 4.55 / 
\plot 0.5 4.55  2.5 4.55 / 
\endpicture} at 0 0

\put{\beginpicture
    \setcoordinatesystem units <.212cm,.212cm>
\multiput{} at 0 0  3 10 /
\plot 0 5  5 0  6 1  1 6 /
\plot 0 5  5 10  7 8  2 3 /
\plot 1 4  6 9 /
\plot 3 2  7 6  4 9 /
\plot 4 1  7 4  3 8 /
\plot 2 7  6 3 /
\setshadegrid span <.3mm>
\vshade 0 5 5 <z,z,,> 1 4 6 <z,z,,> 2 5 5  <z,z,,> 
                      3 4 6 <z,z,,> 4 5 5  <z,z,,> 
                      5 4 6 <z,z,,> 6 5 5 /
\endpicture} at 8 0
\endpicture}
$$
In the partition $\lambda(E)$ on the right, we have shaded the diagonal boxes
in order to stress that they correspond to the boxes of $V$ which contain a bullet. 
The parts of the partition of $V$ seen on the left are
folded (or kinked) at the bullet level 
in order to obtain the central hooks of $\lambda(E)$.

\subsection{Counting isomorphism classes.}
\label{sec-thirteen-two}

We are going to count the number of isomorphism classes of invariant subspaces
$X = (U,V)$ with $bU = 1$  and suitable additional properties.

\bigskip\noindent 
{\bf (1) Height.} 
{\it 
The number of indecomposable objects $X = (U,V)$ with $bU = 1$, and
with $X$ of height $n \ge 1$ is $2^{n-1}$} 
(note that there is no indecomposable object of height $0$).
{\it For $n \ge 0,$  the number of indecomposable objects $X = (U,V)$ in $\Cal S(n)$
with $bU = 1$ is $2^n-1$.}

\begin{proof}
  We count the objects $M(E)$, where $E$ is a subset of $\{1,2,\dots,n\}$
  such that $n$ belongs to $E$. Thus, $E$ is an arbitrary subset of $\{1,2,\dots,n-1\}$.
  For the second assertion, we have to add the numbers $2^m$ with $1\le m \le n-1.$
\end{proof}

	\bigskip\noindent
{\bf (2)} {\bf Height and width.} 
{\it Let $n \ge 1.$
The number of indecomposable objects $X = (U,V)$ with $bU= 1$, 
with $U$ contained in the radical of $V$, and  
with $X$ of height $n$ and width $b$ is $\binom {n-1}{2b-1}$.

The number of indecomposable objects $X = (U,V)$ with $bU = 1$,
with $U$ not contained in the radical of $V$, and  
with $X$ of height $n$ and width $b$ is $\binom {n-1}{2b-2}$.

The number of indecomposable objects $X = (U,V)$ with $bU = 1$, 
and  with $X$ of height $n$ and width $b$ is $\binom {n}{2b-1}$.}

\begin{proof}
  First, we deal with the subspaces $U$ contained in the radical of $V$.
  The objects $(U,V)$ correspond bijectively to the subsets $E$ of $\{1,2,\dots,n\}$
  of cardinality $2b$ (since $U$ is contained in the radical of $V$) and
  containing the element $n$ (since the height is equal to $n$), thus $E\setminus\{n\}$
  is a subset of $\{1,2,\dots,n-1\}$ of cardinality $2b-1$. The number of such
  subsets is $\binom{n-1}{2b-1}$.
  
  \smallskip
  Second, we coonsider the subspaces $U$ not contained in the radical of $V$.
  The objects $(U,V)$ correspond bijectively to the subsets $E$ of $\{1,2,\dots,n\}$
  of cardinality $2b-1$ (since $U$ is not contained in the radical of $V$) and
  containing the element $n$ (since the height is equal to $n$), thus $E\setminus\{n\}$
  is a subset of $\{1,2,\dots,n-1\}$ of cardinality $2b-2$. The number of such
  subsets is $\binom{n-1}{2b-2}$.

  \smallskip
According to the first two assertions, the number 
of indecomposable objects $X = (U,V)$ with $bU = 1$, 
and  with $X$ of height $n$ and width $b$ 
is $\binom{n-1}{2b-1}+ \binom{n-1}{2b-2}
= \binom{n}{2b-1}$.
\end{proof}

\bigskip\noindent
{\bf (3) Length of the global space.}
{\it Let $v \ge 1.$ 
The number of indecomposable objects $X = (U,V)$ with $bU = 1$, and $v = |V|$
is equal to the number of partitions of $v$.}

\begin{proof}
  This follows directly from the bijection 
  $M(E) \leftrightarrow \lambda(E)$ exhibited at the end of Section~\ref{sec-thirteen-one}
\end{proof}

  \bigskip
  Before we proceed, let us recall 
  a well-known counting 
procedure dealing with paths in a rectangular grid. The grid paths
we consider use unit steps going from left to right, as well as unit steps going up.

\medskip
\begin{lemma}
  \label{lem-thirteen-two}
  \begin{itemize}[leftmargin=3em]
    \item[\rm (a)]
        Given the rectangle
        $$  
        {\beginpicture
          \setcoordinatesystem units <.3cm,.3cm>
          \multiput{} at 0 0  7 5 /
          \plot 0 0  0 5   7 5  7 0  0 0 /
          \put{$A$} at -.8 -.5 
          \put{$C$} at 7.8 5.5 
          \put{$a$\strut} at -.7 2.5
          \put{$c$\strut} at 3.5  5.7
          \setdots <1mm>
          \plot 1 0  1 5 /
          \plot 2 0  2 5 /
          \plot 0 1  7 1 /
          \plot 0 4  7 4 /
          \multiput{$\bullet$} at 0 0  7 5 /
          \endpicture}
        $$
        the number of paths from $A$ to $C$ is $\binom{a+c}a = \binom{a+c}c$.
      \item[\rm (b)]
        Given a grid vertex $B$ on the dashed diagonal
        $$  
        {\beginpicture
          \setcoordinatesystem units <.3cm,.3cm>
          \multiput{} at 0 0  7 5 /
          \plot 0 0  0 5   7 5  7 0  0 0 /
          \put{$A$} at -.8 -.5 
          \put{$C$} at 7.8 5.5 
          \setdots <1mm>
          \plot 1 0  1 5 /
          \plot 2 0  2 5 /
          \plot 0 1  7 1 /
          \plot 0 4  7 4 /
          \multiput{$\bullet$} at 0 0  7 5  2 3 /
          \setdashes <1mm>
          \plot 0 5  5 0 /
          \setsolid
          \plot 0 3  2 3  2 5 /
          \put{$B$} at 2.8 3.2
          \multiput{$b$\strut} at -.7 4  1  5.7 /
          \endpicture}
        $$
        the number of paths from $A$ to $C$ via $B$ is the product $\binom a b\binom c b$.
  \end{itemize}
\end{lemma}

\begin{proof}
  (a) Any path from $A$ to $C$ 
has length $a+c$, and precisely $a$ of the steps are going upwards.
(b) It follows from (a) that the number of paths from $A$ to $B$ is $\binom a b$,
and that the number of paths from $B$ to $C$ is $\binom c b$. Since all combinations
are allowed, the total number of paths from $A$ to $C$ via $B$ is the product.
\end{proof}

	\medskip
Since any path from $A$ to $C$ passes through precisely one vertex $B$ on the dashed
diagonal, we obtain the following special case of the Vandermonde formula: {\it For all
natural numbers $a,c,$ we have}
$$
 \binom{a+c}a = \sum_{b\ge 0} \binom ab \binom cb.
$$
	\bigskip
        Now we return to the category $\Cal S(n)$. 
        
	\bigskip\noindent
{\bf (4) Length of the subspace, height.} 
{\it The number of indecomposable objects $X = (U,V)$ with $bU = 1$ and $U$ of length $u$,
and with $X$ of height $n \ge 1$ is $\binom {n-1}{u-1}$.}

{\it The number of indecomposable objects $X = (U,V)$ in $\Cal S(n)$
with $bU= 1$  and $U$ of length $u$
is $\binom {n}{u}$.}

\smallskip
\begin{proof}
  We look at the partitions with $u$ parts and perimeter $n$. Thus, we deal 
  with the grid paths from $A$ to $C$ 
  in the rectangle with sides of length $u-1$ and $n-u$.
  $$  
  {\beginpicture
    \setcoordinatesystem units <.3cm,.3cm>
    \multiput{} at -1 0  7 6 /
    \plot -1 0  -1 6  7 6  7 5  0 5  0 0  -1 0 /
    \setdashes <1mm> 
    \plot  7 5  7 0  0 0 /
    \put{$u$\strut} at -1.7 3
    \put{$n-u+1$\strut} at 3  6.7
    \setdots <1mm>
    \plot 1 0  1 5 /
    \plot 2 0  2 5 /
    \plot 0 1  7 1 /
    \plot 0 4  7 4 /
    \put{$A$} at -.2 -.9 
    \put{$C$} at 8 5.2 
    
    \multiput{$\bullet$} at 0 0  7 5 /
    \endpicture}
  $$
  According to Lemma~\ref{lem-thirteen-two}(a), the number of paths is $\binom {n-1}{u-1}$.
  
  The second assertion, follows by induction.
\end{proof}

  \bigskip\noindent
{\bf (5) Length of the subspace, height and width.}
{\it Let $n \ge 1.$ 
The number of indecomposable objects $X = (U,V)$ with $bU = 1$ and $U$ of length $u$,
with $X$ of height $n$ and width $b$ is $\binom {u-1}{b-1}\binom {n-u}{b-1}$.}

\smallskip
\begin{proof}
  As in (4), we look at the partitions with $u$ parts and perimeter $n$. 
  but we look only at the paths which pass through $B$ (the lower right corner of the box
  $(b,b)$).
  $$  
  {\beginpicture
    \setcoordinatesystem units <.3cm,.3cm>
    \multiput{} at -1 0  7 6 /
    \plot -1 0  -1 6  7 6  7 5  0 5  0 0  -1 0 /
    \plot -1 3  2 3  2 6 /
    \setdashes <1mm> 
    \plot  7 5  7 0  0 0 /
    \put{$u$\strut} at -1.7 3
    \put{$n-u+1$\strut} at 3  6.7
    \setdots <1mm>
    \plot 1 0  1 5 /
    \plot 2 0  2 5 /
    \plot 0 1  7 1 /
    \plot 0 4  7 4 /
    \multiput{$\bullet$} at 0 0  2 3  7 5 /
    \put{$B$} at 2.7 3
    \put{$A$} at -.2 -.9 
    \put{$C$} at 8 5.2 
    \endpicture}
  $$
  According to Lemma~\ref{lem-thirteen-two}(b), the number of paths from $A$ to $C$ via $B$ 
  is $\binom {u-1}{b-1}\binom {n-u}{b-1}$ (the rectangle between $A$ and $B$ has sides 
  $u-b$ and $b-1$; the rectangle between $B$ and $C$ has sides $b-1$ and $n-u+1-b$). 
\end{proof}

\subsection{Insertion:
  The $T$-height sequence.} 
\label{sec-thirteen-three}

We have started to look at 
the objects $X = (U,V)$ in $\Cal S(n)$ with $bU = 1$. This is the proper context for mentioning
the concept of the $T$-height and the
$T$-height sequence of a non-zero element $y\in V$, where $V$
is an object in $\Cal N(n)$, for some $n$, as considered by
Pr\"ufer \cite{P} already in 1923. (Warning: 
The reader should be aware that height (as defined in Section~\ref{sec-one-one}) and $T$-height are completely
different concepts: Whereas the height of $\Lambda y$ is just the length of
$\Lambda y$ and does not depend on the global space $V$, the $T$-height 
strongly depends on the embedding of $y$ in $V$.)

\medskip
\begin{definitions}
  Let $V$ be an object in $\Cal N(n)$, and $y$ a non-zero
  element in $V$: The {\it $T$-height $h(y) = h_V(y)$ 
    of $y$ in $V$} is the largest number $d$ such that $y\in T^dV.$
  The {\it $T$-height sequence $H(y) = H_V(y)$ 
    of $y$ in $V$} is the sequence $h_V(y),h_V(Ty),\dots,h_V(T^{e-1}y)$,
  where $T^{e-1}y \neq 0$ and $T^ey = 0.$
  Obviously, {\it the $T$-height sequence $H(y)$ of a non-zero
    element $y\in V$ is a strictly increasing sequence of numbers in
    $\{0,\dots,n-1\}.$}
  Note that the $T$-height sequence does not depend on the choice of the generator $y$
  of $\Lambda y$.  Hence, if $U=\Lambda y$ is a non-zero cyclic submodule of $V$
  we may write $H(U)=H(y)$.
  We should mention that the sequence $H(x)$ and its relevance 
  were discussed by Kaplansky \cite{Kap}, see p.57 ff (under the name Ulm sequence). 
\end{definitions}

\medskip
In order to be able to phrase these definitions
in terms of the category $\Cal S(n)$, we need the following observation. 

\medskip
We recall from \cite{ARS} that a homomorphism $f\:V \to V'$ in a length category 
is said to be {\it left minimal} provided any
endomorphism $g\:V' \to V'$ with $gf = f$ is an automorphism. And given any morphism 
$f\:U \to V$, there is a direct decomposition $V = V'\oplus V''$ such
that the image of $f$ is contained in $V'$ and $p_{V'}f\:U \to V'$ is left minimal,
where $p_{V'}$ is the projection of $V$ onto $V'$; the map $p_{V'}f\:U \to V'$ is
called a {left minimal version of} $f$. If 
$p_{V'}f\:U \to V'$  and $p_{V_1'}f\:U \to V_1'$ are left minimal versions of $f\:U \to V$,
then there is an isomorphism $h\:V' \to V_1'$ with $p_{V'}f = hp_{V'}$. 

\medskip
\begin{lemma}
  \label{lem-thirteen-three}
  Let $U$ be indecomposable in $\Cal N(n)$ and $X = (U,V)$ an object in
  $\Cal S(n)$. Then $(U,V)$ is indecomposable if and only if the inclusion map
  $U \to X$ is left minimal.
\end{lemma}

\begin{proof}
  First, assume that $(U,V)$ is indecomposable. Let $(U,V')$ be a left miniaml version
  of the inclusion map $U \to V$, say with $V = V'\oplus V''$.
  Then we get a direct decomposition
  $(U,V) = (U,V')\oplus (0,V'')$ in $\Cal S(n)$, thus the indecomposability of $(U,V)$ asserts
  that $V'' = 0,$  thus $V' = V$. This shows that the inclusion map $U \to V$ is left minimal.
  Conversely, assume that the inclusion map $U \to V$ is left minimal. Given any
  direct decomposition $(U,V) = (U',V')\oplus (U'',V''),$ we have $U = U'\oplus U''.$ Since
  we assume that $U$ is indecomposable, we see that one of $U', U''$ is zero, say $U'' = 0,$
  thus we deal with a direct decomposition $V = V'\oplus V''$ such that $U = U' \subseteq V'$.
  Since $U \to V$ is left minimal, $V'' = 0,$ thus $(U,V)$ is indecomposable.
\end{proof}

\medskip
\begin{remark}
  In Lemma~\ref{lem-thirteen-three},
  we need the assumption that $U$ is indecomposable. Namely, a $0$-picket $(0,[m])$ 
  is indecomposable, but $0 \to [m]$ is not left minimal. Also, if $U$ is non-zero and 
  decomposable, then $(U,U)$ is decomposable, but the identity map $U \to U$ is left minimal.
\end{remark}

\bigskip
Let us return to the setting we are interested in. There is given a 
(usually decomposable) object $V$ in $\Cal N(n)$, 
and a non-zero element $y$ in $V$, thus, we consider $U = \Lambda y.$ 
Then there is a direct decomposition
$V = V'\oplus V''$ in $\Cal N(n)$ such that the map $U \to V'$ is a left minimal version of 
the inclusion map $U \to V.$ Then, 
the pair $(U,V')$ is indecomposable in $\Cal S(n)$ and uniquely determined (up to isomorphism).

\medskip
Note that given $V$ in $\Cal N(n)$, and a non-zero element $y\in V$, let $U = \Lambda y$
and let $V = V'\oplus V''$ be a direct decomposition with $U \subseteq V'$ such that 
$U \to V'$ is a left minimal version of the inclusion map $U \to V.$ Then 
neither $V'$ nor $V''$ are usually uniquely determined:
As an example, take $V = [2,1]$ with generators $x_1,x_2$ annihilated
by $T^2$ and $T$, respectively, and take $y = Tx_1$.
Then we have in $\Cal N(n)$ the decompositions $V = \Lambda x_1\oplus \Lambda x_2 =
 \Lambda (x_1+x_2)\oplus \Lambda (Tx_1+x_2)$ and $y = Tx_1 = T(x_1+x_2).$ 

\subsection{A further combinatorial description.}
\label{sec-thirteen-four}

We use the considerations in Section~\ref{sec-thirteen-three}
in order to construct a second description of the
set of indecomposable objects $X = (U,V)$ with $bU = 1$. 

\medskip
Starting with an arbitrary non-empty subset $E$ of $\{1,\dots,n\}$, we may 
apply $H$ to the subspace $U$ in $M(E)$ and write $H(E) = H_{M(E)}(U)$
for the $T$-height sequence of $U$ in $M(E)$.

\medskip
Let us consider as an example the case
$E = \{2,3,5,6,8,9\}$, 
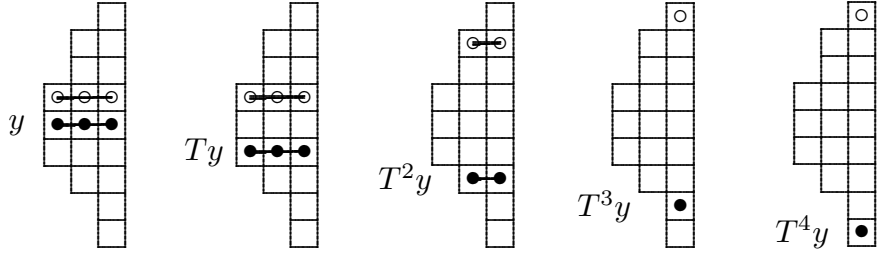
$$  
{\beginpicture
    \setcoordinatesystem units <.3cm,.3cm>
\put{\beginpicture
\multiput{} at 0 0  3 9 /
\plot 2 0  3 0  3 9  2 9  2 0 /
\plot 3 2  1 2  1 8  3 8 /
\plot 3 3  0 3  0 6  3 6 /
\plot 2 1  3 1 /
\plot 0 4  3 4 /
\plot 0 5  3 5  /
\plot 1 7  3 7 /
\multiput{$\bullet$} at .5 4.5  1.5 4.5  2.5 4.5 /
\plot 0.5 4.5  2.5 4.5 / 
\plot 0.5 4.45  2.5 4.55 / 
\plot 0.5 4.55  2.5 4.55 / 
\multiput{$\circ$} at .5 5.5  1.5 5.5  2.5 5.5 /
\plot 0.5 5.5  2.5 5.5 / 
\plot 0.5 5.45  2.5 5.55 / 
\plot 0.5 5.55  2.5 5.55 / 
\put{$y$} at -1 4.5
\endpicture} at 0 0
\put{\beginpicture
\multiput{} at 0 0  3 9 /
\plot 2 0  3 0  3 9  2 9  2 0 /
\plot 3 2  1 2  1 8  3 8 /
\plot 3 3  0 3  0 6  3 6 /
\plot 2 1  3 1 /
\plot 0 4  3 4 /
\plot 0 5  3 5  /
\plot 1 7  3 7 /
\multiput{$\bullet$} at .5 3.5  1.5 3.5  2.5 3.5 /
\plot 0.5 3.5  2.5 3.5 / 
\plot 0.5 3.45  2.5 3.55 / 
\plot 0.5 3.55  2.5 3.55 / 
\multiput{$\circ$} at .5 5.5  1.5 5.5  2.5 5.5 /
\plot 0.5 5.5  2.5 5.5 / 
\plot 0.5 5.45  2.5 5.55 / 
\plot 0.5 5.55  2.5 5.55 / 
\put{$Ty$} at -1.2 3.5
\endpicture} at 6.8 0
\put{\beginpicture
\multiput{} at 0 0  3 9 /
\plot 2 0  3 0  3 9  2 9  2 0 /
\plot 3 2  1 2  1 8  3 8 /
\plot 3 3  0 3  0 6  3 6 /
\plot 2 1  3 1 /
\plot 0 4  3 4 /
\plot 0 5  3 5  /
\plot 1 7  3 7 /
\multiput{$\bullet$} at  1.5 2.5  2.5 2.5 /
\plot 1.5 2.5  2.5 2.5 / 
\plot 1.5 2.45  2.5 2.55 / 
\plot 1.5 2.55  2.5 2.55 / 
\multiput{$\circ$} at  1.5 7.5  2.5 7.5 /
\plot 1.5 7.5  2.5 7.5 / 
\plot 1.5 7.45  2.5 7.55 / 
\plot 1.5 7.55  2.5 7.55 / 
\put{$T^2y$} at -1 2.5
\endpicture} at 14 0
\put{\beginpicture
\multiput{} at 0 0  3 9 /
\plot 2 0  3 0  3 9  2 9  2 0 /
\plot 3 2  1 2  1 8  3 8 /
\plot 3 3  0 3  0 6  3 6 /
\plot 2 1  3 1 /
\plot 0 4  3 4 /
\plot 0 5  3 5  /
\plot 1 7  3 7 /
\multiput{$\bullet$} at 2.5 1.5  /
\multiput{$\circ$} at  2.5 8.5 /
\put{$T^3y$} at -.3 1.5

\endpicture} at 21 0
\put{\beginpicture
\multiput{} at 0 0  3 9 /
\plot 2 0  3 0  3 9  2 9  2 0 /
\plot 3 2  1 2  1 8  3 8 /
\plot 3 3  0 3  0 6  3 6 /
\plot 2 1  3 1 /
\plot 0 4  3 4 /
\plot 0 5  3 5  /
\plot 1 7  3 7 /
\multiput{$\bullet$} at 2.5 0.5  /
\multiput{$\circ$} at  2.5 8.5 /
\put{$T^4y$} at .3 .5
\endpicture} at 28 0

\endpicture}
$$
For $1\le i \le 5$, we have inserted in the global space of $M(E)$
both the element $T^{i-1}y$ (using bullets), as
well as an element $z$ with $T^{H_i}z = T^{i-1}y$ (using small circles). 
We get $H(U) = (H_1,\dots,H_5) = (1,2,5,7,8).$

\medskip
\begin{proposition}
  \label{prop-thirteen-four}
  The map $H\:E\mapsto H(E)$ induces a bijection between
  the non-empty subsets of $\{1,\dots,n\}$ and the non-empty subsets of $\{0,\dots,n-1\}$.
\end{proposition}

\begin{proof}
  It is easy to see that $(H_1,\dots,H_t)$ 
  is the concatenation of
  $b$ intervalls of natural numbers, namely of the intervalls
  $$
  [d_1,\cdots,d_1+e_1-1],\ [d_2+e_1,\dots,d_2+e_2-1],\ \cdots,\ [d_b+e_{b-1}+1,\dots,d_b+e_b-1],
  $$
  in the prescribed order;
  here, the intervalls have length $e_i-e_{i-1}$ (with $e_0 = 0$) and end in $d_i+e_i-1$;
  in-between the intervalls is always a jump by at least 2. 
  (In the example above, the sequence $H(U)$ is cut into the intervalls 
  $[1,2],\ [5],\ [7,8].$) 
  This shows that we can recover from the sequence $(H_1,\dots,H_t),$ 
  first the sequence $e_i$, then the sequence
  $d_i$, thus we recover $E$.
  We see in this way that $H$ is an injective map, thus a bijective map.
\end{proof}

\subsection{The indecomposable objects $X = (U,V)$ with $bU = 1 = pX.$}
\label{sec-thirteen-five}

We have seen in Section~\ref{sec-seven-two} that there is a bijection between, on one hand, 
the isomorphism classes of the indecomposable objects
$X = (U,V)$ in $\Cal S(n)$ with $pX = 1$ and $bU = 1,$ and, on the other hand, the
strongly decreasing partitions bounded by $n$, 
by sending  $X = (U,V)$ to $[V].$ The indecomposable object $X$ corresponding to $\lambda$
has been denoted by $X = C_\lambda.$ 

\medskip
\begin{proposition}
  \label{prop-thirteen-five}
  The number of strongly decreasing partitions $\lambda$ of height $n$
  is the $n$-th Fibonacci number $F_n$.
\end{proposition}

\begin{proof}
  Let $\Pi(n)$ be the set of strongly decreasing partitions $\lambda$ with
  $\lambda_1 = n$. Clearly, $|\Pi(0)         | = 0 = F_0,$ and $|\Pi(1)| = 1 = F_1.$ 
  For $i \ge 2,$ let $\Pi_1(n)$ be the subset of $\Pi(n)$
  consisting of those partitions which end with $1$.
  If $\lambda = [\lambda_1,\dots,\lambda_b]\in \Pi_1(n)$, 
  then $[\lambda_1-2,\dots,\lambda_{b-1}-2]$ belongs to $\Pi(n-2)$; this map provides
  a bijection between $\Pi_1(n)$ and $\Pi(n-2)$. Similarly, if
  $\lambda = [\lambda_1,\dots,\lambda_b]\in \Pi(n)\setminus \Pi_1(n)$, 
  then $[\lambda_1-1,\dots,\lambda_{b}-1]$ belongs to $\Pi(n-1)$; this map provides
  a bijection between $\Pi(n)\setminus \Pi_1(n)$ and $\Pi(n-1)$. Altogether, we see
  that $|\Pi(n)| = |\Pi(n-2)|+|\Pi(n-1)|$.
\end{proof}

\medskip
\begin{corollary}
  \label{cor-thirteen-five}
  The number of indecomposable objects $X=(U,V)\in\Cal S(n)$
  with $pX=1$ and $bU=1$ is $F_{n+2}-1$ where $F_{n+2}$ is the $(n+2)$-nd Fibonacci number.
\end{corollary}

\begin{proof}
  Note $\sum_{i=1}^n F_i=F_{n+2}-F_2=F_{n+2}-1$.
\end{proof}

\subsection{The indecomposable objects $X = (U,V)$ with $bU = 1 = bW.$}
\label{sec-thirteen-six}

\begin{proposition}
  \label{prop-thirteen-six}
  Let $X = (U,V)$ be indecomposable. Then $bU \le 1,$ and $bW\le 1$
  if and only if $X$ is a picket or else a bipicket of the form $B(c_1,0,c_3,0,c_5).$
\end{proposition}

\begin{proof}
  We only have to show: 
  If $bU \le 1$ and $bW \le 1,$ then $bV \le 2.$ Now $bU = |\soc U|$ and $bW = |\soc W|.$ 
  We have $\soc V = \Hom(k,V)$ for $V \in \Cal N$. If
  we apply $\Hom(k,-)$ to the exact 
  sequence $0 \to U \to V \to W \to 0$, we obtain the exact sequence
  $0 \to \soc U \to \soc V \to \soc W.$ 
\end{proof}

\medskip
We see: If $X = (U,V)$ is indecomposable with $bU = 1 = bW$, then either $X$ is a
picket which is neither void nor full, or else it is the bipicket 
$B(c_1,0,c_3,0,c_5)$ (with $U = [c_1+c_3]$ and $W = [c_3+c_5]$, see
Section~\ref{sec-twelve-two}).
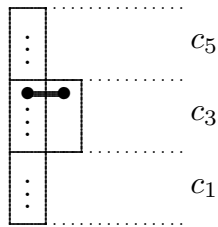
$$
{\beginpicture
   \setcoordinatesystem units <.4cm,.4cm>
\multiput{} at 0 0  2 6 /
\plot 0 0  1 0  1 6  0 6  0 0 /
\plot 0 2  2 2  2 4  0 4 /
\multiput{$\vdots$} at 0.5 1.1  0.5 3.1  0.5 5.1   /
\multiput{$\bullet$} at  0.5 3.6  1.5  3.6 /
\plot  0.5 3.6  1.5  3.6 /
\plot  0.5 3.55  1.5  3.55 /
\plot  0.5 3.65  1.5  3.65 /
\put{$c_1$} [l] at 5 1
\put{$c_3$} [l] at 5 3
\put{$c_5$} [l] at 5 5
\setdots <1mm>
\plot 1 0  5 0 /
\plot 1 2  5 2 /
\plot 1 4  5 4 /
\plot 1 6  5 6 /
\endpicture}
$$

\bigskip
\section{Objects in $\Cal S$ which are not gradable.}
\label{sec-fourteen}

All the examples which have been presented until now, were gradable. 
But we should stress that the concepts of mean and level apply to all non-zero
objects in $\Cal S$, not only for the gradable ones. 

\bigskip
As we have mentioned, all pickets are, of course, gradable, as are all
bipickets, see Proposition~\ref{prop-twelve-two}
Also, it has been shown in \cite{RS1} that all objects in $\Cal S(6)$
are gradable. 

\subsection{Some examples.}
\label{sec-fourteen-one}

Here is a typical non-gradable object $X$ in $\Cal S(7)$:
$$
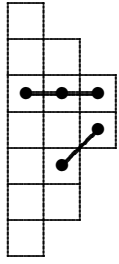
{\beginpicture
    \setcoordinatesystem units <.4cm,.4cm>
\multiput{} at 0 0  3 7 /
\plot 0 0  1 0  1 7  0 7  0 0 /
\plot 0 2  2 2 /
\plot 0 1  2 1  2 6  0 6 /
\plot 0 3  3 3  3 5  0 5 /
\plot 0 4  3 4 /
\multiput{$\bullet$} at 0.5 4.5  1.5 4.5  2.5 4.5  1.5 2.5  2.5 3.5 /
\plot 0.5 4.5  2.5 4.5 /
\plot 0.5 4.45  2.5 4.45 /
\plot 0.5 4.55  2.5 4.55 /
\plot 1.5 2.5  2.5 3.5 /
\plot 1.5 2.46  2.5 3.46 /
\plot 1.5 2.54  2.5 3.54 /
\endpicture}
$$
(to be precise: $X = (U,V)$ with $V$ generated by $v_1,v_2,v_3$ such that 
$T^7v_1 = T^5v_2 = T^2v_3 = 0$, and $U$ generated by $T^2v_1+Tv_2+v_3$ and
$T^3v_2+Tv_4$). We have
$$
 \bpar X = ([5,2],[7,5,2],[5,2]),
$$
thus $\bpar \tau^iX = \bpar X$ for all $i$. In particular, all objects $\tau^i X$ are 
central. More generally, let us consider the following example. 

\medskip
\begin{example}
  Let $n\ge 6$.
  let us consider the following object $X = (U,V)$ in $\Cal S(n)$
  with 
  $$
  \bpar X = ([n-2,2],[n,n-2,2],[n-2,2]),
  $$
  with $V$ generated by $v_1,v_2,v_3$ such that 
  $T^nv_1 = T^{n-2}v_2 = T^2v_3 = 0$, and $U$ generated by 
  $$  
  u_1 = T^2v_1+Tv_2+v_3\quad\text{and}\quad u_2 = T^{n-4}v_2+Tv_4.
  $$ 
  For example, for $n = 7$, we obtain the object $X$ shown above.
\end{example}

\medskip 
\begin{proposition}
  Let $n \ge 7.$ There are no indecomposable objects $\widetilde X$ in
  $\widetilde{\Cal S}(n)$ with $\bpar \widetilde X = ([n-2,2],[n,n-2,2],[n-2,2]).$
\end{proposition}

\begin{proof}
  All the algebras and the modules considered in the proof are $\mathbb Z$-graded,
  and the elements considered are homogeneous elements. 
  To simplify the notation, we will usually omit the use of a tilde.
  \smallskip
  
  We assume that $X = (U,V)$ is an indecomposable (graded!) object in $\widetilde{\Cal S}(n)$ with
  $\bpar X = ([n-2,2],[n,n-2,2],[n-2,2]).$ We decompose $V = V(1)\oplus V(2)\oplus V(3)$
  with $V(1),V(2),V(3)$ of length $n,\ n-2,\ 2$, respectively, with $V(i)$ generated in
  degree $d(i)$, say by $x_i$, for $1\le i \le 3.$ Without loss of generality, we may 
  assume that $d(1) = n.$
  \smallskip
  
  (1) {\it The canonical map $U \to V \to V(2)$ is not surjective.} Because otherwise, there
  is an element $u\in U$ of the form $u = \lambda_1x_1+x_2+\lambda_3x_3$ 
  (with $\lambda_1,\lambda_3 \in \Lambda$). Then $V$ is generated by $x_1, u, x_3$.
  Now $\Lambda u$ must have length $n-2$, thus $V = \Lambda u \oplus (V(1)+V(3))$.
  Since $\Lambda u\subseteq  U$, the modular law shows that $U = \Lambda u\oplus
  ((V(1)+V(3))\cap U),$ therefore $(U,V)$ is decomposable, a contradiction. 
  \smallskip
  
  (2) {\it The canonical map $U \to V \to V(3)$ is surjective.} Otherwise, $U$
  is contained in $\rad V(1)\oplus \rad V(2)\oplus \rad V(3) = \rad V$, and $W = V/U$
  has $bW = 3,$ a contradiction.
  \smallskip
  
  (3) {\it There is no element $u'\in U$ 
    such that $\Lambda u'$ has length $2$ and maps onto $V(3).$}
  Assume there is an element $u'\in U$ of the form $u' = \lambda_1x_1+
  \lambda_2 x_2+ x_3$.  Then $V$ is generated by $x_1, x_2, u'$. As in (1) we see that
  $V = \Lambda u'\oplus (V(1)+ V(2))$ and 
  $U = \Lambda u'\oplus
  ((V(1)+V(2))\cap U),$ therefore $(U,V)$ is decomposable, a contradiction. 
  \smallskip
  
  (4) According to (2) and (3), 
  there is an element $u\in U$ such that $\Lambda u$ has length $n-2$ and maps onto $V(3).$
  According to (1), 
  we can assume that $u = T^2x_1+T^mx_2 + x_3$ for some $m\ge 1$ (changing, if necessary, the
  generators $x_1,x_2$). 
  \smallskip
        
  (5) {\it We must have $m\ge 2.$} Assume, for the contrary, that $m = 1.$ 
  Then, for $m\ge 7$, the global space $V$ with $u$ inserted looks as follows:
  $$  
  
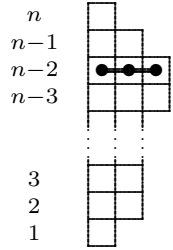
{\beginpicture
    \setcoordinatesystem units <.3cm,.3cm>
    \multiput{} at 0 0  3 9 /
    \plot 0 3.2  0 0  1 0  1 3.2 /
    \plot 0 1  2 1  2 3.2 /
    \plot 0 2  2 2 /
    \plot 0 3  2 3 /
    \plot 0 4.4  0 9  1 9  1 4.4 /
    \plot 0 8  2 8  2 4.4 /
    \plot 0 7  3 7  3 5  0 5 /
    \plot 0 6  3 6 /
    \setdots <1mm>
    \plot 0 3.3  0 4.7 /
    \plot 1 3.3  1 4.7 /
    \plot 2 3.3  2 4.7 /
    \put{$\ssize 1$} at -2 0.5
    \put{$\ssize 2$} at -2 1.5
    \put{$\ssize 3$} at -2 2.5
    \put{$\ssize n-3$} at -2 5.5
    \put{$\ssize n-2$} at -2 6.5
    \put{$\ssize n-1$} at -2 7.5
    \put{$\ssize n$} at -2 8.5
    \multiput{$\bullet$} at 0.5  6.5  1.5 6.5  2.5 6.5 /
    \setsolid
    \plot .5 6.5  2.5 6.5 /
    \plot .5 6.45  2.5 6.45 /
    \plot .5 6.55  2.5 6.55 /
    \endpicture}
  $$
  
  Since the partition of $U$ is $[n-2,2]$, we have $U = \Lambda u\oplus \Lambda u'$,
  where $\Lambda u'$ has length 2. 
  Since  $\Lambda u'$ has length 2, the degree $d'$ of $u'$ has to be $2,\ 3$ or $n-2.$
  Now $d' = n-2$ is impossible, according to (3). Also $d' = 2$ is impossible, since
  in this case both $\soc \Lambda u'$ and $\soc \Lambda u$ are the elements in
  $V$ of degree 1, thus $\Lambda u'\cap \Lambda u \neq 0.$
  It follows that $u'$ has degree 3. In this case, all elements of $V$ of degree at most 3
  belong to $U$ (here, we use $n\ge 7$). 
  As a consequence, $W = V/U$ is annihilated by $T^{n-3}$, thus its partition
  cannot be $[n-2,2].$
  \smallskip
  
  (6) {\it We can assume that $u = T^2x_1+x_3$.} 
  Namely, consider a case where $2\le m \le n-3.$ Since $T^mx_2 \neq 0,$ 
  the homogeneous element $u$ with degree $n-2$ shows that $d(2) = n-2+m,$ 
  thus $T^{m-2}x_2$ has degree $n$. Let $x'_1 = x_1+T^{m-2}x_2$. Then 
  $V$ is generated by $x'_1, x_2,x_3$ and the element $T^2x'_1 +x_3 = 
  T^2(x_1+T^{m-2}x_2)+x_3$ is just the element $u$, and therefore  
  belongs to $U$. Thus, it is sufficient to replace $x_1$ by $x'_1.$ 
  \smallskip
  
  (7) Finally, we are in the situation that $V$ is generated by $x_1,x_2,x_3$ and
  $U$ is generated by the element $u= T^2x_1+x_3$ with $\Lambda u$ of length $n-2$
  and an additional element $u'$ with $\Lambda u'$ of length 2.
  According to (3), $u'$ lies in $\Lambda x_1+\Lambda x_2 + \Lambda Tx_3,$
  thus in $\Lambda T^{n-2}x_1+\Lambda T^{n-4}x_2 + \Lambda Tx_3,$ since $\Lambda u$ has
  length 2. Now $T^{n-2}x_1 = T^{n-4}u$ belongs to $U$, thus we can assume that
  $u'$ belongs to $\Lambda T^{n-4}x_2 + \Lambda Tx_3,$ or even that $u' =  
  T^{n-4}x_2 + Tx_3.$ It follows that $T^{n-3}x_2 = Tu'$ belongs to $U$.
  Since also $T^{n-3}x_1$ belongs to $U$, we see also here that $W = V/U$ is annihilated by
  $T^{n-3}$, thus its partition
  cannot be $[n-2,2].$
\end{proof}

\medskip
\begin{corollary}
  For $n\ge 7$, the objects $X = (U,V)$ in $\Cal S(n)$ with
  $$
  \bpar X = ([n-2,2],[n,n-2,2],[n-2,2])
  $$ 
  presented above, are not gradable.
\end{corollary}

\subsection{These objects belong to homogeneous tubes.}
\label{sec-fourteen-two}

We recall from \cite[Corollary 6.5]{RS2}
that all indecomposable objects in $\Cal S(n)$ for $n\geq 6$
occur in tubes of rank 1, 2, 3, or 6; the tubes of rank 1 are said to be
homogeneous.	

\medskip
\begin{proposition}
  \label{prop-fourteen-two}
  Let $n> 6$. Let $X$ be one of the objects as defined in Section~\ref{sec-fourteen-one}.
  Then $\tau_nX = X.$
\end{proposition}

\begin{proof}
  We work in the category $\Cal N(n)$ of all $\Lambda$-modules, where $\Lambda =
  k[T]/\langle T^n\rangle$.  
  Canonical inclusion maps in $\Cal N(n)$ are denoted by $\mu,$ canonical epimorphisms by
  $\varepsilon$, and we denote the multiplication map by $T$ just by $T$.
  
  If we use the generators $v_1,v_2,-v_3$ of $[n,n-2,2]$,
  we see that $X = (U,V)$ is given by the inclusion map
  $$
  \left[\smallmatrix \mu & 0 \cr
    T & \mu \cr
    -\varepsilon & -T \endsmallmatrix \right]\:[n-2,2] \to [n,n-2,2],
  $$
  or also by the inclusion map
  $$
  X' = \left[\smallmatrix 1-T^{n-6} & 0 & 0\cr
    0 & 1 & 0 \cr
    0 & 0 & 1  \endsmallmatrix \right]\cdot
  \left[\smallmatrix \mu & 0 \cr
    T & \mu \cr
    -\varepsilon & -T \endsmallmatrix \right]\:[n-2,2] \to [n,n-2,2]
  $$
  (obtained from $X$ by applying an automorphism of $[n,n-2,2]$). 
  \smallskip
  
  We are going to form $\tau_n X = \tau_n X' = \Mimo\tau_\Lambda\Cok X'$ (see \cite{RS2}). 
  
  There is the commutative diagram
  $$
  {\beginpicture
    \setcoordinatesystem units <3cm,1cm>
    \multiput{$[n-2,2]$} at 0 0  1 -1.1   2 0 /
    \multiput{$[n]$} at  0.7  1.1  1.3 1.1 /
    \arr{0.3 0.4}{0.5 .9}
    \arr{0.3 -.4}{0.7 -.8}
    \arr{1.5 .9}{1.7 0.4}
    \arr{1.3 -.8}{1.7 -.4}
    \arr{0.9 1.1}{1.1 1.1}
    \put{$\left[\smallmatrix \mu & 0 
        \endsmallmatrix \right]$} at 0.2 .8
    \put{$\ssize 1-T^{n-6}$} at 1 1.4
    \put{$\left[\smallmatrix \varepsilon \cr 0 
        \endsmallmatrix \right]$} at 1.8 .8
    \put{$\left[\smallmatrix T & \mu \cr
        -\varepsilon & -T 
        \endsmallmatrix \right]$} at 0.2 -.95
    \put{$\left[\smallmatrix T & \mu \cr
        \varepsilon & T 
        \endsmallmatrix \right]$} at 1.75 -.95
    \endpicture}
  $$
  since for both paths, the composition from left to right is equal to 
  $\left[\smallmatrix T^2-T^{n-4} & 0 \cr
    0 & 0 
    \endsmallmatrix \right].$
  
  It shows that the following sequence is exact:
  $$
  \CD
  0 @>>> [n-2,2] 
  @>X'>> [n,n-2,2] @>Y>> [n-2,2] @>>> 0
  \endCD
  $$
  with
  $$  Y = 
  \left[\smallmatrix -\varepsilon & T & \mu  \cr
    0 & \varepsilon & T \endsmallmatrix \right]\:[n,n-2,2] \to [n-2,2]
  $$
  (namely, the commutativity shows that the composition $YX'$ is zero; we know that 
  the map $X'$ is a monomorphism, and we easily see that the map $Y$ is an
  epimorphism; finally, the length of the middle term is the sum of the length of the end-terms).
  
  Now $\tau_\Lambda Y$ is given by the matrix
  $$
  \tau_\Lambda Y =
  \left[\smallmatrix  T & \mu  \cr
    \varepsilon & T \endsmallmatrix \right]\:[n-2,2] \to [n-2,2].
  $$
  It remains to form 
  $$ 
  \Mimo  \tau_\Lambda Y = 
  \left[\smallmatrix  \mu & 0 \cr
    T & \mu  \cr
    \varepsilon & T \endsmallmatrix \right]\:[n-2,2] \to [n,n-2,2],
  $$
  but this is just $X$. This completes the proof.
\end{proof}

\medskip
\begin{corollary}
  \label{cor-fourteen-two}
  Each category $\Cal S(n)$ for $n\geq 6$ contains homogeneous tubes.
\end{corollary}

\medskip
\begin{remark}
  Proposition~\ref{prop-fourteen-two} asserts that for $n > 6$, the object
  defined in Section~\ref{sec-fourteen-one} belongs to a homogeneous tube. In contrast, for $n = 6$,
  the object lies in the tube of rank 3 and rationality index 0.
\end{remark}

\subsection{Homogeneous tubes containing graded objects.}
\label{sec-fourteen-three}

\begin{proposition}
  There exist homogeneous tubes of gradable objects in $\Cal S(n)$ 
  if and only if $6$ divides $n$.
\end{proposition}

\begin{proof}
  If $6$ does not devide $n$, then see \cite[Corollary 6.6]{RS2}.

Thus, let us assume that $6$ divides $n$. We want to exhibit objects $M_c(n)$ which 
are homogeneous.
Write $n=6\ell$ and define $M_c(n)=(U,V)$ with $V$ generated by elements $v_1$, $v_2$, $v_3$ of height $6\ell$,
$4\ell$ and $2\ell$, respectively, and let $U$ be generated by $u_1=T^{2\ell}v_1+c_0T^\ell v_2+c_1v_3$ and $u_2=T^{2\ell}v_2+T^\ell v_3$,
thus $\bpar M_c(n)=([4\ell,2\ell],[6\ell,4\ell,2\ell],[4\ell,2\ell])$.  Then $M_c(n)$ is gradable,
in fact, if $M_c(n)=\pi \widetilde M_c(n)$ then
$\tau_{\widetilde{\Cal S}}\widetilde M_c(n)=\widetilde M_c(n)\big[\frac{n-6}6\big]=\widetilde M_c(n)[\ell-1]$,
for $(c_0:c_1)\neq (0:1),(1:0),(1:1)$.  Namely, the representation $\widetilde M_c(n)$
is obtained from $\widetilde M_c(6)$ by replicating the vector spaces $\ell$ times.
Applying the dual of the transpose in $\widetilde{\Cal S}(n)$ results in a representation
obtained from $\tau_{\widetilde{\Cal S}}\widetilde M_c(6)$ by a corresponding replication of its
vector spaces.
It follows that $M_c(n)$ lies on the mouth of a homogeneous tube.
\end{proof}

\subsection{AR-components in $\Cal S(n)$ containing gradable objects.}
\label{sec-fourteen-four}

For $n\geq 6$, every Auslander-Reiten component in $\Cal S(n)$ is an $m$-tube for $m$ a divisor
of six \cite{RS2}.  We have seen in Section~\ref{sec-fourteen-three} that, unless $n$ is a multiple of $6$,
no Auslander-Reiten component in $\Cal S(n)$ which contains
gradable objects is homogeneous. For tubes of rank 2 or 3 we have the following statement.

\medskip 
\begin{proposition}
  Let $n\geq 6$.
  \begin{itemize}[leftmargin=3em]
  \item[(a)] The category $\Cal S(n)$ has $3$-tubes consisting of gradable objects
    if and only if $n$ is even.
  \item[(b)] $\Cal S(n)$ has $2$-tubes consisting of gradable objects
    if and only if $n$ is divisible by $3$.
  \end{itemize}
\end{proposition}

\begin{proof}
  We first show the ``only if'' direction.
Let $X$ be a gradable object in $\Cal S(n)$ which is indecomposable
and not projective. Denote by $m$ its minimal $\tau$-period.
Of course, we know that $m$ is a divisor of $6$, thus $mm' = 6$ for some integer $m'$.
We show that $m'$ is a divisor of $n-6$.

Since $X$ is gradable, there is $\widetilde X$ in $\widetilde {\Cal S}(n)$ with
$\pi\widetilde X = X.$ Since $\tau^m X = X,$ there is an integer $s$ such that
$\widetilde\tau^m \widetilde X=\widetilde X [s]$;
here, $\widetilde \tau$ is the Auslander-Reiten translation in $\widetilde{\Cal S}(n)$
and $[1]$ is the shift-operator in $\widetilde{\Cal S}(n).$ 
It follows that $\widetilde\tau^6 \widetilde X =\widetilde\tau^{mm'} \widetilde X = 
\widetilde X[sm'].$ 

According to \cite{RS2}, the Auslander-Reiten translation in any $\widetilde{\Cal S}(n)$
satisfies $\widetilde\tau^6 M = M[n-6]$ for any
indecomposable non-projective object $M$ in $\widetilde{\Cal S}(n)$, therefore
$\widetilde\tau^6 \widetilde X = \widetilde X[n-6]$.

Altogether, we see that $\widetilde X[sm'] = \widetilde\tau^6 \widetilde X = \widetilde X[n-6]$,
therefore $sm' = n-6,$ thus $m'$ is a divisor of $n-6$.
In particular, if $m=3$ then $n$ must be even and if $m=2$ then $n$ must be divisible by $3$.

For the converse direction, if $n\geq6$ is even, we can consider the picket $(0;[\frac n2])$
which is determined uniquely by its partition type.  The mouth of the $3$-tube which contains
this object has the following shape.
$$\textstyle (0;[\frac n2])\;\stackrel\tau\longmapsfrom \;([\frac n2];[n])\;
\stackrel\tau\longmapsfrom \;([\frac n2];[\frac n2])\;\stackrel\tau\longmapsfrom \;(0;[\frac n2])$$

Similarly, if $n$ is divisible by $3$, the picket $([\frac n3];[\frac{2n}3])$ occurs
on the mouth of a $2$-tube.
$$\textstyle ([\frac n3];[\frac{2n}3])\;\stackrel\tau\longmapsfrom\;
([\frac{2n}3];[n,\frac n3];[\frac{2n}3])\;\stackrel\tau\longmapsfrom\;
([\frac n3];[\frac{2n}3])$$
\end{proof}

\medskip
\begin{corollary}
  Any Auslander-Reiten component in $\Cal S(7)$ which
  contains gradable objects is a $6$-tube.
\end{corollary}

\subsection{Gradability of objects.}
\label{sec-fourteen-five}

As we have mentioned, all objects in $\Cal S(6)$ are gradable, see \cite{RS1}. 
The novice may be surprised by this fact, thus let us consider explicitly three objects 
$X = (U,V)$ in $\Cal S(6)$ which look similar to the objects in $\Cal S(7)$ which have been 
shown to be non-gradable. In all cases, let $V$ be 
generated by $v_1,v_2,v_3$ of height $6,4,2,$ respectively, and $U$ 
generated by $u_1 =  T^2v_1+Tv_2+v_3$ and a further element $u_2.$

\smallskip
(1) Let $u_2 = T^2v_2+v_3.$ 
Then $V$ is generated by  $v'_1 = v_1,\ v'_2 = v_2-Tv_2,\ v'_3 = v_3+T^2v_2,$
and $U$
is generated by $u_1 = T^2v'_1+Tv'_2$ and $u_2 = v'_3,$ thus even decomposable. 

\smallskip
(2) Let $u_2 = T^3v_2+v_3.$ 
Then $V$ is generated by  $v'_1 = v_1,\ v'_2 = v_2-T^2v_2,\ v'_3 = v_3+T^3v_2,$
and $U$
is generated by $u_1 = T^2v'_1+Tv'_2$ and $u_2 = v'_3,$ thus again decomposable. 

\smallskip
(3) Let $u_2 = T^3v_2+Tv_3.$ 
Then $V$ is generated by  $v'_1 = v_1,\ v'_2 = v_2-Tv_2,\ v'_3 = v_3+T^2v_2,$
and $U$
is generated by $u_1 = T^2v'_1+Tv'_2+v'_3$ and $u_2 = Tv'_3,$ thus indecomposable, but 
obviously gradable.

\medskip
Here are the corresponding pictures:
$$
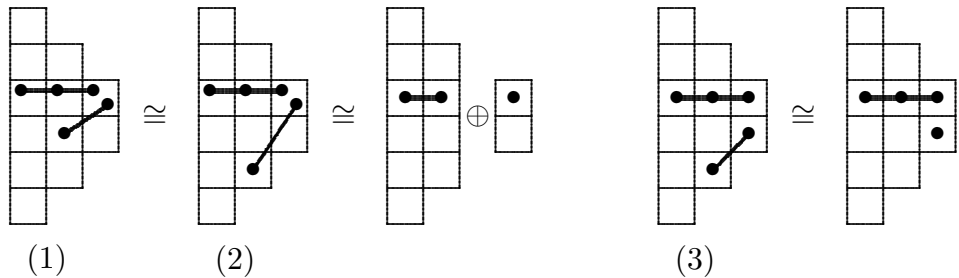
{\beginpicture
    \setcoordinatesystem units <.4cm,.4cm>
\multiput{} at 0 0  4 6 /
\plot 0 0  1 0  1 6  0 6  0 0 /
\plot 0 3  3 3 /
\plot 0 1  2 1  2 5  0 5 /
\plot 0 2  3 2  3 4  0 4 /
\multiput{$\bullet$} at 0.3 3.7  1.3 3.7  2.3 3.7   1.5 2.5  2.7 3.3 /
\plot 0.3 3.7   2.3 3.7 /
\plot 0.3 3.75  2.3 3.75 /
\plot 0.3 3.65  2.3 3.65 /
\plot 1.5 2.5   2.7 3.3  /
\plot 1.5 2.46  2.7 3.26 /
\plot 1.5 2.54  2.7 3.34 /
\put{$\cong$} at 4 3 
\put{(1)} at 1 -1
\endpicture}
\quad
{\beginpicture
    \setcoordinatesystem units <.4cm,.4cm>
\multiput{} at 0 0  4 6 /
\plot 0 0  1 0  1 6  0 6  0 0 /
\plot 0 3  3 3 /
\plot 0 1  2 1  2 5  0 5 /
\plot 0 2  3 2  3 4  0 4 /
\multiput{$\bullet$} at 0.3 3.7  1.3 3.7  2.3 3.7   1.5 1.5  2.7 3.3 /
\plot 0.3 3.7   2.3 3.7 /
\plot 0.3 3.75  2.3 3.75 /
\plot 0.3 3.65  2.3 3.65 /
\plot 1.5 1.5   2.7 3.3  /
\plot 1.5 1.46  2.7 3.26 /
\plot 1.5 1.54  2.7 3.34 /
\put{$\cong$} at 4 3 
\put{(2)} at 1 -1
\endpicture}
\quad
{\beginpicture
    \setcoordinatesystem units <.4cm,.4cm>
\multiput{} at 0 0  4 6 /
\plot 0 0  1 0  1 6  0 6  0 0 /
\plot 0 3  2 3 /
\plot 0 1  2 1  2 5  0 5 /
\plot 0 2  2 2 /
\plot 0 4  2 4 /
\plot 3 2  4 2  4 4  3 4  3 2 /
\plot 3 3  4 3 /
\multiput{$\bullet$} at 0.5 3.5  1.5 3.5  3.5 3.5 /
\plot 0.5 3.5   1.5 3.5  /
\plot 0.5 3.55  1.5 3.55 /
\plot 0.5 3.45  1.5 3.45 /
\put{$\oplus$} at 2.5 3 
\endpicture}
\qquad\qquad
{\beginpicture
    \setcoordinatesystem units <.4cm,.4cm>
\multiput{} at 0 0  4 6 /
\plot 0 0  1 0  1 6  0 6  0 0 /
\plot 0 3  3 3 /
\plot 0 1  2 1  2 5  0 5 /
\plot 0 2  3 2  3 4  0 4 /
\multiput{$\bullet$} at 0.5 3.5  1.5 3.5  2.5 3.5    1.5 1.5  2.5 2.5 /
\plot 0.5 3.5   2.5 3.5 /
\plot 0.5 3.55  2.5 3.55 /
\plot 0.5 3.45  2.5 3.45 /
\plot 1.5 1.5   2.5 2.5 /
\plot 1.5 1.46  2.5 2.46 /
\plot 1.5 1.54  2.5 2.54 /
\put{$\cong$} at 4 3 
\put{(3)} at 1 -1
\endpicture}
\quad
{\beginpicture
    \setcoordinatesystem units <.4cm,.4cm>
\multiput{} at 0 0  3 6 /
\plot 0 0  1 0  1 6  0 6  0 0 /
\plot 0 3  3 3 /
\plot 0 1  2 1  2 5  0 5 /
\plot 0 2  3 2  3 4  0 4 /
\multiput{$\bullet$} at 0.5 3.5  1.5 3.5  2.5 3.5     2.5 2.5 /
\plot 0.5 3.5   2.5 3.5 /
\plot 0.5 3.55  2.5 3.55 /
\plot 0.5 3.45  2.5 3.45 /
\endpicture}
$$
Thus, we have found objects isomorphic to (1), (2) and (3) which are gradable.

\bigskip
\section{Gallery: Some Examples.}
\label{sec-fifteen}

We are going to exhibit some examples of objects in $\Cal S$ which are
worth to be aware of. 

\subsection{The cases $\Cal S(n)$ with $n \le 5$.}
\label{sec-fifteen-one}

In this Section we present the pr-triangles for the representation finite
categories $\Cal S(n)$ where $n\leq 5$. The corresponding Auslander-Reiten quivers
can be found in \cite{RS1}.  

In view of later questions, it seems to be interesting to watch in Sections~\ref{sec-fifteen-two}--(g) and
\ref{sec-fifteen-three}
how the lines $p=1$ (horizontal) and $q=3$ (diagonal) evolve as $n$ increases,
and which objects arise in the regions $2<q<3$ and $3<q<4$.

\medskip
{\bf The case  $n=1$.} (Note that this case is not covered in \cite{RS1}, but is, of course, trivial:
It is the so-called Basis Extension Theorem from linear algebra which provides
the classification of the indecomposable objects in
$\Cal S(1)$.)
$$
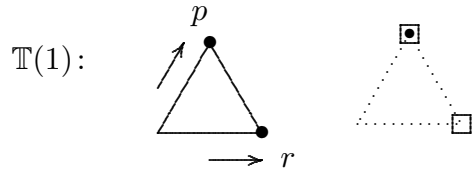
{\beginpicture
   \setcoordinatesystem units <.57735cm,1cm>
\put{\beginpicture
\put{$\mathbb T(1)\:$} at -3 .8
\multiput{} at -1 -1  1 2 /
\plot -1 0  1 0  0 1 -1 0 /
\multiput{$\bullet$} at 0 1  1 0 /
\setsolid 
\arr{0 -.3}{1 -.3}
\put{$r$} at 1.5 -.3
\arr{-1 .5}{-.5 1}
\put{$p$} at -.2 1.3
\endpicture} at 0 0 
\put{\beginpicture
\multiput{} at -1 0  1 0  0 1 /
\setdots <1mm>
\plot -1 0  1 0  0 1  -1 0 /

\put{\beginpicture
   \setcoordinatesystem units <.2cm,.2cm>
  \multiput{} at 0 -1  1 2 /
  \setsolid
  \plot 0 0  1 0  1 1  0 1  0 0 /
  \endpicture} at 1 0
\put{\beginpicture
   \setcoordinatesystem units <.2cm,.2cm>
  \multiput{} at 0 0  1 1 /
  \setsolid
  \plot 0 0  1 0  1 1  0 1  0 0 /
  \put{$\ssize\bullet$} at 0.5 0.5 
  \endpicture} at 0 1
\endpicture} at 5 0 
\endpicture}
$$
	
{\bf The case  $n=2$.} The objects in $\Cal S(2)$ are determined in \cite[(6.2)]{RS1}.
We stress: {\it Each indecomposable object $(U,V)$ in $\Cal S(2)$ lies on the
boundary of $\mathbb T(2)$; in particular, it is a picket.}
$$
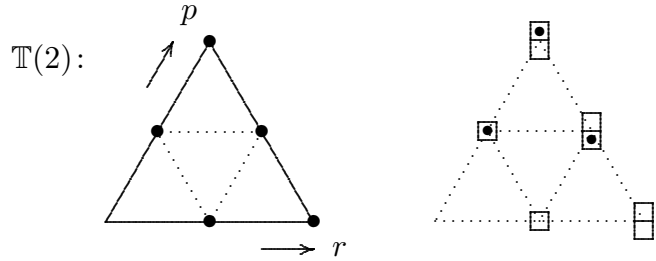
{\beginpicture
   \setcoordinatesystem units <.57735cm,1cm>
\put{\beginpicture
\put{$\mathbb T(2)\:$} at -3 1.8
\multiput{} at -2 0  2 0  0 2 /
\plot -2 0  2 0  0 2 -2 0 /
\setdots <1mm>
\plot 0 0  -1 1  1 1  0 0 /
\multiput{$\bullet$} at 0 0  -1 1  2 0  1 1  0 2 /
\setsolid 
\arr{1 -.3}{2 -.3}
\put{$r$} at 2.5 -.3
\arr{-1.2 1.5}{-.7 2}
\put{$p$} at -.4 2.3
\endpicture} at 0 0 
\put{\beginpicture
\multiput{} at -2 0  2 0  0 2 /
\setdots <1mm>
\plot -2 0  2 0  0 2  -2 0 /
\plot 0 0  -1 1  1 1  0 0 /

\put{\beginpicture
   \setcoordinatesystem units <.2cm,.2cm>
  \multiput{} at 0 0  1 1 /
  \setsolid
  \plot 0 0  1 0  1 1  0 1  0 0 /
  \endpicture} at 0 0
\put{\beginpicture
   \setcoordinatesystem units <.2cm,.2cm>
  \multiput{} at 0 0  1 2 /
  \setsolid
  \plot 0 0  1 0  1 2  0 2  0 0 /
  \plot 0 1  1 1 /
  \endpicture} at 2 0
\put{\beginpicture
   \setcoordinatesystem units <.2cm,.2cm>
  \multiput{} at 0 0  1 1 /
  \setsolid
  \plot 0 0  1 0  1 1  0 1  0 0 /
  \put{$\ssize\bullet$} at 0.5 0.5 
  \endpicture} at -1 1
\put{\beginpicture
   \setcoordinatesystem units <.2cm,.2cm>
  \multiput{} at 0 0  1 2 /
  \setsolid
  \plot 0 0  1 0  1 2  0 2  0 0 /
  \plot 0 1  1 1 /
  \put{$\ssize\bullet$} at 0.5 1.5 
  \endpicture} at 0 2
\put{\beginpicture
   \setcoordinatesystem units <.2cm,.2cm>
  \multiput{} at 0 0  1 2 /
  \setsolid
  \plot 0 0  1 0  1 2  0 2  0 0 /
  \plot 0 1  1 1 /
  \put{$\ssize\bullet$} at 0.5 .5 
  \endpicture} at 1 1
\endpicture} at 7 0 
\endpicture}
$$
	
{\bf The case  $n=3$.} There is a single indecomposable object which is not a picket,
namely $E = E_2^2$ with uwb-vector $\frac{2|2}2,$ thus with pr-vector $(1,1) = z(3).$
{\it Each indecomposable object $(U,V)$ in $\Cal S(3)$ lies on the grid of 
$\mathbb T(3)$ with integeral coefficients.}
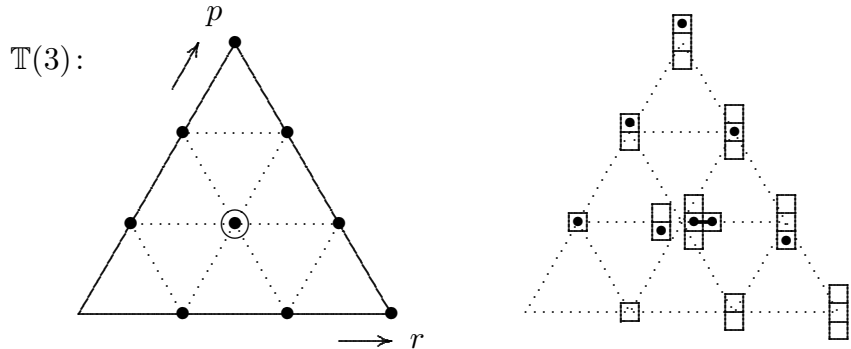
$$
{\beginpicture
   \setcoordinatesystem units <.57735cm,1cm>
\put{\beginpicture
\put{$\mathbb T(3)\:$} at -3.5 2.8
\multiput{} at -3 0  3 0  0 3 /
\plot -3 0  3 0  0 3  -3 0 /
\setdots <1mm>
\plot -2 1  -1 0  0 1  1 0  2 1  -2 1 /
\plot -1 2  0 1  1 2  -1 2 /
\multiput{$\bullet$} at -1 0  1 0  3 0  -2 1  0 1  2 1  -1 2  1 2  0 3 /
\put{$\bigcirc$} at 0 1 
\setsolid 
\arr{2 -.3}{3 -.3}
\put{$r$} at 3.5 -.3
\arr{-1.2 2.5}{-.7 3}
\put{$p$} at -.4 3.3
\endpicture} at 0 0 
\put{\beginpicture
\multiput{} at -3 0  3 0  0 3 /
\setdots <1mm>
\plot -3 0  3 0  0 3  -3 0 /
\plot -2 1  -1 0  0 1  1 0  2 1  -2 1 /
\plot -1 2  0 1  1 2  -1 2 /

\put{\beginpicture
   \setcoordinatesystem units <.2cm,.2cm>
  \multiput{} at 0 0  1 1 /
  \setsolid
  \plot 0 0  1 0  1 1  0 1  0 0 /
  \endpicture} at -1 0
\put{\beginpicture
   \setcoordinatesystem units <.2cm,.2cm>
  \multiput{} at 0 0  1 2 /
  \setsolid
  \plot 0 0  1 0  1 2  0 2  0 0 /
  \plot 0 1  1 1 /
  \endpicture} at 1 0
\put{\beginpicture
   \setcoordinatesystem units <.2cm,.2cm>
  \multiput{} at 0 0  1 3 /
  \setsolid
  \plot 0 0  1 0  1 3  0 3  0 0 /
  \plot 0 1  1 1 /
  \plot 0 2  1 2 /
  \endpicture} at 3 0

\put{\beginpicture
   \setcoordinatesystem units <.2cm,.2cm>
  \multiput{} at 0 0  1 1 /
  \setsolid
  \plot 0 0  1 0  1 1  0 1  0 0 /
  \put{$\ssize\bullet$} at 0.5 0.5 
  \endpicture} at -2 1
\put{\beginpicture
   \setcoordinatesystem units <.2cm,.2cm>
  \multiput{} at 0 0  1 2 /
  \setsolid
  \plot 0 0  1 0  1 2  0 2  0 0 /
  \plot 0 1  1 1 /
  \put{$\ssize\bullet$} at 0.5 1.5 
  \endpicture} at -1 2
\put{\beginpicture
   \setcoordinatesystem units <.2cm,.2cm>
  \multiput{} at 0 0  1 3 /
  \setsolid
  \plot 0 0  1 0  1 3  0 3  0 0 /
  \plot 0 1  1 1 /
  \plot 0 2  1 2 /
  \put{$\ssize\bullet$} at 0.5 2.5 
  \endpicture} at 0 3
\put{\beginpicture
   \setcoordinatesystem units <.2cm,.2cm>
  \multiput{} at 0 0  1 3 /
  \setsolid
  \plot 0 0  1 0  1 3  0 3  0 0 /
  \plot 0 1  1 1 /
  \plot 0 2  1 2 /
  \put{$\ssize\bullet$} at 0.5 1.5 
  \endpicture} at 1 2
\put{\beginpicture
   \setcoordinatesystem units <.2cm,.2cm>
  \multiput{} at 0 0  1 3 /
  \setsolid
  \plot 0 0  1 0  1 3  0 3  0 0 /
  \plot 0 1  1 1 /
  \plot 0 2  1 2 /
  \put{$\ssize\bullet$} at 0.5 0.5 
  \endpicture} at 2 1

\put{\beginpicture
   \setcoordinatesystem units <.2cm,.2cm>
  \multiput{} at 0 0  1 2 /
  \setsolid
  \plot 0 0  1 0  1 2  0 2  0 0 /
  \plot 0 1  1 1 /
  \put{$\ssize\bullet$} at 0.5 .5 
  \endpicture} at -.4 1

\put{\beginpicture
   \setcoordinatesystem units <.2cm,.2cm>
  \multiput{} at 0 0  2 3 /
  \setsolid
  \plot 0 0  1 0  1 3  0 3  0 0 /
  \plot 0 1  2 1  2 2  0 2 /
  \multiput{$\ssize\bullet$} at 0.5 1.5  1.5 1.5 /
  \plot 0.5 1.5  1.5 1.5 /
  \plot 0.5 1.45  1.5 1.45 /
  \plot 0.5 1.55  1.5 1.55 /
  \endpicture} at .4 1 
\endpicture} at 9 0 
\endpicture}
$$

{\bf The case $n = 4$.} The $\tau_n$-orbits of the pickets provide 
all but 3 indecomposable objects, 
the missing ones form the $\tau_n$-orbit of $E_2^2=C_{(3,1)}$. 
$$
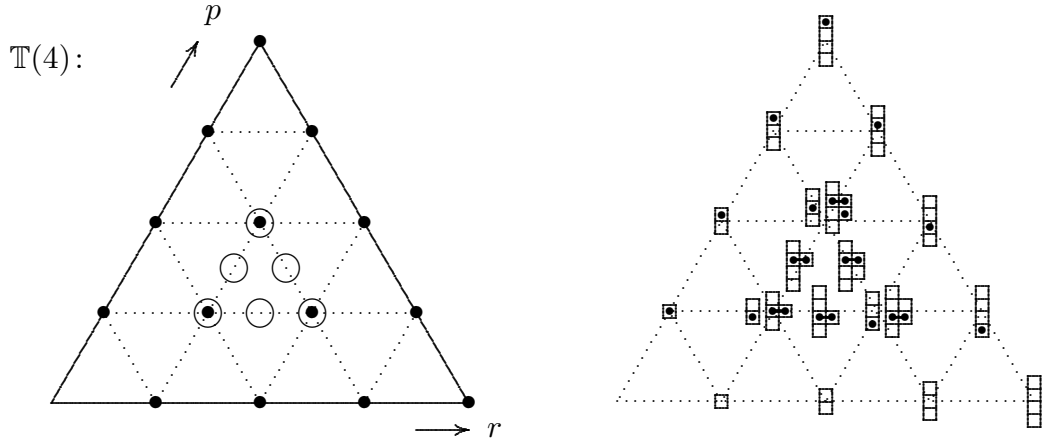
{\beginpicture
   \setcoordinatesystem units <.57735cm,1cm>
\put{\beginpicture
\put{$\mathbb T(4)\:$} at -4 3.8
\multiput{} at -4 0  4 0  0 4 /
\plot -4 0  4 0  0 4  -4 0 /
\setdots <1mm>
\plot -2 0  -3 1  3 1  2 0  -1 3  1 3  -2 0  /
\plot 0 0  -2 2  2 2  0 0 /
\multiput{$\bullet$} at -2 0  0 0  2 0  4 0  -3 1  -1 1  1 1  3 1  -2 2  0 2  2 2  -1 3  1 3  0 4 /
\multiput{$\bigcirc$} at -1 1  0 1  1 1  -.5 1.5  .5 1.5  0 2 / 
\setsolid 
\arr{3 -.3}{4 -.3}
\put{$r$} at 4.5 -.3
\arr{-1.7 3.5}{-1.2 4}
\put{$p$} at -.9 4.3
\endpicture} at 0 0 
\put{\beginpicture
\multiput{} at -4 0  4 0  0 4 /
\setdots <1mm>
\plot -4 0  4 0  0 4  -4 0 /
\plot -2 0  -3 1  3 1  2 0  -1 3  1 3  -2 0  /
\plot 0 0  -2 2  2 2  0 0 /

\put{\beginpicture
   \setcoordinatesystem units <.14cm,.14cm>
  \multiput{} at 0 0  1 1 /
  \setsolid
  \plot 0 0  1 0  1 1  0 1  0 0 /
  \endpicture} at -2 0
\put{\beginpicture
   \setcoordinatesystem units <.14cm,.14cm>
  \multiput{} at 0 0  1 2 /
  \setsolid
  \plot 0 0  1 0  1 2  0 2  0 0 /
  \plot 0 1  1 1 /
  \endpicture} at 0 0
\put{\beginpicture
   \setcoordinatesystem units <.14cm,.14cm>
  \multiput{} at 0 0  1 3 /
  \setsolid
  \plot 0 0  1 0  1 3  0 3  0 0 /
  \plot 0 1  1 1 /
  \plot 0 2  1 2 /
  \endpicture} at 2 0

\put{\beginpicture
   \setcoordinatesystem units <.14cm,.14cm>
  \multiput{} at 0 0  1 1 /
  \setsolid
  \plot 0 0  1 0  1 1  0 1  0 0 /
  \put{$\sssize\bullet$} at 0.5 0.5 
  \endpicture} at -3 1
\put{\beginpicture
   \setcoordinatesystem units <.14cm,.14cm>
  \multiput{} at 0 0  1 2 /
  \setsolid
  \plot 0 0  1 0  1 2  0 2  0 0 /
  \plot 0 1  1 1 /
  \put{$\sssize\bullet$} at 0.5 1.5 
  \endpicture} at -2 2
\put{\beginpicture
   \setcoordinatesystem units <.14cm,.14cm>
  \multiput{} at 0 0  1 3 /
  \setsolid
  \plot 0 0  1 0  1 3  0 3  0 0 /
  \plot 0 1  1 1 /
  \plot 0 2  1 2 /
  \put{$\sssize\bullet$} at 0.5 2.5 
  \endpicture} at -1 3
\put{\beginpicture
   \setcoordinatesystem units <.14cm,.14cm>
  \multiput{} at 0 0  1 3 /
  \setsolid
  \plot 0 0  1 0  1 3  0 3  0 0 /
  \plot 0 1  1 1 /
  \plot 0 2  1 2 /
  \put{$\sssize\bullet$} at 0.5 1.5 
  \endpicture} at -.25 2.15
\put{\beginpicture
   \setcoordinatesystem units <.14cm,.14cm>
  \multiput{} at 0 0  1 3 /
  \setsolid
  \plot 0 0  1 0  1 3  0 3  0 0 /
  \plot 0 1  1 1 /
  \plot 0 2  1 2 /
  \put{$\sssize\bullet$} at 0.5 0.5 
  \endpicture} at .9 1

\put{\beginpicture
   \setcoordinatesystem units <.14cm,.14cm>
  \multiput{} at 0 0  1 2 /
  \setsolid
  \plot 0 0  1 0  1 2  0 2  0 0 /
  \plot 0 1  1 1 /
  \put{$\sssize\bullet$} at 0.5 .5 
  \endpicture} at -1.4 1

\put{\beginpicture
   \setcoordinatesystem units <.14cm,.14cm>
  \multiput{} at 0 0  2 3 /
  \setsolid
  \plot 0 0  1 0  1 3  0 3  0 0 /
  \plot 0 1  2 1  2 2  0 2 /
  \multiput{$\sssize\bullet$} at 0.5 1.5  1.5 1.5 /
  \plot 0.5 1.5  1.5 1.5 /
  \plot 0.5 1.45  1.5 1.45 /
  \plot 0.5 1.55  1.5 1.55 /
  \endpicture} at -.9 1 
\put{\beginpicture
   \setcoordinatesystem units <.14cm,.14cm>
  \multiput{} at 0 0  1 4 /
  \setsolid
  \plot 0 0  1 0  1 4  0 4  0 0 /
  \plot 0 1  1 1 /
  \plot 0 2  1 2 /
  \plot 0 3  1 3 /
  \put{$\sssize\bullet$} at 0.5 3.5 
  \endpicture} at 0 4
\put{\beginpicture
   \setcoordinatesystem units <.14cm,.14cm>
  \multiput{} at 0 0  1 4 /
  \setsolid
  \plot 0 0  1 0  1 4  0 4  0 0 /
  \plot 0 1  1 1 /
  \plot 0 2  1 2 /
  \plot 0 3  1 3 /
  \put{$\sssize\bullet$} at 0.5 2.5 
  \endpicture} at 1 3
\put{\beginpicture
   \setcoordinatesystem units <.14cm,.14cm>
  \multiput{} at 0 0  1 4 /
  \setsolid
  \plot 0 0  1 0  1 4  0 4  0 0 /
  \plot 0 1  1 1 /
  \plot 0 2  1 2 /
  \plot 0 3  1 3 /
  \put{$\sssize\bullet$} at 0.5 1.5 
  \endpicture} at 2 2
\put{\beginpicture
   \setcoordinatesystem units <.14cm,.14cm>
  \multiput{} at 0 0  1 4 /
  \setsolid
  \plot 0 0  1 0  1 4  0 4  0 0 /
  \plot 0 1  1 1 /
  \plot 0 2  1 2 /
  \plot 0 3  1 3 /
  \put{$\sssize\bullet$} at 0.5 .5 
  \endpicture} at 3 1
\put{\beginpicture
   \setcoordinatesystem units <.14cm,.14cm>
  \multiput{} at 0 0  1 4 /
  \setsolid
  \plot 0 0  1 0  1 4  0 4  0 0 /
  \plot 0 1  1 1 /
  \plot 0 2  1 2 /
  \plot 0 3  1 3 /
  \endpicture} at 4 0 
\put{\beginpicture
   \setcoordinatesystem units <.14cm,.14cm>
  \multiput{} at 0 0  2 4 /
  \setsolid
  \plot 0 0  1 0  1 4  0 4  0 0 /
  \plot 0 1  2 1  2 2  0 2 /
  \plot 0 3  1 3 /
  \multiput{$\sssize\bullet$} at 0.5 1.5  1.5 1.5 /
  \plot 0.5 1.5  1.5 1.5 /
  \plot 0.5 1.45  1.5 1.45 /
  \plot 0.5 1.55  1.5 1.55 /
  \endpicture} at 0 1 
\put{\beginpicture
   \setcoordinatesystem units <.14cm,.14cm>
  \multiput{} at 0 0  2 4 /
  \setsolid
  \plot 0 0  1 0  1 4  0 4  0 0 /
  \plot 0 1  2 1  2 3  0 3 /
  \plot 0 2  2 2 /
  \multiput{$\sssize\bullet$} at 0.5 1.5  1.5 1.5 /
  \plot 0.5 1.5  1.5 1.5 /
  \plot 0.5 1.45  1.5 1.45 /
  \plot 0.5 1.55  1.5 1.55 /
  \endpicture} at 1.4 1 
\put{\beginpicture
   \setcoordinatesystem units <.14cm,.14cm>
  \multiput{} at 0 0  2 4 /
  \setsolid
  \plot 0 0  1 0  1 4  0 4  0 0 /
  \plot 0 2  2 2  2 3  0 3 /
  \plot 0 1  1 1 /
  \multiput{$\sssize\bullet$} at 0.5 2.5  1.5 2.5 /
  \plot 0.5 2.5  1.5  2.5 /
  \plot 0.5 2.45  1.5 2.45 /
  \plot 0.5 2.55  1.5 2.55 /
  \endpicture} at -.5 1.5 
\put{\beginpicture
   \setcoordinatesystem units <.14cm,.14cm>
  \multiput{} at 0 0  2 4 /
  \setsolid
  \plot 0 0  1 0  1 4  0 4  0 0 /
  \plot 0 1  2 1  2 3  0 3 /
  \plot 0 2  2 2 /
  \multiput{$\sssize\bullet$} at 0.5 2.5  1.5 2.5 /
  \plot 0.5 2.5  1.5  2.5 /
  \plot 0.5 2.45  1.5 2.45 /
  \plot 0.5 2.55  1.5 2.55 /
  \endpicture} at .5 1.5 
\put{\beginpicture
   \setcoordinatesystem units <.14cm,.14cm>
  \multiput{} at 0 0  2 4 /
  \setsolid
  \plot 0 0  1 0  1 4  0 4  0 0 /
  \plot 0 1  2 1  2 3  0 3 /
  \plot 0 2  2 2 /
  \multiput{$\sssize\bullet$} at 0.5 2.5  1.5 2.5  1.5 1.5 /
  \plot 0.5 2.5  1.5  2.5 /
  \plot 0.5 2.45  1.5 2.45 /
  \plot 0.5 2.55  1.5 2.55 /
  \endpicture} at .25 2.15
\endpicture} at 11 0 
\endpicture}
$$

{\bf The case $n = 5$.}  The $\tau_n$-orbits of the pickets provide 26 indecomposable objects.
There are 24 additional objects, they 
form four $\tau_n$-orbits of cardinality $6$. Three of these orbits contain objects of the form
$C_{\lambda}$, with $\lambda$ a partition (see Sections~\ref{sec-seven-two} and \ref{sec-thirteen-five}):
One orbit contains $C_{(4,1)}$ and  $C_{(5,2)}$,
the second contains $C_{(4,2)}$ and $C_{(5,3,1)}$, the third contains $C_{(3,1)}$ and 
$C_{(5,3)}$. The fourth is the branching orbit, it is the orbit of the object $X$ with
$\bpar X = ([3,1],[5,3,1],[3,2]).$

$$
{\beginpicture
   \setcoordinatesystem units <.57735cm,1cm>
\put{\beginpicture
\put{$\mathbb T(5)\:$} at -4.5 4.8
\multiput{} at -5 0  5 0  0 5 /
\plot -5 0  5 0  0 5  -5 0 /
\setdots <1mm>
\plot -3 0  -4 1  4 1  3 0  -1 4  1 4  -3 0 /
\plot -1 0  -3 2  3 2  1 0  -2 3  2 3  -1 0 /
\multiput{$\bullet$} at -3 0  -1 0  1 0  3 0  5 0  -4 1  -2 1  0 .95  2 1  4 1  -3 2  -1 1.95  1 1.95  3 2  -2 3  0 3  2 3  -1 4  1 4  0 5  /
\multiput{$\bigcirc$} at -2 1  -1 1  1 1  2 1  -1.5 1.5  1.5 1.5  -.5 2.5  .5 2.5  0 3 /
\setshadegrid span <.3mm>
\vshade -1 2 2  <,z,,> 0 1 2 <z,,,> 1 2 2 /
\setsolid 
\arr{4 -.3}{5 -.3}
\put{$r$} at 5.5 -.3
\arr{-1.6 4.5}{-1.1 5}
\put{$p$} at -.8 5.3
\endpicture} at 0 0 
\put{\beginpicture
\multiput{} at -5 0  5 0  0 5 /
\setdots <1mm>
\plot -5 0  5 0  0 5  -5 0 /
\plot -3 0  -4 1  4 1  3 0  -1 4  1 4  -3 0 /
\plot -1 0  -3 2  3 2  1 0  -2 3  2 3  -1 0 /
  \setshadegrid span <.3mm>
  \vshade -1 2 2  <,z,,> 0 1 2 <z,,,> 1 2 2 /
\put{\beginpicture
   \setcoordinatesystem units <.14cm,.14cm>
  \multiput{} at 0 0  1 1 /
  \setsolid
  \plot 0 0  1 0  1 1  0 1  0 0 /
  \endpicture} at -3 0
\put{\beginpicture
   \setcoordinatesystem units <.14cm,.14cm>
  \multiput{} at 0 0  1 2 /
  \setsolid
  \plot 0 0  1 0  1 2  0 2  0 0 /
  \plot 0 1  1 1 /
  \endpicture} at -1 0
\put{\beginpicture
   \setcoordinatesystem units <.14cm,.14cm>
  \multiput{} at 0 0  1 3 /
  \setsolid
  \plot 0 0  1 0  1 3  0 3  0 0 /
  \plot 0 1  1 1 /
  \plot 0 2  1 2 /
  \endpicture} at 1 0

\put{\beginpicture
   \setcoordinatesystem units <.14cm,.14cm>
  \multiput{} at 0 0  1 1 /
  \setsolid
  \plot 0 0  1 0  1 1  0 1  0 0 /
  \put{$\sssize\bullet$} at 0.5 0.5 
  \endpicture} at -4 1
\put{\beginpicture
   \setcoordinatesystem units <.14cm,.14cm>
  \multiput{} at 0 0  1 2 /
  \setsolid
  \plot 0 0  1 0  1 2  0 2  0 0 /
  \plot 0 1  1 1 /
  \put{$\sssize\bullet$} at 0.5 1.5 
  \endpicture} at -3 2
\put{\beginpicture
   \setcoordinatesystem units <.14cm,.14cm>
  \multiput{} at 0 0  1 3 /
  \setsolid
  \plot 0 0  1 0  1 3  0 3  0 0 /
  \plot 0 1  1 1 /
  \plot 0 2  1 2 /
  \put{$\sssize\bullet$} at 0.5 2.5 
  \endpicture} at -2 3

\put{\beginpicture
   \setcoordinatesystem units <.14cm,.14cm>
  \multiput{} at 0 0  1 2 /
  \setsolid
  \plot 0 0  1 0  1 2  0 2  0 0 /
  \plot 0 1  1 1 /
  \put{$\sssize\bullet$} at 0.5 .5 
  \endpicture} at -2.4 1

\put{\beginpicture
   \setcoordinatesystem units <.14cm,.14cm>
  \multiput{} at 0 0  2 3 /
  \setsolid
  \plot 0 0  1 0  1 3  0 3  0 0 /
  \plot 0 1  2 1  2 2  0 2 /
  \multiput{$\sssize\bullet$} at 0.5 1.5  1.5 1.5 /
  \plot 0.5 1.5  1.5 1.5 /
  \plot 0.5 1.45  1.5 1.45 /
  \plot 0.5 1.55  1.5 1.55 /
  \endpicture} at -1.9 1 
\put{\beginpicture
   \setcoordinatesystem units <.14cm,.14cm>
  \multiput{} at 0 0  1 4 /
  \setsolid
  \plot 0 0  1 0  1 4  0 4  0 0 /
  \plot 0 1  1 1 /
  \plot 0 2  1 2 /
  \plot 0 3  1 3 /
  \put{$\sssize\bullet$} at 0.5 3.5 
  \endpicture} at -1 4
\put{\beginpicture
   \setcoordinatesystem units <.14cm,.14cm>
  \multiput{} at 0 0  1 4 /
  \setsolid
  \plot 0 0  1 0  1 4  0 4  0 0 /
  \plot 0 1  1 1 /
  \plot 0 2  1 2 /
  \plot 0 3  1 3 /
  \put{$\sssize\bullet$} at 0.5 2.5 
  \endpicture} at -.2 3.1
\put{\beginpicture
   \setcoordinatesystem units <.14cm,.14cm>
  \multiput{} at 0 0  1 4 /
  \setsolid
  \plot 0 0  1 0  1 4  0 4  0 0 /
  \plot 0 1  1 1 /
  \plot 0 2  1 2 /
  \plot 0 3  1 3 /
  \put{$\sssize\bullet$} at 0.5 .5 
  \endpicture} at 1.85 .95
\put{\beginpicture
   \setcoordinatesystem units <.14cm,.14cm>
  \multiput{} at 0 0  1 4 /
  \setsolid
  \plot 0 0  1 0  1 4  0 4  0 0 /
  \plot 0 1  1 1 /
  \plot 0 2  1 2 /
  \plot 0 3  1 3 /
  \endpicture} at 3 0 
\put{\beginpicture
   \setcoordinatesystem units <.14cm,.14cm>
  \multiput{} at 0 0  2 4 /
  \setsolid
  \plot 0 0  1 0  1 4  0 4  0 0 /
  \plot 0 1  2 1  2 2  0 2 /
  \plot 0 3  1 3 /
  \multiput{$\sssize\bullet$} at 0.5 1.5  1.5 1.5 /
  \plot 0.5 1.5  1.5 1.5 /
  \plot 0.5 1.45  1.5 1.45 /
  \plot 0.5 1.55  1.5 1.55 /
  \endpicture} at -1 .95 
\put{\beginpicture
   \setcoordinatesystem units <.14cm,.14cm>
  \multiput{} at 0 0  2 4 /
  \setsolid
  \plot 0 0  1 0  1 4  0 4  0 0 /
  \plot 0 2  2 2  2 3  0 3 /
  \plot 0 1  1 1 /
  \multiput{$\sssize\bullet$} at 0.5 2.5  1.5 2.5 /
  \plot 0.5 2.5  1.5  2.5 /
  \plot 0.5 2.45  1.5 2.45 /
  \plot 0.5 2.55  1.5 2.55 /
  \endpicture} at -1.5 1.55 
\put{\beginpicture
   \setcoordinatesystem units <.14cm,.14cm>
  \multiput{} at 0 0  1 5 /
  \setsolid
  \plot 0 0  1 0  1 5  0 5  0 0 /
  \plot 0 1  1 1 /
  \plot 0 2  1 2 /
  \plot 0 3  1 3 /
  \plot 0 4  1 4 /
  \put{$\sssize\bullet$} at 0.5 4.5 
  \endpicture} at 0 5
\put{\beginpicture
   \setcoordinatesystem units <.14cm,.14cm>
  \multiput{} at 0 0  1 5 /
  \setsolid
  \plot 0 0  1 0  1 5  0 5  0 0 /
  \plot 0 1  1 1 /
  \plot 0 2  1 2 /
  \plot 0 3  1 3 /
  \plot 0 4  1 4 /
  \put{$\sssize\bullet$} at 0.5 3.5 
  \endpicture} at 1 4 
\put{\beginpicture
   \setcoordinatesystem units <.14cm,.14cm>
  \multiput{} at 0 0  1 5 /
  \setsolid
  \plot 0 0  1 0  1 5  0 5  0 0 /
  \plot 0 1  1 1 /
  \plot 0 2  1 2 /
  \plot 0 3  1 3 /
  \plot 0 4  1 4 /
  \put{$\sssize\bullet$} at 0.5 2.5 
  \endpicture} at 2 3 
\put{\beginpicture
   \setcoordinatesystem units <.14cm,.14cm>
  \multiput{} at 0 0  1 5 /
  \setsolid
  \plot 0 0  1 0  1 5  0 5  0 0 /
  \plot 0 1  1 1 /
  \plot 0 2  1 2 /
  \plot 0 3  1 3 /
  \plot 0 4  1 4 /
  \put{$\sssize\bullet$} at 0.5 1.5 
  \endpicture} at 3 2
\put{\beginpicture
   \setcoordinatesystem units <.14cm,.14cm>
  \multiput{} at 0 0  1 5 /
  \setsolid
  \plot 0 0  1 0  1 5  0 5  0 0 /
  \plot 0 1  1 1 /
  \plot 0 2  1 2 /
  \plot 0 3  1 3 /
  \plot 0 4  1 4 /
  \put{$\sssize\bullet$} at 0.5 0.5 
  \endpicture} at 4 1
\put{\beginpicture
   \setcoordinatesystem units <.14cm,.14cm>
  \multiput{} at 0 0  1 5 /
  \setsolid
  \plot 0 0  1 0  1 5  0 5  0 0 /
  \plot 0 1  1 1 /
  \plot 0 2  1 2 /
  \plot 0 3  1 3 /
  \plot 0 4  1 4 /
  \endpicture} at 5 0 
\put{\beginpicture
   \setcoordinatesystem units <.14cm,.14cm>
  \multiput{} at 0 0  2 5 /
  \setsolid
  \plot 0 0  1 0  1 5  0 5  0 0 /
  \plot 0 2  2 2  2 4  0 4 /
  \plot 0 3  2 3 /
  \plot 0 1  1 1 /
  \multiput{$\sssize\bullet$} at 0.5 3.5  1.5 3.5  1.5 2.5 /
  \plot 0.5 3.5  1.5  3.5 /
  \plot 0.5 3.45  1.5 3.45 /
  \plot 0.5 3.55  1.5 3.55 /
  \endpicture} at -.6 2.45 
\put{\beginpicture
   \setcoordinatesystem units <.14cm,.14cm>
  \multiput{} at 0 0  2 5 /
  \setsolid
  \plot 0 0  1 0  1 5  0 5  0 0 /
  \plot 0 1  2 1  2 4  0 4 /
  \plot 0 3  2 3 /
  \plot 0 2  2 2 /
  \multiput{$\sssize\bullet$} at 0.5 3.5  1.5 3.5  1.5 2.5 /
  \plot 0.5 3.5  1.5  3.5 /
  \plot 0.5 3.45  1.5 3.45 /
  \plot 0.5 3.55  1.5 3.55 /
  \endpicture} at .25 3.1
\put{\beginpicture
   \setcoordinatesystem units <.14cm,.14cm>
  \multiput{} at 0 0  2 5 /
  \setsolid
  \plot 0 0  1 0  1 5  0 5  0 0 /
  \plot 0 1  2 1  2 4  0 4 /
  \plot 0 3  2 3 /
  \plot 0 2  2 2 /
  \multiput{$\sssize\bullet$} at 0.5 3.5  1.5 3.5  1.5 1.5 /
  \plot 0.5 3.5  1.5  3.5 /
  \plot 0.5 3.45  1.5 3.45 /
  \plot 0.5 3.55  1.5 3.55 /
  \endpicture} at .6 2.45
\put{\beginpicture
   \setcoordinatesystem units <.14cm,.14cm>
  \multiput{} at 0 0  2 5 /
  \setsolid
  \plot 0 0  1 0  1 5  0 5  0 0 /
  \plot 0 1  2 1  2 4  0 4 /
  \plot 0 3  2 3 /
  \plot 0 2  2 2 /
  \multiput{$\sssize\bullet$} at 0.5 2.5  1.5 2.5  /
  \plot 0.5 2.5  1.5  2.5 /
  \plot 0.5 2.45  1.5 2.45 /
  \plot 0.5 2.55  1.5 2.55 /
  \endpicture} at 1.5 1.55
\put{\beginpicture
   \setcoordinatesystem units <.14cm,.14cm>
  \multiput{} at 0 0  2 5 /
  \setsolid
  \plot 0 0  1 0  1 5  0 5  0 0 /
  \plot 0 1  2 1  2 4  0 4 /
  \plot 0 3  2 3 /
  \plot 0 2  2 2 /
  \multiput{$\sssize\bullet$} at 0.5 1.5  1.5 1.5 /
  \plot 0.5 1.5  1.5  1.5 /
  \plot 0.5 1.45  1.5 1.45 /
  \plot 0.5 1.55  1.5 1.55 /
  \endpicture} at 2.35 .95
\put{\beginpicture
   \setcoordinatesystem units <.14cm,.14cm>
  \multiput{} at 0 0  2 5 /
  \setsolid
  \plot 0 0  1 0  1 5  0 5  0 0 /
  \plot 0 1  2 1  2 3  0 3 /
  \plot 0 4  1 4 /
  \plot 0 2  2 2 /
  \multiput{$\sssize\bullet$} at 0.5 1.5  1.5 1.5 /
  \plot 0.5 1.5  1.5  1.5 /
  \plot 0.5 1.45  1.5 1.45 /
  \plot 0.5 1.55  1.5 1.55 /
  \endpicture} at .95 .95
\endpicture} at 9 -3 
\endpicture}
$$

$$
{\beginpicture
   \setcoordinatesystem units <.8083cm,1.4cm>
\put{\beginpicture
\put{$\mathbb T(5)$} at -3.5 .5
\put{(shaded part)} at -3.5 0
\multiput{} at -3 3  3 3  0 0 /
\setdots <.5mm>
\plot -3 3  3 3  0 0  -3 3 /
\setdots <1mm>
\plot -1 1  1 1  -1 3  -2 2  2 2  1 3  -1 1 /
\plot -1.5 1.5  1.5 1.5  0 3  -1.5 1.5 /
\plot -.5 .5  .5 .5  -2 3  -2.5 2.5  2.5 2.5  2 3  -.5 .5 /
\multiput{$\bullet$} at -3 3  3 3  0 0 /
\multiput{$\bigcirc$} at  0 0  -3 3  3 3 
            -1 1  -1.5 1.5  -2 2  -1 3  0 3  1 3  2 2  1.5 1.5  1 1 /
\put{$\ssize (p,r)=(1,2)$} at 0 -.3
\put{$\ssize (2,1)$} at -3.5 2.8
\put{$\ssize (2,2)$} at 3.5 2.8
\multiput{$\sssize 1223$} at -3 3.2  3 3.2  0.6 0 /
\multiput{$\sssize 22$} at 0 3.2  -1.9 1.5  1.9 1.5 /
\multiput{$\sssize 3$} at -1 3.2  1 3.2  -2.4 2  -1.4 1  2.4 2  1.4 1 /
\endpicture} at 0 0 
\put{\beginpicture
\multiput{} at -3 3  3 3  0 0 /
\setdots <.5mm>
\plot .3 .3  .6 .6 /
\plot 1.1 1.1  1.3 1.3 /
\plot 1.75 1.75  1.85 1.85 /
\plot 2.45 2.45  2.7 2.7 /
\plot -.3 .3  -.65 .65 /
\plot -1.1 1.1  -1.27 1.27 /
\plot -1.7 1.7  -1.9 1.9 /
\plot -2.4 2.4  -2.7 2.7 /
\plot -1.8 3  -1.5 3 /
\plot -.75 3  -.5 3 /
\plot 1.6 3  2.1 3 /
\plot .8 3  .55 3 /
\put{\beginpicture
   \setcoordinatesystem units <.14cm,.14cm>
  \multiput{} at 0 0  2 5 /
  \setsolid
  \plot 0 0  1 0  1 5  0 5  0 0 /
  \plot 0 4  2 4  2 1  0 1 /
  \plot 0 3  2 3 /
  \plot 0 2  2 2 /
  \multiput{$\sssize\bullet$} at 0.5 3.5  1.5 3.5 /
  \plot 0.5 3.5  1.5  3.5 /
  \plot 0.5 3.45  1.5 3.45 /
  \plot 0.5 3.55  1.5 3.55 /
  \endpicture} at 3.2 3 
\put{\beginpicture
   \setcoordinatesystem units <.14cm,.14cm>
  \multiput{} at 0 0  2 5 /
  \setsolid
  \plot 0 0  1 0  1 5  0 5  0 0 /
  \plot 0 1  2 1  2 4  0 4 /
  \plot 0 2  2 2 /
  \plot 0 3  2 3 /
  \multiput{$\sssize\bullet$} at 0.5 2.5  1.5 2.5  1.5 1.5 /
  \plot 0.5 2.5  1.5  2.5 /
  \plot 0.5 2.45  1.5 2.45 /
  \plot 0.5 2.55  1.5 2.55 /
  \endpicture} at 2.75 3 
\put{\beginpicture
   \setcoordinatesystem units <.14cm,.14cm>
  \multiput{} at 0 0  1 4 /
  \setsolid
  \plot 0 0  1 0  1 4  0 4  0 0 /
  \plot 0 1  1 1 /
  \plot 0 2  1 2 /
  \plot 0 3  1 3 /
  \put{$\sssize\bullet$} at 0.5 1.5 
  \endpicture} at 2.35 3
\put{\beginpicture
   \setcoordinatesystem units <.14cm,.14cm>
  \multiput{} at 0 0  3 6 /
  \setsolid
  \plot 0 0  1 0  1 5  0 5  0 0 /
  \plot 2 1  3 1  3 6  2 6  2 1 /
  \plot 2 5  3 5 /
  \plot 0 4  3 4 /
  \plot 0 3  3 3 /
  \plot 0 2  3 2 /
  \plot 0 1  1 1 /
  \multiput{$\sssize\bullet$} at 0.5 3.5  1.5 3.5  1.5 2.5  2.5 2.5 /
  \plot 0.5 3.5  1.5  3.5 /
  \plot 0.5 3.45  1.5 3.45 /
  \plot 0.5 3.55  1.5 3.55 /
  \plot 1.5 2.5   2.5 2.5 /
  \plot 1.5 2.45  2.5 2.45 /
  \plot 1.5 2.55  2.5 2.55 /
  \endpicture} at 3.75 3  
\put{\beginpicture
   \setcoordinatesystem units <.14cm,.14cm>
  \multiput{} at 0 1  3 6 /
  \setsolid
  \plot 0 1  1 1  1 5  0 5  0 1 /
  \plot 2 1  3 1  3 6  2 6  2 1 /
  \plot 2 5  3 5 /
  \plot 0 4  3 4 /
  \plot 0 3  3 3 /
  \plot 0 2  3 2 /
  \multiput{$\sssize\bullet$} at 0.5 3.5  1.5 3.5  1.5 2.5  2.5 2.5 /
  \plot 0.5 3.5  1.5  3.5 /
  \plot 0.5 3.45  1.5 3.45 /
  \plot 0.5 3.55  1.5 3.55 /
  \plot 1.5 2.5   2.5 2.5 /
  \plot 1.5 2.45  2.5 2.45 /
  \plot 1.5 2.55  2.5 2.55 /
  \endpicture} at 2.15 2.15 
\put{\beginpicture
   \setcoordinatesystem units <.14cm,.14cm>
  \multiput{} at 0 0  3 5 /
  \setsolid
  \plot 0 0  1 0  1 5  0 5  0 0 /
  \plot 2 1  3 1  3 5  2 5  2 1 /
  \plot 0 4  3 4 /
  \plot 0 3  3 3 /
  \plot 0 2  3 2 /
  \plot 0 1  1 1 /
  \multiput{$\sssize\bullet$} at 0.5 3.5  1.5 3.5  1.5 2.5  2.5 2.5 /
  \plot 0.5 3.5  1.5  3.5 /
  \plot 0.5 3.45  1.5 3.45 /
  \plot 0.5 3.55  1.5 3.55 /
  \plot 1.5 2.5   2.5 2.5 /
  \plot 1.5 2.45  2.5 2.45 /
  \plot 1.5 2.55  2.5 2.55 /
  \endpicture} at 1.2 3  
\put{\beginpicture
   \setcoordinatesystem units <.14cm,.14cm>
  \multiput{} at 0 0  2 5 /
  \setsolid
  \plot 0 0  1 0  1 5  0 5  0 0 /
  \plot 0 4  2 4  2 2  0 2 /
  \plot 0 3  2 3 /
  \plot 0 1  1 1 /
  \multiput{$\sssize\bullet$} at 0.5 3.5  1.5 3.5 /
  \plot 0.5 3.5  1.5  3.5 /
  \plot 0.5 3.45  1.5 3.45 /
  \plot 0.5 3.55  1.5 3.55 /
  \endpicture} at 0.25 3 
\put{\beginpicture
   \setcoordinatesystem units <.14cm,.14cm>
  \multiput{} at 0 0  2 5 /
  \setsolid
  \plot 0 0  1 0  1 5  0 5  0 0 /
  \plot 0 1  2 1  2 3  0 3 /
  \plot 0 2  2 2 /
  \plot 0 4  1 4 /
  \multiput{$\sssize\bullet$} at 0.5 2.5  1.5 2.5  1.5 1.5 /
  \plot 0.5 2.5  1.5  2.5 /
  \plot 0.5 2.45  1.5 2.45 /
  \plot 0.5 2.55  1.5 2.55 /
  \endpicture} at -.2 3 
\put{\beginpicture
   \setcoordinatesystem units <.14cm,.14cm>
  \multiput{} at 0 0  2 5 /
  \setsolid
  \plot 0 0  1 0  1 5  0 5  0 0 /
  \plot 0 4  2 4  2 3  0 3 /
  \plot 0 2  1 2 /
  \plot 0 1  1 1 /
  \multiput{$\sssize\bullet$} at 0.5 3.5  1.5 3.5 /
  \plot 0.5 3.5  1.5  3.5 /
  \plot 0.5 3.45  1.5 3.45 /
  \plot 0.5 3.55  1.5 3.55 /
  \endpicture} at -2.75 3
\put{\beginpicture
   \setcoordinatesystem units <.14cm,.14cm>
  \multiput{} at 0 0  2 4 /
  \setsolid
  \plot 0 0  1 0  1 4  0 4  0 0 /
  \plot 0 1  2 1  2 3  0 3 /
  \plot 0 2  2 2 /
  \multiput{$\sssize\bullet$} at 0.5 2.5  1.5 2.5  1.5 1.5 /
  \plot 0.5 2.5  1.5  2.5 /
  \plot 0.5 2.45  1.5 2.45 /
  \plot 0.5 2.55  1.5 2.55 /
  \endpicture} at -3.2 3 
\put{\beginpicture
   \setcoordinatesystem units <.14cm,.14cm>
  \multiput{} at 0 0  1 3 /
  \setsolid
  \plot 0 0  1 0  1 3  0 3  0 0 /
  \plot 0 1  1 1 /
  \plot 0 2  1 2 /
  \put{$\sssize\bullet$} at 0.5 1.5 
  \endpicture} at -3.6 3
\put{\beginpicture
   \setcoordinatesystem units <.14cm,.14cm>
  \multiput{} at 0 0  3 5 /
  \setsolid
  \plot 0 0  1 0  1 5  0 5  0 0 /
  \plot 0 1  2 1  2 4  0 4 /
  \plot 0 2  3 2  3 3  0 3 /
  \multiput{$\sssize\bullet$} at 0.5 3.5  1.5 3.5  1.5 2.5  2.5 2.5 /
  \plot 0.5 3.5   1.5 3.5 /
  \plot 0.5 3.45  1.5 3.45 /
  \plot 0.5 3.55  1.5 3.55 /
  \plot 1.5 2.5   2.5 2.5 /
  \plot 1.5 2.45  2.5 2.45 /
  \plot 1.5 2.55  2.5 2.55 /
  \endpicture} at -2.2 3 
\put{\beginpicture
   \setcoordinatesystem units <.14cm,.14cm>
  \multiput{} at 0 0  3 5 /
  \setsolid
  \plot 0 0  1 0  1 5  0 5  0 0 /
  \plot 0 1  2 1  2 4  0 4 /
  \plot 0 2  3 2  3 3  0 3 /
  \multiput{$\sssize\bullet$} at 0.3 2.7  1.3 2.7  1.7 2.3  2.7 2.3 /
  \plot 0.3 2.7   1.3 2.7 /
  \plot 0.3 2.65  1.3 2.65 /
  \plot 0.3 2.75  1.3 2.75 /
  \plot 1.7 2.3   2.7 2.3 /
  \plot 1.7 2.25  2.7 2.25 /
  \plot 1.7 2.35  2.7 2.35 /
  \endpicture} at -2.1 2.1  
\put{\beginpicture
   \setcoordinatesystem units <.14cm,.14cm>
  \multiput{} at 0 0  3 5 /
  \setsolid
  \plot 0 0  1 0  1 5  0 5  0 0 /
  \plot 0 1  2 1  2 4  0 4 /
  \plot 0 2  3 2  3 4  0 4 /
  \plot 0 3  3 3 /
  \multiput{$\sssize\bullet$} at 0.5 3.5  1.5 3.5  1.5 2.5  2.5 2.5 /
  \plot 0.5 3.5   1.5 3.5 /
  \plot 0.5 3.45  1.5 3.45 /
  \plot 0.5 3.55  1.5 3.55 /
  \plot 1.5 2.5   2.5 2.5 /
  \plot 1.5 2.45  2.5 2.45 /
  \plot 1.5 2.55  2.5 2.55 /
  \endpicture} at -1.1 3 
\put{\beginpicture
   \setcoordinatesystem units <.14cm,.14cm>
  \multiput{} at 0 0  2 5 /
  \setsolid
  \plot 0 0  1 0  1 5  0 5  0 0 /
  \plot 0 2  2 2  2 3  0 3 /
  \plot 0 4  1 4 /
  \plot 0 1  1 1 /
  \multiput{$\sssize\bullet$} at 0.5 2.5  1.5 2.5 /
  \plot 0.5 2.5  1.5  2.5 /
  \plot 0.5 2.45  1.5 2.45 /
  \plot 0.5 2.55  1.5 2.55 /
  \endpicture} at -1.3 1.55 
\put{\beginpicture
   \setcoordinatesystem units <.14cm,.14cm>
  \multiput{} at 0 0  2 4 /
  \setsolid
  \plot 0 0  1 0  1 4  0 4  0 0 /
  \plot 0 1  2 1  2 3  0 3 /
  \plot 0 2  2 2 /
  \multiput{$\sssize\bullet$} at 0.5 2.5  1.5 2.5 /
  \plot 0.5 2.5  1.5  2.5 /
  \plot 0.5 2.45  1.5 2.45 /
  \plot 0.5 2.55  1.5 2.55 /
  \endpicture} at -1.75 1.45 
\put{\beginpicture
   \setcoordinatesystem units <.14cm,.14cm>
  \multiput{} at 0 0  2 5 /
  \setsolid
  \plot 0 0  1 0  1 5  0 5  0 0 /
  \plot 0 1  2 1  2 2  0 2 /
  \plot 0 4  1 4 /
  \plot 0 3  1 3 /
  \multiput{$\sssize\bullet$} at 0.5 1.5  1.5 1.5 /
  \plot 0.5 1.5  1.5  1.5 /
  \plot 0.5 1.45  1.5 1.45 /
  \plot 0.5 1.55  1.5 1.55 /
  \endpicture} at .2 0 
\put{\beginpicture
   \setcoordinatesystem units <.14cm,.14cm>
  \multiput{} at 0 0  2 4 /
  \setsolid
  \plot 0 0  1 0  1 4  0 4  0 0 /
  \plot 0 1  2 1  2 3  0 3 /
  \plot 0 2  2 2 /
  \multiput{$\sssize\bullet$} at 0.5 1.5  1.5 1.5 /
  \plot 0.5 1.5  1.5  1.5 /
  \plot 0.5 1.45  1.5 1.45 /
  \plot 0.5 1.55  1.5 1.55 /
  \endpicture} at -.25 0 
\put{\beginpicture
   \setcoordinatesystem units <.14cm,.14cm>
  \multiput{} at 0 0  3 5 /
  \setsolid
  \plot 0 0  1 0  1 5  0 5  0 0 /
  \plot 0 1  2 1  2 4  0 4 /
  \plot 0 2  3 2  3 3  0 3 /
  \multiput{$\sssize\bullet$} at 0.5 2.5  1.5 2.5  2.5 2.5 /
  \plot 0.5 2.5   2.5 2.5 /
  \plot 0.5 2.45  2.5 2.45 /
  \plot 0.5 2.55  2.5 2.55 /
  \endpicture} at .75 0 
\put{\beginpicture
   \setcoordinatesystem units <.14cm,.14cm>
  \multiput{} at 0 0  3 5 /
  \setsolid
  \plot 0 0  1 0  1 5  0 5  0 0 /
  \plot 0 1  2 1  2 4  0 4 /
  \plot 0 2  3 2  3 3  0 3 /
  \multiput{$\sssize\bullet$} at 0.5 2.5  1.5 2.5  2.5 2.5  1.5 1.5 /
  \plot 0.5 2.5   2.5 2.5 /
  \plot 0.5 2.45  2.5 2.45 /
  \plot 0.5 2.55  2.5 2.55 /
  \endpicture} at -.85 .85
\put{\beginpicture
   \setcoordinatesystem units <.14cm,.14cm>
  \multiput{} at 0 0  3 5 /
  \setsolid
  \plot 0 0  1 0  1 5  0 5  0 0 /
  \plot 0 1  2 1  2 4  0 4 /
  \plot 0 1  3 1  3 3  0 3 /
  \plot 0 2  3 2 /
  \multiput{$\sssize\bullet$} at 0.5 2.5  1.5 2.5  2.5 2.5  2.5 1.5 /
  \plot 0.5 2.5   2.5 2.5 /
  \plot 0.5 2.45  2.5 2.45 /
  \plot 0.5 2.55  2.5 2.55 /
  \endpicture} at .9 .9
\put{\beginpicture
   \setcoordinatesystem units <.14cm,.14cm>
  \multiput{} at 0 0  1 3 /
  \setsolid
  \plot 0 0  1 0  1 3  0 3  0 0 /
  \plot 0 1  1 1 /
  \plot 0 2  1 2 /
  \put{$\sssize\bullet$} at 0.5 0.5 
  \endpicture} at -.65 0
\put{\beginpicture
   \setcoordinatesystem units <.14cm,.14cm>
  \multiput{} at 0 0  2 5 /
  \setsolid
  \plot 0 0  1 0  1 5  0 5  0 0 /
  \plot 0 1  2 1  2 3  0 3 /
  \plot 0 4  1 4 /
  \plot 0 2  2 2 /
  \multiput{$\sssize\bullet$} at 0.5 2.5  1.5 2.5 /
  \plot 0.5 2.5  1.5  2.5 /
  \plot 0.5 2.45  1.5 2.45 /
  \plot 0.5 2.55  1.5 2.55 /
  \endpicture} at 1.25 1.55 
\put{\beginpicture
   \setcoordinatesystem units <.14cm,.14cm>
  \multiput{} at 0 0  2 5 /
  \setsolid
  \plot 0 0  1 0  1 5  0 5  0 0 /
  \plot 0 2  2 2  2 4  0 4 /
  \plot 0 1  1 1 /
  \plot 0 3  2 3 /
  \multiput{$\sssize\bullet$} at 0.5 2.5  1.5 2.5 /
  \plot 0.5 2.5  1.5  2.5 /
  \plot 0.5 2.45  1.5 2.45 /
  \plot 0.5 2.55  1.5 2.55 /
  \endpicture} at 1.7 1.45 
\endpicture} at 7 -2
\endpicture}
$$

\subsection{Some indecomposables in $\Cal S(6).$}
\label{sec-fifteen-two}

The Auslander-Reiten quiver of $\Cal S(6)$ has been discussed in \cite{RS1}.
Some of the components,
in particular the principal component $\Cal P(6)$, are displayed in \cite{RS1}.
The remaining components which contain pickets
or bipickets are exhibited in \cite{S2}. 

We first sketch the pr-triangle for $\Cal S(6)$ and exhibit the objects outside the central hexagon.

$$
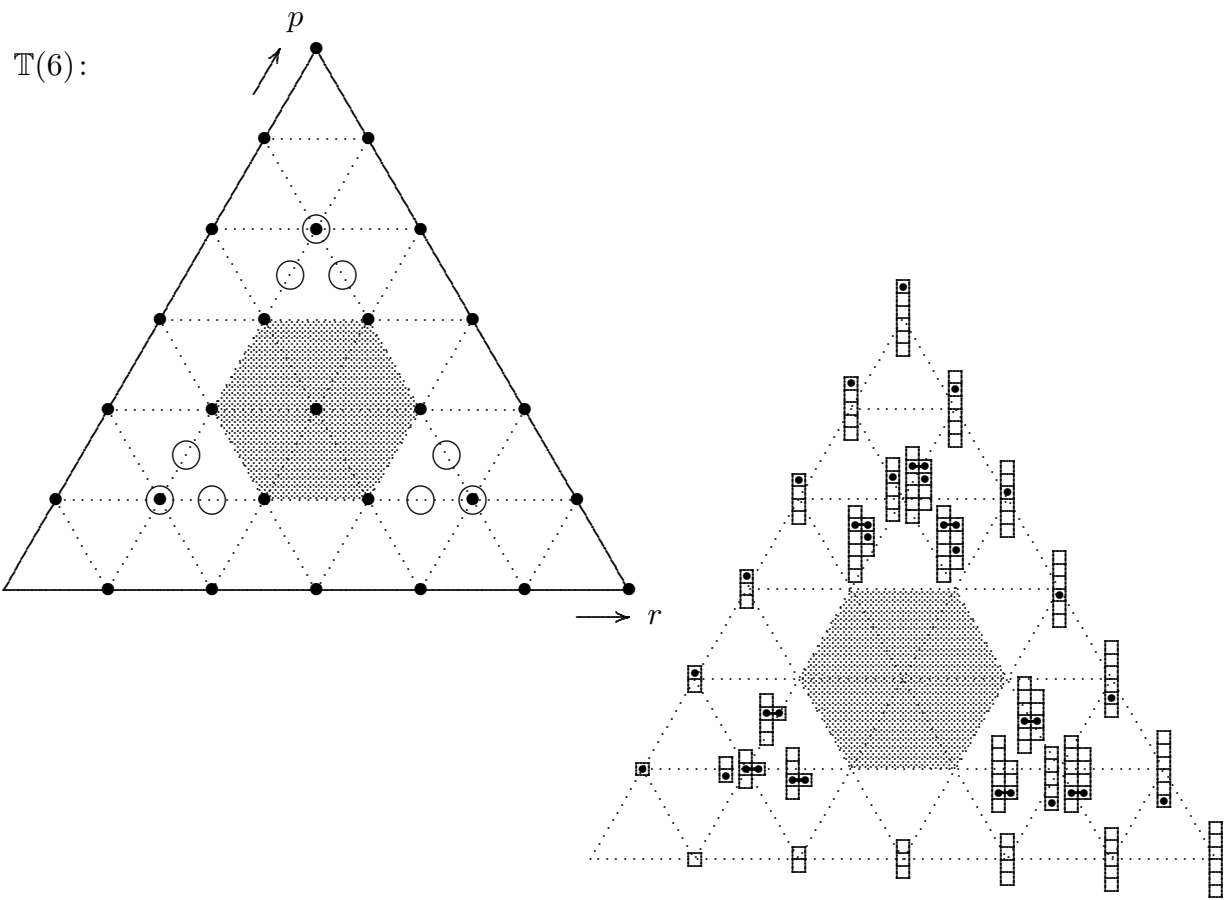
{\beginpicture
   \setcoordinatesystem units <.57735cm,1cm>
\put{\beginpicture
\put{$\mathbb T(6)\:$} at -5 5.8
\multiput{} at -6 0  6 0  0 6 /
\plot -6 0  6 0  0 6  -6 0 /
\setdots <1mm>
\plot -5 1  5 1  4 0  -1 5  1 5  -4 0  -5 1 /
\plot -4 2  4 2  2 0  -2 4  2 4  -2 0  -4 2 /
\plot 0 0  -3 3  3 3  0 0 /
\multiput{$\bullet$} at -4 0  -2 0  0 0  2 0  4 0  6 0  -5 1  -3 1  -1 1  1 1  3 1  5 1  -4 2  -2 2  0 2  2 2  4 2
                     -3 3  -1 3  1 3  3 3  -2 4  0 4  2 4  -1 5  1 5  0 6 /
\multiput{$\bigcirc$} at -3 1  -2 1  2 1  3 1  -2.5 1.5  2.5 1.5  -.5 3.5  .5 3.5  0 4 /
\setshadegrid span <.3mm>
\vshade -2 2 2  <,z,,> -1 1 3 <z,z,,> 1 1 3 <z,,,> 2 2 2 /
\setsolid 
\arr{5 -.3}{6 -.3}
\put{$r$} at 6.5 -.3
\arr{-1.2 5.5}{-.7 6}
\put{$p$} at -.4 6.3
\endpicture} at 0 0 
\put{\beginpicture
\multiput{} at -6 0  6 0  0 6 /
\setdots <1mm>
\plot -5 1  5 1  4 0  -1 5  1 5  -4 0  -5 1 /
\plot -4 2  4 2  2 0  -2 4  2 4  -2 0  -4 2 /
\plot 0 0  -3 3  3 3  0 0 /
\plot -6 0  0 6  6 0  -6 0 /
\setshadegrid span <.3mm>
\vshade -2 2 2  <,z,,> -1 1 3 <z,z,,> 1 1 3 <z,,,> 2 2 2 /
\put{\beginpicture
   \setcoordinatesystem units <.14cm,.14cm>
  \multiput{} at 0 0  1 1 /
  \setsolid
  \plot 0 0  1 0  1 1  0 1  0 0 /
  \endpicture} at -4 0
\put{\beginpicture
   \setcoordinatesystem units <.14cm,.14cm>
  \multiput{} at 0 0  1 2 /
  \setsolid
  \plot 0 0  1 0  1 2  0 2  0 0 /
  \plot 0 1  1 1 /
  \endpicture} at -2 0
\put{\beginpicture
   \setcoordinatesystem units <.14cm,.14cm>
  \multiput{} at 0 0  1 3 /
  \setsolid
  \plot 0 0  1 0  1 3  0 3  0 0 /
  \plot 0 1  1 1 /
  \plot 0 2  1 2 /
  \endpicture} at 0 0

\put{\beginpicture
   \setcoordinatesystem units <.14cm,.14cm>
  \multiput{} at 0 0  1 1 /
  \setsolid
  \plot 0 0  1 0  1 1  0 1  0 0 /
  \put{$\sssize\bullet$} at 0.5 0.5 
  \endpicture} at -5 1
\put{\beginpicture
   \setcoordinatesystem units <.14cm,.14cm>
  \multiput{} at 0 0  1 2 /
  \setsolid
  \plot 0 0  1 0  1 2  0 2  0 0 /
  \plot 0 1  1 1 /
  \put{$\sssize\bullet$} at 0.5 1.5 
  \endpicture} at -4 2
\put{\beginpicture
   \setcoordinatesystem units <.14cm,.14cm>
  \multiput{} at 0 0  1 3 /
  \setsolid
  \plot 0 0  1 0  1 3  0 3  0 0 /
  \plot 0 1  1 1 /
  \plot 0 2  1 2 /
  \put{$\sssize\bullet$} at 0.5 2.5 
  \endpicture} at -3 3

\put{\beginpicture
   \setcoordinatesystem units <.14cm,.14cm>
  \multiput{} at 0 0  1 2 /
  \setsolid
  \plot 0 0  1 0  1 2  0 2  0 0 /
  \plot 0 1  1 1 /
  \put{$\sssize\bullet$} at 0.5 .5 
  \endpicture} at -3.4 1

\put{\beginpicture
   \setcoordinatesystem units <.14cm,.14cm>
  \multiput{} at 0 0  2 3 /
  \setsolid
  \plot 0 0  1 0  1 3  0 3  0 0 /
  \plot 0 1  2 1  2 2  0 2 /
  \multiput{$\sssize\bullet$} at 0.5 1.5  1.5 1.5 /
  \plot 0.5 1.5  1.5 1.5 /
  \plot 0.5 1.45  1.5 1.45 /
  \plot 0.5 1.55  1.5 1.55 /
  \endpicture} at -2.9 1 
\put{\beginpicture
   \setcoordinatesystem units <.14cm,.14cm>
  \multiput{} at 0 0  1 4 /
  \setsolid
  \plot 0 0  1 0  1 4  0 4  0 0 /
  \plot 0 1  1 1 /
  \plot 0 2  1 2 /
  \plot 0 3  1 3 /
  \put{$\sssize\bullet$} at 0.5 3.5 
  \endpicture} at -2 4
\put{\beginpicture
   \setcoordinatesystem units <.14cm,.14cm>
  \multiput{} at 0 0  1 5 /
  \setsolid
  \plot 0 0  1 0  1 5  0 5  0 0 /
  \plot 0 1  1 1 /
  \plot 0 2  1 2 /
  \plot 0 3  1 3 /
  \plot 0 4  1 4 /
  \put{$\sssize\bullet$} at 0.5 3.5 
  \endpicture} at -.2 4.1
\put{\beginpicture
   \setcoordinatesystem units <.14cm,.14cm>
  \multiput{} at 0 0  1 5 /
  \setsolid
  \plot 0 0  1 0  1 5  0 5  0 0 /
  \plot 0 1  1 1 /
  \plot 0 2  1 2 /
  \plot 0 3  1 3 /
  \plot 0 4  1 4 /
  \put{$\sssize\bullet$} at 0.5 .5 
  \endpicture} at 2.85 .9
\put{\beginpicture
   \setcoordinatesystem units <.14cm,.14cm>
  \multiput{} at 0 0  1 4 /
  \setsolid
  \plot 0 0  1 0  1 4  0 4  0 0 /
  \plot 0 1  1 1 /
  \plot 0 2  1 2 /
  \plot 0 3  1 3 /
  \endpicture} at 2 0 
\put{\beginpicture
   \setcoordinatesystem units <.14cm,.14cm>
  \multiput{} at 0 0  2 4 /
  \setsolid
  \plot 0 0  1 0  1 4  0 4  0 0 /
  \plot 0 1  2 1  2 2  0 2 /
  \plot 0 3  1 3 /
  \multiput{$\sssize\bullet$} at 0.5 1.5  1.5 1.5 /
  \plot 0.5 1.5  1.5 1.5 /
  \plot 0.5 1.45  1.5 1.45 /
  \plot 0.5 1.55  1.5 1.55 /
  \endpicture} at -2 .95 
\put{\beginpicture
   \setcoordinatesystem units <.14cm,.14cm>
  \multiput{} at 0 0  2 4 /
  \setsolid
  \plot 0 0  1 0  1 4  0 4  0 0 /
  \plot 0 2  2 2  2 3  0 3 /
  \plot 0 1  1 1 /
  \multiput{$\sssize\bullet$} at 0.5 2.5  1.5 2.5 /
  \plot 0.5 2.5  1.5  2.5 /
  \plot 0.5 2.45  1.5 2.45 /
  \plot 0.5 2.55  1.5 2.55 /
  \endpicture} at -2.5 1.55 
\put{\beginpicture
   \setcoordinatesystem units <.14cm,.14cm>
  \multiput{} at 0 0  1 5 /
  \setsolid
  \plot 0 0  1 0  1 5  0 5  0 0 /
  \plot 0 1  1 1 /
  \plot 0 2  1 2 /
  \plot 0 3  1 3 /
  \plot 0 4  1 4 /
  \put{$\sssize\bullet$} at 0.5 4.5 
  \endpicture} at -1 5
\put{\beginpicture
   \setcoordinatesystem units <.14cm,.14cm>
  \multiput{} at 0 0  1 6 /
  \setsolid
  \plot 0 0  1 0  1 6  0 6  0 0 /
  \plot 0 1  1 1 /
  \plot 0 2  1 2 /
  \plot 0 3  1 3 /
  \plot 0 4  1 4 /
  \plot 0 5  1 5 /
  \put{$\sssize\bullet$} at 0.5 5.5 
  \endpicture} at 0 6
\put{\beginpicture
   \setcoordinatesystem units <.14cm,.14cm>
  \multiput{} at 0 0  1 6 /
  \setsolid
  \plot 0 0  1 0  1 6  0 6  0 0 /
  \plot 0 1  1 1 /
  \plot 0 2  1 2 /
  \plot 0 3  1 3 /
  \plot 0 4  1 4 /
  \plot 0 5  1 5 /
  \put{$\sssize\bullet$} at 0.5 4.5 
  \endpicture} at 1 5 
\put{\beginpicture
   \setcoordinatesystem units <.14cm,.14cm>
  \multiput{} at 0 0  1 6 /
  \setsolid
  \plot 0 0  1 0  1 6  0 6  0 0 /
  \plot 0 1  1 1 /
  \plot 0 2  1 2 /
  \plot 0 3  1 3 /
  \plot 0 4  1 4 /
  \plot 0 5  1 5 /
  \put{$\sssize\bullet$} at 0.5 3.5 
  \endpicture} at 2 4 
\put{\beginpicture
   \setcoordinatesystem units <.14cm,.14cm>
  \multiput{} at 0 0  1 6 /
  \setsolid
  \plot 0 0  1 0  1 6  0 6  0 0 /
  \plot 0 1  1 1 /
  \plot 0 2  1 2 /
  \plot 0 3  1 3 /
  \plot 0 4  1 4 /
  \plot 0 5  1 5 /
  \put{$\sssize\bullet$} at 0.5 2.5 
  \endpicture} at 3 3 
\put{\beginpicture
   \setcoordinatesystem units <.14cm,.14cm>
  \multiput{} at 0 0  1 6 /
  \setsolid
  \plot 0 0  1 0  1 6  0 6  0 0 /
  \plot 0 1  1 1 /
  \plot 0 2  1 2 /
  \plot 0 3  1 3 /
  \plot 0 4  1 4 /
  \plot 0 5  1 5 /
  \put{$\sssize\bullet$} at 0.5 1.5 
  \endpicture} at 4 2
\put{\beginpicture
   \setcoordinatesystem units <.14cm,.14cm>
  \multiput{} at 0 0  1 6 /
  \setsolid
  \plot 0 0  1 0  1 6  0 6  0 0 /
  \plot 0 1  1 1 /
  \plot 0 2  1 2 /
  \plot 0 3  1 3 /
  \plot 0 4  1 4 /
  \plot 0 5  1 5 /
  \put{$\sssize\bullet$} at 0.5 0.5 
  \endpicture} at 5 1
\put{\beginpicture
   \setcoordinatesystem units <.14cm,.14cm>
  \multiput{} at 0 0  1 5 /
  \setsolid
  \plot 0 0  1 0  1 5  0 5  0 0 /
  \plot 0 1  1 1 /
  \plot 0 2  1 2 /
  \plot 0 3  1 3 /
  \plot 0 4  1 4 /
  \endpicture} at 4 0 
\put{\beginpicture
   \setcoordinatesystem units <.14cm,.14cm>
  \multiput{} at 0 0  1 6 /
  \setsolid
  \plot 0 0  1 0  1 6  0 6  0 0 /
  \plot 0 1  1 1 /
  \plot 0 2  1 2 /
  \plot 0 3  1 3 /
  \plot 0 4  1 4 /
  \plot 0 5  1 5 /
  \endpicture} at 6 0 
\put{\beginpicture
   \setcoordinatesystem units <.14cm,.14cm>
  \multiput{} at 0 0  2 6 /
  \setsolid
  \plot 0 0  1 0  1 6  0 6  0 0 /
  \plot 0 2  2 2  2 5  0 5 /
  \plot 0 4  2 4 /
  \plot 0 3  2 3 /
  \plot 0 1  1 1 /
  \multiput{$\sssize\bullet$} at 0.5 4.5  1.5 4.5  1.5 3.5 /
  \plot 0.5 4.5   1.5 4.5 /
  \plot 0.5 4.45  1.5 4.45 /
  \plot 0.5 4.55  1.5 4.55 /
  \endpicture} at -.8 3.5 
\put{\beginpicture
   \setcoordinatesystem units <.14cm,.14cm>
  \multiput{} at 0 0  2 6 /
  \setsolid
  \plot 0 0  1 0  1 6  0 6  0 0 /
  \plot 0 1  2 1  2 5  0 5 /
  \plot 0 4  2 4 /
  \plot 0 3  2 3 /
  \plot 0 2  2 2 /
  \multiput{$\sssize\bullet$} at 0.5 4.5  1.5 4.5  1.5 3.5 /
  \plot 0.5 4.5   1.5 4.5 /
  \plot 0.5 4.45  1.5 4.45 /
  \plot 0.5 4.55  1.5 4.55 /
  \endpicture} at .3 4.15
\put{\beginpicture
   \setcoordinatesystem units <.14cm,.14cm>
  \multiput{} at 0 0  2 6 /
  \setsolid
  \plot 0 0  1 0  1 6  0 6  0 0 /
  \plot 0 1  2 1  2 5  0 5 /
  \plot 0 4  2 4 /
  \plot 0 3  2 3 /
  \plot 0 2  2 2 /
  \multiput{$\sssize\bullet$} at 0.5 4.5  1.5 4.5  1.5 2.5 /
  \plot 0.5 4.5   1.5 4.5 /
  \plot 0.5 4.45  1.5 4.45 /
  \plot 0.5 4.55  1.5 4.55 /
  \endpicture} at .9 3.5
\put{\beginpicture
   \setcoordinatesystem units <.14cm,.14cm>
  \multiput{} at 0 0  2 6 /
  \setsolid
  \plot 0 0  1 0  1 6  0 6  0 0 /
  \plot 0 1  2 1  2 5  0 5 /
  \plot 0 4  2 4 /
  \plot 0 3  2 3 /
  \plot 0 2  2 2 /
  \multiput{$\sssize\bullet$} at 0.5 2.5  1.5 2.5  /
  \plot 0.5 2.5  1.5  2.5 /
  \plot 0.5 2.45  1.5 2.45 /
  \plot 0.5 2.55  1.5 2.55 /
  \endpicture} at 2.45 1.6
\put{\beginpicture
   \setcoordinatesystem units <.14cm,.14cm>
  \multiput{} at 0 0  2 6 /
  \setsolid
  \plot 0 0  1 0  1 6  0 6  0 0 /
  \plot 0 1  2 1  2 5  0 5 /
  \plot 0 4  2 4 /
  \plot 0 3  2 3 /
  \plot 0 2  2 2 /
  \multiput{$\sssize\bullet$} at 0.5 1.5  1.5 1.5 /
  \plot 0.5 1.5  1.5  1.5 /
  \plot 0.5 1.45  1.5 1.45 /
  \plot 0.5 1.55  1.5 1.55 /
  \endpicture} at 3.35 .95
\put{\beginpicture
   \setcoordinatesystem units <.14cm,.14cm>
  \multiput{} at 0 0  2 6 /
  \setsolid
  \plot 0 0  1 0  1 6  0 6  0 0 /
  \plot 0 1  2 1  2 4  0 4 /
  \plot 0 5  1 5 /
  \plot 0 3  2 3 /
  \plot 0 2  2 2 /
  \multiput{$\sssize\bullet$} at 0.5 1.5  1.5 1.5 /
  \plot 0.5 1.5  1.5  1.5 /
  \plot 0.5 1.45  1.5 1.45 /
  \plot 0.5 1.55  1.5 1.55 /
  \endpicture} at 1.95 .95
\endpicture} at 9 -5
\endpicture}
$$

\medskip
Now we are going to exhibit the boundary of the hexagon.

\medskip
{\bf (a) The pickets in the boundary of the hexagon of $\Cal S(6)$.} 
$$
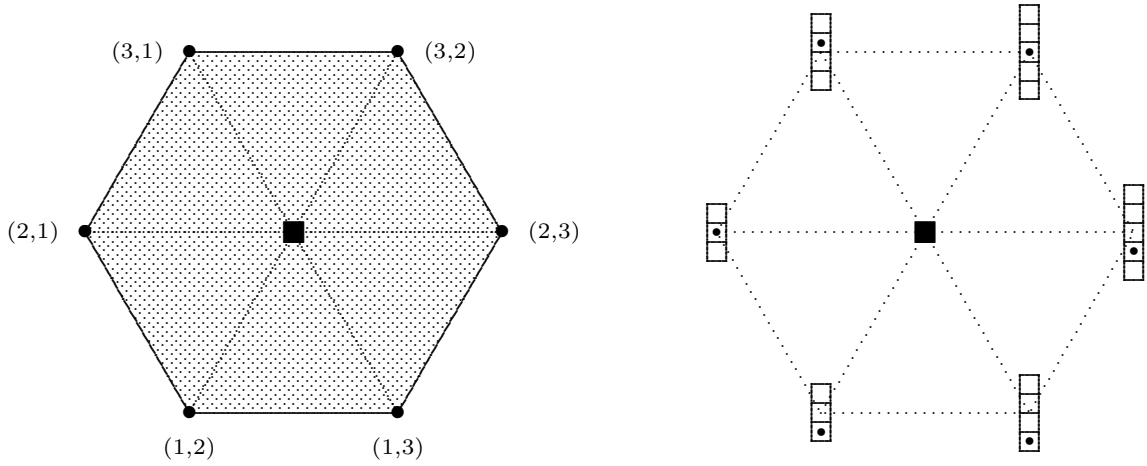
{\beginpicture
   \setcoordinatesystem units <1.1547cm,2cm>
   \multiput{} at -2 0  2 2 /
   \plot -2 1  -1 2  1 2  2 1  1 0  -1 0  -2 1 /
   \multiput{$\bullet$} at -2 1  -1 2  1 2  2 1  1 0  -1 0  -2 1 /
   \put{$\blacksquare$} at 0 1
   \setdots<.5mm>
   \plot -2 1  2 1 /
   \plot -1 0  1 2 /
   \plot -1 2  1 0 /
   \setshadegrid span <.5mm>
   \vshade -2 1 1  <,z,,>  -1 0 2  <z,z,,>  1 0 2 <z,z,,>  2 1 1 /
   \put{$\ss (1,2)$} at -1 -.2
   \put{$\ss (1,3)$} at 1 -.2
   \put{$\ss (2,1)$} at -2.5 1
   \put{$\ss (3,1)$} at -1.5 2
   \put{$\ss (2,3)$} at 2.5 1
   \put{$\ss (3,2)$} at 1.5 2 
\endpicture}
\qquad
\qquad
{\beginpicture
   \setcoordinatesystem units <1.1547cm,2cm>
   \multiput{} at -2 0  2 2 /
   \setdots<1mm>
   \plot -2 1  -1 2  1 2  2 1  1 0  -1 0  -2 1 /
   \put{$\blacksquare$} at 0 1
   \plot -2 1  2 1 /
   \plot -1 0  1 2 /
   \plot -1 2  1 0 /
   \put{\beginpicture
     \setcoordinatesystem units <.21cm,.21cm>
     \multiput{} at 0 0  1 3  /
     \setsolid
     \plot 0 0  0 3  1 3  1 0  0 0 /
     \plot 0 1  1 1 /
     \plot 0 2  1 2 /
     \put{$\sss\bullet$} at .5 .5 
     \endpicture} at -1 0 
   \put{\beginpicture
     \setcoordinatesystem units <.21cm,.21cm>
     \multiput{} at 0 0  1 3  /
     \setsolid
     \plot 0 0  0 3  1 3  1 0  0 0 /
     \plot 0 1  1 1 /
     \plot 0 2  1 2 /
     \put{$\sss\bullet$} at .5 1.5 
     \endpicture} at -2 1 
   \put{\beginpicture
     \setcoordinatesystem units <.21cm,.21cm>
     \multiput{} at 0 0  1 4  /
     \setsolid
     \plot 0 0  0 4  1 4  1 0  0 0 /
     \plot 0 1  1 1 /
     \plot 0 2  1 2 /
     \plot 0 3  1 3 /
     \put{$\sss\bullet$} at .5 .5 
     \endpicture} at 1 0 
   \put{\beginpicture
     \setcoordinatesystem units <.21cm,.21cm>
     \multiput{} at 0 0  1 4  /
     \setsolid
     \plot 0 0  0 4  1 4  1 0  0 0 /
     \plot 0 1  1 1 /
     \plot 0 2  1 2 /
     \plot 0 3  1 3 /
     \put{$\sss\bullet$} at .5 2.5 
     \endpicture} at -1 2 
   \put{\beginpicture
     \setcoordinatesystem units <.21cm,.21cm>
     \multiput{} at 0 0  1 5  /
     \setsolid
     \plot 0 0  0 5  1 5  1 0  0 0 /
     \plot 0 1  1 1 /
     \plot 0 2  1 2 /
     \plot 0 3  1 3 /
     \plot 0 4  1 4 /
     \put{$\sss\bullet$} at .5 1.5 
     \endpicture} at 2 1
   \put{\beginpicture
     \setcoordinatesystem units <.21cm,.21cm>
     \multiput{} at 0 0  1 5  /
     \setsolid
     \plot 0 0  0 5  1 5  1 0  0 0 /
     \plot 0 1  1 1 /
     \plot 0 2  1 2 /
     \plot 0 3  1 3 /
     \plot 0 4  1 4 /
     \put{$\sss\bullet$} at .5 2.5 
     \endpicture} at 1 2 
\endpicture}
$$

\bigskip
{\bf (b) The bipickets in the boundary of the hexagon of $\Cal S(6)$:}  
$$
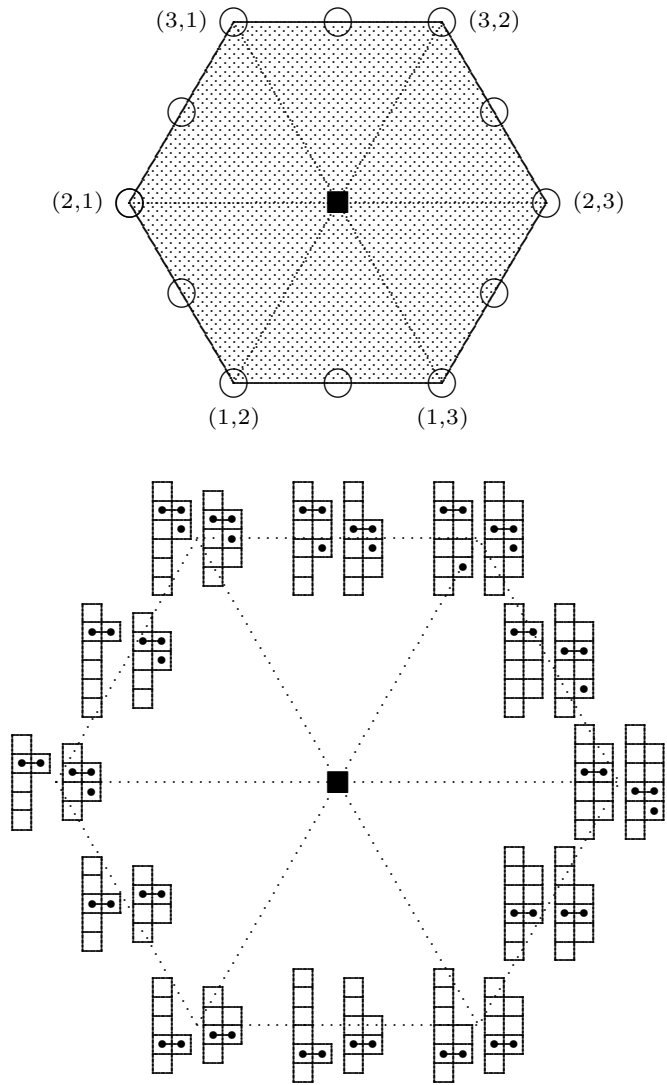
{\beginpicture
   \setcoordinatesystem units <1.1547cm,2cm>
   \multiput{} at -2 0  2 2 /
   \plot -2 1  -1 2  1 2  2 1  1 0  -1 0  -2 1 /
   \multiput{$\bigcirc$} at -2 1  -1.5 1.5  -1 2  0 2  1 2  1.5 1.5  2 1  1.5 .5
                         1 0  0 0  -1 0  -1.5 .5  -2 1 /
   \put{$\blacksquare$} at 0 1
   \setdots<.5mm>
   \plot -2 1  2 1 /
   \plot -1 0  1 2 /
   \plot -1 2  1 0 /
   \setshadegrid span <.5mm>
   \vshade -2 1 1  <,z,,>  -1 0 2  <z,z,,>  1 0 2 <z,z,,>  2 1 1 /
   \put{$\ss (1,2)$} at -1 -.2
   \put{$\ss (1,3)$} at 1 -.2
   \put{$\ss (2,1)$} at -2.5 1
   \put{$\ss (3,1)$} at -1.5 2
   \put{$\ss (2,3)$} at 2.5 1
   \put{$\ss (3,2)$} at 1.5 2 
\endpicture}
$$

\bigskip
$$
{\beginpicture
   \setcoordinatesystem units <1.5588cm, 2.7cm>
   \multiput{} at -2 0  2 2 /
   \setdots<1mm>
   \plot -2 1  -1 2  1 2  2 1  1 0  -1 0  -2 1 /
   \put{$\blacksquare$} at 0 1
   \plot -2 1  2 1 /
   \plot -1 0  1 2 /
   \plot -1 2  1 0 /
   \put{\beginpicture
     \setcoordinatesystem units <.21cm,.21cm>
     \multiput{} at 0 0  2 5  /
     \setsolid
     \plot 0 0  0 5  1 5  1 0  0 0 /
     \plot 0 1  2 1  2 2  0 2  /
     \plot 0 3  1 3 /
     \plot 0 4  1 4 /
     \multiput{$\sss\bullet$} at .5 1.5  1.5 1.5 /
     \plot .5 1.5  1.5 1.5 /
     \endpicture} at -1.18 0 
   \put{\beginpicture
     \setcoordinatesystem units <.21cm,.21cm>
     \multiput{} at 0 0  2 4  /
     \setsolid
     \plot 0 0  0 4  1 4  1 0  0 0 /
     \plot 0 1  2 1  2 3  0 3  /
     \plot 0 2  2 2 /
     \multiput{$\sss\bullet$} at .5 1.5  1.5 1.5 /
     \plot .5 1.5  1.5 1.5 /
     \endpicture} at -.82  0 
   \put{\beginpicture
     \setcoordinatesystem units <.21cm,.21cm>
     \multiput{} at 0 0  2 6  /
     \setsolid
     \plot 0 0  0 6  1 6  1 0  0 0 /
     \plot 0 1  2 1  2 2  0 2  /
     \plot 0 3  1 3 /
     \plot 0 4  1 4 /
     \plot 0 5  1 5 /
     \multiput{$\sss\bullet$} at .5 1.5  1.5 1.5 /
     \plot .5 1.5  1.5 1.5 /
     \endpicture} at -.18 0 
   \put{\beginpicture
     \setcoordinatesystem units <.21cm,.21cm>
     \multiput{} at 0 0  2 5  /
     \setsolid
     \plot 0 0  0 5  1 5  1 0  0 0 /
     \plot 0 1  2 1  2 3  0 3  /
     \plot 0 2  2 2 /
     \plot 0 4  1 4 /
     \multiput{$\sss\bullet$} at .5 1.5  1.5 1.5 /
     \plot .5 1.5  1.5 1.5 /
     \endpicture} at .18 0 
   \put{\beginpicture
     \setcoordinatesystem units <.21cm,.21cm>
     \multiput{} at 0 0  2 6  /
     \setsolid
     \plot 0 0  0 6  1 6  1 0  0 0 /
     \plot 0 1  2 1  2 3  0 3  /
     \plot 0 2  2 2 /
     \plot 0 4  1 4 /
     \plot 0 5  1 5 /
     \multiput{$\sss\bullet$} at .5 1.5  1.5 1.5 /
     \plot .5 1.5  1.5 1.5 /
     \endpicture} at .82 0 
   \put{\beginpicture
     \setcoordinatesystem units <.21cm,.21cm>
     \multiput{} at 0 0  2 5  /
     \setsolid
     \plot 0 0  0 5  1 5  1 0  0 0 /
     \plot 0 1  2 1  2 4  0 4  /
     \plot 0 2  2 2 /
     \plot 0 3  2 3 /
     \multiput{$\sss\bullet$} at .5 1.5  1.5 1.5 /
     \plot .5 1.5  1.5 1.5 /
     \endpicture} at 1.18 0 
   \put{\beginpicture
     \setcoordinatesystem units <.21cm,.21cm>
     \multiput{} at 0 0  2 6  /
     \setsolid
     \plot 0 0  0 6  1 6  1 0  0 0 /
     \plot 0 2  2 2  2 5  0 5  /
     \plot 0 1  1 1 /
     \plot 0 4  2 4 /
     \plot 0 3  2 3 /
     \multiput{$\sss\bullet$} at .5 2.5  1.5 2.5 /
     \plot .5 2.5  1.5 2.5 /
     \endpicture} at 1.32 .5 
   \put{\beginpicture
     \setcoordinatesystem units <.21cm,.21cm>
     \multiput{} at 0 0  2 6  /
     \setsolid
     \plot 0 0  0 6  1 6  1 0  0 0 /
     \plot 0 1  2 1  2 4  0 4  /
     \plot 0 2  2 2 /
     \plot 0 3  2 3 /
     \plot 0 5  1 5 /
     \multiput{$\sss\bullet$} at .5 2.5  1.5 2.5 /
     \plot .5 2.5  1.5 2.5 /
     \endpicture} at 1.68 .5 
   \put{\beginpicture
     \setcoordinatesystem units <.21cm,.21cm>
     \multiput{} at 0 0  2 6  /
     \setsolid
     \plot 0 0  0 6  1 6  1 0  0 0 /
     \plot 0 1  2 1  2 5  0 5  /
     \plot 0 2  2 2 /
     \plot 0 4  2 4 /
     \plot 0 3  2 3 /
     \multiput{$\sss\bullet$} at .5 3.5  1.5 3.5 /
     \plot .5 3.5  1.5 3.5 /
     \endpicture} at 1.82 1 
   \put{\beginpicture
     \setcoordinatesystem units <.21cm,.21cm>
     \multiput{} at 0 0  2 6  /
     \setsolid
     \plot 0 0  0 6  1 6  1 0  0 0 /
     \plot 0 1  2 1  2 5  0 5  /
     \plot 0 2  2 2 /
     \plot 0 3  2 3 /
     \plot 0 4  2 4 /
     \multiput{$\sss\bullet$} at .5 2.5  1.5 2.5  1.5 1.5  /
     \plot .5 2.5  1.5 2.5 /
     \endpicture} at 2.18 1 
   \put{\beginpicture
     \setcoordinatesystem units <.21cm,.21cm>
     \multiput{} at 0 0  2 6  /
     \setsolid
     \plot 0 0  0 6  1 6  1 0  0 0 /
     \plot 0 1  2 1  2 5  0 5  /
     \plot 0 2  2 2 /
     \plot 0 4  2 4 /
     \plot 0 3  2 3 /
     \multiput{$\sss\bullet$} at .5 4.5  1.5 4.5 /
     \plot .5 4.5  1.5 4.5 /
     \endpicture} at 1.32 1.5 
   \put{\beginpicture
     \setcoordinatesystem units <.21cm,.21cm>
     \multiput{} at 0 0  2 6  /
     \setsolid
     \plot 0 0  0 6  1 6  1 0  0 0 /
     \plot 0 1  2 1  2 5  0 5  /
     \plot 0 2  2 2 /
     \plot 0 3  2 3 /
     \plot 0 4  2 4 /
     \multiput{$\sss\bullet$} at .5 3.5  1.5 3.5  1.5 1.5  /
     \plot .5 3.5  1.5 3.5 /
     \endpicture} at 1.68 1.5 
   \put{\beginpicture
     \setcoordinatesystem units <.21cm,.21cm>
     \multiput{} at 0 0  2 6  /
     \setsolid
     \plot 0 0  0 6  1 6  1 0  0 0 /
     \plot 0 1  2 1  2 5  0 5  /
     \plot 0 2  2 2 /
     \plot 0 4  2 4 /
     \plot 0 3  2 3 /
     \multiput{$\sss\bullet$} at .5 4.5  1.5 4.5  1.5 1.5 /
     \plot .5 4.5  1.5 4.5 /
     \endpicture} at .82 2
   \put{\beginpicture
     \setcoordinatesystem units <.21cm,.21cm>
     \multiput{} at 0 0  2 6  /
     \setsolid
     \plot 0 0  0 6  1 6  1 0  0 0 /
     \plot 0 1  2 1  2 5  0 5  /
     \plot 0 2  2 2 /
     \plot 0 3  2 3 /
     \plot 0 4  2 4 /
     \multiput{$\sss\bullet$} at .5 3.5  1.5 3.5  1.5 2.5  /
     \plot .5 3.5  1.5 3.5 /
     \endpicture} at 1.18 2 
   \put{\beginpicture
     \setcoordinatesystem units <.21cm,.21cm>
     \multiput{} at 0 0  2 6  /
     \setsolid
     \plot 0 0  0 6  1 6  1 0  0 0 /
     \plot 0 2  2 2  2 5  0 5  /
     \plot 0 1  1 1 /
     \plot 0 4  2 4 /
     \plot 0 3  2 3 /
     \multiput{$\sss\bullet$} at .5 4.5  1.5 4.5  1.5 2.5 /
     \plot .5 4.5  1.5 4.5 /
     \endpicture} at -.18 2
   \put{\beginpicture
     \setcoordinatesystem units <.21cm,.21cm>
     \multiput{} at 0 0  2 6  /
     \setsolid
     \plot 0 0  0 6  1 6  1 0  0 0 /
     \plot 0 1  2 1  2 4  0 4  /
     \plot 0 2  2 2 /
     \plot 0 3  2 3 /
     \plot 0 5  1 5 /
     \multiput{$\sss\bullet$} at .5 3.5  1.5 3.5  1.5 2.5  /
     \plot .5 3.5  1.5 3.5 /
     \endpicture} at .18 2 
   \put{\beginpicture
     \setcoordinatesystem units <.21cm,.21cm>
     \multiput{} at 0 0  2 6  /
     \setsolid
     \plot 0 0  0 6  1 6  1 0  0 0 /
     \plot 0 3  2 3  2 5  0 5  /
     \plot 0 1  1 1 /
     \plot 0 4  2 4 /
     \plot 0 2  1 2 /
     \multiput{$\sss\bullet$} at .5 4.5  1.5 4.5  1.5 3.5 /
     \plot .5 4.5  1.5 4.5 /
     \endpicture} at -1.18 2
   \put{\beginpicture
     \setcoordinatesystem units <.21cm,.21cm>
     \multiput{} at 0 0  2 5  /
     \setsolid
     \plot 0 0  0 5  1 5  1 0  0 0 /
     \plot 0 1  2 1  2 4  0 4  /
     \plot 0 2  2 2 /
     \plot 0 3  2 3 /
     \multiput{$\sss\bullet$} at .5 3.5  1.5 3.5  1.5 2.5  /
     \plot .5 3.5  1.5 3.5 /
     \endpicture} at -.82 2 
   \put{\beginpicture
     \setcoordinatesystem units <.21cm,.21cm>
     \multiput{} at 0 0  2 6  /
     \setsolid
     \plot 0 0  0 6  1 6  1 0  0 0 /
     \plot 0 4  2 4  2 5  0 5  /
     \plot 0 1  1 1 /
     \plot 0 3  1 3 /
     \plot 0 2  1 2 /
     \multiput{$\sss\bullet$} at .5 4.5  1.5 4.5 /
     \plot .5 4.5  1.5 4.5 /
     \endpicture} at -1.68 1.5
   \put{\beginpicture
     \setcoordinatesystem units <.21cm,.21cm>
     \multiput{} at 0 0  2 5  /
     \setsolid
     \plot 0 0  0 5  1 5  1 0  0 0 /
     \plot 0 2  2 2  2 4  0 4  /
     \plot 0 1  1 1 /
     \plot 0 3  2 3 /
     \multiput{$\sss\bullet$} at .5 3.5  1.5 3.5  1.5 2.5  /
     \plot .5 3.5  1.5 3.5 /
     \endpicture} at -1.32 1.5 
   \put{\beginpicture
     \setcoordinatesystem units <.21cm,.21cm>
     \multiput{} at 0 0  2 5  /
     \setsolid
     \plot 0 0  0 5  1 5  1 0  0 0 /
     \plot 0 3  2 3  2 4  0 4  /
     \plot 0 1  1 1 /
     \plot 0 2  1 2 /
     \multiput{$\sss\bullet$} at .5 3.5  1.5 3.5 /
     \plot .5 3.5  1.5 3.5 /
     \endpicture} at -2.18 1 
   \put{\beginpicture
     \setcoordinatesystem units <.21cm,.21cm>
     \multiput{} at 0 0  2 4  /
     \setsolid
     \plot 0 0  0 4  1 4  1 0  0 0 /
     \plot 0 1  2 1  2 3  0 3  /
     \plot 0 2  2 2 /
     \multiput{$\sss\bullet$} at .5 2.5  1.5 2.5  1.5 1.5  /
     \plot .5 2.5  1.5 2.5 /
     \endpicture} at -1.82 1
   \put{\beginpicture
     \setcoordinatesystem units <.21cm,.21cm>
     \multiput{} at 0 0  2 5  /
     \setsolid
     \plot 0 0  0 5  1 5  1 0  0 0 /
     \plot 0 2  2 2  2 3  0 3  /
     \plot 0 1  1 1 /
     \plot 0 4  1 4 /
     \multiput{$\sss\bullet$} at .5 2.5  1.5 2.5 /
     \plot .5 2.5  1.5 2.5 /
     \endpicture} at -1.68 .5 
   \put{\beginpicture
     \setcoordinatesystem units <.21cm,.21cm>
     \multiput{} at 0 0  2 4  /
     \setsolid
     \plot 0 0  0 4  1 4  1 0  0 0 /
     \plot 0 1  2 1  2 3  0 3  /
     \plot 0 2  2 2 /
     \multiput{$\sss\bullet$} at .5 2.5  1.5 2.5   /
     \plot .5 2.5  1.5 2.5 /
     \endpicture} at -1.32 .5 
\endpicture}
$$

\medskip
The position (rationality index, rank of tube) of each bipicket in $\Cal S(6)$
is specified in Section \ref{sec-fifteen-two}--(g).

	\bigskip
{\bf (c) The tripickets in the boundary of the hexagon of $\Cal S(6)$:}

\smallskip
We picture the 18
\phantomsection{tripickets}
on the boundary of the hexagon.

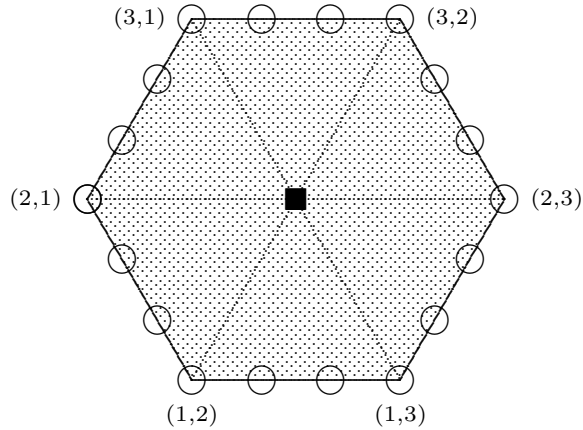
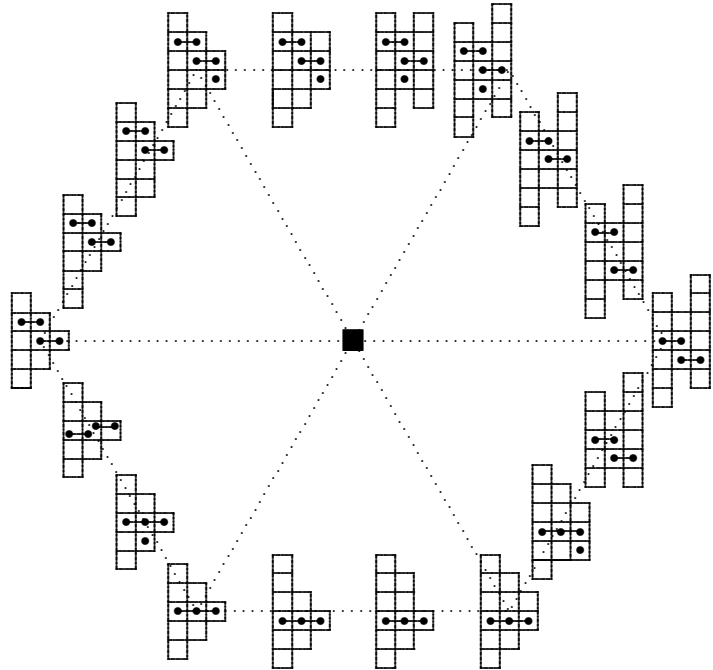
\addcontentsline{lof}{subsection}{The 18 tripickets on the boundary of the hexagon of $\mathcal S(6)$.}

$$
{\beginpicture
   \setcoordinatesystem units <1.1547cm,2cm>
   \multiput{} at -2 0  2 2 /
   \plot -2 1  -1 2  1 2  2 1  1 0  -1 0  -2 1 /
   \multiput{$\bigcirc$} at -2 1  -1.67 1.33  -1.33 1.67  -1 2 -.33 2  .33 2  1 2
                         1.33 1.67  1.67 1.33   2 1  1.67 .67  1.33 .33 
                         1 0  .33 0  -.33 0  -1 0  -1.33 .33  -1.67 .67  -2 1 /
   \put{$\blacksquare$} at 0 1
   \setdots<.5mm>
   \plot -2 1  2 1 /
   \plot -1 0  1 2 /
   \plot -1 2  1 0 /
   \setshadegrid span <.5mm>
   \vshade -2 1 1  <,z,,>  -1 0 2  <z,z,,>  1 0 2 <z,z,,>  2 1 1 /
   \put{$\ss (1,2)$} at -1 -.2
   \put{$\ss (1,3)$} at 1 -.2
   \put{$\ss (2,1)$} at -2.5 1
   \put{$\ss (3,1)$} at -1.5 2
   \put{$\ss (2,3)$} at 2.5 1
   \put{$\ss (3,2)$} at 1.5 2 
\endpicture}
$$
$$
{\beginpicture
   \setcoordinatesystem units <1.73205cm, 3cm>
   \multiput{} at -2 0  2 2 /
   \setdots<1mm>
   \plot -2 1  -1 2  1 2  2 1  1 0  -1 0  -2 1 /
   \put{$\blacksquare$} at 0 1
   \plot -2 1  2 1 /
   \plot -1 0  1 2 /
   \plot -1 2  1 0 /
   \put{\beginpicture
     \setcoordinatesystem units <.21cm,.21cm>
     \multiput{} at 0 0  3 5  /
     \setsolid
     \plot 0 0  0 5  1 5  1 0  0 0 /
     \plot 0 1  2 1  2 4  0 4 /
     \plot 0 2  3 2  3 3  0 3 /
     \multiput{$\sss\bullet$} at .5 2.5  1.5 2.5  2.5 2.5 /
     \plot .5 2.5  2.5 2.5 /
     \endpicture} at -1 0 
   \put{\beginpicture
     \setcoordinatesystem units <.21cm,.21cm>
     \multiput{} at 0 0  3 6  /
     \setsolid
     \plot 0 0  0 6  1 6  1 0  0 0 /
     \plot 0 1  2 1  2 4  0 4 /
     \plot 0 2  3 2  3 3  0 3 /
     \plot 0 5  1 5 /
     \multiput{$\sss\bullet$} at .5 2.5  1.5 2.5  2.5 2.5 /
     \plot .5 2.5  2.5 2.5 /
     \endpicture} at -.33 0 
   \put{\beginpicture
     \setcoordinatesystem units <.21cm,.21cm>
     \multiput{} at 0 0  3 6  /
     \setsolid
     \plot 0 0  0 6  1 6  1 0  0 0 /
     \plot 0 1  2 1  2 5  0 5 /
     \plot 0 2  3 2  3 3  0 3 /
     \plot 0 4  2 4 /
     \multiput{$\sss\bullet$} at .5 2.5  1.5 2.5  2.5 2.5 /
     \plot .5 2.5  2.5 2.5 /
     \endpicture} at .33 0 
   \put{\beginpicture
     \setcoordinatesystem units <.21cm,.21cm>
     \multiput{} at 0 0  3 6  /
     \setsolid
     \plot 0 0  0 6  1 6  1 0  0 0 /
     \plot 0 1  2 1  2 5  0 5 /
     \plot 0 2  3 2  3 4  0 4 /
     \plot 0 3  3 3 /
     \multiput{$\sss\bullet$} at .5 2.5  1.5 2.5  2.5 2.5 /
     \plot .5 2.5  2.5 2.5 /
     \endpicture} at 1 0 
   \put{\beginpicture
     \setcoordinatesystem units <.21cm,.21cm>
     \multiput{} at 0 0  3 6  /
     \setsolid
     \plot 0 0  0 6  1 6  1 0  0 0 /
     \plot 0 1  2 1  2 5  0 5 /
     \plot 0 1  3 1  3 4  0 4 /
     \plot 0 3  3 3 /
     \plot 0 2  3 2 /
     \multiput{$\sss\bullet$} at .5 2.5  1.5 2.5  2.5 2.5  2.5 1.5 /
     \plot .5 2.5  2.5 2.5 /
     \endpicture} at 1.33 .33
   \put{\beginpicture
     \setcoordinatesystem units <.21cm,.21cm>
     \multiput{} at 0 0  3 6  /
     \setsolid
     \plot 0 0  0 5  1 5  1 0  0 0 /
     \plot 2 0  2 6  3 6  3 0  2 0 /
     \plot 0 1  3 1 /
     \plot 0 2  3 2 /
     \plot 0 3  3 3 /
     \plot 0 4  3 4 /
     \plot 2 5  3 5 /
     \multiput{$\sss\bullet$} at .5 2.5  1.5 2.5  1.5 1.5  2.5 1.5 /
     \plot .5 2.5  1.5 2.5 /
     \plot 1.5 1.5  2.5 1.5 /
     \endpicture} at 1.67 .67
   \put{\beginpicture
     \setcoordinatesystem units <.21cm,.21cm>
     \multiput{} at 0 -1  3 6  /
     \setsolid
     \plot 0 -1  0 5  1 5  1 -1  0 -1 /
     \plot 2 0  2 6  3 6  3 0  2 0 /
     \plot 0 0  1 0 /
     \plot 0 1  3 1 /
     \plot 0 2  3 2 /
     \plot 0 3  3 3 /
     \plot 0 4  3 4 /
     \plot 2 5  3 5 /
     \multiput{$\sss\bullet$} at .5 2.5  1.5 2.5  1.5 1.5  2.5 1.5 /
     \plot .5 2.5  1.5 2.5 /
     \plot 1.5 1.5  2.5 1.5 /
     \endpicture} at 2.1 1 
   \put{\beginpicture
     \setcoordinatesystem units <.21cm,.21cm>
     \multiput{} at 0 -1  3 6  /
     \setsolid
     \plot 0 -1  0 5  1 5  1 -1  0 -1 /
     \plot 2 0  2 6  3 6  3 0  2 0 /
     \plot 0 0  1 0 /
     \plot 0 1  3 1 /
     \plot 0 2  3 2 /
     \plot 0 3  3 3 /
     \plot 0 4  3 4 /
     \plot 2 5  3 5 /
     \multiput{$\sss\bullet$} at .5 3.5  1.5 3.5  1.5 1.5  2.5 1.5 /
     \plot .5 3.5  1.5 3.5 /
     \plot 1.5 1.5  2.5 1.5 /
     \endpicture} at 1.67 1.33 
   \put{\beginpicture
     \setcoordinatesystem units <.21cm,.21cm>
     \multiput{} at 0 -1  3 6  /
     \setsolid
     \plot 0 -1  0 5  1 5  1 -1  0 -1 /
     \plot 2 0  2 6  3 6  3 0  2 0 /
     \plot 0 0  1 0 /
     \plot 0 1  3 1 /
     \plot 0 2  3 2 /
     \plot 0 3  3 3 /
     \plot 0 4  3 4 /
     \plot 2 5  3 5 /
     \multiput{$\sss\bullet$} at .5 3.5  1.5 3.5  1.5 2.5  2.5 2.5 /
     \plot .5 3.5  1.5 3.5 /
     \plot 1.5 2.5  2.5 2.5 /
     \endpicture} at 1.25 1.67 
   \put{\beginpicture
     \setcoordinatesystem units <.21cm,.21cm>
     \multiput{} at 0 -1  3 6  /
     \setsolid
     \plot 0 -1  0 5  1 5  1 -1  0 -1 /
     \plot 2 0  2 6  3 6  3 0  2 0 /
     \plot 0 0  1 0 /
     \plot 0 1  3 1 /
     \plot 0 2  3 2 /
     \plot 0 3  3 3 /
     \plot 0 4  3 4 /
     \plot 2 5  3 5 /
     \multiput{$\sss\bullet$} at .5 3.5  1.5 3.5  1.5 2.5  2.5 2.5  1.5 1.5 /
     \plot .5 3.5  1.5 3.5 /
     \plot 1.5 2.5  2.5 2.5 /
     \endpicture} at .83 2 
   \put{\beginpicture
     \setcoordinatesystem units <.21cm,.21cm>
     \multiput{} at 0 -1  3 5  /
     \setsolid
     \plot 0 -1  0 5  1 5  1 -1  0 -1 /
     \plot 2 0  2 5  3 5  3 0  2 0 /
     \plot 0 0  1 0 /
     \plot 0 1  3 1 /
     \plot 0 2  3 2 /
     \plot 0 3  3 3 /
     \plot 0 4  3 4 /
     \multiput{$\sss\bullet$} at .5 3.5  1.5 3.5  1.5 2.5  2.5 2.5  1.5 1.5 /
     \plot .5 3.5  1.5 3.5 /
     \plot 1.5 2.5  2.5 2.5 /
     \endpicture} at .33 2 
   \put{\beginpicture
     \setcoordinatesystem units <.21cm,.21cm>
     \multiput{} at 0 -1  3 5  /
     \setsolid
     \plot 0 -1  0 5  1 5  1 -1  0 -1 /
     \plot 0 4  3 4  3 1  0 1  0 0  2 0  2 4 /
     \plot 0 2  3 2 /
     \plot 0 3  3 3 /
     \multiput{$\sss\bullet$} at .5 3.5  1.5 3.5  1.5 2.5  2.5 2.5  2.5 1.5 /
     \plot .5 3.5  1.5 3.5 /
     \plot 1.5 2.5  2.5 2.5 /
     \endpicture} at -.33 2 
   \put{\beginpicture
     \setcoordinatesystem units <.21cm,.21cm>
     \multiput{} at 0 -1  3 5  /
     \setsolid
     \plot 0 -1  0 5  1 5  1 -1  0 -1 /
     \plot 0 4  2 4  2 0  0 0  /
     \plot 0 1  3 1  3 3  0 3  /
     \plot 0 2  3 2 /
     \multiput{$\sss\bullet$} at .5 3.5  1.5 3.5  1.5 2.5  2.5 2.5  2.5 1.5 /
     \plot .5 3.5  1.5 3.5 /
     \plot 1.5 2.5  2.5 2.5 /
     \endpicture} at -1 2 
   \put{\beginpicture
     \setcoordinatesystem units <.21cm,.21cm>
     \multiput{} at 0 -1  3 5  /
     \setsolid
     \plot 0 -1  0 5  1 5  1 -1  0 -1 /
     \plot 0 4  2 4  2 0  0 0  /
     \plot 0 2  3 2  3 3  0 3  /
     \plot 0 1  2 1 /
     \multiput{$\sss\bullet$} at .5 3.5  1.5 3.5  1.5 2.5  2.5 2.5  /
     \plot .5 3.5  1.5 3.5 /
     \plot 1.5 2.5  2.5 2.5 /
     \endpicture} at -1.33 1.67
   \put{\beginpicture
     \setcoordinatesystem units <.21cm,.21cm>
     \multiput{} at 0 -1  3 5  /
     \setsolid
     \plot 0 -1  0 5  1 5  1 -1  0 -1 /
     \plot 0 4  2 4  2 1  0 1  /
     \plot 0 2  3 2  3 3  0 3  /
     \plot 0 0  1 0 /
     \multiput{$\sss\bullet$} at .5 3.5  1.5 3.5  1.5 2.5  2.5 2.5  /
     \plot .5 3.5  1.5 3.5 /
     \plot 1.5 2.5  2.5 2.5 /
     \endpicture} at -1.67 1.33 
   \put{\beginpicture
     \setcoordinatesystem units <.21cm,.21cm>
     \multiput{} at 0 0  3 5  /
     \setsolid
     \plot 0 0  0 5  1 5  1 0  0 0 /
     \plot 0 4  2 4  2 1  0 1  /
     \plot 0 2  3 2  3 3  0 3  /
     \plot 0 0  1 0 /
     \multiput{$\sss\bullet$} at .5 3.5  1.5 3.5  1.5 2.5  2.5 2.5  /
     \plot .5 3.5  1.5 3.5 /
     \plot 1.5 2.5  2.5 2.5 /
     \endpicture} at -2 1
   \put{\beginpicture
     \setcoordinatesystem units <.21cm,.21cm>
     \multiput{} at 0 0  3 5  /
     \setsolid
     \plot 0 0  0 5  1 5  1 0  0 0 /
     \plot 0 4  2 4  2 1  0 1  /
     \plot 0 2  3 2  3 3  0 3  /
     \multiput{$\sss\bullet$} at .3 2.3  1.3 2.3  1.7 2.7  2.7 2.7   /
     \plot .3 2.3  1.3 2.3  /
     \plot 1.7 2.7  2.7 2.7 /
     \endpicture} at -1.67 .67
   \put{\beginpicture
     \setcoordinatesystem units <.21cm,.21cm>
     \multiput{} at 0 0  3 5  /
     \setsolid
     \plot 0 0  0 5  1 5  1 0  0 0 /
     \plot 0 4  2 4  2 1  0 1  /
     \plot 0 2  3 2  3 3  0 3  /
     \multiput{$\sss\bullet$} at .5 2.5  1.5 2.5  1.5 1.5  2.5 2.5  /
     \plot .5 2.5  2.5 2.5 /
     \endpicture} at -1.33 .33
   \endpicture}
   $$

\bigskip
The 18 objects form 3 orbits under $\Sigma_3$.
The leftmost object in the bottom row occurs on the mouth of the tube of rank 3
and rationality index $\frac12$.  Its left and right neighbors are
on the mouth of the tube of rank 6 and index $\frac34$ and $\frac23$, respectively.

\bigskip
{\bf (d) The tetrapickets in the boundary of the hexagon of $\Cal S(6)$:}
$$
{\beginpicture
   \setcoordinatesystem units <1.1547cm,2cm>
   \multiput{} at -2 0  2 2 /
   \plot -2 1  -1 2  1 2  2 1  1 0  -1 0  -2 1 /
   \multiput{$\bigcirc$} at -1.5 1.5  0 2  1.5 1.5  1.5 .5  0 0  -1.5 .5 /
   \put{$\blacksquare$} at 0 1
   \setdots<.5mm>
   \plot -2 1  2 1 /
   \plot -1 0  1 2 /
   \plot -1 2  1 0 /
   \setshadegrid span <.5mm>
   \vshade -2 1 1  <,z,,>  -1 0 2  <z,z,,>  1 0 2 <z,z,,>  2 1 1 /
   \put{$\ss (1,\frac52)$} at 0 -.2
   \put{$\ss (\frac32,3)$} at 2 .4
   \put{$\ss (\frac32,\frac32)$} at -2 .4
   \put{$\ss (\frac52,1)$} at -2 1.6
   \put{$\ss (\frac52,\frac52)$} at 2 1.6
   \put{$\ss (3,\frac32)$} at 0 2.2
\endpicture}
\qquad\qquad
{\beginpicture
   \setcoordinatesystem units <1.1547cm,2cm>
   \multiput{} at -2 0  2 2 /
   \setdots<1mm>
   \plot -2 1  -1 2  1 2  2 1  1 0  -1 0  -2 1 /
   \put{$\blacksquare$} at 0 1
   \plot -2 1  2 1 /
   \plot -1 0  1 2 /
   \plot -1 2  1 0 /
   \put{\beginpicture
     \setcoordinatesystem units <.21cm,.21cm>
     \multiput{} at 0 0  4 6  /
     \setsolid
     \plot 0 0  0 6  1 6  1 0  0 0 /
     \plot 0 5  2 5  2 1 /
     \plot 0 4  3 4  3 1  0 1 /
     \plot 0 3  4 3  4 2  0 2 /
     \multiput{$\sss\bullet$} at .5 2.5  1.5 2.5  2.5 2.5  3.5 2.5  2.5 1.5 /
     \plot .5 2.5  3.5 2.5 /
     \endpicture} at 0 0
   \put{\beginpicture
     \setcoordinatesystem units <.21cm,.21cm>
     \multiput{} at 0 0  4 7  /
     \setsolid
     \plot 0 1  0 7  1 7  1 1  0 1 /
     \plot 2 0  2 6  3 6  3 0  2 0 /
     \plot 2 1  4 1  4 5  2 5 /
     \plot 0 2  4 2 /
     \plot 0 3  4 3 /
     \plot 0 4  4 4 /
     \plot 0 5  1 5 /
     \plot 0 6  1 6 /
     \multiput{$\sss\bullet$} at .5 2.5  1.5 2.5  1.5 3.5  2.5 3.5  3.5 3.5  /
     \plot .5 2.5  1.5 2.5 /
     \plot 1.5 3.5  3.5 3.5 /
     \endpicture} at 1.5 .5 
   \put{\beginpicture
     \setcoordinatesystem units <.21cm,.21cm>
     \multiput{} at 0 0  4 6  /
     \setsolid
     \plot 1 0  1 6  2 6  2 0  1 0 /
     \plot 1 1  3 1  3 5  1 5 /
     \plot 3 4  0 4  0 3  4 3  4 2  1 2 /
     \multiput{$\sss\bullet$} at .5 3.5  1.5 3.5  2.5 3.5  2.5 2.5 3.5 2.5 /
     \plot 2.5 2.5  3.5 2.5 /
     \plot .5 3.5  2.5 3.5 /
     \endpicture} at -1.5 .5
   \put{\beginpicture
     \setcoordinatesystem units <.21cm,.21cm>
     \multiput{} at 0 0  4 6  /
     \setsolid
     \plot 0 0  0 6  1 6  1 0  0 0 /
     \plot 0 1  2 1  2 5 /
     \plot 0 2  3 2  3 5  0 5 /
     \plot 0 3  4 3  4 4  0 4 /
     \multiput{$\sss\bullet$} at .5 4.5  1.5 4.5  2.5 4.5  1.5 3.5  2.5 3.5 3.5 3.5 /
     \plot 2.5 3.5  3.5 3.5 /
     \plot .5 4.5  2.5 4.5 /
     \endpicture} at -1.5 1.5 
   \put{\beginpicture
     \setcoordinatesystem units <.21cm,.21cm>
     \multiput{} at 0 0  4 7  /
     \setsolid
     \plot 0 1  0 7  1 7  1 1  0 1 /
     \plot 0 2  2 2  2 6  0 6 /
     \plot 3 0  3 6  4 6  4 0  3 0 /
     \plot 3 1  4 1 /
     \plot 3 2  4 2 /
     \plot 0 3  4 3 /
     \plot 0 4  4 4 /
     \plot 0 5  4 5 /
     \multiput{$\sss\bullet$} at .5 5.5  1.5 5.5  1.5 3.5  2.5 3.5  2.5 4.5  3.5 4.5 /
     \plot 1.5 3.5  2.5 3.5 /
     \plot .5 5.5  1.5 5.5 /
     \plot 2.5 4.5  3.5 4.5 /
     \endpicture} at 0 2 
   \put{\beginpicture
     \setcoordinatesystem units <.21cm,.21cm>
     \multiput{} at 0 0  4 7  /
     \setsolid
     \plot 0 0  0 6  1 6  1 0  0 0 /
     \plot 3 1  3 7  4 7  4 1  2 1  2 6  4 6 /
     \plot 0 5  4 5 /
     \plot 0 4  4 4 /
     \plot 0 3  4 3 /
     \plot 0 2  4 2 /
     \plot 0 1  1 1 /
     \multiput{$\sss\bullet$} at .5 4.5  1.5 4.5  1.5 3.5  2.5 3.5  3.5 3.5  2.5 2.5 /
     \plot 1.5 3.5  3.5 3.5 /
     \plot .5 4.5  1.5 4.5 /
     \endpicture} at 1.5 1.5 
\endpicture}$$

\medskip
The objects in position $\pr=(\frac32,\frac32), (\frac32,3), (3,\frac32)$ occur in the
principal component, the remaining ones in the fourth layer in the 6-tube of 
rationality index $1$.

	\bigskip
{\bf (e) The pickets and bipickets in $\Cal S(6)$ 
in the fundamental region $\mathbb F$.}
We denote by $\mathbb F$ the set of all elements $(p,r)$ in $\mathbb T(6)$
with $p \le r \le 6-2p.$ 
This is a fundamental region for the action of $\Sigma_3$ on $\mathbb T(6)$. 
In the second picture, we restrict the attention to the elements $(p,r)$
with $p \ge 1.$
$$
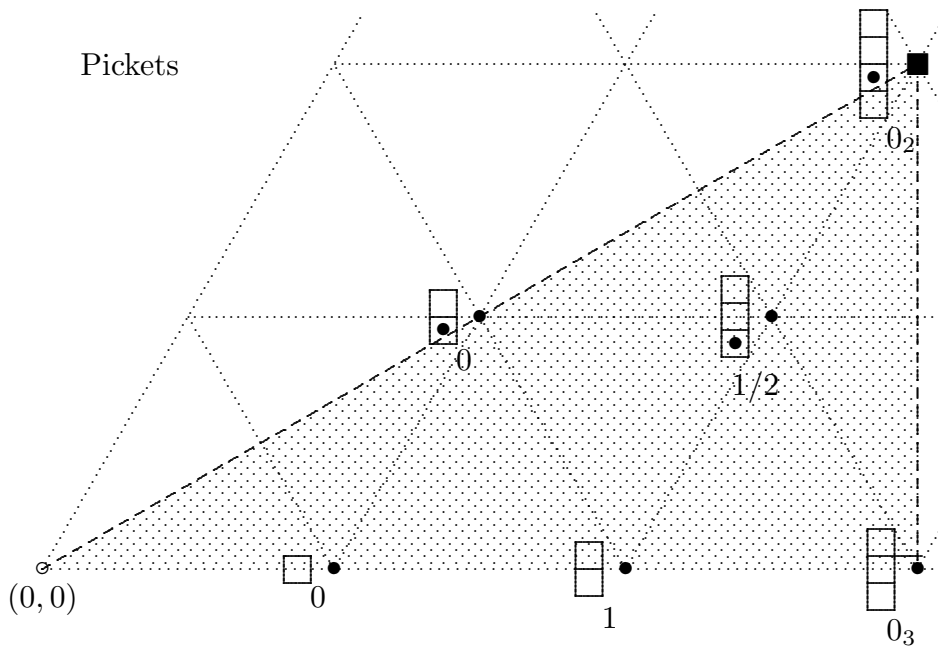
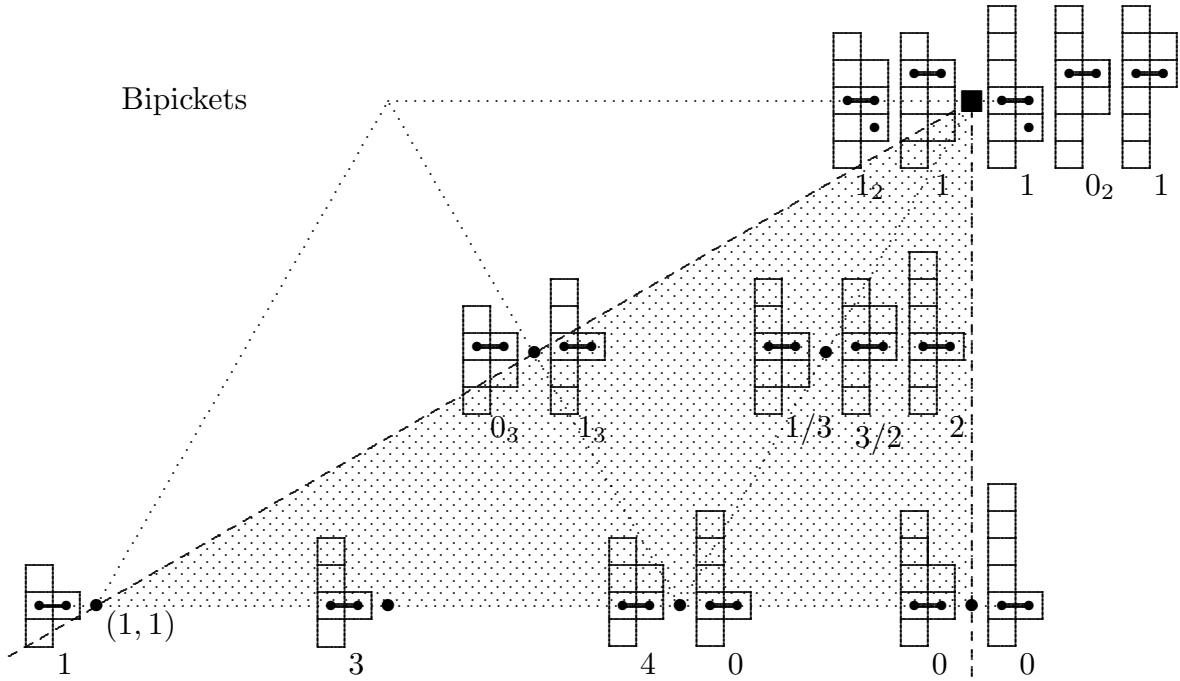
{\beginpicture
   \setcoordinatesystem units <1.6162cm,2.8cm>
\put{Pickets} at -5.4 2 
\multiput{}  at -6 0  0.9 2.2 /
\setdots <1mm> 
\plot -6 0  0 0  0 2  -6 0 /
\setdots <.8mm> 
\plot 0 2  -2 0  -3 1  -4 0  -5 1  0 1 /
\plot 0 0  -2 2  -3 1  -4 2  -3.8 2.2 /
\plot -6 0  -4 2   0 2 /
\plot -2.2 2.2  -2 2  -1.8 2.2 /
\plot -.2 2.2  0 2  .2 2.2 /
\plot  0 2  .2 1.8 /
\plot  .2 0  0 0  .2 .2 /
\plot  .2 1  0 1 /
\setdashes <1mm>
\plot -6 0  0 2  0 0 /
\multiput{$\bullet$} at   -4 0  -2 0  0 0   -3 1  -1 1   /
\multiput{$\circ$} at -6 0  /
\multiput{$\blacksquare$} at 0 2 /
\setshadegrid span <.7mm>
\vshade -6 0 0 <z,z,,> 0 0 2  /

\setsolid 
\put{\beginpicture  
   \setcoordinatesystem units <.3cm,.3cm>
  \multiput{} at 0 0  1 1  /
  \plot 0 0  1 0  1 1  0 1  0 0 /
\endpicture} at  -4.25 0

\put{\beginpicture  
   \setcoordinatesystem units <.3cm,.3cm>
  \multiput{} at 0 0  1 2  /
  \plot 0 0  1 0  1 2  0 2  0 0 /
  \plot 0 1  1 1 /
\endpicture} at  -2.25 0

\put{\beginpicture  
   \setcoordinatesystem units <.3cm,.3cm>
  \multiput{} at 0 0  1 3  /
  \plot 0 0  1 0  1 3  0 3  0 0 /
  \plot 0 1  1 1 /
  \plot 0 2  2 2 /
\endpicture} at  -.25 0

\put{\beginpicture  
   \setcoordinatesystem units <.3cm,.3cm>
  \multiput{} at 0 0  1 2  /
  \plot 0 0  1 0  1 2  0 2  0 0 /
  \plot 0 1  1 1 /
  \put{$\bullet$} at .5 0.5 
\endpicture} at  -3.25 1

\put{\beginpicture  
   \setcoordinatesystem units <.3cm,.3cm>
  \multiput{} at 0 0  1 3  /
  \plot 0 0  1 0  1 3  0 3  0 0 /
  \plot 0 1  1 1 /
  \plot 0 2  1 2 /
  \put{$\bullet$} at .5 0.5 
\endpicture} at  -1.25 1

\put{\beginpicture  
   \setcoordinatesystem units <.3cm,.3cm>
  \multiput{} at 0 0  1 4  /
  \plot 0 0  1 0  1 4  0 4  0 0 /
  \plot 0 1  1 1 /
  \plot 0 2  1 2 /
  \plot 0 3  1 3 /
  \put{$\bullet$} at .5 1.5 
\endpicture} at  -.3 2

\put{\strut$0$} at -4.11 -.12
\put{\strut$1$} at -2.11 -.2
\put{\strut$0_3$} at -.11 -.25
\put{\strut$0$} at -3.11 .82
\put{\strut$1/2$} at -1.11 .72
\put{\strut$0_2$} at -.11 1.7

\put{$(0,0)$} at -6 -.12
\endpicture}
$$

$$
{\beginpicture
   \setcoordinatesystem units <2.8857cm,5cm>

\put{Bipickets} at -2.7 2 
  \multiput{} at -3.3  2  0.2  0.8 /
  \put{$(1,1)$} at -2.85 0.96
\setdots <1mm> 
\plot -3.2 1  0.2 1 /
\plot -3  1  -2 2  -1 1  0 2 /
\plot -2 2  0.2 2 /

\setdashes <1mm>
  \setdashes <1mm>
  \plot -3.3 0.9  0 2  0 0.85 /
\multiput{$\bullet$} at -3 1  -2 1  -1 1  0 1  -1.5 1.5  -.5 1.5 /
\multiput{$\blacksquare$} at 0 2 /
\setshadegrid span <.7mm>
\vshade -3 1 1 <z,z,,> 0 1 2  /

\setsolid 
\put{\beginpicture  
   \setcoordinatesystem units <.3cm,.3cm>
  \multiput{} at 0 0  2 3 /
  \plot 0 0  1 0  1 3  0 3  0 0 /
  \plot 0 1  2 1  2 2  0 2 /
  \multiput{$\ssize \bullet$} at 0.5 1.5  1.5 1.5 /
  \plot  0.5 1.5  1.5 1.5 /
  \plot  0.5 1.45  1.5 1.45 /
  \plot  0.5 1.55  1.5 1.55 /
\endpicture} at  -3.15 1
\put{\beginpicture   
   \setcoordinatesystem units <.3cm,.3cm>
  \multiput{} at 0 0  2 3 /
  \plot 0 0  1 0  1 4  0 4  0 0 / 
  \plot 0 1  2 1  2 2  0 2 /
  \plot 0 3  1 3 /
  \multiput{$\ssize \bullet$} at 0.5 1.5  1.5 1.5 /
  \plot  0.5 1.5  1.5 1.5 /
  \plot  0.5 1.45  1.5 1.45 /
  \plot  0.5 1.55  1.5 1.55 /
\endpicture} at  -2.15 1
\put{\beginpicture  
   \setcoordinatesystem units <.3cm,.3cm>
  \multiput{} at 0 0  2 3 /
  \plot 0 0  1 0  1 4  0 4  0 0 / 
  \plot 0 1  2 1  2 2  0 2 /
  \plot 0 3  1 3 /
  \plot 2 2  2 3  1 3 /
  \multiput{$\ssize \bullet$} at 0.5 1.5  1.5 1.5 /
  \plot  0.5 1.5  1.5 1.5 /
  \plot  0.5 1.45  1.5 1.45 /
  \plot  0.5 1.55  1.5 1.55 /
\endpicture} at  -1.15 1
\put{\beginpicture
   \setcoordinatesystem units <.3cm,.3cm>
  \multiput{} at 0 0  2 3 /
  \plot 0 0  1 0  1 4  0 4  0 0 / 
  \plot 0 1  2 1  2 2  0 2 /
  \plot 0 3  1 3 /
  \plot 0 4  0 5  1 5  1 4 /
  \multiput{$\ssize \bullet$} at 0.5 1.5  1.5 1.5 /
  \plot  0.5 1.5  1.5 1.5 /
  \plot  0.5 1.45  1.5 1.45 /
  \plot  0.5 1.55  1.5 1.55 /
\endpicture} at  -0.85 1

\put{\beginpicture  
   \setcoordinatesystem units <.3cm,.3cm>
  \multiput{} at 0 0  2 3 /
  \plot 0 0  1 0  1 4  0 4  0 0 / 
  \plot 0 1  2 1  2 2  0 2 /
  \plot 0 3  2 3  2 2 /
  \plot 0 4  0 5  1 5  1 3  /
  \multiput{$\ssize \bullet$} at 0.5 1.5  1.5 1.5 /
  \plot  0.5 1.5  1.5 1.5 /
  \plot  0.5 1.45  1.5 1.45 /
  \plot  0.5 1.55  1.5 1.55 /
\endpicture} at  -.15 1
\put{\beginpicture
   \setcoordinatesystem units <.3cm,.3cm>
  \multiput{} at 0 0  2 3 /
  \plot 0 0  1 0  1 4  0 4  0 0 / 
  \plot 0 1  2 1  2 2  0 2 /
  \plot 0 3  1 3 /
  \plot 0 4  0 6  1 6  1 4   /
  \plot 0 5  1 5 /
  \multiput{$\ssize \bullet$} at 0.5 1.5  1.5 1.5 /
  \plot  0.5 1.5  1.5 1.5 /
  \plot  0.5 1.45  1.5 1.45 /
  \plot  0.5 1.55  1.5 1.55 /
\endpicture} at  .15 1

\put{\beginpicture  
   \setcoordinatesystem units <.3cm,.3cm>
  \multiput{} at 0 0  2 4.5 /
  \plot 0 0  1 0  1 4  0 4  0 0 / 
  \plot 0 1  2 1  2 2  0 2 /
  \plot 0 3  2 3  2 2 /
  \multiput{$\ssize \bullet$} at 0.5 2.5  1.5 2.5 /
  \plot  0.5 2.5  1.5 2.5 /
  \plot  0.5 2.45  1.5 2.45 /
  \plot  0.5 2.55  1.5 2.55 /
\endpicture} at  -1.65 1.5
\put{\beginpicture
   \setcoordinatesystem units <.3cm,.3cm>
  \multiput{} at 0 0  2 4.5 /
  \plot 0 0  1 0  1 5  0 5  0 0 / 
  \plot 0 2  2 2  2 3  0 3 /
  \plot 0 1  1 1 /
  \plot 0 4  1 4 /
  \multiput{$\ssize \bullet$} at 0.5 2.5  1.5 2.5 /
  \plot  0.5 2.5  1.5 2.5 /
  \plot  0.5 2.45  1.5 2.45 /
  \plot  0.5 2.55  1.5 2.55 /
\endpicture} at  -1.35 1.5

\put{\beginpicture  
   \setcoordinatesystem units <.3cm,.3cm>
  \multiput{} at 0 0  2 4.5 /
  \plot 0 0  1 0  1 4  0 4  0 0 / 
  \plot 0 1  2 1  2 2  0 2 /
  \plot 0 3  2 3  2 2 /
  \plot 0 4  0 5  1 5  1 4 /
  \multiput{$\ssize \bullet$} at 0.5 2.5  1.5 2.5 /
  \plot  0.5 2.5  1.5 2.5 /
  \plot  0.5 2.45  1.5 2.45 /
  \plot  0.5 2.55  1.5 2.55 /
\endpicture} at  -.65 1.5
\put{\beginpicture
   \setcoordinatesystem units <.3cm,.3cm>
  \multiput{} at 0 0  2 4.5 /
  \plot 0 0  1 0  1 5  0 5  0 0 / 
  \plot 0 2  2 2  2 3  0 3 /
  \plot 0 1  1 1 /
  \plot 0 4  1 4 /
  \plot 1 4  2 4  2 3 /
  \multiput{$\ssize \bullet$} at 0.5 2.5  1.5 2.5 /
  \plot  0.5 2.5  1.5 2.5 /
  \plot  0.5 2.45  1.5 2.45 /
  \plot  0.5 2.55  1.5 2.55 /
\endpicture} at  -.35 1.5
\put{\beginpicture
   \setcoordinatesystem units <.3cm,.3cm>
  \multiput{} at 0 0  2 4.5 /
  \plot 0 0  1 0  1 5  0 5  0 0 / 
  \plot 0 2  2 2  2 3  0 3 /
  \plot 0 1  1 1 /
  \plot 0 4  1 4 /
  \plot 0 5  0 6  1 6  1 5 /
  \multiput{$\ssize \bullet$} at 0.5 2.5  1.5 2.5 /
  \plot  0.5 2.5  1.5 2.5 /
  \plot  0.5 2.45  1.5 2.45 /
  \plot  0.5 2.55  1.5 2.55 /
\endpicture} at  -.12 1.5

\put{\beginpicture  
   \setcoordinatesystem units <.3cm,.3cm>
  \multiput{} at 0 0  2 5 /
  \plot 0 0  1 0  1 5  0 5  0 0 / 
  \plot 0 2  2 2  2 3  0 3 /
  \plot 0 1  1 1 /
  \plot 0 4  1 4 /
  \plot 1 4  2 4  2 3 /
  \plot 1 1  2 1  2 2 /
  \multiput{$\ssize \bullet$} at 0.5 2.5  1.5 2.5  1.5 1.5 /
  \plot  0.5 2.5  1.5 2.5 /
  \plot  0.5 2.45  1.5 2.45 /
  \plot  0.5 2.55  1.5 2.55 /
\endpicture} at  -.38 2
\put{\beginpicture
   \setcoordinatesystem units <.3cm,.3cm>
  \multiput{} at 0 0  2 5 /
  \plot 0 0  1 0  1 5  0 5  0 0 / 
  \plot 0 2  2 2  2 3  0 3 /
  \plot 0 1  1 1 /
  \plot 0 4  1 4 /
  \plot 1 4  2 4   2 3 /
  \plot 1 1  2 1  2 2 /
  \multiput{$\ssize \bullet$} at 0.5 3.5  1.5 3.5 /
  \plot  0.5 3.5  1.5 3.5 /
  \plot  0.5 3.45  1.5 3.45 /
  \plot  0.5 3.55  1.5 3.55 /
\endpicture} at  -.15 2

\put{\beginpicture
   \setcoordinatesystem units <.3cm,.3cm>
  \multiput{} at 0 0  2 5 /
  \plot 0 0  1 0  1 5  0 5  0 0 / 
  \plot 0 2  2 2  2 3  0 3 /
  \plot 0 1  1 1 /
  \plot 0 4  1 4 /
  \plot 1 1  2 1  2 2 /
  \plot 0 5  0 6  1 6  1 5 /
  \multiput{$\ssize \bullet$} at 0.5 2.5  1.5 2.5  1.5 1.5 /
  \plot  0.5 2.5  1.5 2.5 /
  \plot  0.5 2.45  1.5 2.45 /
  \plot  0.5 2.55  1.5 2.55 /
\endpicture} at  .15 2
\put{\beginpicture
   \setcoordinatesystem units <.3cm,.3cm>
  \multiput{} at 0 0  2 5 /
  \plot 0 0  1 0  1 5  0 5  0 0 / 
  \plot 0 1  1 1 /
  \plot 0 3  2 3  2 4  /
  \plot 0 4  1 4 /
  \plot 1 4  2 4   2 3 /
  \plot 0 5  0 6  1 6  1 5 /
  \plot 0 2  2 2  2 3 /
  \multiput{$\ssize \bullet$} at 0.5 3.5  1.5 3.5 /
  \plot  0.5 3.5  1.5 3.5 /
  \plot  0.5 3.45  1.5 3.45 /
  \plot  0.5 3.55  1.5 3.55 /
\endpicture} at  .38 2
\put{\beginpicture
   \setcoordinatesystem units <.3cm,.3cm>
  \multiput{} at 0 0  2 5 /
  \plot 0 0  1 0  1 5  0 5  0 0 / 
  \plot 0 2  1 2  /
  \plot 0 1  1 1 /
  \plot 0 4  1 4 /
  \plot 1 4  2 4   2 3 /
  \plot 0 5  0 6  1 6  1 5 /
  \plot 0 3  2 3 /
  \plot 2 4  2 5  1 5 /
  \multiput{$\ssize \bullet$} at 0.5 3.5  1.5 3.5 /
  \plot  0.5 3.5  1.5 3.5 /
  \plot  0.5 3.45  1.5 3.45 /
  \plot  0.5 3.55  1.5 3.55 /
\endpicture} at  .61 2

\put{\strut$1$} at -3.11 0.88
\put{\strut$3$} at -2.11 0.88
\put{\strut$4$} at -1.11 0.88
\put{\strut$0$} at -.81 0.88
\put{\strut$0$} at -.11 0.88
\put{\strut$0$} at .19 0.88

\put{\strut$0_3$} at -1.6 1.35
\put{\strut$1_3$} at -1.3 1.35
\put{\strut$1/3$} at -.56 1.35
\put{\strut$3/2$} at -.32 1.33
\put{\strut$2$} at -.05 1.35

\put{\strut$1_2$} at -.35 1.83
\put{\strut$1$} at -.1 1.83
\put{\strut$1$} at .19 1.83
\put{\strut$0_2$} at .44 1.83
\put{\strut$1$} at .65 1.83
\endpicture}
$$
The number mentioned below an object $X$ is its rationality index $\gamma$,
this is a number in $\mathbb Q_0^+$ (see Section~\ref{sec-eleven-one}).
Usually, the object $X$ belongs to a tube which has rank $6$.
In case $X$ belongs to a tube $\Cal C$ with rank $2$ or $3$ (this happens only for $\gamma = 0$ 
and $1$),  we write $\gamma_r$ with $r \in \{2,3\}$
in order to indicate that $\Cal C$ has rationality index $\gamma$
and rank $r$. Note that the pickets and bipickets labeled $0$ are those which
belong to $\Cal P(6)$, the remaining pickets and bipickets belong to stable tubes. 
Altogether we see: {\it There are $14$ Auslander-Reiten
components in $\Cal S(6)$ which contain pickets or bipickets, namely the tubes of rank $6,\ 3,\ 2$ with
rationality index $0$ or $1$, and the tubes of rank $6$ with rationality index}
$$
  \frac14,\ \ \frac13,\ \ \frac12,\ \ \frac23,\qquad \frac32,\ \ 2,\ \ 3,\ \ 4.
$$
Most pickets and bipickets which occur in stable tubes
have quasi-length 1, the only exception are the bipickets in the tube of rank 6
and rationality index $1/1$ which have quasi-length 2.
The pickets and bipickets in $\Cal P(n)$ have quasi-length at most 4.

The tubes with rank $6,\ 3,\ 2$  and rationality index $0$ have been presented in detail
in \cite{RS1}, for the remaining tubes containing pickets or bipickets, see \cite{S2}. 

\medskip
\begin{remark}
  There are bipickets in $\Cal S(6)$ with $p = r$ which
  are not self-dual, namely the bipickets $X = ([4],[6,2],[3,1])$ and $\D X = ([3,1],[6,2],[4])$
  (which are dual to each other).
\end{remark}

\medskip
In contrast, {\it for $n \le 5$, all indecomposable objects $X$ with $p = r$ (equivalently,
  with $u = w$) are self-dual.} Namely, according to \cite[(6.5)]{RS1},
all but one indecomposable
objects in $\Cal S(5)$ with $u=w$  are pickets or bipickets. The exception is of course
self-dual. Pickets with $u = w$ are self-dual. And a bipicket $X$ in $\Cal S(5)$ 
satisfies $q \le 4$, thus if $u = v$, then $X$ belongs to the fundamental region $\mathbb F$
of $\mathbb T(6)$ presented above: The 5  bipickets of this kind are all self-dual. 

\bigskip
{\bf (f) The tetrapickets on the triangle $\Delta_{5/4}$:}

\medskip
All tetrapickets in $\Cal S(6)$ occur inside the hexagon pictured below, each on an intersection point
of dotted lines (thus, both pr-coordinates are integer multiples of $\frac14$). 
In Section 15.2 (d) we have already pictured the three tetrapickets on the boundary
of $\Delta_1$.

\medskip
Here we show the tetrapickets on the boundary of $\Delta_{5/4}$, their positions are as indicated.
There are no other objects on the boundary of $\Delta_{5/4}$,
and the only objects outside $\Delta_{5/4}$ are the pickets on $\Delta_0$
and the objects on $\Delta_1$ shown in Section \ref{sec-fifteen-two},  (a)--(d).

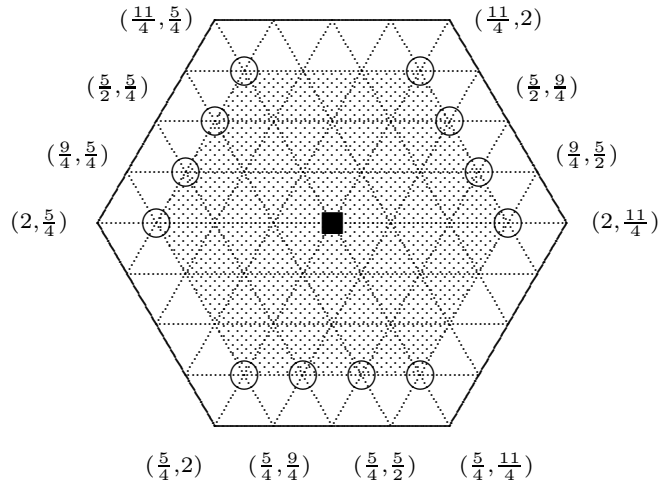
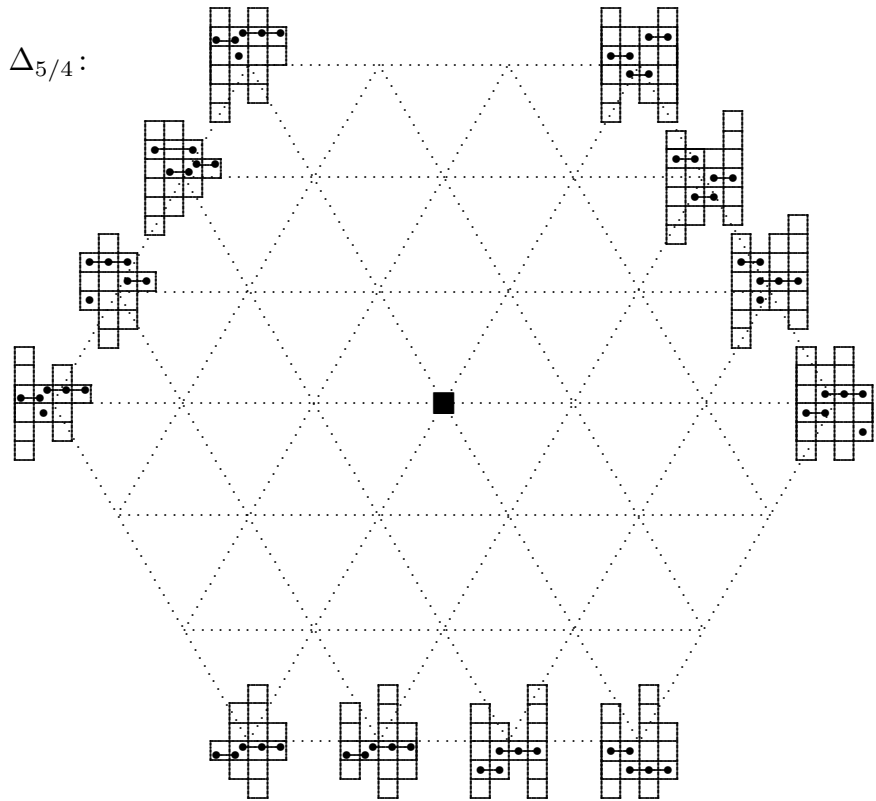
\medskip
$$
{\beginpicture
   \setcoordinatesystem units <1.299cm, 2.25cm>
   \multiput{} at -2 0  2 2 /
   \plot -2 1  -1 2  1 2  2 1  1 0  -1 0  -2 1 /
   \multiput{$\bigcirc$} at   -.75 .25  -.25 .25  .25 .25  .75 .25
                         -1.5 1  -1.25 1.25  -1 1.5  -.75 1.75
                          1.5 1  1.25 1.25  1 1.5  .75 1.75 /
     \put{$\blacksquare$} at 0 1
   \setdots<.5mm>
   \plot -2 1  2 1 /
   \plot -1 0  1 2 /
   \plot -1 2  1 0 /
   \plot 0 0  -1.5 1.5  1.5 1.5  0 0 /
   \plot -.5 0  -1.75 1.25  1.75 1.25  .5 0  -1.25 1.75  1.25 1.75  -.5 0 /
   \plot -1.25 .25  .5 2  1.75 .75  -1.75 .75  -.5 2   1.25 .25   -1.25 .25 /
   \plot -1.5 .5  0 2  1.5 .5  -1.5 .5 /
   \setshadegrid span <.5mm>
   \vshade -1.5 1 1  <,z,,>  -.75 .25 1.75  <z,z,,>  .75 .25 1.75  <z,z,,> 1.5 1 1 /
   \put{$\ss (\frac54,2)$} at -1.35 -.2
   \put{$\ss (\frac54,\frac{11}4)$} at 1.4 -.2
   \put{$\ss (2,\frac54)$} at -2.5 1
   \put{$\ss (\frac{11}4,\frac54)$} at -1.5 2
   \put{$\ss (2,\frac{11}4)$} at 2.5 1
   \put{$\ss (\frac{11}4,2)$} at 1.5 2
   \put{$\ss (\frac54,\frac94)$} at -.45 -.2
   \put{$\ss (\frac54,\frac52)$} at .45 -.2
   \put{$\ss (\frac94,\frac52)$} at 2.17 1.3333
   \put{$\ss (\frac52,\frac94)$} at 1.83 1.6667
   \put{$\ss (\frac52,\frac54)$} at -1.83 1.6667
   \put{$\ss (\frac94,\frac54)$} at -2.17 1.3333
\endpicture}
$$
\medskip
$$
{\beginpicture
   \setcoordinatesystem units <2.165cm,3.75cm>
   \multiput{} at -2.333 0  2.3333 2.3333 /
   \setdots<1mm>
   \plot -1 0  1 0 /
   \plot -2 1  -1 2 /
   \plot 2 1  1 2 /
   \put{$\blacksquare$} at 0 1
   \put{$\Delta_{5/4}\:$} at -2 2 
   \plot -1 0  1 2  -1 2  1 0  2 1  -2 1 -1 0 /
   \plot -.33 0  -1.67 1.33  1.67 1.33  .33 0  -1.33 1.67  1.33 1.67 -.33 0 /
   \plot -.33 2  -1.67 .67  1.67 .67  .33 2  -1.33 .33  1.33 .33  -.33 2 /
   \put{\beginpicture
     \setcoordinatesystem units <.21cm,.21cm>
     \multiput{} at 0 0  4 6  /
     \setsolid
     \plot 1 4  4 4  4 2  0 2  0 3  4 3  /
     \plot 3 5  1 5  1 1  3 1 /
     \plot 2 0  2 6  3 6  3 0  2 0 /
     \multiput{$\sss\bullet$} at .3 2.3  1.3 2.3  1.7 2.7  2.7 2.7  3.7 2.7 /
     \plot .3 2.3  1.3 2.3 /
     \plot 1.7 2.7  3.7 2.7 /
     \endpicture} at -1 0  
   \put{\beginpicture
     \setcoordinatesystem units <.21cm,.21cm>
     \multiput{} at 0 0  4 6  /
     \setsolid
     \plot 0 1  0 5  1 5  1 1  0 1 /
     \plot 0 2  4 2  4 3  0 3 /
     \plot 0 4  1 4 /
     \plot 2 0  2 6  3 6  3 0  2 0 /
     \plot 2 1  4 1  4 4  2 4 /
     \plot 2 5  3 5 /
     \multiput{$\sss\bullet$} at .3 2.3  1.3 2.3  1.7 2.7  2.7 2.7  3.7 2.7 /
     \plot .3 2.3  1.3 2.3 /
     \plot 1.7 2.7  3.7 2.7 /
     \endpicture} at -.3333 0  
   \put{\beginpicture
     \setcoordinatesystem units <.21cm,.21cm>
     \multiput{} at 0 0  4 6  /
     \setsolid
     \plot 0 0  0 5  1 5  1 0  0 0 /
     \plot 0 1  2 1  2 4  0 4 /
     \plot 0 2  4 2  /
     \plot 0 3  4 3 /
     \plot 3 0  3 6  4 6  4 0  3 0 /
     \plot 3 1  4 1 /
     \plot 3 4  4 4 /
     \plot 3 5  4 5 /
     \multiput{$\sss\bullet$} at .5 1.5  1.5 1.5  1.5 2.5  2.5 2.5  3.5 2.5 /
     \plot .5 1.5  1.5 1.5 /
     \plot 1.5 2.5  3.5 2.5 /
     \endpicture} at .3333 0  
   \put{\beginpicture
     \setcoordinatesystem units <.21cm,.21cm>
     \multiput{} at 0 0  4 6  /
     \setsolid
     \plot 0 0  0 5  1 5  1 0  0 0 /
     \plot 2 4  4 4  4 1  0 1  0 3  4 3 /
     \plot 2 0  2 6  3 6  3 0  2 0 /
     \plot 0 2  4 2 /
     \plot 0 4 1 4 /
     \plot 2 5 3 5 /
     \multiput{$\sss\bullet$} at .5 2.5  1.5 2.5  1.5 1.5  2.5 1.5  3.5 1.5 /
     \plot .5 2.5  1.5 2.5 /
     \plot 1.5 1.5  3.5 1.5 /
     \endpicture} at 1 0
   \put{\beginpicture
     \setcoordinatesystem units <.21cm,.21cm>
     \multiput{} at 0 0  4 6  /
     \setsolid
     \plot 0 0  0 6  1 6  1 0  0 0 /
     \plot 4 2  0 2  0 5  4 5  4 1  2 1  2 5 /
     \plot 3 0  3 6  4 6  4 0  3 0 /
     \plot 0 3  4 3  4 4  0 4 /
     \plot 0 1  1 1 /
     \multiput{$\sss\bullet$} at .5 3.5  1.5 3.5  1.5 2.5  2.5 2.5  2.5 4.5  3.5 4.5 /
     \plot .5 3.5  1.5 3.5 /
     \plot 1.5 2.5  2.5 2.5 /
     \plot 2.5 4.5  3.5 4.5 /
     \endpicture} at 1 2
   \put{\beginpicture
     \setcoordinatesystem units <.21cm,.21cm>
     \multiput{} at 0 0  4 7  /
     \setsolid
     \plot 0 0  0 6  1 6  1 0  0 0 /
     \plot 0 1  2 1  2 5  0 5 /
     \plot 0 2  4 2 /
     \plot 0 3  4 3 /
     \plot 0 4  4 4 /
     \plot 0 5  4 5 /
     \plot 3 6  4 6 /
     \plot 3 1  3 7  4 7  4 1  3 1 /
     \multiput{$\sss\bullet$} at .5 4.5  1.5 4.5  1.5 2.5  2.5 2.5  2.5 3.5  3.5 3.5 /
     \plot .5 4.5  1.5 4.5 /
     \plot 1.5 2.5  2.5 2.5 /
     \plot 2.5 3.5  3.5 3.5 /
     \endpicture} at 1.3333 1.6667
   \put{\beginpicture
     \setcoordinatesystem units <.21cm,.21cm>
     \multiput{} at 0 0  4 6  /
     \setsolid
     \plot 0 0  0 6  1 6  1 0  0 0 /
     \plot 0 1  1 1 /
     \plot 0 2  4 2 /
     \plot 0 3  4 3 /
     \plot 0 4  4 4 /
     \plot 0 5  4 5 /
     \plot 2 2  2 6  4 6 /
     \plot 3 1  3 7  4 7  4 1  3 1 /
     \multiput{$\sss\bullet$} at .5 4.5  1.5 4.5  1.5 2.5  1.5 3.5  2.5 3.5  3.5 3.5 /
     \plot .5 4.5  1.5 4.5 /
     \plot 1.5 3.5  3.5 3.5 /
     \endpicture} at 1.6667 1.3333
   \put{\beginpicture
     \setcoordinatesystem units <.21cm,.21cm>
     \multiput{} at 0 0  4 6  /
     \setsolid
     \plot 0 0  0 6  1 6  1 0  0 0 /
     \plot 0 1  1 1 /
     \plot 0 3  4 3  4 5  0 5  0 2  3 2 /
     \plot 2 1  2 6  3 6  3 1  2 1 /
     \plot 0 4  4 4 /
     \multiput{$\sss\bullet$} at .3 4.3  1.3 4.3  1.5 3.5  1.7 4.7  2.7 4.7  3.7 4.7 /
     \plot .3 4.3  1.3 4.3 /
     \plot 1.7 4.7  3.7 4.7 /
     \endpicture} at -1 2
   \put{\beginpicture
     \setcoordinatesystem units <.21cm,.21cm>
     \multiput{} at 0 0  4 6  /
     \setsolid
     \plot 1 6  1 0  0 0  0 6  2 6  2 1  0 1 /
     \plot 0 2  3 2  3 5  0 5 /
     \plot 0 3  4 3  4 4  0 4 /
     \multiput{$\sss\bullet$} at .5 4.5  2.5 4.5  1.3 3.3  2.3 3.3  2.7 3.7  3.7 3.7 /
     \plot .5 4.5  2.5 4.5 /
     \plot 1.3 3.3  2.3 3.3 /
     \plot 2.7 3.7  3.7 3.7 /
     \endpicture} at -1.3333 1.6667
   \put{\beginpicture
     \setcoordinatesystem units <.21cm,.21cm>
     \multiput{} at 0 0  4 6  /
     \setsolid
     \plot 1 0  1 6  2 6  2 0  1 0 /
     \plot 3 2  0 2  0 5  3 5  3 1  1 1  /
     \plot 0 3  4 3  4 4  0 4 /
     \multiput{$\sss\bullet$} at .5 2.5  .5 4.5  1.5 4.5  2.5 4.5  2.5 3.5  3.5 3.5 /
     \plot .5 4.5  2.5 4.5 /
     \plot 2.5 3.5  3.5 3.5 /
     \endpicture} at -1.6667 1.3333
   \put{\beginpicture
     \setcoordinatesystem units <.21cm,.21cm>
     \multiput{} at 0 0  4 6  /
     \setsolid
     \plot 0 0  0 6  1 6  1 0  0 0 /
     \plot 2 0  2 6  3 6  3 0  2 0 /
     \plot 0 4  4 4  4 1  0 1  0 5  3 5 /
     \plot 0 2  4 2  4 3  0 3 /
     \multiput{$\sss\bullet$} at .5 2.5  1.5 2.5  1.5 3.5  2.5 3.5  3.5 3.5  3.5 1.5 /
     \plot .5 2.5  1.5 2.5 /
     \plot 1.5 3.5  3.5 3.5 /
     \endpicture} at 2 1
   \put{\beginpicture
     \setcoordinatesystem units <.21cm,.21cm>
     \multiput{} at 0 0  4 6  /
     \setsolid
     \plot 0 0  0 6  1 6  1 0  0 0 /
     \plot 0 1  1 1 /
     \plot 3 2  0 2  0 4  4 4  4 3  0 3 /
     \plot 2 1  2 5  3 5  3 1  2 1 /
     \plot 0 5  1 5 /
     \multiput{$\sss\bullet$} at .3 3.3  1.3 3.3  1.5 2.5  1.7 3.7  2.7 3.7  3.7 3.7 /
     \plot .3 3.3  1.3 3.3 /
     \plot 1.7 3.7  3.7 3.7 /
     \endpicture} at -2 1
\endpicture}
$$

\bigskip
The objects in position $(\frac54,2)$ and its
rotations have quasi-length 2 in the 6-tube with rationality index $3$;
those in position $(\frac54,\frac94)$ and its rotations on the mouth of the 6-tube
with rationality index $\frac53$;
the modules at $(\frac54,\frac52)$ and its rotations on the mouth of the 6-tube
with rationality index $\frac35$; and the objects at $(\frac54,\frac{11}4)$ and its rotations
have quasi-length 2 in the 6-tube with rationality index $\frac 13$.

\bigskip

{\bf (g) The line $p = 1$ for $\Cal S(6)$:}
$$
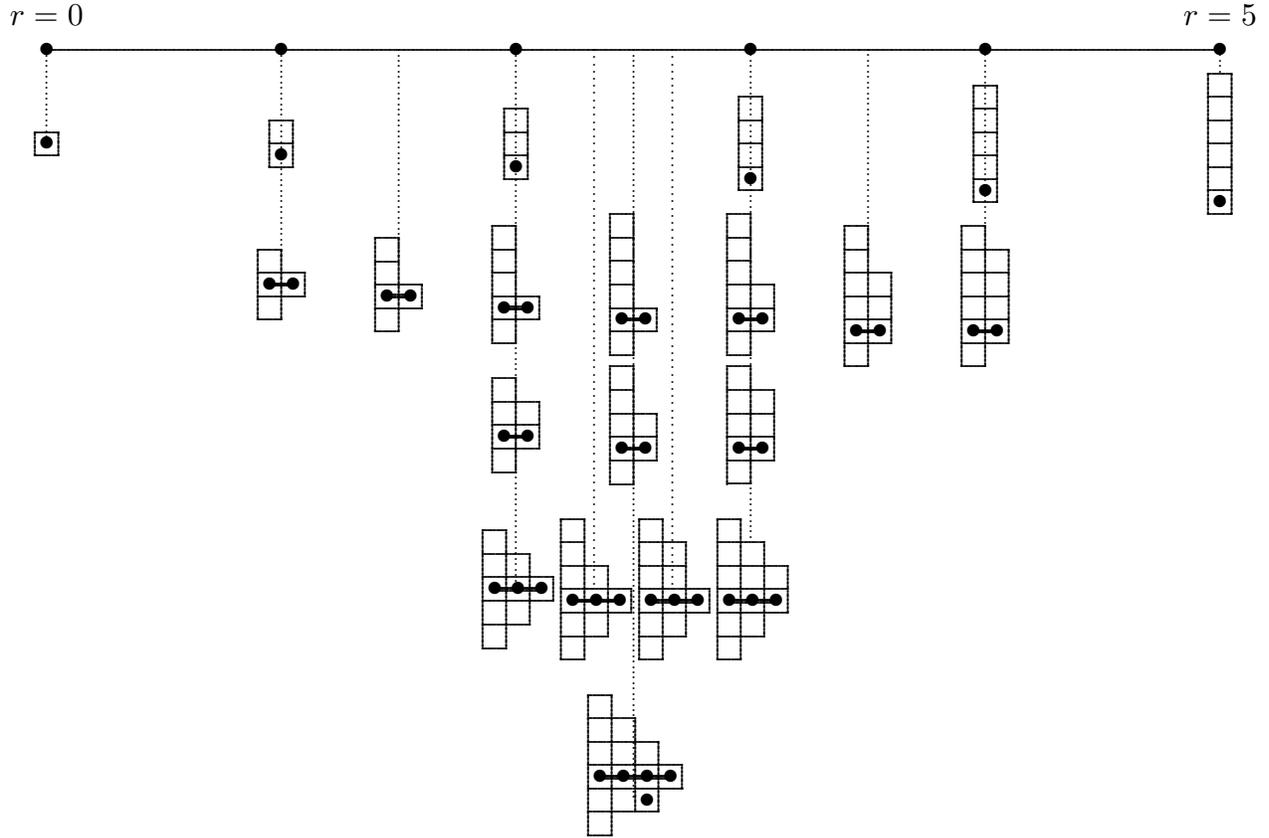
{\beginpicture
    \setcoordinatesystem units <.23cm,.23cm>
\multiput{$\bullet$} at -10 0  0 0  10 0  20 0  30 0  40 0 /
\plot -10 0  40 0 /
\put{$r = 0$} at -10 1.5
\put{$r = 5$} at 40 1.5

\setdots <.6mm>
\plot -10 0  -10 -3 /
\plot 0 0   0 -2.5 /
\plot 0 -6  0 -8 /
\plot 5 0   5 -7.5 /
\plot 10 0  10 -2 /
\plot 10 -6  10 -7 /
\plot 10 -13  10 -13.5 /
\plot 10 -18.7  10 -21 /
\plot 15 0  15 -6.5 /
\plot 15 -19 15 -28 /
\plot 20 0  20 -1.5 /
\plot 20 -6  20 -7 /
\plot 20 -19.2  20 -20.5 /
\plot 25 0  25 -7 /
\plot 30 0   30 -1 /
\plot 30 -7  30 -7.5 /
\plot 40 0   40 -.6 /
\plot 13.333 0  13.333 -21.5 /
\plot 16.667 0  16.667 -20.5 /
\setsolid

\put{\beginpicture
\multiput{} at 0 0  1 1 /
\plot 0 0  1 0  1 1  0 1  0 0  /
\multiput{$\bullet$} at 0.5 0.5   /
\endpicture} at -10 -4

\put{\beginpicture
\multiput{} at 0 0  1 2 /
\plot 0 0  1 0  1 2  0 2  0 0 /
\plot 0 1  1 1 /
\multiput{$\bullet$} at 0.5 0.5   /
\endpicture} at 0 -4

\put{\beginpicture
\multiput{} at 0 0  1 3 /
\plot 0 0  1 0  1 3  0 3  0 0 /
\plot 0 1  1 1 /
\plot 0 2  1 2 /
\multiput{$\bullet$} at 0.5 0.5   /
\endpicture} at 10 -4

\put{\beginpicture
\multiput{} at 0 0  1 4 /
\plot 0 0  1 0  1 4  0 4  0 0 /
\plot 0 1  1 1 /
\plot 0 2  1 2 /
\plot 0 3  1 3 /
\multiput{$\bullet$} at 0.5 0.5   /
\endpicture} at 20 -4

\put{\beginpicture
\multiput{} at 0 0  1 5 /
\plot 0 0  1 0  1 5  0 5  0 0 /
\plot 0 1  1 1 /
\plot 0 2  1 2 /
\plot 0 3  1 3 /
\plot 0 4  1 4 /
\multiput{$\bullet$} at  0.5 0.5   /
\endpicture} at 30 -4

\put{\beginpicture
\multiput{} at 0 0  1 6 /
\plot 0 0  1 0  1 6  0 6  0 0 /
\plot 0 1  1 1 /
\plot 0 2  1 2 /
\plot 0 3  1 3 /
\plot 0 4  1 4 /
\plot 0 5  1 5 /
\multiput{$\bullet$} at 0.5 0.5   /
\endpicture} at 40 -4

\put{\beginpicture
\multiput{} at 0 0  2 3 /
\plot 0 0  1 0  1 3  0 3  0 0 /
\plot 0 1  2 1  2 2  0 2 /
\multiput{$\bullet$} at 0.5 1.5  1.5 1.5   /
\plot 0.5 1.5  1.5 1.5 /
\plot 0.5 1.55  1.5 1.55 /
\plot 0.5 1.45  1.5 1.45 /
\endpicture} at 0 -10

\put{\beginpicture
\multiput{} at 0 0  2 4 /
\plot 0 0  1 0  1 4  0 4  0 0 /
\plot 0 1  2 1  2 2  0 2 /
\plot 0 3  1 3 /
\multiput{$\bullet$} at 0.5 1.5  1.5 1.5   /
\plot 0.5 1.5  1.5 1.5 /
\plot 0.5 1.55  1.5 1.55 /
\plot 0.5 1.45  1.5 1.45 /
\endpicture} at 5 -10

\put{\beginpicture
\multiput{} at 0 0  2 5 /
\plot 0 0  1 0  1 5  0 5  0 0 /
\plot 0 1  2 1  2 2  0 2 /
\plot 0 3  1 3 /
\plot 0 4  1 4 /
\multiput{$\bullet$} at 0.5 1.5  1.5 1.5   /
\plot 0.5 1.5  1.5 1.5 /
\plot 0.5 1.55  1.5 1.55 /
\plot 0.5 1.45  1.5 1.45 /
\endpicture} at 10 -10

\put{\beginpicture
\multiput{} at 0 0  2 6 /
\plot 0 0  1 0  1 6  0 6  0 0 /
\plot 0 1  2 1  2 2  0 2 /
\plot 0 3  1 3 /
\plot 0 4  1 4 /
\plot 0 5  1 5 /
\multiput{$\bullet$} at 0.5 1.5  1.5 1.5   /
\plot 0.5 1.5  1.5 1.5 /
\plot 0.5 1.55  1.5 1.55 /
\plot 0.5 1.45  1.5 1.45 /
\endpicture} at 15 -10 

\put{\beginpicture
\multiput{} at 0 0  2 4 /
\plot 0 0  1 0  1 4  0 4  0 0 /
\plot 0 1  2 1  2 3  0 3 /
\plot 0 2  2 2 /
\multiput{$\bullet$} at 0.5 1.5  1.5 1.5   /
\plot 0.5 1.5  1.5 1.5 /
\plot 0.5 1.55  1.5 1.55 /
\plot 0.5 1.45  1.5 1.45 /
\endpicture} at 10 -16 

\put{\beginpicture
\multiput{} at 0 0  2 5 /
\plot 0 0  1 0  1 5  0 5  0 0 /
\plot 0 1  2 1  2 3  0 3 /
\plot 0 2  2 2 /
\plot 0 4  1 4 /
\multiput{$\bullet$} at 0.5 1.5  1.5 1.5   /
\plot 0.5 1.5  1.5 1.5 /
\plot 0.5 1.55  1.5 1.55 /
\plot 0.5 1.45  1.5 1.45 /
\endpicture} at 15 -16 

\put{\beginpicture
\multiput{} at 0 0  2 6 /
\plot 0 0  1 0  1 6  0 6  0 0 /
\plot 0 1  2 1  2 3  0 3 /
\plot 0 2  2 2 /
\plot 0 4  1 4 /
\plot 0 5  1 5 /
\multiput{$\bullet$} at 0.5 1.5  1.5 1.5   /
\plot 0.5 1.5  1.5 1.5 /
\plot 0.5 1.55  1.5 1.55 /
\plot 0.5 1.45  1.5 1.45 /
\endpicture} at 20 -10 

\put{\beginpicture
\multiput{} at 0 0  2 5 /
\plot 0 0  1 0  1 5  0 5  0 0 /
\plot 0 1  2 1  2 4  0 4 /
\plot 0 2  2 2 /
\plot 0 3  2 3 /
\multiput{$\bullet$} at 0.5 1.5  1.5 1.5   /
\plot 0.5 1.5  1.5 1.5 /
\plot 0.5 1.55  1.5 1.55 /
\plot 0.5 1.45  1.5 1.45 /
\endpicture} at 20 -16

\put{\beginpicture
\multiput{} at 0 0  2 6 /
\plot 0 0  1 0  1 6  0 6  0 0 /
\plot 0 1  2 1  2 4  0 4 /
\plot 0 2  2 2 /
\plot 0 3  2 3 /
\plot 0 5  1 5 /
\multiput{$\bullet$} at 0.5 1.5  1.5 1.5   /
\plot 0.5 1.5  1.5 1.5 /
\plot 0.5 1.55  1.5 1.55 /
\plot 0.5 1.45  1.5 1.45 /
\endpicture} at 25 -10.5

\put{\beginpicture
\multiput{} at 0 0  2 6 /
\plot 0 0  1 0  1 6  0 6  0 0 /
\plot 0 1  2 1  2 5  0 5 /
\plot 0 2  2 2 /
\plot 0 3  2 3 /
\plot 0 4  2 4 /
\multiput{$\bullet$} at 0.5 1.5  1.5 1.5   /
\plot 0.5 1.5  1.5 1.5 /
\plot 0.5 1.55  1.5 1.55 /
\plot 0.5 1.45  1.5 1.45 /
\endpicture} at 30 -10.5


\put{\beginpicture
\multiput{} at 0 0  2 5 /
\plot 0 0  1 0  1 5  0 5  0 0 /
\plot 0 1  2 1  2 4  0 4 /
\plot 0 2  3 2  3 3  0 3 /
\multiput{$\bullet$} at 0.5 2.5  1.5 2.5  2.5 2.5  /
\plot 0.5 2.5  2.5 2.5 /
\plot 0.5 2.55  2.5 2.55 /
\plot 0.5 2.45  2.5 2.45 /
\endpicture} at 10 -23

\put{\beginpicture
\multiput{} at 0 0  2 6 /
\plot 0 0  1 0  1 6  0 6  0 0 /
\plot 0 1  2 1  2 4  0 4 /
\plot 0 2  3 2  3 3  0 3 /
\plot 0 5  1 5 /
\multiput{$\bullet$} at 0.5 2.5  1.5 2.5  2.5 2.5  /
\plot 0.5 2.5  2.5 2.5 /
\plot 0.5 2.55  2.5 2.55 /
\plot 0.5 2.45  2.5 2.45 /
\endpicture} at 13.333 -23

\put{\beginpicture
\multiput{} at 0 0  2 6 /
\plot 0 0  1 0  1 6  0 6  0 0 /
\plot 0 1  2 1  2 5  0 5 /
\plot 0 2  3 2  3 3  0 3 /
\plot 0 4  2 4 /
\multiput{$\bullet$} at 0.5 2.5  1.5 2.5  2.5 2.5  /
\plot 0.5 2.5  2.5 2.5 /
\plot 0.5 2.55  2.5 2.55 /
\plot 0.5 2.45  2.5 2.45 /
\endpicture} at 16.667 -23

\put{\beginpicture
\multiput{} at 0 0  2 6 /
\plot 0 0  1 0  1 6  0 6  0 0 /
\plot 0 1  2 1  2 5  0 5 /
\plot 0 2  3 2  3 4  0 4 /
\plot 0 3  3 3 /
\multiput{$\bullet$} at 0.5 2.5  1.5 2.5  2.5 2.5  /
\plot 0.5 2.5  2.5 2.5 /
\plot 0.5 2.55  2.5 2.55 /
\plot 0.5 2.45  2.5 2.45 /
\endpicture} at 20 -23

\put{\beginpicture
\multiput{} at 0 0  2 6 /
\plot 0 0  1 0  1 6  0 6  0 0 /
\plot 0 1  2 1  2 5  0 5 /
\plot 0 2  4 2  4 3  0 3 /
\plot 2 1  3 1  3 4  0 4 /
\multiput{$\bullet$} at 0.5 2.5  1.5 2.5  2.5 2.5  3.5 2.5   2.5 1.5 /
\plot 0.5 2.5  3.5 2.5 /
\plot 0.5 2.55  3.5 2.55 /
\plot 0.5 2.45  3.5 2.45 /
\endpicture} at 15 -30.5

\endpicture}
$$

\subsection{Objects in $\Cal S$ with $q < 4.$}
\label{sec-fifteen-three}

{\bf (a) Indecomposable objects in $\Cal S$ with $2<qX<3$:}

\medskip
Here are the known indecomposable objects $X$ with $2 < qX < 3.$ 
They have $qX = \frac 52$:
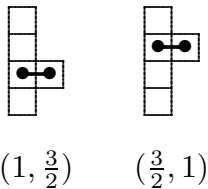
$$
{\beginpicture
    \setcoordinatesystem units <.3cm,.3cm>
\put{\beginpicture
\multiput{} at 0 0  2 4 /
\plot 0 0  1 0  1 4  0 4  0 0 /
\plot 0 1  2 1  2 2  0 2 /
\plot 0 3  1 3 /
\multiput{$\bullet$} at 0.5 1.5  1.5 1.5   /
\plot 0.5 1.5  1.5 1.5 /
\plot 0.5 1.55  1.5 1.55 /
\plot 0.5 1.45  1.5 1.45 /
\endpicture} at 0 0
\put{$(1,\frac32)$} at 0 -4
\put{\beginpicture
\multiput{} at 0 0  2 4 /
\plot 0 0  1 0  1 4  0 4  0 0 /
\plot 0 2  2 2  2 3  0 3 /
\plot 0 1  1 1 /
\multiput{$\bullet$} at 0.5 2.5  1.5 2.5   /
\plot 0.5 2.5  1.5 2.5 /
\plot 0.5 2.55  1.5 2.55 /
\plot 0.5 2.45  1.5 2.45 /
\endpicture} at 5 0
\put{$(\frac32,1)$} at 5 -4
\endpicture}
$$

\bigskip
{\bf (b) The indecomposable objects $X$ in $\Cal S(6)$ with $qX = 3$.}
	\medskip
        
There are fifteen indecomposable objects $X$ with $qX = 3$.

$$
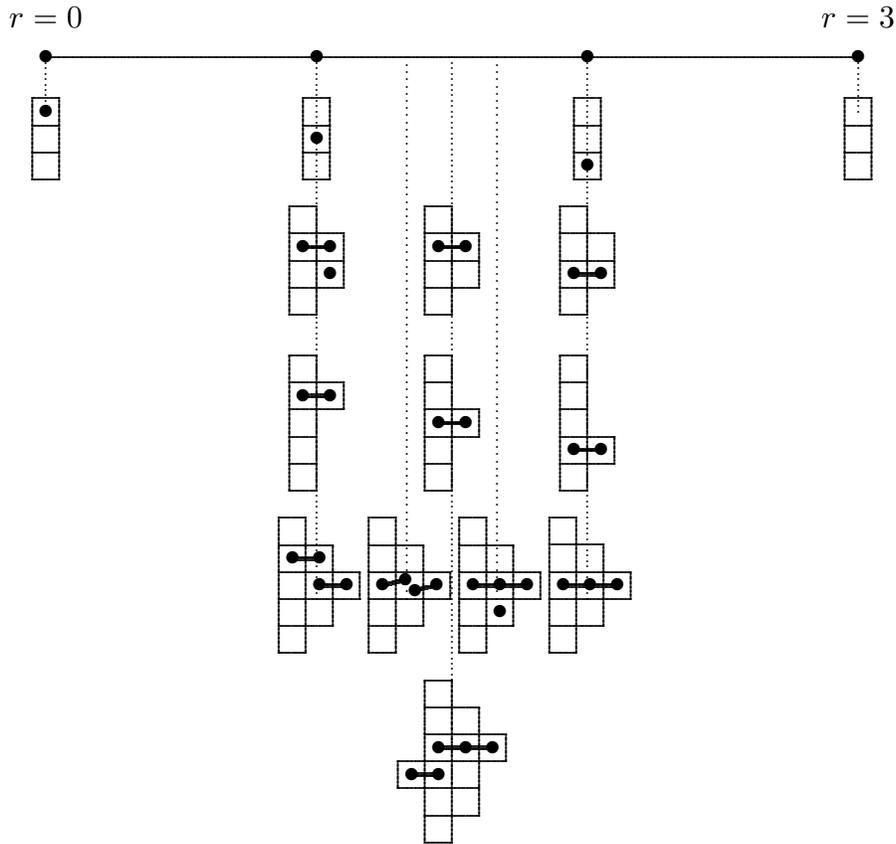
{\beginpicture
    \setcoordinatesystem units <.3cm,.3cm>
\put{$r = 0$} at 0 1.5
\put{$r = 3$} at 30 1.5
\multiput{$\bullet$} at 0 0  10 0  20 0  30 0 /
\plot 0 0  30 0 /

\setdots <.6mm>
\plot 0 0   0 -1 /
\plot 10 0  10 -1 /
\plot 10 -4.5  10 -5.2 /
\plot 10 -10  10 -10.5 /
\plot 10 -16.5  10 -17.5 /
\plot 30 0   30 -1 /
\plot 20 0  20 -1 /
\plot 20 -4.5  20 -5.2 /
\plot 20 -10  20 -10.5 /
\plot 20 -16.5  20 -17.5 /
\plot 15 0  15 -5.2 /
\plot 15 -10  15 -10.5 /
\plot 15 -16.5 15 -22.5 /
\plot 13.333 0  13.333 -17.5 /
\plot 16.667 0  16.667 -17.5 /
\setsolid

\put{\beginpicture
\multiput{} at 0 0  1 3 /
\plot 0 0  1 0  1 3  0 3  0 0 /
\plot 0 1  1 1 /
\plot 0 2  1 2 /
\endpicture} at 30 -3

\put{\beginpicture
\multiput{} at 0 0  1 3 /
\plot 0 0  1 0  1 3  0 3  0 0 /
\plot 0 1  1 1 /
\plot 0 2  1 2 /
\multiput{$\bullet$} at 0.5 0.5   /
\endpicture} at 20 -3
\put{\beginpicture
\multiput{} at 0 0  1 3 /
\plot 0 0  1 0  1 3  0 3  0 0 /
\plot 0 1  1 1 /
\plot 0 2  1 2 /
\multiput{$\bullet$} at 0.5 1.5   /
\endpicture} at 10 -3

\put{\beginpicture
\multiput{} at 0 0  1 3 /
\plot 0 0  1 0  1 3  0 3  0 0 /
\plot 0 1  1 1 /
\plot 0 2  1 2 /
\multiput{$\bullet$} at 0.5 2.5   /
\endpicture} at 0 -3

\put{\beginpicture
\multiput{} at 0 0  2 4 /
\plot 0 0  1 0  1 4  0 4  0 0 /
\plot 0 1  2 1  2 3  0 3 /
\plot 0 2  2 2 /
\multiput{$\bullet$} at 0.5 1.5  1.5 1.5   /
\plot 0.5 1.5  1.5 1.5 /
\plot 0.5 1.55  1.5 1.55 /
\plot 0.5 1.45  1.5 1.45 /
\endpicture} at 20 -7.5 

\put{\beginpicture
\multiput{} at 0 0  2 4 /
\plot 0 0  1 0  1 4  0 4  0 0 /
\plot 0 1  2 1  2 3  0 3 /
\plot 0 2  2 2 /
\multiput{$\bullet$} at 0.5 2.5  1.5 2.5   /
\plot 0.5 2.5  1.5 2.5 /
\plot 0.5 2.55  1.5 2.55 /
\plot 0.5 2.45  1.5 2.45 /
\endpicture} at 15 -7.5

\put{\beginpicture
\multiput{} at 0 0  2 4 /
\plot 0 0  1 0  1 4  0 4  0 0 /
\plot 0 1  2 1  2 3  0 3 /
\plot 0 2  2 2 /
\multiput{$\bullet$} at 0.5 2.5  1.5 2.5  1.5 1.5  /
\plot 0.5 2.5  1.5 2.5 /
\plot 0.5 2.55  1.5 2.55 /
\plot 0.5 2.45  1.5 2.45 /
\endpicture} at 10 -7.5

\put{\beginpicture
\multiput{} at 0 0  2 4 /
\plot 0 0  1 0  1 5  0 5  0 0 /
\plot 0 1  2 1  2 2  0 2 /
\plot 0 3  1 3 /
\plot 0 4  1 4 /
\multiput{$\bullet$} at 0.5 1.5  1.5 1.5   /
\plot 0.5 1.5  1.5 1.5 /
\plot 0.5 1.55  1.5 1.55 /
\plot 0.5 1.45  1.5 1.45 /
\endpicture} at 20 -14

\put{\beginpicture
\multiput{} at 0 0  2 4 /
\plot 0 0  1 0  1 5  0 5  0 0 /
\plot 0 2  2 2  2 3   0 3 /
\plot 0 1  1 1 /
\plot 0 4  1 4 /
\multiput{$\bullet$} at 0.5 2.5  1.5 2.5   /
\plot 0.5 2.5  1.5 2.5 /
\plot 0.5 2.55  1.5 2.55 /
\plot 0.5 2.45  1.5 2.45 /
\endpicture} at 15 -14

\put{\beginpicture
\multiput{} at 0 0  2 4 /
\plot 0 0  1 0  1 5  0 5  0 0 /
\plot 0 3  2 3  2 4  0 4 /
\plot 0 1  1 1 /
\plot 0 2  1 2 /
\multiput{$\bullet$} at 0.5 3.5  1.5 3.5   /
\plot 0.5 3.5  1.5 3.5 /
\plot 0.5 3.55  1.5 3.55 /
\plot 0.5 3.45  1.5 3.45 /
\endpicture} at 10 -14

\put{\beginpicture
\multiput{} at 0 0  2 4 /
\plot 0 0  1 0  1 5  0 5  0 0 /
\plot 0 1  2 1  2 4  0 4 /
\plot 0 2  3 2  3 3  0 3 /
\multiput{$\bullet$} at 0.5 2.5  1.5 2.5  2.5 2.5  /
\plot 0.5 2.5  2.5 2.5 /
\plot 0.5 2.55  2.5 2.55 /
\plot 0.5 2.45  2.5 2.45 /
\endpicture} at 20 -20

\put{\beginpicture
\multiput{} at 0 0  2 4 /
\plot 0 0  1 0  1 5  0 5  0 0 /
\plot 0 1  2 1  2 4  0 4 /
\plot 0 2  3 2  3 3  0 3 /
\multiput{$\bullet$} at 0.5 2.5  1.5 2.5  2.5 2.5  1.5 1.5 /
\plot 0.5 2.5  2.5 2.5 /
\plot 0.5 2.55  2.5 2.55 /
\plot 0.5 2.45  2.5 2.45 /
\endpicture} at 16.667 -20

\put{\beginpicture
\multiput{} at 0 0  2 4 /
\plot 0 0  1 0  1 5  0 5  0 0 /
\plot 0 1  2 1  2 4  0 4 /
\plot 0 2  3 2  3 3  0 3 /
\multiput{$\bullet$} at 0.5 2.5  1.35 2.7  
                         1.7 2.3  2.5 2.5  /
\plot 0.5 2.5  1.3 2.7 /
\plot 0.5 2.55  1.3 2.75 /
\plot 0.5 2.45  1.3 2.65 /

\plot 1.6 2.3  2.5 2.5 /
\plot 1.6 2.35  2.5 2.55 /
\plot 1.6 2.25  2.5 2.45 /
\endpicture} at 13.333 -20

\put{\beginpicture
\multiput{} at 0 0  2 4 /
\plot 0 0  1 0  1 5  0 5  0 0 /
\plot 0 1  2 1  2 4  0 4 /
\plot 0 2  3 2  3 3  0 3 /
\multiput{$\bullet$} at 0.5 3.5  1.5 3.5  1.5 2.5  2.5 2.5  /
\plot 0.5 3.5  1.5 3.5 /
\plot 0.5 3.55  1.5 3.55 /
\plot 0.5 3.45  1.5 3.45 /
\plot 1.5 2.5  2.5 2.5 /
\plot 1.5 2.55  2.5 2.55 /
\plot 1.5 2.45  2.5 2.45 /
\endpicture} at 10 -20

\put{\beginpicture
\multiput{} at 0 0  2 4 /
\plot 0 0  1 0  1 6  0 6  0 0 /
\plot 0 1  2 1  2 5  0 5 /
\plot 0 2  2 2 /
\plot 0 3  3 3  3 4  0 4 /
\plot 0 2 -1 2  -1 3  0 3 /

\multiput{$\bullet$} at 0.5 2.5  -.5 2.5  0.5 3.5 1.5 3.5  2.5 3.5  /
\plot 0.5 2.5  -.5 2.5 /
\plot 0.5 2.55  -.5 2.55 /
\plot 0.5 2.45  -.5 2.45 /

\plot 0.5 3.5  2.5 3.5 /
\plot 0.5 3.55  2.5 3.55 /
\plot 0.5 3.45  2.5 3.45 /

\endpicture} at 15 -27

\endpicture}
$$

\bigskip
{\bf (c)  Indecomposable objects in $\Cal S(9)$ with $q < 4.$}
	\medskip
        
It is easy to modify some indecomposable objects in $\Cal S(6)$ with mean
strictly smaller than 4, in order
to obtain objects $X$ in $\Cal S(n)$ with $n = 7,8,9$ which still satisfy the
condition $qX < 4.$ Here are some objects $X$ in $\Cal S(9)$ with
$vX = 15$ and $bX = 4$, thus $qX = 15/4 < 4;$ they are obvious modifications of
the object $Y = S[6]$ in $\Cal S(6)$.
$$
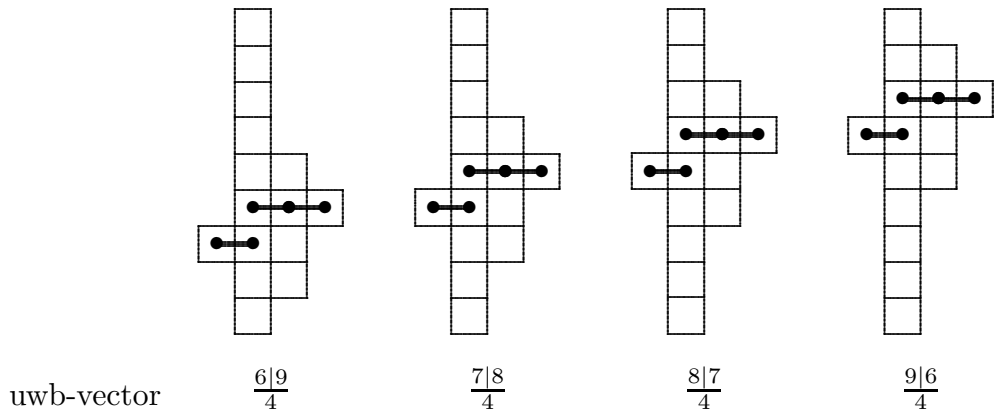
{\beginpicture
    \setcoordinatesystem units <.35cm,.35cm>
\put{\beginpicture

\multiput{} at 0 -2   2 6 /
\put{uwb-vector} at 2 -4.5
\endpicture} at -5 0
\put{\beginpicture

\multiput{} at 0 0   2 9 /
\plot 1 0  2 0  2 6  1 6  1 0 /
\plot 1 1  3 1  3 5  1 5 /
\plot 0 2  3 2 /
\plot 0 2  0 3  4 3  4 4  1 4 /
\multiput{$\bullet$} at 0.5 2.5  1.5 2.5  1.5 3.5  2.5 3.5  3.5 3.5  2.5 3.5 /
\plot .5 2.5   1.5 2.5 /
\plot .5 2.45  1.5 2.45 /
\plot .5 2.55  1.5 2.55 /
\plot 1.5 3.5   3.5 3.5 /
\plot 1.5 3.45  3.5 3.45 /
\plot 1.5 3.55  3.5 3.55 /

\plot 1 6  1 9  2 9  2 6 /
\plot 1 7  2 7 /
\plot 1 8  2 8 /
\put{$\frac{6|9}4$} at 2 -1.5
\endpicture} at 0 0

\put{\beginpicture

\multiput{} at 0 -1   2 8 /
\plot 1 0  2 0  2 6  1 6  1 0 /
\plot 1 1  3 1  3 5  1 5 /
\plot 0 2  3 2 /
\plot 0 2  0 3  4 3  4 4  1 4 /
\multiput{$\bullet$} at 0.5 2.5  1.5 2.5  1.5 3.5  2.5 3.5  3.5 3.5  2.5 3.5 /
\plot .5 2.5   1.5 2.5 /
\plot .5 2.45  1.5 2.45 /
\plot .5 2.55  1.5 2.55 /
\plot 1.5 3.5   3.5 3.5 /
\plot 1.5 3.45  3.5 3.45 /
\plot 1.5 3.55  3.5 3.55 /

\plot 1 6  1 8  2 8  2 6 /
\plot 1 7  2 7 /

\plot 1 0  1 -1  2 -1 2 0 /
\put{$\frac{7|8}4$} at 2 -2.5
\endpicture} at 6 0

\put{\beginpicture

\multiput{} at 0 -2   2 7 /
\plot 1 0  2 0  2 6  1 6  1 0 /
\plot 1 1  3 1  3 5  1 5 /
\plot 0 2  3 2 /
\plot 0 2  0 3  4 3  4 4  1 4 /
\multiput{$\bullet$} at 0.5 2.5  1.5 2.5  1.5 3.5  2.5 3.5  3.5 3.5  2.5 3.5 /
\plot .5 2.5   1.5 2.5 /
\plot .5 2.45  1.5 2.45 /
\plot .5 2.55  1.5 2.55 /
\plot 1.5 3.5   3.5 3.5 /
\plot 1.5 3.45  3.5 3.45 /
\plot 1.5 3.55  3.5 3.55 /

\plot 1 6  1 7  2 7  2 6 /

\plot 1 0  1 -2  2 -2  2 0 /
\plot 1 -1  2 -1 /
\put{$\frac{8|7}4$} at 2 -3.5
\endpicture} at 12 0

\put{\beginpicture

\multiput{} at 0 -2   2 6 /
\plot 1 0  2 0  2 6  1 6  1 0 /
\plot 1 1  3 1  3 5  1 5 /
\plot 0 2  3 2 /
\plot 0 2  0 3  4 3  4 4  1 4 /
\multiput{$\bullet$} at 0.5 2.5  1.5 2.5  1.5 3.5  2.5 3.5  3.5 3.5  2.5 3.5 /
\plot .5 2.5   1.5 2.5 /
\plot .5 2.45  1.5 2.45 /
\plot .5 2.55  1.5 2.55 /
\plot 1.5 3.5   3.5 3.5 /
\plot 1.5 3.45  3.5 3.45 /
\plot 1.5 3.55  3.5 3.55 /

\plot 1 0  1 -3  2 -3  2 0 /
\plot 1 -1  2 -1 /
\plot 1 -2  2 -2 /
\put{$\frac{9|6}4$} at 2 -4.5
\endpicture} at 18 0

\endpicture}
$$

These objects have non-trivial extensions with $Y$. Here are indecomposable objects
in $\Cal S(9)$ which occur in this way (using the first of the four objects presented above):
$$
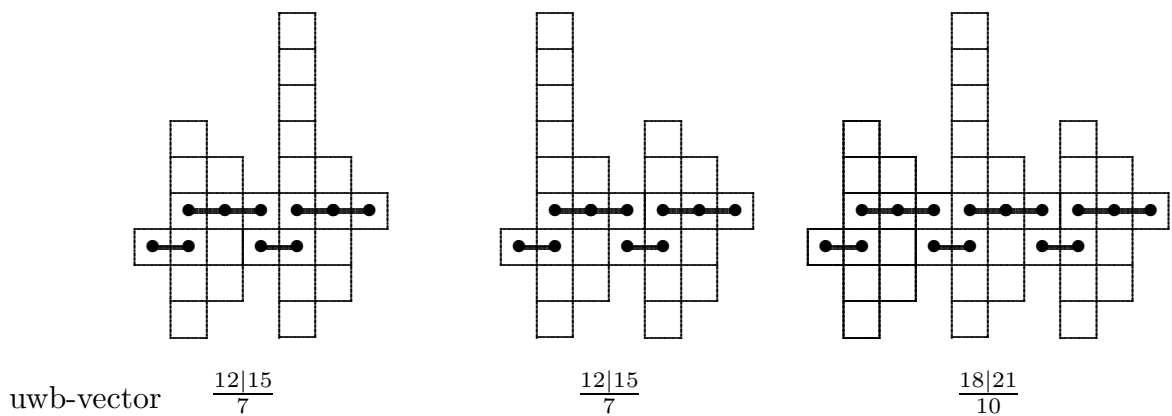
{\beginpicture
    \setcoordinatesystem units <.35cm,.35cm>
\put{\beginpicture

\multiput{} at 0 -2   2 6 /
\put{uwb-vector} at 2 -4.5
\endpicture} at -5 0
\put{\beginpicture

\multiput{} at -2 0   4 9 /
\plot 1 0  2 0  2 6  1 6  1 0 /
\plot 1 1  3 1  3 5  1 5 /
\plot 0 2  3 2 /
\plot 0 2  0 3  4 3  4 4  1 4 /
\multiput{$\bullet$} at 0.5 2.5  1.5 2.5  1.5 3.5  2.5 3.5  3.5 3.5  2.5 3.5 /
\plot .5 2.5   1.5 2.5 /
\plot .5 2.45  1.5 2.45 /
\plot .5 2.55  1.5 2.55 /
\plot 1.5 3.5   3.5 3.5 /
\plot 1.5 3.45  3.5 3.45 /
\plot 1.5 3.55  3.5 3.55 /

\plot 1 6  1 9  2 9  2 6 /
\plot 1 7  2 7 /
\plot 1 8  2 8 /

\plot -2 0  -1 0  -1 6  -2 6  -2 0 /
\plot -2 1  0 1  0 5  -2 5 /
\plot -2 2  0 2 /
\plot -2 3  0 3 /
\plot -2 4  1 4 /
\plot -2 2  -3 2  -3 3  -2 3 /
\put{$\frac{12|15}7$} at 0 -1.5

\multiput{$\bullet$} at -2.5 2.5  -1.5 2.5  -1.5 3.5  -.5 3.5  0.5 3.5  /
\plot -2.5 2.5   -1.5 2.5 /
\plot -2.5 2.45  -1.5 2.45 /
\plot -2.5 2.55  -1.5 2.55 /
\plot -1.5 3.5   .5 3.5 /
\plot -1.5 3.45  .5 3.45 /
\plot -1.5 3.55  .5 3.55 /

\endpicture} at 0 0

\put{\beginpicture

\multiput{} at  0 0   7 9 /
\plot 1 0  2 0  2 6  1 6  1 0 /
\plot 1 1  3 1  3 5  1 5 /
\plot 0 2  3 2 /
\plot 0 2  0 3  4 3  4 4  1 4 /
\multiput{$\bullet$} at 0.5 2.5  1.5 2.5  1.5 3.5  2.5 3.5  3.5 3.5  2.5 3.5 /
\plot .5 2.5   1.5 2.5 /
\plot .5 2.45  1.5 2.45 /
\plot .5 2.55  1.5 2.55 /
\plot 1.5 3.5   3.5 3.5 /
\plot 1.5 3.45  3.5 3.45 /
\plot 1.5 3.55  3.5 3.55 /

\plot 1 6  1 9  2 9  2 6 /
\plot 1 7  2 7 /
\plot 1 8  2 8 /

\multiput{$\bullet$} at 3.5 2.5  4.5 2.5   4.5 3.5  5.5 3.5  6.5 3.5  /
\plot 3.5 2.5   4.5 2.5 /
\plot 3.5 2.45  4.5 2.45 /
\plot 3.5 2.55  4.5 2.55 /
\plot 4.5 3.5   6.5 3.5 /
\plot 4.5 3.45  6.5 3.45 /
\plot 4.5 3.55  6.5 3.55 /

\plot 4 0  4 6  5 6  5 0  4 0 /
\plot 4 1  6 1  6 5  4 5 /
\plot 4 3  7 3  7 4  4 4 /

\plot 3 2  6 2 /

\put{$\frac{12|15}7$} at 3 -1.5
\endpicture} at 10 0

\put{\beginpicture

\multiput{} at -2 0   2 9 /
\plot 1 0  2 0  2 6  1 6  1 0 /
\plot 1 1  3 1  3 5  1 5 /
\plot 0 2  3 2 /
\plot 0 2  0 3  4 3  4 4  1 4 /
\multiput{$\bullet$} at 0.5 2.5  1.5 2.5  1.5 3.5  2.5 3.5  3.5 3.5  2.5 3.5 /
\plot .5 2.5   1.5 2.5 /
\plot .5 2.45  1.5 2.45 /
\plot .5 2.55  1.5 2.55 /
\plot 1.5 3.5   3.5 3.5 /
\plot 1.5 3.45  3.5 3.45 /
\plot 1.5 3.55  3.5 3.55 /

\plot 1 6  1 9  2 9  2 6 /
\plot 1 7  2 7 /
\plot 1 8  2 8 /

\plot -2 0  -1 0  -1 6  -2 6  -2 0 /
\plot -2 1  0 1  0 5  -2 5 /
\plot -2 2  0 2 /
\plot -2 3  0 3 /
\plot -2 4  1 4 /
\plot -2 2  -3 2  -3 3  -2 3 /

\multiput{$\bullet$} at -2.5 2.5  -1.5 2.5  -1.5 3.5  -.5 3.5  0.5 3.5  /
\plot -2.5 2.5   -1.5 2.5 /
\plot -2.5 2.45  -1.5 2.45 /
\plot -2.5 2.55  -1.5 2.55 /
\plot -1.5 3.5   .5 3.5 /
\plot -1.5 3.45  .5 3.45 /
\plot -1.5 3.55  .5 3.55 /

\plot -2 0  -1 0  -1 6  -2 6  -2 0 /
\plot -2 1  0 1  0 5  -2 5 /
\plot -2 2  0 2 /
\plot -2 3  0 3 /
\plot -2 4  1 4 /
\plot -2 2  -3 2  -3 3  -2 3 /

\multiput{$\bullet$} at 3.5 2.5  4.5 2.5   4.5 3.5  5.5 3.5  6.5 3.5  /
\plot 3.5 2.5   4.5 2.5 /
\plot 3.5 2.45  4.5 2.45 /
\plot 3.5 2.55  4.5 2.55 /
\plot 4.5 3.5   6.5 3.5 /
\plot 4.5 3.45  6.5 3.45 /
\plot 4.5 3.55  6.5 3.55 /

\plot 4 0  4 6  5 6  5 0  4 0 /
\plot 4 1  6 1  6 5  4 5 /
\plot 4 3  7 3  7 4  4 4 /

\plot 3 2  6 2 /

\put{$\frac{18|21}{10}$} at 2 -1.5
\endpicture} at 20 0

\endpicture}
$$

Using induction, we can construct indecomposable objects in $\Cal S(9)$
with uwb-vector of the form $\dfrac{6t|6t\!+\!3}{3t\!+\!1}$ for any $t\ge 1,$ thus with
mean $(12t+3)/(3t+1) < 4.$ Of course, the sequence of the corresponding pr-vectors converges 
to $(2,2).$ For $t = 1,2,3,$ we may start with the following objects:
$$
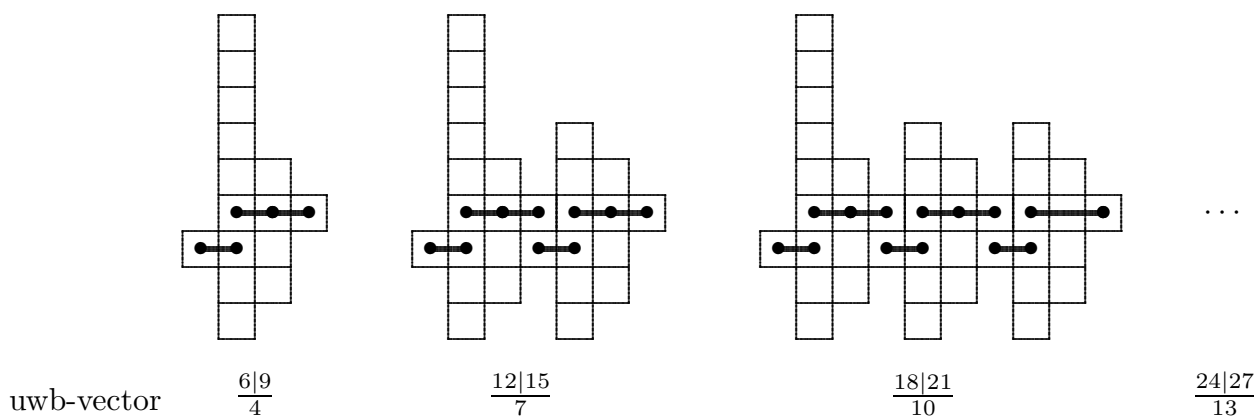
{\beginpicture
    \setcoordinatesystem units <.35cm,.35cm>
\put{\beginpicture
\multiput{} at 0 -2   2 6 /
\put{uwb-vector} at 3 -4.5
\endpicture} at -5 0
\put{\beginpicture

\multiput{} at 0 0   2 9 /
\plot 1 0  2 0  2 6  1 6  1 0 /
\plot 1 1  3 1  3 5  1 5 /
\plot 0 2  3 2 /
\plot 0 2  0 3  4 3  4 4  1 4 /
\multiput{$\bullet$} at 0.5 2.5  1.5 2.5  1.5 3.5  2.5 3.5  3.5 3.5  2.5 3.5 /
\plot .5 2.5   1.5 2.5 /
\plot .5 2.45  1.5 2.45 /
\plot .5 2.55  1.5 2.55 /
\plot 1.5 3.5   3.5 3.5 /
\plot 1.5 3.45  3.5 3.45 /
\plot 1.5 3.55  3.5 3.55 /

\plot 1 6  1 9  2 9  2 6 /
\plot 1 7  2 7 /
\plot 1 8  2 8 /
\put{$\frac{6|9}4$} at 2 -1.5
\endpicture} at 0 0

\put{\beginpicture

\multiput{} at  0 0   7 9 /
\plot 1 0  2 0  2 6  1 6  1 0 /
\plot 1 1  3 1  3 5  1 5 /
\plot 0 2  3 2 /
\plot 0 2  0 3  4 3  4 4  1 4 /
\multiput{$\bullet$} at 0.5 2.5  1.5 2.5  1.5 3.5  2.5 3.5  3.5 3.5  2.5 3.5 /
\plot .5 2.5   1.5 2.5 /
\plot .5 2.45  1.5 2.45 /
\plot .5 2.55  1.5 2.55 /
\plot 1.5 3.5   3.5 3.5 /
\plot 1.5 3.45  3.5 3.45 /
\plot 1.5 3.55  3.5 3.55 /

\plot 1 6  1 9  2 9  2 6 /
\plot 1 7  2 7 /
\plot 1 8  2 8 /

\multiput{$\bullet$} at 3.5 2.5  4.5 2.5   4.5 3.5  5.5 3.5  6.5 3.5  /
\plot 3.5 2.5   4.5 2.5 /
\plot 3.5 2.45  4.5 2.45 /
\plot 3.5 2.55  4.5 2.55 /
\plot 4.5 3.5   6.5 3.5 /
\plot 4.5 3.45  6.5 3.45 /
\plot 4.5 3.55  6.5 3.55 /

\plot 4 0  4 6  5 6  5 0  4 0 /
\plot 4 1  6 1  6 5  4 5 /
\plot 4 3  7 3  7 4  4 4 /

\plot 3 2  6 2 /

\put{$\frac{12|15}7$} at 3 -1.5
\endpicture} at 8 0
\put{\beginpicture

\multiput{} at  0 0   7 9 /
\plot 1 0  2 0  2 6  1 6  1 0 /
\plot 1 1  3 1  3 5  1 5 /
\plot 0 2  3 2 /
\plot 0 2  0 3  4 3  4 4  1 4 /
\multiput{$\bullet$} at 0.5 2.5  1.5 2.5  1.5 3.5  2.5 3.5  3.5 3.5  2.5 3.5 /
\plot .5 2.5   1.5 2.5 /
\plot .5 2.45  1.5 2.45 /
\plot .5 2.55  1.5 2.55 /
\plot 1.5 3.5   3.5 3.5 /
\plot 1.5 3.45  3.5 3.45 /
\plot 1.5 3.55  3.5 3.55 /

\plot 1 6  1 9  2 9  2 6 /
\plot 1 7  2 7 /
\plot 1 8  2 8 /

\multiput{$\bullet$} at 3.5 2.5  4.5 2.5   4.5 3.5  5.5 3.5  6.5 3.5  /
\plot 3.5 2.5   4.5 2.5 /
\plot 3.5 2.45  4.5 2.45 /
\plot 3.5 2.55  4.5 2.55 /
\plot 4.5 3.5   6.5 3.5 /
\plot 4.5 3.45  6.5 3.45 /
\plot 4.5 3.55  6.5 3.55 /

\plot 4 0  4 6  5 6  5 0  4 0 /
\plot 4 1  6 1  6 5  4 5 /
\plot 4 3  7 3  7 4  4 4 /

\plot 3 2  6 2 /

\multiput{$\bullet$} at 6.5 2.5  7.5 2.5   7.5 3.5  9.5 3.5  9.5 3.5  /
\plot 6.5 2.5   7.5 2.5 /
\plot 6.5 2.45  7.5 2.45 /
\plot 6.5 2.55  7.5 2.55 /
\plot 7.5 3.5   9.5 3.5 /
\plot 7.5 3.45  9.5 3.45 /
\plot 7.5 3.55  9.5 3.55 /

\plot 7 0  7 6  8 6  8 0  7 0 /
\plot 7 1  9 1  9 5  9 5 /
\plot 7 3  10 3  10 4  9 4 /

\plot 6 2  9 2 /

\plot 7 4  9 4 /
\plot 7 5  9 5 /
\put{$\frac{18|21}{10}$} at 4.5 -1.5
\endpicture} at 19 0 
\put{\beginpicture
\multiput{} at  0 0   0 9 /
\put{$\cdots$} at 0 3.5 
\put{$\frac{24|27}{13}$} at 0 -1.5
\endpicture} at 27 0

\endpicture} 
$$

\vfill\eject

\centerline{\Gross Sixth part: Questions and comments.}
\addcontentsline{toc}{part}{Sixth part: Questions and comments.}

\section{Questions.}
\label{sec-sixteen}
        
We will formulate a lot of open questions, the reader will
observe that at present there are more questions than results. In case we consider pr-vectors
$(p,r)$ in $\mathbb T(n)$, we write $q = p+r.$

\subsection{First question: finite pr-vectors.}
\label{sec-sixteen-one}

Let us call a pr-vector $(p,r)$ {\it $n$-finite} provided there are only finitely many
indecomposable objects in $\Cal S(n)$ with pr-vector $(p,r),$ and {\it finite,} provided
$(p,r)$ is $n$-finite for all $n\in \mathbb N.$ 
	\medskip
        
According to Theorem~\ref{theoremone}, we know that all pr-vectors $(p,r)$ with $p < 1$ or $r < 1$
are finite. According to Theorem~\ref{theoremtwo}, also the pr-vector $(1,1)$ is finite. The first 
question to be mentioned is: {\it Which additional pr-vectors are finite?} We expect that many,
may-be all vectors $(p,r)$ with $q < 4$ are finite.

\subsection{Second: indecomposable objects $X$ with $qX < 4.$}
\label{sec-sixteen-two}

\phantom m
\vspace{-3mm}

{\bf (1)} Let $X$ be indecomposable with $qX < 3$. Is $X$ in 
$\Cal S(4)$, thus a picket or a bipicket?
	\medskip
        
{\bf (2)} Let $X$ be indecomposable with $qX = 3$. Is $X$ in $\Cal S(6)$, thus
$bX \le 4$ ?
	\newline
The 15 indecomposable objects $X$ in $\Cal S(6)$ with $qX = 3$ are 
exhibited in Section \ref{sec-fifteen-three}--(c).
	\medskip
        
{\bf (3)} Is there an indecomposable object $X$ in $\Cal S$ with $3 < qX < \frac{16}5$ ?\newline
In Appendix~\ref{app-B-six} we present two pentapickets, each satisfies $qX=\frac{16}5$. 
	\medskip
        
{\bf (4)}  Is the height of the indecomposable objects $X$
with $qX < 4$ bounded? Or, at least: Let $\epsilon > 0$.
 Is the height of the indecomposable objects $X$
with $qX \le 4-\epsilon$ bounded?
	\medskip
        
{\bf (5)} Is there a BTh-vector with $q < 4$ ?
	\medskip
        
{\bf (6)} 
The region $p \ge 1,\ r\ge 1, q< 4$ is of great interest. Can one determine
the values $q$ for indecomposable objects in this region? Does there exist an
accumulation point different from 4 ?
$$
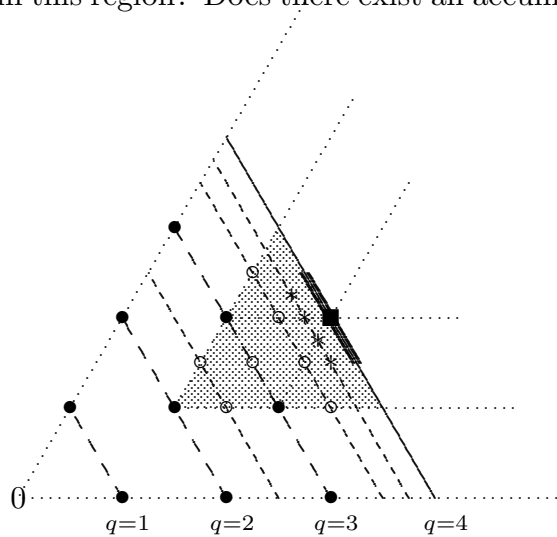
{\beginpicture
   \setcoordinatesystem units <.57735cm,1cm>
\put{} at 5 5 
\multiput{$\bullet$} at 2 0  4 0  6 0  1 1  2 2  3 3  3 1  4 2  5 1 /
\plot 8 0  4 4 /
\plot 5.55 2.5  6.55 1.5 /
\plot 5.45 2.5  6.45 1.5 /
\plot 5.58 2.5  6.58 1.5 /
\plot 5.42 2.5  6.42 1.5 /
\setdots <1mm>
\plot 10.5 0  0 0  5.5 5.5 /
\plot 9.5 1  3 1  6.5 4.5 /
\plot 7.5  3.5  5 1  4 2  8.5 2 /
\put{$\ssize \blacksquare$} at 6 2 
\multiput{$\circ$} at 4 1  3.5 1.5  4.5 1.5  4.5 2.5  5 2  5.5 1.5  6 1 /
\setdashes <1.5mm>
\plot 2 0  1 1 /
\plot 4 0  2 2 /
\plot 6 0  3 3 /

\setdashes <.8mm>
\plot 5 0  2.5 2.5 /
\plot 7 0  3.5 3.5 /
\put{$\ssize q = 1$} at 2.1 -.3 
\put{$\ssize q = 2$} at 4.1 -.3 
\put{$\ssize q = 3$} at 6.1 -.3 
\put{$\ssize q = 4$} at 8.2 -.3 

\plot 7.5 0  3.75  3.75 /
\multiput{$*$} at  6 1.5  5.75 1.75  5.5 2  5.25 2.25 /

\put{$0$} at 0 0 
\setshadegrid span <.35mm>
\vshade 3 1 1  <z,z,z,z> 5 1 3  <z,z,,> 7 1 1 /
\endpicture}
$$
Here, the dashed lines show the pr-vectors with constant mean $q$, equal to $1,\ 2,\ 
5/2,\ 3,\ 7/2,$ and $15/4.$ The four stars $*$ on the line $q = 15/4$ 
indicate the position of the four indecomposables
with mean equal to $15/4$ which have been exhibited at the beginning of Section \ref{sec-fifteen-three}--(c). 
(Of course, there are also indecomposables with $q = 9/3,\ 10/3,$ and $11/3$, see
Section \ref{sec-fifteen-two}--(c) as well as Appendix~\ref{app-B})

\subsection{What is the role of the line $3p+r = 7$ ?}
\label{sec-sixteen-three}

{\bf Question.} 
  Are there  BTh-vectors in the triangle with corners $(1,1),$ $(2,1),$ and 
   $(1,4)$?
	\medskip
        
Here is the triangle. The four marked points represent BTh-families discussed in this paper.
See Section~\ref{sec-eight-one} for the family in $\Cal S(7)$ with uwb-vector $\frac{6|10}4$.
The family in $\Cal S_3(7)$ with uwb-vector $\frac{7|14}5$ is presented in Section~\ref{sec-eight-two}.
In Section~\ref{sec-eight-three} we exhibit the family in $\Cal S(8)$ with uwb-vector $\frac{6|17}5$.
The last point marks the position of the family in $\Cal S(9)$ of uwb-vector
$\frac{6|24}6$ from Section~\ref{sec-seven-four}.
	
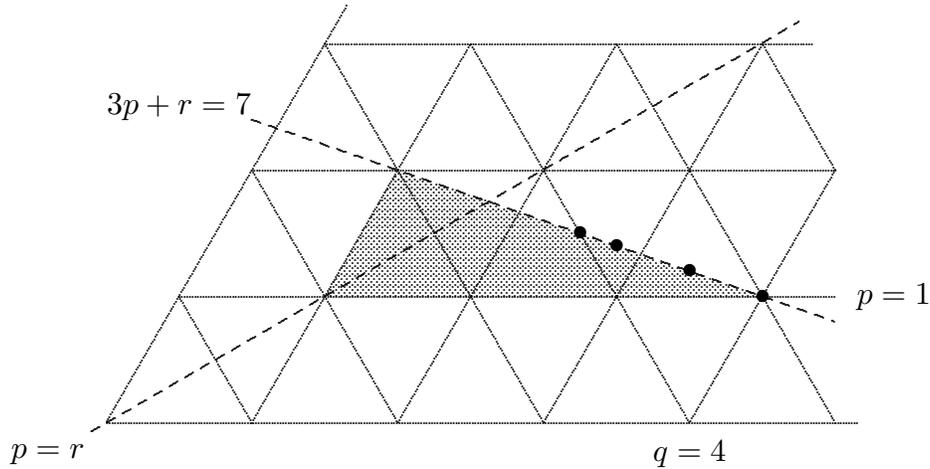
\medskip

$$  
{\beginpicture
   \setcoordinatesystem units <.404cm,.7cm>
   \setcoordinatesystem units <.808cm,1.4cm>
\multiput{} at 0 0  10 3.5  /
\multiput{$\bullet$} at 6.5 1.5  7 1.4  8 1.2  9 1 /
\setdots <.3mm>
\plot 10.3 0  10 0  0 0  3 3  6 0  9 3  10 2 /
\plot 10 0  7 3  4 0  2 2  10 2  /
\plot 10 2  8 0  5 3  2 0  1 1  10 1 /
\plot 3.3 3.3  3 3  9.7 3 /
\setdashes <1mm>
\plot -.2 -.07  9.5 3.175  /
\plot 2 2.4  9 1  10 0.8 /
\put{$3p+r = 7$} at 1 2.5 

\setshadegrid span <.3mm>
\vshade 3 1 1 <z,z,,> 4 1 2    <z,z,,>  9 1 1  /

\put{$q = 4$} at 8 -.25
\put{$p = 1$} at 10.8 1 
\put{$ p = r$} at -.8 -.25
\endpicture}
$$

\subsection{Terra Incognita.}
\label{sec-sixteen-four}
        
Many questions concern the convex closure of the
pr-vectors $(1,1),\ (1,4),\ (2,4),\ (2,2),$
$$  
{\beginpicture
   \setcoordinatesystem units <.404cm,.7cm>
   \setcoordinatesystem units <.808cm,1.4cm>
\multiput{} at 0 0  10 3.5  /
\multiput{$\bullet$} at 3 1  6 2   9 1  10 2 /
\put{$(1,1)$} at 3.2 .75
\put{$(1,4)$} at 9.2 .75
\put{$(2,2)$} at 5.85 2.25
\put{$(2,4)$} at 10.2 2.25

\setdots <.3mm>
\plot 10.3 0  10 0  0 0  3 3  6 0  9 3  10 2 /
\plot 10 0  7 3  4 0  2 2  10 2  /
\plot 10 2  8 0  5 3  2 0  1 1  10 1 /
\plot 3.3 3.3  3 3  9.7 3 /
\setdashes <1mm>
\plot -.2 -.07  9.5 3.175  /

\setshadegrid span <.3mm>
\vshade 3 1 1 <z,z,,> 6 1 2 <z,z,,>  9 1 2  <z,z,,>  10 2 2  /

\put{$\ssize p = 1$} at 11.5 1 
\put{$\ssize p = 2$} at 11.5    2 
\put{$\ssize p = r$} at -.8 -.2
\put{$\ssize r = 4$} at 8 -.4

\setdots <1mm>
\plot 8 0  7.5 -.5 /
\plot 10 1  11 1 / 
\plot 10 2  11 2 /

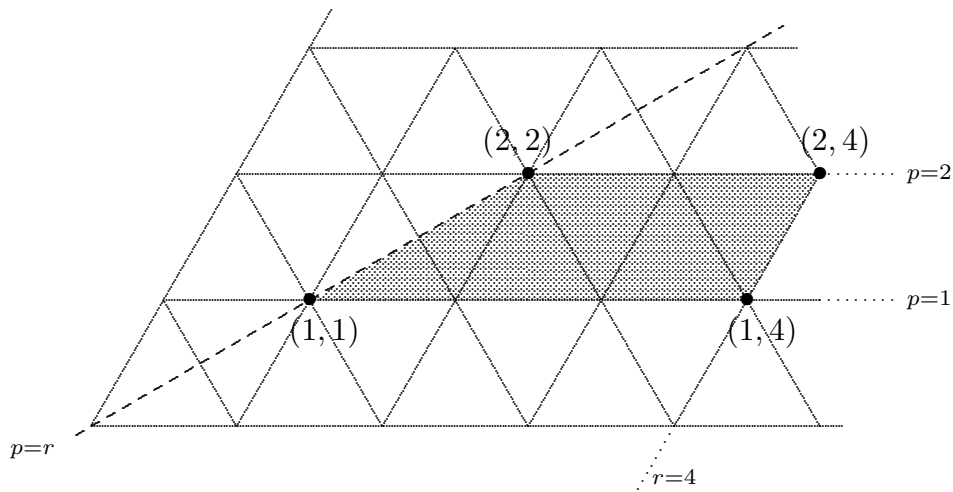
\endpicture}
$$
Note that the corner $(2,2)$ is the pr-vector of the standard family in $\Cal S(6)$,
whereas the corner $(1,4)$ is the pr-vector of the family in $\Cal S(9)$ with $p = 1.$

{\it A pr-vector with boundary distance at least $1$ either lies, up to the symmetries,
  in this trapezoid, or else it is known to be a BTh-vector,} see Theorems~\ref{theoremfour} and
\ref{theoremten}.

\subsection{The set of BTh-vectors.}
\label{sec-sixteen-five}

One of the main aims of the present paper concerns the existence of BTh-vectors.
There are many questions. As we have mentioned already: Is there a BTh-vector with $q < 4$ ?
Is there a BTh-vector with $3p+r < 7$ ? Further questions:
	\medskip
        
{\bf (1)} Is there a BTh-vector for $\Cal S(8)$ with $p = 1$ ? (There is one for $\Cal S(9)$.)
	\medskip
        
{\bf (2)} For any $n$, is the region of all BTh-vectors for $\Cal S(n)$ convex
in $\mathbb T(n)\cap \mathbb Q^2$ ?
\medskip
        
{\bf (3)} Is the
\phantomsection{region}
of all BTh-vectors for $\Cal S(n)$ even 
the convex hull of the BTh-vectors exhibited in
Sections~\ref{sec-six}, \ref{sec-seven}, \ref{sec-eight}?
Here are 
the cases $n = 7,\ 8,\ 9,\ 10,\ 15.$ The dark region is the
convex hull in question. 

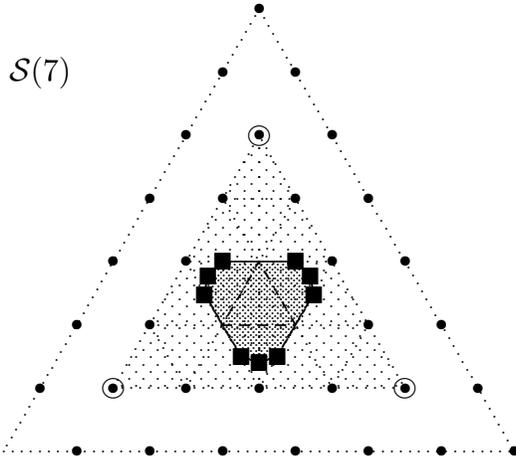
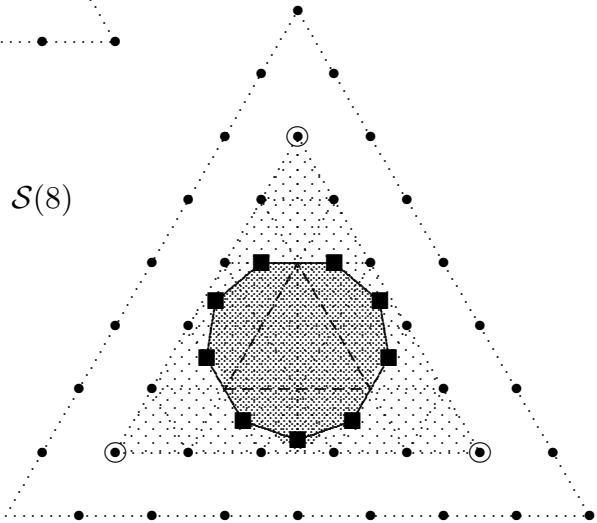
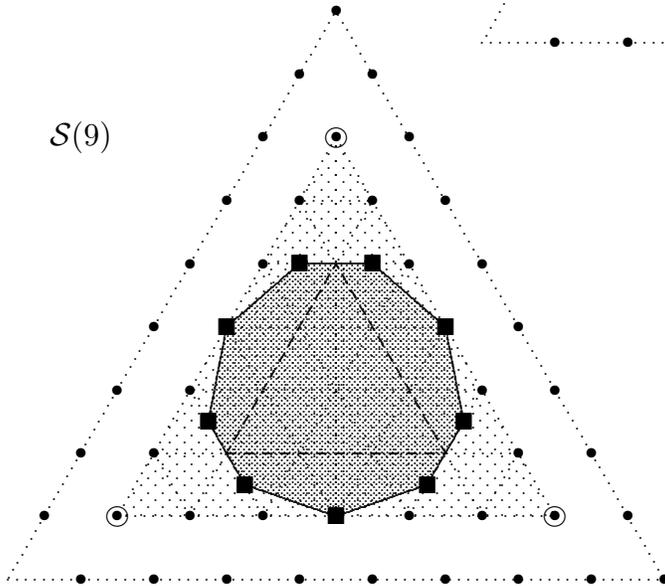
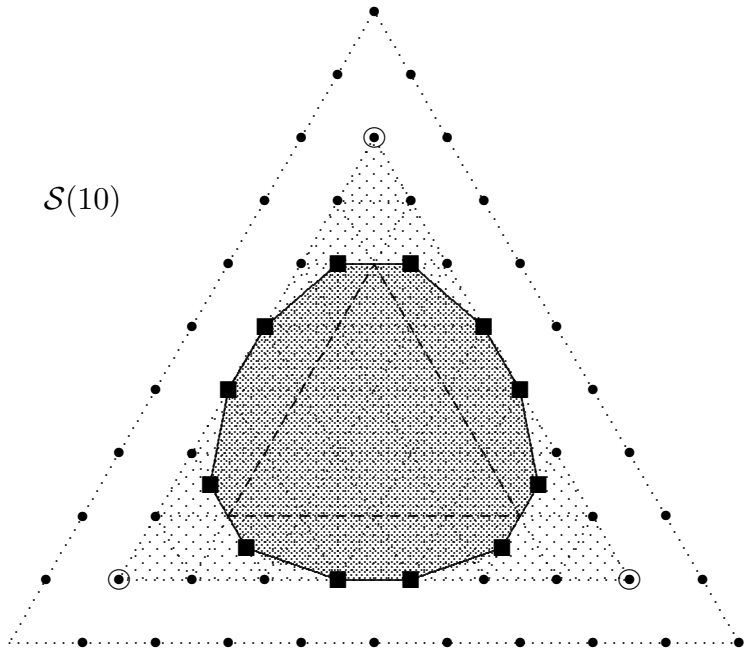
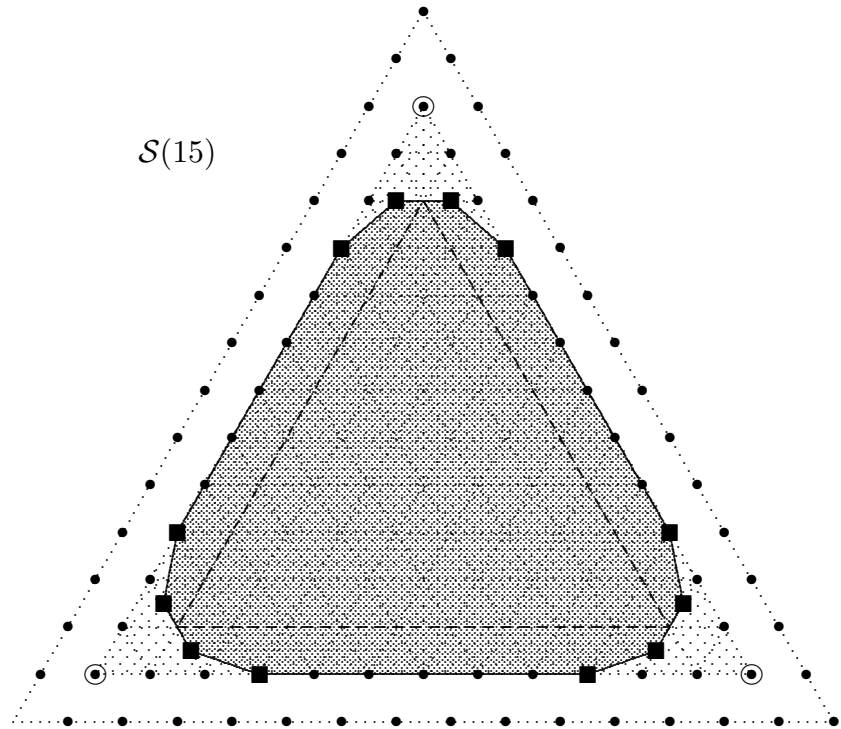
\addcontentsline{lof}{subsection}{Terra incognita.}

$$  
{\beginpicture
   \setcoordinatesystem units  <.404cm,.7cm>

\put{\beginpicture
   \setcoordinatesystem units  <.404cm,.7cm>
\multiput{} at -7 0  7 7  /
\multiput{$\ssize \bullet$} at   -5 0  -3 0  -1 0  1 0  3 0  5 0  7 0
   -6 1  -5 2  -4 3  -3 4  
   -2 5  -1 6   0 7   
    6 1   5 2   4 3   3 4   2 5   1 6    /
\multiput{$\ssize \bullet$} at 
    -4 1  -2 1  0 1  2 1  4 1   -3 2  -2 3  -1 4  0 5    3 2  2 3  1 4    /
\multiput{$\ssize \bigcirc$} at -4 1  4 1  0 5 /
\setdots <1mm>
\plot -7 0  7 0  0 7  -7 0 /
\plot -4 1  4 1  0 5  -4 1 /
\plot -3 2  -2 1  1 4  -1 4  2 1  3 2  -3 2 /
\plot  -2 3  0 1  2 3  -2 3 /
\multiput{$\ssize \blacksquare$} at -.5 1.5  0 1.4   .5 1.5
                                    1.5 2.47  1.4 2.77   1 3
                                    -1 3   -1.4 2.77  -1.5 2.47 /
\setsolid
\plot   -.5 1.5  0 1.4   .5 1.5
                                    1.5 2.5  1.4 2.77   1 3
                                    -1 3   -1.4 2.77  -1.5 2.5  -.5 1.5 /
\setshadegrid span <.3mm>
\vshade -1.5 2.5 2.5   <z,z,,> -1 1.9 3  <z,z,,> 0 1.4 3 <z,z,,>  1 1.9 3  <z,z,,>
    1.5 2.5 2.5 /

\setshadegrid span <.8mm>
\vshade -4  1 1  <z,z,,> 0 1 5  <z,z,,> 4 1 1 /

\put{$\Cal S(7)$} at -6 6
\setdashes <1mm> 
\plot -1 2  1 2  0 3  -1 2 /
\endpicture} at 0 0

\put{\beginpicture
   \setcoordinatesystem units  <.404cm,.7cm>
\multiput{} at -8 0  8 8 /
\multiput{$\ssize \bullet$} at   -6 0  -4 0  -2 0  0 0  2 0  4 0  6 0  8 0 
   -7 1  -6 2  -5 3  -4 4  -3 5  -2 6  -1 7   0 8 
    7 1   6 2   5 3   4 4   3 5   2 6   1 7   /
\multiput{$\ssize \bullet$} at 
    -5 1  -3 1  -1 1  1 1  3 1  5 1  4 2  3 3  2 4  1 5  0 6  -4 2  -3 3  -2 4  -1 5    /
\multiput{$\ssize \bigcirc$} at -5 1  5 1  0 6 /
\setdots <1mm>
\plot -8 0  8 0  0 8  -8 0 /
\plot -5 1  5 1  0 6  -5 1 /
\plot -4 2  -3 1  1 5  -1 5  3 1  4 2   -4 2 /
\plot -3 3  -1 1  2 4  -2 4  1 1  3 3  -3 3 /
\multiput{$\ssize \blacksquare$} at 1.5 1.5  2.5 2.5  1 4  -1 4  -1.5 1.5  -2.5 2.5  
    0  1.2  2.25 3.4  -2.25 3.4 /
\setsolid
\plot 0 1.2  1.5 1.5  2.5 2.5  2.25 3.4  1 4  -1 4  -2.25 3.4  -2.5 2.5  -1.5 1.5  0 1.2   /
\setshadegrid span <.3mm>
\vshade -2.5 2.5 2.5   <z,z,,> -2.25 2.4 3.4  <z,z,,> -1.5 1.4 3.8 <z,z,,> -1 1.4 4 
    <z,z,,> 0 1 4  <z,z,,> 1 1.4 4 <z,z,,> 1.5 1.4 3.8 <z,z,,> 2.25 2.4 3.4  <z,z,,> 2.5 2.5 2.5 /

\setshadegrid span <.8mm>
\vshade -5  1 1  <z,z,,> 0 1 6  <z,z,,> 5 1 1 /

\put{$\Cal S(8)$} at -7 5
\setdashes <1mm> 
\plot -2 2  2 2  0 4  -2 2 /
\endpicture} at 12 -7

\put{\beginpicture
   \setcoordinatesystem units  <.404cm,.7cm>
\multiput{} at -9 0  9 9 /
\multiput{$\ssize \bullet$} at -7 0  -5 0  -3 0  -1 0  1 0  3 0  5 0  7 0  9 0
   -8 1  -7 2  -6 3  -5 4  -4 5  -3 6  -2 7  -1 8  0 9 
    8 1   7 2   6 3   5 4   4 5   3 6   2 7   1 8  /
\multiput{$\ssize \bullet$} at 
    -6 1  -4 1  -2 1  0 1  2 1  4 1  6 1 
     -5 2  -4 3  -3 4  -2 5  -1 6  0 7  1 6  2 5  3 4  4 3  5 2   /
\multiput{$\ssize \bigcirc$} at -6 1  6 1  0 7 /
\setdots <1mm>
\plot -9 0  9 0  0 9  -9 0 /
\plot -6 1  6 1  0 7  -6 1 /
\plot -5 2  -4 1  1 6  -1 6  4 1  5 2  -5 2 /
\plot -4 3  -2 1  2 5  -2 5  2 1  4 3  -4 3  /
\plot -3 4  0 1  3 4  -3 4 /
\multiput{$\ssize \blacksquare$} at 0 1  3 4  -3 4  2.5 1.5  3.5 2.5  -2.5 1.5  -3.5 2.5 
   -1 5  1 5  /
\setsolid
\plot 0 1  2.5 1.5  3.5 2.5  3 4  1 5  -1 5 -3 4  -3.5 2.5  -2.5 1.5  0 1  /
\setshadegrid span <.3mm>
\vshade -3.5  2.5  2.5 <z,z,,> -3 2 4  <z,z,,> -2.5  1.5 4.3 <z,z,,> -1 1.18 5 <z,z,,> 0 1 5 
   <z,z,,> 1 1.18 5  <z,z,,> 2.5  1.5 4.3 <z,z,,> 3 2 4  <z,z,,>  3.5  2.5  2.5 /
\setshadegrid span <.8mm>
\vshade -6 1 1  <z,z,,> 0 1 7  <z,z,,> 6 1 1 /

\put{$\Cal S(9)$} at -7 7 
\setdashes <1mm> 
\plot -3 2  3 2  0 5  -3 2 /

\endpicture} at 0 -15 
\endpicture}
$$

$$
{\beginpicture
   \setcoordinatesystem units  <.404cm,.7cm>
\put{\beginpicture
   \setcoordinatesystem units  <.404cm,.7cm>
\multiput{} at -9 0  11 10 /
\multiput{$\ssize \bullet$} at -7 0  -5 0  -3 0  -1 0  1 0  3 0  5 0  7 0  9 0
   -8 1  -7 2  -6 3  -5 4  -4 5  -3 6  -2 7  -1 8  0 9 
    8 1   7 2   6 3   5 4   4 5   3 6   2 7   1 8  /
\multiput{$\ssize \bullet$} at 
    -6 1  -4 1  -2 1  0 1  2 1  4 1  6 1 
     -5 2  -4 3  -3 4  -2 5  -1 6  0 7  
     11 0  10 1  9 2  8 3  7 4  6 5  5 6  4 7  3 8  2 9  1 10  /
\multiput{$\ssize \bigcirc$} at -6 1  8 1  1 8 /
\setdots <1mm>
\plot -9 0  11 0  1 10  -9 0  /
\plot -6 1  8 1  1 8  -6 1 /
\plot -5 2  -4 1  2 7  0 7  6 1  7 2  -5 2 /
\plot -4 3  -2 1  3 6  -1 6  4 1  6 3 -4 3 /
\plot -3 4  0 1  4 5  -2 5  2 1  5 4  -3 4 /
\multiput{$\ssize \blacksquare$} at  -3 4  -2.5 1.5 -3.5 2.5 
   -2 5  0 1  2 1  4.5 1.5  5.5 2.5  5 4  4 5  2 6  0 6   /
\setsolid
\plot   -2 5  -3 4  -3.5 2.5  -2.5 1.5  0 1  2 1 
      4.5 1.5  5.5 2.5  5 4  4 5  2 6  0 6  -2 5 /
\setshadegrid span <.3mm>
\vshade -3.5  2.5  2.5 <z,z,,> -3 2 4  <z,z,,> -2.5  1.5 4.6 <z,z,,> 
        -2 1.3 5  <z,z,,> 0 1 6 <z,z,,> 2 1 6  <z,z,,> 
        4 1.3 5  <z,z,,> 4.5 1.5 4.6   <z,z,,>  5 2 4  <z,z,,> 5.5 2.5 2.5 /
\setshadegrid span <.8mm>
\vshade -6 1 1  <z,z,,> 1 1 8  <z,z,,> 8 1 1 /

\put{$\Cal S(10)$} at -7 7 
\setdashes <1mm> 
\plot -3 2  5 2  1 6  -3 2 /
\endpicture} at 0 0

\put{\beginpicture
   \setcoordinatesystem units  <.404cm,.7cm>
   \setcoordinatesystem units  <.303cm,.525cm>
\multiput{} at -15 0  15 15 /
\multiput{$\ssize \bullet$} at 0 15 
   -14 1  -13 2  -12 3  -11 4  -10 5  -9 6  -8 7  -7 8  -6 9  -5 10  -4 11  -3 12  -2 13  -1 14 
   -13 0  -11 0  -9 0  
   -7 0  -5 0  -3 0  -1 0  1 0  3 0  5 0  7 0  9 0  11 0  13 0  15 0 
   14 1  13 2  12 3  11 4  10 5  9 6  8 7  7 8  6 9  5 10  4 11  3 12  2 13  1 14 
   /
\multiput{$\ssize \bigcirc$} at -12 1  12 1  0 13 /
\setdots<1mm>
\plot -15 0  15 0  0 15  -15 0 /
\plot -12 1  12 1  0 13 -12 1 /
\plot -10 1  -11 2  11 2  10 1  -1 12  1 12  -10 1 /
\plot -8 1  -10 3  10 3  8 1  -2 11  2 11  -8 1 /
\plot -6 1  -9 4  9 4  6 1  -3 10  3 10  -6 1 /
\plot -4 1  -8 5  8 5  4 1  -4 9  4 9  -4 1 /
\plot -2 1  -7 6  7 6  2 1  -5 8  5 8  -2 1 /
\plot 0 1  -6 7  6 7  0 1 /
\multiput{$\ssize \bullet$} at -12 1  -11 2  -10 3  -9 4  -8 5  -7 6  -6 7  -5 8 
   -4 9  -3 10  -2 11 -1 12  0 13 
   12 1  11 2  10 3  9 4  8 5  7 6  6 7  5 8  4 9  3 10  2 11  1 12 
   -10 1  -8 1  -6 1  -4 1  -2 1  0 1  2 1  4 1  6 1  8 1  10 1 /

\multiput{$\ssize \blacksquare$} at  -9 4  -9.5 2.5  -8.5 1.5  -6 1 
      9 4  9.5 2.5  8.5 1.5  6 1  -3 10  -1 11  1 11  3 10 /
\setsolid
\plot -3 10  -9 4  -9.5 2.5  -8.5 1.5  -6 1  6 1 /
\plot  3 10  9 4  9.5 2.5  8.5 1.5  6 1 /
\plot -3 10  -1 11  1 11  3 10 /
\setshadegrid span <.3mm>
\vshade -9.5 2.5 2.5  <z,z,,> -9 2 4 <z,z,,> -8.5 1.5 4.6 <z,z,,> 
        -6 1 7 <z,z,,> -3 1 10  <z,z,,> 
        -1 1 11  <z,z,,>
        1 1 11   <z,z,,>  3 1 10  <z,z,,> 6 1 7  <z,z,,> 8.5 1.5 4.6 <z,z,,> 
         9 2 4 <z,z,,> 9.5 2.5 2.5  /
\setshadegrid span <.8mm>
\vshade -12 1 1  <z,z,,> 0 1 13  <z,z,,> 12 1 1 /

\put{$\Cal S(15)$} at -9 12 
\setdashes <1mm> 
\plot -9 2  9 2  0 11  -9 2   /
\endpicture} at 7 -12
\endpicture} 
$$

\subsection{Tame and wild.}
\label{sec-sixteen-six}

Starting from a tame family, using expansion we get a wild family. But this increases 
the colevel $r$.

	\medskip
        
        {\bf (1)} Is there a wild family with $q \le 4$ ?
        In Section~\ref{sec-seven-five} we have seen that there is a wild
        family with $p=1$.

	\medskip
        
{\bf (2)} Is there a wild family with pr-vector $(7/5, 14/5)$ ?
          \medskip
          
{\bf (3)} What can we say about the region of rational points in $\mathbb T(n)$ which support
           wild families?

\subsection{Infinite-dimensional objects.}
\label{sec-sixteen-seven}

In this paper (but also in earlier ones) we have restricted the attention to 
finite-dimensional vector spaces. But it is of interest (also for applications)
to deal with arbitrary vector spaces. Thus, in general, the setting should be 
expanded as follows: One should look at pairs $(U,V)$, where $V$ is an arbitary,
not necessarily finite dimensional, vector space with a locally nilpotent operator,
and $U$ is an invariant subspace of $V$, thus $(U,V)$ is a direct limit of objects
in $\Cal S$. 

\subsection{The case $n = 6.$}
\label{sec-sixteen-eight}

Looking at \cite{RS1} and Section~\ref{sec-eleven} of the present paper, one may be inclined to believe
that the case $n=6$ is now well established. But this is definitely not the case.
As usually, such a boundary case between finite and wild type is of special interest,
since one hopes that it sheds light on the general behaviour. Here are some questions.
	\medskip
        
{\bf (1)} Is there a formula (or at least a reasonable estimate)
  for the number of
  indecomposable objects in $\Cal S(6)$ with fixed pr-vector?
  In which way do these numbers increase in the direction towards the center? Of course,
  it is sufficient to consider the central half-lines in $\mathbb T(6)$ which pass through one of the
  five vertices $(0,0),\ (0,1),\ (0,2),\ (0,8/3),$ and $(0,3).$
	\medskip
        
{\bf (2)} For any $t\ge 1,$ what is the support of the indecomposable objects 
   $X$ with fixed width $bX = t$ ?
	\medskip
        
{\bf (3)} The reference space $\mathbb T(6)$ does not provide any clue on
   central objects: They just vanish in the center $(2,2)$ (thus, we call it the black hole).
   What is the structure of the full subcategory of central objects?
   	\medskip
        
{\bf (4)} Are there central objects in $\Cal S(6)$ which are at the 
   boundary of a 6-tube or a 3-tube?
	\medskip
        
{\bf (5)} The reference space $\mathbb T(6)$ seems to provide not much information about the
   rationality index of the objects in $\Cal S(6)$. In which way can one
   combine the combinatorial information provided by $\mathbb T(6)$ with the information provided
   by the rationality index?

   \bigskip
\section{Final comments.}
\label{sec-seventeen}

\subsection{The mystery of boundary distance $1.$}
\label{sec-seventeen-one}

Let us draw the attention to the dichotomy between $p < 1$ and $p \ge 1;$ 
and similarly between $q \le n-1$ and $q > n-1.$
	\medskip
        
We have seen in the paper: 
On the one hand, there do exist (for $n\ge 9$)
infinite families of indecomposable objects $X$ in $\Cal S(n)$
with $qX = n\!-\!1.$ On the other hand, there are 
no non-trivial indecomposable objects $Y$ in $\Cal S(n)$ with 
$n\!-\!1 < qY$ (non-trivial means: $bY \ge 2$).
Thus, starting with $X = (U,V)$, it is impossible to make a proper enlargement 
$V \subseteq V'$ and $U \subseteq U'$ in $\Cal S(n)$ with $bV = bV'$ and keeping
indecomposability ...
	\smallskip
        
In contrast, look at the family $X_c$ with covering modules
$
 \smallmatrix
           &&  2&1  \cr
           1&2&3&3&2&1,\endsmallmatrix
$
such that the map $3\leftarrow 3$ is an isomorphism. There is 
a (unique) indecomposable object $Y$ with covering
$
 \smallmatrix
           &&  2&2  \cr
           1&2&3&3&2&1,\endsmallmatrix
$
and such that any $X_c$ is a submodule of $Y$. Thus, here we deal with
a situation that starting with $X = (U,V)$, we can enlarge
$U$ to $U'$ (inside $V$) and $Y = (U',V)$ is still indecomposable.
	\medskip
        
There are also examples of families $X = (U,V)$ where we
enlarge just $V$, not $U$: Start with $V = [7,4,2]$ and 
$U = [4,2]$ and enlarge $V$ to $V' = [7,5,2]$ ...
But as we have shown this cannot happen in case $qX = n-1$.

\subsection{The plus-construction.}
\label{sec-seventeen-two}

The plus-construction defined in Section~\ref{sec-nine-one} focuses the attention for $n\ge 6$ to 
an important property of the category $\Cal S(n)$: It provides for
any indecomposable object $X$ an arithmetical sequence of 
indecomposable objects which describes a move towards the center of $\mathbb T(n)$.
We may call this a process of concentration or amplification: 
There is the path $X \to X^+ \to X^{+2}\to \cdots$ along a ray
which has to be seen as a sequence of amplifications of $X$ and of 
moves towards the center.

But there is also the dual feature of decentration or shrinking: Namely,
using duality, the sectional path $X \to X^+$ yields
a corresponding sectional path
$\D(X^+) \to \D X,$ and actually we have $\D(X^+) = (\D X)^+$: Thus for 
any indecomposable object $Y$ there is an arithmetical sequence of 
indecomposable objects which moves away from the center of $\mathbb T(n)$ and
shows a shrinking of the objects: This concerns
the path $\cdots \to Y^{+2}\to Y^+ \to Y$ along a coray.

As we have seen in Section~\ref{sec-ten}, there is a difference between the
stable components and the principal component $\Cal P(n)$:
If we deal with a stable component $\Cal C$, the arithmetical sequence
defined by an object $X$, say supported by the half-line $H$, gives rise
to a corresponding arithmetical sequence defined again by an object of 
$\Cal C$ but now supported by the complementary half-line $H'$ of $H$.
Thus, we observe an oscillation behaviour on the central line $H\cup H'$.
In contrast, for the principal component $\Cal P(n)$, 
there are half-lines which occur in 
the half-line support in $\Cal P(n)$, where the complementary half-line does not 
play a role. 

Altogether, we should have the following picture in mind: There are
the half-lines $\mathbb H_\ell(n)$, and for $n\ge 7$ also the half-lines $\mathbb K_s(n)$
which support the concentration or amplification along rays inside $\Cal P(n)$, 
then, inside the
stable components, we observe the oscilation on complementary half-lines, and
finally we look again at the principal component with its arithmetical sequences 
in corays, moving away from the center of $\mathbb T(n)$, as a sign of 
decentration or shrinking.

This picture has some similarity to the flow in the categories $\mod K(m)$ with $m\ge 2$
(here $K(m)$ is the $m$-Kronecker algebra),
where we start in the preprojective component and end in the preinjective component,
passing in-between through the regular components. The regular components 
are of tree class $\mathbf A_\infty$, and we always 
use first a coray (which means decrease), then a ray (which means increase). 
Of course, the fact that $\Cal S(n)$ is a Frobenius category (where projective
objects and injective objects coincide) explains that the exceptional behaviour
of dealing with half-lines without complements is not devoted to two components
(as in the case of $\mod K(m)$), a preprojective and a preinjective one, 
but to the single principle component
which contains the indecomposable projective-injective objects. 
All components of $\Cal S(n)$ with $n\ge 6$ are of tree type 
$\mathbf A_\infty$; and clearly, the principal component combines features of 
a preprojective and a preinjective component.
We are convinced that the categories $\mod K(m)$ and $\Cal S(n)$ are
basic examples in representation theory and their similarities should be kept in mind.

\subsection{Combinatorics versus algebra.}
\label{sec-seventeen-three}

As a summary, we want to stress that the target of our investigations is to 
single out settings which are purely combinatorial. 

Dealing with
finite-dimensional representations of a $k$-algebra, the invariants which 
have to be described may or may not refer to the number of elements of the ground field $k$
or to properties of $k$.
Invariants which do not refer to $k$ may be said
to be {\it combinatorial} invariants; 
the remaining ones we will call {\it algebraic} invariants.
	
The algebraic invariants which we encounter 
usually refer just to the cardinality of $k$, sometimes also to the algebraic closure of $k$, 
but it seems that more information about the structure of $k$, about its characteristic or
about completeness properties never play a role. 

In order to facilitate the reading (say even for undergraduate students), we avoid the use
of field extensions. In particular, the  projective line $\mathbb P^1$ over $k$ is defined just
as the set of its rational points. We hope that the more advanced reader has no problem to
invoke also the non-rational points ... 

The combinatorics which we encounter leads first of all 
to formulae which involve binomial coeffcients
as well as Fibonacci numbers. It should not come as a surprise that also the exceptional
root systems, in particular the root system $\mathbf E_8$, play a prominent role. 
After all, we have seen already in \cite{RS1} and \cite{S1} in which way the root systems 
known from Lie theory come into play. 

Our detailed study of the case $n=6$ in \cite{RS1} as well as in the present paper relies on
the conviction that this tame border case between finite and wild type may be crucial for
a general understanding of the categories $\Cal S(n)$. But is this conviction justified? We
do not know! One should be aware of the fact that 
the categories $\Cal S(n)$ form what is called an even $\mathbf {ADE}$ chain, see for example
\cite{Bi}: The categories $\Cal S(1),\Cal S(2),\Cal S(3),\Cal S(4),\Cal S(5)$ 
are of tree type $\mathbf A_0, \mathbf A_2, \mathbf D_4, \mathbf E_6, \mathbf E_8$, respectively,
the category $\Cal S(6)$ is of tubular type $(6,3,2).$ One may expect that the extension of
this chain to the wild cases may shed some
light on the structure of the categories $\Cal S(n)$ with $n\ge 7.$

\subsection{The category $\Cal S(6)$ and the root system $\mathbf E_8$.}
\label{sec-seventeen-four}

Our detailled study of the category $\Cal S(6)$ in \cite{RS1} and in Section~\ref{sec-eleven} of
the present paper 
relies on the relationship between
$\Cal S(6)$ and the module category of the tubular algebra $\Theta$, and this 
furnishes a relationship between $\Cal S(6)$ and the root system $\mathbf E_8$. 

We need to distinguish between
the objects of $\Cal S(6)$ living in stable tubes and those
belonging to $\Cal P(6)$. 
The objects in $\Cal P(6)$ are known explicitly, thus one
can check any question directly. It is decisive to understand
the objects which live in stable tubes.
For every indecomposable object $X$ which lives in a stable tube, 
we have attached the element $\mathbf r(X)$ which is a root of 
$\chi_\Xi$ or the zero vector. Note that the roots of $\chi_\Xi$ 
form a root system of type $\mathbf E_8.$

This approach to distinguish 
between the objects of $\Cal S(6)$ living in stable tubes and those
belonging to $\Cal P(6)$,
has been used four times in our investigations. In
\cite{RS1}, we have discussed in this way the relationship between $u$ and $v$,
Theorem~\ref{theoremseven} deals in the same way with the relationship between $u$
and $b$. Then there is Theorem~\ref{theoremeight} which focuses the attention first 
to the 12 support lines and second to the nested sequence of standard triangles
which support the category.
To repeat: We always had to distinguish between
the objects of $\Cal S(6)$ living in stable tubes and those
belonging to $\Cal P(6)$. For an object $X$ in a
stable component, we used properties of $\mathbf r(X)$, and 
we could refer to the table provided in Appendix~\ref{app-A}. (Let us stress that we have used this
appendix only for dealing with objects of $\Cal S(6)$ which belong to stable tubes.) 
As we know, the remaining objects in $\Cal P(6)$ can be described very well, and we
have checked the required properties directly.
	\medskip
        
Let us insert a short remark, concerning the separation of $\Cal S(6)$ into $\Cal P(6)$ and
the stable tubes:
It is fascinating to see that $\Cal S(6)$ can be reconstructed by starting from the stable tubes of 
$\mod \Theta$, identifying the stable tubes in
$\Cal T_0$ and in $\Cal T_\infty,$ and then adding the principal component $\Cal P(6)$.
Note that the identified categories $\Cal T_0$ and $\Cal T_\infty$ together with 
$\Cal P(6)$ yield the objects with rationality index $0$.
Looking at the function $\mathbf r$ defined for all objects in stable tubes of $\Cal S(6)$,
it would be proper to attach to an object $X$ which belongs to a stable tube with
rationality index 0 not just the one element $\mathbf r(X)$, but the
{\bf pair} $\{\mathbf r(X),\mathbf r'(X)\}$, where 
$\mathbf r'(X)$ is obtained from $\mathbf r(X)$ by shifting the coefficients one step to the right.  
	\medskip
        
As we have mentioned, we have used known properties of the tubular algebra $\Theta$ in order
to get information about $\Cal S(6)$. But conversely, our investigation of $\Cal S(6)$ 
may provide a new insight into $\mod \Theta$, and this could be helpful for a better understanding 
of the module categories of tubular algebras.
It seems that the relevance of the 12 central lines in $\mathbb T(6)$ which we encounter
when dealing with $\Cal S(6)$, but which are of similar interest for $\mod \Theta$, 
was not yet noticed for $\Theta$. 
For any tubular algebra, there should exist a corresponding set
of planar lines which deserves to be investigated.

\subsection{The spider web for $n = 6.$}
\label{sec-seventeen-five}

Theorem~\ref{theoremeight} presents two essential features for the case $n = 6$,
namely the 12 lines and the nested sequence of standard triangles which provide the 
support of $\Cal S(6)$. There cannot be any doubt about the importance of
the 12 lines. However, there may be questions about the relevance of the triangles. 
Instead of the triangles one may also look at costandard triangles or at 
a corresponding nested sequence of hexagons. Actually, the authors had
numerous discussions which point of view should be preferable. Of course, 
the triangles (the standard and the costandard ones) as well as the hexagons,
all these shapes are definitely of interest, but for simplicity, 
we wanted to focus the first attention 
to a single shape.

In order to be precise: The {\it hexagons} $H_d$ which play a role are given as the
boundary of the convex hull of the $\Sigma_3$-orbit of a vertex of the form $(d,2)$,
where $1\le d < 2$ (it is a subset of the union of the standard triangle $\Delta_d$
and the costandard triangle $\nabla_{4-d}$). The configuration given by 
the 12 lines $\mathbb L(6)$ and the hexagons $H_d$ which support indecomposable objects
in $\Cal S(6)$ seems to look like an impressive spider web. 
	
The first hexagon $H_1$ has been featured prominently at the end of Section~\ref{sec-one-eight}.
Here, all the
vertices on $H_1$ which belong to one of the 12 lines support indecomposable objects
of $\Cal S(6)$. 
But this is no longer true for some of the further hexagons $H_d$ which occur.

Next, one should look at the hexagon $H_{5/4}$ as shown in Sections \ref{sec-fifteen-two}--(f) and
\ref{sec-eight-five}.
Here we encounter 
a missing symmetry with respect to the reflection $(p,r) \mapsto (4-p,p+r-2):$ there is
a unique indecomposable object in $\Cal S(6)$ supported by $(9/4,5/4)$, but there are
two indecomposable objects supported by $(7/4,6/4).$ 
Similarly, looking at $H_{4/3}$, 
Appendix~\ref{app-B-two} shows that the number of tripickets with support $(5/3,5/3)$ is 2, whereas
the number of tripickets with support $(7/3,4/3)$ is 3.
Obviously, one has to be aware that not all symmetries of the hexagon are realised 
in our setting. 

In-between $H_1$ and $H_{5/4}$, a further hexagon is needed, namely $H_{6/5}$
in order to accomodate the 6 pentapickets exhibited in Appendix~\ref{app-B-six}, for example the two
pentapickets
with pr-vector $(8/5,8/5).$ But note that only three vertices of this
hexagon $H_{6/5}$
support indecomposable objects of $\Cal S(6)$, namely the corners of $\Delta_{8/5}.$
This may be a hint that the use of the standard triangles is more appropriate than
that of the hexagons. 

We finally decided to stress the relevance of the standard triangles $\Delta_d$ and to
postpone a study of the hexagons for later investigations.

Some
\phantomsection{hexagons}
$H_d$ (with $d = 1,\ \frac 65,\ \frac 54,\ \frac 43,\ \frac75,\ 
\frac{10}7,\ \frac32$\, ), and the position of the
indecomposables on $H_d$ with  
$d = 1,\ \frac 65,\ \frac 54\,$:

\addcontentsline{lof}{subsection}{The nested hexagons (“spider web”).}

$$  
{\beginpicture
   \setcoordinatesystem units <3.5cm,6.062cm>
\setdots <.4mm>
\multiput{} at  -1 .9  0 2.05 /
\put{$\blacksquare$} at 0 2 
\multiput{$\bullet$} at 1 1  0 1  -1  1  -1.5 1.5  -2 2 
    0.333 1  -0.333 1  -1.333 1.333  -1.667 1.667 /

\plot -1.1  .9  0 2 /
\plot -1 1  -2 2 /

\setsolid
\plot -2 2  -1 1  1 1 /

\plot -1.1  .9  0 2 /
\plot 0 2  0 .9 /
\plot -2 2  0 2 /
\plot -2 1.3333  0 2 /
\plot 0 2  1.1 .9 /

\setdashes <1mm>
\plot 1 1  1.13 1 /
\put{$1$} at 1.18 1 

\setsolid
\plot -1 2  -.5 1.5  .5 1.5 /

\plot -2 2  -1.95  2.05  /

\setsolid
\plot -1 2  -.95 2.05  /
\setdashes <1mm>
\plot .5 1.5  .67 1.5 /
\put{$3/2$} at .79 1.5 

\setsolid
\plot -1.5 2  -1.45 2.05  /
\plot 0.75 1.25  -.75 1.25  -1.5 2 /
\setdashes <1mm>
\plot .75 1.25  1.05 1.25 /
\multiput{$\bullet$} at -1.5 2  -1.25 1.75 -1 1.5
   -0.75 1.25  -.25 1.25   .25 1.25  0.75 1.25 /
\put{$5/4$} at 1.15 1.25

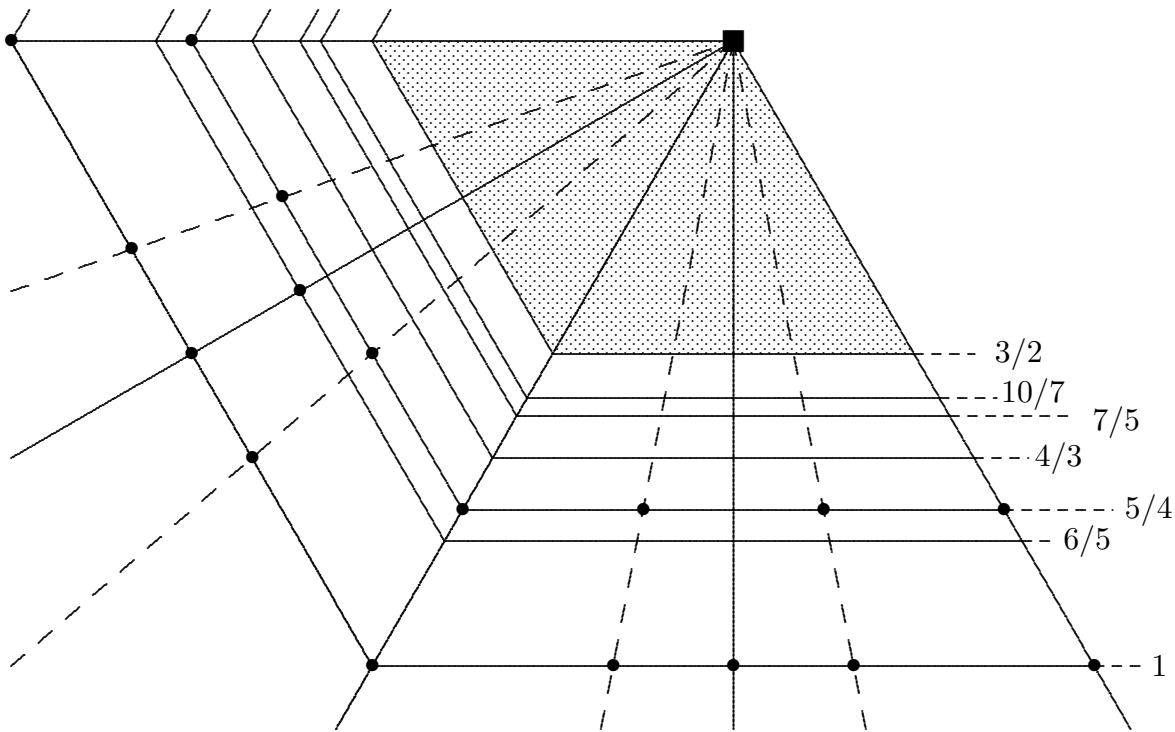
\setsolid
\plot -1.333 2  -1.283  2.05  /
\plot 0.6667 1.333  -.6667  1.333   -1.333 2 /
\setdashes <1mm>
\plot 0.6667 1.333 0.82 1.333 /
\multiput{$\bullet$} at /
\put{$4/3$} at .9 1.333 

\setsolid
\plot -1.2 2  -1.15  2.05  /
\plot 0.6 1.4  -.6  1.4   -1.2 2 /
\setdashes <1mm>
\plot 0.6 1.4  0.95 1.4 /
\multiput{$\bullet$} at /
\put{$7/5$} at 1.06 1.39

\setsolid
\plot -1.1428 2  -1.0923  2.05  /
\plot  0.5714 1.4286  -.5714 1.4286  -1.1428 2  /
\setdashes <1mm>
\plot   0.5714 1.4286   0.73 1.4286 /
\multiput{$\bullet$} at /
\put{$10/7$} at .83 1.435

\setshadegrid span <.5mm>
\vshade -1 2 2  <z,z,,> -.5 1.5 2 <z,z,,> 0 1.5 2 <z,z,,> .5 1.5  1.5 /

\setsolid 
\plot .8 1.2  -.8 1.2  -1.6 2  -1.55 2.05 /
\setdashes <1mm>
\plot .8 1.2  .88 1.2  /
\multiput{$\bullet$} at -1.2 1.6 /
\put{$6/5$} at .98 1.2

\setdashes <2mm>
\plot -2 1.6  0 2 /
\plot -2 1  0 2 /
\plot -.3667 .9  0 2 /
\plot  .3667 .9  0 2 /
\endpicture}
$$

Following the discussion in Section~\ref{sec-seventeen-four},
the spider web of $\Cal S(6)$ should be seen as
a kind of shadow of the root system of $\chi_\Theta$.
It should be possible to draw similar spider webs also for tubular
algebras. 

\subsection{The triangle $\mathbb T(n)$.}
\label{sec-seventeen-six}

The title of the paper refers to three 
invariants: level $p$, mean $q$ and colevel $r$, as well as to the use of these
numbers by invoking the triangle $\mathbb T(n)$. Our visualization of $\mathbb T(n)$
focuses the attention to the level: Level $0$ is the horizontal basis of the triangle,
the increase of the level corresponds to going upwards in the triangle. 
However, in this way the important duality functor $\D$ is not seen well, being 
the reflextion at an inclined line. 

This flaw could be remedied by a rotation of $\mathbb T(n)$:
By drawing the origin $(0,0)$ of the triangle at the bottom, so that going upwards corresponds
to an increase of $q$, see the following picture on the left.
Then the duality $\D$ is seen well as
the reflection through the central vertical line. 
$$
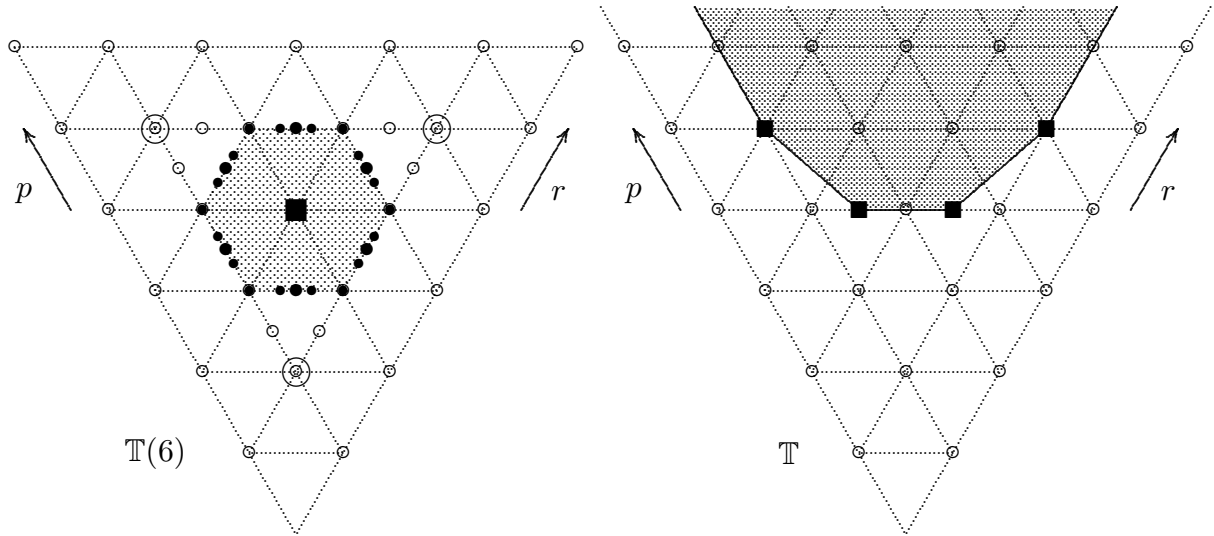
{\beginpicture
   \setcoordinatesystem units <.46187cm,.8cm>
\put{\beginpicture
\multiput{} at -6 0  6 6.5 /
\arr{-4.8 4}{-5.8 5}
\put{$p$} at -5.8 4.2 
\arr{4.8 4}{5.8 5}
\put{$r$} at 5.6 4.2 
\setdots <.5mm> 
\plot -6 6  6 6  0 0  -6 6 /
\plot -3 5  3 5  0 2  -3 5 /
\plot -1 5  1 3  -1 3  1 5  2 4  -2 4  -1 5 /

\plot -5 5  -4 6  -3 5  -2 6  -1 5  0 6  1 5  2 6  3 5  4 6  5 5 /
\plot -5 5  -3 5  -4 4  -2 4  -3 3  -1 3  -2 2  0 2  -1 1  1 1 /
\plot  5 5   3 5   4 4   2 4   3 3   1 3   2 2  0 2   1 1  /
\multiput{$\circ$} at 
     -4 6  -2 6  0 6  2 6  4 6  
     -3 5  -1 5  1 5  3 5 
     -2 4  0 4  2 4
     -1 3  1 3
      0 2   6 6   -6 6 
       -5 5  -4 4  -3 3  -2 2  -1 1  1 1  2 2  3 3  4 4  5 5  /
\multiput{$\circ$} at -.5 2.5  .5 2.5  -2 5  2 5  -2.5 4.5  2.5 4.5 /
  \multiput{$\bigcirc$} at -3 5  3 5  0 2  /

\multiput{$\ssize\bullet$} at -.333 5  .333 5  1.333 3.333  1.667 3.667  -1.333 3.333  -1.667 3.667 
     -2 4  -1 5  -1 3  1 5  1 3  2 4 -.333 3  .333 3 
      1.333 4.667   1.667 4.333 -1.333 4.667  -1.667 4.333 /
\multiput{$\bullet$} at 0 5  0 3  -1.5 4.5  1.5 4.5  -1.5 3.5  1.5 3.5  /
\put{$\blacksquare$} at 0 4
   \setshadegrid span <.4mm>
   \vshade -2 4 4  <,z,,>  -1 3 5   <z,z,,>  1 3 5  <z,z,,>  2 4 4  /
\put{$\mathbb T(6)\:$} at -5 2
\endpicture} at 0 0 
\put{\beginpicture
\multiput{} at -6 0  6 6.5 /
\arr{-4.8 4}{-5.8 5}
\put{$p$} at -5.8 4.2 
\arr{4.8 4}{5.8 5}
\put{$r$} at 5.6 4.2 
\setdots <.5mm> 
\plot -6 6  6 6  0 0  -6 6 /
\plot -3 5  3 5  0 2  -3 5 /
\plot -1 5  1 3  -1 3  1 5  2 4  -2 4  -1 5 /

\plot -5 5  -4 6  -3 5  -2 6  -1 5  0 6  1 5  2 6  3 5  4 6  5 5 /
\plot -5 5  -3 5  -4 4  -2 4  -3 3  -1 3  -2 2  0 2  -1 1  1 1 /
\plot  5 5   3 5   4 4   2 4   3 3   1 3   2 2  0 2   1 1  /

\plot -6.5 6.5  -6 6 /
\plot -4.5 6.5  -4 6  -3.5 6.5 /
\plot -2.5 6.5  -2 6  -1.5 6.5 /
\plot -0.5 6.5   0 6   0.5 6.5 /
\plot  6.5 6.5   6 6 /
\plot  4.5 6.5   4 6   3.5 6.5 /
\plot  2.5 6.5   2 6   1.5 6.5 /
\multiput{$\circ$} at 
     -4 6  -2 6  0 6  2 6  4 6  
     -3 5  -1 5  1 5  3 5 
     -2 4  0 4  2 4
     -1 3  1 3
      0 2   6 6   -6 6 
       -5 5  -4 4  -3 3  -2 2  -1 1  1 1  2 2  3 3  4 4  5 5  /
  
\multiput{$\ssize \blacksquare$} at -1 4  1 4  -3 5  3 5 /
   \setshadegrid span <.3mm>
   \vshade -4.5 6.5 6.5  <,z,,>  -3 5 6.5   <z,z,,>  -1 4 6.5  <z,z,,>  1 4 6.5  <z,z,,>
     3 5 6.5  <z,z,,> 4.5 6.5 6.5 /
\put{$\mathbb T\:$} at -4.5 2 
\setsolid
\plot -4.5 6.5  -3 5  -1 4  1 4  3 5  4.5 6.5 /
\endpicture} at 13 0 
\endpicture}
$$

The rotated triangles $\mathbb T(n)$ seem to be of special interest when
looking at the whole category $\Cal S = \bigcup_n \Cal S(n)$ and the corresponding (now infinite)
triangle $\mathbb T = \bigcup_n \mathbb T(n)$. On the right, we have drawn 
part of the rotated triangle $\mathbb T$,
with the region of suspected BTh-vectors.

\subsection{The rotation of $\mathbb T(n)$.}
\label{sec-seventeen-seven}

Whatever optical presentation of $\mathbb T(n)$ the reader prefers (to
stress $p = 0$ as the basis of $\mathbb T(n)$ as usually done in the paper, or 
else to use the 
cone version with vertical increase of the mean, as mentioned in Section~\ref{sec-seventeen-six}, 
it is the rotation by 120$^\circ$
which really is the outstanding feature of the reference space $\mathbb T(n)$.
This surprising rotation was discovered by Schmidmeier a long time ago, it has been 
the fundament of our description of the Auslander-Reiten translation \cite{RS2}
and now has led us to introduce the set $\mathbb E(n)$, see Section~\ref{sec-three-one}. 
As we have mentioned, it is a main ingredient for nearly all the results 
of the present paper. 
A look at the hexagon pictures in Section~\ref{sec-fifteen} and in Appendix~\ref{app-B}
shows that this hidden symmetry relates objects which 
on a first view seem to be really different. 
Clearly, the rotation 
should be helpful also for any further investigation of $\Cal S(n)$. 

There is the following surprising observation:
{\it If $X$ is a reduced object of $\Cal S(n)$, 
then $pX = n - q(\tau_n^2X)$.} Why is this relationship between
$p$ and $q$ so astonishing? 
Whereas for the calculation of $pX$, one  needs
to know both the global space $VX$ and the subspace $UX$, 
the value of $qY$ only depends on the global space $VY$.
Namely, we have $pX = uX/bX = |U|/bV,$ but $qX = vX/bX = |V|/bV.$
Looking at the triangle $\mathbb T(n)$, there
are three barycentric coordinates $p,\ r,\ \omega,$ with 
$\omega = n-p-q$, and the existence of the rotation $\rho = \tau^2$ 
shows that the three coordinates are on an equal footing,
in sharp contrast to an intuitive interpretation of the 
invariants $p$ and $r$ on the one hand, and $q$ on the other. 

\medskip
{\bf Acknowledgements.}
The
\phantomsection{authors}
\addcontentsline{toc}{section}{Acknowledgements.}
wish to thank the Managing Editor
and the three referees for thoughtful comments which have
led to substantial improvements of the paper.
They are grateful for advice from William Crawley-Boevey
regarding the folding algorithm for partitions
which led to a better understanding of embeddings
with a cyclic submodule in Section~\ref{sec-thirteen}.
The first author wants to thank Piotr Dowbor
for helpful advice when he prepared Appendix~\ref{app-D}.
The second author would like to thank Bernhard Keller for
his guidance regarding the functor $\tau^2$ and Shijie Zhu
for interesting discussions about regions of BTh-vectors.

\medskip
{\bf Conflict of Interest Statement.}
On behalf of all authors, the corresponding author states
that there is no conflict of interest.

	\bigskip
\renewcommand{\refname}{Essential References (for proofs and illustrations)}

\bigskip
\renewcommand{\refname}{Additional References (history, relevance)}

\bigskip\bigskip
\vfill\eject


\def\roote#1#2#3#4#5#6#7#8{\beginpicture \setcoordinatesystem units <.12cm,.18cm>
                                         \multiput{} at 0 0  6 1 /
                                         \put{\kl #1} at 6 0
                                         \put{\kl #2} at 5 0
                                         \put{\kl #3} at 4 0
                                         \put{\kl #4} at 3 0
                                         \put{\kl #5} at 2 0
                                         \put{\kl #6} at 1 0
                                         \put{\kl #7} at 0 0
                                         \put{\kl #8} at 2 1 
                                         \endpicture}
\def\roott#1#2#3#4#5#6#7#8{\beginpicture \setcoordinatesystem units <.12cm,.18cm>
                                         \multiput{} at 0 0  5 1 /
                                         \put{\kl #1} at 0 0
                                         \put{\kl #2} at 1 0
                                         \put{\kl #3} at 2 0
                                         \put{\kl #4} at 3 0
                                         \put{\kl #5} at 4 0
                                         \put{\kl #6} at 5 0
                                         \put{\kl #7} at 1 1
                                         \put{\kl #8} at 2 1 
                                         \endpicture}
\def\uwbvec#1#2#3{$\frac{#1\,|\,#2}{#3}$}
\def\isphi#1#2#3{${\ss \frac{#1}{#2}= #3}$}
\def\ubvec#1#2#3#4#5#6{$#3,#4$}

\newcounter{sectionapp}
\renewcommand\thesection{\Alph{sectionapp}}
\renewcommand\thesubsection{(\Alph{sectionapp}.\arabic{subsection})}
\renewcommand\thetheorem{\Alph{sectionapp}.\arabic{theorem}}

\centerline{\Gross \phantomsection{Appendices}.}
\addcontentsline{toc}{part}{Appendices.}

\setcounter{section}{19}
\setcounter{sectionapp}{1}
\section{The roots of $\mathbf E_8$.}
\label{app-A}

\bigskip
by Markus Schmidmeier

\bigskip
	\bigskip
In the root table Appendix~\ref{app-A-one}, we list for each positive root  of $\mathbf E_8$ some data which
describe the corresponding representations of the algebra $\Theta$ in Section~\ref{sec-eleven-one}.
Some examples are given in Appendix~\ref{app-A-two}.

\subsection{The root table.}
\label{app-A-one}

Consider the following two quivers, $Q_{\Xi}$, which is a restriction of $Q_\Theta$ in
Section~\ref{sec-eleven-one}, and $Q_H$,
which is the quiver of a path algebra of type $\mathbf E_8$. 
$$
\beginpicture\setcoordinatesystem units <.7cm,.7cm>
        \multiput{} at 0 0  5 1 /
        \put{$Q_{\Xi}\:$} at -1 .5
        \put{$\ssize 1$} at 0 0
        \put{$\ssize 2$} at 1 0
        \put{$\ssize 3$} at 2 0
        \put{$\ssize 4$} at 3 0
        \put{$\ssize 5$} at 4 0
        \put{$\ssize 6$} at 5 0
        \put{$\ssize 2'$} at 1 1 
        \put{$\ssize 3'$} at 2 1
        \arr{.7 0}{.3 0}
        \arr{1.7 0}{1.3 0}
        \arr{2.7 0}{2.3 0}
        \arr{3.7 0}{3.3 0}
        \arr{4.7 0}{4.3 0}
        \arr{1.7 1}{1.3 1}
        \arr{1 .7}{1 .3}
        \arr{2 .7}{2 .3}
        \put{$\ssize\alpha$} at 1.5 1.3
        \put{$\ssize\beta$} at 2.3 .5
        \setdots<2pt>
        \plot 1.7 .7  1.3 .3 /
        \endpicture
\qquad\qquad
\beginpicture\setcoordinatesystem units <.7cm,.7cm>
        \multiput{} at 0 0  6 1 /
        \put{$Q_{H}\:$} at -1 .5
        \put{$\ssize 1$} at 0 0 
        \put{$\ssize 2$} at 1 0 
        \put{$\ssize \bar3$} at 2 0 
        \put{$\ssize 3$} at 3 0 
        \put{$\ssize 4$} at 4 0 
        \put{$\ssize 5$} at 5 0 
        \put{$\ssize 6$} at 6 0 
        \put{$\ssize 2'$} at 2 1
        \arr{.7 0}{.3 0}
        \arr{1.7 0}{1.3 0}
        \arr{2.7 0}{2.3 0}
        \arr{3.7 0}{3.3 0}
        \arr{4.7 0}{4.3 0}
        \arr{5.7 0}{5.3 0}
        \arr{2 .7}{2 .3}
        \put{$\ssize\gamma$} at 1.7 .5
        \put{$\ssize\delta$} at 2.5 -.3
        \endpicture
        $$

The indecomposable representations $M$ for $\Xi$ for which the map $(M_\alpha,M_\beta)^t:M_{3'}\to M_{2'}\oplus M_3$ is injective
are in one-to-one correspondence to the indecomposable
representations $N$ for $H$ for which the map $(N_\gamma,N_\delta):N_{2'}\oplus N_3\to N_{\bar 3}$ is surjective.
To see this, form the pushout of the commutative square.
All indecomposable representations for $\Xi$ are obtained from roots for $\mathbf E_8$,
for the indecomposable representations of $\Theta$ see Section~\ref{sec-eleven-one}.

\medskip
In the table we list for each of the 120 positive roots of $\mathbf E_8$:

\medskip
\begin{itemize}[leftmargin=3em]
\item[(1)] The number of the root as listed on \cite{W}.

\medskip
\item[(2)] The root for $\mathbf E_8$.  For roots 5, 13, and 20,
           the corresponding map $(N_\gamma,N_\delta)$ is not onto.

\medskip
\item[(3)] The corresponding root for $\chi_\Xi$.  
           The roots for 5, 13, and 20 are mixed roots for $\chi_\Xi$,
           we abbreviate ``$-1$'' by ``$-$''. (Note that given an indecomposable 
           module $M$ in $'\Cal D$, the root $\mathbf r(M)$ of $\chi_\Xi$ 
           does not have to be positive, but may be negative or mixed.
            Of course, the negative of a mixed root is again mixed.)

\medskip
\item[(4)] For a root $\mathbf r$ for $\chi_\Xi$, the formal uwb-vector:
$$u=\mathbf r_1+\mathbf r_{2'}+\mathbf r_{3'}, \quad w=\sum_{i=1}^6\mathbf r_i-u, \quad b=\mathbf r_3$$

\medskip
\item[(5)] The slope $\phi$ of the corresponding central half-line, presented in the
following format:
$$\frac{u-2b}{w-2b}=\phi$$ 

\medskip\item[(6)] For non-central objects, the type of the corresponding central half-line.

\smallskip
\begin{itemize}[leftmargin=3em]
\item[$\mathbb P$:]
The half-line is parallel to one of the coordinate axes:  The lines
$p=2$, $r=2$, $6-q=2$ have slope $0$, $\infty$, and $-1$ respectively.

\smallskip \item[$\mathbb D$:] The diagonal lines
$p=r$, $r=6-q$, $p=6-q$ have slope $1$, $-2$, $-\frac12$ respectively.

\smallskip\item[$\mathbb H$:]
The object lies on a central half-line in $\mathbb H(6)$,
its slope is $\frac12$, $2$, $-3$, $-\frac32$, $-\frac23$ or $-\frac13$

\smallskip\item[c:] The object is central.
\end{itemize}
\medskip\item[] The symbols $\mathbb D$ and $\mathbb H$
                   carry a subscript, $s$ or $\ell$,
                   depending on whether objects corresponding to the root occur
                   on the short or the long half-line 
                   of the central line.
                   In this case, and also when the symbol is $\mathbb P$,
                   the negative root will give rise
                   to objects on the opposite side
                   of the same central line.
\medskip\item[(7)] Two numbers $r_\Delta,r_\nabla$,
                   defined as follows.
                   $$r_\Delta=\max\{\,2b-u,\,v-4b,\,2b-w\}\quad\text{and}\quad
                   r_\nabla=\max\{\,u-2b,\,4b-v,\,w-2b\}$$
                   are computed from the formal uwb-vector in column (4).
\smallskip\item[]  Suppose $X\in\Cal S(n)$ occurs in a stable tube and has (positive)
                   root $\mathbf r=\mathbf r(X)$.  We can compute
\begin{itemize}[leftmargin=3em]
\smallskip\item[$\bullet$] the in-radius of the standard ($\Delta$-)triangle through $\pr X$
                   as $\displaystyle \frac{r_\Delta}{bX}$;
\smallskip\item[$\bullet$]  the boundary distance as $\displaystyle dX=2-\frac{r_\Delta}{bX}$,
                                so $\pr X$ lies on $\Delta_{dX}$;
\smallskip\item[$\bullet$]  the in-radius of the costandard ($\nabla$-)triangle as
                   $\displaystyle\frac{r_\nabla}{bX}$
                   (so $\pr X$ occurs on $\displaystyle\nabla_{2+r_\nabla/bX}$,
                   using notation from Section~\ref{sec-ten-seven});
\smallskip\item[$\bullet$] the primitive pair as $(u,b)=(2bX-r_\Delta,bX)$
                   (using notation from 11.8).
\smallskip\item[$\bullet$] and the in-radius of the hexagon
                   as $\displaystyle r_H(X)=\frac{\max\{r_\Delta,r_\nabla\}}{bX}$
                   (so $\pr X$ lies on $H_{2-r_H(X)}$, see Section~\ref{sec-seventeen-four}).
                   \end{itemize}
\medskip\item[]                   
                   If $X$ occurs in a stable tube and has (negative) root
                   $\mathbf r(X)=-\mathbf r$, the five above quantities
                   are computed as $\displaystyle \frac{r_\nabla}{bX}$ (in-radius of the standard
                   triangle),
                   $\displaystyle dX=2-\frac{r_\nabla}{bX}$ (boundary distance),
                   $\displaystyle \frac{r_\Delta}{bX}$ (in-radius of the costandard triangle),
                   $(2bX-r_\nabla,bX)$ (primitive pair),
                   and $r_H(X)$ (in-radius of the hexagon), respectively.
\end{itemize}

\newpage
\noindent Table of
\phantomsection{roots}
\addcontentsline{lof}{subsection}{The table of the roots of $\Bbb E_8$.}%
of $\mathbf E_8$ and $\chi_\Xi$ --- Roots 1 to 60

$$\beginpicture
        \setcoordinatesystem units <1.05cm,.58cm>
        \multiput{} at 0.4 0  2 30 /
        \put{1} at 0.4 29 
        \put{\roote10000000} at 1.2 29
        \put{\roott00000100} at 2.1 29
        \put{\uwbvec010} at 3 29
        \put{\isphi010} at 4 29
        \put{$\mathbb P$} at 4.8 29
        \put{\ubvec011133} at 5.3 29
        \put{2} at 0.4 28
        \put{\roote01000000} at 1.2 28
        \put{\roott00001000} at 2.1 28
        \put{\uwbvec010} at 3 28
        \put{\isphi010} at 4 28
        \put{$\mathbb P$} at 4.8 28
        \put{\ubvec011133} at 5.3 28
        \put{3}      at 0.4 27
        \put{\roote00100000} at 1.2 27
        \put{\roott00010000} at 2.1 27
        \put{\uwbvec010}   at 3 27
        \put{\isphi010} at 4 27
        \put{$\mathbb P$}     at 4.8 27
        \put{\ubvec011133} at 5.3 27
        \put{4}      at 0.4 26
        \put{\roote00010000} at 1.2 26
        \put{\roott00100001} at 2.1 26
        \put{\uwbvec101}   at 3 26
        \put{\isphi{-1}{-2}{\frac12}} at 4 26
        \put{$\mathbb H_\ell$}     at 4.8 26
        \put{\ubvec{-1}{-3}2315} at 5.3 26
        \put{5}      at 0.4 25
        \put{\roote00001000} at 1.2 25
        \put{\roott0000000-} at 2.1 25
        \put{\uwbvec{-1}10}   at 3 25
        \put{\isphi{-1}1{-1}} at 4 25
        \put{$\mathbb P$}     at 4.8 25
        \put{\ubvec{-1}01133} at 5.3 25
        \put{6}      at 0.4 24
        \put{\roote00000100} at 1.2 24
        \put{\roott01000000} at 2.1 24
        \put{\uwbvec010}   at 3 24
        \put{\isphi010} at 4 24
        \put{$\mathbb P$}     at 4.8 24
        \put{\ubvec011133} at 5.3 24
        \put{7}      at 0.4 23
        \put{\roote00000010} at 1.2 23
        \put{\roott10000000} at 2.1 23
        \put{\uwbvec100}   at 3 23
        \put{\isphi10\infty} at 4 23
        \put{$\mathbb P$}     at 4.8 23
        \put{\ubvec111133} at 5.3 23
        \put{8}      at 0.4 22
        \put{\roote00000001} at 1.2 22
        \put{\roott00000011} at 2.1 22
        \put{\uwbvec2{-2}0}   at 3 22
        \put{\isphi2{-2}{-1}} at 4 22
        \put{$\mathbb P$}     at 4.8 22
        \put{\ubvec202233} at 5.3 22
        \put{9}      at 0.4 21
        \put{\roote11000000} at 1.2 21
        \put{\roott00001100} at 2.1 21
        \put{\uwbvec020}   at 3 21
        \put{\isphi020} at 4 21
        \put{$\mathbb P$}     at 4.8 21
        \put{\ubvec022233} at 5.3 21
        \put{10}      at 0.4 20
        \put{\roote01100000} at 1.2 20
        \put{\roott00011000} at 2.1 20
        \put{\uwbvec020}   at 3 20
        \put{\isphi020} at 4 20
        \put{$\mathbb P$}     at 4.8 20
        \put{\ubvec022233} at 5.3 20
        \put{11}      at 0.4 19
        \put{\roote00110000} at 1.2 19
        \put{\roott00110001} at 2.1 19
        \put{\uwbvec111}   at 3 19
        \put{\isphi{-1}{-1}1} at 4 19
        \put{$\mathbb D_\ell$}     at 4.8 19
        \put{\ubvec{-1}{-2}1212} at 5.3 19
        \put{12}      at 0.4 18
        \put{\roote00011000} at 1.2 18
        \put{\roott00100000} at 2.1 18
        \put{\uwbvec011}   at 3 18
        \put{\isphi{-2}{-1}2} at 4 18
        \put{$\mathbb H_\ell$}     at 4.8 18
        \put{\ubvec{-2}{-3}2315} at 5.3 18
        \put{13}      at 0.4 17
        \put{\roote00001100} at 1.2 17
        \put{\roott0100000-} at 2.1 17
        \put{\uwbvec{-1}20}   at 3 17
        \put{\isphi{-1}2{-\frac12}} at 4 17
        \put{$\mathbb D_\ell$}     at 4.8 17
        \put{\ubvec{-1}11233} at 5.3 17
        \put{14}      at 0.4 16
        \put{\roote00001001} at 1.2 16
        \put{\roott00000010} at 2.1 16
        \put{\uwbvec1{-1}0}   at 3 16
        \put{\isphi1{-1}{-1}} at 4 16
        \put{$\mathbb P$}     at 4.8 16
        \put{\ubvec101133} at 5.3 16
        \put{15}      at 0.4 15
        \put{\roote00000110} at 1.2 15
        \put{\roott11000000} at 2.1 15
        \put{\uwbvec110}   at 3 15
        \put{\isphi111} at 4 15
        \put{$\mathbb D_s$}     at 4.8 15
        \put{\ubvec122133} at 5.3 15
        \put{16}      at 0.4 14
        \put{\roote11100000} at 1.2 14
        \put{\roott00011100} at 2.1 14
        \put{\uwbvec030}   at 3 14
        \put{\isphi030} at 4 14
        \put{$\mathbb P$}     at 4.8 14
        \put{\ubvec033333} at 5.3 14
        \put{17}      at 0.4 13
        \put{\roote01110000} at 1.2 13
        \put{\roott00111001} at 2.1 13
        \put{\uwbvec121}   at 3 13
        \put{\isphi{-1}0\infty} at 4 13
        \put{$\mathbb P$}     at 4.8 13
        \put{\ubvec{-1}{-1}1112} at 5.3 13
        \put{18}      at 0.4 12
        \put{\roote00111000} at 1.2 12
        \put{\roott00110000} at 2.1 12
        \put{\uwbvec021}   at 3 12
        \put{\isphi{-2}0\infty} at 4 12
        \put{$\mathbb P$}     at 4.8 12
        \put{\ubvec{-2}{-2}2212} at 5.3 12
        \put{19}      at 0.4 11
        \put{\roote00011100} at 1.2 11
        \put{\roott01100000} at 2.1 11
        \put{\uwbvec021}   at 3 11
        \put{\isphi{-2}0\infty} at 4 11
        \put{$\mathbb P$}     at 4.8 11
        \put{\ubvec{-2}{-2}2212} at 5.3 11
        \put{20}      at 0.4 10
        \put{\roote00001110} at 1.2 10
        \put{\roott1100000-} at 2.1 10
        \put{\uwbvec020}   at 3 10
        \put{\isphi020} at 4 10
        \put{$\mathbb P$}     at 4.8 10
        \put{\ubvec022233} at 5.3 10
        \put{21}      at 0.4 9
        \put{\roote00011001} at 1.2 9
        \put{\roott00100011} at 2.1 9
        \put{\uwbvec2{-1}1}   at 3 9
        \put{\isphi0{-3}0} at 4 9
        \put{$\mathbb P$}     at 4.8 9
        \put{\ubvec0{-3}3345} at 5.3 9
        \put{22}      at 0.4 8
        \put{\roote00001101} at 1.2 8
        \put{\roott01000010} at 2.1 8
        \put{\uwbvec100}   at 3 8
        \put{\isphi10\infty} at 4 8
        \put{$\mathbb P$}     at 4.8 8
        \put{\ubvec111133} at 5.3 8
        \put{23}      at 0.4 7
        \put{\roote00011101} at 1.2 7
        \put{\roott01100011} at 2.1 7
        \put{\uwbvec201}   at 3 7
        \put{\isphi0{-2}0} at 4 7
        \put{$\mathbb P$}     at 4.8 7
        \put{\ubvec0{-2}2212} at 5.3 7
        \put{24}      at 0.4 6
        \put{\roote11110000} at 1.2 6
        \put{\roott00111101} at 2.1 6
        \put{\uwbvec131}   at 3 6
        \put{\isphi{-1}1{-1}} at 4 6
        \put{$\mathbb P$}     at 4.8 6
        \put{\ubvec{-1}01112} at 5.3 6
        \put{25}      at 0.4 5
        \put{\roote01111000} at 1.2 5
        \put{\roott00111000} at 2.1 5
        \put{\uwbvec031}   at 3 5
        \put{\isphi{-2}1{-2}} at 4 5
        \put{$\mathbb D_s$}     at 4.8 5
        \put{\ubvec{-2}{-1}2112} at 5.3 5
        \put{26}      at 0.4 4
        \put{\roote00111100} at 1.2 4
        \put{\roott01110000} at 2.1 4
        \put{\uwbvec031}   at 3 4
        \put{\isphi{-2}1{-2}} at 4 4
        \put{$\mathbb D_s$}     at 4.8 4
        \put{\ubvec{-2}{-1}2112} at 5.3 4
        \put{27}      at 0.4 3
        \put{\roote00111001} at 1.2 3
        \put{\roott00110011} at 2.1 3
        \put{\uwbvec201}   at 3 3
        \put{\isphi0{-2}0} at 4 3
        \put{$\mathbb P$}     at 4.8 3
        \put{\ubvec0{-2}2242} at 5.3 3
        \put{28}      at 0.4 2
        \put{\roote00011110} at 1.2 2
        \put{\roott11100000} at 2.1 2
        \put{\uwbvec121}   at 3 2
        \put{\isphi{-1}0\infty} at 4 2
        \put{$\mathbb P$}     at 4.8 2
        \put{\ubvec{-1}{-1}1112} at 5.3 2
        \put{29}      at 0.4 1
        \put{\roote00001111} at 1.2 1
        \put{\roott11000010} at 2.1 1
        \put{\uwbvec200}   at 3 1
        \put{\isphi20\infty} at 4 1
        \put{$\mathbb P$}     at 4.8 1
        \put{\ubvec222233} at 5.3 1
        \put{30}      at 0.4 0
        \put{\roote00111101} at 1.2 0
        \put{\roott01110011} at 2.1 0
        \put{\uwbvec211}   at 3 0
        \put{\isphi0{-1}0} at 4 0
        \put{$\mathbb P$}     at 4.8 0
        \put{\ubvec0{-1}1112} at 5.3 0
        \put{31}        at 6.4 29
        \put{\roote11111000}  at 7.2 29
        \put{\roott00111100}  at 8.1 29
        \put{\uwbvec041}    at 9 29
        \put{\isphi{-2}2{-1}} at 10 29
        \put{$\mathbb P$}     at 10.8 29
        \put{\ubvec{-2}02212} at 11.3 29
        \put{32}        at 6.4 28
        \put{\roote01111100}  at 7.2 28
        \put{\roott01111000}  at 8.1 28
        \put{\uwbvec041}    at 9 28
        \put{\isphi{-2}2{-1}} at 10 28
        \put{$\mathbb P$}     at 10.8 28
        \put{\ubvec{-2}02212} at 11.3 28
        \put{33}        at 6.4 27
        \put{\roote01111001}  at 7.2 27
        \put{\roott00111011}  at 8.1 27
        \put{\uwbvec211}    at 9 27
        \put{\isphi0{-1}0} at 10 27
        \put{$\mathbb P$}     at 10.8 27
        \put{\ubvec0{-1}1142} at 11.3 27
        \put{34}        at 6.4 26
        \put{\roote00012101}  at 7.2 26
        \put{\roott01100010}  at 8.1 26
        \put{\uwbvec111}    at 9 26
        \put{\isphi{-1}{-1}1} at 10 26
        \put{$\mathbb D_\ell$}     at 10.8 26
        \put{\ubvec{-1}{-2}1212} at 11.3 26
        \put{35}        at 6.4 25
        \put{\roote00111110}  at 7.2 25
        \put{\roott11110000}  at 8.1 25
        \put{\uwbvec131}    at 9 25
        \put{\isphi{-1}1{-1}} at 10 25
        \put{$\mathbb P$}     at 10.8 25
        \put{\ubvec{-1}01112} at 11.3 25
        \put{36}        at 6.4 24
        \put{\roote00011111}  at 7.2 24
        \put{\roott11100011}  at 8.1 24
        \put{\uwbvec301}    at 9 24
        \put{\isphi1{-2}{-\frac12}} at 10 24
        \put{$\mathbb D_s$}     at 10.8 24
        \put{\ubvec1{-1}2112} at 11.3 24
        \put{37}        at 6.4 23
        \put{\roote00111111}  at 7.2 23
        \put{\roott11110011}  at 8.1 23
        \put{\uwbvec311}    at 9 23
        \put{\isphi1{-1}{-1}} at 10 23
        \put{$\mathbb P$}     at 10.8 23
        \put{\ubvec101112} at 11.3 23
        \put{38}        at 6.4 22
        \put{\roote11111100}  at 7.2 22
        \put{\roott01111100}  at 8.1 22
        \put{\uwbvec051}    at 9 22
        \put{\isphi{-2}3{-\frac23}} at 10 22
        \put{$\mathbb H_\ell$}     at 10.8 22
        \put{\ubvec{-2}12315} at 11.3 22
        \put{39}        at 6.4 21
        \put{\roote11111001}  at 7.2 21
        \put{\roott00111111}  at 8.1 21
        \put{\uwbvec221}    at 9 21
        \put{$\frac00\quad$ND} at 10 21
        \put{c}     at 10.8 21
        \put{40}        at 6.4 20
        \put{\roote01111110}  at 7.2 20
        \put{\roott11111000}  at 8.1 20
        \put{\uwbvec141}    at 9 20
        \put{\isphi{-1}2{-\frac12}} at 10 20
        \put{$\mathbb D_\ell$}     at 10.8 20
        \put{\ubvec{-1}11212} at 11.3 20
        \put{41}        at 6.4 19
        \put{\roote01111101}  at 7.2 19
        \put{\roott01111011}  at 8.1 19
        \put{\uwbvec221}    at 9 19
        \put{$\frac00\quad$ND} at 10 19
        \put{c}     at 10.8 19
        \put{42}        at 6.4 18
        \put{\roote00112101}  at 7.2 18
        \put{\roott01110010}  at 8.1 18
        \put{\uwbvec121}    at 9 18
        \put{\isphi{-1}0\infty} at 10 18
        \put{$\mathbb P$}     at 10.8 18
        \put{\ubvec{-1}{-1}1112} at 11.3 18
        \put{43}        at 6.4 17
        \put{\roote00012111}  at 7.2 17
        \put{\roott11100010}  at 8.1 17
        \put{\uwbvec211}    at 9 17
        \put{\isphi0{-1}0} at 10 17
        \put{$\mathbb P$}     at 10.8 17
        \put{\ubvec0{-1}1112} at 11.3 17
        \put{44}        at 6.4 16
        \put{\roote01111111}  at 7.2 16
        \put{\roott11111011}  at 8.1 16
        \put{\uwbvec321}    at 9 16
        \put{\isphi10\infty} at 10 16
        \put{$\mathbb P$}     at 10.8 16
        \put{\ubvec111112} at 11.3 16
        \put{45}        at 6.4 15
        \put{\roote11111110}  at 7.2 15
        \put{\roott11111100}  at 8.1 15
        \put{\uwbvec151}    at 9 15
        \put{\isphi{-1}3{-\frac13}} at 10 15
        \put{$\mathbb H_\ell$}     at 10.8 15
        \put{\ubvec{-1}22315} at 11.3 15
        \put{46}        at 6.4 14
        \put{\roote00112111}  at 7.2 14
        \put{\roott11110010}  at 8.1 14
        \put{\uwbvec221}    at 9 14
        \put{$\frac00\quad$ND} at 10 14
        \put{c}     at 10.8 14
        \put{47}        at 6.4 13
        \put{\roote00122101}  at 7.2 13
        \put{\roott01210011}  at 8.1 13
        \put{\uwbvec222}    at 9 13
        \put{\isphi{-2}{-2}1} at 10 13
        \put{$\mathbb D_\ell$}     at 10.8 13
        \put{\ubvec{-2}{-4}2424} at 11.3 13
        \put{48}        at 6.4 12
        \put{\roote11111101}  at 7.2 12
        \put{\roott01111111}  at 8.1 12
        \put{\uwbvec231}    at 9 12
        \put{\isphi010} at 10 12
        \put{$\mathbb P$}     at 10.8 12
        \put{\ubvec011112} at 11.3 12
        \put{49}        at 6.4 11
        \put{\roote01112101}  at 7.2 11
        \put{\roott01111010}  at 8.1 11
        \put{\uwbvec131}    at 9 11
        \put{\isphi{-1}1{-1}} at 10 11
        \put{$\mathbb P$}     at 10.8 11
        \put{\ubvec{-1}01112} at 11.3 11
        \put{50}        at 6.4 10
        \put{\roote00012211}  at 7.2 10
        \put{\roott12100010}  at 8.1 10
        \put{\uwbvec221}    at 9 10
        \put{$\frac00\quad$ND} at 10 10
        \put{c}     at 10.8 10 
        \put{51}        at 6.4 9
        \put{\roote11112101}  at 7.2 9
        \put{\roott01111110}  at 8.1 9
        \put{\uwbvec141}    at 9 9
        \put{\isphi{-1}2{-\frac12}} at 10 9
        \put{$\mathbb D_\ell$}     at 10.8 9
        \put{\ubvec{-1}11212} at 11.3 9
        \put{52}        at 6.4 8
        \put{\roote01122101}  at 7.2 8
        \put{\roott01211011}  at 8.1 8
        \put{\uwbvec232}    at 9 8
        \put{\isphi{-2}{-1}2} at 10 8
        \put{$\mathbb H_\ell$}     at 10.8 8
        \put{\ubvec{-2}{-3}2324} at 11.3 8
        \put{53}        at 6.4 7
        \put{\roote11111111}  at 7.2 7
        \put{\roott11111111}  at 8.1 7
        \put{\uwbvec331}    at 9 7
        \put{\isphi111} at 10 7
        \put{$\mathbb D_s$}     at 10.8 7
        \put{\ubvec122112} at 11.3 7
        \put{54}        at 6.4 6
        \put{\roote01112111}  at 7.2 6
        \put{\roott11111010}  at 8.1 6
        \put{\uwbvec231}    at 9 6
        \put{\isphi010} at 10 6
        \put{$\mathbb P$}     at 10.8 6
        \put{\ubvec011112} at 11.3 6
        \put{55}        at 6.4 5
        \put{\roote00112211}  at 7.2 5
        \put{\roott12110010}  at 8.1 5
        \put{\uwbvec231}    at 9 5
        \put{\isphi010} at 10 5
        \put{$\mathbb P$}     at 10.8 5
        \put{\ubvec011142} at 11.3 5
        \put{56}        at 6.4 4
        \put{\roote00122111}  at 7.2 4
        \put{\roott11210011}  at 8.1 4
        \put{\uwbvec322}    at 9 4
        \put{\isphi{-1}{-2}{\frac12}} at 10 4
        \put{$\mathbb H_\ell$}     at 10.8 4
        \put{\ubvec{-1}{-3}2324} at 11.3 4
        \put{57}        at 6.4 3
        \put{\roote11122101}  at 7.2 3
        \put{\roott01211111}  at 8.1 3
        \put{\uwbvec242}    at 9 3
        \put{\isphi{-2}0\infty} at 10 3
        \put{$\mathbb P$}     at 10.8 3
        \put{\ubvec{-2}{-2}2224} at 11.3 3
        \put{58}        at 6.4 2
        \put{\roote01222101}  at 7.2 2
        \put{\roott01221011}  at 8.1 2
        \put{\uwbvec242}    at 9 2
        \put{\isphi{-2}0\infty} at 10 2
        \put{$\mathbb P$}     at 10.8 2
        \put{\ubvec{-2}{-2}2221} at 11.3 2
        \put{59}        at 6.4 1
        \put{\roote00122211}  at 7.2 1
        \put{\roott12210011}  at 8.1 1
        \put{\uwbvec332}    at 9 1
        \put{\isphi{-1}{-1}1} at 10 1
        \put{$\mathbb D_\ell$}     at 10.8 1
        \put{\ubvec{-1}{-2}1221} at 11.3 1
        \put{60} at 6.4 0
        \put{\roote01122111} at 7.2 0
        \put{\roott11211011} at 8.1 0
        \put{\uwbvec332} at 9 0
        \put{\isphi{-1}{-1}1} at 10 0
        \put{$\mathbb D_\ell$} at 10.8 0
        \put{\ubvec{-1}{-2}1221} at 11.3 0
\endpicture$$        

\newpage
\noindent Table of roots of $\mathbf E_8$ and $\chi_\Xi$ --- Roots 61 to 120

$$\beginpicture
        \setcoordinatesystem units <1.05cm,.58cm>
        \multiput{} at 0.4 0  2 30 /
        \put{61} at 0.4 29 
        \put{\roote11112111} at 1.2 29
        \put{\roott11111110} at 2.1 29
        \put{\uwbvec241} at 3 29
        \put{\isphi020} at 4 29
        \put{$\mathbb P$} at 4.8 29
        \put{\ubvec022212} at 5.3 29
        \put{62} at 0.4 28
        \put{\roote01112211} at 1.2 28
        \put{\roott12111010} at 2.1 28
        \put{\uwbvec241} at 3 28
        \put{\isphi020} at 4 28
        \put{$\mathbb P$} at 4.8 28
        \put{\ubvec022242} at 5.3 28
        \put{63}      at 0.4 27
        \put{\roote11112211} at 1.2 27
        \put{\roott12111110} at 2.1 27
        \put{\uwbvec251}   at 3 27
        \put{\isphi030} at 4 27
        \put{$\mathbb P$}     at 4.8 27
        \put{\ubvec033345} at 5.3 27
        \put{64}      at 0.4 26
        \put{\roote11222101} at 1.2 26
        \put{\roott01221111} at 2.1 26
        \put{\uwbvec252}   at 3 26
        \put{\isphi{-2}1{-2}} at 4 26
        \put{$\mathbb D_s$}     at 4.8 26
        \put{\ubvec{-2}{-1}2121} at 5.3 26
        \put{65}      at 0.4 25
        \put{\roote11122111} at 1.2 25
        \put{\roott11211111} at 2.1 25
        \put{\uwbvec342}   at 3 25
        \put{\isphi{-1}0\infty} at 4 25
        \put{$\mathbb P$}     at 4.8 25
        \put{\ubvec{-1}{-1}1121} at 5.3 25
        \put{66}      at 0.4 24
        \put{\roote01122211} at 1.2 24
        \put{\roott12211011} at 2.1 24
        \put{\uwbvec342}   at 3 24
        \put{\isphi{-1}0\infty} at 4 24
        \put{$\mathbb P$}     at 4.8 24
        \put{\ubvec{-1}{-1}1121} at 5.3 24
        \put{67}      at 0.4 23
        \put{\roote01222111} at 1.2 23
        \put{\roott11221011} at 2.1 23
        \put{\uwbvec342}   at 3 23
        \put{\isphi{-1}0\infty} at 4 23
        \put{$\mathbb P$}     at 4.8 23
        \put{\ubvec{-1}{-1}1121} at 5.3 23
        \put{68}      at 0.4 22
        \put{\roote00123211} at 1.2 22
        \put{\roott12210010} at 2.1 22
        \put{\uwbvec242}   at 3 22
        \put{\isphi{-2}0\infty} at 4 22
        \put{$\mathbb P$}     at 4.8 22
        \put{\ubvec{-2}{-2}2221} at 5.3 22
        \put{69}      at 0.4 21
        \put{\roote12222101} at 1.2 21
        \put{\roott01222111} at 2.1 21
        \put{\uwbvec262}   at 3 21
        \put{\isphi{-2}2{-1}} at 4 21
        \put{$\mathbb P$}     at 4.8 21
        \put{\ubvec{-2}02221} at 5.3 21
        \put{70}      at 0.4 20
        \put{\roote01123211} at 1.2 20
        \put{\roott12211010} at 2.1 20
        \put{\uwbvec252}   at 3 20
        \put{\isphi{-2}1{-2}} at 4 20
        \put{$\mathbb D_s$}     at 4.8 20
        \put{\ubvec{-2}{-1}2121} at 5.3 20
        \put{71}      at 0.4 19
        \put{\roote11222111} at 1.2 19
        \put{\roott11221111} at 2.1 19
        \put{\uwbvec352}   at 3 19
        \put{\isphi{-1}1{-1}} at 4 19
        \put{$\mathbb P$}     at 4.8 19
        \put{\ubvec{-1}01121} at 5.3 19
        \put{72}      at 0.4 18
        \put{\roote01222211} at 1.2 18
        \put{\roott12221011} at 2.1 18
        \put{\uwbvec352}   at 3 18
        \put{\isphi{-1}1{-1}} at 4 18
        \put{$\mathbb P$}     at 4.8 18
        \put{\ubvec{-1}01121} at 5.3 18
        \put{73}      at 0.4 17
        \put{\roote11122211} at 1.2 17
        \put{\roott12211111} at 2.1 17
        \put{\uwbvec352}   at 3 17
        \put{\isphi{-1}1{-1}} at 4 17
        \put{$\mathbb P$}     at 4.8 17
        \put{\ubvec{-1}01124} at 5.3 17
        \put{74}      at 0.4 16
        \put{\roote00123212} at 1.2 16
        \put{\roott12210021} at 2.1 16
        \put{\uwbvec422}   at 3 16
        \put{\isphi0{-2}0} at 4 16
        \put{$\mathbb P$}     at 4.8 16
        \put{\ubvec0{-2}2221} at 5.3 16
        \put{75}      at 0.4 15
        \put{\roote11222211} at 1.2 15
        \put{\roott12221111} at 2.1 15
        \put{\uwbvec362}   at 3 15
        \put{\isphi{-1}2{-\frac12}} at 4 15
        \put{$\mathbb D_\ell$}     at 4.8 15
        \put{\ubvec{-1}11224} at 5.3 15
        \put{76}      at 0.4 14
        \put{\roote11123211} at 1.2 14
        \put{\roott12211110} at 2.1 14
        \put{\uwbvec262}   at 3 14
        \put{\isphi{-2}2{-1}} at 4 14
        \put{$\mathbb P$}     at 4.8 14
        \put{\ubvec{-2}02221} at 5.3 14
        \put{77}      at 0.4 13
        \put{\roote12222111} at 1.2 13
        \put{\roott11222111} at 2.1 13
        \put{\uwbvec362}   at 3 13
        \put{\isphi{-1}2{-\frac12}} at 4 13
        \put{$\mathbb D_\ell$}     at 4.8 13
        \put{\ubvec{-1}11221} at 5.3 13
        \put{78}      at 0.4 12
        \put{\roote01223211} at 1.2 12
        \put{\roott12221010} at 2.1 12
        \put{\uwbvec262}   at 3 12
        \put{\isphi{-2}2{-1}} at 4 12
        \put{$\mathbb P$}     at 4.8 12
        \put{\ubvec{-2}02221} at 5.3 12
        \put{79}      at 0.4 11
        \put{\roote01123212} at 1.2 11
        \put{\roott12211021} at 2.1 11
        \put{\uwbvec432}   at 3 11
        \put{\isphi0{-1}0} at 4 11
        \put{$\mathbb P$}     at 4.8 11
        \put{\ubvec0{-1}1121} at 5.3 11
        \put{80}      at 0.4 10
        \put{\roote01233211} at 1.2 10
        \put{\roott12321011} at 2.1 10
        \put{\uwbvec363}   at 3 10
        \put{\isphi{-3}0\infty} at 4 10
        \put{$\mathbb P$}     at 4.8 10
        \put{\ubvec{-3}{-3}3333} at 5.3 10
        \put{81}      at 0.4 9
        \put{\roote12222211} at 1.2 9
        \put{\roott12222111} at 2.1 9
        \put{\uwbvec372}   at 3 9
        \put{\isphi{-1}3{-\frac13}} at 4 9
        \put{$\mathbb H_\ell$}     at 4.8 9
        \put{\ubvec{-1}22324} at 5.3 9
        \put{82}      at 0.4 8
        \put{\roote01223212} at 1.2 8
        \put{\roott12221021} at 2.1 8
        \put{\uwbvec442}   at 3 8
        \put{$\frac00\quad$ND} at 4 8
        \put{c}     at 4.8 8
        \put{83}      at 0.4 7
        \put{\roote11223211} at 1.2 7
        \put{\roott12221110} at 2.1 7
        \put{\uwbvec272}   at 3 7
        \put{\isphi{-2}3{-\frac23}} at 4 7
        \put{$\mathbb H_\ell$}     at 4.8 7
        \put{\ubvec{-2}12324} at 5.3 7
        \put{84}      at 0.4 6
        \put{\roote11123212} at 1.2 6
        \put{\roott12211121} at 2.1 6
        \put{\uwbvec442}   at 3 6
        \put{$\frac00\quad$ND} at 4 6
        \put{c}     at 4.8 6
        \put{85}      at 0.4 5
        \put{\roote01233212} at 1.2 5
        \put{\roott12321022} at 2.1 5
        \put{\uwbvec543}   at 3 5
        \put{\isphi{-1}{-2}{\frac12}} at 4 5
        \put{$\mathbb H_\ell$}     at 4.8 5
        \put{\ubvec{-1}{-3}2333} at 5.3 5
        \put{86}      at 0.4 4
        \put{\roote11233211} at 1.2 4
        \put{\roott12321111} at 2.1 4
        \put{\uwbvec373}   at 3 4
        \put{\isphi{-3}1{-3}} at 4 4
        \put{$\mathbb H_s$}     at 4.8 4
        \put{\ubvec{-3}{-2}3233} at 5.3 4
        \put{87}      at 0.4 3
        \put{\roote12223211} at 1.2 3
        \put{\roott12222110} at 2.1 3
        \put{\uwbvec282}   at 3 3
        \put{\isphi{-2}4{-\frac12}} at 4 3
        \put{$\mathbb D_\ell$}     at 4.8 3
        \put{\ubvec{-2}22424} at 5.3 3
        \put{88}      at 0.4 2
        \put{\roote11223212} at 1.2 2
        \put{\roott12221121} at 2.1 2
        \put{\uwbvec452}   at 3 2
        \put{\isphi010} at 4 2
        \put{$\mathbb P$}     at 4.8 2
        \put{\ubvec011121} at 5.3 2
        \put{89}      at 0.4 1
        \put{\roote12233211} at 1.2 1
        \put{\roott12322111} at 2.1 1
        \put{\uwbvec383}   at 3 1
        \put{\isphi{-3}2{-\frac32}} at 4 1
        \put{$\mathbb H_s$}     at 4.8 1
        \put{\ubvec{-3}{-1}3233} at 5.3 1
        \put{90}      at 0.4 0
        \put{\roote11233212} at 1.2 0
        \put{\roott12321122} at 2.1 0
        \put{\uwbvec553}   at 3 0
        \put{\isphi{-1}{-1}1} at 4 0
        \put{$\mathbb D_\ell$}     at 4.8 0
        \put{\ubvec{-1}{-2}1233} at 5.3 0
        \put{91}        at 6.4 29
        \put{\roote12223212}  at 7.2 29
        \put{\roott12222121}  at 8.1 29
        \put{\uwbvec462}    at 9 29
        \put{\isphi020} at 10 29
        \put{$\mathbb P$}     at 10.8 29
        \put{\ubvec022221} at 11.3 29
        \put{92}        at 6.4 28
        \put{\roote01234212}  at 7.2 28
        \put{\roott12321021}  at 8.1 28
        \put{\uwbvec453}    at 9 28
        \put{\isphi{-2}{-1}2} at 10 28
        \put{$\mathbb H_\ell$}     at 10.8 28
        \put{\ubvec{-2}{-3}2333} at 11.3 28
        \put{93}        at 6.4 27
        \put{\roote12333211}  at 7.2 27
        \put{\roott12332111}  at 8.1 27
        \put{\uwbvec393}    at 9 27
        \put{\isphi{-3}3{-1}} at 10 27
        \put{$\mathbb P$}     at 10.8 27
        \put{\ubvec{-3}03333} at 11.3 27
        \put{94}        at 6.4 26
        \put{\roote11234212}  at 7.2 26
        \put{\roott12321121}  at 8.1 26
        \put{\uwbvec463}    at 9 26
        \put{\isphi{-2}0\infty} at 10 26
        \put{$\mathbb P$}     at 10.8 26
        \put{\ubvec{-2}{-2}2233} at 11.3 26
        \put{95}        at 6.4 25
        \put{\roote12233212}  at 7.2 25
        \put{\roott12322122}  at 8.1 25
        \put{\uwbvec563}    at 9 25
        \put{\isphi{-1}0\infty} at 10 25
        \put{$\mathbb P$}     at 10.8 25
        \put{\ubvec{-1}{-1}1133} at 11.3 25
        \put{96}        at 6.4 24
        \put{\roote01234312}  at 7.2 24
        \put{\roott13321021}  at 8.1 24
        \put{\uwbvec463}    at 9 24
        \put{\isphi{-2}0\infty} at 10 24
        \put{$\mathbb P$}     at 10.8 24
        \put{\ubvec{-2}{-2}2233} at 11.3 24
        \put{97}        at 6.4 23
        \put{\roote12333212}  at 7.2 23
        \put{\roott12332122}  at 8.1 23
        \put{\uwbvec573}    at 9 23
        \put{\isphi{-1}1{-1}} at 10 23
        \put{$\mathbb P$}     at 10.8 23
        \put{\ubvec{-1}01133} at 11.3 23
        \put{98}        at 6.4 22
        \put{\roote12234212}  at 7.2 22
        \put{\roott12322121}  at 8.1 22
        \put{\uwbvec473}    at 9 22
        \put{\isphi{-2}1{-2}} at 10 22
        \put{$\mathbb D_s$}     at 10.8 22
        \put{\ubvec{-2}{-1}2133} at 11.3 22
        \put{99}        at 6.4 21
        \put{\roote11234312}  at 7.2 21
        \put{\roott13321121}  at 8.1 21
        \put{\uwbvec473}    at 9 21
        \put{\isphi{-2}1{-2}} at 10 21
        \put{$\mathbb D_s$}     at 10.8 21
        \put{\ubvec{-2}{-1}2133} at 11.3 21
        \put{100}        at 6.4 20
        \put{\roote01234322}  at 7.2 20
        \put{\roott23321021}  at 8.1 20
        \put{\uwbvec563}    at 9 20
        \put{\isphi{-1}0\infty} at 10 20
        \put{$\mathbb P$}     at 10.8 20
        \put{\ubvec{-1}{-1}1133} at 11.3 20
        \put{101}        at 6.4 19
        \put{\roote12234312}  at 7.2 19
        \put{\roott13322121}  at 8.1 19
        \put{\uwbvec483}    at 9 19
        \put{\isphi{-2}2{-1}} at 10 19
        \put{$\mathbb P$}     at 10.8 19
        \put{\ubvec{-2}02233} at 11.3 19
        \put{102}        at 6.4 18
        \put{\roote12334212}  at 7.2 18
        \put{\roott12332121}  at 8.1 18
        \put{\uwbvec483}    at 9 18
        \put{\isphi{-2}2{-1}} at 10 18
        \put{$\mathbb P$}     at 10.8 18
        \put{\ubvec{-2}02233} at 11.3 18
        \put{103}        at 6.4 17
        \put{\roote11234322}  at 7.2 17
        \put{\roott23321121}  at 8.1 17
        \put{\uwbvec573}    at 9 17
        \put{\isphi{-1}1{-1}} at 10 17
        \put{$\mathbb P$}     at 10.8 17
        \put{\ubvec{-1}01133} at 11.3 17
        \put{104}        at 6.4 16
        \put{\roote12344212}  at 7.2 16
        \put{\roott12432122}  at 8.1 16
        \put{\uwbvec584}    at 9 16
        \put{\isphi{-3}0\infty} at 10 16
        \put{$\mathbb P$}     at 10.8 16
        \put{\ubvec{-3}{-3}3345} at 11.3 16
        \put{105}        at 6.4 15
        \put{\roote12334312}  at 7.2 15
        \put{\roott13332121}  at 8.1 15
        \put{\uwbvec493}    at 9 15
        \put{\isphi{-2}3{-\frac23}} at 10 15
        \put{$\mathbb H_\ell$}     at 10.8 15
        \put{\ubvec{-2}12333} at 11.3 15
        \put{106}        at 6.4 14
        \put{\roote12234322}  at 7.2 14
        \put{\roott23322121}  at 8.1 14
        \put{\uwbvec583}    at 9 14
        \put{\isphi{-1}2{-\frac12}} at 10 14
        \put{$\mathbb D_\ell$}     at 10.8 14
        \put{\ubvec{-1}11233} at 11.3 14
        \put{107}        at 6.4 13
        \put{\roote12334322}  at 7.2 13
        \put{\roott23332121}  at 8.1 13
        \put{\uwbvec593}    at 9 13
        \put{\isphi{-1}3{-\frac13}} at 10 13
        \put{$\mathbb H_\ell$}     at 10.8 13
        \put{\ubvec{-1}22333} at 11.3 13
        \put{108}        at 6.4 12
        \put{\roote12344312}  at 7.2 12
        \put{\roott13432122}  at 8.1 12
        \put{\uwbvec594}    at 9 12
        \put{\isphi{-3}1{-3}} at 10 12
        \put{$\mathbb H_s$}     at 10.8 12
        \put{\ubvec{-3}{-2}3242} at 11.3 12
        \put{109}        at 6.4 11
        \put{\roote12344322}  at 7.2 11
        \put{\roott23432122}  at 8.1 11
        \put{\uwbvec694}    at 9 11
        \put{\isphi{-2}1{-2}} at 10 11
        \put{$\mathbb D_s$}     at 10.8 11
        \put{\ubvec{-2}{-1}2142} at 11.3 11
        \put{110}        at 6.4 10
        \put{\roote12345312}  at 7.2 10
        \put{\roott13432121}  at 8.1 10
        \put{\uwbvec4{10}4}    at 9 10
        \put{\isphi{-4}2{-2}} at 10 10
        \put{$\mathbb D_s$}     at 10.8 10
        \put{\ubvec{-4}{-2}4242} at 11.3 10
        \put{111}        at 6.4 9
        \put{\roote12345313}  at 7.2 9
        \put{\roott13432132}  at 8.1 9
        \put{\uwbvec684}    at 9 9
        \put{\isphi{-2}0{\infty}} at 10 9
        \put{$\mathbb P$}     at 10.8 9
        \put{\ubvec{-2}{-2}2242} at 11.3 9
        \put{112}        at 6.4 8
        \put{\roote12345322}  at 7.2 8
        \put{\roott23432121}  at 8.1 8
        \put{\uwbvec5{10}4}    at 9 8
        \put{\isphi{-3}2{-\frac32}} at 10 8
        \put{$\mathbb H_s$}     at 10.8 8
        \put{\ubvec{-3}{-1}3242} at 11.3 8
        \put{113}        at 6.4 7
        \put{\roote12345422}  at 7.2 7
        \put{\roott24432121}  at 8.1 7
        \put{\uwbvec5{11}4}    at 9 7
        \put{\isphi{-3}3{-1}} at 10 7
        \put{$\mathbb P$}     at 10.8 7
        \put{\ubvec{-3}03345} at 11.3 7
        \put{114}        at 6.4 6
        \put{\roote12345323}  at 7.2 6
        \put{\roott23432132}  at 8.1 6
        \put{\uwbvec784}    at 9 6
        \put{\isphi{-1}0\infty} at 10 6
        \put{$\mathbb P$}     at 10.8 6
        \put{\ubvec{-1}{-1}1142} at 11.3 6
        \put{115}        at 6.4 5
        \put{\roote12345423}  at 7.2 5
        \put{\roott24432132}  at 8.1 5
        \put{\uwbvec794}    at 9 5
        \put{\isphi{-1}1{-1}} at 10 5
        \put{$\mathbb P$}     at 10.8 5
        \put{\ubvec{-1}01142} at 11.3 5
        \put{116}        at 6.4 4
        \put{\roote12346423}  at 7.2 4
        \put{\roott24432131}  at 8.1 4
        \put{\uwbvec6{10}4} at 9 4
        \put{\isphi{-2}2{-1}} at 10 4
        \put{$\mathbb P$}     at 10.8 4
        \put{\ubvec{-2}02242} at 11.3 4
        \put{117}        at 6.4 3
        \put{\roote12356423}  at 7.2 3
        \put{\roott24532132}  at 8.1 3
        \put{\uwbvec7{10}5} at 9 3
        \put{\isphi{-3}0\infty} at 10 3
        \put{$\mathbb P$}     at 10.8 3
        \put{\ubvec{-3}{-3}3354} at 11.3 3
        \put{118}        at 6.4 2
        \put{\roote12456423}  at 7.2 2
        \put{\roott24542132}  at 8.1 2
        \put{\uwbvec7{11}5} at 9 2
        \put{\isphi{-3}1{-3}} at 10 2
        \put{$\mathbb H_s$}     at 10.8 2
        \put{\ubvec{-3}{-2}3251} at 11.3 2
        \put{119}        at 6.4 1
        \put{\roote13456423}  at 7.2 1
        \put{\roott24543132}  at 8.1 1
        \put{\uwbvec7{12}5}  at 9 1
        \put{\isphi{-3}2{-\frac32}} at 10 1
        \put{$\mathbb H_s$}     at 10.8 1
        \put{\ubvec{-3}{-1}3251} at 11.3 1
        \put{120} at 6.4 0
        \put{\roote23456423} at 7.2 0
        \put{\roott24543232} at 8.1 0
        \put{\uwbvec7{13}5} at 9 0
        \put{\isphi{-3}3{-1}} at 10 0
        \put{$\mathbb P$} at 10.8 0
        \put{\ubvec{-3}03354} at 11.3 0
\endpicture$$        

\vfill\eject

\subsection{Some special roots.}
\label{app-A-two}

{\bf (a) The object $E_2^2$ as an example.}
Consider the row corresponding to root $\mathbf r_{47}$ in the table:

\smallskip
\centerline{47 $\quad$ \roote00122101 $\quad$ \roott01210011 $\quad$ {\uwbvec222} $\quad$
\isphi{-2}{-2}1 $\quad \mathbb D_\ell \quad$ \ubvec{-2}{-4}2424}

\smallskip
The entry after the root number is the corresponding root of $\mathbf E_8$,
then the dimension vector $\mathbf r_{47}$ of
the object $M\in {}'\Cal D$ with $\pi M=E_2^2$,
and the uwb-vector \uwbvec222.  In $\mathbb T(6)$, the object $E_2^2$
can be found on the diagonal line with slope \isphi{-2}{-2}1,
on the long half line of a central line of type $\mathbb D$.
The central distance is given by the in-radius $\frac{r_\Delta}{bE_2^2}=\frac22=1$
of the standard triangle
(so the boundary distance
is $dE_2^2=2-\frac{r_\Delta}{bE_2^2}=1$) and the in-radius
$\frac{r_\nabla}{bE_2^2}=\frac42=2$
of the costandard triangle (so $E_2^2$ occurs on $\nabla_{4}$).
Finally, the corresponding primitive pair is $(2bE_2^2-r_\Delta,bE_2^2)=(2,2)$.

\smallskip The root given by dimension vector $\mathbf h_0+\mathbf r_{47}$ corresponds to
an object $X$ of width $3+2=5$ (a pentapicket), it occurs on the same half-line in $\mathbb T(6)$,
but the central distances are given by the in-radii $\frac{r_\Delta}{bX}=\frac2{5}$ and
$\frac{r_\nabla}{bX}=\frac4{5}$ of the standard
and costandard triangles, respectively.
Hence the boundary distance is $2-\frac25=\frac85$, and the object has primitive pair
$(8,5)$.

\smallskip The smallest realization $M$ of a corresponding negative root
$-\mathbf r_{47}$ within $'\Cal D$
has dimension vector 
$\mathbf h_0+\mathbf h_\infty-\mathbf r_{47}$; we abbreviate $\mathbf h_0+\mathbf h_\infty$ as 
$\mathbf h_1$.
(Note that the vectors $\mathbf h_0-\mathbf r_{47}$ and $\mathbf h_\infty-\mathbf r_{47}$ do not have their
largest entry in position 3, hence cannot be realized by objects in $'\Cal D$,
see Lemma~\ref{lem-eleven-two}.
The object $X=\pi M\in\Cal S(6)$ occurs on the short half line of the same central line;
since $bX=\dim M_3=6-2=4$, $X$ is a tetrapicket,
its central distance is given by the in-radius $\frac{r_\nabla}{bX}=\frac44$
of the standard triangle and the in-radius $\frac{r_\Delta}{bX}=\frac24$ of the costandard triangle.
As the root is negative, the boundary distance is obtained as $dX=2-\frac{r_\nabla}{bX}=2-\frac44=1$
and the primitive pair as $(4,4)$. 

\medskip
{\bf (b) The objects on the costandard triangle $\nabla_{5}$.}
The vertices next to the corners of the triangle $\mathbb T(6)$ support objects of quasi-length $1$ in the
principal component.  All six objects can be realized by roots of $\chi_\Theta$,
they have been pictured in the introduction in Section~\ref{sec-fifteen-two}
and are listed in the order in which they
occur on their $\tau$-orbit.
$$\begin{matrix} S=(0,[1]) & ([1],[1]) & ([1],[6]) & (0,[5]) & ([5],[5]) & ([5],[6]) \cr
\mathbf r_{12} & \mathbf r_4 &  \mathbf r_{45} & \mathbf r_{38} &
      \mathbf h_1-\mathbf r_{119} & \mathbf h_1-\mathbf r_{118}
\end{matrix}$$

\smallskip
By adding to each root either $\mathbf h_0$ or $\mathbf h_\infty$ we obtain the dimension vectors of the 12
tetrapickets pictured in Appendix~\ref{app-B-five}.

\medskip
{\bf (c) The bipickets in the corners of $\Delta_1$.}
The three bipickets occur in the 6-tube of rationality index $\frac 11$, each is given by a root.
We have already considered $E_2^2$ in (a).
$$\begin{matrix} E_2^2=([2],[3,1]) & ([2],[6,4]) & ([5,3],[6,4],[2]) \cr
          \mathbf r_{47} & \mathbf r_{87} & \mathbf h_1-\mathbf r_{110} \cr
          \end{matrix}$$
For each of those three roots $\mathbf r$, the root $\mathbf h_1-\mathbf r$ defines the tetrapicket
on $\Delta_1$, on the lines $q=5$, $r=1$ and $p=1$, respectively.  They are pictured in
Section \ref{sec-fifteen-two}--(d).

\smallskip
Moreover, the roots of the form $\mathbf h_0+\mathbf r$ and $\mathbf h_\infty+\mathbf r$
define the six pentapickets on the vertices of $\Delta_{8/5}$, see Appendix~\ref{app-B-six}.

\medskip 
{\bf (d) The six bipickets on the costandard triangle $\nabla_{7/2}$.}
The objects, which lie on the lines $q=\frac52, r=\frac72$ and $p=\frac72$, respectively, are
listed in the counterclockwise order in which they occur around the center of $\mathbb T(6)$. 
$$\smallmatrix E_2^3=([2],[4,1]) & ([2],[6,3]) & ([3],[6,4]) & ([5,2],[6,4],[3]) & ([5,2],[6,3],[2]) & ([3],[4,1]) \cr
\mathbf r_{52} & \mathbf r_{83} & \mathbf r_{81} & \mathbf h_1-\mathbf r_{108} & \mathbf h_1-\mathbf r_{112} & \mathbf r_{56}\cr
      \endsmallmatrix$$
      
The objects are pictured at the beginning of Section~\ref{sec-fifteen-two}.

\medskip
{\bf (e) The 18 tripickets on the hexagon.}
We specify the roots of the 18 tripickets on the boundary of the hexagon in $\Cal S(6)$,
as presented in Section \ref{sec-fifteen-two}--(c).
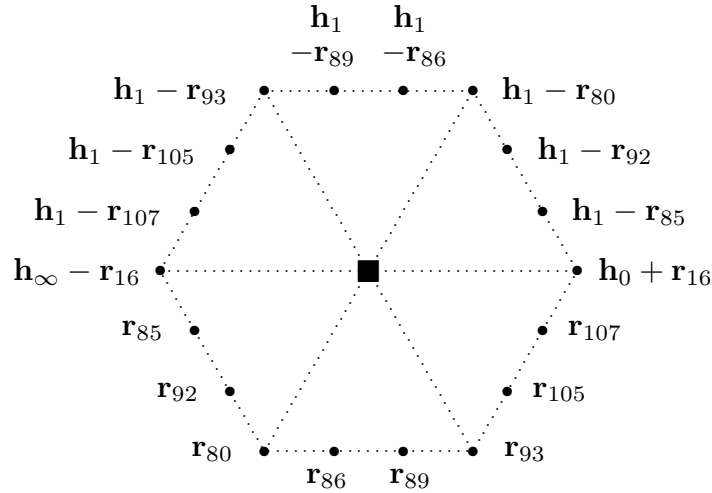
$$
{\beginpicture
   \setcoordinatesystem units <1.155cm,2cm>
   \multiput{} at -2 0  2 2 /
   \setdots<1mm>
   \plot -2 1  -1 2  1 2  2 1  1 0  -1 0  -2 1 /
   \put{$\blacksquare$} at 0 1
   \plot -2 1  2 1 /
   \plot -1 0  1 2 /
   \plot -1 2  1 0 /
   \multiput{$\ssize\bullet$} at -1 0  -.33 0  .33 0  1 0  1.33 .33  1.67 .67  2 1  1.67 1.33  1.33 1.67  1 2
                        .33 2  -.33 2  -1 2  -1.33 1.67  -1.67 1.33  -2 1  -1.67 .67  -1.33 .33 /
   \put{$\mathbf r_{80}$} at -1.5 0
   \put{$\mathbf r_{86}$} at -.4 -.15
   \put{$\mathbf r_{89}$} at .4 -.15
   \put{$\mathbf r_{93}$} at 1.5 0
   \put{$\mathbf r_{105}$} at 1.83 .33
   \put{$\mathbf r_{107}$} at 2.17 .67
   \put{$\mathbf h_0+\mathbf r_{16}$} at 2.75 1
   \put{$\mathbf h_1-\mathbf r_{85}$} at 2.5 1.33
   \put{$\mathbf h_1-\mathbf r_{92}$} at 2.17 1.67
   \put{$\mathbf h_1-\mathbf r_{80}$} at 1.83 2
   \put{$\mathbf h_1$} at .4 2.4
   \put{$\;-\mathbf r_{86}$} at .4 2.2 
   \put{$\mathbf h_1$} at -.4 2.4
   \put{$-\mathbf r_{89}\;$} at -.4 2.2
   \put{$\mathbf h_1-\mathbf r_{93}$} at -1.9 2
   \put{$\mathbf h_1-\mathbf r_{105}$} at -2.27 1.67
   \put{$\mathbf h_1-\mathbf r_{107}$} at -2.6 1.33
   \put{$\mathbf h_\infty-\mathbf r_{16}$} at -2.8 1
   \put{$\mathbf r_{85}$} at -2.17 .67
   \put{$\mathbf r_{92}$} at -1.83 .33
   \endpicture}
$$

\medskip
{\bf (f) Some objects which are not roots.}  In the principal component $\Cal P(6)$, some objects cannot be
realized by roots or radical vectors of $\chi_\Theta$.
Those include the projective module $P'=([6],[6])$: If $P'=\pi M$ for some
object $M\in\widetilde{\Cal S}(6)$ then $\dim M_{i'}=1$ for six consecutive integers $i$, which is not
possible for any root of $\chi_\Theta$.
Similarly, the projective module $P=(0,[6])$: Assume $P=\pi M$ with $\bdim M$ a root of $\chi_\Theta$, then either
$M_1=0$ or else $M_{1'}\neq0$ and neither is possible.
Also, $P[4]=([5,2],[6,6,2],[5,2])$ cannot be realized since any $M\in\widetilde{\Cal S}(6)$ with $P[4]=\pi M$
has support on an interval of length $8$.
 	  
\bigskip
\renewcommand{\refname}{Additional reference for Appendix A}

\vfill\eject


\setcounter{sectionapp}{2}
\setcounter{section}{20}

\section{The tripickets in $\mathcal S(6)$,
  some tetra- and some pentapickets.}
\label{app-B}

by Markus Schmidmeier

\subsection{Overview.}
\label{app-B-one}

All tripickets occur inside or on the boundary of the hexagon of $\Cal S(6)$.
In Section \ref{sec-fifteen-two}--(c), we have already determined the 18 tripickets on the boundary.

 The remaining tripickets reside inside the hexagon.  Several
tripickets will share their position.  For the 57 non-central objects in the interior of the hexagon,
we indicate this multiplicity.

Finally, there are the central objects.  There are 3 objects in $\Cal P(6)$,
5 in the stable non-homogeneous tubes of index 0,
and one on the mouth of each homogeneous tube of index 0.

\medskip
$$
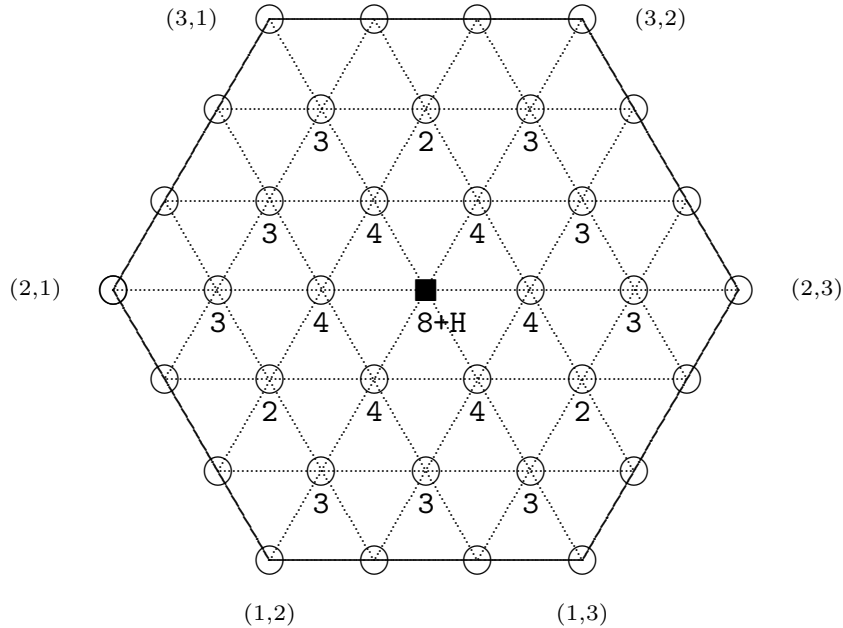
{\beginpicture
   \setcoordinatesystem units <1.732cm, 3cm>
   \multiput{} at -2 0  2 2 /
   \plot -2 1  -1 2  1 2  2 1  1 0  -1 0  -2 1 /
   \multiput{$\bigcirc$} at -2 1  -1.67 1.33  -1.33 1.67  -1 2 -.33 2  .33 2  1 2
                         1.33 1.67  1.67 1.33   2 1  1.67 .67  1.33 .33 
                         1 0  .33 0  -.33 0  -1 0  -1.33 .33  -1.67 .67  -2 1 /
   \multiput{$\bigcirc$} at -.67 .33  0 .33  .67 .33  -1 .67  -.33 .67  .33 .67  1 .67
                         -1.33 1  -.67 1  .67 1  1.33 1
                         -1 1.33  -.33 1.33  .33 1.33  1 1.33  -.67 1.67  0 1.67  .67 1.67 /
   \multiput{{\tt 2}} at -1 .55  0 1.55  1 .55 /
   \multiput{{\tt 3}} at -.67 .21   0 .21  .67 .21  -1.33 .88  1.33 .88  -1 1.21  1 1.21  -.67 1.55  .67 1.55 /
   \multiput{{\tt 4}} at -.33 .55  .33 .55  -.67 .88  .67 .88  -.33 1.21  .33 1.21 /
   \put{{\tt 8+H}} at .1 .88
   \put{$\blacksquare$} at 0 1
   \setdots<.5mm>
   \plot -2 1  2 1 /
   \plot -1 0  1 2 /
   \plot -1 2  1 0 /
   \plot -.33 0  -1.67 1.33  1.67 1.33   .33 0  -1.33 1.67  1.33 1.67  -.33 0 /
   \plot -1.33 .33  .33 2  1.67 .67  -1.67 .67  -.33 2  1.33 .33  -1.33 .33 / 
   \put{$\ss (1,2)$} at -1 -.2
   \put{$\ss (1,3)$} at 1 -.2
   \put{$\ss (2,1)$} at -2.5 1
   \put{$\ss (3,1)$} at -1.5 2
   \put{$\ss (2,3)$} at 2.5 1
   \put{$\ss (3,2)$} at 1.5 2 
\endpicture}
$$

\medskip
In Appendix~\ref{app-B-two}, we picture the non-central objects on the diagonal lines (15).
In Appendix~\ref{app-B-three},
we deal with the non-central tripickets on $\Delta_{4/3}$ which lie on a central line parallel
to one of the axes (18).
In Appendix~\ref{app-B-four},
we present the non-central tripickets in $\Delta_{5/3}$ which lie on a central line parallel
to one of the axes (24).

\medskip
In this appendix, we also exhibit some additional tetrapickets,
in addition to those on the boundary of the hexagon in Section \ref{sec-fifteen-two}--(d)
and on the triangle $\Delta_{5/4}$ in Section \ref{sec-fifteen-two}--(f). 
The tetrapickets on the triangle $\Delta_{5/3}$
and on the lines $q=\frac{13}4$, $p=\frac{11}4$ and $r=\frac{11}4$
are presented in Appendix~\ref{app-B-five}.
Finally, in Appendix~\ref{app-B-six} we picture some pentapickets inside but
close to the boundary of the hexagon.

\medskip
Here, we list the non-homogeneous central tripickets from \cite[(2.3), (8)]{RS1},
the homogeneous ones are pictured in Remark~\ref{rem-two-seven-one} in Section~\ref{sec-two-seven}.
$$\beginpicture
        \setcoordinatesystem units <1.2cm, 1cm>
        \multiput{} at 0 0  9 2.5 /
        \put{\strut tripickets} at 1 .8
        \put{\strut in $\Cal P(6)$} at 1 .4
        \put{\strut (quasi-length 6)} at 1 0
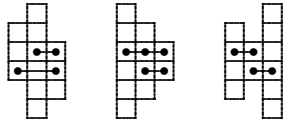
        \put{\strut tripickets in 3-tube} at 5 .8
        \put{\strut of index 0} at 5 .4
        \put{\strut(quasi-length 3)} at 5 0
        \put{\strut tripickets in 2-tube} at 8.5 .8
        \put{\strut of index 0} at 8.5 .4
        \put{\strut(quasi-length 2)} at 8.5 0
   \put{\beginpicture 
     \setcoordinatesystem units <.21cm,.21cm>
     \multiput{} at 0 0  3 8  /
     \setsolid
     \plot 0 2  0 7  1 7  1 2  0 2 /
     \plot 0 6  1 6 /
     \plot 0 5  3 5 /
     \plot 0 4  3 4 /
     \plot 0 3  3 3 /
     \plot 2 1  2 6  3 6  3 1  2 1 /
     \plot 2 2  3 2 /
     \multiput{$\sss\bullet$} at .5 3.5  1.5 3.5  1.5 4.5  2.5 4.5 /
     \plot .5 3.5  1.5 3.5 /
     \plot 1.5 4.5  2.5 4.5 /
     \endpicture} at 0 2
   \put{\beginpicture 
     \setcoordinatesystem units <.21cm,.21cm>
     \multiput{} at 0 0  3 6  /
     \setsolid
     \plot 1 0  1 6  2 6  2 0  1 0 /
     \plot 1 1  3 1  3 5  1 5  /
     \plot 3 2  0 2  0 4  3 4 /
     \plot 0 3  3 3 /
     \multiput{$\sss\bullet$} at .5 3.5  1.5 3.5  2.5 3.5  .5 2.5  2.5 1.5 /
     \plot .5 3.5  2.5 3.5 /
     \endpicture} at 1 2
   \put{\beginpicture
     \setcoordinatesystem units <.21cm,.21cm>
     \multiput{} at 0 0  3 6   /
     \setsolid
     \plot 0 0  0 6  1 6  1 0  0 0 /
     \plot 0 1  2 1  2 5  0 5 /
     \plot 0 2  3 2  3 4  0 4 /
     \plot 0 3  3 3 /
     \multiput{$\sss\bullet$} at .5 3.5  1.5 3.5  1.5 2.5  2.5 2.5 /
     \plot 1.5 2.5  2.5 2.5 /
     \plot .5 3.5  1.5 3.5 /
     \endpicture} at 2 2
   \put{\beginpicture 
     \setcoordinatesystem units <.21cm,.21cm>
     \multiput{} at 0 0  3 6  /
     \setsolid
     \plot 1 0  1 6  2 6  2 0  1 0 /
     \plot 2 5  0 5  0 2  3 2 /
     \plot 0 4  3 4  3 1  1 1 /
     \plot 0 3  3 3 /
     \multiput{$\sss\bullet$} at .5 2.5  2.5 2.5  1.5 3.5  2.5 3.5 /
     \plot .5 2.5  2.5 2.5 /
     \plot 1.5 3.5  2.5 3.5 /
     \endpicture} at 4 2
   \put{\beginpicture
     \setcoordinatesystem units <.21cm,.21cm>
     \multiput{} at 0 0  3 6  /
     \setsolid
     \plot 0 0  0 6  1 6  1 0  0 0 /
     \plot 0 1  2 1  2 5  0 5 /
     \plot 0 2  3 2  3 4  0 4 /
     \plot 0 3  3 3 /
     \multiput{$\sss\bullet$} at .5 3.5  1.5 3.5  2.5 3.5  1.5 2.5  2.5 2.5 /
     \plot 1.5 2.5  2.5 2.5 /
     \plot .5 3.5  2.5 3.5 /
     \endpicture} at 5 2
   \put{\beginpicture
     \setcoordinatesystem units <.21cm,.21cm>
     \multiput{} at 0 0  3 6  /
     \setsolid
     \plot 2 0  2 6  3 6  3 0   2 0 /
     \plot 0 1  1 1  1 5  0 5  0 1 /
     \plot 2 1  3 1 /
     \plot 0 2  3 2 /
     \plot 0 3  3 3 /
     \plot 0 4  3 4 /
     \plot 2 5  3 5 /
     \multiput{$\sss\bullet$} at .5 3.5  1.5 3.5  1.5 2.5  2.5 2.5 /
     \plot 1.5 2.5  2.5 2.5 /
     \plot .5 3.5  1.5 3.5 /
     \endpicture} at 6 2
   \put{\beginpicture 
     \setcoordinatesystem units <.21cm,.21cm>
     \multiput{} at 0 0  3 6  /
     \setsolid
     \plot 0 0  0 6  1 6  1 0  0 0 /
     \plot 2 1  3 1  3 5  2 5  2 1 /
     \plot 0 1  1 1 /
     \plot 0 2  3 2 /
     \plot 0 3  3 3 /
     \plot 0 4  3 4 /
     \plot 0 5  1 5 /
     \multiput{$\sss\bullet$} at .5 3.5  1.5 3.5  1.5 2.5  2.5 2.5 /
     \plot 1.5 2.5  2.5 2.5 /
     \plot .5 3.5  1.5 3.5 /
     \endpicture} at 8 2
   \put{\beginpicture
     \setcoordinatesystem units <.21cm,.21cm>
     \multiput{} at 0 0  3 6  /
     \setsolid
     \plot 1 0  1 6  2 6  2 0  1 0 /
     \plot 2 1  0 1  0 5  2 5  /
     \plot 0 2  3 2  3 4  0 4  /
     \plot 0 3  3 3 /
     \multiput{$\sss\bullet$} at .5 2.5  .5 3.5  1.5 3.5  2.5 3.5 /
     \plot .5 3.5  2.5 3.5 /
     \endpicture} at 9 2
\endpicture
$$

\subsection{The non-central tripickets on the diagonal lines.}
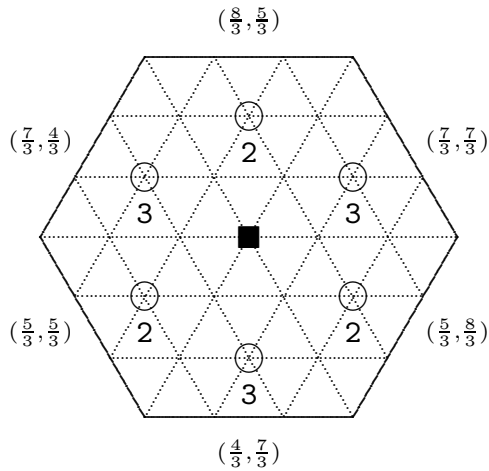
\label{app-B-two}

$$
{\beginpicture
   \setcoordinatesystem units <1.1547cm, 2cm>
   \multiput{} at -2 0  2 2 /
   \plot -2 1  -1 2  1 2  2 1  1 0  -1 0  -2 1 /
   \multiput{$\bigcirc$} at -1 .67  0 1.67  1 .67
                           0 .33  -1 1.33  1 1.33 /
   \multiput{{\tt 2}} at -1 .47  0 1.47  1 .47 /
   \multiput{{\tt 3}} at  0 .13  -1 1.13  1 1.13 /
     \put{$\blacksquare$} at 0 1
   \setdots<.5mm>
   \plot -2 1  2 1 /
   \plot -1 0  1 2 /
   \plot -1 2  1 0 /
   \plot -.33 0  -1.67 1.33  1.67 1.33   .33 0  -1.33 1.67  1.33 1.67  -.33 0 /
   \plot -1.33 .33  .33 2  1.67 .67  -1.67 .67  -.33 2  1.33 .33  -1.33 .33 / 
   \put{$\ss (\frac43,\frac73)$} at 0 -.2
   \put{$\ss (\frac53,\frac83)$} at 2 .47
   \put{$\ss (\frac53,\frac53)$} at -2 .47
   \put{$\ss (\frac73,\frac43)$} at -2 1.53
   \put{$\ss (\frac73,\frac73)$} at 2 1.53
   \put{$\ss (\frac83,\frac53)$} at 0 2.2
\endpicture}
$$

\medskip
The objects pictured on the left (the right) have quasi-length 2 in the 6-tube of index $\frac12$
(of index $\frac 21$); the middle object in the bottom row (and its rotations)
occurs in the principal component $\Cal P(6)$ with quasi-length 4.

\medskip
Note that there are 3 orbits of objects under the action of
the symmetry group $\Sigma_3$:  
One orbit consists of the six objects on the long central half lines
$\mathbb D_\ell$, the second has six objects (the orbits under rotation of the outer
objects pictured at the bottom) on the short central half lines $\mathbb D_s$,
while the third orbit contains only three objects on half lines of type $\mathbb D_s$.

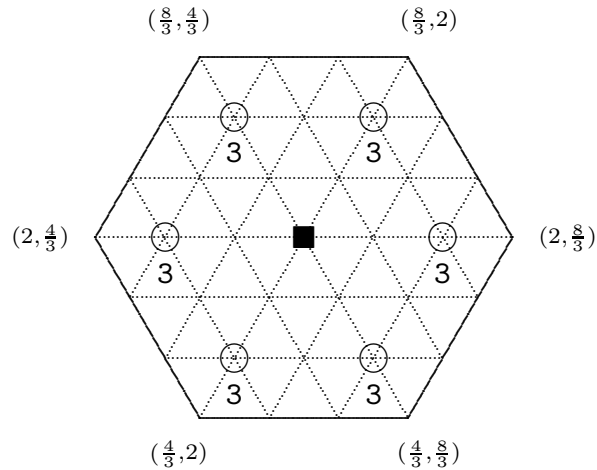
$$
{\beginpicture
   \setcoordinatesystem units <1.4722cm, 2.55cm>
   \multiput{} at -3 -.5  3 3 /
   \setdots<1mm>
   \plot -3 1  -1.5 2.5  1.5 2.5  3 1  1.5 -.5  -1.5 -.5  -3 1 /
   \put{$\blacksquare$} at 0 1
   \plot -3 1  3 1 /
   \plot -1.5 -.5  1.5 2.5 /
   \plot -1.5 2.5  1.5 -.5 /
   \plot -.5  -.5  -2.5 1.5  2.5 1.5  .5 -.5  -2 2  2 2  -.5 -.5 /
   \plot -2 0  2 0  -.5 2.5  -2.5 .5  2.5 .5  .5 2.5  -2 0 /
   \put{\beginpicture
     \setcoordinatesystem units <.21cm,.21cm>
     \multiput{} at 0 0  3 6  /
     \setsolid
     \plot 0 0  0 6  1 6  1 0  0 0 /
     \plot 0 2  2 2  2 5  0 5 /
     \plot 0 3  3 3  3 4  0 4 /
     \plot 0 1  1 1 /
     \multiput{$\sss\bullet$} at .5 3.5  1.5 3.5  2.5 3.5  1.5 2.5 /
     \plot .5 3.5  2.5 3.5 /
     \endpicture} at -1.75 .5
   \put{\beginpicture
     \setcoordinatesystem units <.21cm,.21cm>
     \multiput{} at 0 0  3 6  /
     \setsolid
     \plot 0 0  0 6  1 6  1 0  0 0 /
     \plot 0 1  2 1  2 4  0 4 /
     \plot 0 2  3 2  3 3  0 3 /
     \plot 0 5  1 5 /
     \multiput{$\sss\bullet$} at .4 2.4  1.4 2.7  1.7 2.3  2.7 2.7 /
     \plot .4 2.4  1.4 2.7 /
     \plot 1.7 2.3  2.7 2.7 /
     \endpicture} at -1.25 .5
   \put{\beginpicture
     \setcoordinatesystem units <.21cm,.21cm>
     \multiput{} at 0 0  3 6  /
     \setsolid
     \plot 0 0  0 6  1 6  1 0  0 0 /
     \plot 2 0  2 5  3 5  3 0  2 0 /
     \plot 0 1  3 1 /
     \plot 0 2  3 2 /
     \plot 0 3  3 3 /
     \plot 0 4  1 4 /
     \plot 2 4  3 4 /
     \plot 0 5  1 5 /
     \multiput{$\sss\bullet$} at .5 2.5  1.5 2.5  1.5 1.5  2.5 1.5 /
     \plot .5 2.5  1.5 2.5 /
     \plot 1.5 1.5  2.5 1.5 /
     \endpicture} at 1.25 .5
   \put{\beginpicture
     \setcoordinatesystem units <.21cm,.21cm>
     \multiput{} at 0 0  3 6  /
     \setsolid
     \plot 0 0  0 6  1 6  1 0  0 0 /
     \plot 0 1  3 1  3 5  2 5  2 1  /
     \plot 0 2  3 2 /
     \plot 0 3  3 3 /
     \plot 0 4  3 4 /
     \plot 0 5  1 5 /
     \multiput{$\sss\bullet$} at .5 3.5  1.5 3.5  1.5 1.5  2.5 1.5 /
     \plot .5 3.5  1.5 3.5 /
     \plot 1.5 1.5  2.5 1.5 /
     \endpicture} at 1.75 .5
   \put{\beginpicture
     \setcoordinatesystem units <.21cm,.21cm>
     \multiput{} at 0 0  3 6  /
     \setsolid
     \plot 0 0  0 6  1 6  1 0  0 0 /
     \plot 0 1  2 1  2 5  /
     \plot 0 2  3 2  3 5  0 5 /
     \plot 0 3  3 3  /
     \plot 0 4  3 4 /
     \multiput{$\sss\bullet$} at .5 4.5  1.5 4.5  2.5 4.5  1.5 3.5 /
     \plot .5 4.5  2.5 4.5 /
     \endpicture} at -.25 2
   \put{\beginpicture
     \setcoordinatesystem units <.21cm,.21cm>
     \multiput{} at 0 0  3 6  /
     \setsolid
     \plot 0 0  0 6  1 6  1 0  0 0 /
     \plot 2 1  2 6  3 6  3 1  2 1 /
     \plot 0 1  1 1 /
     \plot 0 2  1 2 /
     \plot 2 2  3 2 /
     \plot 0 3  3 3 /
     \plot 0 4  3 4 /
     \plot 0 5  3 5 /
     \multiput{$\sss\bullet$} at .5 4.5  1.5 4.5  1.5 3.5   2.5 3.5 /
     \plot .5 4.5  1.5 4.5 /
     \plot 1.5 3.5  2.5 3.5 /
     \endpicture} at .25 2 
   \put{\beginpicture 
     \setcoordinatesystem units <.21cm,.21cm>
     \multiput{} at 0 0  3 6  /
     \setsolid
     \plot 0 0  0 6  1 6  1 0  0 0 /
     \plot 0 1  3 1  3 3  0 3 /
     \plot 2 1  2 4  0 4 /
     \plot 0 2  3 2 /
     \plot 0 5  1 5 /
     \multiput{$\sss\bullet$} at .5 2.5  1.5 2.5  2.5 2.5  2.5 1.5 /
     \plot .5 2.5  2.5 2.5 /
     \endpicture} at -.5 0 
   \put{\beginpicture
     \setcoordinatesystem units <.21cm,.21cm>
     \multiput{} at 0 0  3 6  /
     \setsolid
     \plot 0 0  0 6  1 6  1 0  0 0 /
     \plot 0 1  2 1  2 5  0 5 /
     \plot 0 3  3 3  3 4  0 4 /
     \plot 0 2  2 2 /
     \multiput{$\sss\bullet$} at .5 3.5  1.5 3.5  2.5 3.5  /
     \plot .5 3.5  2.5 3.5 /
     \endpicture} at 0 0 
   \put{\beginpicture
     \setcoordinatesystem units <.21cm,.21cm>
     \multiput{} at 0 0  3 6  /
     \setsolid
     \plot 0 0  0 6  1 6  1 0  0 0 /
     \plot 0 1  2 1  2 5  0 5 /
     \plot 0 2  3 2  3 3  0 3 /
     \plot 0 4  2 4 /
     \multiput{$\sss\bullet$} at .5 2.5  1.5 2.5  2.5 2.5  1.5 1.5 /
     \plot .5 2.5  2.5 2.5 /
     \endpicture} at .5 0
   \put{\beginpicture  
     \setcoordinatesystem units <.21cm,.21cm>
     \multiput{} at 0 0  3 6  /
     \setsolid
     \plot 0 0  0 6  1 6  1 0  0 0 /
     \plot 0 1  2 1  2 5  0 5 /
     \plot 0 3  3 3  3 4  0 4 /
     \plot 0 2  2 2 /
     \multiput{$\sss\bullet$} at .4 3.4  1.4 3.7  1.7 3.4  2.7 3.7 /
     \plot .4 3.4  1.4 3.7 /
     \plot 1.7 3.4  2.7 3.7 /
     \endpicture} at -2 1.5
   \put{\beginpicture
     \setcoordinatesystem units <.21cm,.21cm>
     \multiput{} at 0 0  3 6  /
     \setsolid
     \plot 0 0  0 6  1 6  1 0  0 0 /
     \plot 0 1  2 1  2 5  0 5 /
     \plot 0 2  3 2  3 3  0 3 /
     \plot 0 4  2 4 /
     \multiput{$\sss\bullet$} at .5 4.5  1.5 4.5  1.5 2.5   2.5 2.5 /
     \plot .5 4.5  1.5 4.5 /
     \plot 1.5 2.5  2.5 2.5 /
     \endpicture} at -1.5 1.5
   \put{\beginpicture
     \setcoordinatesystem units <.21cm,.21cm>
     \multiput{} at 0 0  3 6  /
     \setsolid
     \plot 0 0  0 6  1 6  1 0  0 0 /
     \plot 0 2  2 2  2 5 /
     \plot 0 3  3 3  3 5  0 5 /
     \plot 0 4  3 4  /
     \plot 0 1  1 1 /
     \multiput{$\sss\bullet$} at .5 4.5  1.5 4.5  2.5 4.5  1.5 3.5 /
     \plot .5 4.5  2.5 4.5 /
     \endpicture} at -1 1.5 
   \put{\beginpicture  
     \setcoordinatesystem units <.21cm,.21cm>
     \multiput{} at 0 0  3 6  /
     \setsolid
     \plot 0 0  0 6  1 6  1 0  0 0 /
     \plot 2 1  2 6  3 6  3 1  2 1 /
     \plot 0 1  1 1 /
     \plot 0 2  3 2 /
     \plot 0 3  3 3 /
     \plot 0 4  3 4 /
     \plot 0 5  3 5 /
     \multiput{$\sss\bullet$} at .5 4.5  1.5 4.5  1.5 2.5   2.5 2.5 /
     \plot .5 4.5  1.5 4.5 /
     \plot 1.5 2.5  2.5 2.5 /
     \endpicture} at 1 1.5
   \put{\beginpicture
     \setcoordinatesystem units <.21cm,.21cm>
     \multiput{} at 0 0  3 8  /
     \setsolid
     \plot 0 0  0 6  1 6  1 0  0 0 /
     \plot 2 2  2 8  3 8  3 2  2 2 /
     \plot 0 1  1 1 /
     \plot 0 2  1 2 /
     \plot 0 3  3 3 /
     \plot 0 4  3 4 /
     \plot 0 5  3 5 /
     \plot 2 6  3 6 /
     \plot 2 7  3 7 /
     \multiput{$\sss\bullet$} at .5 4.5  1.5 4.5  1.5 3.5   2.5 3.5 /
     \plot .5 4.5  1.5 4.5 /
     \plot 1.5 3.5  2.5 3.5 /
     \endpicture} at 1.5 1.5 
   \put{\beginpicture
     \setcoordinatesystem units <.21cm,.21cm>
     \multiput{} at 0 0  3 6  /
     \setsolid
     \plot 0 5  2 5  2 0  0 0  0 6  1 6  1 0  /
     \plot 0 1  3 1  3 4  0 4 /
     \plot 0 2  3 2 /
     \plot 0 3  3 3 /
     \multiput{$\sss\bullet$} at .5 3.5  1.5 3.5  2.5 3.5  1.5 2.5  /
     \plot .5 3.5  2.5 3.5 /
     \endpicture} at 2 1.5
   \endpicture}
   $$

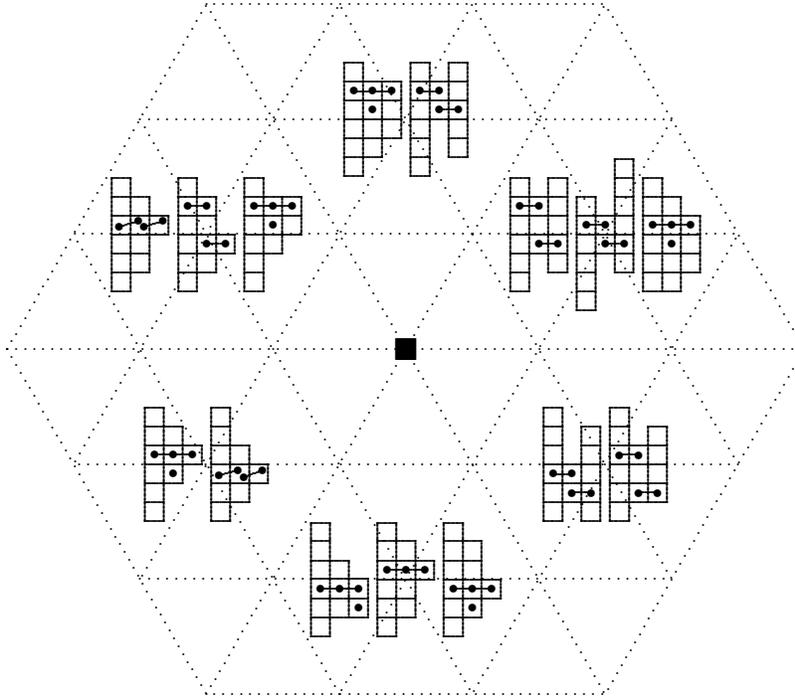
\subsection{The non-central tripickets in $\Delta_{4/3}$ on axis-parallel lines.}
\label{app-B-three}

$$
{\beginpicture
   \setcoordinatesystem units <1.1547cm, 2cm>
   \multiput{} at -2 0  2 2 /
   \plot -2 1  -1 2  1 2  2 1  1 0  -1 0  -2 1 /
   \multiput{$\bigcirc$} at -.67 .33  -1.33 1  -.67 1.67  .67 1.67  1.33 1  .67 .33 /
   \multiput{{\tt 3}} at   -.67 .13  -1.33 .8  -.67 1.47  .67 1.47  1.33 .8  .67 .13 /
     \put{$\blacksquare$} at 0 1
   \setdots<.5mm>
   \plot -2 1  2 1 /
   \plot -1 0  1 2 /
   \plot -1 2  1 0 /
   \plot -.33 0  -1.67 1.33  1.67 1.33   .33 0  -1.33 1.67  1.33 1.67  -.33 0 /
   \plot -1.33 .33  .33 2  1.67 .67  -1.67 .67  -.33 2  1.33 .33  -1.33 .33 / 
   \put{$\ss (\frac43,2)$} at -1.2 -.2 
   \put{$\ss (\frac43,\frac83)$} at 1.2 -.2 
   \put{$\ss (2,\frac43)$} at -2.53 1
   \put{$\ss (\frac83,\frac43)$} at -1.2 2.2 
   \put{$\ss (2,\frac83)$} at 2.53 1
   \put{$\ss (\frac83,2)$} at 1.2 2.2 
\endpicture}
$$

\medskip
At each position, there are three tripickets,
together they form three $\Sigma_3$-orbits.
In the fundamental domain, the object pictured on the left is on
the mouth of the 6-tube of index $\frac14$, the object in the middle has
quasi-length 3 in the 6-tube of index $\frac11$, and the object on the right lies
on the mouth of the 6-tube of index $\frac25$.

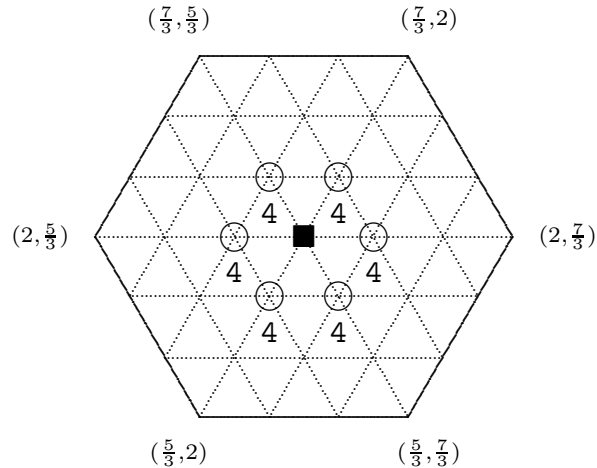
\medskip

$$
{\beginpicture
   \setcoordinatesystem units <1.4722cm, 2.55cm>
   \multiput{} at -3 -.5  3 3 /
   \setdots<1mm>
   \plot -3 1  -1.5 2.5  1.5 2.5  3 1  1.5 -.5  -1.5 -.5  -3 1 /
   \put{$\blacksquare$} at 0 1
   \plot -3 1  3 1 /
   \plot -1.5 -.5  1.5 2.5 /
   \plot -1.5 2.5  1.5 -.5 /
   \plot -.5  -.5  -2.5 1.5  2.5 1.5  .5 -.5  -2 2  2 2  -.5 -.5 /
   \plot -2 0  2 0  -.5 2.5  -2.5 .5  2.5 .5  .5 2.5  -2 0 /
   \put{\beginpicture
     \setcoordinatesystem units <.21cm,.21cm>
     \multiput{} at 0 0  3 6  /
     \setsolid
     \plot 0 0  0 6  1 6  1 0  0 0 /
     \plot 0 2  2 2  2 5  0 5 /
     \plot 0 3  3 3  3 4  0 4 /
     \plot 0 1  1 1 /
     \multiput{$\sss\bullet$} at .5 3.5  1.5 3.5  2.5 3.5 /
     \plot .5 3.5  2.5 3.5 /
     \endpicture} at -1.5 0
   \put{\beginpicture
     \setcoordinatesystem units <.21cm,.21cm>
     \multiput{} at 0 0  3 6  /
     \setsolid
     \plot 0 0  0 6  1 6  1 0  0 0 /
     \plot 0 1  2 1  2 4  0 4 /
     \plot 0 3  3 3  3 2  0 2 /
     \plot 0 5  1 5 /
     \multiput{$\sss\bullet$} at .5 2.5  1.5 2.5  2.5 2.5  1.5 1.5 /
     \plot .5 2.5  2.5 2.5 /
     \endpicture} at -1 0
   \put{\beginpicture
     \setcoordinatesystem units <.21cm,.21cm>
     \multiput{} at 0 0  3 5  /
     \setsolid
     \plot 0 0  0 5  1 5  1 0  0 0 /
     \plot 0 2  3 2 /
     \plot 0 3  3 3  3 1  0 1 /
     \plot 0 4  2 4  2 1 /
     \multiput{$\sss\bullet$} at .5 2.5  1.5 2.5  2.5 2.5  2.5 1.5 /
     \plot .5 2.5  2.5 2.5 /
     \endpicture} at -.5 0
   \put{\beginpicture
     \setcoordinatesystem units <.21cm,.21cm>
     \multiput{} at 0 0  3 6  /
     \setsolid
     \plot 0 0  0 6  1 6  1 0  0 0 /
     \plot 0 1  2 1  2 5  0 5 /
     \plot 0 2  3 2  3 4  0 4 /
     \plot 0 3  3 3 /
     \multiput{$\sss\bullet$} at .5 3.5  1.5 3.5  2.5 3.5  /
     \plot .5 3.5  2.5 3.5 /
     \endpicture} at .5 0
   \put{\beginpicture
     \setcoordinatesystem units <.21cm,.21cm>
     \multiput{} at 0 0  3 6  /
     \setsolid
     \plot 0 0  0 6  1 6  1 0  0 0 /
     \plot 0 1  3 1  3 3  0 3 /
     \plot 2 1  2 5  0 5 / 
     \plot 0 2  3 2 /
     \plot 0 4  2 4 /
     \multiput{$\sss\bullet$} at .5 2.5  1.5 2.5  2.5 2.5  2.5 1.5 /
     \plot .5 2.5  2.5 2.5 /
     \endpicture} at 1 0
   \put{\beginpicture
     \setcoordinatesystem units <.21cm,.21cm>
     \multiput{} at 0 0  3 6  /
     \setsolid
     \plot 0 0  0 6  1 6  1 0  0 0 /
     \plot 0 1  2 1  2 5  0 5 /
     \plot 0 2  3 2  3 4  0 4 /
     \plot 0 3  3 3 /
     \multiput{$\sss\bullet$} at .5 2.5  1.5 2.5  2.5 2.5  1.5 1.5 /
     \plot .5 2.5  2.5 2.5 /
     \endpicture} at 1.5 0
   \put{\beginpicture
     \setcoordinatesystem units <.21cm,.21cm>
     \multiput{} at 0 -1  3 6  /
     \setsolid
     \plot 0 0  0 6  1 6  1 0  0 0 /
     \plot 2 -1  2 5  3 5  3 -1  2 -1 /
     \plot 2 0  3 0 /
     \plot 0 1  3 1 /
     \plot 0 2  3 2 /
     \plot 0 3  3 3 /
     \plot 0 4  1 4 /
     \plot 2 4  3 4 /
     \plot 0 5  1 5 /
     \multiput{$\sss\bullet$} at .5 1.5  1.5 1.5  1.5 2.5  2.5 2.5 /
     \plot .5 1.5  1.5 1.5 /
     \plot 1.5 2.5  2.5 2.5 /
     \endpicture} at 1.5 1
   \put{\beginpicture
     \setcoordinatesystem units <.21cm,.21cm>
     \multiput{} at 0 0  3 6  /
     \setsolid
     \plot 0 0  0 6  1 6  1 0  0 0 /
     \plot 2 0  2 5  3 5  3 0  2 0 /
     \plot 0 1  3 1 /
     \plot 0 2  3 2 /
     \plot 0 3  3 3 /
     \plot 0 4  3 4 /
     \plot 0 5  1 5 /
     \multiput{$\sss\bullet$} at  1.5 3.5  2.5 3.5  1.5 1.5  .5 1.5 /
     \plot 1.5 3.5  2.5 3.5 /
     \plot .5 1.5  1.5 1.5 /
     \endpicture} at 2 1
   \put{\beginpicture
     \setcoordinatesystem units <.21cm,.21cm>
     \multiput{} at 0 0  3 6  /
     \setsolid
     \plot 0 0  0 6  1 6  1 0  0 0 /
     \plot 2 1  2 6  3 6  3 1  2 1 /
     \plot 0 1  1 1 /
     \plot 0 2  3 2 /
     \plot 0 3  3 3 /
     \plot 0 4  3 4 /
     \plot 0 5  3 5 /
     \multiput{$\sss\bullet$} at .5 3.5  1.5 3.5  1.5 2.5  2.5 2.5 /
     \plot .5 3.5  1.5 3.5 /
     \plot 1.5 2.5  2.5 2.5 /
     \endpicture} at 2.5 1
   \put{\beginpicture
     \setcoordinatesystem units <.21cm,.21cm>
     \multiput{} at 0 -1  3 6  /
     \setsolid
     \plot 0 -1  0 5  1 5  1 -1  0 -1 /
     \plot 2 0  2 6  3 6  3 0  2 0 /
     \plot 0 0  1 0 /
     \plot 0 1  1 1 /
     \plot 2 1  3 1 /
     \plot 0 2  3 2 /
     \plot 0 3  3 3 /
     \plot 0 4  3 4 /
     \plot 2 5  3 5 /
     \multiput{$\sss\bullet$} at 2.5 2.5  1.5 2.5  1.5 3.5  .5 3.5   /
     \plot 1.5 3.5  .5 3.5 /
     \plot 2.5 2.5  1.5 2.5 /
     \endpicture} at .5 2
   \put{\beginpicture
     \setcoordinatesystem units <.21cm,.21cm>
     \multiput{} at 0 0  3 6  /
     \setsolid
     \plot 0 0  0 6  1 6  1 0  0 0 /
     \plot 2 1  2 6  3 6  3 1  2 1 /
     \plot 0 1  1 1 /
     \plot 0 2  3 2 /
     \plot 0 3  3 3 /
     \plot 0 4  3 4 /
     \plot 0 5  3 5 /
     \multiput{$\sss\bullet$} at .5 4.5  1.5 4.5  2.5 3.5  1.5 3.5 /
     \plot 1.5 3.5  2.5 3.5 /
     \plot .5 4.5  1.5 4.5 /
     \endpicture} at 1 2
   \put{\beginpicture
     \setcoordinatesystem units <.21cm,.21cm>
     \multiput{} at 0 0  3 6  /
     \setsolid
     \plot 0 0  0 6  1 6  1 0  0 0 /
     \plot 2 0  2 5  3 5  3 0  2 0 /
     \plot 0 1  3 1 /
     \plot 0 2  3 2 /
     \plot 0 3  3 3 /
     \plot 0 4  3 4 /
     \plot 0 5  1 5 /
     \multiput{$\sss\bullet$} at  1.5 3.5  2.5 3.5  1.5 1.5  .5 2.5  1.5 2.5 /
     \plot 1.5 3.5  2.5 3.5 /
     \plot .5 2.5  1.5 2.5 / 
     \endpicture} at 1.5 2
   \put{\beginpicture
     \setcoordinatesystem units <.21cm,.21cm>
     \multiput{} at 0 0  3 6  /
     \setsolid
     \plot 0 0  0 6  1 6  1 0  0 0 /
     \plot 0 1  2 1  2 5  0 5 /
     \plot 0 2  3 2  3 4  0 4 /
     \plot 0 3  3 3 /
     \multiput{$\sss\bullet$} at .5 4.5  1.5 4.5  1.5 3.5  2.5 3.5 /
     \plot 1.5 3.5  2.5 3.5 /
     \plot .5 4.5  1.5 4.5 /
     \endpicture} at -1.5 2
   \put{\beginpicture
     \setcoordinatesystem units <.21cm,.21cm>
     \multiput{} at 0 0  3 6  /
     \setsolid
     \plot 0 0  0 6  1 6  1 0  0 0 /
     \plot 0 1  2 1  2 5 /
     \plot 0 2  2 2 /
     \plot 0 3  3 3  3 5  0 5 /
     \plot 0 4  3 4 /
     \multiput{$\sss\bullet$} at .5 4.5  1.5 4.5  2.5 4.5  1.5 3.5 /
     \plot .5 4.5  2.5 4.5 /
     \endpicture} at -1 2
   \put{\beginpicture
     \setcoordinatesystem units <.21cm,.21cm>
     \multiput{} at 0 0  3 6  /
     \setsolid
     \plot 0 0  0 6  1 6  1 0  0 0 /
     \plot 0 1  2 1  2 5  0 5 /
     \plot 0 2  3 2  3 4  0 4 /
     \plot 0 3  3 3 /
     \multiput{$\sss\bullet$} at .4  3.4  1.4 3.7  1.7 3.4  2.7 3.7  2.5 2.5 /
     \plot .4 3.4  1.4 3.7 /
     \plot 1.7 3.4  2.7 3.7 /
     \endpicture} at -.5 2
   \put{\beginpicture
     \setcoordinatesystem units <.21cm,.21cm>
     \multiput{} at 0 0  3 6  /
     \setsolid
     \plot 0 0  0 6  1 6  1 0  0 0 /
     \plot 0 1  2 1  2 4  0 4 /
     \plot 0 3  3 3  3 2  0 2 /
     \plot 0 5  1 5 /
     \multiput{$\sss\bullet$} at .5 3.5  1.5 3.5  2.5 2.5  1.5 2.5 /
     \plot .5 3.5  1.5 3.5 /
     \plot 1.5 2.5  2.5 2.5 /
     \endpicture} at -2.5 1
   \put{\beginpicture
     \setcoordinatesystem units <.21cm,.21cm>
     \multiput{} at 0 0  3 6  /
     \setsolid
     \plot 0 0  0 6  1 6  1 0  0 0 /
     \plot 0 2  2 2  2 5  0 5 /
     \plot 0 3  3 3  3 4  0 4 /
     \plot 0 1  1 1 /
     \multiput{$\sss\bullet$} at .4 3.4  1.4 3.7  1.7 3.3  2.7 3.6 /
     \plot .4 3.4  1.4 3.7 /
     \plot 1.7 3.3  2.7 3.6 /
     \endpicture} at -2 1
   \put{\beginpicture
     \setcoordinatesystem units <.21cm,.21cm>
     \multiput{} at 0 0  3 5  /
     \setsolid
     \plot 0 0  0 5  1 5  1 0  0 0 /
     \plot 0 2  3 2  3 4  0 4 /
     \plot 0 3  3 3  /
     \plot 0 1  2 1  2 4  /
     \multiput{$\sss\bullet$} at .5 3.5  1.5 3.5  2.5 2.5  1.5 2.5 /
     \plot .5 3.5  1.5 3.5 /
     \plot 1.5 2.5  2.5 2.5 /
     \endpicture} at -1.5 1
\endpicture}
$$

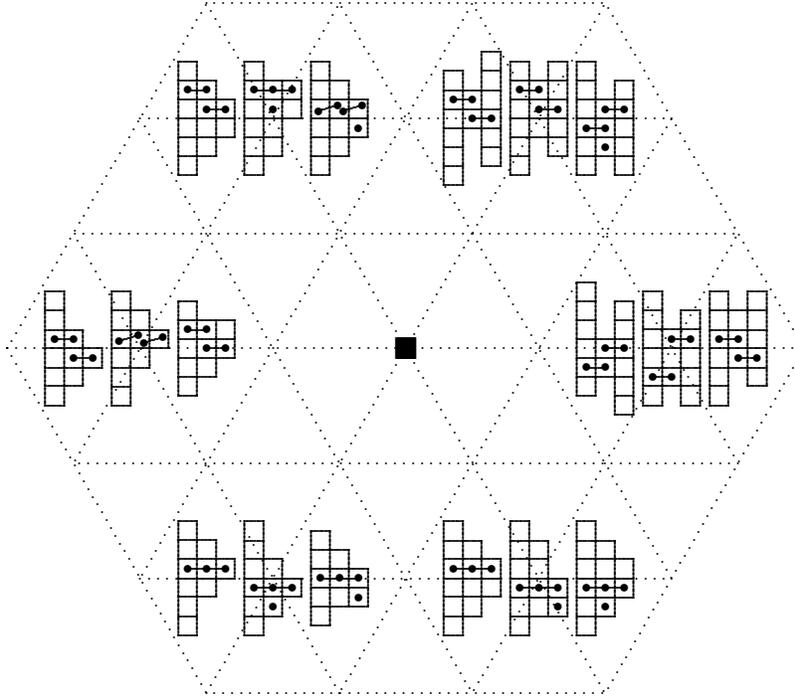
\subsection{The non-central tripickets in $\Delta_{5/3}$ on axis-parallel lines.}
\label{app-B-four}

$$
{\beginpicture
   \setcoordinatesystem units <1.1547cm, 2cm>
   \multiput{} at -2 0  2 2 /
   \plot -2 1  -1 2  1 2  2 1  1 0  -1 0  -2 1 /
   \multiput{$\bigcirc$} at -.33 .67  -.67 1  -.33 1.33  .33 1.33  .67 1  .33 .67 /
   \multiput{{\tt 4}} at -.33 .47  -.67 .8  -.33 1.13  .33 1.13  .67 .8  .33 .47 /
     \put{$\blacksquare$} at 0 1
   \setdots<.5mm>
   \plot -2 1  2 1 /
   \plot -1 0  1 2 /
   \plot -1 2  1 0 /
   \plot -.33 0  -1.67 1.33  1.67 1.33   .33 0  -1.33 1.67  1.33 1.67  -.33 0 /
   \plot -1.33 .33  .33 2  1.67 .67  -1.67 .67  -.33 2  1.33 .33  -1.33 .33 / 
   \put{$\ss (\frac53,2)$} at -1.2 -.2 
   \put{$\ss (\frac53,\frac73)$} at 1.2 -.2 
   \put{$\ss (2,\frac53)$} at -2.53 1
   \put{$\ss (\frac73,\frac53)$} at -1.2 2.2 
   \put{$\ss (2,\frac73)$} at 2.53 1
   \put{$\ss (\frac73,2)$} at 1.2 2.2 
\endpicture}
$$

\medskip
There are four tripickets at each position, they form four  $\Sigma_3$-orbits.
Among the objects in the fundamental domain, 
the first occurs on the mouth of the 6-tube of index $\frac15$,
the second on the mouth of the 6-tube of index $\frac51$.  (Thus the first two orbits contain
all the modules on the mouth of those two tubes.)  The third object
occurs on the 5th layer in $\Cal P(6)$; the last in the tube of index $\frac 16$.

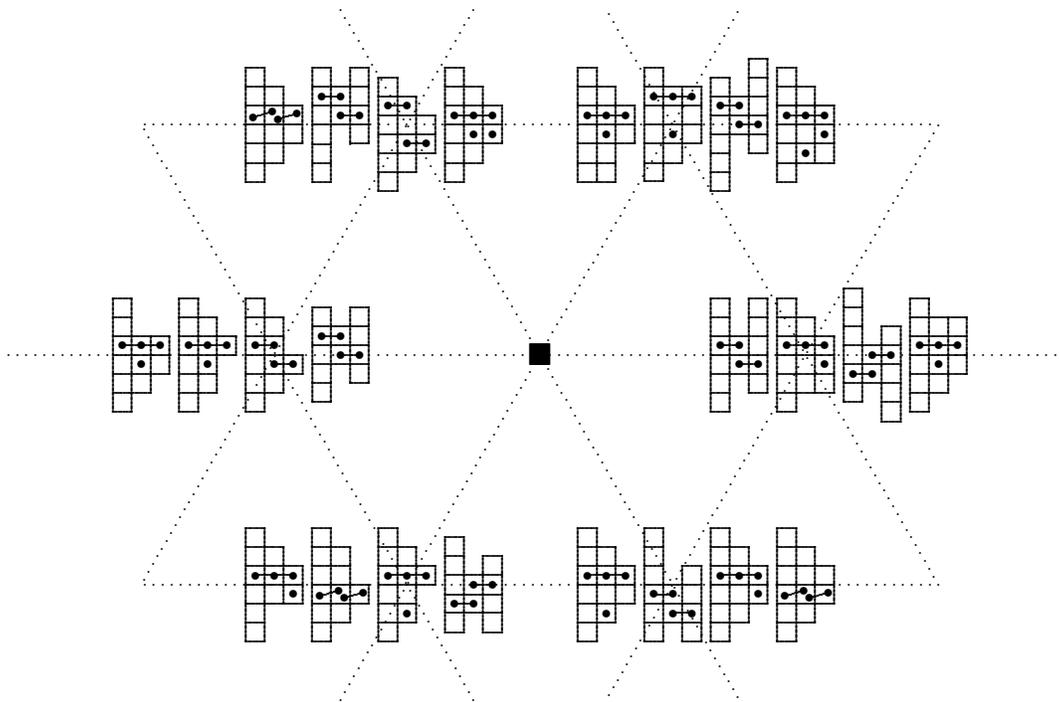
$$
{\beginpicture
   \setcoordinatesystem units <1.4722cm, 2.55cm>
   \multiput{} at -3 -.5  3 3 /
   \setdots<1mm>
   \put{$\blacksquare$} at 0 1
   \plot -4 1  4 1 /
   \plot -.5 2.5  -3 0  3 0  .5 2.5 /
   \plot -.5 -.5  -3 2  3 2  .5 -.5 /
   \plot -1.5 -.5  1.5 2.5 /
   \plot -1.5 2.5  1.5 -.5 /
   \put{\beginpicture
     \setcoordinatesystem units <.21cm,.21cm>
     \multiput{} at 0 0  3 6  /
     \setsolid
     \plot 0 0  0 6  1 6  1 0  0 0 /
     \plot 0 2  3 2  3 4  0 4 /
     \plot 2 2  2 5  0 5 /
     \plot 0 1  1 1 /
     \plot 0 3  3 3 /
     \multiput{$\sss\bullet$} at .5 3.5  1.5 3.5  2.5 3.5  2.5 2.5  /
     \plot .5 3.5  2.5 3.5 /
     \endpicture} at -2 0
   \put{\beginpicture
     \setcoordinatesystem units <.21cm,.21cm>
     \multiput{} at 0 0  3 6  /
     \setsolid
     \plot 0 0  0 6  1 6  1 0  0 0 /
     \plot 0 1  2 1  2 5  0 5 /
     \plot 0 2  3 2  3 3  0 3 /
     \plot 0 4  2 4 /
     \multiput{$\sss\bullet$} at .4 2.4  1.4 2.7  1.7 2.3  2.7 2.6 /
     \plot .4 2.4  1.4 2.7 /
     \plot 1.7 2.3  2.7 2.6 /
     \endpicture} at -1.5 0
   \put{\beginpicture
     \setcoordinatesystem units <.21cm,.21cm>
     \multiput{} at 0 0  3 6  /
     \setsolid
     \plot 0 0  0 6  1 6  1 0  0 0 /
     \plot 0 1  2 1  2 5  0 5 /
     \plot 0 3  3 3  3 4  0 4 /
     \plot 0 2  2 2 /
     \multiput{$\sss\bullet$} at .5 3.5  1.5 3.5  2.5 3.5  1.5 1.5  /
     \plot .5 3.5  2.5 3.5 /
     \endpicture} at -1 0
   \put{\beginpicture
     \setcoordinatesystem units <.21cm,.21cm>
     \multiput{} at 0 0  3 5  /
     \setsolid
     \plot 0 0  0 5  1 5  1 0  0 0 /
     \plot 2 0  2 4  3 4  3 0  2 0 /
     \plot 0 1  3 1 /
     \plot 0 2  3 2 /
     \plot 0 3  3 3 /
     \plot 0 4  1 4 /
     \multiput{$\sss\bullet$} at .5 1.5  1.5 1.5  1.5 2.5  2.5 2.5 /
     \plot .5 1.5  1.5 1.5 /
     \plot 1.5 2.5  2.5 2.5 /
     \endpicture} at -.5 0
   \put{\beginpicture
     \setcoordinatesystem units <.21cm,.21cm>
     \multiput{} at 0 0  3 6  /
     \setsolid
     \plot 0 0  0 6  1 6  1 0  0 0 /
     \plot 0 1  2 1  2 5  0 5 /
     \plot 0 2  3 2  3 4  0 4 /
     \plot 0 3  3 3 /
     \multiput{$\sss\bullet$} at .5 3.5  1.5 3.5  2.5 3.5  1.5 1.5 /
     \plot .5 3.5  2.5 3.5 /
     \endpicture} at .5 0
   \put{\beginpicture
     \setcoordinatesystem units <.21cm,.21cm>
     \multiput{} at 0 0  3 6  /
     \setsolid
     \plot 0 0  0 6  1 6  1 0  0 0 /
     \plot 2 0  2 4  3 4  3 0  2 0 /
     \plot 0 1  3 1 /
     \plot 0 2  3 2 /
     \plot 0 3  3 3 /
     \plot 0 4  1 4 /
     \plot 0 5  1 5 /
     \multiput{$\sss\bullet$} at .5 2.5  1.5 2.5  1.5 1.5  2.5 1.5 /
     \plot .5 2.5  1.5 2.5 /
     \plot 1.5 1.5  2.5 1.5 /
     \endpicture} at 1 0
   \put{\beginpicture
     \setcoordinatesystem units <.21cm,.21cm>
     \multiput{} at 0 0  3 6  /
     \setsolid
     \plot 0 0  0 6  1 6  1 0  0 0 /
     \plot 0 1  2 1  2 5  0 5 /
     \plot 0 2  3 2  3 4  0 4 /
     \plot 0 3  3 3 /
     \multiput{$\sss\bullet$} at .5 3.5  1.5 3.5  2.5 3.5  2.5 2.5 /
     \plot .5 3.5  2.5 3.5 /
     \endpicture} at 1.5 0
   \put{\beginpicture
     \setcoordinatesystem units <.21cm,.21cm>
     \multiput{} at 0 0  3 6  /
     \setsolid
     \plot 0 0  0 6  1 6  1 0  0 0 /
     \plot 0 1  2 1  2 5  0 5 /
     \plot 0 2  3 2  3 4  0 4 /
     \plot 0 3  3 3  /
     \multiput{$\sss\bullet$} at .4 2.4  1.4  2.7  1.7 2.3  2.7 2.6 /
     \plot .4 2.4  1.4 2.7 /
     \plot 1.7 2.3  2.7 2.6 /
     \endpicture} at 2 0
   \put{\beginpicture
     \setcoordinatesystem units <.21cm,.21cm>
     \multiput{} at 0 0  3 6  /
     \setsolid
     \plot 0 0  0 6  1 6  1 0  0 0 /
     \plot 2 1  2 6  3 6  3 1  2 1 /
     \plot 0 1  1 1 /
     \plot 0 2  3 2 /
     \plot 0 3  3 3 /
     \plot 0 4  3 4 /
     \plot 0 5  1 5 /
     \plot 2 5  3 5 /
     \multiput{$\sss\bullet$} at .5 3.5  1.5 3.5  1.5 2.5  2.5 2.5 /
     \plot .5 3.5  1.5 3.5 /
     \plot 1.5 2.5  2.5 2.5 /
     \endpicture} at 1.5 1
   \put{\beginpicture
     \setcoordinatesystem units <.21cm,.21cm>
     \multiput{} at 0 0  3 6  /
     \setsolid
     \plot 0 0  0 6  1 6  1 0  0 0 /
     \plot 0 1  3 1  3 4  0 4 /
     \plot 2 1  2 5  0 5 /
     \plot 0 3  3 3  /
     \plot 0 2  3 2 /
     \multiput{$\sss\bullet$} at .5 3.5  1.5 3.5  2.5 3.5  2.5 2.5 /
     \plot .5 3.5  2.5 3.5 /
     \endpicture} at 2 1
   \put{\beginpicture
     \setcoordinatesystem units <.21cm,.21cm>
     \multiput{} at 0 0  3 7  /
     \setsolid
     \plot 0 1  0 7  1 7  1 1  0 1 /
     \plot 2 0  2 5  3 5  3 0  2 0 /
     \plot 2 1  3 1 /
     \plot 0 2  3 2 /
     \plot 0 3  3 3 /
     \plot 0 4  3 4 /
     \plot 0 5  1 5 /
     \plot 0 6  1 6 /
     \multiput{$\sss\bullet$} at .5 2.5  1.5 2.5   1.5 3.5  2.5 3.5 /
     \plot 1.5 3.5  2.5 3.5 /
     \plot .5 2.5  1.5 2.5 /
     \endpicture} at 2.5 1
   \put{\beginpicture
     \setcoordinatesystem units <.21cm,.21cm>
     \multiput{} at 0 0  3 6  /
     \setsolid
     \plot 0 0  0 6  1 6  1 0  0 0 /
     \plot 0 1  2 1  2 5 /
     \plot 0 2  3 2  3 5  0 5 /
     \plot 0 3  3 3 /
     \plot 0 4  3 4 /
     \multiput{$\sss\bullet$} at .5 3.5  1.5 3.5  2.5 3.5  1.5 2.5 /
     \plot .5 3.5  2.5 3.5 /
     \endpicture} at 3 1
   \put{\beginpicture
     \setcoordinatesystem units <.21cm,.21cm>
     \multiput{} at 0 0  3 6  /
     \setsolid
     \plot 0 5  2 5  2 0  0 0  0 6  1 6  1 0 /
     \plot 0 1  2 1 /
     \plot 0 2  3 2  3 4  0 4 /
     \plot 0 3  3 3 /
     \multiput{$\sss\bullet$} at .5 3.5  1.5 3.5  2.5 3.5  1.5 2.5 /
     \plot .5 3.5  2.5 3.5 /
     \endpicture} at .5 2
   \put{\beginpicture
     \setcoordinatesystem units <.21cm,.21cm>
     \multiput{} at 0 0  3 6  /
     \setsolid
     \plot 0 0  0 6  1 6  1 0  0 0 /
     \plot 0 1  2 1  2 5 /
     \plot 0 2  3 2  3 5  0 5 /
     \plot 0 3  3 3 /
     \plot 0 4  3 4 /
     \multiput{$\sss\bullet$} at .5 4.5  1.5 4.5  2.5 4.5  1.5 2.5  /
     \plot .5 4.5  2.5 4.5 /
     \endpicture} at 1 2
   \put{\beginpicture  
     \setcoordinatesystem units <.21cm,.21cm>
     \multiput{} at 0 0  3 7  /
     \setsolid
     \plot 0 0  0 6  1 6  1 0  0 0 /
     \plot 0 1  1 1 /
     \plot 0 2  1 2 /
     \plot 2 2  2 7  3 7  3 2  2 2 /
     \plot 0 3  3 3  /
     \plot 0 4  3 4 /
     \plot 0 5  3 5 /
     \plot 2 6  3 6 /
     \multiput{$\sss\bullet$} at .5 4.5  1.5 4.5  1.5 3.5  2.5 3.5 /
     \plot 1.5 3.5  2.5 3.5 /
     \plot .5 4.5  1.5 4.5 /
     \endpicture} at 1.5 2
   \put{\beginpicture
     \setcoordinatesystem units <.21cm,.21cm>
     \multiput{} at 0 0  3 6  /
     \setsolid
     \plot 0 0  0 6  1 6  1 0  0 0 /
     \plot 0 1  3 1  3 4  0 4 /
     \plot 2 1  2 5  0 5 /
     \plot 0 2  3 2 /
     \plot 0 3  3 3 /
     \plot 0 4  3 4 /
     \multiput{$\sss\bullet$} at .5 3.5  1.5 3.5  2.5 3.5  1.5 1.5  2.5 2.5 /
     \plot .5 3.5  2.5 3.5 /
     \endpicture} at 2 2 
   \put{\beginpicture
     \setcoordinatesystem units <.21cm,.21cm>
     \multiput{} at 0 0  3 6  /
     \setsolid
     \plot 0 0  0 6  1 6  1 0  0 0 /
     \plot 0 1  2 1  2 5  0 5 /
     \plot 0 2  3 2  3 4  0 4 /
     \plot 0 3  3 3 /
     \multiput{$\sss\bullet$} at .4 3.4  1.4 3.7  1.7 3.3  2.7 3.6  /
     \plot .4 3.4  1.4 3.7 /
     \plot 1.7 3.3  2.7 3.6 /
     \endpicture} at -2 2 
   \put{\beginpicture
     \setcoordinatesystem units <.21cm,.21cm>
     \multiput{} at 0 0  3 6  /
     \setsolid
     \plot 0 0  0 6  1 6  1 0  0 0 /
     \plot 0 1  1 1 /
     \plot 0 2  1 2 /
     \plot 2 2  2 6  3 6  3 2  2 2 /
     \plot 0 3  3 3 /
     \plot 0 4  3 4 /
     \plot 0 5  3 5 /
     \multiput{$\sss\bullet$} at .5 4.5  1.5 4.5  1.5 3.5  2.5 3.5 /
     \plot .5 4.5  1.5 4.5 /
     \plot 1.5 3.5  2.5 3.5 /
     \endpicture} at -1.5 2
   \put{\beginpicture
     \setcoordinatesystem units <.21cm,.21cm>
     \multiput{} at 0 0  3 7  /
     \setsolid
     \plot 0 0  0 6  1 6  1 0  0 0 /
     \plot 0 1  2 1  2 5  0 5 /
     \plot 0 2  3 2  3 4  0 4 /
     \plot 0 3  3 3 /
     \multiput{$\sss\bullet$} at .5 4.5  1.5 4.5   1.5 2.5  2.5 2.5 /
     \plot .5 4.5  1.5 4.5 /
     \plot 1.5 2.5  2.5 2.5 /
     \endpicture} at -1 2 
   \put{\beginpicture
     \setcoordinatesystem units <.21cm,.21cm>
     \multiput{} at 0 0  3 6  /
     \setsolid
     \plot 0 0  0 6  1 6  1 0  0 0 /
     \plot 0 1  2 1  2 5  0 5 /
     \plot 0 2  3 2  3 4  0 4 /
     \plot 0 3  3 3 /
     \multiput{$\sss\bullet$} at .5 3.5  1.5 3.5  2.5 3.5  1.5 2.5  2.5 2.5 /
     \plot .5 3.5  2.5 3.5 /
     \endpicture} at -.5 2
   \put{\beginpicture
     \setcoordinatesystem units <.21cm,.21cm>
     \multiput{} at 0 0  3 6  /
     \setsolid
     \plot 0 0  0 6  1 6  1 0  0 0 /
     \plot 0 1  2 1  2 4 /
     \plot 0 2  3 2  3 4  0 4 /
     \plot 0 3  3 3 /
     \plot 0 5  1 5 /
     \multiput{$\sss\bullet$} at .5 3.5  1.5 3.5  2.5 3.5  1.5 2.5 /
     \plot .5 3.5  2.5 3.5 /
     \endpicture} at -3 1
   \put{\beginpicture
     \setcoordinatesystem units <.21cm,.21cm>
     \multiput{} at 0 0  3 6  /
     \setsolid
     \plot 0 0  0 6  1 6  1 0  0 0 /
     \plot 0 1  2 1  2 5  0 5 /
     \plot 0 3  3 3  3 4  0 4 /
     \plot 0 2  2 2 /
     \multiput{$\sss\bullet$} at .5 3.5  1.5 3.5  2.5 3.5  1.5 2.5 /
     \plot .5 3.5  2.5 3.5 /
     \endpicture} at -2.5 1
   \put{\beginpicture
     \setcoordinatesystem units <.21cm,.21cm>
     \multiput{} at 0 0  3 6  /
     \setsolid
     \plot 0 0  0 6  1 6  1 0  0 0 /
     \plot 0 1  2 1  2 5  0 5 /
     \plot 0 3  3 3  3 2  0 2 /
     \plot 0 4  2 4 /
     \multiput{$\sss\bullet$} at .5 3.5  1.5 3.5  1.5 2.5  2.5 2.5 /
     \plot .5 3.5  1.5 3.5 /
     \plot 1.5 2.5  2.5 2.5 /
     \endpicture} at -2 1
   \put{\beginpicture
     \setcoordinatesystem units <.21cm,.21cm>
     \multiput{} at 0 0  3 5  /
     \setsolid
     \plot 0 0  0 5  1 5  1 0  0 0 /
     \plot 2 1  2 5  3 5  3 1  2 1 /
     \plot 0 1  1 1 /
     \plot 0 2  3 2 /
     \plot 0 3  3 3 /
     \plot 0 4  3 4 /
     \multiput{$\sss\bullet$} at .5 3.5  1.5 3.5  1.5 2.5  2.5 2.5 /
     \plot .5 3.5  1.5 3.5 /
     \plot 1.5 2.5  2.5 2.5 /
     \endpicture} at -1.5 1
\endpicture}
$$


\subsection{Some tetrapickets in $\Cal S(6)$.}
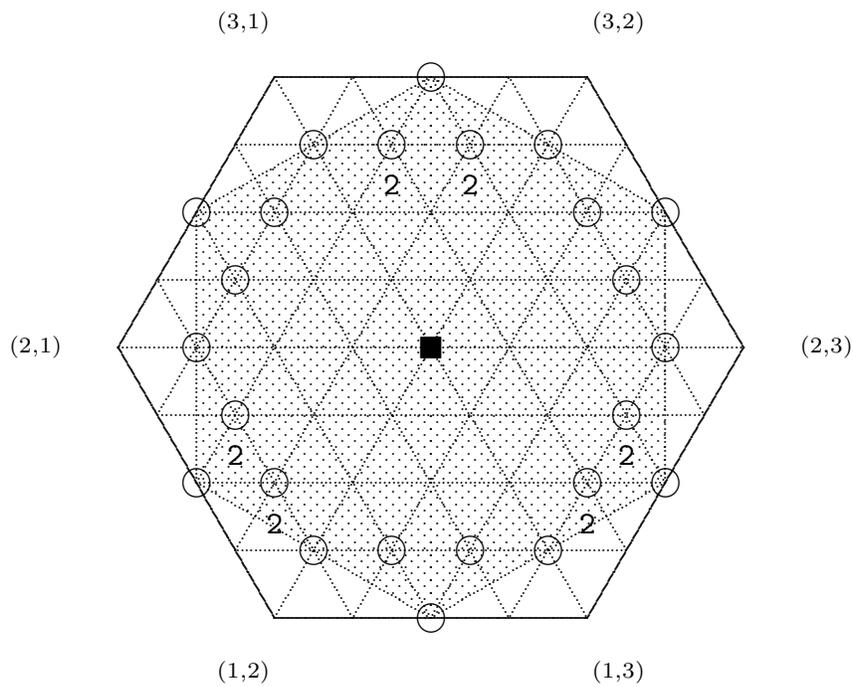
\label{app-B-five}

$$
{\beginpicture
   \setcoordinatesystem units <1.732cm, 3cm>
   \multiput{} at -2 0  2 2 /
   \plot -2 1  -1 2  1 2  2 1  1 0  -1 0  -2 1 /
   \multiput{$\bigcirc$} at  0 0  -.75 .25  -.25 .25  .25 .25  .75 .25
                         -1.5 1.5  -1.5 1  -1.25 1.25  -1 1.5  -.75 1.75
                         1.5 1.5  1.5 1  1.25 1.25  1 1.5  .75 1.75 
                         -1 .5  -1.25 .75  -.25 1.75  .25 1.75  1 .5  1.25 .75 /
   \multiput{$\bigcirc$} at  0 2  -1.5 0.5  1.5 0.5  /
     \put{$\blacksquare$} at 0 1
     \multiput{{\tt 2}} at -1 .35  -1.25 .6  -.25 1.6  .25 1.6  1 .35  1.25 .6 /
   \setdots<.5mm>
   \plot -2 1  2 1 /
   \plot -1 0  1 2 /
   \plot -1 2  1 0 /
   \plot 0 0  -1.5 1.5  1.5 1.5  0 0 /
   \plot -.5 0  -1.75 1.25  1.75 1.25  .5 0  -1.25 1.75  1.25 1.75  -.5 0 /
   \plot -1.25 .25  .5 2  1.75 .75  -1.75 .75  -.5 2   1.25 .25   -1.25 .25 /
   \plot -1.5 .5  0 2  1.5 .5  -1.5 .5 /
   \setshadegrid span <.7mm>
   \vshade -1.5 .5 1.5  <,z,,>  0 0 2   <z,z,,>  1.5 .5 1.5 /
   \setdots <.5mm>
   \plot -1.5 .5  0 0  1.5 .5  1.5 1.5  0 2  -1.5 1.5  -1.5 0.5 /
   \put{$\ss (1,2)$} at -1.2 -.2 
   \put{$\ss (1,3)$} at 1.2 -.2 
   \put{$\ss (2,1)$} at -2.53 1
   \put{$\ss (3,1)$} at -1.2 2.2 
   \put{$\ss (2,3)$} at 2.53 1
   \put{$\ss (3,2)$} at 1.2 2.2 
\endpicture}
$$

\medskip
All tetrapickets in $\Cal S(6)$ occur inside the hexagon, each on an intersection point
of dotted lines in the above diagram.
In Section \ref{sec-fifteen-two}--(d) we have already pictured the six tetrapickets on the boundary
of the hexagon.

\medskip
The 12 tetrapickets on the boundary of $\Delta_{5/4}$ are presented in Section \ref{sec-fifteen-two}--(f).
There are no other objects on the boundary of $\Delta_{5/4}$,
and the only objects outside $\Delta_{5/4}$ are the pickets on $\Delta_0$ and the objects on $\Delta_1$
shown in Section~\ref{sec-fifteen-two}.

\medskip
There are 12 additional tetrapickets on the lines $q=\frac{13}4$, $p=\frac{11}4$
and $r=\frac{11}4$ (namely at the points marked by ``2'').
They lie on the half-lines connecting the center
with the positions of the boundary pickets in the $\Sigma_3$-orbit of the simple object $S$.
At each position,  two tetrapickets occur;
each dimension vector is the sum of the root for the corresponding boundary picket
and one of the two radical vectors $\mathbf h_0$ or $\mathbf h_\infty$ for $\chi_\Theta$.

\medskip
The objects lie on two orbits under the action of $\Sigma_3$.
The module pictured in the bottom row on the left occurs on the mouth of the 6-tube of index $\frac 61$,
its neighbor has quasi-length 7 in the principal component.

\medskip
$$
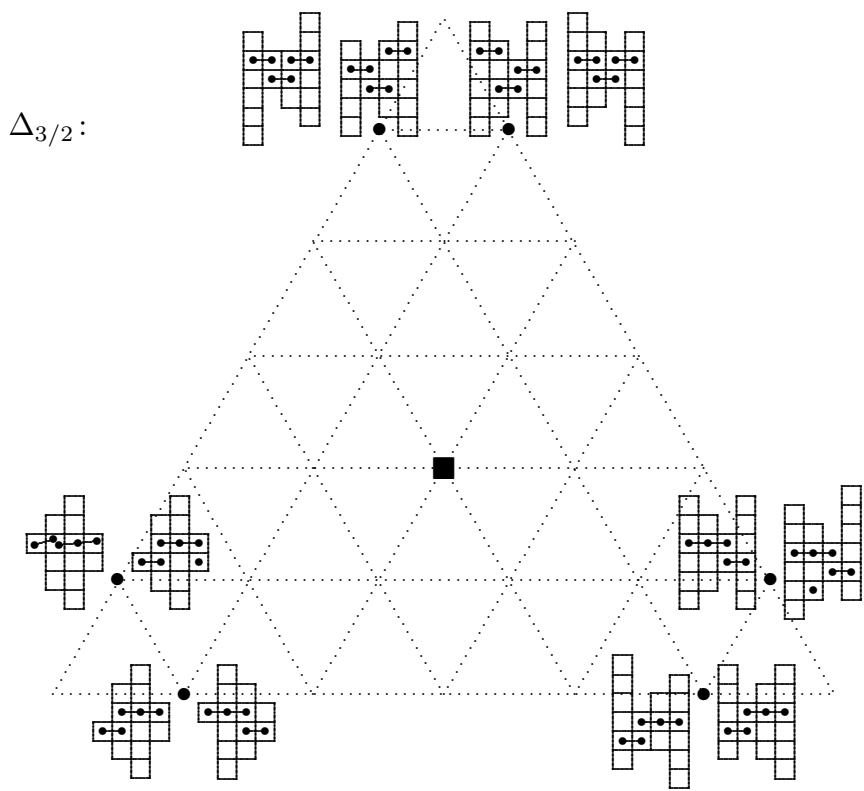
{\beginpicture
   \setcoordinatesystem units <2.165cm,3.75cm>
   \multiput{} at -2.333 .33  2.3333 2.3333 /
   \setdots<1mm>
   \put{$\blacksquare$} at 0 1
   \put{$\Delta_{3/2}\:$} at -2 2 
   \plot -2 .33  2 .33  0 2.33  -2 .33 /
   \plot -1.33 .33  .33 2  -.33 2  1.33 .33 1.67 .67  -1.67 .67 -1.33 .33 /
   \plot -.67 .33  .67 1.67  -.67 1.67  .67 .33  1.33 1  -1.33 1  -.67 .33 /
   \plot 0 .33  -1 1.33  1 1.33  0 .33 /
   \multiput{$\bullet$} at -1.33 .33  .33 2  -.33 2  1.33 .33  1.67 .67  -1.67 .67 /
   \put{\beginpicture
     \setcoordinatesystem units <.21cm,.21cm>
     \multiput{} at 0 0  4 6  /
     \setsolid
     \plot 1 4  4 4  4 2  0 2  0 3  4 3  /
     \plot 3 5  1 5  1 1  3 1 /
     \plot 2 0  2 6  3 6  3 0  2 0 /
     \multiput{$\sss\bullet$} at .5 2.5  1.5 2.5  1.5 3.5  2.5 3.5  3.5 3.5  /
     \plot .5 2.5  1.5 2.5 /
     \plot 1.5 3.5  3.5 3.5 /
     \endpicture} at -1.6 .25
   \put{\beginpicture
     \setcoordinatesystem units <.21cm,.21cm>
     \multiput{} at 0 0  4 6  /
     \setsolid
     \plot 1 0  1 6  2 6  2 0  1 0 /
     \plot 1 1  3 1  3 5  1 5 /
     \plot 1 2  4 2  4 4  0 4  0 3  4 3 /
     \multiput{$\sss\bullet$} at .5 3.5  1.5 3.5  2.5 3.5  2.5 2.5  3.5 2.5 /
     \plot .5 3.5  2.5 3.5 /
     \plot 2.5 2.5  3.5 2.5 /
     \endpicture} at -1.06 .25
   \put{\beginpicture
     \setcoordinatesystem units <.21cm,.21cm>
     \multiput{} at 0 0  4 6  /
     \setsolid
     \plot 2 0  2 6  3 6  3 0  2 0 /
     \plot 3 1  1 1  1 5  3 5 /
     \plot 1 2  4 2  4 4  0 4  0 3  4 3  /
     \multiput{$\sss\bullet$} at .4 3.4  1.4 3.7  1.7 3.4  2.7 3.5  3.7 3.6 /
     \plot .4 3.4  1.4 3.7 /
     \plot 1.7 3.4  3.7 3.6 /
     \endpicture} at -1.94 .75
   \put{\beginpicture
     \setcoordinatesystem units <.21cm,.21cm>
     \multiput{} at 0 0  4 6  /
     \setsolid
     \plot 2 0  2 6  3 6  3 0  2 0 /
     \plot 3 1  1 1  1 5  3 5 /
     \plot 1 4  4 4  4 2  0 2  0 3  4 3 /
     \multiput{$\sss\bullet$} at .5 2.5  1.5 2.5  1.5 3.5  2.5 3.5  3.5 3.5  3.5 2.5 /
     \plot .5 2.5  1.5 2.5 /
     \plot 1.5 3.5  3.5 3.5 /
     \endpicture} at -1.4 .75
   \put{\beginpicture
     \setcoordinatesystem units <.21cm,.21cm>
     \multiput{} at 0 0  4 7  /
     \setsolid
     \plot 0 1  0 7  1 7  1 1  0 1 /
     \plot 3 0  3 6  4 6  4 0  3 0 /
     \plot 4 5  2 5  2 2 /
     \plot 3 1  4 1 /
     \plot 0 2  4 2 /
     \plot 0 3  4 3 /
     \plot 0 4  4 4 /
     \plot 0 5  1 5 /
     \plot 0 6  1 6 /
     \multiput{$\sss\bullet$} at .5 2.5  1.5 2.5  1.5 3.5  2.5 3.5  3.5 3.5 /
     \plot .5 2.5  1.5 2.5 /
     \plot 1.5 3.5  3.5 3.5 /
     \endpicture} at 1.06 .25
   \put{\beginpicture
     \setcoordinatesystem units <.21cm,.21cm>
     \multiput{} at 0 0  4 6  /
     \setsolid
     \plot 0 1  0 6  1 6  1 1  0 1 /
     \plot 3 0  3 6  4 6  4 0  3 0 /
     \plot 4 1  2 1  2 5  4 5 /
     \plot 0 2  4 2 /
     \plot 0 3  4 3 /
     \plot 0 4  4 4 /
     \plot 0 5  1 5 /
     \multiput{$\sss\bullet$} at .5 2.5  1.5 2.5  1.5 3.5  2.5 3.5  3.5 3.5 /
     \plot .5 2.5  1.5 2.5 /
     \plot 1.5 3.5  3.5 3.5 /
     \endpicture} at 1.6 .25
   \put{\beginpicture
     \setcoordinatesystem units <.21cm,.21cm>
     \multiput{} at 0 0  4 7  /
     \setsolid
     \plot 0 0  0 6  1 6  1 0  0 0 /
     \plot 0 1  2 1  2 5  0 5  /
     \plot 3 1  3 7  4 7  4 1  3 1 /
     \plot 0 2  4 2 /
     \plot 0 3  4 3 /
     \plot 0 4  4 4 /
     \plot 3 5  4 5 /
     \plot 3 6  4 6 /
     \multiput{$\sss\bullet$} at 1.5 1.5  .5 3.5  1.5 3.5  2.5 3.5  2.5 2.5  3.5 2.5 /
     \plot .5 3.5  2.5 3.5 /
     \plot 2.5 2.5  3.5 2.5 /
     \endpicture} at 1.94 .75
   \put{\beginpicture
     \setcoordinatesystem units <.21cm,.21cm>
     \multiput{} at 0 0  4 6  /
     \setsolid
     \plot 0 0  0 6  1 6  1 0  0 0  /
     \plot 0 1  2 1  2 5  0 5 /
     \plot 3 0  3 6  4 6  4 0  3 0 /
     \plot 3 1  4 1 /
     \plot 0 2  4 2 /
     \plot 0 3  4 3 /
     \plot 0 4  4 4 /
     \plot 3 5  4 5 /
     \multiput{$\sss\bullet$} at .5 3.5  1.5 3.5  2.5 3.5  2.5 2.5  3.5 2.5 /
     \plot .5 3.5  2.5 3.5 /
     \plot 2.5 2.5  3.5 2.5 /
     \endpicture} at 1.4 .75
   \put{\beginpicture
     \setcoordinatesystem units <.21cm,.21cm>
     \multiput{} at 0 0  4 7  /
     \setsolid
     \plot 0 0  0 6  1 6  1 0  0 0 /
     \plot 3 1  3 7  4 7  4 1  3 1 /
     \plot 4 2  2 2  2 5  /
     \plot 0 1  1 1 /
     \plot 0 2  1 2 /
     \plot 0 3  4 3 /
     \plot 0 4  4 4 /
     \plot 0 5  4 5 /
     \plot 3 6  4 6 /
     \multiput{$\sss\bullet$} at .5 4.5  1.5 4.5  1.5 3.5  2.5 3.5  2.5 4.5 3.5 4.5 /
     \plot .5 4.5  1.5 4.5 /
     \plot 1.5 3.5  2.5 3.5 /
     \plot 2.5 4.5  3.5 4.5 /
     \endpicture} at -.83 2.15
   \put{\beginpicture
     \setcoordinatesystem units <.21cm,.21cm>
     \multiput{} at 0 0  4 6  /
     \setsolid
     \plot 0 0  0 5  1 5  1 0  0 0 /
     \plot 3 0  3 6  4 6  4 0  3 0 /
     \plot 4 1  2 1  2 5  4 5 /
     \plot 0 1  1 1 /
     \plot 0 2  4 2 /
     \plot 0 3  4 3 /
     \plot 0 4  4 4 /
     \multiput{$\sss\bullet$} at .5 3.5  1.5 3.5  1.5 2.5  2.5 2.5  2.5 4.5  3.5 4.5 /
     \plot .5 3.5  1.5 3.5 /
     \plot 1.5 2.5  2.5 2.5 /
     \plot 2.5 4.5  3.5 4.5 /
     \endpicture} at -.33 2.15
   \put{\beginpicture
     \setcoordinatesystem units <.21cm,.21cm>
     \multiput{} at 0 0  4 6  /
     \setsolid
     \plot 0 0  0 6  1 6  1 0  0 0 /
     \plot 3 0  3 6  4 6  4 0  3 0 /
     \plot 0 1  2 1  2 5  0 5 /
     \plot 3 1  4 1 /
     \plot 0 2  4 2 /
     \plot 0 3  4 3 /
     \plot 0 4  4 4 /
     \plot 3 5  4 5 /
     \multiput{$\sss\bullet$} at .5 4.5  1.5 4.5  1.5 2.5  2.5 2.5  2.5 3.5  3.5 3.5 /
     \plot .5 4.5  1.5 4.5 /
     \plot 1.5 2.5  2.5 2.5 /
     \plot 2.5 3.5  3.5 3.5 /
     \endpicture} at .33 2.15
   \put{\beginpicture
     \setcoordinatesystem units <.21cm,.21cm>
     \multiput{} at 0 0  4 7  /
     \setsolid
     \plot 3 0  3 6  4 6  4 0  3 0 /
     \plot 0 1  0 7  1 7  1 1  0 1 /
     \plot 3 1  4 1 /
     \plot 0 2  2 2  2 6  0 6 /
     \plot 3 2  4 2 /
     \plot 0 3  4 3 /
     \plot 0 4  4 4 /
     \plot 0 5  4 5 /
     \multiput{$\sss\bullet$} at .5 4.5  1.5 4.5  1.5 3.5  2.5 3.5  2.5 4.5  3.5 4.5 /
     \plot .5 4.5  1.5 4.5 /
     \plot 1.5 3.5  2.5 3.5 /
     \plot 2.5 4.5  3.5 4.5 /
     \endpicture} at .83 2.15
\endpicture}
$$


\subsection{Six pentapickets.}
\label{app-B-six}

Here, we present six pentapickets, they occur pairwise at the vertices of the triangle
$\Delta_{8/5}$.  Each dimension vector is
the sum of the root for $E_2^2$ or one of its rotations and one of the two radical vectors.

\medskip
$$
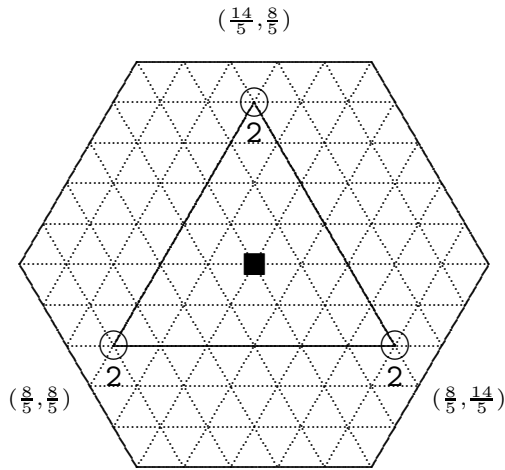
{\beginpicture
   \setcoordinatesystem units <1.299cm, 2.25cm>
   \multiput{} at -2 0  2 2 /
   \plot -2 1  -1 2  1 2  2 1  1 0  -1 0  -2 1 /
   \multiput{$\bigcirc$} at  -1.2 .6  1.2 .6  0 1.8 /
   \plot -1.2 .6  1.2 .6  0 1.8  -1.2 .6 /
     \put{$\blacksquare$} at 0 1
     \multiput{{\tt 2}} at -1.2 .45  1.2 .45  0 1.65 /
   \setdots<.5mm>
   \plot -2 1  2 1 /
   \plot -1 0  1 2 /
   \plot -1 2  1 0 /
   \plot -.6 0  -1.8 1.2  1.8 1.2  .6 0  -1.2 1.8  1.2 1.8  -.6 0 /
   \plot -.2 0  -1.6 1.4  1.6 1.4  .2 0  -1.4 1.6  1.4 1.6  -.2 0 /
   \plot -.6 2  -1.8 .8  1.8 .8  .6 2  -1.2 .2  1.2 .2  -.6 2 /
   \plot -.2 2  -1.6 .6  1.6 .6  .2 2  -1.4 .4  1.4 .4  -.2 2 /
   \put{$\ss (\frac85,\frac85)$} at -1.83 .333 
   \put{$\ss (\frac85,\frac{14}5)$} at 1.83 .333
   \put{$\ss (\frac{14}5,\frac85)$} at 0 2.2
\endpicture}
$$

\medskip
The objects form one $\Sigma_3$-orbit consisting of modules of quasi-length two
in the 6-tubes of index $\frac 14$ (left) or $\frac 41$ (right).

\medskip
$$
{\beginpicture
   \setcoordinatesystem units <1.1547cm,2cm>
   \multiput{} at -3 0  3 3 /
   \setdots<1mm>
   \plot -3 0  3 0  0 3 -3 0 /
   \put{$\blacksquare$} at 0 1
   \put{$\Delta_{8/5}\:$} at -3 2.5 
   \plot -2 0  -2.5 .5  2.5 .5  2 0  -.5 2.5  .5 2.5  -2 0 /
   \plot -1 0  -2 1  2 1  1 0  -1 2  1 2  -1 0 /
   \plot 0 0  -1.5 1.5  1.5 1.5  0 0 /
   \multiput{$\bullet$} at -3 0  3 0  0 3  /
   \put{\beginpicture
     \setcoordinatesystem units <.21cm,.21cm>
     \multiput{} at 0 0  5 6  /
     \setsolid
     \plot 0 1  0 5  1 5  1 1  0 1 /
     \plot 2 0  2 6  3 6  3 0  2 0 /
     \plot 2 1  3 1 /
     \plot 2 5  4 5  4 2  0 2 /
     \plot 0 3  5 3  5 4  0 4 /
     \multiput{$\sss\bullet$} at .4 3.4  1.4 3.7  1.5 2.5  1.7 3.3  2.7 3.4  3.7 3.5  4.7 3.6 /
     \plot .4 3.4   1.4 3.7 /
     \plot 1.7 3.3 4.7 3.6 /
     \endpicture} at -3.6 0 
   \put{\beginpicture
     \setcoordinatesystem units <.21cm,.21cm>
     \multiput{} at 0 0  5 6  /
     \setsolid
     \plot 0 0  0 6  1 6  1 0  0 0 /
     \plot 0 1  2 1  2 4  0 4 /
     \plot 3 1  3 5  4 5  4 1  3 1 /
     \plot 3 4  5 4  5 2  0 2 0 3  5 3 /
     \plot 0 5  1 5 /
     \multiput{$\sss\bullet$} at .5 3.5  1.5 3.5  1.4 2.4 2.4 2.7  2.7 2.3  3.7 2.5  4.7 2.7 /
     \plot .5 3.5  1.5 3.5 /
     \plot 1.4 2.4  2.4 2.7 /
     \plot 2.7 2.3  4.7 2.7 /
     \endpicture} at -2.4 0
   \put{\beginpicture
     \setcoordinatesystem units <.21cm,.21cm>
     \multiput{} at 0 0  5 7  /
     \setsolid
     \plot 0 1  0 6  1 6  1 1  0 1 /
     \plot 2 0  2 6  3 6  3 0  2 0 /
     \plot 4 1  4 7  5 7  5 1  4 1 /
     \plot 2 1  3 1 /
     \plot 0 2  5 2 /
     \plot 0 3  5 3 /
     \plot 0 4  5 4 /
     \plot 0 5  3 5 /
     \plot 4 5  5 5 /
     \plot 4 6  5 6 /
     \multiput{$\sss\bullet$} at .5 2.5  1.5 2.5  1.5 3.5  2.5 3.5  3.5 3.5  3.5 2.5  4.5 2.5 /
     \plot .5 2.5  1.5 2.5 /
     \plot 1.5 3.5  3.5 3.5 /
     \plot 3.5 2.5  4.5 2.5 /
     \endpicture} at 2.4 0
   \put{\beginpicture
     \setcoordinatesystem units <.21cm,.21cm>
     \multiput{} at 0 0  5 6  /
     \setsolid
     \plot 0 0  0 6  1 6  1 0  0 0 /
     \plot 0 1  2 1  2 5  0 5 /
     \plot 3 0  3 6  4 6  4 0  3 0 /
     \plot 3 1  5 1  5 5  3 5 /
     \plot 0 2  5 2 /
     \plot 0 3  5 3 /
     \plot 0 4  5 4 /
     \multiput{$\sss\bullet$} at .5 3.5  1.5 3.5  2.5 3.5  2.5 2.5  3.5 2.5  4.5 2.5  4.5 1.5 /
     \plot .5 3.5  2.5 3.5 /
     \plot 2.5 2.5  4.5 2.5 /
     \endpicture} at 3.6 0 
   \put{\beginpicture
     \setcoordinatesystem units <.21cm,.21cm>
     \multiput{} at 0 0  5 6  /
     \setsolid
     \plot 0 0  0 6  1 6  1 0  0 0 /
     \plot 0 1  4 1  4 5  0 5 /
     \plot 2 0  2 6  3 6  3 0  2 0 /
     \plot 0 2  5 2  5 4  0 4 /
     \plot 0 3  5 3 /
     \multiput{$\sss\bullet$} at .5 3.5  1.5 3.5  1.5 2.5  1.5 4.5  2.5 4.5  3.5 4.5  3.5 3.5  4.5 3.5 /
     \plot .5 3.5  1.5 3.5 /
     \plot 1.5 4.5  3.5 4.5 /
     \plot 3.5 3.5  4.5 3.5 /
     \endpicture} at -.6 3 
   \put{\beginpicture
     \setcoordinatesystem units <.21cm,.21cm>
     \multiput{} at 0 0  5 7  /
     \setsolid
     \plot 2 0  2 6  3 6  3 0  2 0 /
     \plot 0 1  0 6  1 6  1 1  0 1 /
     \plot 4 1  4 7  5 7  5 1  4 1 /
     \plot 2 1  3 1 /
     \plot 4 2  5 2 /
     \plot 0 3  5 3 /
     \plot 0 4  5 4 /
     \plot 0 5  5 5 /
     \plot 0 2  3 2 /
     \plot 4 6  5 6 /
     \multiput{$\sss\bullet$} at 2.5 4.5  3.5 4.5  3.7 3.3  4.5 3.5  .5 4.5  1.5 4.5  1.5 3.5  3.3 3.7 /
     \plot 2.5 4.5  3.5 4.5 /
     \plot 3.7 3.3  4.5 3.5 /
     \plot .5 4.5  1.5 4.5 /
     \plot 1.5 3.5  3.3 3.7 / 
     \endpicture} at .6 3
\endpicture}
$$

\smallskip
The line $q=\frac{16}5$ contains only the two mentioned pentapickets.
This completes our description of the indecomposable objects in $\Cal S(6)$ with
$q\le \frac{13}4$.
$$
{\beginpicture
   \setcoordinatesystem units <1cm,.5cm>
\put{$q$} at -.7 0
\put{$1$} at 1 0
\put{$2$} at 2 0
\put{$\frac52$} at 3 0
\put{$3$} at 4 0
\put{$\frac{16}5$} at 5 0
\put{$\frac{13}4$} at 6 0
\put{$\cdots$} at 7 0
\put{$\#\ind$} at -.7 -1
\put{$2$} at 1 -1
\put{$4$} at 2 -1
\put{$2$} at 3 -1
\put{$15$} at 4 -1
\put{$2$} at 5 -1
\put{$6$} at 6 -1
\put{$\cdots$} at 7 -1
\plot 0.2 0.5  0.2 -1.5 /
\plot -1.4 -.5  7.5 -.5 /
\endpicture}
$$

\vfill\eject

\setcounter{sectionapp}{3}
\setcounter{section}{21}

\section{Applications.}
\label{app-C}

\bigskip
by Markus Schmidmeier

\bigskip
	\bigskip
	\medskip
        Invariant subspaces are ubiquiteous in pure and applied mathematics and in applications.
        We give a few examples as a reminder.

\subsection{Invariant subspaces in undergraduate mathematics.}
\label{app-C-one}

        Invariant subspaces occur frequently in undergraduate mathematics although the student may
        not always be aware of it.

        \smallskip
        Consider for example the method of Undetermined Coefficients
        to deal with a linear higher-order non-homogeneous differential equation with constant coefficients.
        For illustration, in the equation
        \begin{equation*}
        Ly=e^{3t}+te^t, \quad\text{where}\quad Ly=y^{(4)}-5y^{(3)}+9y''-7y'+2y,\tag{DE}
        \end{equation*}
        the characteristic polynomial for the linear operator $L$ has roots 1, 1, 1 and 2
        (so $ L=(\frac d{dt}-1)^3(\frac d{dt}-2)$).  Consider the first term on the right hand side
        of (DE).  Clearly,  $L$ acts
        invertibly on the $\frac d{dt}$-invariant subspace $\langle e^{3t}\rangle$,
        hence a particular solution
        for $Ly=e^{3t}$ has the form $A\,e^{3t}$ for some coefficient
        $A\in\mathbb R$.
        
        Let us now deal with the second term, $te^t$, on the right hand side of (DE).
        For $n\in \mathbb N_0$,
        denote by $\Eig_1(n)$ the $(n+1)$-dimensional vector space
        $\langle t^ne^t,\ldots,t^0e^t\rangle$;
        the operator $T=(\frac d{dt}-1)$
        acts nilpotently with kernel $\Eig_1(0)$ and image $\Eig_1(n-1)$.
        Hence for $n=4$, the second term on the right hand side of the differential equation,
        $te^t$, occurs in the image of $(\frac d{dt}-1)^3$ acting on $V=\Eig_1(4)$.
        In other words, we are dealing with the picket
        $$\textstyle V=\langle t^4e^t,\ldots,e^t\rangle,\quad
        T=(\frac d{dt}-1),\quad U=\langle te^t,e^t\rangle$$
        of type $([2],[5])$ where, 
        by construction,
        the subspace generator, $te^t$, has $T$-height $3$ in $V$ (see Section~\ref{sec-thirteen-two}).
        Since the operator $L$ is the product of $T^3$ with
        the operator $(\frac d{dt}-2)$ (which is acting
        invertibly on $V$ and which commutes with $T$), the 
        Ansatz $L(B_1\, t^4e^t+B_0\, t^3e^t)=te^t$ will lead to a particular solution for $Ly=te^t$.
        
        Using superposition, we see that
        there is a particular solution to (DE) in the form
        $y_p(t)=A\,e^{3t}+B_1\,t^4e^t+B_0\,t^3e^t$.

        The instructor can avoid exposing students to invariant subspaces at this early stage.
        By comparison, the proposed method in the textbook on differential equations,
        \cite[6.4.1]{C},
        involves filling in a 5-column table...

\subsection{Linear time-invariant dynamical systems.}
\label{app-C-two}

In the theory of linear time-invariant control systems, the state-space representation of the system 
consists of (finite-dimensional complex) vector spaces $U'$, $V$, $W'$,
called the input space, the state space and the output space;
a linear operator $A:V\to V$;  and linear maps $B:U'\to V$, $C:V\to W'$.
$$\beginpicture
    \setcoordinatesystem units <1cm,1cm>
    \multiput{} at 0 -.25  2 .8 /
    \put{$U'$} at 0 0
    \put{$V$} at 1 0
    \put{$W'$} at 2 0
    \arr{.3 0}{.7 0}
    \arr{1.3 0}{1.7 0}
    \put{$\ssize B$} at .5 -.25
    \put{$\ssize C$} at 1.5 -.25
    \put{$\ssize A$} at 1.4 .7
    \circulararc 270 degrees from 1.2 .2   center at 1 .4
    \arr{0.78 .23}{0.8 .2}
    \put{(DS)} at -2 0
    \endpicture
$$
The corresponding dynamical system is given by the vector equations
$$\dot x(t) = A\,x(t)+B\,u(t),\qquad w(t)=C\,x(t).$$
In general, the classification of control systems is considered infeasible.
For example, the above quiver has infinite representation type even
if the relation that $A^2=0$ is imposed.

Invariant subspaces occur naturally in control theory, see e.g.\ \cite[Theorem 1.24ff]{KS}.
In particular for $T=A$, the {\it controllable subspace} $U=\sum_{i\geq0}\Im T^i\,B$
of the state space $V$ is $T$-invariant, and the dynamical system (DS) is
said to be {\it completely controllable} provided $U=V$.
Focusing on the state space and its controllable subspace, we are dealing with a
representation
$$\beginpicture
\setcoordinatesystem units <1cm,1cm>
\multiput{} at 0 0  2 .8 /
\put{$U$} at 0 0
\put{$V$} at 2 0
\arr{.3 0}{1.7 0}
\put{$\ssize{\text{incl}}$} at 1 .25
\put{$\ssize T$} at 2.4 .7
\put{$\ssize T|U$} at -.4 .7
\circulararc 270 degrees from 2.2 .2   center at 2 .4
\circulararc 270 degrees from .2 .2   center at 0 .4
\arr{1.78 .23}{1.8 .2}
\arr{-.22 .23}{-.2 .2}
\endpicture
$$
of the quiver $Q$ in Section~\ref{sec-two-three} with commutativity relations,
but with no condition on nilpotency.

The state space $V$ decomposes as the direct sum of its generalized eigen spaces,
and this state space decomposition induces a direct sum decomposition
of the embedding $(U,V)$.  Without loss of generality we may hence assume that
the operator has only one eigen value, and that this eigen value is zero.
Thus, we are dealing with an object $(U,V)$ in $\mathcal S$ and our classification results apply.

\medskip
Let us now consider the special case where the input space $U'$ has dimension one
(and where the map $B:U'\to V$ is the inclusion),
then the dynamical system is called a {\it single input (SI) system.}  In this case,
we are dealing with an embedding $(U,V)$ where the controllable
subspace $U=\sum_{i\geq0} T^i(U')$ is cyclic.
If the operator acts nilpotently
such systems have been considered in Section~\ref{sec-thirteen}.

For example, we have seen in (3) in Section~\ref{sec-thirteen-two} that for fixed
dimension $v=\dim V$ of the state space $V$, the number of such indecomposable systems
is equal to the number of partitions of $v$.  Which of these systems are completely
controllable?  Only one, namely the picket $([v],[v])$, satisfies
the required condition $U=V$.

\subsection{The Algebraic Riccati Equation.}
\label{app-C-three}

        In Optimal Control Theory, the multivariable Linear-Quadratic Regulator problem gives rise to the matrix equation
        \begin{equation*}
        A+BX+XC+XDX=0\tag{ARE}
        \end{equation*}
        in the unknown $n\times n$-matrix $X$.  Solutions to (ARE) correspond to invariant subspaces, as follows.
        Let $\Cal H=\left(\smallmatrix C & D\\ -A & -B\endsmallmatrix\right)$, then the computation
        $$\Cal H\cdot\left(\smallmatrix I\\ X\endsmallmatrix\right)
               =\left(\smallmatrix C+DX\\ -A-BX\endsmallmatrix\right)
               =\left(\smallmatrix C+DX\\ XC+XDX\endsmallmatrix\right)
               =\left(\smallmatrix I\\ X\endsmallmatrix\right)\cdot (C+DX)$$
        shows that the $n$-dimensional subspace $U$ of $k^{2n}$
        spanned by the columns of $\left(\smallmatrix I\\X\endsmallmatrix\right)$
        is $\Cal H$-invariant.  Conversely, any $\Cal H$-invariant $n$-dimensional subspace $U$ of $k^{2n}$
        with the extra property that the projection $U\to k^n, (u_i)_{i=1}^{2n}\mapsto (u_i)_{i=1}^n,$ is onto,
        gives rise to a solution of (ARE).

        \smallskip
        Wikipedia \cite{W2} writes, ``The algebraic Riccati equation
        determines the solution of [...] two of the most fundamental
        problems in control theory.''

\subsection{Complexes in Topological Data Analysis.}
\label{app-C-four}

If a simplicial complex $\Cal X$ in topology is represented algebraically as a
differential complex (that is, a graded
square-zero linear operator $T:X_*\to X_*$), then a simplicial filtration of $\Cal X$
indexed by a poset $\Cal P$ with maximal element $*$
is represented as a system $(X_i,\leq)$ of $T$-invariant subspaces of $X_*$, indexed by $\Cal P$.
The graded subspaces $X_i$ satisfy $X_i\leq X_j$ whenever $i\leq j$ holds in $\Cal P$.

       \medskip
In Topological Data Analysis, the persistence diagram (``bar code'') records the evolution of the
topological space (typically given by the presence of homology) as the filtration
parameter varies over the index set $\Cal P$. In the single parameter case, the poset
$\Cal P$ is totally ordered and the category of systems of $T$-invariant subspaces
has discrete representation type (recall that $T^2=0$).
In fact, the poset representations form the torsionless objects in a module category
and the homology functor provides a correspondence with representations of a linear quiver
\cite{RZ}. 

\medskip
In multi-parameter persistence, the poset $\Cal P$ needs no longer be totally ordered.
Then two complications arise.  The category of poset representations may be of wild
representation type, and indecomposable non-projective representations may no longer be
determined uniquely by homology.
Multi-parameter persistence is topic of recent research, see \cite{O} for an introduction.

\subsection{The invariant subspace problem in analysis.}
\label{app-C-five}
        
        According to \cite{W1},
        the invariant subspace problem in functional analysis is a partially unresolved problem
        asking whether every bounded operator on a complex Banach space
        sends some non-trivial closed subspace
        to itself.
        Many variants of the problem have been solved,
        by restricting the class of bounded operators considered
        or by specifying a particular class of Banach spaces.

        \smallskip
        At this time,
        the problem may still be open for Hilbert spaces, although
        several solutions have been proposed since 2023.

	\bigskip\bigskip
	\bigskip
\renewcommand{\refname}{Additional references for Appendix C}

\vfill\eject

\setcounter{sectionapp}{4}
\setcounter{section}{22}

\section{Covering Theory, a Report.}
\label{app-D}

\bigskip
by Claus Michael Ringel

\bigskip
\centerline{Reprint from the Izmir Lectures (2014), slightly revised}

\bigskip\bigskip

Covering theory provides an important method to construct indecomposable representations
of a quiver with cycles using a suitable covering of the quiver. 
Let us stress that in this appendix we 
deal with quivers and their representations in general: Relations which are satisfied
by some given representations will not play a role. 
The quivers which can be handled, are the orbit quivers 
$Q = \widetilde Q/G$, where 
$G$ is a group of automorphisms of a quiver $\widetilde Q$. 
In the application which we have in mind,
namely $\Cal S(n)$, the group $G$ is just the additive group $\mathbb Z$, and the representations of
$\widetilde Q$ are the graded (that means $\mathbb Z$-graded) $kQ$-modules.  
Parallel to the development of covering theory by Bongartz,
Gabriel and Riedtmann \cite{BG, G, Rm},
the (equivalent) theory of 
dealing with group-graded algebras and the corresponding graded modules
was developed by Gordan and Green \cite{GG}. In the language of covering theory, the forgetful
functor from the category of graded modules to the category of ungraded modules is called
the push-down functor $\pi_\lambda\:\mod k\widetilde Q \to \mod k Q$. The covering
theory was introduced in order to deal with representation-finite algebras. 
It later was extended by Dowbor and Skowro\'nski \cite{DS} to finite-dimensional algebras
in general.
	\medskip
        
From now on, we change the perspective (and thus the notation). The quiver which we will work
with and which is the start of all the considerations will 
be denoted by $Q$ (and not by $\widetilde Q$ as in the previous sentences): We start
with a quiver $Q$ and with a group $G$ of automorphisms of $Q$.

\subsection{Locally fin.-dim.\ representations of locally finite quivers.}
\label{app-D-one}
        
A quiver $Q$ is said to be {\it locally finite} provided any vertex $x$ is head or target
of only finitely many arrows. A representation $M$ of a quiver $Q$ is said to be {\it locally
finite-dimensional} provided all the vector spaces $M_x$ are finite-dimensional. Let us denote
by $\Mod kQ$ the category of locally finite-dimensional representations of $Q$ (and by $\MOD kQ$
the category of all the representations of $Q$). 
	\medskip
        
Given a (not necessarily finite) direct sum $M = \bigoplus_{i\in I} M_i$ with
indecomposable modules $M_i$, 
we say that any indecomposable module occurs (in this direct sum decomposition)
with {\it finite multiplicity}, provided for
any module $N$, the number of indices $i\in I$ such that $M_i$ is isomorphic to $N$, is finite. 
\medskip
        
\begin{theorem}[Dowbor-Skowro\'nski]
  \label{thm-D-one-one}
  Let $Q$ be a locally finite quiver. 
  Any indecomposable locally finite-dimensional representation has local endomorphism
  ring. Any locally finite-dimensional representation
  is the direct sum of indecomposable representations, each one occurring with finite multiplicity.
\end{theorem}

\medskip
\begin{remark}
  \label{rem-D-one-one}
  Recall that the
  theorem of Krull-Remak-Schmidt-Azuyama asserts that for a direct sum $\bigoplus_i M_i$ 
  of modules $M_i$ with {\bf local} endomorphism rings, the number of indices $i\in I$ with 
  $M_i$ isomorphic to a fixed module $N$ does not depend on the decomposition. According to
  Theorem, we can use this result.
\end{remark}

\begin{proof}[Proof of Theorem]
First, let $M$ be an indecomposable locally finite-dimensional representation of $Q$.
We show that for any endomorphism $f = (f_x)_{x\in Q_0}$ of $M$ either all the maps $f_x$
are nilpotent (in this case $f$ is said to be {\it locally nilpotent}) 
or else that all non-zero maps $f_x$ are automorphisms. Recall that given a 
finite-dimensional vector space
$V$ and an endomorphism $\phi$ of $V$, there is a (unique!) direct decomposition $V = V'\oplus V''$
of vector spaces such that $\phi(V')\subseteq V', \phi(V'') \subseteq V''$ so that the 
restriction $\phi' = \phi|V'$ is bijective, whereas the restriction $\phi'' = \phi|V''$ is 
nilpotent. 

Looking at the vector space endomorphism $f_x$ of $M_x$, we obtain in this
way a direct decomposition $M_x = M'_x\oplus M''_x$ such that 
$f_x(M_x')\subseteq M_x'$ and $f_x(M_x'') \subseteq M_x''$, with $f_x' = f_x|M_x'$ bijective, 
and  $f_x'' = f_x|M_x''$ nilpotent. Let $\alpha\:x \to y$ be an arrow of $Q$, thus there is given 
$M_\alpha\:M_x \to M_y.$ With respect to the direct decompositions $M_x = M_x' \oplus M_x''$
we can write  $f_x$ in matrix form 
$\left[\smallmatrix (f_x')^t & 0 \cr 0 & (f''_x)^t \endsmallmatrix\right]$. Similarly, 
we use the direct decompositions $M_x = M_x' \oplus M_x''$
and $M_y = M'_y \oplus M''_y$ in order to write $M_\alpha$ in matrix form $M_\alpha = 
\left[\smallmatrix A & B \cr C & D\endsmallmatrix\right]$. Since $f$ is an endomorphism of $M$,
there is the commutativity condition $M_\alpha f_x = f_y M_\alpha$, thus also 
$M_\alpha f_x^t =
f_y^tM_\alpha$ for all $t\ge 0$. In terms of matrices, this means that
$$
 \left[\smallmatrix A & B \cr C & D\endsmallmatrix\right]
 \left[\smallmatrix (f_x')^t & 0 \cr 0 & (f''_x)^t \endsmallmatrix\right] = 
\left[\smallmatrix (f'_y)^t  & 0 \cr 0 &  (f''_y)^t\endsmallmatrix\right]\left[\smallmatrix A & B \cr C & D\endsmallmatrix\right],
$$
in particular, we have
$$
 B(f''_x)^t = (f'_y)^tB,\quad\text{and}\quad  C(f'_x)^t = (f''_y)^tC,
$$
 for all $t\ge 0$. Since $f''_x$ and $f''_y$ are nilpotent,
 the maps $B(f''_x)^t$ and $(f''_y)^tC$ are zero
for $t$ large. Since $f'_x$ and $f'_y$ are invertible, we conclude that $B = 0,\ C = 0.$

This shows that we have vector space decompositions $M_x = M'_x \oplus M''_x$ such that
for any arrow $\alpha\:x\to y$, the map $M_\alpha$ maps $M'_x$ into $M'_y$ and $M''_x$ into $M''_y,$
thus $M = M'\oplus M''$ is a direct decomposition of $kQ$-modules. By assumption $M$ is indecomposable,
thus either $M = M'$ or $M = M''$. In the first case, all non-zero maps $f_x$ are automorphisms,
in the second case, $f$ is locally nilpotent. 

It follows as usual that the set of locally nilpotent endomorphisms of $M$ is an ideal in 
the endomorphism ring $\End(M)$ of $M$ and that this ideal is 
the unique maximal ideal of $\End(M)$, thus $\End(M)$ is a local ring.
	\medskip
        
Now let us consider arbitrary locally finite-dimensional representations $M$ of $Q$,
where $Q$ is a locally finite quiver. Without loss of generality, we may assume that $Q$
is connected, thus clearly $Q_0$ is a countable set. Since the assertion of the theorem
is well-known for finite quivers, we assume that $Q$ is infinite, thus we can label the
vertices as an infinite sequence $x(1),x(2),\dots$ and we put
$\Cal X(t) = \{x(1),x(2),\dots,x(t)\}.$

If $\Cal X$ is a set of vertices of $Q$, we say that $M\in \Mod kQ$ 
is {\it $\Cal X$-indecomposable} provided
for any direct decomposition $M = M' \oplus M''$, we have $M'_x = 0$ for all $x \in \Cal X$
or $M''_x = 0$ for all $x\in \Cal X$. Note that we do not require that $M_x \neq 0$, not even that
$M\neq 0$, thus a direct summand of an $\Cal X$-indecomposable module is $\Cal X$-indecomposable. 

A finite direct decomposition $M = \bigoplus_{i=1}^m M(i)$ is called an {\it $\Cal X$-decomposition}
provided all the $M(i)$ are $\Cal X$-indecomposable. 
\medskip

\begin{lemma}
  \label{lem-D-one}
  If $\Cal X$ is a finite set, then any
  $M\in \Mod kQ$ has an $\Cal X$-decomposition.
\end{lemma}

\medskip
Note that the direct summands $M(i)$ in an
$\Cal X$-decomposition are usually not unique, not even up to isomorphism. 
\smallskip
        
\begin{proof}[Proof of the Lemma]
  Either $M$ is already $\Cal X$-indecomposable, then
  nothing has to be done. Otherwise, there is a direct decomposition
  $M = M'\oplus M''$ with $M'_x \neq 0$ for some $x\in \Cal X$ and $M''_y \neq 0$
  for some $y\in \Cal X.$ Next we consider $M'$ and $M''$ separately, and so on.
  Why does this process stops? Let $d(M) = \sum_{x\in \Cal X} \dim M_x,$
  this is a finite number (since $\Cal X$ is finite and $M$ is locally
  finite-dimensional). We have $d(M) = d(M')+d(M'')$ and by assumption,
  $d(M') \neq 0,$ and $d(M'') \neq 0,$ thus $d(M') < d(M)$ and $d(M'') < d(M).$
  Thus, we see that after a finite number of steps, the process has to
  stop: This means that the corresponding direct summands are
  $\Cal X$-indecomposable.
\end{proof}

\bigskip 
Now we start with a module $M\in \Mod kQ$ and want to decompose it. This is done inductively, 
looking at the vertices $x(1), x(2),\dots .$  
The direct summands $M(v)$ of $M$ obtained in step $t$ will be indexed by 
a set $I(t)$ of sequences $v = (v_1,\dots,v_t)$ of natural numbers. In step 1, take an $\Cal X(1)$-decomposition
$$
 M = \bigoplus_{v_1 = 1}^{m(1)} M(v_1) = \bigoplus_{I(1)} M(v)
$$
where $I(1)$ is the set of numbers $1\le v_1 \le m(1)$. 
Assume we have constructed already 
an $\Cal X(t)$-decomposition $M = \bigoplus_{v} M(v_1,\dots,v_t)$, then we take for any module
$M(v_1,\dots,v_t)$ an 
$\{x(t+1)\}$-decomposition 
$$
 M(v_1,\dots,v_t) = \bigoplus_{v_{t+1}= 1}^{m(v_1,\dots,v_t)} M(v_1,\dots,v_t,v_{t+1}).
$$
Of course, since $M(v_1,\dots,v_t)$ is $\Cal X(t)$-indecomposable, 
all the modules $M(v_1,\dots,v_t,v_{t+1})$ are $\Cal X(t+1)$-indecomposable.
Thus, we obtain in this
way an $\Cal X(t+1)$-decomposition 
$$
 M =  
  \bigoplus_{v\in I(t+1)} M(v)
$$
with $I(t+1)$ the set of sequences
$(v_1,\dots,v_{t+1})$ such that 
$(v_1,\dots,v_{t})$ belongs to $I(t)$ and $1 \le v_{t+1} \le m(v_1,\dots,v_t).$

Let $I$ be the set of infinite sequences $v =(v_1,v_2,\dots)$ such that 
$(v_1,\dots,v_{t})$ belongs to $I(t)$, for all $t$. For $v\in I$, we define the module 
$M(v)$ as $M(v) = \bigcap_t M(v_1,\dots,v_t).$ Note that the restriction of $M(v)$
to the full subquiver with vertices in $\Cal X(t)$ is equal to the restriction of $M(v_1,\dots,v_t)$
to this subquiver. This shows that $M(v)$ is $\Cal X(t)$-indecomposable, for all $t$ and
that $M = \bigoplus_I M(v).$ Since $M(v)$ is $\Cal X(t)$-indecomposable for all $t$, we see that
$M(v)$ is either zero or else indecomposable. If we denote by $I'$ the set of indices $v$
such that $M(v) \neq 0$, then $M = \bigoplus_{I'} M(v)$ is a direct decomposition
with indecomposable direct summands. 

It remains to stress that multiplicities have to be finite: If $M, N$ are locally finite-dimensional
$kQ$-modules and  $N_x \neq 0$ for some vertex $x$, then
in any direct decomposition of $M$, the number of direct summands 
which are isomorphic to $N$ is bounded by $\dim M_x.$
\end{proof}

\medskip
\begin{remark}
  \label{rem-D-one-two}
  In the proof given above, we had to single out at the end the indices $v\in I$
  with $M(v) = 0,$ thus replacing the index set $I$ by $I'$. 
  One may wonder whether one can avoid this.
  Given a module $M$ in $\Mod kQ$, we have used $\{x\}$-decompositions $M = \bigoplus M^{(i)}$.
  Of course, in case $M_x = 0$, we may require to use as decomposition just the trivial
  one $M = M$, and for 
  $M_x \neq 0$, we may require that $M^{(i)}_x \neq 0$
  for all $i$. In general, this will reduce the size of $I$, but still some of the 
  summands $M(v)$ with $v\in I$ may be zero.

  As an example, consider the quiver of type $\mathbb A_\infty$ with arrows
  from right to left, and the following representation $M$:
  $$\beginpicture
  \setcoordinatesystem units <1cm,.7cm>
  \multiput{} at 0 0  6 1 /
  \multiput{$k^2$} at 0 0  2 0  4 0 /
  \multiput{$\left[\smallmatrix 0 & 0 \\ 1 & 0 \endsmallmatrix\right]$} at 1 .5  3 .5 /
  \arr{1.5 0}{.5 0}
  \arr{3.5 0}{2.5 0}
  \arr{5.5 0}{4.5 0}
  \put{$\cdots$} at 6 0 
  \endpicture
  $$
  In the first decomposition $M = M(1)\oplus M(2)$, we may assume that 
  one of the direct summands, say $M(1)$,
  is simple projective. In the second step, we can decompose $M(2) = M(2,1) \oplus
  M(2,2)$ with $M(2,1)$ of length 2, then, in the next step, $M(2,2) = M(2,2,1) \oplus M(2,2,2)$ with
  $M(2,2,1)$ of length 2, and so on. In this way, we obtain as index set $I'$
  the set of sequences of the form $(2,\dots,2,1,1,\dots)$ starting with $s\ge 0$ entries equal
  to $2$, all others equal to $1$, with $M(1,1,\dots)$ the simple projective module,
  all other modules $M(2,\dots,2,1,1,\dots)$ of length $2$. But the set $I$ contains in
  addition the constant sequence $(2,2,\dots)$ with $M(2,2,\dots) = 0.$
\end{remark}

\medskip
        
\begin{remark}
  \label{rem-D-one-three}
  We have seen that given an indecomposable
  locally finite-dimensional representation $M$ of a locally
  finite quiver, any non-invertible endomorphism $f$ of $M$ is locally nilpotent. But note
  that {\it $f$ does not have to be nilpotent.} Here is an example:
  We consider the quiver $Q$ presented on the left, and the representation $M$ of $Q$ outlined
  on the right:
  $$
  
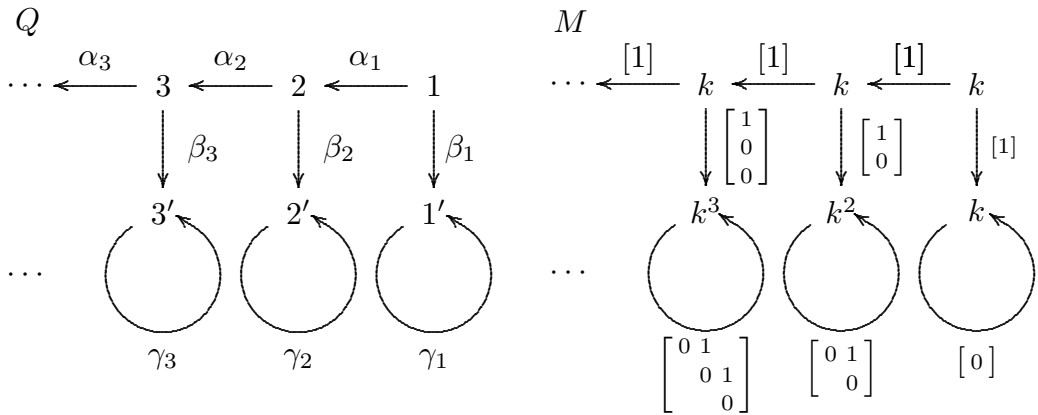
{\beginpicture
    \setcoordinatesystem units <1.5cm,1.4cm>
    \put{\beginpicture
      \multiput{} at 0 0.5  3 -3 /
      \put{$Q$} at 0 .5
      \put{$1$} at 3 0
      \put{$2$} at 2 0
      \put{$3$} at 1 0
      \put{$1'$} at 3 -1
      \put{$2'$} at 2 -1
      \put{$3'$} at 1 -1
      \multiput{$\cdots$} at 0 0  0 -1.5 /
      \arr{2.8 0}{2.2 0}
      \arr{1.8 0}{1.2 0}
      \arr{0.8 0}{0.2 0}
      \arr{3 -.2}{3 -.8}
      \arr{2 -.2}{2 -.8}
      \arr{1 -.2}{1 -.8}
      \circulararc -300 degrees from 3.2 -1.1 center at 3 -1.5
      \arr{3.2 -1.1}{3.1 -1.03}
      \circulararc -300 degrees from 2.2 -1.1 center at 2 -1.5
      \arr{2.2 -1.1}{2.1 -1.03}
      \circulararc -300 degrees from 1.2 -1.1 center at 1 -1.5
      \arr{1.2 -1.1}{1.1 -1.03}
      \put{$\alpha_1$} at 2.5 0.2
      \put{$\alpha_2$} at 1.5 0.2
      \put{$\alpha_3$} at 0.5 0.2
      \put{$\beta_1$} at 3.2 -.5
      \put{$\beta_2$} at 2.3 -.5
      \put{$\beta_3$} at 1.3 -.5
      
      \put{$\gamma_1$} at 3 -2.17
      \put{$\gamma_2$} at 2 -2.17
      \put{$\gamma_3$} at 1 -2.17
      \endpicture} at 0 0 
    \put{\beginpicture
      \multiput{} at 0 0.5  3 -3 /
      \put{$M$} at 0 .5
      \put{$k$} at 3 0
      \put{$k$} at 2 0
      \put{$k$} at 1 0
      \put{$k$} at 3 -1
      \put{$k^2$} at 2 -1
      \put{$k^3$} at 1 -1
      \put{$[1]$} at 2.5 0.2
      \put{$[1]$} at 1.5 0.2
      \put{$[1]$} at 0.5 0.2
      \put{$[1]$} at 2.5 0.2
      \put{$\ss[1]$} at 3.2 -.5
      \put{$\left[\smallmatrix 1 \cr 0 \endsmallmatrix \right]$} at 2.3 -.5
      \put{$\left[\smallmatrix 1 \cr 0\cr 0 \endsmallmatrix \right]$} at 1.3 -.5
      
      \put{$\left[\smallmatrix 0 \endsmallmatrix \right]$} at 3 -2.2
      \put{$\left[\smallmatrix 0 &1 \cr & 0  \endsmallmatrix \right]$} at 2 -2.25
      \put{$\left[\smallmatrix 0 &1\cr & 0 & 1\cr & & 0 \endsmallmatrix \right]$} at 1 -2.3
      \multiput{$\cdots$} at 0 0  0 -1.5 /
      \arr{2.8 0}{2.2 0}
      \arr{1.8 0}{1.2 0}
      \arr{0.8 0}{0.2 0}
      \arr{3 -.2}{3 -.8}
      \arr{2 -.2}{2 -.8}
      \arr{1 -.2}{1 -.8}
      \circulararc -300 degrees from 3.2 -1.1 center at 3 -1.5
      \arr{3.2 -1.1}{3.1 -1.03}
      \circulararc -300 degrees from 2.2 -1.1 center at 2 -1.5
      \arr{2.2 -1.1}{2.1 -1.03}
      \circulararc -300 degrees from 1.2 -1.1 center at 1 -1.5
      \arr{1.2 -1.1}{1.1 -1.03}
      \endpicture} at 4 0 
    \endpicture}
  $$
  The endomorphism $f$ of $M$ 
  which we are interested in, is defined by $f_i = 0$ and $f_{i'}$
  being the nilpotent $i\times i$ Jordan matrix (as $M_{\gamma_i}$)
  for all $i\in \mathbb N_1.$ Of course, $f$ is locally 
  nilpotent, but not nilpotent. 
\end{remark}

\subsection{Automorphisms of a quiver $Q$ which operate freely on $Q_0$.}
\label{app-D-two}
        
Let $Q$ be a quiver and $g$ an automorphism of $Q$. 
We say that $g$ {\it operates freely} on $Q_0$ provided given a vertex $x$ of $Q$ and a
natural number $s$ with $g^sx = x$, we have $g^s = 1.$ Of course, then $G$ acts also freely on the
arrow set $Q_1$ (namely, if $\alpha\:x \to y$ is an arrow with $g(\alpha) = \alpha$, then $g(x) = x$,
thus $g = 1$).

If $M$ is a representation of $Q$, there is a representation $M^g$ defined as follows:
$(M^g)_x = M_{gx}$ for any vertex $x$, and $(M^g)_\alpha = M_{g\alpha}$ for any arrow $\alpha$.
\medskip
        
\begin{lemma}
  \label{lem-D-two}
  Assume that $g$ is an automorphism of $Q$ which acts freely on $Q_0$ 
  and has infinite order. If $M\neq 0$
  is a finite-dimensional representation of $Q$,
  then the modules $M$ and $M^g$ have different support.
\end{lemma}

\medskip
In particular, the module $M^g$ is not isomorphic to $M$.
\smallskip
        
\begin{proof}
  Let $x$ be a vertex of $Q$ which belongs to the support of $M$. For any $i\in \mathbb Z$, let
  $x_i = g^ix.$ Since $g$ operates freely on $Q_0$, we obtain in this way  infinitely many 
  pairwise different elements of $Q_0$. Since $M$ is finite-dimensional, there is $i$ maximal 
  such that $x_i$ belongs to the support of $M$. Then 
  $M^g_{x_i} = M_{gx_i} = M_{x_{i+1}} = 0$ shows that $x_i$ does not belong to the support of $M^g.$
  Thus, $M$ and $M^g$ do not have the same support.
\end{proof}

\medskip
\begin{example}
  Let us present 
  a typical example of a quiver $Q$ with an automorphism $g$ of finite order $p \ge 2$
  acting freely on $Q_0$ and an indecomposable representation $M$ of $Q$ such that $M^g = M$
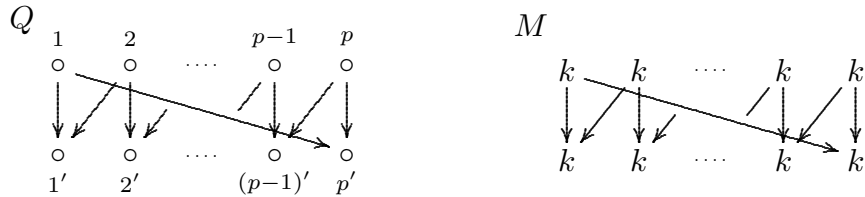
$$
{\beginpicture
\setcoordinatesystem units <.8cm,1cm>
\put{\beginpicture
\multiput{$\circ$} at 0 0  0 1  1 0  1 1  3 0  3 1  4 0  4 1 /
\arr{0 0.8}{0 0.2}
\arr{1 0.8}{1 0.2}
\arr{3 0.8}{3 0.2}
\arr{4 0.8}{4 0.2}
\arr{0.8 0.8}{0.2 0.2}
\arr{3.8 0.8}{3.2 0.2}
\arr{1.5 0.5}{1.2 0.2}
\plot 2.8 0.8  2.5 0.5 /
\arr{0.25 0.9}{3.75 0.1}
\setdots <1mm>
\plot 1.8 0  2.2 0 /
\plot 1.8 1  2.2 1 /
\put{$\ssize 1$} at 0 1.3 
\put{$\ssize 1'$} at 0 -.3 
\put{$\ssize 2$} at 1 1.3 
\put{$\ssize 2'$} at 1 -.3 
\put{$\ssize p-1$} at 3 1.3 
\put{$\ssize (p-1)'$} at 3 -.3 
\put{$\ssize p$} at 4 1.3 
\put{$\ssize p'$} at 4 -.3 
\put{$Q$} at -.5 1.5 
\endpicture} at 0 0
\put{\beginpicture
\multiput{$k$} at 0 0  0 1  1 0  1 1  3 0  3 1  4 0  4 1 /
\arr{0 0.8}{0 0.2}
\arr{1 0.8}{1 0.2}
\arr{3 0.8}{3 0.2}
\arr{4 0.8}{4 0.2}
\arr{0.8 0.8}{0.2 0.2}
\arr{3.8 0.8}{3.2 0.2}
\arr{1.5 0.5}{1.2 0.2}
\plot 2.8 0.8  2.5 0.5 /
\arr{0.25 0.9}{3.75 0.1}
\setdots <1mm>
\plot 1.8 0  2.2 0 /
\plot 1.8 1  2.2 1 /
\put{$M$} at -.5 1.5 
\multiput{} at 0 1.3  0 -.3 /
\endpicture} at 7 0
\endpicture}
$$
Here, $M$ is defined by 
$M_x = k$ for all vertices $x$ and $M_\alpha = 1_k$ for all arrows $\alpha$.
The automorphism $g$ of $Q$ is defined by $g(i) = i\!+\!1$ and $g(i') = (i\!+\!1)'$
(modulo $p$). 
\end{example}

\subsection{Groups operating freely on a quiver.}
\label{app-D-three}
        
Assume now that there is given an automorphism group $G$ of a quiver $Q$ which operates freely on 
$Q_0.$
We denote by $Q/G$ the orbit quiver: If $x$ is a vertex of $Q$, let $\pi(x)$ be the
$G$-orbit of $x$, if $\alpha\:x \to y$ is an arrow, let $\pi(\alpha)\:\pi(x) \to \pi(y)$ be the 
$G$-orbit of $\alpha$, thus $\pi\:Q \to Q/G$ is a morphism of quivers. 

Our interest concerns the functor $\pi_\lambda\:\mod kQ \to \mod k(Q/G)$ which is
defined as follows: If $N$ is a finite-dimensional
representation of $Q$, then 
$$
 (\pi_\lambda N)_z = \bigoplus_{x\in \pi^{-1}(z)} N_x,\quad
 (\pi_\lambda N)_\gamma = \bigoplus_{\alpha\in \pi^{-1}(\gamma)} N_\alpha,
$$
for all vertices $z$ and all arrows $\gamma$ of $Q$;
this functor $\pi_\lambda$ is usually called the {\it push-down functor.}

\medskip
\begin{theorem}[Gabriel]
  Let $Q$ be a locally finite quiver and $G$ a group of automorphisms of $Q$
  which acts freely on $Q_0$. Let $N$ be a finite-dimensional indecomposable representation of $Q$
  such that the representations $N^g$ for $g\in G$ are pairwise non-isomorphic.
  Then $\pi_\lambda N$ is an indecomposable representation of $Q/G$. If $N'$ is 
  a finite-dimensional indecomposable representation of $Q$ 
  such that the representations $\pi_\lambda N$ and $\pi_\lambda N'$ are isomorphic, then there is
  $g\in G$ such that $N'$ is isomorphic to $N^g$.
\end{theorem}

\begin{proof}
We need a further functor, namely  
$\pi_.\:\mod k(Q/G) \to \Mod kQ$ which is defined as follows: If $M$ is
a finite-dimensional $k(Q/G)$-module, then
$$
 (\pi_.M)_x = M_{\pi(x)},\quad (\pi_.M)_\alpha = M_{\pi(\alpha)},
$$
for all vertices $x$ and all arrows $\alpha$ of $Q$; this functor is called the {\it pull-up
functor} (note that by definition $\pi_.M$ is locally finite-dimensional, thus in $\Mod kQ$).
Of course, for any element $g\in G$, we have $(\pi_.M)^g = \pi_.M.$
	\medskip
        
Given a finite-dimensional $kQ$-module $N$, let us consider $\pi_.\pi_\lambda N.$ For $x$
a vertex of $Q$, we have
$$
 (\pi_.\pi_\lambda N)_x = (\pi_\lambda N)_{\pi(x)} = \bigoplus_{y\in \pi^{-1}\pi(x)} N_y.
$$
The set $\pi^{-1}\pi(x)$ is by definition the $G$-orbit of $x$, thus 
$$
 (\pi_.\pi_\lambda N)_x = \bigoplus\nolimits_{g\in G}N_{g(x)} = \bigoplus\nolimits_{g\in G}N^g_{x} =
 \Bigl(\bigoplus\nolimits_{g\in G}N^g\Bigr)_{x}.
$$
It follows that 
$$
 \pi_.\pi_\lambda N = \bigoplus\nolimits_{g\in G}N^g.
$$
Since by assumption the modules $N^g$ are pairwise non-isomorphic (and indecomposable),
we have obtained in this way a direct decomposition of $\pi_.\pi_\lambda N$ 
using pairwise non-isomorphic modules with local endomorphism rings. 

Let us show now that $\pi_\lambda N$ is indecomposable. Thus, assume there is given a direct
decomposition $\pi_\lambda N = M \oplus M'$ of $k(Q/G)$-modules.
Then
$$
  \bigoplus\nolimits_{g\in G}N^g =  \pi_.\pi_\lambda N = \pi_. M \oplus \pi_.M'.
$$
According to D.1 and the Krull-Remak-Schmidt-Azumaya theorem, there is a 
subset $H \subseteq G$ such that $\pi_.M$ is isomorphic to 
$\bigoplus\nolimits_{h\in H}N^h.$ If $M\neq 0$, then $H$ is not empty.
Since $(\pi_.M)^g = \pi_.M$ for all $g \in G$, it follows that $\pi_.M$ is also
isomorphic to $\bigoplus\nolimits_{h\in H}N^{hg},$ for all $g\in G$, thus there is
a direct decomposition of $\pi_.M$ into indecomposable modules such that one of the summands 
is isomorphic to $N$.   
Assume now that both $M, M'$ are non-zero. Then there is
a direct decomposition of $\pi_.M\oplus \pi_.M' = \pi_.\pi_\lambda N$ 
into indecomposable modules such that two of the summands 
are isomorphic to $N$.  This contradicts 
the theorem of Krull-Remak-Schmidt-Azuyama. 

Finally, assume that $N, N'$ are finite-dimensional 
indecomposable representations of $Q$ such that 
the representations $\pi_\lambda N$ and $\pi_\lambda N'$ are isomorphic.
Then $\pi_\lambda N$ is isomorphic to $\bigoplus\nolimits_{g\in G}N^g$
as well as to $\bigoplus\nolimits_{g\in G}(N')^g$. Again using the assumption that
the modules $N^g$ are pairwise non-isomorphic (and indecomposable) and that $N'$
is indecomposable, the theorem of Krull-Remak-Schmidt-Azumaya implies that 
$N'$ is isomorphic to some $N^g.$ 
This completes the proof.
\end{proof}

\medskip
\begin{corollary}
  \label{cor-D-three}
  Let $Q$ be a locally finite quiver and $G$ a torsionfree
  group of automorphisms of $Q$ which acts freely on $Q_0$. Then $\pi_\lambda$
  provides an injective map from the set of $G$-orbits of isomorphism classes of indecomposable
  $kQ$-modules to the set of isomorphism classes of indecomposable $k(Q/G)$-modules.
\end{corollary}

\medskip
        
Let us rephrase what it means that the map given by $\pi_\lambda$ is injective: 
\begin{itemize}[leftmargin=3em]
\item[(a)] If $M$ is a finite-dimensional indecomposable
representation of $Q$, then $\pi_\lambda M$ is an indecomposable representation of $Q/G$. 
\item[(b)] If $M, M'$ are 
finite-dimensional indecomposable representations of $Q$, with
$\pi_\lambda M, \pi_\lambda M'$ being
isomorphic representations of $Q/G$, then 
there is $g\in G$ such that $M'$ is isomorphic to $M^g$.
\end{itemize}
\medskip
        
Corollary is a direct consequence of the theorem using Lemma~\ref{lem-D-two}
in Appendix~\ref{app-D-two}.

\medskip
\begin{remark}
  The map given by $\pi_\lambda$ is injective, but usually not surjective.
  A typical example is provided by the quiver $Q$ of type 
  $\mathbb A_\infty^\infty$ with bipartite orientation
$$
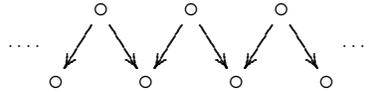
{\beginpicture
\setcoordinatesystem units <.5cm,.8cm>
\multiput{$\circ$} at 0 0  1 1  2 0  3 1  4 0  5 1  6 0  /
\arr{0.8 0.8}{0.2 0.2}
\arr{1.2 0.8}{1.8 0.2}
\arr{2.8 0.8}{2.2 0.2}
\arr{3.2 0.8}{3.8 0.2}
\arr{4.8 0.8}{4.2 0.2}
\arr{5.2 0.8}{5.8 0.2}
\setdots <1mm>
\plot -1 .5  -.4 0.5 /
\plot 6.4 .5  7 0.5 /
\endpicture}
$$
and the shift automorphism $g$ so that there are precisely two $g$-orbits, the
sources and the sinks. Then $Q/G$ is the Kronecker quiver. Obviously, the 
indecomposable representation $M$ of $Q/G$
with $M_x = k$ for both vertices and $M_\alpha = 1_k$
for both arrows is not isomorphic to a representation of the form $\pi_\lambda N$.
\end{remark}

\bigskip\bigskip 
\renewcommand{\refname}{Additional references for the Izmir Lectures, thus Appendix D}

\bigskip\bigskip
Address of the authors:

\bigskip
{\baselineskip=1pt
\rmk
C.\ M.\ Ringel\phantom(\par
Fakult\"at f\"ur Mathematik, Universit\"at Bielefeld \phantom(\par
POBox 100131, D-33501 Bielefeld, Germany  \phantom(\par
ringel@math.uni-bielefeld.de\phantom(\par
ORCID: 0000-0001-7434-0436
}

\bigskip
{\baselineskip=1pt
\rmk
M.\ Schmidmeier\phantom(\par
Department for Mathematics and Statistics, Florida Atlantic University \phantom(\par
777 Glades Road, Boca Raton, Florida 33431, United States of America \phantom(\par
markusschmidmeier@gmail.com\phantom(\par
ORCID: 0000-0003-3365-6666
}

\begin{thebibliography}{9999}
\addcontentsline{toc}{section}{References.}
\bibitem{ARS} M.\ Auslander, I.\ Reiten, S.O.\ Smal\o,
  {\it Representation Theory of Artin Algebras,}
   Cambridge studies in advanced mathematics {\bf 36}, Cambridge 1995.
\bibitem{AS} M.\ Auslander, S.\ O.\ Smal\o.
               {\it Almost split sequences in subcategories,} J.\ Algebra {\bf 69} (1981), 426--454.
\bibitem{DLS} P.\ Dowbor, H.\ Lenzing, A.\ Skowronski.
             {\it Galois coverings by algebras of locally
             support-finite categories,}
             in: Representation Theory I, Proc.\ Conf.\ Ottawa (1984), 91--93.
\bibitem{R1} C.\ M.\ Ringel, {\it Representations of K-species and bimodules,}
           J.\ Algebra {\bf 41} (1976), 269--302.
\bibitem{R2} C.\ M.\ Ringel, {\it Tame Algebras and Integral Quadratic Forms,}
           Springer Lect.\ Notes Math.\ {\bf 1099} (1984).
\bibitem{R3} C.\ M.\ Ringel, {\it Covering Theory,}
           Izmir Lectures (2014), 1--8. Reprinted (slightly revised) as Appendix~\ref{app-D} below. 
\bibitem{RS1} C.\ M.\ Ringel, M.\ Schmidmeier,
             {\it Invariant subspaces of nilpotent linear operators.\ I,}
             J.\ reine angew.\ Math.\ {\bf 614} (2008), 1--52.
\bibitem{RS2} C.\ M.\ Ringel, M.\ Schmidmeier,
             {\it The Auslander-Reiten translation in submodule categories,}
             Trans.\ Amer.\ Math.\ Soc.\ {\bf 360} (2008), 691--716.
\bibitem{S1} M.\ Schmidmeier, {\it Bounded submodules of modules,}
            J. Pure Appl.\ Alg.\ {\bf 203} (2005), 45--82.
\bibitem{S2} M.\ Schmidmeier,
            {\it Invariant subspaces of nilpotent linear operators.\ II,}
            in preparation.
\bibitem{S3} M.\ Schmidmeier, {\it Systems of submodules and an isomorphism
            problem for Auslander-Reiten quivers,}
  Bull.\ Belg.\ Math.\ Soc.\ Simon Stevin {\bf 15} (2008), 523--546.
\bibitem{Str} A.\ Straub, {\it Core partitions into distinct parts and an analog of Euler's theorem,}
            European J.\ Combin.\  {\bf 57} (2016), 40--49.
\end{thebibliography}

\begin{thebibliography}{9999}
\setcounter{enumiv}{12}
\bibitem{Bi} Birep, {\it The ADE chain workshop,} \\
   {\tt\small https://www.math.uni-bielefeld.de/$\sim$sek/ADE/,} aufgerufen am 11.\ 05.\ 2024.
\bibitem{Bh} G.\ Birkhoff, {\it Subgroups of abelian groups,}
           Proc.\ Lond.\ Math.\ Soc.\ (II) {\bf 38} (1934), 385--401.
\bibitem{Bob} G.\ Bobinski,
   {\it Algebras with irreducible module varieties III: Birkhoff varieties,}
   International Mathematics Research Notices,
   Vol.\ 2021, No.\ 4, pp.\ 2497--2525.
\bibitem{GKKP} N.\ Gao, J.\ K\"ulshammer, S.\ Kvamme, Ch.\ Psaroudakis,
   {\it A functorial approach to monomorphism categories II: Indecomposables,}
   (2023), 46pp., arXiv:2303.07753
\bibitem{Kap} I.\ Kaplansky, {\it Infinite Abelian Groups,}
           Ann Arbor, 1954.
\bibitem{Kva} S.\ Kvamme, {\it An introduction to monomorphism categories,}
             manuscript (2023), 42pp.
\bibitem{KSS} J.\ Kosakowska, M.\ Schmidmeier, M.\ Schreiner,
           {\it Abelian $p$-groups with a p-bounded factor or a $p$-bounded subgroup,} 
        arXiv:2312.01451.
\bibitem{P}  H.\ Pr\"ufer, {\it Untersuchungen \"uber die Zerlegbarkeit
            der abz\"ahlbaren prim\"aren Abelschen Gruppen,}
            Math.\ Z.\ {\bf 17} (1923), 35--61.
\bibitem{Z} Interview mit G.\ M.\ Ziegler, in {\it Ein Moment f\"ur Mensch und 
   Mathematik,} Hrg. C.\ Cederbaum, Ph.\ von Homeyer. Freiburg im Breisgau (2007),
   102--111.
\end{thebibliography}

\begin{thebibliography}{9999}
\setcounter{enumiv}{21}
\bibitem{W} Wikipedia article on $\mathbf E_8$ (2023-10-15), {\ttk https://en.wikipedia.org/wiki/E8\_(mathematics)}
\end{thebibliography}

\begin{thebibliography}{9999}
\setcounter{enumiv}{22}
\bibitem{C} C.\ Constanda, {\it Differential equations,} 2nd ed., Springer 2017.      
\bibitem{KS} H.\ Kwakernaak, R.\ Sivan, {\it Linear Optimal Control Systems,}
            Wiley-Interscience [John Wiley \& Sons], New York-London-Sydney, 1972.
\bibitem{O}  S.\ Y.\ Oudot, {\it Persistence theory: From quiver representations
            to data analysis,} Mathematical Surveys and Monographs {\bf 209,}
            AMS 2015.
\bibitem{RZ} C.\ M.\ Ringel and P.\ Zhang,
            {\it Representations of quivers over the algebra of dual numbers,}
            J.\ Algebra {\bf 475} (2017), 327--360.
\bibitem{W1} Wikipedia, {\it Invariant subspace problem,} recalled 2023-09-09.
\bibitem{W2} Wikipedia, {\it Algebraic Riccati equation,} recalled 2023-09-12.
\end{thebibliography}

\begin{thebibliography}{9999}
\setcounter{enumiv}{28}
\bibitem{BG} K.\ Bongartz, P.\ Gabriel, {\it Covering Spaces in Representation-Theory,}
    Invent.\ Math.\ {\bf 65} (1981/82), pp 331--378.
\bibitem{D} P.\ Dowbor, {\it Properties of G-atoms and full
    Galois covering reduction to stabilizers,} Colloq.\ Math.\ {\bf 83} (2000), 231--265.
\bibitem{DS} P.\ Dowbor, A.\ Skowro\'nski, {\it Galois coverings of
  representation-infinite algebras}, Comment.\ Math.\ Helvetici {\bf 62} (1987), 311--337.
\bibitem{G} P.\ Gabriel, {\it The universal cover of a representation-finite algebra,}
  in: Lecture Notes in Math., vol.\ {\bf 903,} Springer-Verlag, Berlin/New York (1981), 68--105.
\bibitem{GG} R.\ Gordon, E.\ L.\ Green, 
  {\it Representation theory of graded Artin algebras,} J.\ Algebra
  {\bf 76} (1982), 138--152.
\bibitem{Rm} Chr.\ Riedtmann, {\it Algebren, Darstellungsköcher, \"Uberlagerungen und zur\"uck,}
  Comm.\ Math.\ Helv.\ {\bf 55} (1980), 199--224.
\end{thebibliography}
\end{document}